\documentclass[oneside,ngerman,american,reqno]{amsart}
\usepackage[T1]{fontenc}
\usepackage[utf8]{inputenc}
\usepackage{geometry}
\usepackage{color}
\usepackage{babel}
\usepackage{mathrsfs}
\usepackage{mathtools}
\usepackage{enumitem}
\usepackage{amstext}
\usepackage{amsthm}
\usepackage{amssymb}
\usepackage[unicode=true,pdfusetitle,
 bookmarks=true,bookmarksnumbered=false,bookmarksopen=false,
 breaklinks=false,pdfborder={0 0 1},backref=false,colorlinks=false]
 {hyperref}

\makeatletter
\numberwithin{equation}{section}
\theoremstyle{plain}
\newtheorem{thm}{\protect\theoremname}[section]
  \theoremstyle{plain}
  \newtheorem{cor}[thm]{\protect\corollaryname}
  \theoremstyle{remark}
  \newtheorem*{rem*}{\protect\remarkname}
  \theoremstyle{definition}
  \newtheorem{defn}[thm]{\protect\definitionname}
  \theoremstyle{remark}
  \newtheorem{rem}[thm]{\protect\remarkname}
  \theoremstyle{plain}
  \newtheorem{lem}[thm]{\protect\lemmaname}
  \theoremstyle{definition}
  \newtheorem{example}[thm]{\protect\examplename}
 \newlist{casenv}{enumerate}{4}
 \setlist[casenv]{leftmargin=*,align=left,widest={iiii}}
 \setlist[casenv,1]{label={{\itshape\ \casename} \arabic*.},ref=\arabic*}
 \setlist[casenv,2]{label={{\itshape\ \casename} \roman*.},ref=\roman*}
 \setlist[casenv,3]{label={{\itshape\ \casename\ \alph*.}},ref=\alph*}
 \setlist[casenv,4]{label={{\itshape\ \casename} \arabic*.},ref=\arabic*}
  \theoremstyle{plain}
  \newtheorem{assumption}[thm]{\protect\assumptionname}

\usepackage[T1]{fontenc}
\usepackage{dsfont}

\let\emptyset\varnothing

\makeatletter
\@namedef{subjclassname@1991}{2010 Mathematics Subject Classification}
\@namedef{subjclassname@2000}{2010 Mathematics Subject Classification}
\@namedef{subjclassname@2010}{2010 Mathematics Subject Classification}
\makeatother

\theoremstyle{remark}
\newtheorem{remx}[thm]{Remark}
\renewenvironment{rem}
  {\pushQED{\qed}\remx}
  {\popQED\endremx}

\theoremstyle{remark}
\newtheorem*{remxx}{Remark}
\renewenvironment{rem*}
  {\pushQED{\qed}\remxx}
  {\popQED\endremxx}

\theoremstyle{definition}
\newtheorem{examplex}[thm]{Example}
\renewenvironment{example}
  {\pushQED{\qed}\examplex}
  {\popQED\endexamplex}

\theoremstyle{definition}
\newtheorem{defnx}[thm]{Definition}
\renewenvironment{defn}
  {\pushQED{\qed}\defnx}
  {\popQED\enddefnx}

\theoremstyle{plain}
\newtheorem{lemx}[thm]{Lemma}
\renewenvironment{lem}
  {\pushQED{\qed}\lemx}
  {\popQED\endlemx}

\theoremstyle{plain}
\newtheorem{corx}[thm]{Corollary}
\renewenvironment{cor}
  {\pushQED{\qed}\corx}
  {\popQED\endcorx}

\theoremstyle{definition}
\newtheorem{thmx}[thm]{Theorem}
\renewenvironment{thm}
  {\pushQED{\qed}\thmx}
  {\popQED\endthmx}

\theoremstyle{definition}
\newtheorem{innercustomthm}{Theorem}
\newenvironment{customthm}[1]
  {\pushQED{\qed}\renewcommand\theinnercustomthm{#1}\innercustomthm}
  {\popQED\endinnercustomthm}

\theoremstyle{plain}
\newtheorem{assumptionx}[thm]{Assumption}
\renewenvironment{assumption}
  {\pushQED{\qed}\assumptionx}
  {\popQED\endassumptionx}

\DeclareRobustCommand{\SkipTocEntry}[5]{}

\makeatletter
\newcommand\manuallabel[1]{\def\@currentlabel{#1}}
\makeatletter

\usepackage[ddmmyyyy,hhmmss]{datetime}
\newdateformat{daymonthyeardate}{%
  \THEDAY.\THEMONTH.\THEYEAR}
\usepackage{fancyhdr}
\makeatletter

\newsavebox{\@linebox}
\savebox{\@linebox}[3em][t]{\parbox[t]{3em}{%
  \@tempcnta\@ne\relax
  \loop{\underline{\scriptsize\the\@tempcnta}}\\
    \advance\@tempcnta by \@ne\ifnum\@tempcnta<55\repeat}}

\makeatother

\let\origsection\section
\renewcommand{\section}[1]{\origsection{#1}\sectionmark{#1}} 
\renewcommand\sectionmark[1]{\markleft{\thesection: #1}}
\makeatother

  \addto\captionsamerican{\renewcommand{\assumptionname}{Assumption}}
  \addto\captionsamerican{\renewcommand{\corollaryname}{Corollary}}
  \addto\captionsamerican{\renewcommand{\definitionname}{Definition}}
  \addto\captionsamerican{\renewcommand{\examplename}{Example}}
  \addto\captionsamerican{\renewcommand{\lemmaname}{Lemma}}
  \addto\captionsamerican{\renewcommand{\remarkname}{Remark}}
  \addto\captionsamerican{\renewcommand{\theoremname}{Theorem}}
  \addto\captionsngerman{\renewcommand{\assumptionname}{Annahme}}
  \addto\captionsngerman{\renewcommand{\corollaryname}{Korollar}}
  \addto\captionsngerman{\renewcommand{\definitionname}{Definition}}
  \addto\captionsngerman{\renewcommand{\examplename}{Beispiel}}
  \addto\captionsngerman{\renewcommand{\lemmaname}{Lemma}}
  \addto\captionsngerman{\renewcommand{\remarkname}{Bemerkung}}
  \addto\captionsngerman{\renewcommand{\theoremname}{Theorem}}
  \providecommand{\assumptionname}{Assumption}
  \providecommand{\corollaryname}{Corollary}
  \providecommand{\definitionname}{Definition}
  \providecommand{\examplename}{Example}
  \providecommand{\lemmaname}{Lemma}
  \providecommand{\remarkname}{Remark}
 \addto\captionsamerican{\renewcommand{\casename}{Case}}
 \addto\captionsngerman{\renewcommand{\casename}{Fall}}
 \providecommand{\casename}{Case}
\providecommand{\theoremname}{Theorem}

\begin{document}

\title{Embeddings of Decomposition Spaces}

\author{Felix Voigtlaender}


\selectlanguage{ngerman}%

\address{Lehrstuhl A für Mathematik\foreignlanguage{american}{}\\
\foreignlanguage{american}{RWTH Aachen University}\\
\foreignlanguage{american}{52056 Aachen}\\
\foreignlanguage{american}{Germany}}
\selectlanguage{american}%

\email{felix@voigtlaender.xyz}
\selectlanguage{ngerman}%

\curraddr{Katholische Universität Eichstätt–Ingolstadt\\
Lehrstuhl Wissenschaftliches Rechnen\\
Ostenstraße 26\\
85072 Eichstätt\foreignlanguage{american}{}\\
\foreignlanguage{american}{Germany}}

\maketitle
\selectlanguage{american}%
\global\long\def\vertiii#1{{\left\vert \kern-0.25ex  \left\vert \kern-0.25ex  \left\vert #1\right\vert \kern-0.25ex  \right\vert \kern-0.25ex  \right\vert }}
\global\long\def\essup{\esssup}
\global\long\def\with{\,:\,}
\global\long\def\DecompSp#1#2#3{{\mathcal{D}\left({#1},L^{#2},{#3}\right)}}
\global\long\def\FourierDecompSp#1#2#3{{\mathcal{D}_{\mathcal{F}}\left({#1},L^{#2},{#3}\right)}}
\global\long\def\BAPUFourierDecompSp#1#2#3#4{{\mathcal{D}_{\mathcal{F},{#4}}\left({#1},L^{#2},{#3}\right)}}
\global\long\def\R{\mathbb{R}}
\global\long\def\Compl{\mathbb{C}}
\global\long\def\N{\mathbb{N}}
\global\long\def\Z{\mathbb{Z}}
\global\long\def\CalQ{\mathcal{Q}}
\global\long\def\CalP{\mathcal{P}}
\global\long\def\CalB{\mathcal{B}}
\global\long\def\CalR{\mathcal{R}}
\global\long\def\CalS{\mathcal{S}}
\global\long\def\CalD{\mathcal{D}}
\global\long\def\BesovInhom#1#2#3{B_{#3}^{#1,#2}}
\global\long\def\BesovHom#1#2#3{\dot{B}_{#3}^{#1,#2}}
\global\long\def\dimension{d}
\global\long\def\ModSpace#1#2#3{M_{#3}^{#1,#2}}
\global\long\def\AlphaModSpace#1#2#3#4{M_{#3,#4}^{#1,#2}}
\global\long\def\ShearletSmoothness#1#2#3{S_{#3}^{#1,#2}}
\global\long\def\GL{\mathrm{GL}}
\global\long\def\CalO{\mathcal{O}}
\global\long\def\Coorbit{\mathrm{Co}}
\global\long\def\d{\operatorname{d}}
\global\long\def\TestFunctionSpace#1{C_{c}^{\infty}\left(#1\right)}
\global\long\def\DistributionSpace#1{\mathcal{D}'\left(#1\right)}
\global\long\def\SpaceTestFunctions#1{Z\left(#1\right)}
\global\long\def\SpaceReservoir#1{Z'\left(#1\right)}
\global\long\def\Schwartz{\mathcal{S}}
\global\long\def\Fourier{\mathcal{F}}
\global\long\def\supp{\operatorname{supp}}
\global\long\def\dist{\operatorname{dist}}
\global\long\def\Xhookrightarrow#1{\xhookrightarrow{#1}}
\global\long\def\Xmapsto#1{\xmapsto{#1}}
\global\long\def\Indicator{{\mathds{1}}}
\global\long\def\identity{\operatorname{id}}
\global\long\def\mybullet{\bullet}
\markleft{1: Introduction}\global\long\def\UpperExpo#1{#1^{\left(\uparrow\right)}}
\global\long\def\LowerExpo#1{#1^{\left(\downarrow\right)}}
\global\long\def\SignedUpperExpo#1{#1^{\left(\uparrow,\pm\right)}}

\begin{abstract}
A unifying framework for a wide class of smoothness spaces used in
harmonic analysis—Besov spaces, ($\alpha$)-modulation spaces, and
many wavelet-type coorbit spaces—is provided by the theory of \emph{decomposition
spaces}. In this paper we consider the following question: Given
two decomposition spaces, is there a continuous inclusion relation
(an \emph{embedding}) between the two?

A decomposition space $\DecompSp{\CalQ}pY$ is determined by three
ingredients: a covering $\CalQ=\left(Q_{i}\right)_{i\in I}$ of the
frequency domain, an integrability exponent $p\in\left(0,\infty\right]$
and a sequence space $Y\subset\Compl^{I}$ on the index set $I$.
Given these, the decomposition space norm of a distribution $g$
is calculated by decomposing the frequency content of $g$ into different
parts according to the covering $\CalQ$. Each ``frequency localized
piece'' is measured in the $L^{p}$ norm, and the contributions of
the individual parts are globally aggregated using the sequence space
$Y$. More formally, we have ${\left\Vert g\right\Vert _{\DecompSp{\CalQ}pY}=\left\Vert \left(\left\Vert \Fourier^{-1}\left(\varphi_{i}\cdot\widehat{g}\right)\right\Vert _{L^{p}}\right)_{i\in I}\right\Vert _{Y}}$,
where $\left(\varphi_{i}\right)_{i\in I}$ is a suitable partition
of unity subordinate to $\CalQ$.

Up to now, the theory of decomposition spaces was mainly used to simplify
the introduction of new spaces, for instance to justify the independence
of the resulting space from the chosen partition of unity $\left(\varphi_{i}\right)_{i\in I}$.
We will show that identifying a given function space $X$ as a decomposition
space offers much greater advantages: Once this is done, assertions
about the (non)-existence of embeddings $X\hookrightarrow V$ or $V\hookrightarrow X$,
for other decomposition spaces $V$, come almost for free.

More precisely, we establish readily verifiable criteria which ensure
the existence of an embedding $\DecompSp{\CalQ}{p_{1}}Y\hookrightarrow\DecompSp{\CalP}{p_{2}}Z$,
mostly concentrating on the case where $Y=\ell_{w}^{q_{1}}\left(I\right)$
and $Z=\ell_{v}^{q_{2}}\left(J\right)$. Since in this case, the two
decomposition spaces depend only on the quantities $\CalQ,\CalP$,
$w,v$ and $p_{1},p_{2},q_{1},q_{2}$, it should be possible to decide
the existence of the embedding only in terms of these ingredients.
It is not at all clear, however, precisely \emph{which} conditions
on these quantities ensure or prevent the existence of such an embedding.
We will see—under suitable assumptions on $\CalQ,\CalP$—that the
relevant sufficient conditions are $p_{1}\leq p_{2}$ and finiteness
of a nested norm of the form
\[
\left\Vert \left(\left\Vert \left(\alpha_{i}\beta_{j}\cdot v_{j}/w_{i}\right)_{i\in I_{j}}\right\Vert _{\ell^{t}}\right)_{j\in J}\right\Vert _{\ell^{s}}\,,\quad\text{with}\quad I_{j}=\left\{ i\in I\with Q_{i}\cap P_{j}\neq\emptyset\right\} \quad\text{for }j\in J\,.
\]
The sets $I_{j}$ only depend on the coverings $\CalQ=\left(Q_{i}\right)_{i\in I}$
and $\CalP=\left(P_{j}\right)_{j\in J}$. Further, the exponents $t,s$
and the weights $\alpha,\beta$ also only depend on the quantities
used to define the decomposition spaces.  In a nutshell, although
knowledge of Fourier analysis is required to define and understand
decomposition spaces, \emph{no such knowledge is required if one just
wants to apply the embedding results presented in this article}. Instead,
one only has to study the geometric properties of the involved coverings,
so that one can decide the finiteness of certain sequence space norms
defined in terms of the coverings.

Finally, we show—under suitable assumptions on the two coverings $\CalQ,\CalP$—that
the developed criteria are \emph{sharp}. In fact, for almost arbitrary
coverings and certain ranges of $p_{1},p_{2}$, our criteria yield
a \emph{complete characterization} of the existence of an embedding.
The same holds for \emph{arbitrary} values of $p_{1},p_{2}$ under
more stringent assumptions on the two coverings $\CalQ,\CalP$.

As a further application of our proof techniques, we show a \emph{rigidity
result}, namely that—for $\left(p_{1},q_{1}\right)\neq\left(2,2\right)$—two
decomposition spaces $\DecompSp{\CalQ}{p_{1}}{\ell_{w}^{q_{1}}}$
and $\DecompSp{\CalP}{p_{2}}{\ell_{v}^{q_{2}}}$ can only \emph{coincide}
if their ``ingredients'' are equivalent, that is, if $p_{1}=p_{2}$
and $q_{1}=q_{2}$ and if the coverings $\CalQ,\CalP$ and the weights
$w,v$ are equivalent in a suitable sense.

We illustrate the power and convenience of the resulting embedding
theory by applications to $\alpha$-modulation spaces and Besov spaces.
We will see that all known results about embeddings between these
types of spaces are special cases of our general approach; in some
cases, our criteria even show that the previously known embedding
results can be improved considerably. Due to the length of the present
paper, we postpone further applications—for instance to shearlet smoothness
spaces and wavelet-type coorbit spaces, including shearlet coorbit
spaces—to a later contribution.

Finally, it is worth pointing out that our criteria allow for $p_{1}\neq p_{2}$,
so that one can ``trade smoothness for integrability'' in a certain
sense, similar to the classical Sobolev embedding theorems. For $\CalQ\neq\CalP$,
sharp embeddings of this generality have previously not been known,
even for the comparatively well-studied special case of $\alpha$-modulation
spaces.
\end{abstract}

\newpage{}

\label{toc}\tableofcontents{}

\section{Introduction}

Our aim in this paper is to establish a general and easy-to-use framework
for proving and disproving the existence of embeddings between decomposition
spaces.

The structure of this introduction is as follows: First, we illustrate
the relevance of a general framework for embeddings between decomposition
spaces by discussing related work. As we will see, the existing literature
shows that there is quite some interest in embeddings between different
special classes of decomposition spaces. All of these results can
be obtained with ease using the framework presented in this paper.
Furthermore, recent work in the theory of (wavelet type) coorbit spaces
shows that certain decomposition spaces defined using ``exotic, non-classical''
coverings are of interest. This further increases the need for a general,
comprehensive framework for such embeddings.

Next, we fix the notations used in the remainder of the paper. Most
of these notational conventions are standard and are only given here
to fix for instance the exact normalization of the Fourier transform
which we use. There are, however, some non-standard notations which
are crucial in the remainder of the paper, namely the \textbf{upper/lower
(signed) conjugate exponents} $\UpperExpo p,\LowerExpo p$, and $\SignedUpperExpo p$.

Having introduced our notation, we give a quick overview of the results
obtained in this paper, including a brief overview over the most important
notions that are needed to formulate these results.

Finally, we close the introduction with an explanation of the general
structure of the paper.

\subsection{Motivation and related work}

\label{subsec:IntroductionRelatedWork}Let us start at the beginning:
Decomposition spaces were first introduced—in a very general setting—by
Feichtinger and Gröbner in \cite{DecompositionSpaces1,DecompositionSpaces2},
with the aim of giving a unifying framework for \emph{Besov spaces},
\emph{modulation spaces} and certain \emph{Wiener amalgam spaces}.
More precisely, both the scale of Besov spaces $B_{s}^{p,q}\left(\R^{d}\right)$
(including the \emph{anisotropic} Besov spaces \cite{BownikAnisotropicBesovSpaces})
and the scale of modulation spaces $M_{s}^{p,q}\left(\R^{d}\right)$
arise as special decomposition spaces, and thus are amenable to the
theory developed in this article. Furthermore, the modulation spaces
arise as the inverse image under the Fourier transform of certain
Wiener amalgam spaces (originally introduced in \cite{FeichtingerWienerSpaces,FeichtingerWienerInterpolation},
where they are called \emph{Wiener type spaces}). Precisely, $M_{s}^{p,q}\left(\R^{d}\right)=\Fourier^{-1}\left[W\left(\smash{\Fourier L^{p}},\ell_{w_{s}}^{q}\right)\right]$
for a suitable weight $w_{s}$. We emphasize, however, that not all
decomposition spaces can be obtained in this way from Wiener amalgam
spaces.

With the general framework of decomposition spaces established, Gröbner
introduced in his PhD thesis \cite{GroebnerAlphaModulationSpaces}
the so-called \emph{$\alpha$-modulation spaces}  $\AlphaModSpace pqs{\alpha}\left(\R^{\dimension}\right)$
for different values of $\alpha\in\left[0,1\right]$. These spaces
satisfy (up to canonical identifications)
\[
\AlphaModSpace pqs{\alpha}\left(\R^{\dimension}\right)=\DecompSp{\smash{\CalQ^{\left(\alpha\right)}}}p{\ell_{u^{\left(s,\alpha\right)}}^{q}}
\]
for a so-called \emph{$\alpha$-covering} $\CalQ^{\left(\alpha\right)}$
of $\R^{\dimension}$ and a suitable weight $u^{\left(s,\alpha\right)}$.
The idea underlying the $\alpha$-coverings $\CalQ^{\left(\alpha\right)}$
is that they ``geometrically interpolate'' the dyadic covering used
to construct (inhomogeneous) Besov spaces and the uniform covering
with which modulation spaces are defined. Thus, $\alpha$-modulation
spaces can be considered as intermediate spaces between modulation-
and Besov spaces. Note, however, that they can \emph{not} be obtained
by complex interpolation between those spaces, as shown in \cite{AlphaModulationNotInterpolation}.

In his thesis, Gröbner obtained certain sufficient and certain necessary
conditions for the existence of the embedding
\begin{equation}
\AlphaModSpace{p_{1}}{q_{1}}{s_{1}}{\alpha_{1}}\left(\R^{\dimension}\right)\hookrightarrow\AlphaModSpace{p_{2}}{q_{2}}{s_{2}}{\alpha_{2}}\left(\R^{\dimension}\right),\label{eq:IntroductionAlphaModulationEmbedding}
\end{equation}
but did not achieve a complete characterization. More recently, Sugimoto
and Tomita \cite{ModulationBesovEmbedding} characterized the embeddings
between modulation spaces and Besov spaces, which amounts to the above
embedding for $\left(\alpha_{1},\alpha_{2}\right)\in\left\{ \left(0,1\right),\left(1,0\right)\right\} $.
For general $\alpha_{1},\alpha_{2}$, the work of Gröbner was continued
by Toft and Wahlberg \cite{ToftWahlbergAlphaModulationEmbeddings},
shortly before the existence of the embedding~(\ref{eq:IntroductionAlphaModulationEmbedding})
was characterized completely—at least for $\left(p_{1},q_{1}\right)=\left(p_{2},q_{2}\right)$—by
Han and Wang \cite{HanWangAlphaModulationEmbeddings}. Some of their
techniques were used in my PhD thesis \cite{VoigtlaenderPhDThesis}
on which large parts of the present paper are based. We will see in
the last part of the present paper that our results greatly generalize
those of Han and Wang, since they characterize the embedding~(\ref{eq:IntroductionAlphaModulationEmbedding})
\emph{completely}, also for $\left(p_{1},q_{1}\right)\neq\left(p_{2},q_{2}\right)$.
We mention that the same characterization of the embeddings between
$\alpha$-modulation spaces that we develop here has been obtained
independently by Guo et al. (see \cite{GuoAlphaModulationEmbeddingCharacterization}),
but without developing a general framework for the embeddings of decomposition
spaces; in fact, their preprint appeared on the arXiv only a few days
after the first version of the present paper.

Borup and Nielsen also contributed heavily to the theory of decomposition
spaces. In \cite{BorupNielsenDecomposition}, they introduced a more
restrictive class of decomposition spaces than those originally defined
by Feichtinger and Gröbner. It is (essentially) this class of decomposition
spaces that we will consider in the present paper. Furthermore, Borup
and Nielsen introduced the class of structured coverings and gave
a general construction of Banach frames and atomic decompositions
for decomposition spaces \cite{BorupNielsenDecomposition}, which
they first formulated (see \cite{BorupNielsenAlphaModulationSpaces})
in the context of $\alpha$-modulation spaces; see also \cite{NielsenAlphaModulationONB}
for a construction of an orthonormal basis compatible to $\alpha$-modulation
spaces. Related constructions of Banach frames and atomic decompositions
for $\alpha$-modulation spaces were obtained by Dahlke et al.~\cite{QuotientCoorbitTheoryAndAlphaModulationSpaces},
by Fornasier~\cite{FornasierFramesForAlphaModulation}, by Speckbacher
et al.~\cite{SpeckbacherAlphaModulationTransform}, and by the present
author~\cite{StructuredBanachFrames1}. These Banach frames for the
$\alpha$-modulation spaces are obtained by dilating, modulating,
and translating a fixed prototype function, where the parameter $\alpha$
determines the dependence of the dilation factor on the frequency
parameter; see \cite[Section 3]{FornasierFramesForAlphaModulation}
for more details. Because of this special structure, Fornasier calls
these frames ``$\alpha$-Gabor-wavelet frames'' \cite{FornasierFramesForAlphaModulation}.
In the context of $\alpha$-modulation spaces, I would also like to
mention the recent paper \cite{KatoAlphaModulationSobolev} in which
Kato characterizes the embeddings between $\alpha$-modulation spaces
and Sobolev spaces.

Additionally, Borup and Nielsen studied boundedness properties of
certain pseudo-differential operators on $\alpha$-modulation spaces;
see \cite{BorupNielsenPsiDOsOnMultivariateAlphaModulationSpaces}.
In this context, it is also worth mentioning the recent paper \cite{AlphaModulationFourierMultiplier},
in which boundedness properties of certain Fourier multipliers on
$\alpha$-modulation spaces are studied.

A further application of decomposition spaces is the paper \cite{Labate_et_al_Shearlet}
in which Labate et al.\@ define so-called \emph{shearlet smoothness
spaces} and study (among other things) the existence of embeddings
between these spaces and (inhomogeneous) Besov spaces. Note, however,
that their definition of Besov spaces does not coincide with the usual
one, so that their results need to be taken with a grain of salt.
Using the results in the present paper, the existence of the embeddings
studied in \cite{Labate_et_al_Shearlet} can be \emph{characterized
completely} (see \cite[Theorem 6.4.3]{VoigtlaenderPhDThesis}), whereas
Labate et al.\@ only give sufficient criteria. The shearlet smoothness
spaces introduced in \cite{Labate_et_al_Shearlet} have been further
studied in \cite{AlphaShearletSparsity}, where it was shown that
membership of a function/distribution $f$ in the shearlet smoothness
space is equivalent to sparsity of $f$ with respect to a (sufficiently
nice, compactly supported) shearlet frame. A related construction
of a bandlimited shearlet frame with ``nice'' dual has been given
in \cite{GrohsBandlimitedShearletsWithNiceDuals}.

Another result illustrating the richness of decomposition spaces is
the work \cite{FuehrVoigtlaenderCoorbitSpacesAsDecompositionSpaces}
of Führ and the present author, in which it is shown for a large class
of (generalized) \emph{wavelet-type coorbit spaces} that each space
from this class of coorbit spaces is canonically isomorphic to a certain
decomposition space. Precisely, if $H\subset\GL\left(\R^{\dimension}\right)$
is a so-called \emph{admissible dilation group}, we have (up to canonical
identifications) that
\begin{equation}
{\rm Co}\left(L_{m}^{p,q}\left(\R^{\dimension}\rtimes H\right)\right)=\DecompSp{\CalQ_{H}}p{\ell_{u}^{q}},\label{eq:IntroductionCoorbitAsDecomposition}
\end{equation}
where $\CalQ_{H}$ is a so-called \emph{induced covering} of the \emph{dual
orbit} $\CalO=H^{T}\xi_{0}\subset\R^{\dimension}$ of $H$. Here,
$\xi_{0}\in\R^{\dimension}$ is chosen such that $\CalO$ is open
and of full measure; existence of such $\xi_{0}$ is part of the definition
of an admissible dilation group. Originally, this result was only
proven in the Banach setting $p,q\in\left[1,\infty\right]$, but in
my PhD thesis \cite{VoigtlaenderPhDThesis}, this was extended to
the full range $p,q\in\left(0,\infty\right]$. We refer to \cite{FeichtingerVoigtlaender}
for a survey relating the results in \cite{FuehrVoigtlaenderCoorbitSpacesAsDecompositionSpaces}
with the findings in the present paper. A related result—which characterizes
coorbit spaces in terms of \emph{phase-space} coverings instead of
\emph{frequency} space coverings—is obtained in \cite{RomeroCoorbitSpacesAndPhaseSpaceCovers};
see also \cite{RomeroDoerflerFramesAdaptedToCover} for a related
construction of frames. A further construction of frames based on
decomposition methods is given in \cite{FornasierQuasiOrthogonalDecompositions}.

Using the isomorphism (\ref{eq:IntroductionCoorbitAsDecomposition}),
many properties of the coorbit space ${\rm Co}\left(L_{m}^{p,q}\left(\R^{\dimension}\rtimes H\right)\right)$
can be deduced which are not obvious from the coorbit-space point
of view. In particular, combined with the embedding results in this
paper, one can obtain nontrivial embeddings between wavelet-type coorbit
spaces with respect to \emph{different} dilation groups, and also
embeddings between these coorbit spaces and classical function spaces
like Besov spaces. For some results in this direction, see \cite[Sections 6.3 and 6.5]{VoigtlaenderPhDThesis}.
Using the results from \cite{DecompositionIntoSobolev}, one can also
obtain embeddings into the classical Sobolev spaces $W^{k,q}\left(\R^{\dimension}\right)$;
see in particular \cite[Section 7]{DecompositionIntoSobolev} and
\cite{ShearletCoorbitIntoSobolev}. Due to the considerable length
of the present paper, we postpone a more detailed discussion of these
examples to a later contribution. We remark, however, that the results
obtained in \cite{VoigtlaenderPhDThesis,DecompositionIntoSobolev}
are very different from those obtained by Dahlke et al.\@ in their
investigation of shearlet coorbit spaces \cite{Dahlke_etal_sh_coorbit1,Dahlke_etal_sh_coorbit2,DahlkeShearletArbitraryDimension,DahlkeShearletCoorbitEmbeddingsInHigherDimensions,DahlkeToeplitzShearletTransform}:
For example in \cite[Theorem 4.7]{Dahlke_etal_sh_coorbit2}, the authors
consider embeddings of \emph{a subspace} $\mathcal{SCC}_{p,r}$ of
a shearlet coorbit space into \emph{a sum} of \emph{homogeneous} Besov
spaces $\BesovHom pp{\sigma_{1}}\left(\R^{2}\right)+\BesovHom pp{\sigma_{2}}\left(\R^{2}\right)$,
while in \cite{VoigtlaenderPhDThesis}, embeddings of the \emph{whole}
shearlet coorbit space into a \emph{single inhomogeneous} Besov space
are considered.

Finally, we remark that the embedding results presented in this paper
can be used to decide whether a given decomposition space $\DecompSp{\CalQ}p{\ell_{w}^{q}}$
is invariant under \emph{dilation} by a matrix $A\in\GL\left(\R^{\dimension}\right)$
or by \emph{modulation} with $\xi\in\R^{\dimension}$. Indeed, if
we define $\CalP:=\left(A^{-T}\left(Q_{i}-\xi\right)\right)_{i\in I}$
where $\CalQ=\left(Q_{i}\right)_{i\in I}$, then a direct calculation
shows that $\left\Vert M_{\xi}\left(f\circ A\right)\right\Vert _{\DecompSp{\CalQ}p{\ell_{w}^{q}}}\asymp\left\Vert f\right\Vert _{\DecompSp{\CalP}p{\ell_{w}^{q}}}$,
where $M_{\xi}g\left(x\right)=e^{2\pi i\left\langle x,\xi\right\rangle }\,g\left(x\right)$
denotes the modulation of $g$, and where the implied constant depends
only on $A$ and $p$. In other words, the map $f\mapsto M_{\xi}\left(f\circ A\right)$
acts boundedly on $\DecompSp{\CalQ}p{\ell_{w}^{q}}$ if and only if
the embedding ${\DecompSp{\CalQ}p{\ell_{w}^{q}}\hookrightarrow\DecompSp{\CalP}p{\ell_{w}^{q}}}$
holds; the existence of this embedding can in many cases be decided
using the criteria developed in the present article. Since this article
is already quite long, we refrain from discussing this further, and
refer the interested reader to \cite[Section 6.5]{VoigtlaenderPhDThesis}.
We would like to mention, however, that this type of argument has
been used in \cite[Section 8]{WavePacketSpaces} to prove that the
so-called \emph{wave packet smoothness spaces} (a scale of spaces
which strictly includes the scale of $\alpha$-modulation spaces)
are invariant under dilation with arbitrary invertible matrices $A\in\GL\left(\R^{2}\right)$.
Furthermore, a combination of the sketched argument with the identification
of wavelet-type coorbit spaces as decomposition spaces (equation\ (\ref{eq:IntroductionCoorbitAsDecomposition}))
was used in \cite[Section 6.5]{VoigtlaenderPhDThesis} to determine
the group of matrices which act boundedly by dilation on a certain
class of shearlet-type coorbit spaces.

Overall, the discussed body of literature shows that decomposition
spaces appear in a variety of contexts, from $\alpha$-modulation
spaces over the characterization of signals which can be sparsely
represented using shearlets, to the theory of wavelet-type coorbit
spaces. The embedding theory developed in the present paper yields
novel results in each of the considered examples, for example when
comparing these spaces to more ``classical'' function spaces like
Besov spaces, or when studying properties like dilation invariance
of these spaces. Furthermore, the theory developed here greatly simplifies
the task of comparing all the resulting spaces with each other, for
instance to consider embeddings between $\alpha$-modulation spaces
and shearlet smoothness spaces.

\subsection{Notation}

\label{subsec:Notation}We let $\N:=\left\{ n\in\Z\with n\geq1\right\} $
and $\N_{0}:=\left\{ 0\right\} \cup\N$. For $n\in\N_{0}$, we write
$\underline{n}:=\left\{ 1,\dots,n\right\} $. In particular, $\underline{0}=\emptyset$
is the empty set. We use the notation $X\subset Y$ to denote that
$X$ is a subset of $Y$, where possibly $X=Y$.

For $a\in\R$, we denote the positive part of $a$ by $a_{+}:=\max\left\{ a,0\right\} $.
We write $\left|x\right|$ for the euclidean norm of $x\in\R^{\dimension}$.
The \emph{open} euclidean ball of radius $r>0$ around $x\in\R^{\dimension}$
is denoted by $B_{r}\left(x\right)$ and the corresponding closed
ball is denoted by $\overline{B_{r}}\left(x\right)$. 

For $A,B\subset\R^{\dimension}$, we write $A-B$ for the \textbf{Minkowski
difference} (or \textbf{algebraic difference})
\[
A-B:=\left\{ a-b\with a\in A,\,b\in B\right\} \subset\R^{\dimension}\,,
\]
which should not be confused with the set-theoretic difference $A\setminus B=\left\{ a\in A\with a\notin B\right\} $.

For a subset $M\subset V$ of a vector space $V$, we write $\left\langle M\right\rangle $
for the \textbf{span} of $M$. If $\left(x_{i}\right)_{i\in I}$ is
a family in $V$, we also write $\left\langle x_{i}\with i\in I\right\rangle $
for $\left\langle \left\{ x_{i}\with i\in I\right\} \right\rangle $.
If $W$ is a \textbf{subspace} of $V$, we write $W\leq V$.

A \textbf{quasi-norm} $\left\Vert \mybullet\right\Vert :X\to\left[0,\infty\right)$
on a vector space $X$ has to satisfy the same properties as a norm,
with the exception that the triangle inequality is replaced by the
\textbf{quasi-triangle inequality} 
\[
\left\Vert x+y\right\Vert \leq C\cdot\left(\left\Vert x\right\Vert +\left\Vert y\right\Vert \right)\qquad\forall\:x,y\in X,
\]
where the \textbf{triangle constant} $C\geq1$ of $X$ is independent
of $x,y\in X$. Cauchy sequences, convergence of sequences, completeness,
and the topology induced by the quasi-norm are defined as for normed
vector spaces. A complete quasi-normed space is called a \textbf{quasi-Banach
space}. The two main differences between normed vector spaces and
quasi-normed vector spaces are that a quasi-norm $\left\Vert \mybullet\right\Vert $
need not be continuous with respect to the topology that it generates,
and that the topological dual $X'$ of a quasi-normed space need not
separate the points of $X$: It can even happen that $X'=\left\{ 0\right\} $,
see for example \cite[Section 1.47]{RudinFA}. The most important
example of quasi-Banach spaces are the Lebesgue spaces $L^{p}\left(\mu\right)$
for $p\in\left(0,1\right)$.

We remark that the possible discontinuity of a quasi-norm $\left\Vert \mybullet\right\Vert $
is not as serious as it first appears: Indeed, the \emph{Aoki-Rolewicz
Theorem} (cf.\@ \cite[Chapter 2, Theorem 1.1]{ConstructiveApproximation})
shows that $X$ admits an equivalent \textbf{$r$-norm} $\vertiii{\mybullet}$
for some $r\in\left(0,1\right]$, which simply means that $\vertiii{x+y}^{r}\leq\vertiii x^{r}+\vertiii y^{r}$
for all $x,y\in X$. This easily implies that $d\left(x,y\right):=\vertiii{x-y}^{r}$
defines a metric on $X$ which induces the same topology as $\left\Vert \mybullet\right\Vert $.
Furthermore, the metric $d$ is translation invariant (in the sense
$d\left(x+z,y+z\right)=d\left(x,y\right)$), and if $\left(X,\left\Vert \mybullet\right\Vert \right)$
is complete, then so is $\left(X,d\right)$. Consequently, in the
terminology of \cite[Section 1.8]{RudinFA}, each quasi-Banach space
is an $F$-space. Being an $F$-space is sufficient for most classical
theorems from functional analysis to hold. In particular, the \textbf{closed
graph theorem} remains valid for $F$-spaces (see for instance \cite[Theorem 2.15]{RudinFA}),
and thus for quasi-Banach spaces.

Given a linear operator $T:X\to Y$ acting on quasi-normed spaces
$\left(X,\left\Vert \mybullet\right\Vert _{X}\right)$ and $\left(Y,\left\Vert \mybullet\right\Vert _{Y}\right)$,
we denote the \textbf{operator norm} of $T$ (which really is a quasi-norm,
unless $Y$ is a normed vector space) by
\[
\vertiii T:=\vertiii T_{X\to Y}:=\sup_{x\in X,\left\Vert x\right\Vert _{X}\leq1}\left\Vert Tx\right\Vert _{Y}\in\left[0,\infty\right]\,.
\]

For the \textbf{interior} of a subset $M$ of a topological space
$X$, we write $M^{\circ}$. The \textbf{closure} of $M$ is denoted
by $\overline{M}$. Given an open subset $U\subset X$, we write $C_{c}\left(U\right)$
for the space of all continuous functions $f:X\to\Compl$ which have
compact support $\supp f\subset U$, where $\supp f:=\overline{\left\{ x\in X\with f\left(x\right)\neq0\right\} }$.
We denote by $C_{0}\left(\R^{\dimension}\right)$ the space of all
continuous functions $f:\R^{\dimension}\to\Compl$ that \textbf{vanish
at infinity}, i.e., that satisfy $f\left(x\right)\to0$ as $\left|x\right|\to\infty$.

For a function $f:G\to S$, for some group $G$ and some set $S$,
we write
\[
L_{x}f:\:G\to S,\,y\mapsto f\left(\smash{x^{-1}}y\right)\qquad\text{ and }\qquad f^{\vee}:\:G\to S,\,y\mapsto f\left(\smash{y^{-1}}\right)
\]
for each $x\in G$. For $G=\R^{\dimension}$, we occasionally also
write $T_{x}:=L_{x}$.

For a matrix $h\in\R^{k\times\dimension}$, we write $h^{T}$ for
the \textbf{transpose} of $h$. Furthermore, for $f:\R^{\dimension}\to\Compl$
and $\xi\in\R^{\dimension}$, we define the \textbf{modulation} of
$f$ by $\xi$ as 
\[
M_{\xi}f:\R^{\dimension}\to\Compl,x\mapsto e^{2\pi i\left\langle x,\xi\right\rangle }\cdot f\left(x\right),
\]
where $\left\langle x,y\right\rangle =\sum_{j=1}^{\dimension}x_{j}y_{j}$
denotes the \textbf{standard scalar product} on $\R^{\dimension}$.
For $f:\R^{\dimension}\to\Compl$ and $h\in\GL\left(\R^{\dimension}\right)$,
we denote the \textbf{dilation} of $f$ by $h^{T}$ by
\[
D_{h}f:=f\circ h^{T}\,,\quad\text{where as usual}\quad\left(f\circ g\right)\left(x\right)=f\left(g\left(x\right)\right).
\]
In case of $h=a\cdot\identity$ for some $a\in\R^{\ast}=\R\setminus\left\{ 0\right\} $,
we also write $D_{a}:=D_{h}$.

We will frequently use \textbf{multiindex notation}: For $\alpha\in\N_{0}^{d}$
and a (sufficiently smooth) function $f:U\subset\R^{\dimension}\to\Compl$,
we let $\partial^{\alpha}f:=\frac{\partial^{\alpha_{1}}}{\partial x_{1}^{\alpha_{1}}}\cdots\frac{\partial^{\alpha_{\dimension}}}{\partial x_{\dimension}^{\alpha_{\dimension}}}f$,
and $x^{\alpha}:=x_{1}^{\alpha_{1}}\cdots x_{\dimension}^{\alpha_{\dimension}}$
for $x\in\R^{\dimension}$.

For the \textbf{Fourier transform}, we use the normalization
\[
\left(\Fourier f\right)\left(\xi\right):=\widehat{f}\left(\xi\right):=\int_{\R^{\dimension}}f\left(x\right)\cdot e^{-2\pi i\left\langle x,\xi\right\rangle }\,\d x\qquad\forall\:\xi\in\R^{\dimension}
\]
for $f\in L^{1}\left(\R^{\dimension}\right)$. As is well known, $\Fourier$
uniquely extends to a unitary automorphism of $L^{2}\left(\R^{\dimension}\right)$,
where the inverse is the unique extension of the operator $\Fourier^{-1}$,
given on $L^{1}\left(\R^{\dimension}\right)$ by
\[
\left(\Fourier^{-1}g\right)\left(x\right):=\int_{\R^{\dimension}}g\left(\xi\right)\cdot e^{2\pi i\left\langle x,\xi\right\rangle }\,\d\xi\qquad\forall\:x\in\R^{\dimension}.
\]

We will make heavy use of the formalism of \textbf{(tempered) distributions},
which we briefly recall here: For an open set $U\subset\R^{\dimension}$,
we denote by $\TestFunctionSpace U$ the space of all $C^{\infty}$-functions
$f:\R^{\dimension}\to\Compl$ with $\supp f\subset U$. With the topology
described in \cite[Sections 6.3--6.6]{RudinFA}, $\TestFunctionSpace U$
becomes a complete locally convex space, whose (topological) dual
space we denote by $\DistributionSpace U$, the space of \textbf{distributions}
on $U$. The most convenient description of $\DistributionSpace U$
is given by \cite[Theorem 6.8]{RudinFA}: A linear functional $\phi:\TestFunctionSpace U\to\Compl$
belongs to $\DistributionSpace U$ iff for every compact set $K\subset U$,
there are $N\in\N$ and $C>0$ with $\left|\phi\left(g\right)\right|\leq C\cdot\max_{\left|\alpha\right|\leq N}\left\Vert \partial^{\alpha}g\right\Vert _{L^{\infty}}$
for all $g\in\TestFunctionSpace U$ with $\supp g\subset K$. We will
often write $\left\langle \phi,g\right\rangle _{\CalD'}:=\left\langle \phi,g\right\rangle :=\phi\left(g\right)$
for the \emph{bilinear} pairing between $\DistributionSpace U$ and
$\TestFunctionSpace U$. We equip $\DistributionSpace U$ with the
weak-$\ast$-topology, so that a net $\left(\phi_{\alpha}\right)_{\alpha}$
in $\DistributionSpace U$ satisfies $\phi_{\alpha}\to\phi$ iff $\left\langle \phi_{\alpha},g\right\rangle _{\CalD'}\to\left\langle \phi,g\right\rangle _{\CalD'}$
for all $g\in\TestFunctionSpace U$.

We will also need the notion of the \textbf{support of a distribution}
$\phi\in\DistributionSpace U$: We say that $\phi$ vanishes on an
open set $V\subset U$ if $\left\langle \phi,g\right\rangle _{\CalD'}=0$
for all $g\in\TestFunctionSpace V$. With this, we define 
\[
\supp\phi:=U\setminus\bigcup\left\{ V\subset U\with V\text{ open and }\phi\text{ vanishes on }V\right\} .
\]
It is well-known (see for instance  \cite[Proposition 9.2]{FollandRA})
that $\phi$ vanishes on $U\setminus\supp\phi$.

Next, for $f\in C^{\infty}\left(\R^{\dimension};\Compl\right)$, $N\in\N_{0}$
and $\alpha\in\N_{0}^{\dimension}$, write 
\[
\left\Vert f\right\Vert _{N,\alpha}:=\sup_{x\in\R^{n}}\left(1+\left|x\right|\right)^{N}\left|\partial^{\alpha}f\left(x\right)\right|\in\left[0,\infty\right]\,.
\]
Then the space $\Schwartz\left(\R^{\dimension}\right)$ of \textbf{Schwartz
functions} is simply
\[
\Schwartz\left(\R^{\dimension}\right):=\left\{ f\in C^{\infty}\left(\R^{\dimension};\Compl\right)\with\forall\,N\in\N_{0},\alpha\in\N_{0}^{\dimension}\,:\,\left\Vert f\right\Vert _{N,\alpha}<\infty\right\} \,,
\]
which we equip with the topology generated by the family of seminorms
$\big(\left\Vert \bullet\right\Vert _{N,\alpha}\big)_{N\in\N_{0},\alpha\in\N_{0}^{\dimension}}$;
see also \cite[Section 8.1]{FollandRA}. We denote the topological
dual space of $\Schwartz\left(\R^{\dimension}\right)$ by $\Schwartz'\left(\R^{\dimension}\right)$,
the space of \textbf{tempered distributions}. As for ordinary distributions,
we write $\left\langle \phi,g\right\rangle _{\Schwartz'}:=\left\langle \phi,g\right\rangle :=\phi\left(g\right)$
for the bilinear pairing between $\Schwartz'\left(\R^{\dimension}\right)$
and $\Schwartz\left(\R^{\dimension}\right)$, and we equip $\Schwartz'\left(\R^{\dimension}\right)$
with the weak-$\ast$-topology. It is not hard to see that each tempered
distribution $\phi\in\Schwartz'\left(\R^{\dimension}\right)$ can
be restricted to a distribution $\phi|_{U}\in\DistributionSpace U$
for $U\subset\R^{\dimension}$ open, by setting $\left\langle \phi|_{U},g\right\rangle _{\CalD'}:=\left\langle \phi,g\right\rangle _{\Schwartz'}$
for $g\in\TestFunctionSpace U\subset\Schwartz\left(\R^{\dimension}\right)$.
We will sometimes abuse notation and not distinguish between $\phi$
and $\phi|_{U}$, especially in case of $U=\R^{\dimension}$.

Finally, we extend the Fourier transform as usual to $\Schwartz'\left(\R^{\dimension}\right)$:
Since $\Fourier:\Schwartz\left(\R^{\dimension}\right)\to\Schwartz\left(\R^{\dimension}\right)$
is continuous (see \cite[Corollary 8.23]{FollandRA}), we have $\Fourier\phi:=\widehat{\phi}:=\phi\circ\Fourier\in\Schwartz'\left(\R^{\dimension}\right)$
for every $\phi\in\Schwartz'\left(\R^{\dimension}\right)$. Put differently,
this definition simply means $\left\langle \Fourier\phi,\,g\right\rangle _{\Schwartz'}=\left\langle \phi,\Fourier g\right\rangle _{\Schwartz'}$.

The \textbf{Lebesgue measure} of a (measurable) subset $M\subset\R^{\dimension}$
will be denoted by $\lambda\left(M\right)$. The \textbf{cardinality}
of a set $M$ (which will either be a nonnegative integer or $\infty$)
is denoted by $\left|M\right|$.

If $X$ is a set (whose choice is fixed and clear from the context)
and if $M\subset X$, we write $\Indicator_{M}$ for the \textbf{characteristic
function} (or \textbf{indicator function}) 
\[
\Indicator_{M}:X\to\left\{ 0,1\right\} ,x\mapsto\begin{cases}
0, & \text{if }x\notin M,\\
1, & \text{if }x\in M.
\end{cases}
\]

In the remainder of the paper, the\textbf{ weighted $\ell^{q}$ spaces}
$\ell_{w}^{q}\left(I\right)$ will play a crucial role: If $I$ is
a set, if $q\in\left(0,\infty\right]$, and if $w=\left(w_{i}\right)_{i\in I}$
is a weight on $I$ (i.e., $w_{i}\in\left(0,\infty\right)$ for all
$i\in I$), then 
\[
\ell_{w}^{q}\left(I\right):=\left\{ x=\left(x_{i}\right)_{i\in I}\in\Compl^{I}\with\left(w_{i}\cdot x_{i}\right)_{i\in I}\in\ell^{q}\left(I\right)\right\} ,
\]
which we equip with the (quasi)-norm $\left\Vert \left(x_{i}\right)_{i\in I}\right\Vert _{\ell_{w}^{q}}:=\left\Vert \left(w_{i}\cdot x_{i}\right)_{i\in I}\right\Vert _{\ell^{q}}$.
We denote by $\ell_{0}\left(I\right)$ the space of all \textbf{finitely
supported (complex valued) sequences} over the index set $I$. Note
that if $I=\emptyset$, then $\ell_{0}\left(I\right)=\ell_{w}^{q}\left(I\right)=\left\{ 0\right\} $.

For $p\in\left[1,\infty\right]$, we let $p'\in\left[1,\infty\right]$
denote the usual \textbf{conjugate exponent} of $p$, i.e.\@ such
that $\frac{1}{p}+\frac{1}{p'}=1$. For $p\in\left(0,1\right)$, however,
we define the conjugate exponent of $p$ as $p':=\infty$. Note that
$\frac{1}{p}+\frac{1}{p'}=1$ is thus \emph{not} fulfilled for $p\in\left(0,1\right)$.

Furthermore, we define the \textbf{upper conjugate exponent} $\UpperExpo p$
of $p\in\left(0,\infty\right]$ by $\UpperExpo p:=\max\left\{ p,p'\right\} $,
whereas the \textbf{lower conjugate exponent} is given by $\LowerExpo p:=\min\left\{ p,p'\right\} $.
Finally, the \textbf{signed upper conjugate exponent} $\SignedUpperExpo p$
is defined by
\[
\SignedUpperExpo p:=\begin{cases}
p, & \text{if }p\geq2,\\
\frac{p}{p-1}, & \text{if }0<p<2\text{ and }p\neq1,\\
\infty, & \text{if }p=1,
\end{cases}
\]
which ensures that $\frac{1}{\UpperExpo p}=\min\left\{ p^{-1},1-p^{-1}\right\} $.

Our notation and nomenclature for these exponents is designed to be
mnemonic: $\UpperExpo p$ is the \emph{upper} conjugate exponent,
since we are taking the \emph{maximum} of $p$ and its conjugate $p'$.
For the same reason, we choose an \emph{upwards} pointing arrow as
the exponent of $\UpperExpo p$. Our conventions for $\LowerExpo p$
are explained similarly.

Finally, note that we have $\frac{1}{\UpperExpo p}=\min\left\{ \smash{\frac{1}{p}},\smash{\frac{1}{p'}}\right\} $,
with $\frac{1}{p'}=1-\frac{1}{p}$ for $p\in\left[1,\infty\right]$,
but $\frac{1}{p'}=0$ for $p\in\left(0,1\right)$. Thus, if we were
to define $p'$ such that $\frac{1}{p'}=1-\frac{1}{p}$ for all $p\in\left(0,\infty\right]$,
we would arrive at $\frac{1}{\UpperExpo p}=\min\left\{ p^{-1},1-p^{-1}\right\} $.
Note that this identity \emph{does} hold for $\SignedUpperExpo p$
instead of $\UpperExpo p$. Furthermore, note that $1/\SignedUpperExpo p=1-p^{-1}<0$
for $p\in\left(0,1\right)$, so that $\SignedUpperExpo p$ is signed—in
contrast to $\UpperExpo p,\LowerExpo p\in\left(0,\infty\right]$.
These observations explain the notation and nomenclature of the \emph{signed}
upper conjugate exponent $\SignedUpperExpo p$.

\addtocontents{toc}{\SkipTocEntry}

\subsection*{A comment on the writing style and the length of the paper}

The reader might object that the length of the present paper is excessive.
But this length is due to two main factors:
\begin{itemize}[leftmargin=0.6cm]
\item The paper has three different main goals which could in principle
each be handled in a separate paper. Namely,

\begin{itemize}
\item we establish \emph{sufficient} conditions for the existence of embeddings,
\item we establish \emph{necessary} conditions for the existence of embeddings,
\item we simplify these conditions to the point were they are practically
usable.
\end{itemize}
But since these three aims are highly interwoven and dependent on
each other, I decided to instead write one larger contribution.
\item Compared with many other papers, the writing style of this paper is
rather detailed. Furthermore, I try to explain the ideas used in the
proofs, instead of just giving a formally correct but unenlightening
proof. Thus, even though the \emph{number of pages} might be larger
than for many other papers, it is hoped that the required \emph{reading
time} is not as high as one would expect, simply because the reader
does not have to work as hard.
\end{itemize}

\subsection{Overview of new results}

\label{subsec:IntroductionOverview}In this subsection, we give a
flavor of the framework for embeddings between decomposition spaces
which will be developed in the present paper. To this end, we first
need to clarify the setting in which we work, by introducing a few
central concepts.

\subsubsection{Concepts related to decomposition spaces}

\label{subsec:IntroductionDecompositionSpaceConcepts}The main ingredient
for defining a decomposition space $\DecompSp{\CalQ}pY$ is a covering
$\CalQ=\left(Q_{i}\right)_{i\in I}$ of an open subset $\CalO$ of
the frequency domain $\R^{\dimension}$. The idea for computing the
(quasi)-norm $\left\Vert g\right\Vert _{\DecompSp{\CalQ}pY}$ of a
distribution $g$ is to decompose the frequency content $\widehat{g}$
according to the covering $\CalQ$, i.e., to consider the ``frequency
localized pieces'' $g_{i}:=\Fourier^{-1}\left(\varphi_{i}\cdot\widehat{g}\right)$,
where $\left(\varphi_{i}\right)_{i\in I}$ is a (suitable) partition
of unity subordinate to $\CalQ$. To obtain the \emph{global} (quasi)-norm
from the frequency localized pieces $g_{i}$, one uses the sequence
space $Y\subset\Compl^{I}$ to aggregate the contributions of the
$g_{i}$, by defining $\left\Vert g\right\Vert _{\DecompSp{\CalQ}pY}:=\left\Vert \left(\left\Vert g_{i}\right\Vert _{L^{p}}\right)_{i\in I}\right\Vert _{Y}$.
In the following paragraphs, we make this rough description more precise,
and we introduce assumptions on the frequency-covering $\CalQ$ and
on the sequence space $Y$ which ensure that $\DecompSp{\CalQ}pY$
is a well-defined quasi-Banach space.

\medskip{}

Although most results in this paper apply in greater generality, for
this introduction we will restrict ourselves to the setting of \textbf{almost
structured coverings}. Above all, such a covering $\CalQ=\left(Q_{i}\right)_{i\in I}$
of an open set $\CalO\subset\R^{\dimension}$ is required to have
the \textbf{finite overlap property}: If we let
\[
i^{\ast}:=\left\{ \ell\in I\with Q_{\ell}\cap Q_{i}\neq\emptyset\right\} 
\]
denote the \textbf{set of neighbors} of $Q_{i}$, the finite overlap
property simply means that the number of neighbors of each set is
uniformly bounded, i.e.\@ that the following constant is finite:
\[
N_{\CalQ}:=\sup_{i\in I}\left|\,i^{\ast}\right|\,.
\]

A further property of almost structured coverings is that one can
associate to $\CalQ=\left(Q_{i}\right)_{i\in I}$ a family of (invertible)
affine transformations $\left(\xi\mapsto T_{i}\,\xi+b_{i}\right)_{i\in I}$
such that—essentially—every set $Q_{i}$ is obtained as the affine
image of a fixed set $Q\subset\R^{\dimension}$. The actual definition
allows for a slightly larger degree of freedom. Precisely, $\CalQ$
is called an almost structured covering of the open subset $\CalO$
of the frequency space $\R^{\dimension}$, if there are \emph{open,
bounded} subsets $Q_{i}'\subset\R^{\dimension}$ and invertible affine
transformations as above such that the following hold:
\begin{enumerate}
\item $\CalQ$ has the finite overlap property and $Q_{i}=T_{i}Q_{i}'+b_{i}$
for all $i\in I$.
\item The following expression (then a constant) is finite:
\[
C_{\CalQ}:=\sup_{i\in I}\,\sup_{\ell\in i^{\ast}}\,\left\Vert T_{i}^{-1}T_{\ell}\,\right\Vert \,.
\]
\item For each $i\in I$, there is an open set $\emptyset\neq P_{i}'\subset\R^{\dimension}$
satisfying $\overline{P_{i}'}\subset Q_{i}'$ and also the following:

\begin{enumerate}
\item The sets $\left\{ P_{i}'\with i\in I\right\} $ and $\left\{ Q_{i}'\with i\in I\right\} $
are finite.
\item We have $\CalO=\bigcup_{i\in I}\left(T_{i}P_{i}'+b_{i}\right)=\bigcup_{i\in I}Q_{i}$.
\end{enumerate}
\end{enumerate}
Given such a covering $\CalQ$, the associated \textbf{decomposition
space (quasi)-norm} of an element $g$ of the \textbf{reservoir} $Z'\left(\CalO\right)$
(which is a certain distribution space, defined below) is
\begin{equation}
\left\Vert g\right\Vert _{\DecompSp{\CalQ}pY}:=\left\Vert \left(\left\Vert \Fourier^{-1}\left(\varphi_{i}\cdot\widehat{g}\right)\right\Vert _{L^{p}}\right)_{i\in I}\right\Vert _{Y}\in\left[0,\infty\right]\,,\label{eq:IntroductionDecompositionSpaceNorm}
\end{equation}
for a suitable \emph{partition of unity} $\Phi=\left(\varphi_{i}\right)_{i\in I}$
subordinate to $\CalQ$ and a suitable \emph{sequence space} $Y\subset\Compl^{I}$.
The corresponding decomposition space is then simply
\[
\DecompSp{\CalQ}pY:=\left\{ g\in Z'\left(\CalO\right)\with\left\Vert g\right\Vert _{\DecompSp{\CalQ}pY}<\infty\right\} \:.
\]

To make the assumptions concerning the partition of unity $\Phi$
more precise, we mention that we require $\Phi$ to be a so-called
\textbf{$L^{p}$-BAPU}, a concept whose precise definition is immaterial
for this introduction; it suffices to know that we have in particular
$\varphi_{i}\in\TestFunctionSpace{Q_{i}}$ and $\sum_{i\in I}\varphi_{i}\equiv1$
on $\CalO$. Nevertheless, an important property of almost structured
coverings is that for each such covering, there exists a family $\Phi=\left(\varphi_{i}\right)_{i\in I}$
which is an $L^{p}$-BAPU simultaneously for all $p\in\left(0,\infty\right]$,
see Theorem~\ref{thm:AlmostStructuredAdmissibleAdmitsBAPU}.

Finally, to ensure that the (quasi)-norm in equation~(\ref{eq:IntroductionDecompositionSpaceNorm})
is well-defined (i.e., that different choices of the BAPU yield equivalent
(quasi)-norms), we require $Y$ to be \textbf{$\CalQ$-regular}. This
means the following:
\begin{enumerate}
\item $Y$ is a quasi-Banach space and a subspace of the space $\Compl^{I}$
of all sequences over $I$.
\item $Y$ is \textbf{solid}, i.e., if $y=\left(y_{i}\right)_{i\in I}\in Y$
and $x=\left(x_{i}\right)_{i\in I}\in\Compl^{I}$ satisfy $\left|x_{i}\right|\leq\left|y_{i}\right|$
for all $i\in I$, then $x\in Y$ with $\left\Vert x\right\Vert _{Y}\leq\left\Vert y\right\Vert _{Y}$.
\item The \textbf{$\CalQ$-clustering map} $\Gamma_{Q}:Y\to Y,x\mapsto x^{\ast}$
with $x_{i}^{\ast}:=\sum_{\ell\in i^{\ast}}x_{\ell}$ is well-defined
and bounded.
\end{enumerate}
As we will see, given these assumptions, the decomposition space (quasi)-norm
(\ref{eq:IntroductionDecompositionSpaceNorm}) is well-defined, with
different BAPUs giving rise to equivalent quasi-norms.

The most common type of $\CalQ$-regular sequence spaces that we will
consider are the weighted $\ell^{q}$ spaces $\ell_{w}^{q}\left(I\right)$
that we defined in Section~\ref{subsec:Notation}. The space $\ell_{w}^{q}\left(I\right)$
is indeed $\CalQ$-regular (see Lemma~\ref{lem:ModeratelyWeightedSpacesAreRegular})
if $w$ is \textbf{$\CalQ$-moderate}, i.e., if the following expression
(then a constant) is finite:
\[
C_{w,\CalQ}:=\sup_{i\in I}\,\sup_{\ell\in i^{\ast}}\frac{w_{i}}{w_{\ell}}\,.
\]

On first sight, our definition of a $\CalQ$-regular sequence space
might appear weaker than what was required by Feichtinger and Gröbner
in the original introduction of decomposition spaces, since in \cite[Definition 2.4]{DecompositionSpaces1}
it is also required that the coordinate evaluation mappings $Y\to\Compl,\left(x_{i}\right)_{i\in I}\mapsto x_{j}$
are bounded. But in fact, we will see in Lemma~\ref{lem:SolidSequenceSpaceEmbedsIntoWeightedLInfty}
that this property is an automatic consequence of the solidity of
$Y$. In other words, if $Y$ is a solid Banach sequence space, then
it is a \textbf{solid BK-space} in the sense of \cite[Definition 2.4]{DecompositionSpaces1}.
The only remaining difference is that in \cite{DecompositionSpaces1},
it is assumed that $Y$ contains all finitely supported sequences;
for greater generality we omitted this assumption.

\medskip{}

To complete the definition of the decomposition space $\DecompSp{\CalQ}pY$,
we still have to define the \textbf{reservoir} $Z'\left(\CalO\right)$,
whose definition is inspired by Triebel's book \cite{TriebelFourierAnalysisAndFunctionSpaces}.
First, we let $Z\left(\CalO\right):=\Fourier\left(\TestFunctionSpace{\CalO}\right)$,
and we equip this space with the unique topology that makes the Fourier
transform $\Fourier:\TestFunctionSpace{\CalO}\to Z\left(\CalO\right)$
into a homeomorphism. The reservoir $Z'\left(\CalO\right)$ is then
simply defined to be the topological dual space of $Z\left(\CalO\right)$,
equipped with the weak-$\ast$-topology. Similar to the case of tempered
distributions, we can then define a \textbf{Fourier transform on $Z'\left(\CalO\right)$}
by duality, i.e.,
\[
\Fourier:Z'\left(\CalO\right)\to\DistributionSpace{\CalO},\phi\mapsto\Fourier\phi:=\widehat{\phi}:=\phi\circ\Fourier\,.
\]
Then $\Fourier:Z'\left(\CalO\right)\to\DistributionSpace{\CalO}$
is a homeomorphism, with inverse $\Fourier^{-1}:\DistributionSpace{\CalO}\to Z'\left(\CalO\right),\phi\mapsto\phi\circ\Fourier^{-1}$.

We emphasize that with this choice of the reservoir $Z'\left(\CalO\right)$,
the quasi-norm $\left\Vert \mybullet\right\Vert _{\DecompSp{\CalQ}pY}$
as defined in equation~(\ref{eq:IntroductionDecompositionSpaceNorm})
really makes sense: For $g\in Z'\left(\CalO\right)$, we have $\widehat{g}\in\DistributionSpace{\CalO}$.
Because of $\varphi_{i}\in\TestFunctionSpace{\CalO}$, this implies
that $\varphi_{i}\cdot\widehat{g}\in\Schwartz'\left(\R^{\dimension}\right)$
is a (tempered) distribution with compact support. By the Paley-Wiener
theorem, we thus see that $g_{i}:=\Fourier^{-1}\left(\varphi_{i}\cdot\widehat{g}\right)$
is a well-defined smooth function $g_{i}\in C^{\infty}\left(\R^{\dimension}\right)$,
so that $\left\Vert g_{i}\right\Vert _{L^{p}}\in\left[0,\infty\right]$
is well-defined. Thus, $\left\Vert g\right\Vert _{\DecompSp{\CalQ}pY}=\left\Vert \left(\left\Vert g_{i}\right\Vert _{L^{p}}\right)_{i\in I}\right\Vert _{Y}\in\left[0,\infty\right]$
is well-defined, with the convention that $\left\Vert g\right\Vert _{\DecompSp{\CalQ}pY}=\infty$
if $g_{i}\notin L^{p}\left(\R^{\dimension}\right)$ for some $i\in I$,
or if $\left(\left\Vert g_{i}\right\Vert _{L^{p}}\right)_{i\in I}\notin Y$.

At this point the reader might wonder why we chose to define the decomposition
space $\DecompSp{\CalQ}pY$ as a subspace of the uncommon space $Z'\left(\CalO\right)$
instead of the well-known standard space $\Schwartz'\left(\R^{\dimension}\right)$,
which is for example used in \cite{BorupNielsenDecomposition}. There
are two principal reasons for this: First, in case of $\CalO\subsetneq\R^{\dimension}$,
the quasi-norm $\left\Vert \mybullet\right\Vert _{\DecompSp{\CalQ}pY}$
from equation~(\ref{eq:IntroductionDecompositionSpaceNorm}) is not
positive definite on $\Schwartz'\left(\R^{\dimension}\right)$, so
that one would need to factor out the subspace of all $g\in\Schwartz'\left(\R^{\dimension}\right)$
with $\widehat{g}\equiv0$ on $\CalO$, i.e., with $\left\langle \widehat{g},\varphi\right\rangle _{\Schwartz'}=0$
for all $\varphi\in\Schwartz\left(\R^{\dimension}\right)$ with $\supp\varphi\subset\CalO$.
The resulting quotient space is then just as exotic as the space $Z'\left(\CalO\right)$
itself. Secondly, and more importantly, even in case of $\CalO=\R^{\dimension}$,
if one defines $\DecompSp{\CalQ}pY$ as a subspace of $\Schwartz'\left(\R^{\dimension}\right)$,
one occasionally obtains an \emph{incomplete} space, as we will see
in Example~\ref{exa:BorupNielsenDecompositionSpaceIncomplete}. This
does not happen, however, if one uses the reservoir $Z'\left(\CalO\right)$,
as shown in Theorem~\ref{thm:DecompositionSpaceComplete}.

Still, most working analysts have acquired significant intuition on
which operations (like differentiation and multiplication with suitable
smooth functions) are permitted on the distribution spaces $\Schwartz'\left(\R^{\dimension}\right)$
or $\DistributionSpace{\CalO}$, while this is not true for the space
$\SpaceReservoir{\CalO}$. For this—and other—reasons, we will often
find it convenient to work with the so-called \textbf{Fourier-side
decomposition space} 
\[
\FourierDecompSp{\CalQ}pY:=\left\{ g\in\DistributionSpace{\CalO}\with\left\Vert g\right\Vert _{\FourierDecompSp{\CalQ}pY}:=\left\Vert \Fourier^{-1}g\right\Vert _{\DecompSp{\CalQ}pY}<\infty\right\} \,.
\]
We emphasize that the Fourier transform $\Fourier:\DecompSp{\CalQ}pY\to\FourierDecompSp{\CalQ}pY$
is an isometric isomorphism, so that \emph{it does not really matter
which of the spaces we consider; all results (in particular embeddings)
for the Fourier-side spaces translate to results about the space-side
spaces and vice versa}.

\subsubsection{Sufficient conditions for embeddings}

\label{subsec:IntroductionSufficientConditions}In the following,
we always assume that we are given two almost structured coverings
\[
\CalQ=\left(Q_{i}\right)_{i\in I}=\left(\vphantom{P_{j}'}T_{i}Q_{i}'+b_{i}\right)_{i\in I}\qquad\text{ and }\qquad\CalP=\left(P_{j}\right)_{j\in J}=\left(S_{j}P_{j}'+c_{j}\right)_{j\in J}
\]
of the open subsets $\CalO\subset\R^{\dimension}$ and $\CalO'\subset\R^{\dimension}$,
respectively. Furthermore, we assume that $Y\subset\Compl^{I}$ and
$Z\subset\Compl^{J}$ are $\CalQ$-regular and $\CalP$-regular, respectively.
Finally, we fix families $\Phi=\left(\varphi_{i}\right)_{i\in I}$
and $\Psi=\left(\psi_{j}\right)_{j\in J}$ which are $L^{p}$-BAPUs
simultaneously for all $p\in\left(0,\infty\right]$ for $\CalQ$ and
$\CalP$, respectively.

Then, we are interested in \emph{sufficient} conditions for the existence
of an embedding of the form
\begin{equation}
\FourierDecompSp{\CalQ}{p_{1}}Y\hookrightarrow\FourierDecompSp{\CalP}{p_{2}}Z.\label{eq:IntroductionEmbeddingAbstract}
\end{equation}
Part of the problem is to formulate an appropriate notion of an embedding
for this setting, especially if $\CalO\neq\CalO'$. Necessary conditions
are the subject of the next subsection; there, we will see that the
sufficient criteria given in this subsection frequently yield \emph{complete
characterizations} of the existence of (\ref{eq:IntroductionEmbeddingAbstract}).

In this introduction, we will only state a few of our results which
convey the general flavor. For a more detailed discussion, we refer
to the remainder of the paper; in particular we refer to Section~\ref{sec:SummaryOfEmbeddingResults},
where we summarize the most useful embedding results.

In order to state our results, we need the notion of \textbf{almost
subordinateness} of two coverings: Recall from above the notation
$j^{\ast}=\left\{ \ell\in J\with P_{\ell}\cap P_{j}\neq\emptyset\right\} $
for the set of all (indices of) neighbors of the set $P_{j}$. More
generally, for $L\subset J$, let $L^{\ast}:=\bigcup_{\ell\in L}\ell^{\ast}\subset J$
and define inductively $L^{0\ast}:=L$ and $L^{\left(n+1\right)\ast}:=\left(L^{n\ast}\right)^{\ast}$.
Now, we say that $\CalQ$ is \textbf{almost subordinate to $\CalP$},
if there is some $n\in\N_{0}$ with the following property:
\[
\forall\,i\in I\,\exists\:j_{i}\in J:\qquad Q_{i}\subset P_{j_{i}}^{n\ast}\qquad\text{where }\qquad P_{j}^{n\ast}:=\bigcup_{\ell\in j^{n\ast}}P_{\ell}.
\]
The notation $Q_{i}^{n\ast}$ (or $P_{j}^{n\ast}$) will be used over
and over in this paper and is thus worth remembering. Roughly speaking,
the definition of almost subordinateness formalizes the intuition
that $\CalQ$ is ``finer than $\CalP$''. Note in particular that
we have $\CalO\subset\CalO'$ if $\CalQ$ is almost subordinate to
$\CalP$.

Finally, we need the notion of \textbf{nested sequence spaces}, which
we only introduce informally in this introduction. Given a solid sequence
space $X\subset\Compl^{K}$, an exponent $q\in\left(0,\infty\right]$
and for each $k\in K$ some set $M_{k}$, as well as a weight $u=\left(u_{k,m}\right)_{k\in K,m\in M_{k}}$,
we define the space $X\left(\left[\ell_{u}^{q}\left(M_{k}\right)\right]_{k\in K}\right)$
as the space of all complex-valued sequences $x=\left(x_{m}\right)_{m\in M}$,
with $M:=\bigcup_{k\in K}M_{k}$, for which the (quasi)-norm
\[
\left\Vert x\right\Vert _{X\left(\left[\ell_{u}^{q}\left(M_{k}\right)\right]_{k\in K}\right)}:=\left\Vert \left(\left\Vert \left(u_{k,m}\cdot x_{m}\right)_{m\in M_{k}}\right\Vert _{\ell^{q}\left(M_{k}\right)}\right)_{k\in K}\right\Vert _{X}
\]
is finite.

Using these notions, we can finally state our first embedding result.
It might look daunting at first, but as we will see in Corollary~\ref{cor:IntroductionFineIntoCoarseSimplified}
below, its assumptions can be greatly simplified in most cases.
\begin{thm}
\label{thm:IntroductionFineIntoCoarse}(special case of Corollary~\ref{cor:EmbeddingFineIntoCoarse})

Assume that $\CalQ$ is almost subordinate to $\CalP$, and for each
$j\in J$ define
\begin{equation}
I_{j}:=\left\{ i\in I\with Q_{i}\cap P_{j}\neq\emptyset\right\} \,.\label{eq:IntroductionIntersectionSet}
\end{equation}
If $p_{1}\leq p_{2}$ and if the embedding
\begin{equation}
Y\hookrightarrow Z\left(\left[\ell_{\left|\det T_{i}\right|^{p_{1}^{-1}-p_{2}^{-1}}}^{\LowerExpo{p_{2}}}\left(I_{j}\right)\right]_{j\in J}\right)\label{eq:IntroductionQIntoPSequenceSpaceEmbedding}
\end{equation}
holds, then the map $\iota$ is well-defined and bounded, where
\[
\iota:\FourierDecompSp{\CalQ}{p_{1}}Y\to\FourierDecompSp{\CalP}{p_{2}}Z,f\mapsto\sum_{i\in I}\varphi_{i}f\,.
\]
Furthermore, $\iota f\in\DistributionSpace{\CalO'}$ is an extension
of $f\in\FourierDecompSp{\CalQ}{p_{1}}Y\subset\DistributionSpace{\CalO}$.
In particular, if $\CalO=\CalO'$, then $\iota f=f$ for all $f\in\FourierDecompSp{\CalQ}{p_{1}}Y$.
\end{thm}

The main requirement of the preceding theorem is the validity of the
embedding~(\ref{eq:IntroductionQIntoPSequenceSpaceEmbedding}). Of
course, verifying this condition in general can be difficult. But
if $Y=\ell_{w}^{q_{1}}\left(I\right)$ and $Z=\ell_{v}^{q_{2}}\left(J\right)$
are weighted $\ell^{q}$ spaces, we will see in a moment that these
conditions can be simplified considerably.

Even further simplifications are possible if $w$ and $\CalQ$ are
\textbf{relatively $\CalP$-moderate}: For $w$, this means that the
constant
\[
C_{w,\CalQ,\CalP}:=\sup_{j\in J}\,\sup_{i,\ell\in I_{j}}\,\frac{w_{i}}{w_{\ell}}
\]
is finite. Roughly speaking, this condition states that the weight
$w$ is ``uniformly essentially constant'' on the collection of
all ``small'' sets $Q_{i}$ which intersect the same ``large''
set $P_{j}$. Finally, $\CalQ=\left(T_{i}Q_{i}'+b_{i}\right)_{i\in I}$
is called relatively $\CalP$-moderate if the weight $\left(\left|\det T_{i}\right|\right)_{i\in I}$
is relatively $\CalP$-moderate. Roughly speaking, this means that
any two ``small'' sets $Q_{i},Q_{\ell}$ which intersect the same
``large'' set $P_{j}$ have to be of essentially the same measure.
\begin{cor}
\label{cor:IntroductionFineIntoCoarseSimplified}(special case of
Corollary~\ref{cor:EmbeddingFineInCoarseSimplification})

Assume that $\CalQ$ is almost subordinate to $\CalP$ and that the
weights $w=\left(w_{i}\right)_{i\in I}$ and $v=\left(v_{j}\right)_{j\in J}$
are moderate with respect to $\CalQ$ and $\CalP$, respectively.
Finally, assume $Y=\ell_{w}^{q_{1}}\left(I\right)$ and $Z=\ell_{v}^{q_{2}}\left(J\right)$.

Then, validity of the embedding~(\ref{eq:IntroductionQIntoPSequenceSpaceEmbedding})
is equivalent to finiteness of
\begin{equation}
\left\Vert \left(v_{j}\cdot\left\Vert \left(w_{i}^{-1}\cdot\left|\det T_{i}\right|^{p_{1}^{-1}-p_{2}^{-1}}\right)_{i\in I_{j}}\right\Vert _{\ell^{\LowerExpo{p_{2}}\cdot\left(q_{1}/\LowerExpo{p_{2}}\right)'}}\right)_{j\in J}\right\Vert _{\ell^{q_{2}\cdot\left(q_{1}/q_{2}\right)'}}\,.\label{eq:IntroductionQIntoPSimplified}
\end{equation}
Furthermore, if $w$ and $\CalQ$ are relatively $\CalP$-moderate
and if $\CalO=\CalO'$, then finiteness of (\ref{eq:IntroductionQIntoPSimplified})
is equivalent to finiteness of
\begin{equation}
\left\Vert \left(\frac{v_{j}}{w_{i_{j}}}\cdot\left|\det\smash{T_{i_{j}}}\right|^{\frac{1}{p_{1}}-\left(\smash{\frac{1}{\LowerExpo{p_{2}}}}-\frac{1}{q_{1}}\right)_{+}-\frac{1}{p_{2}}}\cdot\left|\det S_{j}\right|^{\left(\smash{\frac{1}{\LowerExpo{p_{2}}}}-\frac{1}{q_{1}}\right)_{+}}\right)_{j\in J}\right\Vert _{\ell^{q_{2}\cdot\left(q_{1}/q_{2}\right)'}},\label{eq:IntroductionQIntoPRelativelyModerateCondition}
\end{equation}
where for each $j\in J$, an arbitrary $i_{j}\in I$ with $P_{j}\cap Q_{i_{j}}\neq\emptyset$
can be selected.
\end{cor}

We remark that the exponent $p\cdot\left(q/p\right)'$ uses the convention
$p\cdot\left(q/p\right)'=\infty$ if $p=\infty$ or if $q\leq p<\infty$.
In the remaining cases, we use the usual arithmetic in $\left(0,\infty\right]$
to evaluate the product, see Lemma~\ref{lem:EmbeddingBetweenWeightedSequenceSpaces}.

Note that the above conditions—in particular condition~(\ref{eq:IntroductionQIntoPRelativelyModerateCondition})—are
usually not hard to check. One simply needs to verify finiteness of
the $\ell^{q}$-norm of a \emph{single}, explicitly computable sequence.

\medskip{}

In the above results, we always assumed $\CalQ$ to be almost subordinate
to $\CalP$. In the following result, we consider the ``reverse''
case in which $\CalP$ is almost subordinate to $\CalQ$.
\begin{thm}
\label{thm:IntroductionCoarseIntoFine}(special case of Corollary~\ref{cor:EmbeddingCoarseIntoFine})\vspace{0.1cm}

Assume that $\CalP$ is almost subordinate to $\CalQ$, and for each
$i\in I$ define
\[
J_{i}:=\left\{ j\in J\with P_{j}\cap Q_{i}\neq\emptyset\right\} \,.\vspace{-0.1cm}
\]
If $p_{1}\leq p_{2}$ and if we have\vspace{-0.1cm}
\begin{equation}
Y\left(\vphantom{\ell_{u}^{M}}\left[\vphantom{\ell_{u}^{M}}\smash{\ell_{u}^{\UpperExpo{p_{1}}}}\!\left(J_{i}\right)\right]_{i\in I}\right)\hookrightarrow Z\qquad\text{where}\qquad u_{i,j}:=\begin{cases}
\left|\det S_{j}\right|^{p_{2}^{-1}-1}\cdot\left|\det T_{i}\right|^{1-p_{1}^{-1}}, & \text{if }p_{1}<1,\\
\vphantom{\rule{0.1cm}{0.55cm}}\left|\det S_{j}\right|^{p_{2}^{-1}-p_{1}^{-1}}, & \text{if }p_{1}\geq1,
\end{cases}\label{eq:IntroductionPIntoQEmbeddingAssumption}
\end{equation}
then we have $\CalO'\subset\CalO$ and the restriction map $\iota$
is well-defined and bounded, where
\[
\iota:\FourierDecompSp{\CalQ}{p_{1}}Y\subset\DistributionSpace{\CalO}\to\FourierDecompSp{\CalP}{p_{2}}Z\subset\DistributionSpace{\CalO'},f\mapsto f|_{\TestFunctionSpace{\CalO'}}\:.\qedhere
\]
\end{thm}

As above, the preceding assumptions simplify considerably if $Y,Z$
are weighted $\ell^{q}$ spaces:
\begin{cor}
\label{cor:IntroductionCoarseIntoFineSimplified}(special case of
Corollary~\ref{cor:EmbeddingCoarseIntoFineSimplification} and Remark~\ref{rem:SufficientCoarseIntoFineSimplification})

Assume that $\CalP$ is almost subordinate to $\CalQ$ and that the
weights $w=\left(w_{i}\right)_{i\in I}$ and $v=\left(v_{j}\right)_{j\in J}$
are moderate with respect to $\CalQ$ and $\CalP$, respectively.
Finally, assume $Y=\ell_{w}^{q_{1}}\left(I\right)$ and $Z=\ell_{v}^{q_{2}}\left(J\right)$.

Then, condition~(\ref{eq:IntroductionPIntoQEmbeddingAssumption})
is equivalent to finiteness of
\begin{equation}
\begin{cases}
\,\left\Vert \left(w_{i}^{-1}\cdot\left\Vert \left(\left|\det S_{j}\right|^{p_{1}^{-1}-p_{2}^{-1}}\cdot v_{j}\right)_{j\in J_{i}}\right\Vert _{\ell^{q_{2}\cdot\left(\UpperExpo{p_{1}}/q_{2}\right)'}}\right)_{i\in I}\right\Vert _{\ell^{q_{2}\cdot\left(q_{1}/q_{2}\right)'}} & \text{if }p_{1}\geq1,\\
\vphantom{\rule{0.1cm}{0.9cm}}\left\Vert \left(\left|\det T_{i}\right|^{p_{1}^{-1}-1}\cdot w_{i}^{-1}\cdot\left\Vert \left(\left|\det S_{j}\right|^{1-p_{2}^{-1}}\cdot v_{j}\right)_{j\in J_{i}}\right\Vert _{\ell^{q_{2}\cdot\left(\UpperExpo{p_{1}}/q_{2}\right)'}}\right)_{i\in I}\right\Vert _{\ell^{q_{2}\cdot\left(q_{1}/q_{2}\right)'}} & \text{if }p_{1}<1.
\end{cases}\label{eq:IntroductionPIntoQSimplified}
\end{equation}

Furthermore, if $v$ and $\CalP$ are relatively $\CalQ$-moderate,
and if $\CalO=\CalO'$, then finiteness of (\ref{eq:IntroductionPIntoQSimplified})
is equivalent to finiteness of
\begin{equation}
\left\Vert \left(\frac{v_{j_{i}}}{w_{i}}\cdot\left|\det\smash{S_{j_{i}}}\right|^{\frac{1}{p_{1}}-\left(\frac{1}{q_{2}}-\smash{\frac{1}{\SignedUpperExpo{p_{1}}}}\right)_{+}-\frac{1}{p_{2}}}\cdot\left|\det T_{i}\right|^{\left(\frac{1}{q_{2}}-\smash{\frac{1}{\SignedUpperExpo{p_{1}}}}\right)_{+}}\right)_{\!\!i\in I}\right\Vert _{\ell^{q_{2}\cdot\left(q_{1}/q_{2}\right)'}},\label{eq:IntroductionPIntoQRelativelyModerateCondition}
\end{equation}
where for each $i\in I$, an arbitrary $j_{i}\in J$ with $Q_{i}\cap P_{j_{i}}\neq\emptyset$
can be selected.
\end{cor}

We emphasize once more that the preceding two theorems are merely
the most accessible of our embedding results. More involved—but more
flexible—results are also derived. In particular, we will obtain a
variant of the preceding theorems were neither $\CalQ$ needs to be
almost subordinate to $\CalP$, nor vice versa. Instead, it suffices
if one can write $\CalO\cap\CalO'=A\cup B$, where $\CalQ$ is almost
subordinate to $\CalP$ ``near $A$'' and $\CalP$ is almost subordinate
to $\CalQ$ ``near $B$''. Since the precise statement is slightly
involved, we don't spell it out in this introduction. Instead, we
refer the interested reader to Theorem~\ref{thm:SummaryMixedSubordinateness}
below.

Finally, there are also results which do not need any subordinateness
at all (see Theorem~\ref{thm:NoSubordinatenessWithConsiderationOfOverlaps}).
The preceding theorems, however, have the advantage that they are
(relatively) easy to apply and that they are fairly sharp—as we will
see now.

\subsubsection{Necessary conditions for embeddings}

\label{subsec:IntroductionNecessaryConditions}In the previous subsection,
we presented various \emph{sufficient} conditions for the existence
of an embedding between two decomposition spaces. In the present generality,
to my knowledge, these are the best known results—simply because no
other results are known which apply in this general setting. Nevertheless,
it is crucial to know how \emph{sharp} these sufficient conditions
are.

As we will see in this subsection, the results of Corollaries \ref{cor:IntroductionFineIntoCoarseSimplified}
and \ref{cor:IntroductionCoarseIntoFineSimplified} \emph{are indeed
reasonably sharp}: We will see for $p_{2}\in\left(0,2\right]\cup\left\{ \infty\right\} $
that the (sufficient) conditions given in Theorem~\ref{thm:IntroductionFineIntoCoarse}
and Corollary~\ref{cor:IntroductionFineIntoCoarseSimplified} are
also \emph{necessary} conditions for the existence of the stated embedding.
Likewise, for $p_{1}\in\left[2,\infty\right]$, the (sufficient) conditions
given in Theorem~\ref{thm:IntroductionCoarseIntoFine} and Corollary~\ref{cor:IntroductionCoarseIntoFineSimplified}
are necessary for the existence of the embedding. Thus, for these
ranges of $p_{2}$ or $p_{1}$, we achieve a \emph{complete characterization}
of the existence of the embeddings.

Finally, in the case of weighted $\ell^{q}$ spaces $Y=\ell_{w}^{q_{1}}\left(I\right)$
and $Z=\ell_{v}^{q_{2}}\left(J\right)$, this characterization extends
to \emph{all} values of $p_{2}$, \emph{as long as} $\CalQ$ and $w$
are relatively $\CalP$-moderate, or to all values of $p_{1}$, if
$\CalP$ and $v$ are relatively $\CalQ$-moderate.

Let us now turn to precise statements. To ensure flexibility, we assume
in the following that we are given a set $K\subset\CalO\cap\CalO'$,
for which the identity map
\begin{equation}
\iota:\left(\,\smash{\CalD_{K}},\smash{\left\Vert \mybullet\right\Vert _{\FourierDecompSp{\CalQ}{p_{1}}Y}}\,\right)\to\FourierDecompSp{\CalP}{p_{2}}Z,f\mapsto f\quad\!\text{with}\quad\!\CalD_{K}\!:=\!\left\{ f\!\in\!\TestFunctionSpace{\R^{\dimension}}:\supp f\!\subset\!K\!\right\} \label{eq:IntroductionNecessaryEmbeddingQIntoPAssumption}
\end{equation}
is well-defined and bounded. In most practical applications, one will
choose $K=\CalO\cap\CalO'$.

Based on the boundedness of $\iota$, we will derive various conditions
which are thus necessary conditions for the existence of the embedding.
For the first such condition, note that Theorems \ref{thm:IntroductionFineIntoCoarse}
and \ref{thm:IntroductionCoarseIntoFine} always included the assumption
$p_{1}\leq p_{2}$. This assumption is unavoidable:
\begin{thm}
\label{thm:IntroductionSimplestNecessaryCriterion}(special case of
Lemma~\ref{lem:SimpleNecessaryCondition})

If the map $\iota$ from equation~(\ref{eq:IntroductionNecessaryEmbeddingQIntoPAssumption})
is well-defined and bounded and if there are $i\in I$ and $j\in J$
satisfying $K^{\circ}\cap Q_{i}\cap P_{j}\neq\emptyset$ and $\delta_{i}\in Y$,
then $p_{1}\leq p_{2}$ and $\delta_{j}\in Z$.

In case of $p_{1}=p_{2}$, we even have
\begin{equation}
\left\Vert \delta_{j}\right\Vert _{Z}\lesssim\left\Vert \delta_{i}\right\Vert _{Y},\label{eq:IntroductionSimplestNecessaryCriterionNormEstimate}
\end{equation}
where the implied constant is \emph{independent} of $i,j$.
\end{thm}

\begin{rem*}
If $Y=\ell_{w}^{q_{1}}\left(I\right)$ and $Z=\ell_{v}^{q_{2}}\left(J\right)$,
then (\ref{eq:IntroductionSimplestNecessaryCriterionNormEstimate})
simply means $v_{j}\lesssim w_{i}$ if $K^{\circ}\cap Q_{i}\cap P_{j}\neq\emptyset$.
\end{rem*}
For completely general (almost structured) coverings $\CalQ,\CalP$,
it is hard to say more. Thus, for our next result, we will assume
that $\CalQ$—or at least a subfamily $\CalQ_{I_{0}}=\left(Q_{i}\right)_{i\in I_{0}}$
with $I_{0}\subset I$—is almost subordinate to $\CalP$, i.e.\@
that there is a fixed $n\in\N_{0}$ and for each $i\in I_{0}$ some
$j_{i}\in J$ satisfying $Q_{i}\subset P_{j_{i}}^{n\ast}$. In this
case, much more can be said:
\begin{thm}
\label{thm:IntroductionNecessaryFineIntoCoarse}(special case of
Theorem~\ref{thm:BurnerNecessaryConditionFineInCoarse})

Let $\emptyset\neq I_{0}\subset I$ and assume that $\iota$ as in
equation~(\ref{eq:IntroductionNecessaryEmbeddingQIntoPAssumption})
is bounded, with $K:=\bigcup_{i\in I_{0}}Q_{i}$. Furthermore, assume
that $\CalQ_{I_{0}}:=\left(Q_{i}\right)_{i\in I_{0}}$ is almost subordinate
to $\CalP$.

If there is some $i_{0}\in I_{0}$ with $\delta_{i_{0}}\in Y$, then
we have $p_{1}\leq p_{2}$, and we have a bounded embedding
\begin{equation}
\ell_{0}\left(I_{0}\right)\cap Y\hookrightarrow Z\left(\left[\,\vphantom{T^{j}}\smash{\ell_{\left|\det T_{i}\right|^{p_{1}^{-1}-p_{2}^{-1}}}^{p_{2}}}\left(I_{j}\cap I_{0}\right)\,\right]_{j\in J}\right)\,.\label{eq:IntroductionNecessaryQIntoPDiscreteEmbedding}
\end{equation}
Here, $\ell_{0}\left(I_{0}\right)$ is the space of all sequences
on $I$ with a finite support which is contained in $I_{0}$.
\end{thm}

\begin{rem*}
For $I_{0}=I$, this embedding coincides with the (sufficient) condition
given in Theorem~\ref{thm:IntroductionFineIntoCoarse}, with the
exception that the ``inner'' norm in Theorem~\ref{thm:IntroductionFineIntoCoarse}
is $\ell^{\LowerExpo{p_{2}}}$, whereas here, it is $\ell^{p_{2}}$
(and we have to restrict to $\ell_{0}\left(I_{0}\right)=\ell_{0}\left(I\right)$,
which is immaterial in most cases).

But for $p_{2}\in\left(0,2\right]$, we have $p_{2}=\LowerExpo{p_{2}}$,
so that we get a \emph{complete characterization} of the existence
of the embedding $\FourierDecompSp{\CalQ}{p_{1}}Y\hookrightarrow\FourierDecompSp{\CalP}{p_{2}}Z$,
at least if $\CalQ$ is almost subordinate to $\CalP$. Furthermore,
for $p_{2}=\infty$, we will see in Theorem~\ref{thm:BurnerNecessaryConditionFineInCoarse}
that if $\FourierDecompSp{\CalQ}{p_{1}}Y\hookrightarrow\FourierDecompSp{\CalP}{p_{2}}Z$
holds, then condition~(\ref{eq:IntroductionNecessaryQIntoPDiscreteEmbedding})
is satisfied with the ``inner norm'' $\ell^{1}=\ell^{\LowerExpo{p_{2}}}$
instead of $\ell^{p_{2}}$, so that the complete characterization
holds for all $p_{2}\in\left(0,2\right]\cup\left\{ \infty\right\} $.

Finally, in case of $Y=\ell_{w}^{q_{1}}\left(I\right)$ and $Z=\ell_{v}^{q_{2}}\left(J\right)$,
the validity of the embedding~(\ref{eq:IntroductionNecessaryQIntoPDiscreteEmbedding})
can again be reformulated in a form that is easier to verify, similar
to condition~(\ref{eq:IntroductionQIntoPSimplified}): One simply
has to replace $\LowerExpo{p_{2}}$ by $p_{2}$ and $I_{j}$ by $I_{j}\cap I_{0}$
everywhere.
\end{rem*}
Of course, the difference between $\LowerExpo{p_{2}}$ (in condition~(\ref{eq:IntroductionQIntoPSequenceSpaceEmbedding}))
and $p_{2}$ (in condition~(\ref{eq:IntroductionNecessaryQIntoPDiscreteEmbedding}))
is somewhat unsatisfactory. In general, I doubt that it can be removed.
Under suitable hypothesis, however, it can:
\begin{thm}
\label{thm:IntroductionNecessaryFineIntoCoarseRelativelyModerate}(special
case of Theorem~\ref{thm:NecessaryConditionForModerateCoveringFineInCoarse})

Assume that $\CalQ$ is almost subordinate to $\CalP$ and that $\iota$
as in equation~(\ref{eq:IntroductionNecessaryEmbeddingQIntoPAssumption})
is bounded, with $K=\CalO=\CalO'$. Let $Y=\ell_{w}^{q_{1}}\left(I\right)$
and $Z=\ell_{v}^{q_{2}}\left(J\right)$, with $w,v$ moderate with
respect to $\CalQ$ and $\CalP$, respectively. Finally, assume that
$\CalQ$ and $w$ are relatively $\CalP$-moderate.

Then, the quantity in equation~(\ref{eq:IntroductionQIntoPRelativelyModerateCondition})
is finite.
\end{thm}

\begin{rem*}
In view of Corollary~\ref{cor:IntroductionFineIntoCoarseSimplified},
we have thus achieved a \emph{complete characterization} of the existence
of the embedding $\FourierDecompSp{\CalQ}{p_{1}}{\ell_{w}^{q_{1}}}\hookrightarrow\FourierDecompSp{\CalP}{p_{2}}{\ell_{v}^{q_{2}}}$,
assuming that $\CalQ$ is almost subordinate to $\CalP$ and that
$\CalQ$ and $w$ are relatively $\CalP$-moderate. A variant of the
above theorem remains true if only a subfamily $\CalQ_{I_{0}}$ of
$\CalQ$ is almost subordinate to $\CalP$; but for the sake of succinctness,
we omit this more general theorem in this introduction.
\end{rem*}
Above, we assumed $\CalQ$ to be almost subordinate to $\CalP$. As
in the previous subsection, we now consider the ``reverse'' case
in which $\CalP$—or at least some subfamily $\CalP_{J_{0}}$ of $\CalP$—is
almost subordinate to $\CalQ$. As above, one can obtain a satisfying
necessary condition (which is very similar to the sufficient condition
from Theorem~\ref{thm:IntroductionCoarseIntoFine}) under this assumption:
\begin{thm}
\label{thm:IntroductionNecessaryCoarseIntoFine}(special case of
Theorem~\ref{thm:BurnerNecessaryConditionCoarseInFine})

Let $\emptyset\neq J_{0}\subset J$ and assume that the map $\iota$
from equation~(\ref{eq:IntroductionNecessaryEmbeddingQIntoPAssumption})
is bounded, with $K:=\bigcup_{j\in J_{0}}P_{j}$. Furthermore, assume
that $\CalP_{J_{0}}:=\left(P_{j}\right)_{j\in J_{0}}$ is almost subordinate
to $\CalQ$.

Then the map
\begin{equation}
\eta\!:\!\ell_{0}\!\left(J_{0}\right)\cap Y\!\left(\!\left[\,\vphantom{T^{j}}\smash{\ell_{\left|\det S_{j}\right|^{p_{2}^{-1}-p_{1}^{-1}}}^{p_{1}}}\!\left(J_{i}\cap\!J_{0}\right)\right]_{i\in I}\right)\!\hookrightarrow Z,\left(x_{j}\right)_{j\in J_{0}}\mapsto\!\left(x_{j}\right)_{j\in J}\text{ with }x_{j}=0\text{ for }j\!\in\!J\!\setminus\!J_{0}\label{eq:IntroductionNecessaryPIntoQDiscreteEmbedding}
\end{equation}
is well-defined and bounded.
\end{thm}

\begin{rem*}
For $p_{1}\in\left[2,\infty\right]$, we have $\UpperExpo{p_{1}}=p_{1}$
and the weight $u$ defined in Theorem~\ref{thm:IntroductionCoarseIntoFine}
(equation~(\ref{eq:IntroductionPIntoQEmbeddingAssumption})) satisfies
$u_{i,j}=\left|\det S_{j}\right|^{p_{2}^{-1}-p_{1}^{-1}}$. Therefore,
assuming $J_{0}=J$, the condition given above coincides with the
sufficient condition~(\ref{eq:IntroductionPIntoQEmbeddingAssumption})
from Theorem~\ref{thm:IntroductionCoarseIntoFine} (up to the intersection
with $\ell_{0}\left(J_{0}\right)=\ell_{0}\left(J\right)$), so that
we achieve a \emph{complete characterization} of the embedding $\FourierDecompSp{\CalQ}{p_{1}}Y\hookrightarrow\FourierDecompSp{\CalP}{p_{2}}Z$,
as long as $p_{1}\in\left[2,\infty\right]$ and as long as $\CalP$
is almost subordinate to $\CalQ$.
\end{rem*}
Finally, one obtains a complete characterization for all $p_{1}\in\left(0,\infty\right]$
if one assumes slightly more:
\begin{thm}
\label{thm:IntroductionNecessaryCoarseIntoFineRelativelyModerate}(special
case of Theorem~\ref{thm:NecessaryConditionForModerateCoveringCoarseInFine})

Assume that $\CalP$ is almost subordinate to $\CalQ$ and that $\iota$
as in equation~(\ref{eq:IntroductionNecessaryEmbeddingQIntoPAssumption})
is bounded, with $K=\CalO=\CalO'$. Furthermore, assume that $Y=\ell_{w}^{q_{1}}\left(I\right)$
and $Z=\ell_{v}^{q_{2}}\left(J\right)$, with $w,v$ moderate with
respect to $\CalQ$ and $\CalP$, respectively. Finally, assume that
$\CalP$ and $v$ are relatively $\CalQ$-moderate.

Then, the quantity in equation~(\ref{eq:IntroductionPIntoQRelativelyModerateCondition})
is finite.
\end{thm}

\begin{rem*}
In view of Corollary~\ref{cor:IntroductionCoarseIntoFineSimplified},
we have thus achieved a \emph{complete characterization} of the existence
of the embedding $\FourierDecompSp{\CalQ}{p_{1}}{\ell_{w}^{q_{1}}}\hookrightarrow\FourierDecompSp{\CalP}{p_{2}}{\ell_{v}^{q_{2}}}$,
assuming that $\CalP$ is almost subordinate to $\CalQ$ and that
$\CalP$ and $v$ are relatively $\CalQ$-moderate. Again, a variant
of this statement remains true if only a subfamily $\CalP_{J_{0}}$
of $\CalP$ is almost subordinate to and relatively moderate with
respect to $\CalQ$.
\end{rem*}
This concludes the overview of the sufficient conditions and necessary
conditions for the existence of embeddings between decomposition spaces
which are derived in this paper.

However, using the same techniques, we will prove one more stunning
result which illustrates the ``rigidity'' of decomposition spaces:
The mapping $\left(p,q,w,\CalQ\right)\mapsto\FourierDecompSp{\CalQ}p{\ell_{w}^{q}}$
is essentially injective. To explain this, we first have to introduce
a suitable notion of \emph{equivalence} for coverings: Inspired by
\cite[Definition 3.3]{DecompositionSpaces1}, we say that two coverings
$\CalQ,\CalP$ of the same set $\CalO\subset\R^{\dimension}$ are
\textbf{weakly equivalent} if
\[
\sup_{i\in I}\left|\left\{ j\in J\with P_{j}\cap Q_{i}\neq\emptyset\right\} \right|<\infty\qquad\text{ and }\qquad\sup_{j\in J}\left|\left\{ i\in I\with Q_{i}\cap P_{j}\neq\emptyset\right\} \right|<\infty.
\]
If $\CalQ,\CalP$ consist only of open, connected sets, one can show
(see Corollary~\ref{cor:WeakSubordinationImpliesSubordinationIfConnected})
that weak equivalence of $\CalQ,\CalP$ already implies that $\CalQ,\CalP$
are \textbf{equivalent}, i.e.\@ that $\CalQ$ is almost subordinate
to $\CalP$ and vice versa.

Now, we state our ``rigidity''-result, which is proven below as
Theorem~\ref{thm:NecessaryCriterionForCoincidenceOfDecompositionSpaces}:
\begin{thm}
\label{thm:EqualDecompositionSpacesOnlyForEquivalentCoverings}Let
$\CalQ,\CalP$ be two almost structured coverings of the open set
$\CalO\subset\R^{\dimension}$ and let $Y\subset\Compl^{I}$ and $Z\subset\Compl^{J}$
be $\CalQ$-regular and $\CalP$-regular, respectively, with $\ell_{0}\left(I\right)\leq Y$
and with $\ell_{0}\left(J\right)\leq Z$. If
\[
\FourierDecompSp{\CalQ}{p_{1}}Y=\FourierDecompSp{\CalP}{p_{2}}Z\quad\text{ for certain }\quad p_{1},p_{2}\in\left(0,\infty\right],
\]
then the following hold:

\begin{enumerate}
\item We have $p_{1}=p_{2}=:p$.
\item In case of $p\neq2$, the coverings $\CalQ$ and $\CalP$ are weakly
equivalent.
\item Furthermore, if $Y=\ell_{w}^{q_{1}}\left(I\right)$ and $Z=\ell_{v}^{q_{2}}\left(J\right)$
for certain $q_{1},q_{2}\in\left(0,\infty\right]$ and certain weights
$w,v$, which are $\CalQ$-moderate and $\CalP$-moderate, respectively,
then the following hold:

\begin{enumerate}
\item We have $\left(p_{1},q_{1}\right)=\left(p_{2},q_{2}\right)=:\left(p,q\right)$.
\item We have $w_{i}\asymp v_{j}$ if $Q_{i}\cap P_{j}\neq\emptyset$.
\item If $\left(p,q\right)\neq\left(2,2\right)$, then $\CalQ$ and $\CalP$
are weakly equivalent.\qedhere
\end{enumerate}
\end{enumerate}
\end{thm}

It is not possible to remove the assumption $\left(p,q\right)\neq\left(2,2\right)$
in the last claim, since an easy application of Plancherel's theorem
shows $\FourierDecompSp{\CalQ}2{\ell_{w}^{2}}=\FourierDecompSp{\CalP}2{\ell_{v}^{2}}$
for any two almost structured coverings $\CalQ,\CalP$, as long as
$w_{i}\asymp v_{j}$ if $Q_{i}\cap P_{j}\neq\emptyset$; see Lemma~\ref{lem:DecompositionSpacesHilbertCase}
for the precise statement. In this sense, the above theorem is best
possible.

Finally, we also show that the decomposition spaces $\FourierDecompSp{\CalQ}{p_{1}}{\ell_{w}^{q_{1}}}$
and $\FourierDecompSp{\CalP}{p_{2}}{\ell_{v}^{q_{2}}}$ are ``fundamentally
different'' if the coverings $\CalQ,\CalP$ cover two \emph{distinct}
sets $\CalO,\CalO'$, and if $\left(p_{1},p_{2},q_{1},q_{2}\right)\neq\left(2,2,2,2\right)$;
for the details, we refer to Theorem~\ref{thm:NoncoincidenceForDifferentOrbits}.

\subsection{Structure of the paper}

\label{subsec:IntroductionStructure}For a quick overview of the structure
of the paper, we refer to the table of contents on page \pageref{toc}.
Here, we give a more detailed overview of the structure and content
of the paper:

We begin our exposition in Section~\ref{sec:Coverings} with a detailed
discussion of different types of coverings. In particular, we recall
the class of \textbf{structured admissible coverings} as introduced
by Borup and Nielsen\cite{BorupNielsenDecomposition}. We generalize
this class of coverings to the two classes of \textbf{semi-structured
admissible coverings} and \textbf{almost structured coverings}. This
last-named class has the advantage of being comparatively nonrestrictive–-in
fact, every specific covering which we will consider is almost structured—while
still possessing all properties needed to obtain well-defined decomposition
spaces and to admit satisfactory criteria for the existence of embeddings.

Additionally, in Section~\ref{sec:Coverings}, we also provide a
brief overview of the possible relations between two coverings, including
\textbf{weak subordinateness} and \textbf{almost subordinateness};
we also discuss the relations between these two concepts. Note that
these notions have already been introduced by Feichtinger and Gröbner
in their seminal work \cite{DecompositionSpaces1,DecompositionSpaces2}
on decomposition spaces. The final—new—notion which we introduce and
explore is that of \textbf{relative moderateness} of two coverings.

\medskip{}

Next, in Section~\ref{sec:DecompositionSpaces}, we begin our study
of decomposition spaces, by giving a detailed formal definition of
these spaces, including a definition of the partitions of unity which
are suitable for defining decomposition spaces. Then we continue by
proving in detail that these spaces are well-defined (i.e., independent
of the chosen partition of unity) and complete.

Strictly speaking, none of these results are new or particularly interesting.
Nevertheless, they seem important enough to include them here, especially
since the (slightly different) definition of decomposition spaces
given in \cite{BorupNielsenDecomposition} yields \emph{incomplete}
spaces in general. Furthermore, well-definedness of the decomposition
space $\DecompSp{\CalQ}pY$ in the quasi-Banach regime $p\in\left(0,1\right)$
is \emph{not} covered by the original paper \cite{DecompositionSpaces1}
on decomposition spaces, since $L^{p}\left(\R^{\dimension}\right)$
is not a Banach-convolution-module over $L^{1}\left(\R^{\dimension}\right)$
for $p<1$. Instead, one has to rely on certain convolution relations
for band-limited $L^{p}$-functions, as treated in \cite{TriebelTheoryOfFunctionSpaces}.
We recall these convolution relations in detail, since they are crucial
not only for the well-definedness of decomposition spaces, but also
for establishing our sufficient criteria for the embeddings between
decomposition spaces.

\medskip{}

In our overview of embedding results (see in particular Theorems \ref{thm:IntroductionFineIntoCoarse}
and \ref{thm:IntroductionCoarseIntoFine}), we saw that the existence
of embeddings between decomposition spaces can be reduced (to a certain
extent) to the existence of embeddings between certain \textbf{nested
sequence spaces}. Thus, a good understanding of these nested sequence
spaces and their embeddings is crucial. We will obtain this understanding
in Section~\ref{sec:NestedSequenceSpaces}, where we study these
spaces in detail. In particular, the results derived in that section
will allow us to pass from the abstract criteria in Theorems \ref{thm:IntroductionFineIntoCoarse}
and \ref{thm:IntroductionCoarseIntoFine} to the concrete ones in
Corollaries \ref{cor:IntroductionFineIntoCoarseSimplified} and \ref{cor:IntroductionCoarseIntoFineSimplified}—a
huge improvement when it comes to the applicability of our criteria.

This kind of simplification is a recurring pattern in the present
paper: In order to simplify proofs, we formulate our sufficient/necessary
criteria in terms of embeddings between nested sequence spaces. Usually,
it would be very painful to verify these embeddings directly. But
using the results of Section~\ref{sec:NestedSequenceSpaces}, one
can then—in a second step—obtain user-friendly versions of these criteria.
Thus, although Section~\ref{sec:NestedSequenceSpaces} is somewhat
technical, it is crucial for the overall goal of the paper.

\medskip{}

In Section~\ref{sec:SufficientConditions}, we properly begin our
investigation of embeddings between decomposition spaces. We start
by estimating the $L^{p}$-norm of a sum $\sum_{i\in I}f_{i}$ of
functions which have ``almost disjoint'' frequency support. Using
this estimate and a ``dual'' version of it, we then derive Theorem~\ref{thm:NoSubordinatenessWithConsiderationOfOverlaps},
which provides a very general sufficient criterion for the existence
of embeddings between two decomposition spaces. However, since the
prerequisites of this theorem are usually hard to verify directly,
we derive several simplified versions of this theorem, which essentially
correspond to Theorems \ref{thm:IntroductionFineIntoCoarse} and \ref{thm:IntroductionCoarseIntoFine}
from above. Finally, we prove Corollary~\ref{cor:MixedSubordinateness},
which applies if we can write $\CalO\cap\CalO'=A\cup B$ in such a
way that $\CalQ$ is almost subordinate to $\CalP$ ``near $A$''
and vice versa ``near $B$''.

\medskip{}

Next, in Section~\ref{sec:NecessaryConditions}, we investigate the
sharpness of the sufficient conditions developed in Section~\ref{sec:SufficientConditions}.
The basic idea is to ``test'' the embedding $\iota:\FourierDecompSp{\CalQ}{p_{1}}Y\hookrightarrow\FourierDecompSp{\CalP}{p_{2}}Z$
using suitably crafted functions. Depending on the precise assumptions
(for instance, $\CalQ$ is almost subordinate to $\CalP$ or vice
versa), the construction varies. Since these constructions—and the
proofs that they yield the desired necessary criteria—are slightly
involved, we refrain from describing them in this introduction. Instead,
we refer the reader to the beginning of Section~\ref{sec:NecessaryConditions}
itself, where the idea is sketched.

The derivation of the necessary criteria itself is divided into a
number of subsections: In subsection \ref{subsec:ElementaryNecessaryConditions},
we begin by proving very elementary conditions: We will see that the
identity map $\iota$ from above can only be bounded if $p_{1}\leq p_{2}$
and if $\left\Vert \delta_{j}\right\Vert _{Z}\lesssim\left\Vert \delta_{i}\right\Vert _{Y}$
for all $i\in I$ and $j\in J$ with $Q_{i}\cap P_{j}\neq\emptyset$.
Although these criteria lack the power of those derived later, they
have the crucial advantage that they impose essentially no assumptions
on (the relation between) the coverings $\CalQ,\CalP$. In particular,
no subordinateness is required. Furthermore, the proofs of these results
allow us to illustrate the general proof techniques for deriving necessary
criteria, without getting bogged down by additional technical difficulties.

Next, in Section~\ref{subsec:CoincidenceOfDecompositionSpaces},
we prove Theorem~\ref{thm:EqualDecompositionSpacesOnlyForEquivalentCoverings};
that is, we show that $\FourierDecompSp{\CalQ}{p_{1}}{\ell_{w}^{q_{1}}}=\FourierDecompSp{\CalP}{p_{2}}{\ell_{v}^{q_{2}}}$
can (essentially) only hold if the ``ingredients'' of the two spaces
are equivalent. The proof of this result is technically much more
demanding than the proofs from the previous subsection; in fact, many
of the proof techniques needed in Subsection~\ref{subsec:ImprovedNecessaryConditions}
are already introduced here.

Subsection~\ref{subsec:ImprovedNecessaryConditions} can be considered
the heart of Section~\ref{sec:NecessaryConditions}: Here, we provide
proofs of Theorems~\ref{thm:IntroductionNecessaryFineIntoCoarse}
and \ref{thm:IntroductionNecessaryCoarseIntoFine}, i.e., we assume
that (a subfamily of) $\CalQ$ is almost subordinate to $\CalP$ (or
vice versa) and we show that the \emph{sufficient} conditions
\begin{equation}
\begin{split}p_{1}\leq p_{2} & \quad\text{and}\quad Y\hookrightarrow Z\left(\left[\,\vphantom{T_{\ell}^{j}}\smash{\ell_{\left|\det T_{i}\right|^{p_{1}^{-1}-p_{2}^{-1}}}^{\LowerExpo{p_{2}}}}\!\left(I_{j}\right)\,\right]_{j\in J}\right)\vphantom{\ell_{\left|\det T_{i}\right|^{p_{1}^{-1}-p_{2}^{-1}}}^{\LowerExpo{p_{2}}}}\\
\text{or}\quad p_{1}\leq p_{2} & \quad\text{and}\quad Y\left(\left[\,\vphantom{T_{j}^{\ell}}\smash{\ell_{u}^{\UpperExpo{p_{1}}}}\!\left(J_{i}\right)\,\right]_{i\in I}\right)\hookrightarrow Z\vphantom{\sum^{T}}
\end{split}
\label{eq:IntroductionSufficientConditionSummary}
\end{equation}
from Theorems \ref{thm:IntroductionFineIntoCoarse} and \ref{thm:IntroductionCoarseIntoFine}
are—slightly modified—also \emph{necessary} for the existence of the
embedding. The main modification needed is that the ``special exponents''
$\LowerExpo{p_{2}}$ or $\UpperExpo{p_{1}}$ are replaced by $p_{2}$
or $p_{1}$, respectively. As noted in the remarks after Theorems~\ref{thm:IntroductionNecessaryFineIntoCoarse}
and \ref{thm:IntroductionNecessaryCoarseIntoFine}, this yields a
\emph{complete characterization} of the existence of the desired embedding
for $p_{2}\in\left(0,2\right]$, or for $p_{1}\in\left[2,\infty\right]$,
respectively.

Of course, the difference between the exponents $\LowerExpo{p_{2}}$
or $\UpperExpo{p_{1}}$ for the sufficient conditions and the exponents
$p_{2}$ or $p_{1}$ for the necessary conditions is unsatisfactory.
In general, I doubt that this gap can be closed; but under the assumption
$p_{1}=p_{2}$, we will see in Subsection \ref{subsec:NecessaryForP1EqualP2}
that this gap can at least be \emph{narrowed}: Note that we always
have $\LowerExpo{p_{2}}\leq2$ and $\UpperExpo{p_{1}}\geq2$. In Subsection
\ref{subsec:NecessaryForP1EqualP2}, we will show —for $p_{1}=p_{2}$—that
the sufficient conditions (\ref{eq:IntroductionSufficientConditionSummary})
are necessary for the existence of the embedding if $\LowerExpo{p_{2}}$
and $\UpperExpo{p_{1}}$ are replaced by $2$. Thus, although our
necessary criteria are not able to reach the exponents $\LowerExpo{p_{2}}$
and $\UpperExpo{p_{1}}$ in general,  we see that no sufficient criterion
similar to condition~(\ref{eq:IntroductionSufficientConditionSummary})
will ever be able to replace $\LowerExpo{p_{2}}$ by an exponent larger
than $2$, or $\UpperExpo{p_{1}}$ by an exponent smaller than $2$.

Finally, in Subsection \ref{subsec:RelativelyModerateCase}, we make
the additional (rather restrictive) assumption that $\CalQ$ is relatively
$\CalP$-moderate (or vice versa). Under this assumption, we achieve
two important goals simultaneously:
\begin{itemize}[leftmargin=0.7cm]
\item We close the offending gap between $\LowerExpo{p_{2}}$ and $p_{2}$
or between $\UpperExpo{p_{1}}$ and $p_{1}$, respectively, thereby
obtaining a \emph{complete characterization} for \emph{all} values
of $p_{2},p_{1}$.\vspace{0.15cm}
\item We are able to simplify our criteria even further, to the point where
only the finiteness of a suitable $\ell^{q}$ norm of a \emph{single}
sequence needs to be checked, see Theorem~\ref{thm:IntroductionNecessaryFineIntoCoarseRelativelyModerate}
and equation~(\ref{eq:IntroductionQIntoPRelativelyModerateCondition})
or Theorem~\ref{thm:IntroductionNecessaryCoarseIntoFineRelativelyModerate}
and equation~(\ref{eq:IntroductionPIntoQRelativelyModerateCondition}).
\end{itemize}
\medskip{}

In view of the large number of sufficient or necessary criteria for
embeddings between decomposition spaces, we summarize our results
in Section~\ref{sec:SummaryOfEmbeddingResults}. In contrast to the
preceding sections, where we strove for maximal generality, our aim
in this section is ease of applicability, even if this makes our results
slightly less general.

\emph{Readers who are mainly interested in applying our embedding
results are thus encouraged to skip directly to }Section~\emph{\ref{sec:SummaryOfEmbeddingResults}},
possibly after familiarizing themselves with the basic concepts of
decomposition spaces outlined in Sections \ref{sec:Coverings}–\ref{sec:DecompositionSpaces}.
In addition to providing a summary of our embedding results, Section~\ref{sec:SummaryOfEmbeddingResults}
also contains a rough ``user's guide'', which should help in deciding
for the correct theorem to apply and in verifying the prerequisites
of the respective theorem.

\medskip{}

In addition to this user's guide, the reader might also want to skim
Section~\ref{sec:Applications} in order to get a feeling for how
our results can be applied. In that section, we \emph{completely characterize}
the existence of the embeddings $\AlphaModSpace{p_{1}}{q_{1}}{s_{1}}{\alpha_{1}}\left(\R^{\dimension}\right)\hookrightarrow\AlphaModSpace{p_{2}}{q_{2}}{s_{2}}{\alpha_{2}}\left(\R^{\dimension}\right)$
between different $\alpha$-modulation spaces, thereby indicating
the power and ease of use of our results. Note that these results
greatly extend the state of the art: Before the present paper, the
most comprehensive knowledge about embeddings of $\alpha$-modulation
spaces had been obtained by Han and Wang \cite{HanWangAlphaModulationEmbeddings};
but they only considered the cases $\alpha_{1}\neq\alpha_{2}$ with
$\left(p_{1},q_{1}\right)=\left(p_{2},q_{2}\right)$, or $\left(p_{1},q_{1}\right)\neq\left(p_{2},q_{2}\right)$
with $\alpha_{1}=\alpha_{2}$, whereas our criteria apply to arbitrary
choices of $\alpha_{1},\alpha_{2},p_{1},p_{2},q_{1},q_{2}$. For completeness,
we mention the paper \cite{GuoAlphaModulationEmbeddingCharacterization}
in which the authors independently derive the same characterization
of embeddings for $\alpha$-modulation spaces as in the present paper.
The preprint of the paper \cite{GuoAlphaModulationEmbeddingCharacterization}
appeared on the arXiv almost simultaneously with the first preprint
of the present paper.

The embeddings between $\alpha$-modulation spaces are particularly
straightforward, since the associated covering $\CalQ^{\left(\alpha_{1}\right)}$
is almost subordinate to and relatively moderate with respect to the
covering $\CalQ^{\left(\alpha_{2}\right)}$ if $\alpha_{1}\leq\alpha_{2}$.
As an example where the two coverings in question are more ``incompatible'',
we also consider embeddings between homogeneous and inhomogeneous
Besov spaces, thereby extending previous results of Triebel \cite{TriebelTheoryOfFunctionSpaces2}.

\medskip{}

The role of Section~\ref{sec:DecompositionSpacesAsSpacesOfTemperedDistributions}
is somewhat special: As we mentioned above—and as we will see in Example~\ref{exa:BorupNielsenDecompositionSpaceIncomplete}—using
the space of tempered distributions as the ``reservoir'' for defining
decomposition spaces can lead to \emph{incomplete spaces}, even in
case of $\bigcup_{i\in I}Q_{i}=\CalO=\R^{\dimension}$. Therefore,
we used the reservoirs $\DistributionSpace{\CalO}$ and $Z'\left(\CalO\right)=\Fourier^{-1}\left(\DistributionSpace{\CalO}\right)$
to define the spaces $\FourierDecompSp{\CalQ}pY$ and $\DecompSp{\CalQ}pY$,
respectively.

But in order to compare our decomposition spaces to classical function
spaces like Besov- and modulation spaces, it is preferable to have
a criterion which ensures that $\DecompSp{\CalQ}pY\hookrightarrow\Schwartz'\left(\R^{\dimension}\right)$.
Hence, in Section~\ref{sec:DecompositionSpacesAsSpacesOfTemperedDistributions},
we derive \emph{sufficient} criteria for this to be the case.

\section{Different classes of coverings and their relations}

\label{sec:Coverings}In this section, we recall the definitions of
various existing classes of coverings and introduce two new types
of coverings, the \textbf{semi-structured admissible coverings} and
the \textbf{almost structured coverings}. The first of these two classes
will be the main class of coverings that we will consider in the remainder
of this paper.

Further, for later use, it will be convenient to have a clear terminology
for describing the relations between two coverings, for instance for
the fact that one is finer than the other. This terminology and certain
auxiliary results are also developed in the present section.

\subsection{Admissible and (semi/almost)-structured coverings}

We begin with the basic concept of \textbf{admissibility} of a covering,
which goes back to the original inception of (general) decomposition
spaces by Feichtinger and Gröbner in \cite{DecompositionSpaces1}.
Admissibility is based on the notion of \textbf{neighbors} of a set
of a covering, which we introduce first.
\begin{defn}
\label{defn:IndexCluster}(cf.\@ \cite[Definition 2.3]{DecompositionSpaces1})

Let $X\neq\emptyset$ be a set and assume that $\CalQ=\left(Q_{i}\right)_{i\in I}$
is a family of subsets of $X$. For a subset $J\subset I$, we define
the \textbf{(index)-cluster} of $J$ (also called the \textbf{set
of neighbors} of $J$) as
\[
J^{\ast}:=\left\{ i\in I\with\exists j\in J:\,Q_{i}\cap Q_{j}\neq\emptyset\right\} .
\]
Inductively, we set $J^{0\ast}:=J$ and $J^{\left(n+1\right)\ast}:=\left(J^{n\ast}\right)^{\ast}$
for $n\in\N_{0}$. We also set $i^{k\ast}:=\left\{ i\right\} ^{k\ast}$
for $i\in I$ and $k\in\N_{0}$. Furthermore, for any subset $J\subset I$,
we define
\[
Q_{J}:=\bigcup_{j\in J}Q_{j}
\]
and we also introduce the shortcuts $Q_{i}^{k\ast}:=Q_{i^{k\ast}}$
and $Q_{i}^{\ast}:=Q_{i^{\ast}}$ for $i\in I$ and $k\in\N_{0}$.

If $\left(\varphi_{i}\right)_{i\in I}$ is a family of functions with
values in some common topological vector space $Y$, we define
\[
\varphi_{J}:=\sum_{i\in J}\varphi_{i}
\]
for every set $J\subset I$ for which the series converges pointwise.
For brevity, we also set $\varphi_{i}^{k\ast}:=\varphi_{i^{k\ast}}$
and $\varphi_{i}^{\ast}:=\varphi_{i^{\ast}}$ for $i\in I$ and $k\in\N_{0}$.
\end{defn}

\begin{rem}
\label{rem:QChainCharacterizationOfCluster}

\begin{enumerate}[leftmargin=0.75cm]
\item If we are considering two different families $\CalQ=\left(Q_{i}\right)_{i\in I}$
and $\CalP=\left(P_{j}\right)_{j\in J}$ with (possibly) nonempty
intersection $I\cap J\neq\emptyset$, we will use notation similar
to $i^{\ast_{\mathcal{Q}}}$ and $i^{\ast_{\mathcal{P}}}$ for $i\in I$
or $i\in J$ to indicate the family which is used to form the respective
cluster.
\item By induction on $n\in\mathbb{N}_{0}$, it is easy to see for $i,j\in I$
that $i\in j^{n\ast}$ is equivalent to $j\in i^{n\ast}$ and that
this holds if and only if there is a \textbf{$\CalQ$-chain} (cf.\@
\cite[Definition 2.3]{DecompositionSpaces1}) of length $n$ from
$j$ to $i$. Here, a finite sequence $\left(i_{\ell}\right)_{\ell=0,\dots,n}$
in $I$ is called a $\CalQ$-chain of length $n$ from $i_{0}$ to
$i_{n}$ if $Q_{i_{\ell-1}}\cap Q_{i_{\ell}}\neq\emptyset$ holds
for all $\ell\in\underline{n}$.
\end{enumerate}
For further properties of admissible coverings and of the cluster
sets $J^{\ast}$, see \cite[Lemma 2.1]{DecompositionSpaces1}.
\end{rem}

Now, we can define the notion of an admissible covering and of moderate
weights.
\begin{defn}
\label{defn:AdmissibleCoveringModerateWeight}(cf.\@ \cite[Definitions 2.1, 3.1 and 3.2]{DecompositionSpaces1}
and \cite[Definition 6]{BorupNielsenDecomposition})

Let $X\neq\emptyset$ be a set and assume that $\CalQ=\left(Q_{i}\right)_{i\in I}$
is a family of subsets of $X$. We say that $\mathcal{Q}$ is an \textbf{admissible
covering} of $X$, if

\begin{enumerate}
\item $\CalQ$ is a covering of $X$, i.e.\@ $X=\bigcup_{i\in I}Q_{i}$
and
\item $Q_{i}\neq\emptyset$ for all $i\in I$, and
\item the constant $N_{\CalQ}:=\sup_{i\in I}\left|i^{\ast}\right|$ is finite.
\end{enumerate}
Let $u=\left(u_{i}\right)_{i\in I}$ be a sequence of positive numbers
$u_{i}>0$. We say that $u$ is a \textbf{$\CalQ$-moderate weight}
if the following expression (then a constant) is finite:
\[
C_{u,\CalQ}:=\sup_{i\in I}\:\sup_{\ell\in i^{\ast}}\:\frac{u_{i}}{u_{\ell}}\:.\qedhere
\]
\end{defn}

\begin{rem*}
It is important to observe that the admissibility condition above
is dependent on the way in which the covering is indexed using the
set $I$. For example, if one set is ``repeated infinitely often'',
i.e.\@ if $\emptyset\neq Q=Q_{i}$ holds for all $i\in I_{0}\subset I$,
where $I_{0}$ is infinite, then this implies $i^{\ast}\supset I_{0}$
for all $i\in I_{0}$, so that $\CalQ$ is \emph{not} admissible.

Furthermore, it is worth noting that every admissible covering is
\emph{of finite height}. This means that the cardinality $\left|\left\{ i\in I\with x\in Q_{i}\right\} \right|$
is bounded uniformly with respect to $x\in X$. More precisely, we
have $\left|\left\{ i\in I\with x\in Q_{i}\right\} \right|\leq N_{\CalQ}$
for all $x\in X$. Equivalently, the function $\sum_{i\in I}\Indicator_{Q_{i}}$
is bounded by $N_{\CalQ}$.
\end{rem*}
As we will see in Section~\ref{sec:DecompositionSpaces}, admissibility
of the covering $\mathcal{Q}$ is the most basic requirement needed
to ensure well-definedness of the decomposition spaces $\FourierDecompSp{\CalQ}pY$.
By well-definedness, we mean that the resulting space does not depend
on the chosen partition of unity $\Phi=\left(\varphi_{i}\right)_{i\in I}$.
Of course, this can only hold under suitable assumptions on $\Phi$.
More precisely, we will assume that $\left\Vert \mathcal{F}^{-1}\varphi_{i}\right\Vert _{L^{1}}$
is uniformly bounded, which ensures well-definedness (at least for
$p\geq1$). Thus, one might ask whether such a partition of unity
$\Phi$ always exists.

In the following lemma, we develop a first set of \emph{necessary}
conditions that the covering $\CalQ$ has to fulfill if such a partition
of unity exists. Note though that we just assume existence of a partition
of unity consisting of $C_{c}$ functions, i.e., we do not yet assume
finiteness of $\sup_{i\in I}\left\Vert \mathcal{F}^{-1}\varphi_{i}\right\Vert _{L^{1}}$.
\begin{lem}
\label{lem:PartitionCoveringNecessary}Let $\CalQ=\left(Q_{i}\right)_{i\in I}$
be an \emph{admissible} covering of an open subset $\CalO\subset X$
of a topological space $X$. We say that $\CalQ$ \textbf{admits a
partition of unity} if there are functions $\varphi_{i}:X\to\Compl$,
$i\in I$, with

\begin{enumerate}
\item $\varphi_{i}\in C_{c}\left(\CalO\right)$,
\item $\varphi_{i}\equiv0$ outside of $Q_{i}$,
\item $\sum_{i\in I}\varphi_{i}\equiv1$ on $\CalO$.
\end{enumerate}
In this case, the following hold:

\begin{enumerate}
\item For any subset $M\subset I$, we have
\begin{equation}
\varphi_{M}\equiv1\quad\text{ on }Q_{i}\text{ for each }i\in I\text{ with }i^{\ast}\subset M.\label{eq:ClusteredPartitionYields1}
\end{equation}
\item We have $\overline{Q_{i}}\subset\CalO$ for each $i\in I$. Furthermore,
$\overline{Q_{i}}$ is compact.
\item The family of topological interiors $\CalQ^{\circ}=\left(Q_{i}^{\circ}\right)_{i\in I}$
covers $\CalO$.
\item The family $\left(\overline{Q_{i}}\right)_{i\in I}$ is locally finite
in $\CalO$ (but not necessarily in $X$). In particular, $\left(\varphi_{i}\right)_{i\in I}$
is a locally finite partition of unity, in the sense that $\left(\supp\varphi_{i}\right)_{i\in I}$
is locally finite in $\CalO$.\qedhere
\end{enumerate}
\end{lem}

\begin{rem*}
Observe that $\varphi_{M}=\sum_{i\in M}\varphi_{i}$ is well-defined
even if $M\subset I$ is an infinite set, since at every point $x\in X$,
at most $N_{\CalQ}$ summands do not vanish.
\end{rem*}
\begin{proof}
Ad (1): Since we have $\sum_{j\in I}\varphi_{j}\equiv1$ on $\CalO\supset Q_{i}$,
it suffices to show $\varphi_{j}\equiv0$ on $Q_{i}$ for $j\notin M$.
But if $\varphi_{j}\left(x\right)\neq0$ holds for some $x\in Q_{i}$,
we get $x\in Q_{j}$ (since $\varphi_{j}$ vanishes outside of $Q_{j}$)
and hence $x\in Q_{i}\cap Q_{j}\neq\emptyset$, i.e.\@ $j\in i^{\ast}\subset M$.
This establishes equation~(\ref{eq:ClusteredPartitionYields1}).

\medskip{}

Ad (2): By equation~(\ref{eq:ClusteredPartitionYields1}), we have
$\sum_{j\in i^{\ast}}\varphi_{j}=\varphi_{i^{\ast}}\equiv1$ on $Q_{i}$.
But $\varphi_{i^{\ast}}\in C_{c}\left(\CalO\right)$, since $i^{\ast}$
is finite, because $\CalQ$ is admissible. Thus, $Q_{i}\subset\supp\varphi_{i^{\ast}}\subset\CalO$,
whence $\overline{Q_{i}}\subset\supp\varphi_{i^{\ast}}\subset\CalO$
is compact.

\medskip{}

Ad (3): By assumption, $\varphi_{i}^{-1}\left(\Compl^{\ast}\right)\subset Q_{i}$
for all $i\in I$ and thus $\varphi_{i}^{-1}\left(\Compl^{\ast}\right)\subset Q_{i}^{\circ}$
by continuity of $\varphi_{i}$. But we have $\sum_{i\in I}\varphi_{i}\equiv1$
on $\CalO$ and hence
\[
\CalO\subset\bigcup_{i\in I}\varphi_{i}^{-1}\left(\Compl^{\ast}\right)\subset\bigcup_{i\in I}Q_{i}^{\circ}\subset\bigcup_{i\in I}Q_{i}=\CalO.
\]

\medskip{}

Ad (4): Let $x\in\CalO$ be arbitrary. Then there is some $i\in I$
with $\varphi_{i}\left(x\right)\neq0$, since $\sum_{i\in I}\varphi_{i}\equiv1$
on $\CalO$. Hence, $U_{x}:=\varphi_{i}^{-1}\left(\Compl^{\ast}\right)\subset Q_{i}^{\circ}\subset\CalO$
is an open neighborhood of $x$. Let $I_{x}:=i^{\ast}$. For any index
$j\in I$ with $U_{x}\cap\overline{Q_{j}}\neq\emptyset$, we have
$\emptyset\neq U_{x}\cap\overline{Q_{j}}\subset Q_{i}^{\circ}\cap\overline{Q_{j}}$.
But this entails $Q_{i}^{\circ}\cap Q_{j}\neq\emptyset$, since otherwise
$Q_{j}$ would be contained in the closed(!) set $\left(Q_{i}^{\circ}\right)^{c}$.
Thus, $\emptyset\neq Q_{i}^{\circ}\cap Q_{j}\subset Q_{i}\cap Q_{j}$
and hence $j\in i^{\ast}=I_{x}$.

Since $I_{x}=i^{\ast}$ is finite by admissibility of $\CalQ$ and
since $U_{x}$ is an open neighborhood of $x$, we see that $\left(\,\overline{Q_{j}}\,\right)_{j\in I}$
is a locally finite family in $\CalO$. Since $\varphi_{i}$ vanishes
outside of $Q_{i}$, we have $\supp\varphi_{i}\subset\overline{Q_{i}}$
for all $i\in I$, so that the last claim follows trivially.
\end{proof}
As a complement to these \emph{necessary} conditions, we would like
to have \emph{sufficient} conditions which ensure existence of a suitable
partition of unity $\left(\varphi_{i}\right)_{i\in I}$ subordinate
to a covering $\CalQ$. To this end, we introduce several more restrictive
classes of coverings: \textbf{semi-structured coverings}, \textbf{almost
structured coverings}, and \textbf{structured coverings}. The general
idea of these coverings is—essentially—that each set $Q_{i}$ should
be an affine image $Q_{i}=T_{i}Q+b_{i}$ of a \emph{fixed} set $Q\subset\R^{\dimension}$.
This is the main requirement for a \emph{structured} covering. For
many practical coverings, however, this assumptions is somewhat strict:
For example, the usual dyadic covering consisting of the dyadic annuli
$\left(B_{2^{n+1}}\left(0\right)\setminus\overline{B_{2^{n-1}}\left(0\right)}\right)_{n\in\mathbb{N}}$
and one ball $B_{4}\left(0\right)$ covering the low frequencies is
\emph{not} of this form, since $B_{4}\left(0\right)$ is convex, while
the dyadic annuli are not. Thus, the notions of almost structured
coverings and semi-structured coverings slightly relax the assumption
$Q_{i}=T_{i}Q+b_{i}$ for all $i\in I$, by allowing the set $Q$
to vary—in a controlled way—with $i\in I$.

As we will see later (see Theorem~\ref{thm:AlmostStructuredAdmissibleAdmitsBAPU}),
one can always construct a suitable partition $\left(\varphi_{i}\right)_{i\in I}$
subordinate to any almost structured covering.

We remark that the notion of structured admissible coverings was first
introduced (for the case $\CalO=\R^{\dimension}$) by Borup and Nielsen
in \cite[Definition 7]{BorupNielsenDecomposition} and then slightly
generalized by Führ and myself in \cite[Definition 13]{FuehrVoigtlaenderCoorbitSpacesAsDecompositionSpaces}
to the definition presented here.
\begin{defn}
\label{defn:DifferentTypesOfCoverings}Let $\emptyset\neq\CalO\subset\R^{\dimension}$
be open and let $I\neq\emptyset$ be an index set. We say that a family
$\CalQ=\left(Q_{i}\right)_{i\in I}$ of subsets $Q_{i}\subset\CalO$
is a \textbf{semi-structured covering} of $\CalO$ if for each $i\in I$,
there are $T_{i}\in\GL\left(\R^{\dimension}\right)$ and $b_{i}\in\R^{\dimension}$
and an \emph{open} subset $Q_{i}'\subset\R^{\dimension}$ with
\[
Q_{i}=T_{i}Q_{i}'+b_{i}
\]
and such that the following properties are fulfilled:

\begin{enumerate}
\item $\CalQ$ is an admissible covering of $\CalO$.
\item The sets $\left(Q_{i}'\right)_{i\in I}$ are uniformly bounded, i.e.\@
the following expression (then a constant) is finite:
\[
R_{\CalQ}:=\sup_{i\in I}\:\sup_{x\in Q_{i}'}\:\left|x\right|\,.
\]
\item For ``neighboring'' indices $i,\ell\in I$, the transformations
$T_{i}\mybullet+b_{i}$ and $T_{\ell}\mybullet+b_{\ell}$ are ``uniformly
compatible,'' i.e.\@ the following expression (then a constant)
is finite:
\[
C_{\CalQ}:=\sup_{i\in I}\:\sup_{\ell\in i^{\ast}}\:\left\Vert T_{i}^{-1}T_{\ell}\right\Vert \,.
\]
\end{enumerate}
If these conditions are satisfied, the family $\left(T_{i},Q_{i}',b_{i}\right)_{i\in I}$
is called a \textbf{(semi-structured) parametrization} of the covering
$\CalQ$.

The semi-structured parametrization $\left(T_{i},Q_{i}',b_{i}\right)_{i\in I}$
(and by abuse of language, the covering $\CalQ$) is called \textbf{tight}
if we additionally have the following:

\begin{enumerate}[resume]
\item There is some $\varepsilon>0$ and for each $i\in I$ some $c_{i}\in\R^{\dimension}$
such that $B_{\varepsilon}\left(c_{i}\right)\subset Q_{i}'$.
\end{enumerate}
The semi-structured parametrization $\left(T_{i},Q_{i}',b_{i}\right)_{i\in I}$
(and again, the covering $\CalQ$) is called \textbf{almost structured}
if for each $i\in I$, there is an \emph{open} set $Q_{i}''\subset\R^{\dimension}$
such that the following hold:

\begin{enumerate}
\item We have $\overline{Q_{i}''}\subset Q_{i}'$ for all $i\in I$.
\item The family $\left(T_{i}Q_{i}''+b_{i}\right)_{i\in I}$ is an admissible
covering of $\CalO$.
\item The sets $\left\{ Q_{i}'\with i\in I\right\} $ and $\left\{ Q_{i}''\with i\in I\right\} $
are finite.
\end{enumerate}
Finally, an almost structured parametrization (or an almost structured
covering) as above is called \textbf{structured }if the set $\left\{ Q_{i}'\with i\in I\right\} $
only consists of one element.
\end{defn}

\begin{rem}
\label{rem:CoveringTypesRemark}

\begin{enumerate}[leftmargin=0.6cm]
\item The adjective ``tight'' from above should be understood as the
opposite of ``loose,'' in the following sense: Given a semi-structured
parametrization $\left(T_{i},Q_{i}',b_{i}\right)_{i\in I}$ of a covering
$\CalQ$, and any sequence $\left(\alpha_{i}\right)_{i\in I}$ of
positive numbers bounded from below, we get an alternative parametrization
$\left(S_{i},Q_{i}'',b_{i}\right)_{i\in I}$ of $\CalQ$ by setting
$S_{i}:=\alpha_{i}T_{i}$ and $Q_{i}'':=\alpha_{i}^{-1}\cdot Q_{i}'$.
But if the sequence $\left(\alpha_{i}\right)_{i\in I}$ is unbounded,
then (by uniform boundedness of the sets $Q_{i}'$), the sets $Q_{i}''$
get arbitrarily small, which is precisely what is excluded by the
definition of a \emph{tight} parametrization. Thus, a tight parametrization
excludes ``loose'' parametrizations like $\left(S_{i},Q_{i}'',b_{i}\right)_{i\in I}$
where the $S_{i}$ are ``unnecessarily large'', while the $Q_{i}''$
are ``unnecessarily small''. A consequence of having a tight parametrization
is that $\lambda\left(Q_{i}\right)\asymp\left|\det T_{i}\right|$,
while for non-tight parametrizations, one only has $\lambda\left(Q_{i}\right)\lesssim\left|\det T_{i}\right|$,
see Lemma~\ref{cor:SemiStructuredDifferenceSetsMeasureEstimate}.
\item In the remainder of the paper, the phrase ``let $\CalQ=\left(Q_{i}\right)_{i\in I}=\left(T_{i}Q_{i}'+b_{i}\right)_{i\in I}$
be a {[}tight{]} semi-structured covering of $\CalO$'' will always
implicitly entail that $\left(T_{i},Q_{i}',b_{i}\right)_{i\in I}$
is a {[}tight{]} parametrization of $\CalQ$. The same holds for ``almost
structured'' or ``structured'' instead of ``semi-structured''.
\item It is interesting to note that any two \emph{tight} parametrizations
$\left(T_{i},Q_{i}',b_{i}\right)_{i\in I}$ and $\left(S_{i},Q_{i}'',c_{i}\right)_{i\in I}$
of a common covering $\CalQ$ are equivalent, in the sense that
\[
\sup_{i\in I}\left[\left\Vert T_{i}^{-1}S_{i}\right\Vert +\left\Vert S_{i}^{-1}T_{i}\right\Vert \right]<\infty\,.
\]
Indeed, by symmetry it suffices to consider $\left\Vert S_{i}^{-1}T_{i}\right\Vert $.
By assumption, $T_{i}Q_{i}'+b_{i}=Q_{i}=S_{i}Q_{i}''+c_{i}$, and
thus $S_{i}^{-1}T_{i}Q_{i}'\subset Q_{i}''+e_{i}$, with $e_{i}:=S_{i}^{-1}\left(c_{i}-b_{i}\right)$.
But there is some $\varepsilon>0$ with $B_{\varepsilon}\left(c_{i}\right)\subset Q_{i}'$
for all $i\in I$, which easily implies $B_{2\varepsilon}\left(0\right)\subset Q_{i}'-Q_{i}'$.
Since we also have $Q_{i}''\subset B_{R}\left(0\right)$ for all $i\in I$
and some $R>0$, we finally see
\[
2\varepsilon\cdot S_{i}^{-1}T_{i}\left(B_{1}\left(0\right)\right)\subset S_{i}^{-1}T_{i}Q_{i}'-S_{i}^{-1}T_{i}Q_{i}'\subset Q_{i}''-Q_{i}''\subset B_{2R}\left(0\right),
\]
so that $\left\Vert S_{i}^{-1}T_{i}\right\Vert \leq R/\varepsilon$
for all $i\in I$.\vspace{0.1cm}
\item If $\CalQ=\left(T_{i}Q_{i}'+b_{i}\right)_{i\in I}$ is a \emph{tight}
semi-structured covering, we set
\[
\varepsilon_{\CalQ}:=\sup\left\{ \varepsilon>0\with\forall\,i\in I\,\exists\,c_{i}\in\R^{\dimension}:\,B_{\varepsilon}\left(c_{i}\right)\subset Q_{i}'\right\} .
\]
As we will see below, this supremum is always attained.
\item Strictly speaking, the constants $R_{\CalQ}$, $C_{\CalQ}$ and $\varepsilon_{\CalQ}$
depend on the choice of the parametrization $\left(T_{i},Q_{i}',b_{i}\right)_{i\in I}$
of $\CalQ$. As above, we will usually suppress this dependence.
\item Every almost structured covering is tight, since the family $\left\{ Q_{i}'\with i\in I\right\} $
of open sets is finite for such a covering $\CalQ=\left(T_{i}Q_{i}'+b_{i}\right)_{i\in I}$.
\item For brevity, we will use the following convention: If a theorem states
that for a given {[}tight{]} semi-structured (or almost-structured)
covering $\CalQ$, there is a constant $C=C\!\left(\CalQ,\dots\right)>0$
{[}or $C=C\!\left(\CalQ,\varepsilon_{\CalQ},\dots\right)${]} with
a certain property (depending on $\CalQ$), then the following is
meant: Given $N,R,K>0$ {[}and $\varepsilon>0${]}, there is a constant
$C_{0}=C_{0}\left(N,R,K,\dots\right)>0$ {[}respectively a constant
$C_{0}=C_{0}\left(N,R,K,\varepsilon,\dots\right)>0${]} such that
this particular property holds—with $C_{0}$ instead of $C$—for every
{[}tight{]} semi-structured (or almost-structured) covering $\CalQ$
with $N_{\CalQ}\leq N$, $R_{\CalQ}\leq R$, $C_{\CalQ}\leq K$ {[}and
$\varepsilon_{\CalQ}\geq\varepsilon${]}.

Note that this means that the constant $C$ may \emph{not} depend
on other properties of the covering $\CalQ$ than on the quantities
$N_{\CalQ},R_{\CalQ},C_{\CalQ}$ {[}and $\varepsilon_{\CalQ}${]},
except possibly for those appearing in the ``$\dots$'' part of
$C\!\left(\CalQ,\dots\right)$.

As an example, note that we have $Q_{i}\subset T_{i}\left(B_{R_{\CalQ}}\left(0\right)\right)+b_{i}$
for all $i\in I$. Furthermore, for $\ell\in i^{\ast}$, \emph{Hadamard's
inequality} $\left|\det A\right|\leq\left\Vert A\right\Vert ^{\dimension}$
for $A\in\R^{\dimension\times\dimension}$ (see \cite[Section 75]{RieszFunctionalAnalysis})
yields
\[
\left|\det T_{\ell}\right|=\left|\det T_{i}\right|\cdot\left|\det\left(T_{i}^{-1}T_{\ell}\right)\right|\leq\left|\det T_{i}\right|\cdot\left\Vert T_{i}^{-1}T_{\ell}\right\Vert ^{\dimension}\leq C_{\CalQ}^{\dimension}\cdot\left|\det T_{i}\right|\,,
\]
and hence
\[
\quad\lambda\left(Q_{\ell}\right)\leq\left|\det T_{\ell}\right|\cdot\lambda\left(B_{R_{\CalQ}}\left(0\right)\right)\leq\lambda\left(B_{1}\left(0\right)\right)\cdot R_{\CalQ}^{\dimension}\cdot\left|\det T_{\ell}\right|\leq\lambda\left(B_{1}\left(0\right)\right)\cdot\left(C_{\CalQ}R_{\CalQ}\right)^{\dimension}\cdot\left|\det T_{i}\right|.
\]
Using the estimate $\left|i^{\ast}\right|\leq N_{\CalQ}$, we finally
get
\[
\lambda\left(Q_{i^{\ast}}\right)\leq\sum_{\ell\in i^{\ast}}\lambda\left(Q_{\ell}\right)\leq\lambda\left(B_{1}\left(0\right)\right)\cdot N_{\CalQ}\cdot\left(C_{\CalQ}R_{\CalQ}\right)^{\dimension}\cdot\left|\det T_{i}\right|\,.
\]
With the above convention, we have thus shown for any semi-structured
covering $\CalQ=\left(T_{i}Q_{i}'+b_{i}\right)_{i\in I}$ that there
is a constant $C=C\!\left(\CalQ,\dimension\right)>0$ such that $\lambda\left(Q_{i^{\ast}}\right)\leq C\cdot\left|\det T_{i}\right|$
for all $i\in I$.
\item \label{enu:CoveringEpsilonIsAttained}Finally, we prove that the supremum
in the definition of $\varepsilon_{\CalQ}$ is attained. To this end,
note that finiteness of $R_{\CalQ}$ implies that each set $Q_{i}'$
is bounded, so that $\overline{Q_{i}'}\subset\R^{\dimension}$ is
compact.

Now, for each fixed $i\in I$ and arbitrary $n\in\N$, let $\varepsilon_{n}:=\left(1-\frac{1}{2n}\right)\varepsilon_{\CalQ}$.
Then there is some $c_{n}\in\R^{\dimension}$ satisfying $c_{n}\in B_{\varepsilon_{n}}\left(c_{n}\right)\subset Q_{i}'\subset\overline{Q_{i}'}$.
For some subsequence, we get $c_{n_{k}}\to c\in\overline{Q_{i}'}$
by compactness of $\overline{Q_{i}'}$.

Finally, let $x\in B_{\varepsilon_{\CalQ}}\left(c\right)$ be arbitrary,
and set $\delta:=\frac{1}{2}\left(\varepsilon_{\CalQ}-\left|x-c\right|\right)>0$.
We have
\[
\left|x-c_{n_{k}}\right|\xrightarrow[k\to\infty]{}\left|x-c\right|<\delta+\left|x-c\right|<\varepsilon_{\CalQ}\:,
\]
and thus $\left|x-c_{n_{k}}\right|<\delta+\left|x-c\right|<\varepsilon_{n_{k}}$
for sufficiently large $k\in\mathbb{N}$. Hence, $x\in B_{\varepsilon_{n_{k}}}\left(c_{n_{k}}\right)\subset Q_{i}'$.
We have thus shown $B_{\varepsilon_{\CalQ}}\left(c\right)\subset Q_{i}'$.\qedhere
\end{enumerate}
\end{rem}

As our next technical result, we show that the ``normalization''
$T_{i}^{-1}\left(Q_{i}-b_{i}\right)\subset B_{R}\left(0\right)$—which
by definition holds for $R=R_{\CalQ}$—remains essentially valid if
the set $Q_{i}$ is replaced by a neighboring set $Q_{j}$, as long
as one is willing to enlarge the ball $B_{R}\left(0\right)$. This
inclusion will become important for establishing well-definedness
of the decomposition space $\FourierDecompSp{\CalQ}pY$ in the quasi-Banach
regime $p\in\left(0,1\right)$. In this case, we will see (see Theorem~\ref{thm:QuasiBanachConvolution})
that the usual estimate $\left\Vert \mathcal{F}^{-1}\left(fg\right)\right\Vert _{L^{p}}\leq\left\Vert \mathcal{F}^{-1}f\right\Vert _{L^{1}}\cdot\left\Vert \mathcal{F}^{-1}g\right\Vert _{L^{p}}$
needs to be replaced by
\[
\left\Vert \mathcal{F}^{-1}\left(fg\right)\right\Vert _{L^{p}}\leq\left[\lambda\left(K-L\right)\right]^{\frac{1}{p}-1}\cdot\left\Vert \mathcal{F}^{-1}f\right\Vert _{L^{p}}\cdot\left\Vert \mathcal{F}^{-1}g\right\Vert _{L^{p}},
\]
where $\supp f\subset K$ and $\supp g\subset L$. Thus, it will be
a convenient consequence of the following lemma that we can estimate
the Lebesgue measure $\lambda\left(\,\overline{Q_{i}}-\overline{Q_{j}}\,\right)$
by a constant multiple of $\left|\det T_{i}\right|$; see Corollary~\ref{cor:SemiStructuredDifferenceSetsMeasureEstimate}.
We remark that the issues related to the constant $\left[\lambda\left(\,\overline{Q_{i}}-\overline{Q_{j}}\,\right)\right]^{\frac{1}{p}-1}$
in contrast to $\left|\det T_{i}\right|^{\frac{1}{p}-1}$ are somewhat
neglected in the treatment of Borup and Nielsen\cite{BorupNielsenDecomposition}.
\begin{lem}
\label{lem:SemiStructuredNormalizationNeighboring}Let $\emptyset\neq\CalO\subset\R^{\dimension}$
be open and let $\CalQ=\left(Q_{i}\right)_{i\in I}=\left(T_{i}Q_{i}'+b_{i}\right)_{i\in I}$
be a semi-structured covering of $\CalO$. Then we have
\[
Q_{j}\subset T_{j}\left(\overline{B_{R}}\left(0\right)\right)+b_{j}\subset T_{i}\left[\overline{B_{\left(2C_{\CalQ}+1\right)^{\ell}R}}\left(0\right)\right]+b_{i}\qquad\forall\:i\in I,\,\ell\in\N_{0},\text{ and }j\in i^{\ell\ast}\,,
\]
as long as $R>0$ is chosen such that $Q_{i}'\subset\overline{B_{R}}\left(0\right)$
holds for all $i\in I$.
\end{lem}

\begin{proof}
Set $K:=C_{\CalQ}$ and let $R>0$ with $\bigcup_{i\in I}Q_{i}'\subset\overline{B_{R}}\left(0\right)$.
We first prove
\begin{equation}
T_{j}\left(\overline{B_{R}}\left(0\right)\right)+b_{j}\subset T_{i}\left(\overline{B_{\left(2K+1\right)R}}\left(0\right)\right)+b_{i}\qquad\forall\:i\in I\text{ and }j\in i^{\ast}\:.\label{eq:SemiStructuredNormalizationInductionStart}
\end{equation}
To see this, note that $j\in i^{\ast}$ implies $\left\Vert \smash{T_{i}^{-1}T_{j}}\right\Vert \leq K$,
as well as
\[
\emptyset\neq Q_{i}\cap Q_{j}\subset\left[T_{i}\left(\overline{B_{R}}\left(0\right)\right)+b_{i}\right]\cap\left[T_{j}\left(\overline{B_{R}}\left(0\right)\right)+b_{j}\right],
\]
which yields $c_{1},c_{2}\in\overline{B_{R}}\left(0\right)$ satisfying
$T_{i}c_{1}+b_{i}=T_{j}c_{2}+b_{j}$. Rearranging results in
\[
b_{j}-b_{i}=T_{i}c_{1}-T_{j}c_{2}=T_{i}\left(c_{1}-T_{i}^{-1}T_{j}c_{2}\right),
\]
with
\[
\left|c_{1}-T_{i}^{-1}T_{j}c_{2}\right|\leq\left|c_{1}\right|+\left\Vert \smash{T_{i}^{-1}}T_{j}\right\Vert \left|c_{2}\right|\leq R+KR=\left(K+1\right)R,
\]
which implies $T_{i}^{-1}\left(b_{j}-b_{i}\right)=c_{1}-T_{i}^{-1}T_{j}c_{2}\in\overline{B_{\left(K+1\right)R}}\left(0\right)$.

Now, let $b\in\overline{B_{R}}\left(0\right)$ be arbitrary and set
$x:=T_{i}^{-1}\left(T_{j}b+\left(b_{j}-b_{i}\right)\right)$. On the
one hand,
\[
\left|x\right|\leq\left|\smash{T_{i}^{-1}}T_{j}b\right|+\left|\smash{T_{i}^{-1}}\left(b_{j}-b_{i}\right)\right|\leq\left\Vert \smash{T_{i}^{-1}}T_{j}\right\Vert \left|b\right|+\left(K+1\right)R\leq\left(2K+1\right)R.
\]
On the other hand, we have $T_{i}x+b_{i}=T_{j}b+\left(b_{j}-b_{i}\right)+b_{i}=T_{j}b+b_{j}$,
and thus
\[
T_{j}b+b_{j}=T_{i}x+b_{i}\in T_{i}\left(\overline{B_{\left(2K+1\right)R}}\left(0\right)\right)+b_{i},
\]
which proves the claimed inclusion~(\ref{eq:SemiStructuredNormalizationInductionStart}).

We can now establish the general claim by induction on $\ell\in\N_{0}$.
The base case $\ell=0$ is a direct consequence of $Q_{i}=T_{i}Q_{i}'+b_{i}\subset T_{i}\left(\overline{B_{R}}\left(0\right)\right)+b_{i}$,
since $i^{0\ast}=\left\{ i\right\} $.

For the induction step, we note that the assumption $j\in i^{\left(\ell+1\right)\ast}$
yields some $k\in i^{\ell\ast}$ with $j\in k^{\ast}$. Let us set
$R':=\left(2K+1\right)R\geq R$, so that $\bigcup_{i\in I}Q_{i}'\subset\overline{B_{R}}\left(0\right)\subset\overline{B_{R'}}\left(0\right)$
holds. By applying the induction hypothesis for $R'$ instead of $R$,
we obtain
\[
T_{k}\left(\overline{B_{R'}}\left(0\right)\right)+b_{k}\subset T_{i}\left(\overline{B_{\left(2K+1\right)^{\ell}R'}}\left(0\right)\right)+b_{i}=T_{i}\left(\overline{B_{\left(2K+1\right)^{\ell+1}R}}\left(0\right)\right)+b_{i}.
\]
But since $j\in k^{\ast}$, we can apply equation~(\ref{eq:SemiStructuredNormalizationInductionStart}),
which results in
\[
Q_{j}\subset T_{j}\left(\overline{B_{R}}\left(0\right)\right)+b_{j}\subset T_{k}\left(\overline{B_{\left(2K+1\right)R}}\left(0\right)\right)+b_{k}=T_{k}\left(\overline{B_{R'}}\left(0\right)\right)+b_{k}\subset T_{i}\left(\overline{B_{\left(2K+1\right)^{\ell+1}R}}\left(0\right)\right)+b_{i}.
\]
This completes the induction step.
\end{proof}
As a corollary, we obtain (a generalization of) the estimate for the
Lebesgue measure of the difference set $\overline{Q_{i}}-\overline{Q_{j}}$
that was announced above.
\begin{cor}
\label{cor:SemiStructuredDifferenceSetsMeasureEstimate}Let $\emptyset\neq\CalO\subset\R^{\dimension}$
be open and let $\CalQ=\left(Q_{i}\right)_{i\in I}=\left(T_{i}Q_{i}'+b_{i}\right)_{i\in I}$
be a semi-structured covering of $\CalO$.

\begin{enumerate}
\item There is a constant $C_{1}=C_{1}\left(\dimension,R_{\CalQ}\right)>0$,
such that
\[
\lambda\left(Q_{i}\right)\leq\lambda\left(\,\overline{Q_{i}}\,\right)\leq C_{1}\cdot\left|\det T_{i}\right|\qquad\forall\,i\in I\,.
\]
\item Conversely, if $\CalQ$ is tight, there is a constant $C_{2}=C_{2}\left(\dimension,\varepsilon_{\CalQ}\right)>0$
satisfying
\[
C_{2}\cdot\left|\det T_{i}\right|\leq\lambda\left(Q_{i}\right)\leq\lambda\left(\,\overline{Q_{i}}\,\right)\qquad\forall\,i\in I\,.
\]
\item Finally, for arbitrary $n\in\N$, there is a constant $C_{3}=C_{3}\left(\CalQ,n,\dimension\right)>0$,
such that 
\[
\max_{j\in i^{n\ast}}\lambda\left(\,\overline{Q_{i}^{n\ast}}-\overline{Q_{j}^{n\ast}}\,\right)\leq C_{3}\cdot\left|\det T_{i}\right|\qquad\forall\,i\in I\,.\qedhere
\]
\end{enumerate}
\end{cor}

\begin{rem*}
For a tight semi-structured covering, we thus get $\lambda\left(Q_{i}\right)\asymp\lambda\left(\,\overline{Q_{i}}\,\right)\asymp\left|\det T_{i}\right|$
uniformly with respect to $i\in I$.
\end{rem*}
\begin{proof}
Let $R\geq R_{\CalQ}$. We first observe $Q_{i}\subset\overline{Q_{i}}=T_{i}\overline{Q_{i}'}+b_{i}\subset T_{i}\left(\,\vphantom{B_{R}}\smash{\overline{B_{R}\left(0\right)}}\,\right)+b_{i}$,
and hence
\[
\lambda\left(Q_{i}\right)\leq\lambda\left(\,\overline{Q_{i}}\,\right)\leq\left|\det T_{i}\right|\cdot\lambda\left(\,\smash{\overline{B_{R}\left(0\right)}}\,\right),
\]
so that we can set $C_{1}:=\lambda\left(\smash{\overline{B_{R}\left(0\right)}}\right)=v_{\dimension}\cdot R^{\dimension}$,
where $v_{\dimension}:=\lambda\left(B_{1}\left(0\right)\right)$ denotes
the measure of the $\dimension$-dimensional Euclidean unit ball.

Conversely, if $\CalQ$ is tight, let $0<\varepsilon\leq\varepsilon_{\CalQ}$.
By part~(\ref{enu:CoveringEpsilonIsAttained}) of Remark~\ref{rem:CoveringTypesRemark},
there is for each $i\in I$ some $c_{i}\in\R^{\dimension}$ with $B_{\varepsilon}\left(c_{i}\right)\subset Q_{i}'$.
Hence, $Q_{i}=T_{i}Q_{i}'+b_{i}\supset T_{i}\left(B_{\varepsilon}\left(c_{i}\right)\right)+b_{i}$,
which yields
\[
\lambda\left(\overline{Q_{i}}\right)\geq\lambda\left(Q_{i}\right)\geq\left|\det T_{i}\right|\cdot\lambda\left(B_{\varepsilon}\left(c_{i}\right)\right)=v_{\dimension}\cdot\varepsilon^{\dimension}\cdot\left|\det T_{i}\right|,
\]
so that we can set $C_{2}:=v_{\dimension}\cdot\varepsilon^{\dimension}$.

It remains to establish the estimate regarding the difference set
$\overline{Q_{i}^{n\ast}}-\overline{Q_{j}^{n\ast}}$ for $j\in i^{n\ast}$.
Note that $j^{n\ast}\subset i^{2n\ast}$. Hence, $\overline{Q_{i}^{n\ast}}\subset\overline{Q_{i}^{2n\ast}}$
and $\overline{Q_{j}^{n\ast}}\subset\overline{Q_{i}^{2n\ast}}$.

Now, let $K\geq C_{\CalQ}$ and $R\geq R_{\CalQ}$. Then, Lemma~\ref{lem:SemiStructuredNormalizationNeighboring}
yields for any $\ell\in i^{2n\ast}$ the inclusion
\[
Q_{\ell}\subset T_{i}\left[\overline{B_{\left(2K+1\right)^{2n}R}}\left(0\right)\right]+b_{i}.
\]
Since $\ell\in i^{2n\ast}$ was arbitrary, we get $\overline{Q_{i}^{2n\ast}}\subset T_{i}\left[\overline{B_{\left(2K+1\right)^{2n}R}}\left(0\right)\right]+b_{i}$,
and hence
\begin{align*}
\overline{Q_{i}^{n\ast}}-\overline{Q_{j}^{n\ast}}\subset\overline{Q_{i}^{2n\ast}}-\overline{Q_{i}^{2n\ast}} & \subset\left(T_{i}\left[\overline{B_{\left(2K+1\right)^{2n}R}}\left(0\right)\right]+b_{i}\right)-\left(T_{i}\left[\overline{B_{\left(2K+1\right)^{2n}R}}\left(0\right)\right]+b_{i}\right)\\
 & \subset T_{i}\left[\overline{B_{\left(2K+1\right)^{2n}R}}\left(0\right)-\overline{B_{\left(2K+1\right)^{2n}R}}\left(0\right)\right]\\
 & \subset T_{i}\left[\overline{B_{\left(2K+1\right)^{2n}2R}}\left(0\right)\right].
\end{align*}
But this implies
\[
\lambda\left(\overline{Q_{i}^{n\ast}}-\overline{Q_{j}^{n\ast}}\right)\leq\left|\det T_{i}\right|\cdot\lambda\left(\overline{B_{\left(2K+1\right)^{2n}2R}}\left(0\right)\right),
\]
so that the choice $C_{3}=v_{\dimension}\cdot\left[\left(2K+1\right)^{2n}2R\right]^{\dimension}$
is possible.
\end{proof}
As the final result in this subsection, we use the inclusion from
Lemma~\ref{lem:SemiStructuredNormalizationNeighboring} to show that
the class of semi-structured coverings is closed under forming clusters,
i.e.\@ that $\CalQ^{k\ast}:=\left(Q_{i}^{k\ast}\right)_{i\in I}$
is again a semi-structured covering of $\CalO$ if $\CalQ$ is. This
property will frequently be useful, for instance in the proof of Lemma~\ref{lem:IntersectionCountForModerateCoverings}
below. Note that it is \emph{not} true in general that $\CalQ^{k\ast}$
is (almost) structured if $\CalQ$ is. This is one of the main reasons
for considering semi-structured coverings and not just almost structured
coverings.
\begin{lem}
\label{lem:SemiStructuredClusterInvariant}Assume that $\CalQ=\left(T_{i}Q_{i}'+b_{i}\right)_{i\in I}$
is a {[}tight{]} semi-structured covering of the open set $\emptyset\neq\CalO\subset\R^{\dimension}$.
Then the following are true:

\begin{enumerate}
\item We have $\left|i^{k\ast}\right|\leq N_{\CalQ}^{k}$ for all $i\in I$
and $k\in\N_{0}$. This even holds if $\CalQ$ is any admissible family.
\item The family of $k$-clusters $\CalQ^{k\ast}:=\left(Q_{i}^{k\ast}\right)_{i\in I}$
is a {[}tight{]} semi-structured covering of $\CalO$ satisfying $N_{\CalQ^{k\ast}}\leq N_{\CalQ}^{2k+1}$.
\item There are suitable open sets $P_{i}'\subset\R^{\dimension}$ for $i\in I$
so that $\left(T_{i},P_{i}',b_{i}\right)_{i\in I}$ is a {[}tight{]}
parametrization of $\CalQ^{k\ast}=\left(T_{i}P_{i}'+b_{i}\right)_{i\in I}$,
with
\[
R_{\CalQ^{k\ast}}\leq\left(2C_{\CalQ}+1\right)^{k}R_{\CalQ}\qquad\text{ and }\qquad C_{\CalQ^{k\ast}}\leq C_{\CalQ}^{2k+1},
\]
and (in the tight case) with $\varepsilon_{\CalQ^{k\ast}}\geq\varepsilon_{\CalQ}$.
\item If $u=\left(u_{i}\right)_{i\in I}$ is a $\CalQ$-moderate weight,
then
\[
u_{\ell}\leq C_{u,\CalQ}^{k}\cdot u_{i}\qquad\forall\:k\in\N_{0},\,i\in I\text{ and }\ell\in i^{k\ast}.
\]
Furthermore, $u$ is $\CalQ^{k\ast}$-moderate with $C_{u,\CalQ^{k\ast}}\leq C_{u,\CalQ}^{2k+1}$.\qedhere
\end{enumerate}
\end{lem}

\begin{proof}
We first show by induction on $k\in\N_{0}$ the following: If $u=\left(u_{i}\right)_{i\in I}$
is $\CalQ$-moderate, then 
\begin{equation}
\left|\smash{i^{k\ast}}\right|\leq N_{\CalQ}^{k}\quad\forall\:i\in I\quad\text{and further}\quad\frac{u_{i}}{u_{j}}\leq C_{u,\CalQ}^{k}\,\,\text{ and }\,\,\left\Vert \smash{T_{i}^{-1}}T_{j}\right\Vert \leq C_{\CalQ}^{k}\quad\forall\:i\in I\text{ and }j\in i^{k\ast}\,.\label{eq:SemiStructuredClusterInvariantAuxiliaryResult}
\end{equation}
The case $k=0$ is trivial, since $i^{0\ast}=\left\{ i\right\} $.
For the induction step, first note
\[
\left|\smash{i^{\left(k+1\right)\ast}}\right|=\left|\vphantom{\bigcup}\smash{\bigcup_{j\in i^{k\ast}}j^{\ast}}\right|\leq\sum_{j\in i^{k\ast}}\left|j^{\ast}\right|\leq N_{\CalQ}\cdot\left|\smash{i^{k\ast}}\right|\leq N_{\CalQ}^{k+1},
\]
where we used the induction hypothesis in the last step. For the second
and last estimates, let $j\in i^{\left(k+1\right)\ast}$ be arbitrary.
Thus, there is some $\ell\in i^{k\ast}$ with $j\in\ell^{\ast}$.
By induction hypothesis, this yields
\[
\frac{u_{i}}{u_{j}}=\frac{u_{i}}{u_{\ell}}\cdot\frac{u_{\ell}}{u_{j}}\leq C_{u,\CalQ}^{k}\cdot C_{u,\CalQ}=C_{u,\CalQ}^{k+1}.
\]
Similarly,
\[
\left\Vert T_{i}^{-1}T_{j}\right\Vert =\left\Vert T_{i}^{-1}T_{\ell}T_{\ell}^{-1}T_{j}\right\Vert \leq\left\Vert T_{i}^{-1}T_{\ell}\right\Vert \cdot\left\Vert T_{\ell}^{-1}T_{j}\right\Vert \leq C_{\CalQ}^{k}\cdot C_{\CalQ}=C_{\CalQ}^{k+1}.
\]
This completes the proof of equation~(\ref{eq:SemiStructuredClusterInvariantAuxiliaryResult}).

\medskip{}

Now, let $i\in I$ and $j\in i^{\ast_{\CalQ^{k\ast}}}$, i.e.\@ with
$Q_{i}^{k\ast}\cap Q_{j}^{k\ast}\neq\emptyset$. This implies $Q_{m}\cap Q_{\ell}\neq\emptyset$
for suitable $m\in i^{k\ast},\ell\in j^{k\ast}$, and hence $j\in\ell^{k\ast}\subset m^{\left(k+1\right)\ast}\subset i^{\left(2k+1\right)\ast}$.
Therefore, we have shown $i^{\ast_{\CalQ^{k\ast}}}\subset i^{\left(2k+1\right)\ast_{\CalQ}}$,
which yields $\left|\smash{i^{\ast_{\CalQ^{k\ast}}}}\right|\leq\left|\smash{i^{\left(2k+1\right)\ast_{\CalQ}}}\right|\leq N_{\CalQ}^{2k+1}$,
so that $\CalQ^{k\ast}$ is admissible with $N_{\CalQ^{k\ast}}\leq N_{\CalQ}^{2k+1}$.

As shown above, we also have $u_{i}/u_{j}\leq C_{u,\CalQ}^{2k+1}$
for $j\in i^{\ast_{\CalQ^{k\ast}}}\subset i^{\left(2k+1\right)\ast_{\CalQ}}$.
Thus, $u$ is $\CalQ^{k\ast}$-moderate with $C_{u,\CalQ^{k\ast}}\leq C_{u,\CalQ}^{2k+1}$.

Now, set $P_{i}':=T_{i}^{-1}\left(Q_{i}^{k\ast}-b_{i}\right)$ for
each $i\in I$. Thus, $T_{i}P_{i}'+b_{i}=Q_{i}^{k\ast}$ and $P_{i}'\neq\emptyset$
is open. We now show that $\CalQ^{k\ast}=\left(T_{i}P_{i}'+b_{i}\right)_{i\in I}$
obeys all requirements from Definition~\ref{defn:DifferentTypesOfCoverings}
for a semi-structured covering.

\begin{enumerate}
\item As seen above, $\CalQ^{k\ast}$ is admissible. Now note that $Q_{i}\neq\emptyset$
implies $i\in i^{k\ast}$ and hence 
\[
Q_{i}\subset Q_{i}^{k\ast}=\bigcup_{j\in i^{k\ast}}Q_{j}\subset\CalO
\]
for all $i\in I$, so that $\CalQ^{k\ast}$ covers $\CalO$.
\item Lemma~\ref{lem:SemiStructuredNormalizationNeighboring} yields
\[
Q_{i}^{k\ast}=\bigcup_{j\in i^{k\ast}}Q_{j}\subset T_{i}\left(\overline{B_{\left(2C_{\CalQ}+1\right)^{k}R_{\CalQ}}}\left(0\right)\right)+b_{i}
\]
for all $i\in I$, and thus $R_{\CalQ^{k\ast}}\leq\left(2C_{\CalQ}+1\right)^{k}R_{\CalQ}<\infty$.
\item For $i\in I$ and $j\in i^{\ast_{\CalQ^{k\ast}}}\subset i^{\left(2k+1\right)\ast_{\CalQ}}$,
equation~(\ref{eq:SemiStructuredClusterInvariantAuxiliaryResult})
yields $\left\Vert \smash{T_{i}^{-1}}T_{j}\right\Vert \leq C_{\CalQ}^{2k+1}$.
But this implies $C_{\CalQ^{k\ast}}\leq C_{\CalQ}^{2k+1}<\infty$.
\end{enumerate}
Finally, if $\CalQ$ is a tight covering, there is a family $\left(c_{i}\right)_{i\in I}\in\left(\R^{\dimension}\right)^{I}$
with $B_{\varepsilon_{\CalQ}}\left(c_{i}\right)\subset Q_{i}'$ for
all $i\in I$. This yields
\[
B_{\varepsilon_{\CalQ}}\left(c_{i}\right)\subset Q_{i}'=T_{i}^{-1}\left(Q_{i}-b_{i}\right)\subset T_{i}^{-1}\left(\smash{Q_{i}^{k\ast}}-b_{i}\right)=P_{i}'
\]
for all $i\in I$, so that $\CalQ^{k\ast}$ is a tight semi-structured
covering with $\varepsilon_{\CalQ^{k\ast}}\geq\varepsilon_{\CalQ}$.
\end{proof}

\subsection{Relations between coverings}

The main goal of this paper is to develop criteria for the existence
of embeddings $\FourierDecompSp{\CalQ}{p_{1}}Y\hookrightarrow\FourierDecompSp{\CalP}{p_{2}}Z$.
For most of these criteria, in particular for necessary conditions,
we will need to impose certain restrictions on the (geometric) relation
between the two coverings $\CalQ,\CalP$. In this subsection, we introduce
a convenient language for describing such relations. Furthermore,
we derive a few consequences of these relations.

We begin by summarizing all possible relations that we consider. We
remark that the notions of (weak/almost) subordinateness of coverings
were already introduced by Feichtinger and Gröbner, cf.\@ \cite[Definition 3.3]{DecompositionSpaces1}.
\begin{defn}
\label{def:RelativeIndexClustersSubordinateCoveringsModerateCoverings}(cf.\@
\cite[Definition 3.3]{DecompositionSpaces1}) Assume that $\CalQ=\left(Q_{i}\right)_{i\in I}$
and $\CalP=\left(P_{j}\right)_{j\in J}$ are two families of subsets
of $\R^{\dimension}$. Then

\begin{enumerate}
\item For $i\in I$ we define the \textbf{$\CalP$-index-cluster} around
$i$ (or the \textbf{set of $\CalP$-neighbors} of $i$) as
\[
J_{i}:=\left\{ j\in J\with P_{j}\cap Q_{i}\neq\emptyset\right\} .
\]
The \textbf{$\CalQ$-index-cluster} around $j\in J$ is defined analogously
and denoted by $I_{j}$.
\item We say that $\CalQ$ is \textbf{weakly subordinate} to $\CalP$ if
the constant 
\[
N\left(\CalQ,\CalP\right):=\sup_{i\in I}\left|J_{i}\right|
\]
 is finite. For equivalent conditions, see \cite[Definition 3.3 and Proposition 3.5]{DecompositionSpaces1}.

We say that $\CalQ$ and $\CalP$ are \textbf{weakly equivalent} if
$\CalQ$ is weakly subordinate to $\CalP$ and $\CalP$ is weakly
subordinate to $\CalQ$.
\item We say that $\CalQ$ is \textbf{almost subordinate} to $\CalP$ if
there is a constant $k=k\left(\CalQ,\CalP\right)\in\N_{0}$ such that
each set $Q_{i}$ is contained in some $P_{j}^{k\ast}$, i.e.\@ if
\[
\forall\:i\in I\,\exists\:j_{i}\in J:\qquad Q_{i}\subset P_{j_{i}}^{k\ast}.
\]
If we can take $k\left(\CalQ,\CalP\right)=0$, we say that $\CalQ$
is \textbf{subordinate} to $\CalP$.

We say that $\CalQ$ and $\CalP$ are \textbf{equivalent} if $\CalQ$
is almost subordinate to $\CalP$ and $\CalP$ is almost subordinate
to $\CalQ$.
\item A weight $u=\left(u_{i}\right)_{i\in I}$ is called \textbf{relatively
$\CalP$-moderate} if there is a constant $C_{u,\CalQ,\CalP}>0$ with
\[
u_{i}\leq C_{u,\CalQ,\CalP}\cdot u_{\ell}
\]
for all $j\in J$ and all $i,\ell\in I$ with $Q_{i}\cap P_{j}\neq\emptyset\neq Q_{\ell}\cap P_{j}$.
\item Now, let $\CalQ=\left(T_{i}Q_{i}'+b_{i}\right)_{i\in I}$ be a semi-structured
covering of some open set $\emptyset\neq\CalO\subset\R^{\dimension}$.
We say that $\CalQ$ is \textbf{relatively $\CalP$-moderate} if there
is a constant $C_{{\rm mod}}\left(\CalQ,\CalP\right)>0$ satisfying
\[
\forall\:j\in J\,\forall\:i,\ell\in I_{j}:\qquad\left|\det\left(\smash{T_{i}^{-1}}T_{\ell}\right)\right|\leq C_{{\rm mod}}\left(\CalQ,\CalP\right).
\]
Since the determinant is multiplicative, an equivalent assumption
is that the weight $\left(\left|\det T_{i}\right|\right)_{i\in I}$
is relatively $\CalP$-moderate.\qedhere
\end{enumerate}
\end{defn}

\begin{rem*}
We remark that the notations $J_{i}$ and $I_{j}$ introduced above
are strictly speaking ambiguous, at least in the case $I=J$, i.e.\@
if the \emph{same} index set is used for $\CalQ$ and $\CalP$. Nevertheless,
the context will always reveal what is meant.

\medskip{}

The condition of relative $\CalP$-moderateness of a weight $u$ means
that if two sets $Q_{i},Q_{\ell}$ of the covering $\CalQ$ are ``close''
to each other \emph{measured with respect to $\CalP$} (i.e., they
intersect the same $P_{j}$), then $u_{i},u_{\ell}$ are of similar
size. Analogously, $\CalQ$ is relatively moderate with respect to
$\CalP$ if the determinants $\left|\det T_{i}\right|,\left|\det T_{\ell}\right|$
are of comparable size if $Q_{i},Q_{\ell}$ are close to each other
measured with respect to $\CalP$.

It seems that relative moderateness is a novel concept. Its significance
will become clear in Subsection \ref{subsec:RelativelyModerateCase},
where we employ it to show that our criteria yield a \emph{complete
characterization} of the existence of the embedding $\FourierDecompSp{\CalQ}{p_{1}}{\ell_{w}^{q_{1}}}\hookrightarrow\FourierDecompSp{\CalP}{p_{2}}{\ell_{v}^{q_{2}}}$
as long as $\CalQ$ and $w$ are relatively $\CalP$-moderate (or
as long as $\CalP$ and $v$ are relatively $\CalQ$-moderate). A
further (related) application is Lemma~\ref{lem:IntersectionCountForModerateCoverings}
below, where a convenient estimate of $\left|J_{i}\right|$ is developed,
subject to suitable moderateness assumptions.

\medskip{}

Historically, one of the earliest examples of the use of equivalent
coverings in the study of function spaces occurred in the definition
of Triebel-Lizorkin spaces. Nowadays, these spaces are most commonly
defined using smooth partitions of unity subordinate to a covering
of $\R^{d}$ consisting of annuli; see for instance \cite[Section 2.3.1]{TriebelTheoryOfFunctionSpaces}.
The original definition by Lizorkin \cite{LizorkinCovering1,LizorkinCovering2},
however, was based on cubes, where essentially each annulus is partitioned
into a constant number of such cubes; see \cite[Section 2.5.4]{TriebelTheoryOfFunctionSpaces}
for more details.
\end{rem*}
The notions introduced in the above definition are of course not independent.
In the next lemmata, we explore the connections between these concepts.
\begin{lem}
\label{lem:SubordinatenessImpliesWeakSubordination}Assume that $\CalQ=\left(Q_{i}\right)_{i\in I}$
and $\CalP=\left(P_{j}\right)_{j\in J}$ are two families of nonempty
subsets of $\R^{\dimension}$. Then the following hold:

\begin{enumerate}
\item If $\CalQ$ is almost subordinate to $\CalP$ and if $\CalP$ is admissible,
then $\CalQ$ is weakly subordinate to $\CalP$ with 
\[
N\left(\CalQ,\CalP\right)\leq N_{\CalP}^{k\left(\CalQ,\CalP\right)+1}.
\]
If $\CalQ,\CalP$ are coverings of $\CalO,\CalO'\subset\R^{\dimension}$,
respectively, then we also have $\CalO\subset\CalO'$.
\item \label{enu:SubordinatenessIndexIsArbitraryIfItIntersects}If $Q_{i}\subset P_{j}^{k\ast}$
holds for some $i\in I$, $j\in J$ and $k\in\N_{0}$, then $J_{i}\subset\ell^{\left(2k+2\right)\ast}$
holds for all $\ell\in J_{i}$. In particular, $Q_{i}\subset P_{\ell}^{\left(2k+2\right)\ast}$
for all $\ell\in J_{i}$.
\item If $\CalQ$ is almost subordinate to $\CalP$, then $Q_{i}\subset P_{j}^{\left(2k\left(\CalQ,\CalP\right)+2\right)\ast}$
holds for all $i\in I$ and all $j\in J_{i}$.\qedhere
\end{enumerate}
\end{lem}

\begin{rem*}
The inclusion $Q_{i}\subset P_{\ell}^{\left(2k+2\right)\ast}$ for
all $\ell\in J_{i}$ from part~(\ref{enu:SubordinatenessIndexIsArbitraryIfItIntersects})
of the Lemma will turn out to be very useful for the proofs of our
embedding results in Sections~\ref{sec:SufficientConditions} and
\ref{sec:NecessaryConditions}: By only assuming $Q_{i}\subset P_{j}^{k\ast}$
for \emph{some} $j\in J$, we can already conclude $Q_{i}\subset P_{\ell}^{m\ast}$
(for suitable $m>k$) for \emph{all} $\ell\in J$ which intersect
$Q_{i}$ nontrivially.
\end{rem*}
\begin{proof}
Ad (1): Choose $k:=k\left(\CalQ,\CalP\right)$ so that for each $i\in I$
there is some $j_{i}\in J$ with $Q_{i}\subset P_{j_{i}}^{k\ast}$.
Now let $j\in J_{i}$ for some $i\in I$. This implies $\emptyset\neq P_{j}\cap Q_{i}\subset P_{j}\cap P_{j_{i}}^{k\ast}$,
and hence $j\in j_{i}^{\left(k+1\right)\ast}$. Thus, $J_{i}\subset j_{i}^{\left(k+1\right)\ast}$,
so that Lemma~\ref{lem:SemiStructuredClusterInvariant} yields $\left|J_{i}\right|\leq\big|j_{i}^{\left(k+1\right)\ast}\big|\leq N_{\CalP}^{k+1}$.
Since this holds for all $i\in I$, we get 
\[
N\left(\CalQ,\CalP\right)=\sup_{i\in I}\left|J_{i}\right|\leq N_{\CalP}^{k+1}<\infty\:.
\]

Finally, assume that $\CalQ,\CalP$ are coverings of $\CalO,\CalO'$,
respectively. Then
\[
\CalO=\bigcup_{i\in I}Q_{i}\subset\bigcup_{i\in I}P_{j_{i}}^{k\ast}\subset\bigcup_{j\in J}P_{j}\subset\CalO'.
\]

\medskip{}

Ad (2): Fix $\ell\in J_{i}$ and let $m\in J_{i}$ be arbitrary. This
implies 
\[
\emptyset\neq Q_{i}\cap P_{s}\subset P_{j}^{k\ast}\cap P_{s}\quad\text{for}\quad s\in\left\{ m,\ell\right\} \,,
\]
and thus $j\in\ell^{\left(k+1\right)\ast}$, as well as $m\in j^{\left(k+1\right)\ast}\subset\ell^{\left(2k+2\right)\ast}$.
We have thus shown $J_{i}\subset\ell^{\left(2k+2\right)\ast}$. Finally,
note that the inclusions $Q_{i}\subset P_{j}^{k\ast}\subset\bigcup_{m\in J}P_{m}$
and $J_{i}\subset\ell^{\left(2k+2\right)\ast}$ imply
\[
Q_{i}\subset\bigcup_{m\in J}\left(P_{m}\cap Q_{i}\right)\subset\bigcup_{m\in J_{i}}P_{m}\subset P_{\ell}^{\left(2k+2\right)\ast}\qquad\forall\:\ell\in J_{i}\,.
\]

\medskip{}

Ad (3): This is a special case of part~(\ref{enu:SubordinatenessIndexIsArbitraryIfItIntersects}).
\end{proof}
It was observed by Feichtinger and Gröbner in \cite[Proposition 3.6]{DecompositionSpaces1}
that weak subordinateness implies almost subordinateness if we impose
certain \emph{connectivity assumptions}. This will turn out to be
very convenient for verifying almost subordinateness for concrete
examples.

The statement of the following lemma is very close to that of \cite[Proposition 3.6]{DecompositionSpaces1}.
The only difference is that a few unnecessary assumptions (like connectedness
of the space $X$) have been removed, and the statement of the lemma
has been made more quantitative. The proof is still the same as that
of \cite[Proposition 3.6]{DecompositionSpaces1} and is hence omitted.
\begin{lem}
\label{lem:WeaksubordinatenessImpliesSubordinatenessIfConnected}(cf.\@
\cite[Proposition 3.6]{DecompositionSpaces1}) Let $\CalQ=\left(Q_{i}\right)_{i\in I}$
and $\CalP=\left(P_{j}\right)_{j\in J}$ be families of subsets of
a topological space $X$. Let $i\in I$ and assume that $Q_{i}$ is
path-connected with $Q_{i}\subset\bigcup_{j\in J}P_{j}$ and so that
\[
J_{i}:=\left\{ j\in J\with P_{j}\cap Q_{i}\neq\emptyset\right\} 
\]
is finite. Furthermore, assume that $P_{j}$ is open for every $j\in J_{i}$.\vspace{0.1cm}

Let $r:=\left|J_{i}\right|$. We then have $J_{i}\subset j^{r\ast}$
and in particular $Q_{i}\subset P_{j}^{r\ast}$ for every $j\in J_{i}$.
\end{lem}

The above lemma immediately yields the following corollary.
\begin{cor}
\label{cor:WeakSubordinationImpliesSubordinationIfConnected}Let $\CalQ=\left(Q_{i}\right)_{i\in I}$
and $\CalP=\left(P_{j}\right)_{j\in J}$ be families of subsets of
$\R^{\dimension}$ such that each $Q_{i}$ is path-connected and such
that each $P_{j}$ is open.

Assume that $\CalQ$ is weakly subordinate to $\CalP$ with $\bigcup_{i\in I}Q_{i}\subset\bigcup_{j\in J}P_{j}$.
Then $\CalQ$ is almost subordinate to $\CalP$ with $k\left(\CalQ,\CalP\right)\leq N\left(\CalQ,\CalP\right)$.
\end{cor}

One of the fundamental tools that we will use again and again in the
remainder of the paper is the following \textbf{disjointization lemma}
for admissible coverings that was developed by Feichtinger and Gröbner
in \cite[Lemma 2.9]{DecompositionSpaces1}.
\begin{lem}
\label{lem:DisjointizationPrinciple}(cf.\@ \cite[Lemma 2.9]{DecompositionSpaces1})
Let $\CalO\neq\emptyset$, and let $\CalQ=\left(Q_{i}\right)_{i\in I}$
be an admissible covering of $\CalO$.

Then, for any $m\in\N_{0}$, there exists a finite partition $I=\biguplus_{r=1}^{r_{0}}I^{\left(r\right)}$
with $Q_{i}^{m\ast}\cap Q_{j}^{m\ast}=\emptyset$ for all $i,j\in I^{\left(r\right)}$
with $i\neq j$ and all $r\in\left\{ 1,\dots,r_{0}\right\} $. In
fact, one can choose $r_{0}=N_{\CalQ}^{2m+1}$.
\end{lem}

Note that the exact choice of $r_{0}$ from above is not explicitly
stated in \cite[Lemma 2.9]{DecompositionSpaces1}. For the sake of
completeness, we thus present a proof which is based on the following
slightly more general lemma. The proof, however, is still similar
to that of \cite[Lemma 2.9]{DecompositionSpaces1}.
\begin{lem}
\label{lem:RelationDisjointification}Let $X\neq\emptyset$ be a set
and let $\sim\subset X\times X$ be a relation which is reflexive
and symmetric, but not necessarily transitive. For $x\in X$, let
$\left[x\right]:=\left\{ y\in X\with y\sim x\right\} $, and assume
that the maximal cardinality $N:=\sup_{x\in X}\left|\left[x\right]\right|$
is finite.

Then there is a partition $X=\biguplus_{\ell=1}^{N}X_{\ell}$ with
$x\not\sim y$ for all $x,y\in X_{\ell}$ with $x\neq y$ and arbitrary
$\ell\in\underline{N}$.
\end{lem}

\begin{proof}
Let $X_{1}\subset X$ be maximal with the following property: For
$x,y\in X_{1}$ with $x\neq y$, we have $x\not\sim y$. Existence
of such a set is an easy consequence of Zorn's Lemma.

Now, if $X_{1},\dots,X_{m}$ are already constructed for some $m\in\underline{N-1}$,
let $X_{m+1}\subset X\setminus\left(X_{1}\cup\dots\cup X_{m}\right)$
be maximal with the same property as above. Again, existence follows
from Zorn's Lemma.

Clearly, $\left(X_{\ell}\right)_{\ell\in\underline{N}}$ are pairwise
disjoint and satisfy the property that $x\not\sim y$ for $x,y\in X_{\ell}$
with $x\neq y$ and arbitrary $\ell\in\underline{N}$.

It remains to show $X=\bigcup_{\ell=1}^{N}X_{\ell}$. Suppose that
this fails. Then there is some $x\in X\setminus\bigcup_{\ell=1}^{N}X_{\ell}$.
Thus, $X_{\ell}\cup\left\{ x\right\} \subset X\setminus\left(X_{1}\cup\dots\cup X_{\ell-1}\right)$
is a strict superset of $X_{\ell}$ for arbitrary $\ell\in\underline{N}$.
By maximality of $X_{\ell}$, we see that there must be some $y_{\ell}\in X_{\ell}\cup\left\{ x\right\} $
with $x\neq y_{\ell}$ and $x\sim y_{\ell}$. Note $y_{\ell}\in X_{\ell}$
since $x\neq y_{\ell}$. By disjointness of the $\left(X_{\ell}\right)_{\ell\in\underline{N}}$,
we see that $x,y_{1},\dots,y_{N}$ are pairwise distinct. Hence,
\[
N\geq\left|\left[x\right]\right|\geq\left|\left\{ x,y_{1},\dots,y_{N}\right\} \right|=N+1,
\]
a contradiction.
\end{proof}

\begin{proof}[Proof of Lemma~\ref{lem:DisjointizationPrinciple}]
For $i,j\in I$, write $i\sim j$ if and only if $Q_{i}^{m\ast}\cap Q_{j}^{m\ast}\neq\emptyset$.
Note that $Q_{i}^{m\ast}\supset Q_{i}\neq\emptyset$ for all $i\in I$,
by definition of an admissible covering. This shows that $\sim$ is
reflexive. Symmetry of $\sim$ is clear.

Finally, note (in the notation of Lemma~\ref{lem:RelationDisjointification})
that $\left[i\right]=i^{\ast_{\CalQ^{m\ast}}}$ satisfies
\[
\left|\left[i\right]\right|=\left|i^{\ast_{\CalQ^{m\ast}}}\right|\leq N_{\CalQ^{m\ast}}\leq N_{\CalQ}^{2m+1}\,,
\]
thanks to Lemma~\ref{lem:SemiStructuredClusterInvariant}. Application
of Lemma~\ref{lem:RelationDisjointification} completes the proof.
\end{proof}
The next lemma shows that restricting attention solely to coverings
for which one is relatively moderate with respect to the other would
prevent one from handling the case where the two covered sets are
distinct.
\begin{lem}
\label{lem:DifferentOrbitsPreventModerateness}Let $\emptyset\neq\CalO,\CalO'\subset\R^{\dimension}$
be open and assume that $\CalQ=\left(Q_{i}\right)_{i\in I}$ is an
admissible covering of $\CalO$ and that $\mathcal{P}=\left(P_{j}\right)_{j\in J}$
is an admissible covering of $\CalO'$.

Finally, assume $\CalO'\cap\partial\CalO\neq\emptyset$ and that $\CalQ,\CalP$
admit partitions of unity; see Lemma~\ref{lem:PartitionCoveringNecessary}.
Then the following hold:

\begin{enumerate}
\item $\CalP$ is \emph{not} weakly subordinate to $\CalQ$. Even more,
we have $\left|I_{j}\right|=\infty$ for some $j\in J$.
\item If $\CalQ=\left(T_{i}Q_{i}'+b_{i}\right)_{i\in I}$ is a \emph{tight}
semi-structured admissible covering of $\CalO$ and if $\CalQ$ is
almost subordinate to $\CalP$, then $\CalQ$ is \emph{not} relatively
$\CalP$-moderate.\qedhere
\end{enumerate}
\end{lem}

\begin{rem*}
It is instructive to note that $\CalO'\cap\partial\CalO\neq\emptyset$
is \emph{always} satisfied if $\CalO\subset\R^{\dimension}$ is dense
(for example, if $\CalO$ is of full measure) with $\CalO\subsetneq\CalO'$.
In particular, this holds if $\CalP$ is the \emph{in}homogeneous
Besov covering of $\CalO'=\R^{\dimension}$ and if $\CalQ$ denotes
the homogeneous Besov covering of $\CalO=\R^{\dimension}\setminus\left\{ 0\right\} $.
These coverings will be formally introduced in Definitions \ref{def:InhomogeneousBesovCovering}
and \ref{def:HomogeneousBesovCovering}, respectively.
\end{rem*}
\begin{proof}
Choose $x\in\CalO'\cap\partial\CalO$. By Lemma~\ref{lem:PartitionCoveringNecessary},
there is some $j\in J$ with $x\in P_{j}^{\circ}$. Because of $x\in\partial\CalO$,
there is a sequence $\left(x_{n}\right)_{n\in\N}\in\CalO^{\N}$ with
$x_{n}\xrightarrow[n\rightarrow\infty]{}x$. Note that we have $x_{n}\in P_{j}^{\circ}$
for $n$ sufficiently large, so that we can assume $x_{n}\in P_{j}^{\circ}$
for all $n\in\N$.

For each $n\in\N$, choose $i_{n}\in I$ with $x_{n}\in Q_{i_{n}}$.
Assume towards a contradiction that there is some $i\in I$ with $i_{n}=i$
for infinitely many $n\in\N$. Since $\overline{Q_{i}}\subset\CalO$
holds by Lemma~\ref{lem:PartitionCoveringNecessary}, this implies
$x\in\overline{Q_{i}}\subset\CalO=\CalO^{\circ}$ in contradiction
to $x\in\partial\CalO$. Thus, restricting to a subsequence, we can
assume that the $\left(i_{n}\right)_{n\in\N}$ are pairwise distinct.
Using $x_{n}\in Q_{i_{n}}\cap P_{j}^{\circ}$, we conclude that $I_{j}\supset\left\{ i_{n}\with n\in\N\right\} $
is infinite, so that $\CalP$ is not weakly subordinate to $\CalQ$.
This establishes the first claim.

\medskip{}

Now assume towards a contradiction that $\CalQ=\left(T_{i}Q_{i}'+b_{i}\right)_{i\in I}$
is a tight semi-structured covering and also almost subordinate to
$\CalP$ and relatively $\CalP$-moderate. Using part~(\ref{enu:SubordinatenessIndexIsArbitraryIfItIntersects})
of Lemma~\ref{lem:SubordinatenessImpliesWeakSubordination} and the
definition of relative $\CalP$-moderateness, we see that there is
some $\ell\in\N$ with $Q_{i_{n}}\subset P_{j}^{\ell\ast}$ for all
$n\in\N$ and some $C>0$ with
\[
\left|\det\left(\smash{T_{i_{n}}^{-1}}T_{i_{m}}\right)\right|\leq C\qquad\forall\:n,m\in\N\,,
\]
because of $Q_{i_{n}}\cap P_{j}\neq\emptyset$. Since $\CalQ$ is
tight, Corollary~\ref{cor:SemiStructuredDifferenceSetsMeasureEstimate}
yields a constant $c=c\left(\CalQ,\varepsilon_{\CalQ},d\right)>0$
satisfying
\[
\lambda\left(Q_{i_{n}}\right)\geq c\cdot\left|\det T_{i_{n}}\right|\geq\frac{c}{C}\cdot\left|\det T_{i_{1}}\right|=:c'\qquad\forall\:n\in\N\,.
\]

But admissible coverings are of finite height; explicitly, we have
$\sum_{n\in\N}\Indicator_{Q_{i_{n}}}\leq\sum_{i\in I}\Indicator_{Q_{i}}\leq N_{\CalQ}$.
Furthermore, we saw above that $Q_{i_{n}}\subset P_{j}^{\ell\ast}$
for all $n\in\N$. All in all, we conclude
\[
\sum_{n\in\N}\Indicator_{Q_{i_{n}}}\leq N_{\CalQ}\cdot\Indicator_{P_{j}^{\ell\ast}}\:,\qquad\text{and thus}\qquad\infty>N_{\CalQ}\cdot\lambda\big(P_{j}^{\ell\ast}\big)\geq\sum_{n\in\N}\lambda\left(Q_{i_{n}}\right)\geq\sum_{n\in\N}c'=\infty\,,
\]
a contradiction. Here, the finite union $P_{j}^{\ell\ast}=\bigcup_{m\in j^{\ell\ast}}\:P_{m}$
is of finite measure, by compactness of each $\overline{P_{m}}$,
see Lemma~\ref{lem:PartitionCoveringNecessary}.
\end{proof}
Using similar techniques as in the proof above, we will now establish
an easy way to estimate the cardinality $\left|I_{j}\right|$ of the
set of $\CalQ$-neighbors of $P_{j}$. This method of estimating $\left|I_{j}\right|$
is implicitly used in \cite{HanWangAlphaModulationEmbeddings} for
the concrete setting of $\alpha$-modulation spaces, but not stated
explicitly. In the present generality, it seems to be a new observation
which will be of great use to us in Subsection \ref{subsec:RelativelyModerateCase}.

It is important to note that we assume $\CalQ$ to be relatively $\CalP$-moderate
and (for the lower estimate) almost subordinate to $\CalP$. Without
these assumptions, it is easy to see that an estimate as in equation~(\ref{eq:IntersectionCountGivenByDeterminantsForModerateCoverings})
below fails in general.
\begin{lem}
\label{lem:IntersectionCountForModerateCoverings}Let $\emptyset\neq\CalO,\CalO'\subset\R^{\dimension}$
be open and assume that the families $\CalQ=\left(Q_{i}\right)_{i\in I}=\left(T_{i}Q_{i}'+b_{i}\right)_{i\in I}$
and $\CalP=\left(P_{j}\right)_{j\in J}=\left(\smash{S_{j}P_{j}'+c_{j}}\right)_{j\in J}$
are \emph{tight} semi-structured coverings of $\CalO$ and $\CalO'$,
respectively.

Finally, let $I_{0}\subset I$, $J_{0}\subset J$ and assume that

\begin{enumerate}
\item \label{enu:IntersectionCountSubordinateness}$\CalQ_{I_{0}}:=\left(Q_{i}\right)_{i\in I_{0}}$
is almost subordinate to $\CalP$.
\item $\CalQ_{I_{0}}$ is relatively $\CalP_{J_{0}}$-moderate, with $\CalP_{J_{0}}:=\left(P_{j}\right)_{j\in J_{0}}$.
\item \label{enu:IntersectionCountCoarseCoveringCoveredByFine}There is
some $r\in\N_{0}$ and some $C_{0}>0$ such that
\begin{equation}
\lambda\left(P_{j}\right)\leq C_{0}\cdot\lambda\bigg(\bigcup_{i\in I_{0}\cap I_{j}}Q_{i}^{r\ast}\bigg)\qquad\forall\:j\in J_{0}\,\text{ with }\,1\leq\left|I_{0}\cap I_{j}\right|<\infty\,.\label{eq:IntersectionCountModerateTechnicalCondition}
\end{equation}
\end{enumerate}
Then there are positive constants 
\[
C_{1}=C_{1}\left(\dimension,\CalQ,\varepsilon_{\CalQ},\CalP,k\left(\smash{\CalQ_{I_{0}}},\CalP\right),C_{{\rm mod}}\left(\smash{\CalQ_{I_{0}}},\CalP_{J_{0}}\right)\right)
\]
and
\[
C_{2}=C_{2}\left(\dimension,C_{0},r,\CalQ,\CalP,\varepsilon_{\CalP},C_{{\rm mod}}\left(\smash{\CalQ_{I_{0}}},\CalP_{J_{0}}\right)\right)
\]
with
\begin{equation}
C_{1}^{-1}\cdot\left|I_{0}\cap I_{j}\right|\leq\frac{\left|\det S_{j}\right|}{\left|\det T_{i}\right|}\leq C_{2}\cdot\left|I_{0}\cap I_{j}\right|\qquad\forall\:i\in I_{0}\cap I_{j}\qquad\forall\,j\in J_{0}.\label{eq:IntersectionCountGivenByDeterminantsForModerateCoverings}
\end{equation}
Finally,
\begin{enumerate}
\item Tightness of $\CalP$ and assumption (\ref{enu:IntersectionCountCoarseCoveringCoveredByFine})
are only needed for the right estimate in equation~(\ref{eq:IntersectionCountGivenByDeterminantsForModerateCoverings}).
\item Tightness of $\CalQ$ and assumption (\ref{enu:IntersectionCountSubordinateness})
are only needed for the left estimate in equation~(\ref{eq:IntersectionCountGivenByDeterminantsForModerateCoverings}).\qedhere
\end{enumerate}
\end{lem}

\begin{rem*}
(1) The moderateness assumption means that there is a constant $C=C_{{\rm mod}}\left(\smash{\CalQ_{I_{0}}},\CalP_{J_{0}}\right)$
with $\left|\det T_{i}\right|\leq C\cdot\left|\det T_{\ell}\right|$
for all $i,\ell\in I_{0}$ for which there is some $j\in J_{0}$ with
$Q_{i}\cap P_{j}\neq\emptyset\neq Q_{\ell}\cap P_{j}$.

\medskip{}

(2) It is worth noting that estimate~(\ref{eq:IntersectionCountModerateTechnicalCondition})
from above is trivially satisfied (with $C_{0}=1$, $r=0$ and arbitrary
$J_{0}\subset J$) if we have $\CalO'\subset\CalO$ (up to a set of
measure zero) and $I_{0}=I$, because in that case we have (up to
a set of measure zero) 
\[
P_{j}\subset\CalO'\subset\CalO=\bigcup_{i\in I}Q_{i}=\bigcup_{i\in I_{0}}Q_{i}^{0\ast}\:,
\]
which easily implies $P_{j}\subset\bigcup_{i\in I_{0}\cap I_{j}}Q_{i}^{0\ast}$
(up to a set of measure zero) and thus $\lambda\left(P_{j}\right)\leq\lambda\big(\bigcup_{i\in I_{0}\cap I_{j}}Q_{i}^{0\ast}\big)$.

\medskip{}

(3) Finally, it is important to observe that the statement in equation~(\ref{eq:IntersectionCountGivenByDeterminantsForModerateCoverings})
is void for $j\in J_{0}$ with $I_{0}\cap I_{j}=\emptyset$.
\end{rem*}
\begin{proof}
We first prove the right estimate, assuming that $\CalP$ is tight,
and assuming equation~(\ref{eq:IntersectionCountModerateTechnicalCondition}).
Since $\CalP$ is a tight semi-structured covering, Corollary~\ref{cor:SemiStructuredDifferenceSetsMeasureEstimate}
yields a constant $c_{1}=c_{1}\left(\dimension,\varepsilon_{\CalP}\right)>0$
such that $\lambda\left(P_{j}\right)\geq c_{1}^{-1}\cdot\left|\det S_{j}\right|$
for all $j\in J$. Likewise, Lemma~\ref{lem:SemiStructuredClusterInvariant}
shows that $\CalQ^{r\ast}$ is also a semi-structured covering with
$R_{\CalQ^{r\ast}}\leq\left(2C_{\CalQ}+1\right)^{r}\cdot R_{\CalQ}$;
therefore, Corollary~\ref{cor:SemiStructuredDifferenceSetsMeasureEstimate}
shows $\lambda\left(Q_{i}^{r\ast}\right)\leq C_{1}\cdot\left|\det T_{i}\right|$
for all $i\in I$ and a suitable constant $C_{1}=C_{1}\left(\dimension,r,\CalQ\right)$.
Finally, set $K:=C_{{\rm mod}}\left(\smash{\CalQ_{I_{0}}},\CalP_{J_{0}}\right)$,
so that $\left|\det\left(\smash{T_{i}^{-1}}T_{\ell}\right)\right|\leq K$
holds for all $j\in J_{0}$ and $i,\ell\in I_{j}\cap I_{0}$.

Choose an arbitrary $j\in J_{0}$ with $I_{0}\cap I_{j}\neq\emptyset$,
noting that equation~(\ref{eq:IntersectionCountGivenByDeterminantsForModerateCoverings})
is trivially satisfied in case of $I_{0}\cap I_{j}=\emptyset$. Fix
any $i_{0}\in I_{0}\cap I_{j}$. This yields $\left|\det T_{i}\right|\leq K\cdot\left|\det T_{i_{0}}\right|$
for all $i\in I_{j}\cap I_{0}$. If $\left|I_{0}\cap I_{j}\right|$
is infinite, the upper estimate in equation~(\ref{eq:IntersectionCountGivenByDeterminantsForModerateCoverings})
is trivial. Hence, we can assume that $I_{0}\cap I_{j}$ is finite,
i.e.\@ $1\leq\left|I_{0}\cap I_{j}\right|<\infty$. Using equation~(\ref{eq:IntersectionCountModerateTechnicalCondition}),
we see
\[
c_{1}^{-1}\cdot\left|\det S_{j}\right|\leq\lambda\left(P_{j}\right)\leq C_{0}\cdot\lambda\biggl(\,\bigcup_{i\in I_{0}\cap I_{j}}Q_{i}^{r\ast}\biggr)\leq C_{0}C_{1}\cdot\!\sum_{i\in I_{0}\cap I_{j}}\left|\det T_{i}\right|\leq C_{0}C_{1}\cdot K\cdot\left|\det T_{i_{0}}\right|\cdot\left|I_{0}\cap I_{j}\right|.
\]
Since $i_{0}\in I_{0}\cap I_{j}$ was arbitrary, the upper estimate
in equation~(\ref{eq:IntersectionCountGivenByDeterminantsForModerateCoverings})
is established.

\medskip{}

We now prove the left estimate, assuming that $\CalQ$ is tight and
that $\CalQ_{I_{0}}$ is almost subordinate to $\CalP$. As above,
let $K:=C_{{\rm mod}}\left(\smash{\CalQ_{I_{0}}},\CalP_{J_{0}}\right)$,
so that $\left|\det\left(\smash{T_{i}^{-1}}T_{\ell}\right)\right|\leq K$
for all $j\in J_{0}$ and $i,\ell\in I_{j}\cap I_{0}$. By tightness
of $\CalQ$, Corollary~\ref{cor:SemiStructuredDifferenceSetsMeasureEstimate}
yields a constant $c_{2}=c_{2}\left(\dimension,\varepsilon_{\CalQ}\right)>0$
with $\lambda\left(Q_{i}\right)\geq c_{2}^{-1}\cdot\left|\det T_{i}\right|$
for all $i\in I$. Lemma~\ref{lem:SubordinatenessImpliesWeakSubordination}
shows that $Q_{i}\subset P_{j}^{m\ast}$ for all $i\in I_{0}$ and
$j\in J_{i}$, where we defined $m:=2\cdot k\left(\smash{\CalQ_{I_{0}}},\CalP\right)+2$.
Next, Lemma~\ref{lem:SemiStructuredClusterInvariant} shows that
$\CalP^{m\ast}$ is a semi-structured covering with $R_{\CalP^{m\ast}}\leq\left(1+2C_{\CalP}\right)^{m}\cdot R_{\CalP}$,
so that Corollary~\ref{cor:SemiStructuredDifferenceSetsMeasureEstimate}
yields a constant $C_{2}=C_{2}\left(\CalP,\dimension,m\right)>0$
satisfying $\lambda\left(P_{j}^{m\ast}\right)\leq C_{2}\cdot\left|\det S_{j}\right|$
for all $j\in J$.

Finally, let $j\in J_{0}$ be arbitrary with $I_{0}\cap I_{j}\neq\emptyset$,
noting that equation~(\ref{eq:IntersectionCountGivenByDeterminantsForModerateCoverings})
is trivially satisfied in case of $I_{0}\cap I_{j}=\emptyset$. Let
$i_{0}\in I_{0}\cap I_{j}$ be arbitrary, and let $\Gamma\subset I_{0}\cap I_{j}$
be an arbitrary finite subset of $I_{0}\cap I_{j}$. For any $i\in\Gamma\subset I_{0}\cap I_{j}$,
we have $i\in I_{0}$ and $j\in J_{i}$, so that $Q_{i}\subset P_{j}^{m\ast}$,
as seen above. Furthermore, the admissibility of $\CalQ$ entails
$\sum_{i\in I}\Indicator_{Q_{i}}\leq N_{\CalQ}$. By combining the
last two observations, we get $\sum_{i\in\Gamma}\Indicator_{Q_{i}}\left(\xi\right)\leq N_{\CalQ}\cdot\Indicator_{P_{j}^{m\ast}}\left(\xi\right)$
for all $\xi\in\R^{\dimension}$. Therefore,
\begin{align*}
\frac{1}{c_{2}K}\cdot\left|\det T_{i_{0}}\right|\cdot\left|\Gamma\right|\leq\sum_{i\in\Gamma}c_{2}^{-1}\cdot\left|\det T_{i}\right|\leq\sum_{i\in\Gamma}\lambda\left(Q_{i}\right)=\int_{\R^{\dimension}}\sum_{i\in\Gamma}\Indicator_{Q_{i}}\left(\xi\right)\,\d\xi & \leq N_{\CalQ}\cdot\int_{\R^{\dimension}}\Indicator_{P_{j}^{m\ast}}\left(\xi\right)\,\d\xi\\
 & \leq N_{\CalQ}\cdot C_{2}\cdot\left|\det S_{j}\right|.
\end{align*}
Since $\Gamma\subset I_{0}\cap I_{j}$ was an arbitrary finite subset,
we see that $I_{0}\cap I_{j}$ is finite and that the lower bound
in equation~(\ref{eq:IntersectionCountGivenByDeterminantsForModerateCoverings})
is satisfied.
\end{proof}

\section{(Fourier-side) decomposition spaces}

\label{sec:DecompositionSpaces}Now that we understand the relevant
types of coverings, we are in a position to properly start our analysis
of decomposition spaces. In this section, we will define these spaces
and derive their most important basic properties, i.e.\@ well-definedness
and completeness.

We note that these are nontrivial issues: Completeness only holds
in general when the right ``reservoir'' of functions is used. We
will see that, in general, the space $\Schwartz'\left(\R^{\dimension}\right)$
of tempered distributions is \emph{not} a suitable choice, even for
$\CalO=\R^{\dimension}$. Furthermore, for $p\in\left(0,1\right)$,
independence of $\FourierDecompSp{\CalQ}pY$ from the chosen partition
of unity is not as straightforward as for $p\in\left[1,\infty\right]$,
since Young's convolution inequality $L^{1}\ast L^{p}\hookrightarrow L^{p}$
fails in the quasi-Banach regime $p\in\left(0,1\right)$. Instead,
we have
\begin{equation}
\left\Vert \Fourier^{-1}\left(fg\right)\right\Vert _{L^{p}}\lesssim\left\Vert \Fourier^{-1}f\right\Vert _{L^{p}}\cdot\left\Vert \Fourier^{-1}g\right\Vert _{L^{p}}\label{eq:ConvolutionQuasiBanachAbstract}
\end{equation}
if both $f,g$ are compactly supported. Somewhat unexpectedly, the
implied constant in (\ref{eq:ConvolutionQuasiBanachAbstract}) depends
on the measure of the \emph{algebraic difference }(or \emph{Minkowski
difference})
\[
K-L=\left\{ k-\ell\with k\in K,\ell\in L\right\} ,
\]
where $\supp f\subset K$ and $\supp g\subset L$. Both of these issues
are somewhat neglected in the standard reference\cite{BorupNielsenDecomposition}
for Fourier-analytic decomposition spaces.

The structure of this section is as follows: In the first subsection,
we briefly present a treatment of a variant of Young's inequality
for convolution in $L^{p}\left(\R^{\dimension}\right)$ in the quasi-Banach
regime $p\in\left(0,1\right)$. This is based on Triebel's book \cite{TriebelTheoryOfFunctionSpaces}.
We repeat some of the arguments given there, since we need to know
the explicit constant arising in equation~(\ref{eq:ConvolutionQuasiBanachAbstract}),
depending on the supports of $f,g$. This constant is not stated explicitly
by Triebel.

The remaining subsections \ref{subsec:DecompositionDefinition}–\ref{subsec:DecompositionCompleteness}
are devoted, respectively, to the definition, well-definedness and
completeness of the (Fourier side) decomposition space $\FourierDecompSp{\CalQ}pY$
and its ``space-side'' counterpart. While defining these spaces,
we also introduce the notion of \textbf{$L^{p}$-BAPUs}, i.e.\@ of
the class of partitions of unity $\left(\varphi_{i}\right)_{i\in I}$
which are suitable to define these spaces. In particular, we note
that every almost-structured covering admits an $L^{p}$-BAPU.

\subsection{Convolution relations for \texorpdfstring{$L^{p}$, $p\in\left(0,1\right)$}{Lᵖ, p∈(0,1)}}

\label{subsec:QuasiBanachConvolution}As noted above, in this subsection,
we will establish the inequality
\begin{equation}
\left\Vert \smash{\Fourier^{-1}\left(f\cdot g\right)}\right\Vert _{L^{p}}\leq C\cdot\left\Vert \smash{\Fourier^{-1}f}\right\Vert _{L^{p}}\cdot\left\Vert \smash{\Fourier^{-1}g}\right\Vert _{L^{p}}\label{eq:QuasiBanachConvolutionNoDetails}
\end{equation}
for compactly supported functions/distributions $f,g$ and $p\in\left(0,1\right)$.
Our treatment is largely identical to that of Triebel in \cite{TriebelTheoryOfFunctionSpaces},
but we still include the proofs, since the dependency of the constant
$C$ on the supports of $f,g$ will be relevant to us. This dependence
is not explicitly stated by Triebel.

We first note that Young's inequality $\left\Vert f\ast g\right\Vert _{L^{p}}\leq\left\Vert f\right\Vert _{L^{1}}\cdot\left\Vert g\right\Vert _{L^{p}}$
fails completely for $p\in\left(0,1\right)$, even if both $f,g$
have compact Fourier support. This is shown in the following example:
\begin{example}
\label{ex:YoungInequalityFailsForQuasiBanach}Let us define 
\[
g_{1}:\R\to\R,\xi\mapsto\max\left\{ 0,1-\left|\xi\right|\right\} =\begin{cases}
0, & \text{if }\left|\xi\right|\geq1,\\
1+\xi, & \text{if }-1\leq\xi\leq0,\\
1-\xi, & \text{if }0\leq\xi\leq1,
\end{cases}
\]
and
\[
g_{2}:\R\to\R,\xi\mapsto\begin{cases}
0, & \text{if }\left|\xi\right|\geq1,\\
\left(\xi-1\right)^{2}\cdot\left(\xi+1\right)^{2}=\xi^{4}-2\xi^{2}+1, & \text{if }\left|\xi\right|\leq1.
\end{cases}
\]
Then $g_{1},g_{2}$ are supported in $\left[-1,1\right]$. Hence,
if we set $f_{j}:=\Fourier^{-1}g_{j}$ for $j\in\left\{ 1,2\right\} $,
then $f_{1},f_{2}$ are bandlimited.

A straightforward, but tedious calculation shows
\[
f_{1}\left(x\right)=-\frac{1}{2}\frac{\cos\left(2\pi x\right)-1}{\pi^{2}x^{2}}\qquad\forall\,\,x\in\R\setminus\left\{ 0\right\} ,
\]
as well as
\[
f_{2}\left(x\right)=\frac{3}{2}\cdot\frac{\sin\left(2\pi x\right)}{\pi^{5}x^{5}}-3\cdot\frac{\cos\left(2\pi x\right)}{\pi^{4}x^{4}}-2\cdot\frac{\sin\left(2\pi x\right)}{\pi^{3}x^{3}}\qquad\forall\,\,x\in\R\setminus\left\{ 0\right\} .
\]
Being (inverse) Fourier transforms of $L^{1}$-functions, $f_{1},f_{2}$
are both bounded. Furthermore, $f_{1}$ decays (at least) like $x^{-2}$
and $f_{2}$ decays (at least) like $\left|x\right|^{-3}$ at $\pm\infty$.
This implies $f_{1}\in L^{p}\left(\R\right)$ for all $p>\frac{1}{2}$,
as well as $f_{2}\in L^{p}\left(\R\right)$ for $p>\frac{1}{3}$.
In particular, $f_{1},f_{2}\in L^{1}\left(\R\right)$.

Another calculation using the convolution theorem, i.e.\@
\[
h:=f_{1}\ast f_{2}=\Fourier^{-1}\left(\,\smash{\widehat{f_{1}}}\cdot\smash{\widehat{f_{2}}}\,\right)=\Fourier^{-1}\left(g_{1}\cdot g_{2}\right)
\]
leads to
\[
h\left(x\right)=-\frac{15}{4}\frac{\cos\left(2\pi x\right)-1}{\pi^{6}x^{6}}-6\frac{\sin\left(2\pi x\right)}{\pi^{5}x^{5}}+3\frac{\cos\left(2\pi x\right)+\frac{1}{2}}{\pi^{4}x^{4}}+\frac{1/2}{\pi^{2}x^{2}}
\]
for $x\in\R\setminus\left\{ 0\right\} $. Since all terms except for
the last one decay strictly faster than $x^{-2}$ as $\left|x\right|\to\infty$,
we get $\left|\left(f_{1}\ast f_{2}\right)\left(x\right)\right|=\left|h\left(x\right)\right|\asymp x^{-2}$
for large $\left|x\right|$.

Thus, $f_{1}\ast f_{2}\notin L^{p}\left(\R\right)$ for $\frac{1}{3}<p\leq\frac{1}{2}$,
although $f_{1}\in L^{1}\left(\R\right)$ and $f_{2}\in L^{p}\left(\R\right)$
for these values of $p$. This shows that the usual form of Young's
inequality fails for $p<1$, even if we assume the factors of the
convolution to be bandlimited.
\end{example}

In the example above, one should note that $f_{1}\ast f_{2}\in L^{p}\left(\R\right)$
holds for $\frac{1}{2}<p<1$, which is exactly the range of $p\in\left(0,1\right)$,
for which $f_{1},f_{2}\in L^{p}\left(\R\right)$ holds. This is indicative
of the kind of convolution relation that we are after (see also equation~(\ref{eq:QuasiBanachConvolutionNoDetails})).

For technical reasons, we start our derivation of the convolution
relations for bandlimited $L^{p}$-functions with $p\in\left(0,1\right)$
by showing that we can always approximate such functions by Schwartz
functions in a suitable way. This will allow us to restrict our attention
to Schwartz functions for most of our proofs.
\begin{lem}
\label{lem:BandlimitedDistributionApproximation}Let $\Omega\subset\R^{\dimension}$
be compact and assume that $f\in\Schwartz'\left(\R^{\dimension}\right)$
is a tempered distribution with compact Fourier support $\supp\widehat{f}\subset\Omega$.

Then $f$ is given by (integration against) a smooth function $g\in C^{\infty}\left(\R^{\dimension}\right)$
with polynomially bounded derivatives of all orders.

Furthermore, there is a sequence of Schwartz functions $\left(g_{n}\right)_{n\in\N}$
with the following properties:

\begin{enumerate}
\item $\supp\widehat{g_{n}}\subset B_{1/n}\left(\Omega\right)$, where $B_{1/n}\left(\Omega\right)$
is the $\frac{1}{n}$-neighborhood of $\Omega$, given by
\[
B_{1/n}\left(\Omega\right)=\left\{ \xi\in\R^{\dimension}\with\dist\left(\xi,\Omega\right)<n^{-1}\right\} .
\]
\item $\left|g_{n}\left(x\right)\right|\leq\left|g\left(x\right)\right|$
for all $x\in\R^{\dimension}$, and
\item $g_{n}\left(x\right)\xrightarrow[n\to\infty]{}g\left(x\right)$ for
all $x\in\R^{\dimension}$.
\end{enumerate}
In particular, $\left\Vert g_{n}-f\right\Vert _{L^{p}}=\left\Vert g_{n}-g\right\Vert _{L^{p}}\xrightarrow[n\to\infty]{}0$
for any $p\in\left(0,\infty\right)$ for which $f\in L^{p}\left(\R^{\dimension}\right)$
is true.
\end{lem}

\begin{rem*}
In the following, we will always identify $f$ with its ``smooth
version'' $g$, i.e.\@ we will write $f\left(x\right)$ instead
of $g\left(x\right)$ for $x\in\R^{\dimension}$. In particular, it
thus makes sense to ask whether $f\in L^{p}\left(\R^{\dimension}\right)$,
or to calculate $\left\Vert f\right\Vert _{L^{p}}\in\left[0,\infty\right]$,
since $f=g$ is a pointwise defined (even smooth) function.
\end{rem*}
\begin{proof}
(based upon the proof of \cite[Theorem 1.4.1]{TriebelTheoryOfFunctionSpaces})
In case of $\Omega=\emptyset$, we have $f=0$, so that all claims
are trivial (with $g_{n}\equiv0$). Hence, we can assume $\Omega\neq\emptyset$
in the following. That $f$ is given by (integration against) a smooth
(even analytic) function $g\in C^{\infty}\left(\R^{\dimension}\right)$
with polynomially bounded derivatives of all orders is a consequence
of the \emph{Paley-Wiener Theorem} (cf.\@ \cite[Theorem 7.23]{RudinFA}).
Hence, we only need to establish existence of the sequence $\left(g_{n}\right)_{n\in\N}$.

To this end, choose $\psi\in C_{c}^{\infty}\left(B_{1}\left(0\right)\right)$
with $0\leq\psi\leq1$ and $\psi\left(0\right)=1$ and set $\widetilde{\varphi}:=\mathcal{F}^{-1}\psi\in\Schwartz\left(\R^{\dimension}\right)$.
Note that
\[
\gamma:=\widetilde{\varphi}\left(0\right)=\int_{\R^{\dimension}}\psi\left(\xi\right)\,\d\xi>0,
\]
since $\psi$ is continuous and nonnegative with $\psi\not\equiv0$.
Finally, define $\varphi:=\widetilde{\varphi}/\gamma$ and observe
$\varphi\left(0\right)=1$, as well as
\[
\supp\widehat{\varphi}=\supp\left(\Fourier\widetilde{\varphi}\right)=\supp\psi\subset B_{1}\left(0\right).
\]
Now, since the Fourier transform $\widehat{\varphi}=\psi/\gamma$
is nonnegative, $\left|\varphi\right|$ attains its global maximum
at the origin, since
\begin{align*}
\left|\varphi\left(x\right)\right|=\left|\left(\Fourier^{-1}\widehat{\varphi}\right)\left(x\right)\right|=\frac{1}{\gamma}\left|\left(\Fourier^{-1}\psi\right)\left(x\right)\right| & \leq\frac{1}{\gamma}\cdot\int_{\R^{\dimension}}\left|\psi\left(\xi\right)\cdot e^{2\pi i\left\langle x,\xi\right\rangle }\right|\,\d\xi\\
 & =\frac{1}{\gamma}\cdot\int_{\R^{\dimension}}\psi\left(\xi\right)\,\d\xi=\frac{1}{\gamma}\cdot\gamma=1=\varphi\left(0\right)\qquad\forall\,\,x\in\R^{\dimension}\,.
\end{align*}

Let us now set
\[
g_{n}:\R^{\dimension}\to\Compl,x\mapsto g\left(x\right)\cdot\varphi\left(\frac{x}{n}\right)
\]
for $n\in\N$. Since $g$ is a smooth function with polynomially bounded
derivatives, the Leibniz rule $\partial^{\alpha}\left(fg\right)=\sum_{\beta\leq\alpha}\binom{\alpha}{\beta}\cdot\partial^{\beta}f\cdot\partial^{\alpha-\beta}g$
easily implies that all derivatives of $g_{n}$ decay faster than
$\left(1+\left|\mybullet\right|\right)^{-N}$ for any $N\in\N$, so
that $g_{n}$ is indeed a Schwartz function. As seen above, $\left|\varphi\left(x\right)\right|\leq1$
for all $x\in\R^{\dimension}$, which immediately yields
\[
\left|g_{n}\left(x\right)\right|=\left|g\left(x\right)\cdot\varphi\left(\frac{x}{n}\right)\right|\leq\left|g\left(x\right)\right|.
\]
Furthermore, the continuity of $\varphi$, together with $\varphi\left(0\right)=1$
shows
\[
g_{n}\left(x\right)=g\left(x\right)\cdot\varphi\left(x/n\right)\xrightarrow[n\to\infty]{}g\left(x\right)\cdot\varphi\left(0\right)=g\left(x\right).
\]

Finally, the convolution theorem (cf.\@ \cite[Theorem 7.19(e)]{RudinFA}
and see \cite[Definition 7.18]{RudinFA} for the definition of the
convolution of a tempered distribution and a Schwartz function) yields
\[
\widehat{g_{n}}=\widehat{g}\ast\widehat{\varphi\left(\frac{\mybullet}{n}\right)}=n^{\dimension}\cdot\left[\widehat{g}\ast\left(\widehat{\varphi}\left(n\mybullet\right)\right)\right].
\]
Since $\widehat{\varphi}=\psi/\gamma$ is supported in $B_{r}\left(0\right)$
for some $r\in\left(0,1\right)$, we see that $\widehat{\varphi}\left(n\mybullet\right)$
has support in $B_{r/n}\left(0\right)$.

By definition of convolution for tempered distributions, we have
\[
\left[\widehat{g}\ast\left(\widehat{\varphi}\left(n\mybullet\right)\right)\right]\left(\xi\right)=\left\langle \widehat{g},\,L_{\xi}\left(\left[\widehat{\varphi}\left(n\mybullet\right)\right]^{\vee}\right)\right\rangle _{\Schwartz'}\qquad\forall\,\,\xi\in\R^{\dimension}\,,
\]
with $\theta^{\vee}\left(\xi\right)=\theta\left(-\xi\right)$. But
since $\left[\widehat{\varphi}\left(n\mybullet\right)\right]^{\vee}$is
supported in $B_{r/n}\left(0\right)$, we see that $L_{\xi}\left(\left[\widehat{\varphi}\left(n\mybullet\right)\right]^{\vee}\right)$
is supported in $B_{r/n}\left(\xi\right)$. Thus, $\left[\widehat{g}\ast\widehat{\varphi}\left(n\mybullet\right)\right]\left(\xi\right)\neq0$
can only happen if $\supp\widehat{g}\cap B_{r/n}\left(\xi\right)\neq\emptyset$,
i.e.\@ for $\xi\in B_{r/n}\left(\Omega\right)$, since $\supp\widehat{g}=\supp\widehat{f}\subset\Omega$.
Hence,
\[
\supp\widehat{g_{n}}\subset\overline{B_{r/n}\left(\Omega\right)}\subset B_{1/n}\left(\Omega\right).
\]

The additional claim regarding the $L^{p}$-convergence is a direct
consequence of the dominated convergence theorem, using the pointwise
convergence, together with $\left|g_{n}\right|\leq\left|g\right|$.
\end{proof}
Using the preceding approximation result, we will now show that bandlimited
functions which are integrable to a power $p$ are automatically integrable
to every power $q\geq p$. This roughly reflects the fact that these
functions are always (uniformly) continuous, so that the only obstruction
to integrability is the decay at infinity. This embedding into $L^{q}$-spaces
with larger $q$ will be a central ingredient for the proof of the
convolution relations in the quasi-Banach regime $p\in\left(0,1\right)$.

We remark that the corollary below is a special case of \cite[1.4.1(3)]{TriebelTheoryOfFunctionSpaces}.
\begin{cor}
\label{cor:BandlimitedEmbedding}Let $\Omega\subset\R^{\dimension}$
be compact and assume that $f\in\mathcal{S}'\left(\R^{\dimension}\right)$
is a tempered distribution with compact Fourier support $\supp\widehat{f}\subset\Omega$.
Then the following hold:

\begin{enumerate}
\item If $f\in L^{p}\left(\R^{\dimension}\right)$ holds for some $p\in\left(0,2\right]$,
then
\begin{equation}
\left\Vert f\right\Vert _{L^{q}}\leq\left[\lambda\left(\Omega\right)\right]^{\frac{1}{p}-\frac{1}{q}}\cdot\left\Vert f\right\Vert _{L^{p}}\label{eq:BandlimitedHigherIntegrabilityQuasiBanach}
\end{equation}
holds for all $q\in\left[p,\infty\right]$.
\item There is a constant $K=K\left(\Omega\right)>0$ such that
\begin{equation}
\left\Vert f\right\Vert _{L^{q}}\leq K\cdot\left\Vert f\right\Vert _{L^{p}}\label{eq:BandlimitedHigherIntegrabilityBanach}
\end{equation}
holds for all $p\in\left[1,\infty\right]$, $q\in\left[p,\infty\right]$
and $f\in L^{p}\left(\R^{\dimension}\right)$ with $\supp\widehat{f}\subset\Omega$.\qedhere
\end{enumerate}
\end{cor}

\begin{proof}
We first observe that it suffices to establish the stated estimates
for Schwartz functions, since the Fatou property for $L^{q}\left(\R^{\dimension}\right)$
implies, using the sequence $\left(g_{n}\right)_{n\in\N}$ given by
Lemma~\ref{lem:BandlimitedDistributionApproximation}, that
\[
\left\Vert f\right\Vert _{L^{q}}\leq\liminf_{n\to\infty}\left\Vert g_{n}\right\Vert _{L^{q}}\leq\liminf_{n\to\infty}\left[\lambda\left(B_{1/n}\left(\Omega\right)\right)\right]^{\frac{1}{p}-\frac{1}{q}}\cdot\left\Vert g_{n}\right\Vert _{L^{p}}\leq\left[\lambda\left(\Omega\right)\right]^{\frac{1}{p}-\frac{1}{q}}\cdot\left\Vert f\right\Vert _{L^{p}}.
\]
This proves equation~(\ref{eq:BandlimitedHigherIntegrabilityQuasiBanach})
in the general case. Observe that the last step used $\left|g_{n}\right|\leq\left|f\right|$
as well as $\lambda\left(B_{1/n}\left(\Omega\right)\right)\to\lambda\left(\Omega\right)$.
This last convergence is a consequence of the continuity of the (Lebesgue)
measure from above (cf.\@ \cite[Theorem 1.8(d)]{FollandRA}) and
of the identity $\Omega=\bigcap_{n\in\N}B_{1/n}\left(\Omega\right)$,
which holds because $\Omega$ is compact.

Validity of equation~(\ref{eq:BandlimitedHigherIntegrabilityBanach})
in the general case is derived similarly, with $K=K\big(\,\overline{B_{1}\left(\Omega\right)}\,\big)$
instead of $K=K\left(\Omega\right)$.

Hence, we can assume $f\in\Schwartz\left(\R^{\dimension}\right)$
with $\supp\widehat{f}\subset\Omega$ in the following. Observe that
by regularity of the Lebesgue measure, there is for each $\varepsilon>0$
some open set $U_{\varepsilon}\supset\Omega$ with $\lambda\left(U_{\varepsilon}\right)<\lambda\left(\Omega\right)+\varepsilon$.
Furthermore, by the $C^{\infty}$ Urysohn Lemma (see for example \cite[Lemma 8.18]{FollandRA})
there is some $\gamma_{\varepsilon}\in C_{c}^{\infty}\left(U_{\varepsilon}\right)$
with $0\leq\gamma_{\varepsilon}\leq1$ and $\gamma_{\varepsilon}\equiv1$
on a neighborhood of $\Omega$. Set $\psi_{\varepsilon}:=\Fourier^{-1}\gamma_{\varepsilon}\in\Schwartz\left(\R^{\dimension}\right)$.

We first observe that the Hausdorff-Young inequality (cf.\@ \cite[Theorem 8.30]{FollandRA})
yields 
\begin{equation}
\left\Vert \psi_{\varepsilon}\right\Vert _{L^{p'}}\leq\left\Vert \gamma_{\varepsilon}\right\Vert _{L^{p}}\leq\left\Vert \Indicator_{U_{\varepsilon}}\right\Vert _{L^{p}}=\left[\lambda\left(U_{\varepsilon}\right)\right]^{1/p}\leq\left[\lambda\left(\Omega\right)+\varepsilon\right]^{1/p}\label{eq:BandlimitedEmbeddingIntoHigherLPHausdorffYoungApplication}
\end{equation}
for all $p\in\left[1,2\right]$, where $p'\in\left[2,\infty\right]$
is the conjugate exponent, i.e.\@ $p'=\frac{p}{p-1}$ for $p\in\left(1,2\right]$
and $p'=\infty$ for $p=1$.

Using the support assumption on $\widehat{f}$, we get $\widehat{f}=\widehat{f}\cdot\gamma_{\varepsilon}$,
and hence 
\[
f=\Fourier^{-1}\widehat{f}=\Fourier^{-1}\left[\,\smash{\widehat{f}}\cdot\gamma_{\varepsilon}\right]=f\ast\psi_{\varepsilon}.
\]
For the proof of equation~(\ref{eq:BandlimitedHigherIntegrabilityQuasiBanach}),
we first note that it suffices to establish the case $q=\infty$;
indeed, in the remaining case $p\leq q<\infty$, we get
\[
\left\Vert f\right\Vert _{L^{q}}^{q}=\int_{\R^{\dimension}}\left|f\right|^{p}\cdot\left|f\right|^{q-p}\,\d x\leq\left\Vert f\right\Vert _{L^{\infty}}^{q-p}\cdot\int_{\R^{\dimension}}\left|f\right|^{p}\,\d x\leq\left[\lambda\left(\Omega\right)\right]^{\frac{q-p}{p}}\cdot\left\Vert f\right\Vert _{L^{p}}^{q-p}\cdot\left\Vert f\right\Vert _{L^{p}}^{p},
\]
where the last step made use of the case $q=\infty$. Taking $q$-th
roots completes the proof of equation~(\ref{eq:BandlimitedHigherIntegrabilityQuasiBanach}).

For the case $q=\infty$, we distinguish two sub-cases:

\textbf{Case 1}: If $p\in\left[1,2\right]$, Hölder's inequality implies
\begin{align}
\left|f\left(x\right)\right|=\left|\left(f\ast\psi_{\varepsilon}\right)\left(x\right)\right|=\left|\int_{\R^{\dimension}}f\left(y\right)\cdot\psi_{\varepsilon}\left(x-y\right)\,\d y\right| & \leq\left\Vert f\right\Vert _{L^{p}}\cdot\left\Vert \psi_{\varepsilon}\left(x-\mybullet\right)\right\Vert _{L^{p'}}\label{eq:BandlimitedHigherIntegrabilityHoelder}\\
 & \leq\left[\lambda\left(\Omega\right)+\varepsilon\right]^{1/p}\cdot\left\Vert f\right\Vert _{L^{p}}\nonumber 
\end{align}
for all $x\in\R^{\dimension}$. Here, equation~(\ref{eq:BandlimitedEmbeddingIntoHigherLPHausdorffYoungApplication})
was used in the last step. Since $\varepsilon>0$ was arbitrary, this
is nothing but equation~(\ref{eq:BandlimitedHigherIntegrabilityQuasiBanach})
for $q=\infty$ and $p\in\left[1,2\right]$.

\textbf{Case 2}: $p\in\left(0,1\right)$. Here, we cannot use Hölder's
inequality as in the last case. Instead, we apply a ``flop''. In
the present context, this means that we will derive an estimate of
the form $Q\leq C\cdot Q^{r}$ for some exponent $r\in\left(0,1\right)$
with $Q:=\sup_{y\in\R^{\dimension}}\left|f\left(y\right)\right|$.
By rearranging, this yields $Q\leq C^{\frac{1}{1-r}}$, which will
imply the desired estimate. Here it is important to note that $Q$
is indeed a finite quantity because of $f\in\Schwartz\left(\R^{\dimension}\right)$.

For the execution of this plan, note
\begin{align*}
\left|f\left(x\right)\right|=\left|\int_{\R^{\dimension}}f\left(y\right)\cdot\psi_{\varepsilon}\left(x-y\right)\,\d y\right| & \leq\int_{\R^{\dimension}}\left|f\left(y\right)\right|^{p}\cdot\left|f\left(y\right)\right|^{1-p}\cdot\left|\psi_{\varepsilon}\left(x-y\right)\right|\,\d y\\
 & \leq\left\Vert \psi_{\varepsilon}\right\Vert _{L^{\infty}}\cdot\left\Vert f\right\Vert _{L^{p}}^{p}\cdot\sup_{y\in\R^{\dimension}}\left|f\left(y\right)\right|^{1-p}\\
 & \leq\left(\lambda\left(\Omega\right)+\varepsilon\right)\cdot\left\Vert f\right\Vert _{L^{p}}^{p}\cdot Q^{1-p}.
\end{align*}
Here, we again used equation~(\ref{eq:BandlimitedEmbeddingIntoHigherLPHausdorffYoungApplication})
(with $p=1$) in the last step. Furthermore, it is important to observe
that we can indeed interchange the supremum and the power $1-p$ because
of $1-p>0$.

Taking the supremum over $x\in\R^{\dimension}$ on the left-hand side
yields
\[
\sup_{x\in\R^{\dimension}}\left|f\left(x\right)\right|\leq\left(\lambda\left(\Omega\right)+\varepsilon\right)\cdot\left\Vert f\right\Vert _{L^{p}}^{p}\cdot\left(\sup_{x\in\R^{\dimension}}\left|f\left(x\right)\right|\right)^{1-p}.
\]
Rearranging, taking $p$-th roots and letting $\varepsilon\downarrow0$
completes the proof of equation~(\ref{eq:BandlimitedHigherIntegrabilityQuasiBanach})
for $p\in\left(0,1\right)$ and $q=\infty$.

\medskip{}

It remains to establish estimate~(\ref{eq:BandlimitedHigherIntegrabilityBanach}).
To this end, we simply note that Hölder's inequality implies as in
equation~(\ref{eq:BandlimitedHigherIntegrabilityHoelder}) that $\left|f\left(x\right)\right|\leq\left\Vert f\right\Vert _{L^{p}}\cdot\left\Vert \psi_{1}\left(x-\mybullet\right)\right\Vert _{L^{p'}}\leq K\cdot\left\Vert f\right\Vert _{L^{p}}$,
with $K:=\max\left\{ 1,\left\Vert \psi_{1}\right\Vert _{L^{1}},\left\Vert \psi_{1}\right\Vert _{L^{\infty}}\right\} $.
The last estimate used $\left\Vert \psi\right\Vert _{L^{r}}\leq\max\left\{ \left\Vert \psi\right\Vert _{L^{1}},\left\Vert \psi\right\Vert _{L^{\infty}}\right\} $
for all $r\in\left[1,\infty\right]$, cf.\@ \cite[Proposition 6.10]{FollandRA}.
This establishes the claim for $q=\infty$.

In general, for $1\leq p\leq q<\infty$, we have
\[
\left\Vert f\right\Vert _{L^{q}}^{q}=\int_{\R^{\dimension}}\left|f\right|^{p}\cdot\left|f\right|^{q-p}\,\d x\leq\left\Vert f\right\Vert _{L^{\infty}}^{q-p}\cdot\left\Vert f\right\Vert _{L^{p}}^{p}\leq K^{q-p}\cdot\left\Vert f\right\Vert _{L^{p}}^{q-p}\cdot\left\Vert f\right\Vert _{L^{p}}^{p}.
\]
Rearranging yields $\left\Vert f\right\Vert _{L^{q}}\leq K^{1-\frac{p}{q}}\cdot\left\Vert f\right\Vert _{L^{p}}$.
Because of $K\geq1$ and $0\leq1-\frac{p}{q}\leq1$, we have $K^{1-\frac{p}{q}}\leq K$,
which completes the proof.
\end{proof}
Now, we can finally establish the convolution relation for bandlimited
functions in $L^{p}\left(\R^{\dimension}\right)$ for $p\in\left(0,1\right)$.
The result is essentially taken from Triebel \cite[Proposition 1.5.1]{TriebelTheoryOfFunctionSpaces};
but Triebel does not state the form of the constant $\left[\lambda\left(Q-\Omega\right)\right]^{\frac{1}{p}-1}$
explicitly in the statement of the theorem. This constant, however,
will be important for us; see for instance the definition of the weight
$v$ (for $q_{k}<1$) in Theorem~\ref{thm:NoSubordinatenessWithConsiderationOfOverlaps}.
\begin{thm}
\label{thm:QuasiBanachConvolution}(cf.\@ \cite[Proposition 1.5.1]{TriebelTheoryOfFunctionSpaces})
Let $\Omega,Q\subset\R^{\dimension}$ be compact and $p\in\left(0,1\right]$.
Furthermore, let $\psi\in L^{1}\left(\R^{\dimension}\right)$ with
$\supp\psi\subset Q$ and such that $\Fourier^{-1}\psi\in L^{p}\left(\R^{\dimension}\right)$.

For each $f\in L^{p}\left(\R^{\dimension}\right)\cap\Schwartz'\left(\R^{\dimension}\right)$
with Fourier support $\supp\widehat{f}\subset\Omega$, we have $\widehat{f}\in C_{0}\left(\R^{\dimension}\right)$,
and
\[
\Fourier^{-1}\left(\,\smash{\psi\cdot\widehat{f}}\,\,\right)=\left(\smash{\Fourier^{-1}\psi}\right)\ast f\in L^{p}\left(\smash{\R^{\dimension}}\right)
\]
with (quasi)-norm estimate
\begin{align*}
\left\Vert \Fourier^{-1}\left(\,\smash{\psi\cdot\widehat{f}}\,\,\right)\right\Vert _{L^{p}}=\left\Vert \left(\smash{\Fourier^{-1}\psi}\right)\ast f\right\Vert _{L^{p}} & \leq\left\Vert \left|\smash{\Fourier^{-1}\psi}\right|\ast\left|f\right|\right\Vert _{L^{p}}\\
 & \leq\left[\lambda\left(Q-\Omega\right)\right]^{\frac{1}{p}-1}\cdot\left\Vert \smash{\Fourier^{-1}\psi}\right\Vert _{L^{p}}\cdot\left\Vert f\right\Vert _{L^{p}},
\end{align*}
where
\[
Q-\Omega=\left\{ q-\omega\with q\in Q\text{ and }\omega\in\Omega\right\} 
\]
is the \textbf{algebraic difference }of $Q$ and $\Omega$, which
is compact and hence measurable.

Finally, we have the pointwise estimate
\begin{equation}
\left|\left[\left(\smash{\Fourier^{-1}\psi}\right)\ast f\right]\left(x\right)\right|\leq\left[\left|\smash{\Fourier^{-1}\psi}\right|\ast\left|f\right|\right]\left(x\right)\leq\left[\lambda\left(Q-\Omega\right)\right]^{\frac{1}{p}-1}\left[\int_{\R^{\dimension}}\left|\left(\smash{\Fourier^{-1}\psi}\right)\left(y\right)\cdot f\left(x-y\right)\right|^{p}\,\d y\right]^{\frac{1}{p}}\label{eq:QuasiBanachConvolutionPointwise}
\end{equation}
for all $x\in\R^{\dimension}$.
\end{thm}

\begin{rem*}
We observe that $\psi\cdot\widehat{f}\in\Schwartz'\left(\R^{\dimension}\right)$
is well-defined, even though $\psi$ might not be a smooth function.
This is because Corollary~\ref{cor:BandlimitedEmbedding}—together
with $p\leq1$—yields $f\in L^{1}\left(\R^{\dimension}\right)$ and
thus $\widehat{f}\in C_{0}\left(\R^{\dimension}\right)$. Hence, $\psi\cdot\widehat{f}\in L^{1}\left(\R^{\dimension}\right)\hookrightarrow\Schwartz'\left(\R^{\dimension}\right)$
because $\psi\in L^{1}\left(\R^{\dimension}\right)$.
\end{rem*}
\begin{proof}
The Riemann-Lebesgue lemma yields $\Fourier^{-1}\psi\in C_{0}\left(\R^{\dimension}\right)$
because $\psi\in L^{1}\left(\R^{\dimension}\right)$. Thus, we get
$\Fourier^{-1}\psi\in L^{\infty}\left(\R^{\dimension}\right)\cap L^{p}\left(\R^{\dimension}\right)\subset L^{q}\left(\R^{\dimension}\right)$
for all $q\in\left[p,\infty\right]$, cf.\@ \cite[Proposition 6.10]{FollandRA}.
For the same range of $q$, Corollary~\ref{cor:BandlimitedEmbedding}
yields $f\in L^{q}\left(\R^{\dimension}\right)$. In particular, $f\in L^{1}\left(\R^{\dimension}\right)$
and hence $\widehat{f}\in C_{0}\left(\R^{\dimension}\right)$.

Because $p\in\left(0,1\right]$, we get $f,\Fourier^{-1}\psi\in L^{1}\left(\R^{\dimension}\right)$
and thus also $\left(\Fourier^{-1}\psi\right)\ast f\in L^{1}\left(\R^{\dimension}\right)$
with Fourier transform
\[
\Fourier\left[\left(\smash{\Fourier^{-1}}\psi\right)\ast f\right]=\psi\cdot\widehat{f}.
\]
Further, $\psi\in L^{1}\left(\R^{\dimension}\right)$ together with
$\widehat{f}\in C_{0}\left(\R^{\dimension}\right)$ yields $\psi\cdot\widehat{f}\in L^{1}\left(\R^{\dimension}\right)$.
By Fourier inversion, this implies as claimed that $\Fourier^{-1}\left(\,\smash{\psi\cdot\widehat{f}}\,\,\right)=\left(\smash{\Fourier^{-1}}\psi\right)\ast f$.

We first observe that the pointwise estimate~(\ref{eq:QuasiBanachConvolutionPointwise})
implies the remaining claims. Indeed, equation~(\ref{eq:QuasiBanachConvolutionPointwise}),
together with Fubini's theorem, yields 
\begin{align*}
\left\Vert \left(\smash{\Fourier^{-1}\psi}\right)\ast f\right\Vert _{L^{p}}^{p} & \leq\left\Vert \left|\smash{\Fourier^{-1}}\psi\right|\ast\left|f\right|\right\Vert _{L^{p}}^{p}\\
 & \leq\left[\lambda\left(Q-\Omega\right)\right]^{1-p}\cdot\int_{\R^{\dimension}}\int_{\R^{\dimension}}\left|\left(\Fourier^{-1}\psi\right)\left(y\right)\right|^{p}\cdot\left|f\left(x-y\right)\right|^{p}\,\d y\,\d x\\
 & =\left[\lambda\left(Q-\Omega\right)\right]^{1-p}\cdot\int_{\R^{\dimension}}\left|\left(\Fourier^{-1}\psi\right)\left(y\right)\right|^{p}\cdot\int_{\R^{\dimension}}\left|f\left(x-y\right)\right|^{p}\,\d x\,\d y\\
 & =\left[\lambda\left(Q-\Omega\right)\right]^{1-p}\cdot\left\Vert \smash{\Fourier^{-1}}\psi\right\Vert _{L^{p}}^{p}\cdot\left\Vert f\right\Vert _{L^{p}}^{p},
\end{align*}
so that taking $p$-th roots completes the proof.

In order to prove equation~(\ref{eq:QuasiBanachConvolutionPointwise}),
fix $x\in\R^{\dimension}$. We will write $g^{\vee}\left(y\right)=g\left(-y\right)$
for arbitrary functions $g:\R^{\dimension}\to\Compl$. Using $f\in L^{\infty}\left(\R^{\dimension}\right)$
and $\Fourier^{-1}\psi\in L^{p}\left(\R^{\dimension}\right)\cap L^{1}\left(\R^{\dimension}\right)$,
we see
\[
F_{x}:=\left(\Fourier^{-1}\psi\right)\cdot L_{x}\,f^{\vee}\in L^{1}\left(\smash{\R^{\dimension}}\right)\cap L^{p}\left(\smash{\R^{\dimension}}\right),
\]
so that the Fourier transform $\Fourier F_{x}\in C_{0}\left(\R^{\dimension}\right)$
is well-defined with
\begin{align*}
\left(\Fourier F_{x}\right)\left(\xi\right) & =\int_{\R^{\dimension}}\widehat{\psi}\left(-y\right)\cdot f\left(x-y\right)\cdot e^{-2\pi i\left\langle y,\xi\right\rangle }\,\d y=\int_{\R^{\dimension}}\widehat{\psi}\left(z\right)\cdot f\left(x+z\right)\cdot e^{2\pi i\left\langle z,\xi\right\rangle }\,\d z\\
 & =\int_{\R^{\dimension}}\psi\left(y\right)\cdot\Fourier\left[z\mapsto f\left(x+z\right)\cdot e^{2\pi i\left\langle z,\xi\right\rangle }\right]\left(y\right)\,\d y\qquad\forall\,\,\xi\in\R^{\dimension}\,,
\end{align*}
where the last step used $\psi,f\in L^{1}\left(\R^{\dimension}\right)$
and that $\int\widehat{f}g=\int f\widehat{g}$ for $f,g\in L^{1}\left(\R^{\dimension}\right)$,
cf.\@ \cite[Lemma 8.25]{FollandRA}.

By elementary properties of the Fourier transform,
\[
\Fourier\left[z\mapsto f\left(x+z\right)\cdot e^{2\pi i\left\langle z,\xi\right\rangle }\right]\left(y\right)=\Fourier\left[M_{\xi}L_{-x}f\right]\left(y\right)=\left(\,L_{\xi}M_{x}\widehat{f}\,\,\right)\left(y\right),
\]
which yields
\[
\left(\Fourier F_{x}\right)\left(\xi\right)=\int_{\R^{\dimension}}\psi\left(y\right)\cdot\left(\,\smash{M_{x}\widehat{f}}\,\,\right)\left(y-\xi\right)\,\d y=\left[\psi\ast\left(\,\smash{M_{x}\widehat{f}}\,\,\right)^{\vee}\right]\left(\xi\right)
\]
and hence\footnote{Note that the distribution $\widehat{f}$ is assumed to satisfy $\supp\widehat{f}\subset\Omega$.
Because of $f\in L^{1}\left(\R^{\dimension}\right)$, $\widehat{f}$
is given by integration against a continuous bounded function. It
is then easy to see that this continuous function also has support
in $\Omega$. This justifies the calculation in equation~(\ref{eq:QuasiBanachConvolutionSupport}).}
\begin{align}
\supp\left(\Fourier F_{x}\right) & \subset\supp\psi+\supp\left(\left[\,\smash{M_{x}\widehat{f}}\,\,\right]^{\vee}\right)\subset Q+\left(-\Omega\right)=Q-\Omega\,.\label{eq:QuasiBanachConvolutionSupport}
\end{align}

Using Corollary~\ref{cor:BandlimitedEmbedding} and $0<p\leq1$,
we finally arrive at
\begin{align*}
\left|\left[\left(\smash{\Fourier^{-1}}\psi\right)\ast f\right]\left(x\right)\right|\leq\left(\left|\smash{\Fourier^{-1}}\psi\right|\ast\left|f\right|\right)\left(x\right) & =\int_{\R^{\dimension}}\left|\left(\smash{\Fourier^{-1}}\psi\right)\left(y\right)\cdot f\left(x-y\right)\right|\,\d y\\
 & =\left\Vert F_{x}\right\Vert _{L^{1}}\leq\left[\lambda\left(Q-\Omega\right)\right]^{\frac{1}{p}-\frac{1}{1}}\cdot\left\Vert F_{x}\right\Vert _{L^{p}},
\end{align*}
which is nothing but the pointwise estimate~(\ref{eq:QuasiBanachConvolutionPointwise}).
\end{proof}

\subsection{Definition of decomposition spaces}

\label{subsec:DecompositionDefinition}In this subsection, we will
formally define the (Fourier side) decomposition space $\FourierDecompSp{\CalQ}pY$
and its ``space side'' version. To this end, we first introduce
the class of partitions of unity which will turn out to be suitable
for the definition of decomposition spaces. We begin with the case
$p\in\left[1,\infty\right]$, since we will have to place more restrictive
assumptions on the covering $\CalQ$ in the quasi-Banach regime $p\in\left(0,1\right)$;
see Definition~\ref{defn:QuasiBanachBAPUAndLpDecompositionCovering}.
\begin{defn}
\label{defn:BanachBAPUUndDecompositionCovering}(inspired by \cite[Definition 2.2]{DecompositionSpaces1})

Let $\emptyset\neq\CalO\subset\R^{\dimension}$ be open and let $\CalQ=\left(Q_{i}\right)_{i\in i}$
be an admissible covering of $\CalO$. A family $\Phi=\left(\varphi_{i}\right)_{i\in I}$
of functions on $\CalO$ is called an \textbf{$L^{p}$-bounded admissible
partition of unity} (\textbf{$L^{p}$-BAPU}) for $\CalQ$ for all
$1\leq p\leq\infty$, if

\begin{enumerate}
\item $\varphi_{i}\in\TestFunctionSpace{\CalO}$ for all $i\in I$,
\item $\sum_{i\in I}\varphi_{i}\left(\xi\right)=1$ for all $\xi\in\CalO$,
\item $\varphi_{i}\left(\xi\right)=0$ for all $\xi\in\R^{\dimension}\setminus Q_{i}$
for all $i\in I$, and
\item the following expression (then a constant) is finite:
\[
C_{\CalQ,\Phi,p}:=\sup_{i\in I}\left\Vert \Fourier^{-1}\varphi_{i}\right\Vert _{L^{1}}\,.
\]
\end{enumerate}
We say that an admissible covering $\CalQ$ of $\CalO$ is an \textbf{$L^{p}$-decomposition
covering} of $\CalO$ for all $1\leq p\leq\infty$ if there is an
$L^{p}$-BAPU $\Phi$ for $\CalQ$.
\end{defn}

\begin{rem*}
The term ``$L^{p}$-bounded'' used above does \emph{not} refer to
the fact that the $L^{p}$-norm of the $\varphi_{i}$ is uniformly
bounded, but to the fact that $\left(\varphi_{i}\right)_{i\in I}$
forms a uniformly bounded family of $L^{p}$ Fourier-multipliers,
as a consequence of Young's inequality $L^{1}\ast L^{p}\hookrightarrow L^{p}$.

Clearly, the constant $C_{\CalQ,\Phi,p}$ does \emph{not} actually
depend on $\CalQ$ and $p$; but below, we will introduce the concept
of an $L^{p}$-BAPU also for $p\in\left(0,1\right)$ and in this case,
the similarly defined constant will depend on $\CalQ$ and $p$. For
consistency, we write $C_{\CalQ,\Phi,p}$ also for $p\in\left[1,\infty\right]$.

Finally, we mention that the term ``BAPU'' goes back to Feichtinger
and Gröbner \cite{DecompositionSpaces1}. Note, however, that in \cite[Definition 2.2]{DecompositionSpaces1},
a BAPU is not required to be smooth. Up to this difference, our terminology
is a special case of \cite[Definition 2.2]{DecompositionSpaces1};
precisely, each $\Fourier L^{1}$-BAPU in the sense of \cite[Definition 2.2]{DecompositionSpaces1}
that consists of smooth functions is also an $L^{p}$-BAPU for all
$1\leq p\leq\infty$ in the sense of the definition above. We need
to impose the smoothness of the $\varphi_{i}$ to ensure compatibility
with the reservoir $\DistributionSpace{\CalO}$ that we will use to
define the (Fourier-side) decomposition spaces, see Definition~\ref{def:FourierSideDecompositionSpace};
in \cite{DecompositionSpaces1}, a different reservoir is used.
\end{rem*}
For the quasi-Banach regime $p\in\left(0,1\right)$, the following
definition will turn out to be suitable. Note that we assume the covering
$\CalQ$ to be semi-structured, whereas for $p\in\left[1,\infty\right]$
all we needed was an admissible covering.
\begin{defn}
\label{defn:QuasiBanachBAPUAndLpDecompositionCovering}(cf.\@ \cite[Definition 2]{BorupNielsenDecomposition})

Let $0<p<1$, let $\emptyset\neq\CalO\subset\R^{\dimension}$ be open
and assume that $\CalQ=\left(Q_{i}\right)_{i\in I}=\left(T_{i}Q_{i}'+b_{i}\right)_{i\in I}$
is a semi-structured covering of $\CalO$. We say that a family $\Phi=\left(\varphi_{i}\right)_{i\in I}$
is an \textbf{$L^{p}$-bounded admissible partition of unity} ($L^{p}$-\textbf{BAPU})
for $\CalQ$, if

\begin{enumerate}
\item $\varphi_{i}\in\TestFunctionSpace{\CalO}$ for all $i\in I$,
\item $\sum_{i\in I}\varphi_{i}\left(\xi\right)=1$ for all $\xi\in\CalO$,
\item $\varphi_{i}\left(\xi\right)=0$ for all $\xi\in\R^{\dimension}\setminus Q_{i}$
for all $i\in I$, and
\item the following expression (then a constant) is finite: 
\[
C_{\CalQ,\Phi,p}:=\sup_{i\in I}\left(\left|\det T_{i}\right|^{\frac{1}{p}-1}\cdot\left\Vert \smash{\Fourier^{-1}}\varphi_{i}\right\Vert _{L^{p}}\right)\,.
\]
\end{enumerate}
We say that $\CalQ$ is an \textbf{$L^{p}$-decomposition covering}
of $\CalO$ if there is an $L^{p}$-BAPU $\Phi$ for $\CalQ$.
\end{defn}

\begin{rem*}
We will see in Corollary~\ref{cor:LpBAPUsAreAlsoLqBAPUsForLargerq}
that every $L^{p}$-BAPU for $\CalQ$ is automatically an $L^{q}$-BAPU
for $\CalQ$ for all $q\in\left[p,\infty\right]$.
\end{rem*}
Before we finally give a formal definition of the decomposition space
$\FourierDecompSp{\CalQ}pY$, we first clarify our assumptions on
the space $Y$.
\begin{defn}
\label{defn:QRegularSequenceSpace}(cf.\@ \cite[Definition 2.5]{DecompositionSpaces1})

\begin{enumerate}[leftmargin=0.8cm]
\item Let $I$ be an index set. A quasi-normed vector space $\left(Y,\left\Vert \mybullet\right\Vert _{Y}\right)$
which is a subspace of $\Compl^{I}$ is called \textbf{solid}, or
a \textbf{solid sequence space over $I$}, if the following holds:
If $x=\left(x_{i}\right)_{i\in I}\in\Compl^{I}$ and $y=\left(y_{i}\right)_{i\in I}\in Y$
are arbitrary with $\left|x_{i}\right|\leq\left|y_{i}\right|$ for
all $i\in I$, then $x\in Y$ and $\left\Vert x\right\Vert _{Y}\leq\left\Vert y\right\Vert _{Y}$.\vspace{0.2cm}
\item Let $\CalQ=\left(Q_{i}\right)_{i\in I}$ be a covering of a set $X\neq\emptyset$.
We say that solid sequence space $\left(Y,\left\Vert \mybullet\right\Vert _{Y}\right)$
over $I$ is \textbf{$\CalQ$-regular}, if the following hold:

\begin{enumerate}
\item $Y$ is complete, i.e.\@ a quasi-Banach space,
\item $Y$ is \textbf{invariant under $\CalQ$-clustering}, i.e.\@, the
\textbf{$\CalQ$-clustering map}
\[
\Gamma_{\CalQ}:Y\to Y,\left(c_{i}\right)_{i\in I}\mapsto\left(c_{i}^{\ast}\right)_{i\in I}:=\left(\,\vphantom{\sum}\smash{\sum_{j\in i^{\ast}}}c_{j}\,\right)_{i\in I}\vphantom{\sum_{j\in i^{\ast}}}
\]
is well-defined and bounded.\qedhere
\end{enumerate}
\end{enumerate}
\end{defn}

\begin{rem*}
(1). We will see in Lemma~\ref{lem:SolidSequenceSpaceEmbedsIntoWeightedLInfty}
that every solid sequence space $Y\subset\Compl^{I}$ automatically
satisfies $Y\hookrightarrow\Compl^{I}$, i.e., each of the coordinate
evaluation mappings is continuous. Hence, each solid Banach sequence
space is a \textbf{solid BK-space} in the sense of \cite[Definition 2.4]{DecompositionSpaces1}.

\medskip{}

(2). If all of the sets $i^{\ast}$ are finite (which holds if $\CalQ$
is admissible), then the closed graph theorem (cf.\@ \cite[Theorem 2.15]{RudinFA}),
together with completeness of $Y$ and with the continuous embedding
$Y\hookrightarrow\Compl^{I}$ and the fact that the clustering map
is continuous with respect to the (Hausdorff!) product topology on
$\Compl^{I}$, imply that \emph{the clustering map $\Gamma_{\CalQ}$
is bounded iff it is well-defined}.

\medskip{}

(3) The most important example of $\CalQ$-regular sequence spaces
that we will consider are weighted $\ell^{q}$-spaces. As we will
see in Lemma~\ref{lem:ModeratelyWeightedSpacesAreRegular}, $\ell_{u}^{q}\left(I\right)$
is $\CalQ$-regular, whenever $u=\left(u_{i}\right)_{i\in I}$ is
$\CalQ$-moderate.
\end{rem*}
\begin{lem}
\label{lem:SolidSequenceSpaceEmbedsIntoWeightedLInfty}Let $I\neq\emptyset$
be a set and let $\left(Y,\left\Vert \mybullet\right\Vert _{Y}\right)$
be a solid sequence space over $I$. For $i\in I$, set
\[
u_{i}:=\begin{cases}
\left\Vert \delta_{i}\right\Vert _{Y}, & \text{if }\delta_{i}\in Y,\\
1, & \text{if }\delta_{i}\notin Y,
\end{cases}\qquad\text{where}\qquad\left(\delta_{i}\right)_{j}=0\text{ for }i\neq j,\text{ and }\left(\delta_{i}\right)_{i}=1\,.
\]
Then $Y\hookrightarrow\ell_{u}^{\infty}\left(I\right)\hookrightarrow\Compl^{I}$.
More precisely, we even have $\left\Vert x\right\Vert _{\ell_{u}^{\infty}}\leq\left\Vert x\right\Vert _{Y}$
for all $x=\left(x_{i}\right)_{i\in I}\in Y$.
\end{lem}

\begin{proof}
Let $x=\left(x_{i}\right)_{i\in I}\in Y$ be arbitrary and let $i\in I$.
In case of $x_{i}\neq0$, we have $\left|\smash{\left(\delta_{i}\right)_{j}}\vphantom{\sum}\right|\leq\left|x_{j}\right|/\left|x_{i}\right|$
for all $j\in I$. By solidity of $Y$, this implies $\delta_{i}\in Y$
with 
\[
\left\Vert \delta_{i}\right\Vert _{Y}\leq\left\Vert \frac{x}{\left|x_{i}\right|}\right\Vert _{Y}=\frac{\left\Vert x\right\Vert _{Y}}{\left|x_{i}\right|},
\]
whence $u_{i}\cdot\left|x_{i}\right|=\left\Vert \delta_{i}\right\Vert _{Y}\cdot\left|x_{i}\right|\leq\left\Vert x\right\Vert _{Y}$.
In case of $x_{i}=0$, we trivially have $u_{i}\cdot\left|x_{i}\right|=0\leq\left\Vert x\right\Vert _{Y}$.

Since $i\in I$ was arbitrary, we conclude $\left\Vert x\right\Vert _{\ell_{u}^{\infty}}\leq\left\Vert x\right\Vert _{Y}<\infty$,
which completes the proof.
\end{proof}
For later use, we state the following result which connects iterated
applications of the clustering map and summing over the clustered
index sets $i^{\ell\ast}$.
\begin{lem}
\label{lem:HigherOrderClusteringMap}Let $\CalQ=\left(Q_{i}\right)_{i\in I}$
be an admissible covering of a set $X\neq\emptyset$ and let $\Gamma_{\CalQ}:\Compl^{I}\to\Compl^{I},c\mapsto c^{\ast}$
denote the $\CalQ$-clustering map from Definition~\ref{defn:QRegularSequenceSpace}.

For  any $\ell\in\N$, we have
\begin{equation}
\left(\Gamma_{\CalQ}^{\ell}\,c\right)_{i}\geq\sum_{j\in i^{\ell\ast}}c_{j}\quad\text{ for all }i\in I\text{ and all sequences }\left(c_{i}\right)_{i\in I}\text{ with nonnegative terms.}\label{eq:HigherOrderClustering}
\end{equation}

In particular, if $Y\subset\Compl^{I}$ is $\CalQ$-regular and $\ell\in\N$
is arbitrary, then the \textbf{$\ell$-fold clustering map}
\[
\Theta_{\ell}:Y\to Y\quad\text{with}\quad\left(\Theta_{\ell}\,c\right)_{i}:=\vphantom{\sum_{j}^{A}}\sum_{j\in i^{\ell\ast}}c_{j}
\]
is well-defined and bounded with $\vertiii{\Theta_{\ell}}_{Y\to Y}\leq\vertiii{\Gamma_{\CalQ}}_{Y\to Y}^{\ell}$.
\end{lem}

\begin{proof}
We first prove estimate~(\ref{eq:HigherOrderClustering}) by induction
on $\ell\in\N$. For $\ell=1$, we have equality by definition of
$\Gamma_{\CalQ}$.

For the induction step, fix for each $j\in i^{\left(\ell+1\right)\ast}$
some $k_{j}\in i^{\ell\ast}$ with $j\in k_{j}^{\ast}$. Since each
term $c_{j}$ is nonnegative, this yields
\begin{align*}
\sum_{j\in i^{\left(\ell+1\right)\ast}}c_{j} & =\sum_{k\in i^{\ell\ast}}\:\sum_{\substack{j\in i^{\left(\ell+1\right)\ast}\\
\text{with }k_{j}=k
}
}c_{j}\\
\left({\scriptstyle \text{since }j\in k_{j}^{\ast}=k^{\ast}\text{ if }k_{j}=k}\right) & \leq\sum_{k\in i^{\ell\ast}}\:\sum_{j\in k^{\ast}}c_{j}=\sum_{k\in i^{\ell\ast}}\left(\Gamma_{\CalQ}\,c\right)_{k}\\
 & \overset{\left(\ast\right)}{\leq}\left(\Gamma_{\CalQ}^{\ell}\Gamma_{\CalQ}\,c\right)_{i}=\left(\Gamma_{\CalQ}^{\ell+1}\,c\right)_{i}\,.
\end{align*}
Here, we used the induction hypothesis (with the nonnegative(!) sequence
$\Gamma_{\CalQ}\,c$ instead of $c$) at $\left(\ast\right)$.

To prove boundedness of $\Theta_{\ell}$, let $c=\left(c_{i}\right)_{i\in I}\in Y$
be arbitrary and set $d:=\left(\left|c_{i}\right|\right)_{i\in I}$.
Observe $d\in Y$ with $\left\Vert d\right\Vert _{Y}=\left\Vert c\right\Vert _{Y}$,
since $Y$ is solid. As seen above, we have
\[
\left|\left(\Theta_{\ell}\,c\right)_{i}\right|\leq\sum_{j\in i^{\ell\ast}}\left|c_{j}\right|\leq\left(\Gamma_{\CalQ}^{\ell}\,d\right)_{i}\qquad\forall\,\,i\in I\,.
\]
Because of $\Gamma_{\CalQ}^{\ell}d\in Y$, solidity of $Y$ implies
$\Theta_{\ell}\,c\in Y$ with
\[
\left\Vert \Theta_{\ell}\,c\right\Vert _{Y}\leq\left\Vert \smash{\Gamma_{\CalQ}^{\ell}}\,d\right\Vert _{Y}\leq\vertiii{\smash{\Gamma_{\CalQ}^{\ell}}}_{Y\to Y}\left\Vert d\right\Vert _{Y}\leq\vertiii{\Gamma_{\CalQ}}_{Y\to Y}^{\ell}\cdot\left\Vert c\right\Vert _{Y}<\infty.\qedhere
\]
\end{proof}
Now, we are in a position to define the (Fourier-side) decomposition
spaces.
\begin{defn}
\label{def:FourierSideDecompositionSpace}Let $\emptyset\neq\CalO\subset\R^{\dimension}$
be an open set and let $p\in\left(0,\infty\right]$. Let $\CalQ=\left(Q_{i}\right)_{i\in I}$
be an $L^{p}$-decomposition covering of $\CalO$ with $L^{p}$-BAPU
$\Phi=\left(\varphi_{i}\right)_{i\in I}$, and let $Y\subset\Compl^{I}$
be $\CalQ$-regular.

For $f\in\DistributionSpace{\CalO}$, define
\[
\left\Vert f\right\Vert _{\FourierDecompSp{\CalQ}pY}:=\left\Vert f\right\Vert _{\CalD_{\Fourier,\Phi}\left(\CalQ,L^{p},Y\right)}:=\left\Vert \left(\left\Vert \Fourier^{-1}\left(\varphi_{i}f\right)\right\Vert _{L^{p}}\right)_{i\in I}\right\Vert _{Y}\in\left[0,\infty\right],
\]
with the convention that for a family $c=\left(c_{i}\right)_{i\in I}$
with $c_{i}\in\left[0,\infty\right]$, the expression $\left\Vert c\right\Vert _{Y}$
is to be read as $\infty$ if $c_{i}=\infty$ for some $i\in I$ or
if $c\notin Y$.

Define the \textbf{(Fourier-side) decomposition space} $\FourierDecompSp{\CalQ}pY$
with respect to the covering $\CalQ$, integrability exponent $p$
and global component $Y$ as
\[
\FourierDecompSp{\CalQ}pY:=\left\{ f\in\DistributionSpace{\CalO}\with\left\Vert f\right\Vert _{\FourierDecompSp{\CalQ}pY}<\infty\right\} .\qedhere
\]
\end{defn}

\begin{rem*}
We remark that $\varphi_{i}f$ is a distribution on $\CalO$ with
compact support in $\CalO$, which thus extends to a (tempered) distribution
on $\R^{\dimension}$. By Lemma~\ref{lem:BandlimitedDistributionApproximation},
this implies that $\mathcal{F}^{-1}\left(\varphi_{i}f\right)\in\Schwartz'\left(\R^{\dimension}\right)$
is given by (integration against) a smooth function. Thus, it makes
sense to write $\left\Vert \Fourier^{-1}\left(\varphi_{i}f\right)\right\Vert _{L^{p}}$,
with the caveat that this expression could be infinite.

We finally remark that the notations $\left\Vert \mybullet\right\Vert _{\FourierDecompSp{\CalQ}pY}$
and $\FourierDecompSp{\CalQ}pY$ suppress the family $\left(\varphi_{i}\right)_{i\in I}$
used to define the (quasi)-norm above. We will see below (see Theorem~\ref{cor:DecompositionSpaceWellDefined})
that the resulting space is independent of the chosen $L^{p}$-BAPU,
with equivalent quasi-norms for different choices, so that this is
justified.
\end{rem*}
For completeness, we also define ``space-side'' decomposition spaces.
To this end, we first introduce the reservoir $\SpaceReservoir{\CalO}$
which we will use for these spaces.
\begin{defn}
\label{def:UnifiedReservoir}For $\emptyset\neq\CalO\subset\R^{\dimension}$
open, we define
\[
\SpaceTestFunctions{\CalO}:=\Fourier\left(\TestFunctionSpace{\CalO}\right):=\left\{ \smash{\widehat{f}}\,\with\,f\in\TestFunctionSpace{\CalO}\right\} \subset\Schwartz\left(\smash{\R^{\dimension}}\right)
\]
and endow this space with the unique topology that makes the Fourier
transform
\[
\Fourier:\TestFunctionSpace{\CalO}\to\SpaceTestFunctions{\CalO}
\]
a homeomorphism.

We equip the topological dual space $\SpaceReservoir{\CalO}:=\left[\SpaceTestFunctions{\CalO}\right]'$
of $\SpaceTestFunctions{\CalO}$ with the weak-$\ast$-topology, i.e.,
with the topology of pointwise convergence on $\SpaceTestFunctions{\CalO}$.

Finally, as on the Schwartz space, we extend the Fourier transform
by duality to $\SpaceReservoir{\CalO}$, i.e.\@ we define
\begin{equation}
\Fourier:\SpaceReservoir{\CalO}\to\DistributionSpace{\CalO},f\mapsto f\circ\Fourier.\label{eq:FourierTransformOnUnifiedReservoir}
\end{equation}
As usual, we write $\widehat{f}:=\Fourier f$ for $f\in\SpaceReservoir{\CalO}$.
\end{defn}

\begin{rem*}
Since $\Fourier:\TestFunctionSpace{\CalO}\to\SpaceTestFunctions{\CalO}$
is a linear homeomorphism, the Fourier transform as defined in equation~(\ref{eq:FourierTransformOnUnifiedReservoir})
is easily seen to be a linear homeomorphism as well.
\end{rem*}
Now, we are in a position to define the \emph{space-side} decomposition
spaces.
\begin{defn}
\label{def:SpaceSideDecompositionSpaces}Let $\emptyset\neq\CalO\subset\R^{\dimension}$
be an open set and let $p\in\left(0,\infty\right]$. Let $\CalQ=\left(Q_{i}\right)_{i\in I}$
be an $L^{p}$-decomposition covering of $\CalO$ with $L^{p}$-BAPU
$\left(\varphi_{i}\right)_{i\in I}$, and let $Y\subset\Compl^{I}$
be $\CalQ$-regular.

For $f\in\SpaceReservoir{\CalO}$, set
\[
\left\Vert f\right\Vert _{\DecompSp{\CalQ}pY}:=\left\Vert \smash{\widehat{f}}\,\right\Vert _{\FourierDecompSp{\CalQ}pY}=\left\Vert \left(\left\Vert \Fourier^{-1}\left(\,\varphi_{i}\,\smash{\widehat{f}}\,\,\right)\right\Vert _{L^{p}}\right)_{i\in I}\right\Vert _{Y}\in\left[0,\infty\right]\,,
\]
and define the \textbf{space-side decomposition space} $\DecompSp{\CalQ}pY$
with respect to the covering $\CalQ$, integrability exponent $p$
and global component $Y$ as 
\[
\DecompSp{\CalQ}pY:=\left\{ f\in\SpaceReservoir{\CalO}\,\with\,\left\Vert f\right\Vert _{\DecompSp{\CalQ}pY}<\infty\right\} .\qedhere
\]
\end{defn}

\begin{rem}
\label{rem:TemperedDistributionsAsReservoirIncomplete}Since the Fourier
transform $\Fourier:\SpaceReservoir{\CalO}\to\DistributionSpace{\CalO}$
defined in equation~(\ref{eq:FourierTransformOnUnifiedReservoir})
is an isomorphism, it is clear that this Fourier transform restricts
to an (isometric) isomorphism
\[
\Fourier:\DecompSp{\CalQ}pY\to\FourierDecompSp{\CalQ}pY.
\]
Hence, for most purposes, it does not matter whether one considers
the ``space-side'' or the ``Fourier-side'' version of these spaces.
In particular, one has $\DecompSp{\CalQ}pY\hookrightarrow\DecompSp{\CalP}qZ$
if and only if $\FourierDecompSp{\CalQ}pY\hookrightarrow\FourierDecompSp{\CalP}qZ$.
But since the space $\DistributionSpace{\CalO}$ is a widely known
standard space, whereas $\SpaceReservoir{\CalO}$ is not, we will
generally prefer to work with the ``Fourier-side'' spaces. Note
in particular that most working analysts have acquired significant
intuition regarding which operations (like differentiation and multiplication
with functions from $C^{\infty}\left(\CalO\right)$) are permitted
on elements of $\DistributionSpace{\CalO}$, while this is not true
for the space $\SpaceReservoir{\CalO}$.

The principal reasons for using the space $\DistributionSpace{\CalO}$
instead of the space $\Schwartz'\left(\R^{\dimension}\right)$ for
the definition of the Fourier-side decomposition space $\FourierDecompSp{\CalQ}pY$
are the following:

\begin{enumerate}[leftmargin=0.6cm]
\item We want to allow the case $\CalO\subsetneq\R^{\dimension}$. If we
were to use the space $\Schwartz'\left(\R^{\dimension}\right)$, the
decomposition space (quasi)-norm would \emph{not} be positive definite,
or we would have to factor out a certain subspace of $\Schwartz'\left(\R^{\dimension}\right)$.
This is for example done in the usual definition of \emph{homogeneous}
Besov spaces, which are subspaces of $\Schwartz'\left(\R^{\dimension}\right)/\mathcal{P}$,
where $\mathcal{P}$ is the space of polynomials. Here, it seems more
natural to use the space $\DistributionSpace{\CalO}$.
\item In case of $\CalO=\R^{\dimension}$, one could use $\Schwartz'\left(\R^{\dimension}\right)$
as the reservoir. For example, Borup and Nielsen (in \cite{BorupNielsenDecomposition})
define their decomposition spaces as
\[
\CalD_{\Schwartz'}\left(\CalQ,L^{p},Y\right):=\left\{ f\in\Schwartz'\left(\smash{\R^{\dimension}}\right)\with\left\Vert f\right\Vert _{\DecompSp{\CalQ}pY}<\infty\right\} .
\]
But as we will see below (see Example~\ref{exa:BorupNielsenDecompositionSpaceIncomplete}),
this does in general \emph{not} yield a complete space, even for the
uniform covering $\CalQ$ of $\R^{\dimension}$ and $Y=\ell_{u}^{1}$
with a $\CalQ$-moderate weight $u$.

Nevertheless, in Section~\ref{sec:DecompositionSpacesAsSpacesOfTemperedDistributions},
we will develop criteria which yield the continuous embedding $\FourierDecompSp{\CalQ}pY\hookrightarrow\Schwartz'\left(\R^{\dimension}\right)$,
and thus also $\DecompSp{\CalQ}pY\hookrightarrow\Schwartz'\left(\R^{\dimension}\right)$.
If these conditions are satisfied, it is easy to see (at least for
$\CalO=\R^{\dimension}$) that the space $\CalD_{\Schwartz'}\left(\CalQ,L^{p},Y\right)$
from above is complete, since it coincides (up to trivial identifications)
with $\DecompSp{\CalQ}pY$.

\end{enumerate}
We finally remark that the notations $\SpaceTestFunctions{\CalO}$
and $\SpaceReservoir{\CalO}$ are inspired by Triebel's book \cite{TriebelFourierAnalysisAndFunctionSpaces};
see in particular \cite[Sections 2.2.1-2.2.4]{TriebelFourierAnalysisAndFunctionSpaces}.
Triebel also defines spaces very similar to $\DecompSp{\CalQ}pY$,
but restricts to the case in which the covering $\CalQ$ consists
of (closed) rectangles with sides parallel to the coordinate axes.
This is due to the fact that he also considers spaces of Triebel-Lizorkin
type as opposed to the spaces of Besov type that we consider here.
\end{rem}

\subsection{Well-definedness of decomposition spaces}

\label{subsec:DecompositionWellDefined}Our goal in this subsection
is to show that the decomposition space $\FourierDecompSp{\CalQ}pY$
is independent of the chosen $L^{p}$-BAPU $\left(\varphi_{i}\right)_{i\in I}$
(for $\CalQ$-regular global components $Y$). For later use, we will
actually show a slightly stronger statement, namely that one can use
an arbitrary \textbf{$L^{p}$-bounded control system} (defined below)
to obtain an equivalent (quasi)-norm.

The proofs of these results depend crucially on certain facts about
Fourier multipliers. For the range $p\in\left[1,\infty\right]$, Young's
inequality will be sufficient. But in the quasi-Banach regime $p\in\left(0,1\right)$,
we have to resort to Theorem~\ref{thm:QuasiBanachConvolution}. To
make application of this theorem more convenient, we note the following
special case.
\begin{cor}
\label{cor:QuasiBanachConvolutionSemiStructured}Let $p_{0}\in\left(0,1\right]$
and let $\CalQ=\left(Q_{i}\right)_{i\in I}=\left(T_{i}Q_{i}'+b_{i}\right)_{i\in I}$
be an $L^{p_{0}}$-decomposition covering of the open set $\emptyset\neq\CalO\subset\R^{\dimension}$.

For each $n\in\N_{0}$, there is a constant $C=C\left(\CalQ,n,\dimension,p_{0}\right)>0$
with the following property: If $p\in\left[p_{0},1\right]$ and $i\in I$
and furthermore

\begin{enumerate}
\item $f\in L^{1}\left(\R^{\dimension}\right)$ with $\supp f\subset\overline{Q_{i}^{n\ast}}$
and $\Fourier^{-1}f\in L^{p}\left(\R^{\dimension}\right)$,\vspace{0.1cm}
\item and $g\in\DistributionSpace{\CalO}$ with $\supp g\subset\overline{Q_{i}^{n\ast}}$
and $\Fourier^{-1}g\in L^{p}\left(\R^{\dimension}\right)$,
\end{enumerate}
\noindent then $g\in C_{0}\left(\R^{\dimension}\right)$, and $\Fourier^{-1}\left(fg\right)\in L^{p}\left(\R^{\dimension}\right)$
with
\[
\left\Vert \smash{\Fourier^{-1}}\left(fg\right)\right\Vert _{L^{p}}\leq C\cdot\left|\det T_{i}\right|^{\frac{1}{p}-1}\cdot\left\Vert \smash{\Fourier^{-1}}f\right\Vert _{L^{p}}\cdot\left\Vert \smash{\Fourier^{-1}}g\right\Vert _{L^{p}}.\qedhere
\]

\end{cor}

\begin{proof}
Lemma~\ref{lem:PartitionCoveringNecessary} implies that $\overline{Q_{i}^{n\ast}}\subset\bigcup_{j\in i^{n\ast}}\overline{Q_{j}}\subset\CalO$
is compact, since $i^{n\ast}$ is finite. Hence, $g\in\DistributionSpace{\CalO}$
is a distribution with compact support in $\CalO$, which means that
$g$ extends to a tempered distribution $g\in\Schwartz'\left(\R^{\dimension}\right)$,
cf.\@ \cite[Theorem 6.24(d) and Example 7.12(a)]{RudinFA}. Hence,
we have $\Fourier^{-1}g\in\Schwartz'\left(\R^{\dimension}\right)\cap L^{p}\left(\R^{\dimension}\right)$
with $\supp\widehat{\Fourier^{-1}g}=\supp g\subset\overline{Q_{i}^{n\ast}}$.
Therefore, Corollary~\ref{cor:BandlimitedEmbedding} shows $\Fourier^{-1}g\in L^{1}\left(\R^{\dimension}\right)$,
so that we get $g=\Fourier\Fourier^{-1}g\in C_{0}\left(\R^{\dimension}\right)$
by the Riemann-Lebesgue lemma. Thus, $f\cdot g\in L^{1}\left(\R^{\dimension}\right)$.
Furthermore, Theorem~\ref{thm:QuasiBanachConvolution} shows
\[
\Fourier^{-1}\left(fg\right)=\Fourier^{-1}\left(f\cdot\widehat{\Fourier^{-1}g}\right)\in L^{p}\left(\smash{\R^{\dimension}}\right)
\]
with
\[
\left\Vert \Fourier^{-1}\left(fg\right)\right\Vert _{L^{p}}=\left\Vert \Fourier^{-1}\left(f\cdot\widehat{\Fourier^{-1}g}\right)\right\Vert _{L^{p}}\leq\left[\lambda\left(\overline{Q_{i}^{n\ast}}-\overline{Q_{i}^{n\ast}}\right)\right]^{\frac{1}{p}-1}\cdot\left\Vert \smash{\Fourier^{-1}}f\right\Vert _{L^{p}}\cdot\left\Vert \smash{\Fourier^{-1}}g\right\Vert _{L^{p}}.
\]
But Corollary~\ref{cor:SemiStructuredDifferenceSetsMeasureEstimate}
yields a constant $C=C\left(\CalQ,n,\dimension\right)>0$ with $\lambda\left(\,\overline{Q_{i}^{n\ast}}-\overline{Q_{i}^{n\ast}}\,\right)\leq C\cdot\left|\det T_{i}\right|$
for all $i\in I$. Without loss of generality, $C\geq1$. Recall $0<p_{0}\leq p\leq1$,
and hence $0\leq\frac{1}{p}-1\leq\frac{1}{p}\leq\frac{1}{p_{0}}$.
All in all, we conclude
\begin{align*}
\left\Vert \Fourier^{-1}\left(fg\right)\right\Vert _{L^{p}} & \leq\left[\lambda\left(\,\overline{Q_{i}^{n\ast}}-\overline{Q_{i}^{n\ast}}\,\right)\right]^{\frac{1}{p}-1}\cdot\left\Vert \smash{\Fourier^{-1}}f\right\Vert _{L^{p}}\cdot\left\Vert \smash{\Fourier^{-1}}g\right\Vert _{L^{p}}\\
 & \leq C^{\frac{1}{p}-1}\cdot\left|\det T_{i}\right|^{\frac{1}{p}-1}\cdot\left\Vert \smash{\Fourier^{-1}}f\right\Vert _{L^{p}}\cdot\left\Vert \smash{\Fourier^{-1}}g\right\Vert _{L^{p}}\\
\left({\scriptstyle \text{since }C\geq1}\right) & \leq C^{1/p_{0}}\cdot\left|\det T_{i}\right|^{\frac{1}{p}-1}\cdot\left\Vert \smash{\Fourier^{-1}}f\right\Vert _{L^{p}}\cdot\left\Vert \smash{\Fourier^{-1}}g\right\Vert _{L^{p}},
\end{align*}
which completes the proof.
\end{proof}
As noted above, instead of just proving independence of $\FourierDecompSp{\CalQ}pY$
of the chosen $L^{p}$-BAPU, we will establish a slightly stronger
claim which will use the notion of an \textbf{$L^{p}$-bounded (control)
system}. This concept is based on that of a \textbf{bounded control
system} as introduced in \cite[Definition 2.6]{DecompositionSpaces1}.
\begin{defn}
\label{def:BoundedControlSystem}Let $p\in\left(0,\infty\right]$,
let $\emptyset\neq\CalO\subset\R^{\dimension}$ be open, and let $\CalQ=\left(Q_{i}\right)_{i\in I}$
be an admissible covering of $\CalO$. For $p\in\left(0,1\right)$,
assume additionally that $\CalQ=\left(T_{i}Q_{i}'+b_{i}\right)_{i\in I}$
is semi-structured.

A family $\Gamma=\left(\gamma_{i}\right)_{i\in I}$ of functions $\gamma_{i}\in\TestFunctionSpace{\CalO}$
is called an \textbf{$L^{p}$-bounded system} for $\CalQ$ if the
following conditions are satisfied:

\begin{enumerate}
\item There is some $\ell=\ell_{\Gamma,\CalQ}\in\N_{0}$ with $\gamma_{i}\equiv0$
on $\R^{\dimension}\setminus Q_{i}^{\ell\ast}$ for all $i\in I$.
\item The following expression (then a constant) is finite:
\[
C_{\CalQ,\Gamma,p}:=\begin{cases}
\sup_{i\in I}\left\Vert \Fourier^{-1}\gamma_{i}\right\Vert _{L^{1}}, & \text{for }p\in\left[1,\infty\right],\\
\vphantom{{\displaystyle \sum^{n}}}\sup_{i\in I}\left|\det T_{i}\right|^{\frac{1}{p}-1}\left\Vert \Fourier^{-1}\gamma_{i}\right\Vert _{L^{p}}, & \text{for }p\in\left(0,1\right).
\end{cases}
\]
\end{enumerate}
If furthermore $\gamma_{i}\equiv1$ on $Q_{i}$ for all $i\in I$,
we say that $\Gamma$ is an \textbf{$L^{p}$-bounded }\textbf{\emph{control}}\textbf{
system} for $\CalQ$.
\end{defn}

\begin{rem*}
We remark once more that the terminology ``$L^{p}$-bounded'' does
\emph{not} refer to the fact that the $\Fourier L^{p}$-norms $\left\Vert \Fourier^{-1}\gamma_{i}\right\Vert _{L^{p}}$
are bounded, but to the fact that the $\left(\gamma_{i}\right)_{i\in I}$
form a uniformly bounded family of $L^{p}$-Fourier multipliers, at
least for $p\in\left[1,\infty\right]$. For $p\in\left(0,1\right)$,
this is only true restricted to the space of $L^{p}$-functions with
Fourier support near $Q_{i}$; see Corollary~\ref{cor:QuasiBanachConvolutionSemiStructured}.
\end{rem*}
\begin{rem}
\label{rem:ClusteredBAPUYIeldsBoundedControlSystem}If $\ell\in\N_{0}$
is fixed and for each $i\in I$, some set $M_{i}\subset i^{\ell\ast}$
is selected, then the family $\Gamma=\left(\gamma_{i}\right)_{i\in I}:=\left(\varphi_{M_{i}}\right)_{i\in I}$
is an $L^{p}$-bounded system for $\CalQ$ if $\left(\varphi_{i}\right)_{i\in I}$
is an $L^{p}$-BAPU for $\CalQ$. To see this, distinguish two cases:

\textbf{Case 1}: If $p\in\left[1,\infty\right]$, simply note (using
the triangle inequality for $L^{1}$ and the estimate $\left|i^{\ell\ast}\right|\leq N_{\CalQ}^{\ell}$
from Lemma~\ref{lem:SemiStructuredClusterInvariant}) that
\[
\left\Vert \Fourier^{-1}\gamma_{i}\right\Vert _{L^{1}}\leq\sum_{j\in M_{i}}\left\Vert \Fourier^{-1}\varphi_{j}\right\Vert _{L^{1}}\leq\left|M_{i}\right|\cdot C_{\CalQ,\Phi,p}\leq\left|i^{\ell\ast}\right|\cdot C_{\CalQ,\Phi,p}\leq N_{\CalQ}^{\ell}\cdot C_{\CalQ,\Phi,p}.
\]

\textbf{Case 2}: For $p\in\left(0,1\right)$, we first note that the
weight $\left(\left|\det T_{i}\right|\right)_{i\in I}$ is $\CalQ$-moderate.
Indeed, Hadamard's inequality (see for instance  \cite[Section 75]{RieszFunctionalAnalysis})
implies $\left|\det A\right|\leq\left\Vert A\right\Vert ^{\dimension}$
for all $A\in\R^{\dimension\times\dimension}$; therefore,
\begin{equation}
\sup_{i\in I}\,\sup_{j\in i^{\ast}}\frac{\left|\det T_{j}\right|}{\left|\det T_{i}\right|}=\sup_{i\in I}\,\sup_{j\in i^{\ast}}\left|\det\left(T_{i}^{-1}T_{j}\right)\right|\leq\sup_{i\in I}\,\sup_{j\in i^{\ast}}\left\Vert T_{i}^{-1}T_{j}\right\Vert ^{\dimension}\leq C_{\CalQ}^{\dimension}.\label{eq:DeterminantIsModerate}
\end{equation}
Now, using the estimate $\left\Vert \sum_{j=1}^{n}f_{j}\right\Vert _{L^{p}}\leq n^{\frac{1}{p}-1}\cdot\sum_{j=1}^{n}\left\Vert f_{j}\right\Vert _{L^{p}}$
(cf.\@ \cite[Exercise 1.1.5(c)]{GrafakosClassical}) and the bound
$\left|M_{i}\right|\leq\left|i^{\ell\ast}\right|\leq N_{\CalQ}^{\ell}$
from Lemma~\ref{lem:SemiStructuredClusterInvariant}, we conclude
\begin{align*}
\left|\det T_{i}\right|^{\frac{1}{p}-1}\cdot\left\Vert \Fourier^{-1}\gamma_{i}\right\Vert _{L^{p}} & \leq N_{\CalQ}^{\ell\left(\frac{1}{p}-1\right)}\cdot\sum_{j\in i^{\ell\ast}}\left|\det T_{i}\right|^{\frac{1}{p}-1}\cdot\left\Vert \Fourier^{-1}\varphi_{j}\right\Vert _{L^{p}}\\
 & \leq N_{\CalQ}^{\ell\left(\frac{1}{p}-1\right)}\cdot\sum_{j\in i^{\ell\ast}}\left[C_{\CalQ}^{d\left(\frac{1}{p}-1\right)}\cdot\left|\det T_{j}\right|^{\frac{1}{p}-1}\cdot\left\Vert \Fourier^{-1}\varphi_{j}\right\Vert _{L^{p}}\right]\\
 & \leq N_{\CalQ}^{\ell\left(\frac{1}{p}-1\right)}C_{\CalQ}^{\dimension\left(\frac{1}{p}-1\right)}\cdot\left|i^{\ell\ast}\right|\cdot C_{\CalQ,\Phi,p}\\
 & \leq N_{\CalQ}^{\ell/p}\cdot C_{\CalQ}^{\dimension\left(\frac{1}{p}-1\right)}\cdot C_{\CalQ,\Phi,p}\qquad\forall\,i\in I.
\end{align*}

All in all, this yields (for both cases) an estimate of the form $C_{\CalQ,\Gamma,p}\leq C\left(\CalQ,C_{\CalQ,\Phi,p},\dimension,p,\ell\right)$.

Finally, if $M_{i}\supset i^{\ast}$ holds for all $i\in I$, then
Lemma~\ref{lem:PartitionCoveringNecessary} shows $\varphi_{M_{i}}\equiv1$
on $Q_{i}$ for all $i\in I$, so that $\Gamma$ is an $L^{p}$-bounded
\emph{control} system for $\CalQ$.
\end{rem}

With these notions, we can now state the following result which allows
to use an $L^{p}$-bounded (control) system instead of an $L^{p}$-BAPU
to determine the quasi-norm on $\FourierDecompSp{\CalQ}pY$ (up to
constant factors). We remark that the statement of the theorem is
essentially the same as \cite[Corollary 2.5]{DecompositionSpaces1},
but the case $p\in\left(0,1\right)$ is of course not covered by the
setting considered in \cite{DecompositionSpaces1}. Roughly, the theorem
shows the following:
\begin{enumerate}[leftmargin=0.6cm]
\item Instead of a $\CalQ$-BAPU $\Phi=\left(\varphi_{i}\right)_{i\in I}$,
one can use any $L^{p}$-bounded control system $\Gamma=\left(\gamma_{i}\right)_{i\in I}$
to calculate the (quasi)-norm $\left\Vert \mybullet\right\Vert _{\FourierDecompSp{\CalQ}pY}$.
Essentially, this means that instead of requiring $\supp\varphi_{i}\subset Q_{i}$,
it is enough to require $\supp\gamma_{i}\subset Q_{i}^{k\ast}$ for
a fixed $k$. Likewise, instead of $\sum_{i\in I}\varphi_{i}\equiv1$
on $\CalQ$, it suffices if $\gamma_{i}\equiv1$ on $Q_{i}$ for each
$i\in I$.
\item Given a (not necessarily disjoint) ``partition'' $I=\bigcup_{r=1}^{r_{0}}I^{\left(r\right)}$
of the index set $I$, we can calculate the (quasi)-norm $\left\Vert f\right\Vert _{\FourierDecompSp{\CalQ}pY}$
``localized'' to each of the sets $Q^{\left(r\right)}=\bigcup_{i\in I^{\left(r\right)}}Q_{i}$
and then aggregate the individual contributions. This will become
useful in connection with the disjointization lemma (Lemma~\ref{lem:DisjointizationPrinciple}).
\end{enumerate}
Regarding the notation in the following theorem, if $Y\subset\Compl^{I}$
is a solid sequence space on $I$ and if $J\subset I$ is an arbitrary
subset, we define the \textbf{restricted sequence space} $Y|_{J}$
as
\begin{equation}
Y|_{J}:=\left\{ c=\left(c_{\ell}\right)_{\ell\in J}\in\Compl^{J}\with\tilde{c}\in Y\right\} ,\label{eq:DefinitionRestrictedSequenceSpace}
\end{equation}
where $\tilde{c}=\left(c_{i}\right)_{i\in I}$ is the \textbf{trivial
extension} of $c=\left(c_{\ell}\right)_{\ell\in J}$ onto $I$, i.e.\@
$c_{i}:=0$ for $i\in I\setminus J$. As expected, we let $\left\Vert c\right\Vert _{Y|_{J}}:=\left\Vert \tilde{c}\right\Vert _{Y}$.
With this definition, $Y|_{J}\subset\Compl^{J}$ is a solid sequence
space on $J$. Furthermore, if $Y$ is complete, then so is $Y|_{J}$,
since $c\mapsto\widetilde{c}$ defines an isometric isomorphism between
$Y|_{J}$ and a closed subspace of $Y$. Closedness of this subspace
is a consequence of Lemma~\ref{lem:SolidSequenceSpaceEmbedsIntoWeightedLInfty}.
\begin{thm}
\label{thm:BoundedControlSystemEquivalentQuasiNorm}Let $\emptyset\neq\CalO\subset\R^{\dimension}$
be open, let $\CalQ=\left(Q_{i}\right)_{i\in I}$ be an $L^{p}$-decomposition
covering of $\CalO$ for some $p\in\left(0,\infty\right]$, and let
$Y\subset\Compl^{I}$ be $\CalQ$-regular. Furthermore, let $r_{0}\in\N$
and assume $I=\bigcup_{r=1}^{r_{0}}I^{\left(r\right)}$ for certain
subsets $I^{\left(r\right)}\subset I$. Finally, let $\Gamma=\left(\gamma_{i}\right)_{i\in I}$
be an $L^{p}$-bounded system for $\CalQ$.

For $f\in\DistributionSpace{\CalO}$, define
\[
\left\Vert f\right\Vert _{\Gamma,\left(I^{\left(r\right)}\right)_{r},L^{p},Y}:=\sum_{r=1}^{r_{0}}\left\Vert \left(\left\Vert \Fourier^{-1}\left(\gamma_{i}f\right)\right\Vert _{L^{p}}\right)_{i\in I^{\left(r\right)}}\right\Vert _{Y|_{I^{\left(r\right)}}}\in\left[0,\infty\right].
\]
If $\Gamma_{\CalQ}:Y\to Y,c\mapsto c^{\ast}$ denotes the clustering
map, then there is a positive constant 
\[
C=\begin{cases}
C\left(\CalQ,p,r_{0},\dimension,\ell_{\Gamma,\CalQ},\vertiii{\Gamma_{\CalQ}}_{Y\to Y}\right), & \text{if }p\in\left(0,1\right),\\
C\left(N_{\CalQ},p,r_{0},\ell_{\Gamma,\CalQ},\vertiii{\Gamma_{\CalQ}}_{Y\to Y}\right), & \text{if }p\in\left[1,\infty\right]
\end{cases}
\]
with 
\[
\left\Vert f\right\Vert _{\Gamma,\left(I^{\left(r\right)}\right)_{r},L^{p},Y}\leq C\cdot C_{\CalQ,\Gamma,p}\cdot\left\Vert f\right\Vert _{\CalD_{\Fourier,\Phi}\left(\CalQ,L^{p},Y\right)}\qquad\forall\,\,f\in\DistributionSpace{\CalO},
\]
for each $L^{p}$-BAPU $\Phi=\left(\varphi_{i}\right)_{i\in I}$ for
$\CalQ$.

\medskip{}

Conversely, if $\Gamma$ is an $L^{p}$-bounded \emph{control} system
for $\CalQ$ and if $C_{0}\geq1$ denotes a triangle constant for
$Y$, then there is a positive constant 
\[
C'=\begin{cases}
C'\left(\CalQ,p,\dimension,C_{0},r_{0},\ell_{\Gamma,\CalQ}\right), & \text{if }p\in\left(0,1\right),\\
C'\left(C_{0},r_{0}\right), & \text{if }p\in\left[1,\infty\right]
\end{cases}
\]
with
\[
\left\Vert f\right\Vert _{\CalD_{\Fourier,\Phi}\left(\CalQ,L^{p},Y\right)}\leq C'\cdot C_{\CalQ,\Phi,p}\cdot\left\Vert f\right\Vert _{\Gamma,\left(I^{\left(r\right)}\right)_{r},L^{p},Y}\text{ for all }f\in\DistributionSpace{\CalO}
\]
for each $L^{p}$-BAPU $\Phi=\left(\varphi_{i}\right)_{i\in I}$ for
$\CalQ$.
\end{thm}

\begin{proof}
Let $\ell:=\ell_{\Gamma,\CalQ}$. Let $\Phi=\left(\varphi_{i}\right)_{i\in I}$
be an $L^{p}$-BAPU for $\CalQ$. Using Lemma~\ref{lem:PartitionCoveringNecessary},
we get $\varphi_{i}^{\left(\ell+1\right)\ast}\equiv1$ on $Q_{i}^{\ell\ast}$
for all $i\in I$. Because of $\gamma_{i}\equiv0$ on $\R^{\dimension}\setminus Q_{i}^{\ell\ast}$,
we conclude
\[
\gamma_{i}=\varphi_{i}^{\left(\ell+1\right)\ast}\gamma_{i}=\sum_{j\in i^{\left(\ell+1\right)\ast}}\gamma_{i}\,\varphi_{j}\,.
\]

We begin by proving the first estimate. To this end, let $f\in\DistributionSpace{\CalO}$
with $\left\Vert f\right\Vert _{\CalD_{\Fourier,\Phi}\left(\CalQ,L^{p},Y\right)}<\infty$
(otherwise, the estimate is trivial). Since $L^{p}\left(\R^{\dimension}\right)$
is a quasi-normed vector space and since Lemma~\ref{lem:SemiStructuredClusterInvariant}
yields the uniform bound $\left|\smash{i^{\left(\ell+1\right)\ast}}\right|\leq N_{\CalQ}^{\ell+1}$,
there is a constant $C_{1}=C_{1}\left(N_{\CalQ},\ell,p\right)>0$
with
\[
\left\Vert \Fourier^{-1}\left(\gamma_{i}\,f\right)\right\Vert _{L^{p}}=\vphantom{\sum_{j\in i^{\left(\ell+1\right)\ast}}}\left\Vert \,\smash{\sum_{j\in i^{\left(\ell+1\right)\ast}}}\vphantom{\sum}\Fourier^{-1}\left(\gamma_{i}\,\varphi_{j}\,f\right)\,\right\Vert _{L^{p}}\leq C_{1}\sum_{j\in i^{\left(\ell+1\right)\ast}}\left\Vert \Fourier^{-1}\left(\gamma_{i}\,\varphi_{j}\,f\right)\right\Vert _{L^{p}}\qquad\forall\,i\in I.
\]
There are now two cases. For $p\in\left[1,\infty\right]$, Young's
inequality $L^{1}\ast L^{p}\hookrightarrow L^{p}$ yields
\[
\left\Vert \Fourier^{-1}\left(\gamma_{i}\,\varphi_{j}\,f\right)\right\Vert _{L^{p}}\leq\left\Vert \Fourier^{-1}\gamma_{i}\right\Vert _{L^{1}}\cdot\left\Vert \Fourier^{-1}\left(\varphi_{j}\,f\right)\right\Vert _{L^{p}}\leq C_{\CalQ,\Gamma,p}\cdot\left\Vert \Fourier^{-1}\left(\varphi_{j}\,f\right)\right\Vert _{L^{p}}.
\]
For $p\in\left(0,1\right)$, we invoke Corollary~\ref{cor:QuasiBanachConvolutionSemiStructured},
which yields a constant $C_{2}=C_{2}\left(\dimension,\ell,p,\CalQ\right)>0$
with
\[
\left\Vert \Fourier^{-1}\left(\gamma_{i}\,\varphi_{j}\,f\right)\right\Vert _{L^{p}}\leq C_{2}\cdot\left|\det T_{i}\right|^{\frac{1}{p}-1}\cdot\left\Vert \Fourier^{-1}\gamma_{i}\right\Vert _{L^{p}}\cdot\left\Vert \Fourier^{-1}\left(\varphi_{j}\,f\right)\right\Vert _{L^{p}}\leq C_{2}C_{\CalQ,\Gamma,p}\cdot\left\Vert \Fourier^{-1}\left(\varphi_{j}\,f\right)\right\Vert _{L^{p}},
\]
where $\CalQ=\left(T_{i}Q_{i}'+b_{i}\right)_{i\in I}$ is semi-structured
by definition of $L^{p}$-decomposition coverings for $p\in\left(0,1\right)$,
see Definition~\ref{defn:QuasiBanachBAPUAndLpDecompositionCovering}.
In the estimate above, we also used the inclusions $\supp\gamma_{i}\subset\overline{Q_{i}^{\ell\ast}}\subset\overline{Q_{i}^{\left(\ell+1\right)\ast}}$
and $\supp\varphi_{j}\subset\overline{Q_{j}}\subset\overline{Q_{i}^{\left(\ell+1\right)\ast}}$
for all $j\in i^{\left(\ell+1\right)\ast}$. Thus, if we set $C_{2}:=1$
for $p\in\left[1,\infty\right]$, the estimate above is valid for
all $i\in I$, $j\in i^{\left(\ell+1\right)\ast}$, and $p\in\left(0,\infty\right]$.

With the ``higher order clustering map'' $\Theta_{\ell+1}$ from
Lemma~\ref{lem:HigherOrderClusteringMap}, and with the solidity
of $Y$, we get
\begin{align*}
\left\Vert \left(\left\Vert \Fourier^{-1}\left(\gamma_{j}f\right)\right\Vert _{L^{p}}\right)_{i\in I^{\left(r\right)}}\right\Vert _{Y|_{I^{\left(r\right)}}} & \leq C_{1}C_{2}C_{\CalQ,\Gamma,p}\cdot\left\Vert \left(\,\smash{\sum_{j\in i^{\left(\ell+1\right)\ast}}}\vphantom{\sum}\left\Vert \Fourier^{-1}\left(\varphi_{j}f\right)\right\Vert _{L^{p}}\,\right)_{i\in I}\right\Vert _{Y}\vphantom{\sum_{j\in i^{\left(\ell+1\right)\ast}}}\\
 & =C_{1}C_{2}C_{\CalQ,\Gamma,p}\cdot\left\Vert \Theta_{\ell+1}\left[\left(\left\Vert \Fourier^{-1}\left(\varphi_{i}f\right)\right\Vert _{L^{p}}\right)_{i\in I}\right]\right\Vert _{Y}\\
 & \leq C_{1}C_{2}C_{\CalQ,\Gamma,p}\cdot\vertiii{\Theta_{\ell+1}}_{Y\to Y}\cdot\left\Vert \left(\left\Vert \Fourier^{-1}\left(\varphi_{i}f\right)\right\Vert _{L^{p}}\right)_{i\in I}\right\Vert _{Y}\\
 & =C_{1}C_{2}C_{\CalQ,\Gamma,p}\cdot\vertiii{\Theta_{\ell+1}}_{Y\to Y}\cdot\left\Vert f\right\Vert _{\CalD_{\Fourier,\Phi}\left(\CalQ,L^{p},Y\right)}.
\end{align*}
Summing over $r\in\underline{r_{0}}$ completes the proof of the first
estimate, since $\vertiii{\Theta_{\ell+1}}_{Y\to Y}\leq\vertiii{\Gamma_{\CalQ}}_{Y\to Y}^{\ell+1}$
by Lemma~\ref{lem:HigherOrderClusteringMap}.

\medskip{}

Now, let $f\in\DistributionSpace{\CalO}$ with $\left\Vert f\right\Vert _{\Gamma,\left(I^{\left(r\right)}\right)_{r},L^{p},Y}<\infty$
and assume that $\Gamma$ is an $L^{p}$-bounded control system for
$\CalQ$, i.e.\@ that $\gamma_{i}\equiv1$ on $Q_{i}$ for all $i\in I$.
This implies $\varphi_{i}=\varphi_{i}\,\gamma_{i}$. Thus, in case
of $p\in\left[1,\infty\right]$, Young's inequality yields
\[
\left\Vert \Fourier^{-1}\left(\varphi_{i}\,f\right)\right\Vert _{L^{p}}=\left\Vert \Fourier^{-1}\left(\varphi_{i}\,\gamma_{i}\,f\right)\right\Vert _{L^{p}}\leq\left\Vert \Fourier^{-1}\varphi_{i}\right\Vert _{L^{1}}\cdot\left\Vert \Fourier^{-1}\left(\gamma_{i}\,f\right)\right\Vert _{L^{p}}\leq C_{2}C_{\CalQ,\Phi,p}\cdot\left\Vert \Fourier^{-1}\left(\gamma_{i}\,f\right)\right\Vert _{L^{p}}
\]
with $C_{2}:=1$. In case of $p\in\left(0,1\right)$, we can use the
same constant $C_{2}=C_{2}\left(\dimension,\ell,p,\CalQ\right)>0$
provided by Corollary~\ref{cor:QuasiBanachConvolutionSemiStructured}
as above to conclude
\begin{align*}
\left\Vert \Fourier^{-1}\left(\varphi_{i}\,f\right)\right\Vert _{L^{p}}=\left\Vert \Fourier^{-1}\left(\varphi_{i}\,\gamma_{i}\,f\right)\right\Vert _{L^{p}} & \leq C_{2}\cdot\left|\det T_{i}\right|^{\frac{1}{p}-1}\cdot\left\Vert \Fourier^{-1}\varphi_{i}\right\Vert _{L^{p}}\cdot\left\Vert \Fourier^{-1}\left(\gamma_{i}\,f\right)\right\Vert _{L^{p}}\\
 & \leq C_{2}C_{\CalQ,\Phi,p}\cdot\left\Vert \Fourier^{-1}\left(\gamma_{i}\,f\right)\right\Vert _{L^{p}},
\end{align*}
so that this estimate holds for all $i\in I$ and $p\in\left(0,\infty\right]$.

Let $c_{i}:=\left\Vert \Fourier^{-1}\left(\gamma_{i}f\right)\right\Vert _{L^{p}}$
for $i\in I$. Since $I=\bigcup_{r=1}^{r_{0}}I^{\left(r\right)}$
and because of $c_{i}\geq0$ for all $i\in I$, we have
\[
0\leq c_{i}\leq\sum_{r=1}^{r_{0}}\left(c\cdot\Indicator_{I^{\left(r\right)}}\right)_{i}\text{ for all }i\in I.
\]
Since $Y$ is a solid quasi-normed sequence space, there is thus a
constant $C_{3}=C_{3}\left(C_{0},r_{0}\right)>0$ (recall that $C_{0}$
is a triangle constant for $Y$) with
\begin{align*}
\left\Vert f\right\Vert _{\CalD_{\Fourier,\Phi}\left(\CalQ,L^{p},Y\right)}=\left\Vert \left(\left\Vert \Fourier^{-1}\left(\varphi_{i}f\right)\right\Vert _{L^{p}}\right)_{i\in I}\right\Vert _{Y} & \leq C_{2}C_{\CalQ,\Phi,p}\cdot\left\Vert c\right\Vert _{Y}\\
 & \leq C_{2}C_{\CalQ,\Phi,p}\cdot\left\Vert \sum_{r=1}^{r_{0}}c\cdot\Indicator_{I^{\left(r\right)}}\right\Vert _{Y}\\
 & \leq C_{2}C_{3}C_{\CalQ,\Phi,p}\cdot\sum_{r=1}^{r_{0}}\left\Vert c\cdot\Indicator_{I^{\left(r\right)}}\right\Vert _{Y}\\
 & =C_{2}C_{3}C_{\CalQ,\Phi,p}\cdot\left\Vert f\right\Vert _{\Gamma,\left(I^{\left(r\right)}\right)_{r},L^{p},Y}<\infty.\qquad\qquad\qedhere
\end{align*}
\end{proof}
Using the theorem above, well-definedness of $\FourierDecompSp{\CalQ}pY$,
i.e.\@ independence of the chosen $L^{p}$-BAPU is an easy consequence:
\begin{cor}
\label{cor:DecompositionSpaceWellDefined}Let $\emptyset\neq\CalO\subset\R^{\dimension}$
be open and let $p\in\left(0,\infty\right]$. Assume that $\CalQ=\left(Q_{i}\right)_{i\in I}$
is an $L^{p}$-decomposition covering of $\CalO$ with $L^{p}$-BAPUs
$\left(\varphi_{i}\right)_{i\in I},\left(\psi_{i}\right)_{i\in I}$
and that $Y$ is $\CalQ$-regular.

Then we have
\[
\left\Vert \left(\left\Vert \Fourier^{-1}\left(\varphi_{i}f\right)\right\Vert _{L^{p}}\right)_{i\in I}\right\Vert _{Y}\asymp\left\Vert \left(\left\Vert \Fourier^{-1}\left(\psi_{i}f\right)\right\Vert _{L^{p}}\right)_{i\in I}\right\Vert _{Y}
\]
for all $f\in\DistributionSpace{\CalO}$, where the implied constants
are independent of $f$. Especially, the left-hand side is finite
iff the right-hand side is.
\end{cor}

\begin{proof}
By Remark~\ref{rem:ClusteredBAPUYIeldsBoundedControlSystem}, we
know that $\Gamma:=\left(\varphi_{i}^{\ast}\right)_{i\in I}$ yields
an $L^{p}$-bounded control system for $\CalQ$. By Theorem~\ref{thm:BoundedControlSystemEquivalentQuasiNorm},
this implies
\[
\left\Vert \mybullet\right\Vert _{\CalD_{\Fourier,\Phi}\left(\CalQ,L^{p},Y\right)}\asymp\left\Vert \mybullet\right\Vert _{\Gamma,\left(I\right),L^{p},Y}\asymp\left\Vert \mybullet\right\Vert _{\CalD_{\Fourier,\Psi}\left(\CalQ,L^{p},Y\right)}
\]
for $\Phi=\left(\varphi_{i}\right)_{i\in I}$ and $\Psi=\left(\psi_{i}\right)_{i\in I}$,
where the implied constants only depend on $p,\dimension,\CalQ,\Phi,\Psi$
and $Y$.
\end{proof}
Note that for $\FourierDecompSp{\CalQ}pY$ to be defined at all, we
always needed to assume that $\CalQ$ is an $L^{p}$-decomposition
covering, i.e.\@ that there is \emph{some} $L^{p}$-BAPU for $\CalQ$.
It is thus important to have sufficient criteria for this to hold.

As we noted after the definition of almost structured coverings, the
most important reason for their introduction is that they always posses
$L^{p}$-BAPUs. This was observed in \cite[Proposition 1]{BorupNielsenDecomposition}
for the case of \emph{structured} coverings of $\CalO=\R^{\dimension}$.
As a slight generalization, the following result was shown in \cite[Theorem 2.8]{DecompositionIntoSobolev}:
\begin{thm}
\label{thm:AlmostStructuredAdmissibleAdmitsBAPU}Let $\CalQ=\left(Q_{i}\right)_{i\in I}=\left(T_{i}Q_{i}'+b_{i}\right)_{i\in I}$
be an \emph{almost structured} covering of the open set $\emptyset\neq\CalO\subset\R^{\dimension}$.

Then $I$ is countably infinite and there is a family $\Phi=\left(\varphi_{i}\right)_{i\in I}$
in $\TestFunctionSpace{\CalO}$ such that $\Phi$ is an $L^{p}$-BAPU
for $\CalQ$ simultaneously for each $p\in\left(0,\infty\right]$.
In particular, every almost structured covering is an $L^{p}$-decomposition
covering, for arbitrary $p\in\left(0,\infty\right]$.
\end{thm}

Now, we know that the Fourier-side decomposition spaces $\FourierDecompSp{\CalQ}pY$—and
thus also the space-side spaces $\DecompSp{\CalQ}pY$—are well-defined
for arbitrary $L^{p}$-decomposition coverings $\CalQ$, in particular
for almost structured $\CalQ$, if $Y$ is $\CalQ$-regular. In the
next subsection, we close our investigation of the basic properties
of decomposition spaces by showing completeness of $\FourierDecompSp{\CalQ}pY$.

\subsection{Completeness of decomposition spaces}

\label{subsec:DecompositionCompleteness}We will now show that the
space $\FourierDecompSp{\CalQ}pY$ is indeed a quasi-Banach space.
For the proof, we first establish an equivalent condition for completeness
of quasi-normed vector spaces. For a normed vector space $X$, it
is well-known (cf.\@ \cite[Theorem 5.1]{FollandRA}) that $X$ is
complete if and only if ``absolute convergence'' of a series in
$X$ implies convergence of the series. For quasi-normed vector spaces,
we have the following replacement:
\begin{lem}
\label{lem:AbsoluteConvergenceForQuasiBanachSpaces}Let $\left(X,\left\Vert \mybullet\right\Vert \right)$
be a quasi-normed vector space and let $C\geq1$ be a triangle constant
for $\left\Vert \mybullet\right\Vert $. Then the following are true:

\begin{enumerate}[leftmargin=0.7cm]
\item If $x_{n}\to x$, then $\left\Vert x\right\Vert \leq C\cdot\liminf_{n\to\infty}\left\Vert x_{n}\right\Vert $.
\item We have $\left\Vert \sum_{i=1}^{n}x_{i}\right\Vert \leq\sum_{i=1}^{n}C^{i}\left\Vert x_{i}\right\Vert $
for all $n\in\N$ and $x_{1},\dots,x_{n}\in X$.
\item $X$ is complete if there is some $M>1$ with the following property:
\[
\quad\text{The series }\sum_{n=1}^{\infty}x_{n}\text{ converges (in \ensuremath{X}) for each }\left(x_{n}\right)_{n\in\N}\in X^{\N}\text{ with }\left\Vert x_{n}\right\Vert \leq M^{-n}\text{ for all }n\in\N.
\]
\item Conversely, if $X$ is complete and if $\left(x_{n}\right)_{n\in\N}\in X^{\N}$
satisfies $\sum_{n=1}^{\infty}C^{n}\left\Vert x_{n}\right\Vert <\infty$,
then $\sum_{n=1}^{\infty}x_{n}$ converges (in $X$), and we have
\[
\left\Vert \,\smash{\sum_{n=1}^{\infty}}\,\vphantom{\sum}x_{n}\,\right\Vert \leq C\cdot\sum_{n=1}^{\infty}C^{n}\left\Vert x_{n}\right\Vert .\qedhere
\]
\end{enumerate}
\end{lem}

\begin{rem*}
Quasi-norms are in general \emph{not} continuous, i.e.\@ $\left\Vert x_{n}-x\right\Vert \xrightarrow[n\to\infty]{}0$
does in general \emph{not} imply $\left\Vert x_{n}\right\Vert \xrightarrow[n\to\infty]{}\left\Vert x\right\Vert $.
Thus, the following proof has to avoid using this property.
\end{rem*}
\begin{proof}
Ad (1): We have $\left\Vert x\right\Vert \leq C\cdot\left[\left\Vert x-x_{n}\right\Vert +\left\Vert x_{n}\right\Vert \right]$.
Taking the $\liminf_{n\to\infty}$ proves the claim.

\medskip{}

Ad (2): We prove the claim by induction on $n\in\N$. For $n=1$,
the claim is trivial, since we have $C\geq1$. For the induction step,
note
\begin{align*}
\vphantom{\sum_{i=1}^{n+1}}\left\Vert \,\smash{\sum_{i=1}^{n+1}}\,\vphantom{\sum}x_{i}\,\right\Vert  & =\left\Vert \,x_{1}+\smash{\sum_{i=2}^{n+1}}\,\vphantom{\sum}x_{i}\,\right\Vert \leq C\left[\left\Vert x_{1}\right\Vert +\vphantom{\sum_{i=2}^{n+1}}\left\Vert \,\smash{\sum_{i=2}^{n+1}}\,\vphantom{\sum}x_{i}\,\right\Vert \right]\\
 & \overset{\left(\ast\right)}{\leq}C\left[\,\left\Vert x_{1}\right\Vert +\smash{\sum_{i=2}^{n+1}}\,\vphantom{\sum}C^{i-1}\left\Vert x_{i}\right\Vert \,\right]\vphantom{\sum_{i=2}^{n+1}}=\sum_{i=1}^{n+1}C^{i}\left\Vert x_{i}\right\Vert .
\end{align*}
At $\left(\ast\right)$, we implicitly used that the induction hypothesis
yields 
\begin{equation}
\vphantom{\sum_{i=2}^{n+1}}\left\Vert \,\smash{\sum_{i=2}^{n+1}}\,\vphantom{\sum}x_{i}\,\right\Vert =\left\Vert \,\smash{\sum_{\ell=1}^{n}}\,\vphantom{\sum}x_{\ell+1}\,\right\Vert \leq\sum_{\ell=1}^{n}C^{\ell}\left\Vert x_{\ell+1}\right\Vert =\sum_{i=2}^{n+1}C^{i-1}\left\Vert x_{i}\right\Vert .\label{eq:IteratedQuasiTriangleInequalityImplicitArgument}
\end{equation}
In the remainder of the proof, we will use similar manipulations without
further comment.

\medskip{}

Ad (3): Assume that $X$ satisfies the stated property and let $\left(y_{n}\right)_{n\in\N}$
be a Cauchy sequence in $X$. Using a trivial induction, we can choose
a strictly increasing sequence $\left(N_{n}\right)_{n\in\N}$ satisfying
\[
\forall\,m,\ell\geq N_{n}:\qquad\left\Vert y_{m}-y_{\ell}\right\Vert \leq M^{-n}.
\]
Define $x_{n}:=y_{N_{n+1}}-y_{N_{n}}$. Because of $N_{n+1}>N_{n}$,
this yields $\left\Vert x_{n}\right\Vert \leq M^{-n}$ for all $n\in\N$,
so that the assumption shows that
\[
x:=\sum_{n=1}^{\infty}x_{n}=\lim_{K\to\infty}\sum_{n=1}^{K}x_{n}=\lim_{K\to\infty}\sum_{n=1}^{K}\left(y_{N_{n+1}}-y_{N_{n}}\right)=\lim_{K\to\infty}\left(y_{N_{K+1}}-y_{N_{1}}\right)
\]
exists in $X$. We conclude $y_{N_{K+1}}\to z:=x+y_{N_{1}}$ as $K\to\infty$.
But this yields
\[
\left\Vert y_{n}-z\right\Vert \leq C\cdot\left[\left\Vert y_{n}-y_{N_{n+1}}\right\Vert +\left\Vert y_{N_{n+1}}-z\right\Vert \right]\xrightarrow[n\to\infty]{}0,
\]
so that $X$ is complete.

\medskip{}

Ad (4): Set $y_{n}:=\sum_{i=1}^{n}x_{i}$ for $n\in\N$. Using the
second part of the lemma, we get for $N\geq M\geq M_{0}$ that
\begin{align*}
\left\Vert y_{N}-y_{M}\right\Vert =\left\Vert \,\smash{\sum_{i=M+1}^{N}}\,\vphantom{\sum}x_{i}\,\right\Vert \vphantom{\sum_{i=M+1}^{N}} & \leq\sum_{i=M+1}^{N}C^{i-M}\left\Vert x_{i}\right\Vert \\
 & \leq C^{-M}\cdot\sum_{i=M+1}^{\infty}C^{i}\left\Vert x_{i}\right\Vert \leq\sum_{i=M_{0}+1}^{\infty}C^{i}\left\Vert x_{i}\right\Vert \xrightarrow[M_{0}\to\infty]{}0,
\end{align*}
so that $\left(y_{n}\right)_{n\in\N}$ is Cauchy and thus convergent
to some $y\in X$ by completeness. But this means $y=\sum_{n=1}^{\infty}x_{n}$.

Finally, the first and second part of the lemma yield
\[
\vphantom{\sum_{n=1}^{\infty}}\left\Vert \,\smash{\sum_{n=1}^{\infty}}\vphantom{\sum}\,x_{n}\,\right\Vert \leq C\cdot\liminf_{N\to\infty}\left\Vert \,\smash{\sum_{n=1}^{N}}\,\vphantom{\sum}x_{n}\,\right\Vert \leq C\cdot\liminf_{N\to\infty}\sum_{n=1}^{N}C^{n}\left\Vert x_{n}\right\Vert =C\cdot\sum_{n=1}^{\infty}C^{n}\left\Vert x_{n}\right\Vert .\qedhere
\]
\end{proof}
Now, we can show that $\FourierDecompSp{\CalQ}pY$ is a quasi-Banach
space which embeds continuously into $\DistributionSpace{\CalO}$.
\begin{thm}
\label{thm:DecompositionSpaceComplete}Let $\emptyset\neq\CalO\subset\R^{\dimension}$
be open and let $p\in\left(0,\infty\right]$. Let $\CalQ=\left(Q_{i}\right)_{i\in I}$
be an $L^{p}$-decomposition covering of $\CalO$ and let $Y\subset\Compl^{I}$
be $\CalQ$-regular.

Then $\FourierDecompSp{\CalQ}pY$ is a quasi-Banach space which embeds
continuously into $\DistributionSpace{\CalO}$. The triangle constant
$C$ for $\FourierDecompSp{\CalQ}pY$ satisfies $C\leq C^{\left(p\right)}C_{Y}$,
where $C_{Y}$ is a triangle constant for $Y$ and $C^{\left(p\right)}$
is a triangle constant for $L^{p}$.

If $Y$ is a Banach space (instead of a quasi-Banach space) and if
$p\in\left[1,\infty\right]$, then $\FourierDecompSp{\CalQ}pY$ is
a Banach space, i.e.\@ $\left\Vert \,\mybullet\right\Vert _{\FourierDecompSp{\CalQ}pY}$
is a norm.
\end{thm}

\begin{proof}
It is clear that $\FourierDecompSp{\CalQ}pY$ is closed under multiplication
with complex scalars and that $\left\Vert \mybullet\right\Vert _{\FourierDecompSp{\CalQ}pY}$
is homogeneous.

Let $\left(\varphi_{i}\right)_{i\in I}$ be an $L^{p}$-BAPU for
$\CalQ$ which is used to define $\left\Vert \mybullet\right\Vert _{\FourierDecompSp{\CalQ}pY}$.
For $f,g\in\FourierDecompSp{\CalQ}pY$ and $i\in I$, we have 
\[
\Fourier^{-1}\left[\varphi_{i}\left(f+g\right)\right]=\Fourier^{-1}\left(\varphi_{i}f\right)+\Fourier^{-1}\left(\varphi_{i}g\right)\in L^{p}\left(\smash{\R^{\dimension}}\right)
\]
with
\[
\left\Vert \Fourier^{-1}\left[\varphi_{i}\left(f+g\right)\right]\right\Vert _{L^{p}}\leq C^{\left(p\right)}\cdot\left[\left\Vert \Fourier^{-1}\left(\varphi_{i}f\right)\right\Vert _{L^{p}}+\left\Vert \Fourier^{-1}\left(\varphi_{i}g\right)\right\Vert _{L^{p}}\right].
\]
By solidity of $Y$, this implies $\left(\left\Vert \Fourier^{-1}\left[\varphi_{i}\left(f+g\right)\right]\right\Vert _{L^{p}}\right)_{i\in I}\in Y$
with
\begin{align*}
\left\Vert f+g\right\Vert _{\FourierDecompSp{\CalQ}pY} & =\left\Vert \left(\left\Vert \Fourier^{-1}\left[\varphi_{i}\left(f+g\right)\right]\right\Vert _{L^{p}}\right)_{i\in I}\right\Vert _{Y}\\
 & \leq C^{\left(p\right)}\cdot\left\Vert \left(\left\Vert \Fourier^{-1}\left(\varphi_{i}f\right)\right\Vert _{L^{p}}+\left\Vert \Fourier^{-1}\left(\varphi_{i}g\right)\right\Vert _{L^{p}}\right)_{i\in I}\right\Vert _{Y}\\
 & \leq C^{\left(p\right)}C_{Y}\cdot\left[\left\Vert \left(\left\Vert \Fourier^{-1}\left(\varphi_{i}f\right)\right\Vert _{L^{p}}\right)_{i\in I}\right\Vert _{Y}+\left\Vert \left(\left\Vert \Fourier^{-1}\left(\varphi_{i}g\right)\right\Vert _{L^{p}}\right)_{i\in I}\right\Vert _{Y}\right]\\
 & =C^{\left(p\right)}C_{Y}\cdot\left[\left\Vert f\right\Vert _{\FourierDecompSp{\CalQ}pY}+\left\Vert g\right\Vert _{\FourierDecompSp{\CalQ}pY}\right]<\infty.
\end{align*}
This shows that $\FourierDecompSp{\CalQ}pY$ is a vector space and
that $\left\Vert \mybullet\right\Vert _{\FourierDecompSp{\CalQ}pY}$
is a quasi-norm with triangle constant $C\leq C^{\left(p\right)}C_{Y}$
(at elast once we know that $\left\Vert \mybullet\right\Vert _{\FourierDecompSp{\CalQ}pY}$
is positive definite, which we show below). If $Y$ is a Banach space
and if $p\in\left[1,\infty\right]$, we can take $C^{\left(p\right)}=C_{Y}=1$,
so that $\left\Vert \mybullet\right\Vert _{\FourierDecompSp{\CalQ}pY}$
is a genuine norm.

\medskip{}

Now, let us prove $\FourierDecompSp{\CalQ}pY\hookrightarrow\DistributionSpace{\CalO}$.
To this end, let $K\subset\CalO$ be an arbitrary compact set. Lemma~\ref{lem:PartitionCoveringNecessary}
shows that $\left(Q_{i}^{\circ}\right)_{i\in I}$ covers $\CalO$.
Since $K\subset\CalO$ is compact, there are finitely many $i_{1},\dots,i_{n}\in I$
with $K\subset\bigcup_{\ell=1}^{n}Q_{i_{\ell}}^{\circ}\subset\bigcup_{\ell=1}^{n}Q_{i_{\ell}}$.
The set $I_{K}:=\left\{ i_{1},\dots,i_{n}\right\} ^{\ast}\subset I$
is finite and Lemma~\ref{lem:PartitionCoveringNecessary} implies
$\varphi_{I_{K}}\equiv1$ on $K$.

Now, choose $u=\left(u_{i}\right)_{i\in I}$ as in Lemma~\ref{lem:SolidSequenceSpaceEmbedsIntoWeightedLInfty},
so that $Y\hookrightarrow\ell_{u}^{\infty}\left(I\right)$. Let us
set $C_{K}:=\min_{i\in I_{K}}u_{i}>0$. For arbitrary $\varphi\in\TestFunctionSpace{\CalO}$
with $\supp\varphi\subset K$ and $f\in\FourierDecompSp{\CalQ}pY\subset\DistributionSpace{\CalO}$,
we now have
\[
\left|\left\langle f,\varphi\right\rangle _{\CalD'}\right|=\left|\vphantom{\sum}\left\langle f,\smash{\sum_{i\in I_{K}}}\varphi_{i}\varphi\right\rangle _{\CalD'}\right|=\left|\vphantom{\sum}\smash{\sum_{i\in I_{K}}}\left\langle \varphi_{i}f,\varphi\right\rangle _{\Schwartz'}\right|\leq\sum_{i\in I_{K}}\left|\left\langle \Fourier^{-1}\left(\varphi_{i}f\right),\widehat{\varphi}\right\rangle _{\Schwartz'}\right|.
\]
There are now two cases: For $p\in\left[1,\infty\right]$, we can
apply Hölder's inequality to conclude
\[
\left|\left\langle \Fourier^{-1}\left(\varphi_{i}f\right),\widehat{\varphi}\right\rangle _{\Schwartz'}\right|\leq\left\Vert \Fourier^{-1}\left(\varphi_{i}f\right)\right\Vert _{L^{p}}\cdot\left\Vert \widehat{\varphi}\right\Vert _{L^{p'}}.
\]
In case of $p\in\left(0,1\right)$, we observe $\supp\widehat{\Fourier^{-1}\left(\varphi_{i}f\right)}=\supp\left(\varphi_{i}f\right)\subset\overline{Q_{i}}$,
so that we can apply Corollary~\ref{cor:BandlimitedEmbedding} with
$q=\infty$ to conclude
\[
\left\Vert \Fourier^{-1}\left(\varphi_{i}f\right)\right\Vert _{L^{\infty}}\leq\left[\lambda\left(\,\overline{Q_{i}}\,\right)\right]^{\frac{1}{p}}\cdot\left\Vert \Fourier^{-1}\left(\varphi_{i}f\right)\right\Vert _{L^{p}}\,.
\]
Hence,
\[
\left|\left\langle \Fourier^{-1}\left(\varphi_{i}f\right),\widehat{\varphi}\right\rangle _{\Schwartz'}\right|\leq\left[\lambda\left(\,\overline{Q_{i}}\,\right)\right]^{\frac{1}{p}}\cdot\left\Vert \widehat{\varphi}\right\Vert _{L^{1}}\cdot\left\Vert \Fourier^{-1}\left(\varphi_{i}f\right)\right\Vert _{L^{p}}\leq C_{K}'\cdot\left\Vert \widehat{\varphi}\right\Vert _{L^{1}}\cdot\left\Vert \Fourier^{-1}\left(\varphi_{i}f\right)\right\Vert _{L^{p}}
\]
for all $i\in I_{K}$, with $C_{K}':=\max_{i\in I_{K}}\left[\lambda\left(\overline{\,Q_{i}}\,\right)\right]^{1/p}$.

In any of the two cases, there is thus an exponent $r\in\left[1,\infty\right]$
and some constant $C_{K}'>0$ satisfying
\begin{align}
\left|\left\langle f,\varphi\right\rangle _{\CalD'}\right|\leq C_{K}'\cdot\left\Vert \widehat{\varphi}\right\Vert _{L^{r}}\cdot\sum_{i\in I_{K}}\left\Vert \Fourier^{-1}\left(\varphi_{i}f\right)\right\Vert _{L^{p}} & \leq\frac{C_{K}'}{C_{K}}\cdot\left\Vert \widehat{\varphi}\right\Vert _{L^{r}}\cdot\sum_{i\in I_{K}}\left[u_{i}\cdot\left\Vert \Fourier^{-1}\left(\varphi_{i}f\right)\right\Vert _{L^{p}}\right]\nonumber \\
 & \leq\frac{C_{K}'}{C_{K}}\cdot\left\Vert \widehat{\varphi}\right\Vert _{L^{r}}\cdot\left\Vert \left(\left\Vert \Fourier^{-1}\left(\varphi_{i}f\right)\right\Vert _{L^{p}}\right)_{i\in I}\right\Vert _{\ell_{u}^{\infty}}\cdot\left|I_{K}\right|\nonumber \\
 & \leq\frac{C_{K}'\cdot\left|I_{K}\right|}{C_{K}}\cdot\left\Vert \widehat{\varphi}\right\Vert _{L^{r}}\cdot\left\Vert f\right\Vert _{\FourierDecompSp{\CalQ}pY}.\label{eq:DecompositionSpaceEmbedsInDistributions}
\end{align}
Here, the term $\left\Vert \widehat{\varphi}\right\Vert _{L^{r}}$
is finite because of $\varphi\in\TestFunctionSpace{\CalO}\subset\Schwartz\left(\R^{\dimension}\right)$.
We recall from above the assumption $\supp\varphi\subset K$, where
$K\subset\CalO$ was an \emph{arbitrary} compact subset of $\CalO$.
In particular, the above estimate proves $f\equiv0$ as an element
of $\DistributionSpace{\CalO}$ if $\left\Vert f\right\Vert _{\FourierDecompSp{\CalQ}pY}=0$.
Hence, $\left\Vert \mybullet\right\Vert _{\FourierDecompSp{\CalQ}pY}$
is positive definite.

Now, if $\left(f_{n}\right)_{n\in\N}$ is a sequence in $\FourierDecompSp{\CalQ}pY$
converging in $\FourierDecompSp{\CalQ}pY$ to $f\in\FourierDecompSp{\CalQ}pY$,
then the above estimate easily implies $\left|\left\langle f_{n}-f,\varphi\right\rangle _{\CalD'}\right|\xrightarrow[n\to\infty]{}0$
for all $\varphi\in\TestFunctionSpace{\CalO}$, and hence $f_{n}\to f$
with convergence in $\DistributionSpace{\CalO}$. This establishes
the continuous embedding $\FourierDecompSp{\CalQ}pY\hookrightarrow\DistributionSpace{\CalO}$.
We note that it suffices to consider sequences to prove this continuity,
since the topology of the quasi-normed vector space $\FourierDecompSp{\CalQ}pY$
is first countable.

\medskip{}

It remains to prove the completeness of $\FourierDecompSp{\CalQ}pY$.
To this end, note that Lemma~\ref{lem:AbsoluteConvergenceForQuasiBanachSpaces}
shows that it suffices to consider a sequence $\left(f_{n}\right)_{n\in\N}$
in $\FourierDecompSp{\CalQ}pY$ with $\left\Vert f_{n}\right\Vert _{\FourierDecompSp{\CalQ}pY}\leq\left(2C_{Y}C^{\left(p\right)}\right)^{-n}$
for all $n\in\N$ and to show that the sequence $g_{N}:=\sum_{n=1}^{N}f_{n}$
converges to some $g\in\FourierDecompSp{\CalQ}pY$.

Let us set 
\[
\theta_{i}^{\left(n\right)}:=\left\Vert \Fourier^{-1}\left(\varphi_{i}f_{n}\right)\right\Vert _{L^{p}}
\]
for $i\in I$ and $n\in\N$ and note $\theta^{\left(n\right)}:=\left(\smash{\theta_{i}^{\left(n\right)}}\right)_{i\in I}\in Y$
with $\left\Vert \smash{\theta^{\left(n\right)}}\right\Vert _{Y}=\left\Vert f_{n}\right\Vert _{\FourierDecompSp{\CalQ}pY}\leq\left(2C_{Y}C^{\left(p\right)}\right)^{-n}$.
Since $Y$ is complete and because of
\[
\sum_{n=N}^{\infty}C_{Y}^{n-N+1}\left\Vert \left(\smash{C^{\left(p\right)}}\right)^{n}\theta^{\left(n\right)}\right\Vert _{Y}\leq\sum_{n=N}^{\infty}\left(C_{Y}\smash{C^{\left(p\right)}}\right)^{n}\left\Vert \smash{\theta^{\left(n\right)}}\right\Vert _{Y}\leq\sum_{n=N}^{\infty}2^{-n}<\infty,
\]
Lemma~\ref{lem:AbsoluteConvergenceForQuasiBanachSpaces} implies
$\sum_{n=N}^{\infty}\left(C^{\left(p\right)}\right)^{n}\theta^{\left(n\right)}\in Y$
for all $N\in\N$, and also 
\begin{equation}
\left\Vert \,\smash{\sum_{n=N}^{\infty}}\,\vphantom{\sum}\left(\smash{C^{\left(p\right)}}\right)^{n}\theta^{\left(n\right)}\,\right\Vert _{Y}\leq C_{Y}\cdot\sum_{n=N}^{\infty}C_{Y}^{n-N+1}\left\Vert \left(\smash{C^{\left(p\right)}}\right)^{n}\theta^{\left(n\right)}\right\Vert _{Y}\leq C_{Y}^{2-N}\sum_{n=N}^{\infty}2^{-n}\xrightarrow[N\to\infty]{}0.\label{eq:DecompositionSpaceCompleteYRescaledSeries}
\end{equation}
Note that implicitly an argument as in equation~(\ref{eq:IteratedQuasiTriangleInequalityImplicitArgument})
was used here. In particular, since $Y\hookrightarrow\ell_{u}^{\infty}\left(I\right)$
for a suitable weight $u$ (see Lemma~\ref{lem:SolidSequenceSpaceEmbedsIntoWeightedLInfty}),
we get $\sum_{n=1}^{\infty}\left(C^{\left(p\right)}\right)^{n}\theta_{i}^{\left(n\right)}<\infty$
for all $i\in I$.

Now, we note that equation~(\ref{eq:DecompositionSpaceEmbedsInDistributions})
implies for each compact $K\subset\CalO$ and each $\varphi\in\TestFunctionSpace{\CalO}$
with $\supp\varphi\subset K$ the estimate
\[
\sum_{n=1}^{\infty}\left|\left\langle f_{n},\varphi\right\rangle _{\CalD'}\right|\leq\frac{C_{K}'\cdot\left|I_{K}\right|}{C_{K}}\cdot\left\Vert \widehat{\varphi}\right\Vert _{L^{r}}\cdot\sum_{n=1}^{\infty}\left\Vert f_{n}\right\Vert _{\FourierDecompSp{\CalQ}pY}<\infty.
\]
Thus, $g\left(\varphi\right):=\sum_{n=1}^{\infty}\left\langle f_{n},\varphi\right\rangle _{\CalD'}$
converges for every $\varphi\in\TestFunctionSpace{\CalO}$. By \cite[Theorem 6.17]{RudinFA},
this implies $g\in\DistributionSpace{\CalO}$. It remains to show
$g\in\FourierDecompSp{\CalQ}pY$ and $\left\Vert g_{N}-g\right\Vert _{\FourierDecompSp{\CalQ}pY}\to0$
as $N\to\infty$.

We recall from \cite[Theorem 7.23]{RudinFA} that the inverse Fourier
transform $\Fourier^{-1}\left(\varphi_{i}\left[g-g_{N}\right]\right)$
is given by
\begin{align}
\left[\Fourier^{-1}\left(\varphi_{i}\left[g-g_{N}\right]\right)\right]\left(x\right) & =\left\langle \varphi_{i}\left[g-g_{N}\right],e_{x}\right\rangle =\left\langle g-g_{N},\,\varphi_{i}e_{x}\right\rangle _{\CalD'}\nonumber \\
 & =\sum_{n=N+1}^{\infty}\left\langle f_{n},\varphi_{i}e_{x}\right\rangle _{\CalD'}=\sum_{n=N+1}^{\infty}\left\langle \varphi_{i}f_{n},e_{x}\right\rangle =\sum_{n=1}^{\infty}\left[\Fourier^{-1}\left(\varphi_{i}f_{n+N}\right)\right]\left(x\right)\label{eq:DecompositionCompletenessPointwiseSeries}
\end{align}
for all $x\in\R^{\dimension}$, with $e_{z}\left(y\right)=e^{2\pi i\left\langle z,y\right\rangle }$
for $y,z\in\R^{\dimension}$.

But above, we saw $\sum_{n=1}^{\infty}\left[\left(C^{\left(p\right)}\right)^{n}\cdot\left\Vert \Fourier^{-1}\left(\varphi_{i}f_{n}\right)\right\Vert _{L^{p}}\right]=\sum_{n=1}^{\infty}\left[\left(C^{\left(p\right)}\right)^{n}\cdot\theta_{i}^{\left(n\right)}\right]<\infty$
for all $i\in I$, which in particular implies 
\[
\sum_{n=1}^{\infty}\left[\left(\smash{C^{\left(p\right)}}\right)^{n}\cdot\left\Vert \Fourier^{-1}\left(\varphi_{i}f_{n+N}\right)\right\Vert _{L^{p}}\right]=\left(\smash{C^{\left(p\right)}}\right)^{-N}\cdot\sum_{\ell=N+1}^{\infty}\left[\left(\smash{C^{\left(p\right)}}\right)^{\ell}\cdot\left\Vert \mathcal{F}^{-1}\left(\varphi_{i}f_{\ell}\right)\right\Vert _{L^{p}}\right]<\infty\quad\forall\,i\in I\,.
\]
By Lemma~\ref{lem:AbsoluteConvergenceForQuasiBanachSpaces}, this
yields $\sum_{n=1}^{\infty}\Fourier^{-1}\left(\varphi_{i}f_{n+N}\right)\in L^{p}\left(\R^{\dimension}\right)$,
since $C^{\left(p\right)}$ is a triangle constant for $L^{p}$. Furthermore,
the same lemma shows
\begin{align*}
0\leq\gamma_{i}^{\left(N\right)}:=\left\Vert \smash{\sum_{n=1}^{\infty}}\,\vphantom{\sum}\Fourier^{-1}\left(\varphi_{i}\,f_{n+N}\right)\right\Vert _{L^{p}} & \leq C^{\left(p\right)}\cdot\sum_{n=1}^{\infty}\left[\left(\smash{C^{\left(p\right)}}\right)^{n}\cdot\left\Vert \Fourier^{-1}\left(\varphi_{i}\,f_{n+N}\right)\right\Vert _{L^{p}}\right]\\
 & =\sum_{n=N+1}^{\infty}\left[\left(\smash{C^{\left(p\right)}}\right)^{n-N+1}\cdot\left\Vert \Fourier^{-1}\left(\varphi_{i}\,f_{n}\right)\right\Vert _{L^{p}}\right]\\
 & \leq\sum_{n=N+1}^{\infty}\left[\left(\smash{C^{\left(p\right)}}\right)^{n}\cdot\left\Vert \Fourier^{-1}\left(\varphi_{i}\,f_{n}\right)\right\Vert _{L^{p}}\right]\\
 & =\sum_{n=N+1}^{\infty}\left[\left(\smash{C^{\left(p\right)}}\right)^{n}\cdot\theta_{i}^{\left(n\right)}\right].
\end{align*}
Note that convergence in $L^{p}\left(\R^{\dimension}\right)$ implies
convergence almost everywhere for a subsequence. Together with equation~(\ref{eq:DecompositionCompletenessPointwiseSeries}),
this implies $\Fourier^{-1}\left(\varphi_{i}\left[g-g_{N}\right]\right)=\sum_{n=1}^{\infty}\Fourier^{-1}\left(\varphi_{i}f_{n+N}\right)\in L^{p}\left(\R^{\dimension}\right)$
and thus also $\gamma_{i}^{\left(N\right)}=\left\Vert \Fourier^{-1}\left(\varphi_{i}\left[g-g_{N}\right]\right)\right\Vert _{L^{p}}$.

But above, we saw $\sum_{n=N+1}^{\infty}\left(C^{\left(p\right)}\right)^{n}\theta^{\left(n\right)}\in Y$.
By solidity of $Y$ and using equation~(\ref{eq:DecompositionSpaceCompleteYRescaledSeries}),
we conclude $\gamma^{\left(N\right)}:=\left(\smash{\gamma_{i}^{\left(N\right)}}\right)_{i\in I}\in Y$
with
\[
\left\Vert g-g_{N}\right\Vert _{\FourierDecompSp{\CalQ}pY}=\left\Vert \gamma^{\left(N\right)}\right\Vert _{Y}\leq\left\Vert \,\smash{\sum_{n=N+1}^{\infty}}\,\vphantom{\sum}\left(\smash{C^{\left(p\right)}}\right)^{n}\cdot\theta^{\left(n\right)}\,\right\Vert _{Y}\vphantom{\sum_{n=N+1}^{\infty}}\xrightarrow[N\to\infty]{\text{eq. }\eqref{eq:DecompositionSpaceCompleteYRescaledSeries}}0.
\]
In particular, we get $g=\left(g-g_{N}\right)+g_{N}\in\FourierDecompSp{\CalQ}pY$,
and the above estimate proves $g=\lim_{N\to\infty}g_{N}$ with convergence
in $\FourierDecompSp{\CalQ}pY$. As noted above, this completes the
proof. 
\end{proof}
We end this section with an example which shows that the completeness
proved above is not trivial. More precisely, we show that the completeness
may fail if the reservoir $\DistributionSpace{\CalO}$ is replaced
by $\Schwartz'\left(\R^{\dimension}\right)$.
\begin{example}
\label{exa:BorupNielsenDecompositionSpaceIncomplete}In the following,
we provide a specific example showing that the space $\CalD_{\Schwartz'}\left(\CalQ,L^{p},Y\right)$
as defined in Remark~\ref{rem:TemperedDistributionsAsReservoirIncomplete}
is in general \emph{not} complete. Specifically, we will take $p=1$
and $Y=\ell_{u}^{1}$. The covering which we will use is the ``uniform
covering'' which is usually used to define modulation spaces; but
in our case, the weight $u$ decays (for $\xi\to\infty$) much faster
than the weights used to define the classical modulation spaces.

Let $I:=\Z$, $T_{i}:=\text{id}_{\R}$ and $b_{i}:=i$ for $i\in\Z$.
Furthermore, define $Q:=\left(-\frac{3}{4},\frac{3}{4}\right)$ and
$P:=\left(-\frac{5}{8},\frac{5}{8}\right)$, as well as 
\[
Q_{i}:=T_{i}Q+b_{i}=\left(i-\tfrac{3}{4},i+\tfrac{3}{4}\right).
\]
It is then easy to see $\bigcup_{i\in I}\left(T_{i}P+b_{i}\right)=\R$
and that $\xi\in Q_{i}\cap Q_{j}\neq\emptyset$ implies $i-\tfrac{3}{4}<\xi<j+\tfrac{3}{4}$,
and hence $i-j<\frac{6}{4}<2$. Since $i-j\in\Z$, we conclude $i-j\leq1$.
By symmetry we get $\left|i-j\right|\leq1$ and thus $i^{\ast}\subset\left\{ i-1,i,i+1\right\} $.
This shows that $\CalQ=\left(Q_{i}\right)_{i\in I}$ is a structured
admissible covering of $\R$.

Now consider the weight $u_{i}:=10^{-i}$ for $i\in\Z$ and note because
of the estimate 
\[
\frac{u_{i}}{u_{j}}=10^{j-i}\leq10^{\left|j-i\right|}\leq10,
\]
which is valid for all $i\in I$ and $j\in i^{\ast}\subset\left\{ i-1,i,i+1\right\} $,
that the weight $u$ is $\CalQ$-moderate.

Theorem~\ref{thm:AlmostStructuredAdmissibleAdmitsBAPU} guarantees
the existence of an $L^{p}$-BAPU $\left(\varphi_{i}\right)_{i\in I}$
for $\CalQ$. Note that for $i\in I$ we have
\[
\bigcup_{j\in I\setminus\left\{ i\right\} }Q_{j}\subset\left(-\infty,\left(i-1\right)+\tfrac{3}{4}\right)\cup\left(\left(i+1\right)-\tfrac{3}{4},\infty\right)=\left(-\infty,i-\tfrac{1}{4}\right)\cup\left(i+\tfrac{1}{4},\infty\right).
\]
Together with $\varphi_{j}\left(\xi\right)=0$ for $\xi\in\R\setminus Q_{j}$
and $\sum_{i\in I}\varphi_{i}\left(\xi\right)=1$ for all $\xi\in\R$,
this implies $\varphi_{i}\left(\xi\right)=1$ for $\xi\in\left[i-\frac{1}{4},i+\frac{1}{4}\right]$
for all $i\in I=\Z$.

Now choose a nonnegative function $\psi\in\TestFunctionSpace{\left(-\frac{1}{4},\frac{1}{4}\right)}\setminus\left\{ 0\right\} $
and define $f_{n}:=\sum_{j=1}^{n}4^{j}\cdot L_{j}\psi$ for $n\in\N$.
Because of 
\[
\supp\left(L_{n}\psi\right)\subset\left(n-\tfrac{1}{4},n+\tfrac{1}{4}\right)\subset\biggl(\bigcup_{j\in I\setminus\left\{ n\right\} }Q_{j}\biggr)^{c}
\]
it is easy to see that 
\[
\varphi_{i}\cdot L_{n}\psi=\begin{cases}
0, & \text{if }i\neq n,\\
L_{i}\psi, & \text{if }i=n
\end{cases}\qquad\qquad\forall\,\,i,n\in\Z\,.
\]
For $n\geq m\geq m_{0}$ we thus get
\begin{align*}
\left\Vert f_{n}-f_{m}\right\Vert _{\FourierDecompSp{\CalQ}1{\ell_{u}^{1}}} & =\sum_{i\in\Z}10^{-i}\left\Vert \Fourier^{-1}\left(\,\varphi_{i}\cdot\smash{\sum_{j=m+1}^{n}}\vphantom{\sum}\,4^{j}\cdot L_{j}\psi\,\right)\right\Vert _{L^{1}}\vphantom{\sum_{j=m+1}^{n}}\\
 & =\sum_{i=m+1}^{n}10^{-i}4^{i}\cdot\left\Vert \Fourier^{-1}\left(L_{i}\psi\right)\right\Vert _{L^{1}}\leq\left\Vert \Fourier^{-1}\psi\right\Vert _{L^{1}}\cdot\sum_{i=m_{0}+1}^{\infty}\left(\frac{4}{10}\right)^{i}\xrightarrow[m_{0}\rightarrow\infty]{}0,
\end{align*}
so that $\left(\Fourier^{-1}f_{n}\right)_{n\in\N}$ is a Cauchy sequence
in $\CalD_{\Schwartz'}\left(\CalQ,L^{1},\ell_{u}^{1}\right)$, since
$\Fourier:\CalD_{\Schwartz'}\left(\CalQ,L^{1},\ell_{u}^{1}\right)\to\FourierDecompSp{\CalQ}1{\ell_{u}^{1}}$
is isometric (but not necessarily surjective).

Let us assume that there is some $f\in\CalD_{\Schwartz'}\left(\CalQ,L^{1},\ell_{u}^{1}\right)\subset\Schwartz'\left(\R\right)$
with $\Fourier^{-1}f_{n}\xrightarrow[n\rightarrow\infty]{\CalD_{\Schwartz'}\left(\CalQ,L^{1},\ell_{u}^{1}\right)}f$.
Using the continuity of 
\[
\Fourier:\CalD_{\Schwartz'}\left(\smash{\CalQ,L^{1},\ell_{u}^{1}}\right)\to\FourierDecompSp{\CalQ}1{\ell_{u}^{1}}
\]
and the continuous embedding $\FourierDecompSp{\CalQ}1{\ell_{u}^{1}}\hookrightarrow\DistributionSpace{\smash{\R}}$
(see Theorem~\ref{thm:DecompositionSpaceComplete}), we conclude
\[
\left\langle \smash{\widehat{f}}\,,g\right\rangle _{\CalD'}=\lim_{n\rightarrow\infty}\left\langle f_{n},g\right\rangle _{\CalD'}\qquad\text{for all }g\in\TestFunctionSpace{\smash{\R}}.
\]
Furthermore, by definition of the topology on $\Schwartz\left(\R\right)$,
by \cite[Proposition 5.15]{FollandRA} and because of $\widehat{f}\in\Schwartz'\left(\R\right)$
(since $f\in\Schwartz'\left(\R\right)$), there exists some $N\in\N$
and some $C>0$ such that 
\[
\left|\left\langle \smash{\widehat{f}}\,,g\right\rangle _{\Schwartz'}\right|\leq C\cdot\sup_{\substack{\alpha\in\N_{0}\\
\alpha\leq N
}
}\:\sup_{\xi\in\R}\left(1+\left|\xi\right|\right)^{N}\cdot\left|\left(\partial^{\alpha}g\right)\left(\xi\right)\right|\qquad\forall\,g\in\Schwartz\left(\R\right)\,.
\]

For $n\in\N$ and $g_{n}:=L_{n}\psi$ we have $\supp g_{n}\subset\left(n-\tfrac{1}{4},n+\tfrac{1}{4}\right)\subset\left[0,n+1\right]$
and thus
\[
\sup_{\substack{\alpha\in\N_{0}\\
\alpha\leq N
}
}\:\sup_{\xi\in\R}\left(1+\left|\xi\right|\right)^{N}\cdot\left|\left(\partial^{\alpha}g_{n}\right)\left(x\right)\right|\leq\left(n+2\right)^{N}\cdot\sup_{\substack{\alpha\in\N_{0}\\
\alpha\leq N
}
}\:\sup_{\xi\in\R}\left|\left(\partial^{\alpha}\psi\right)\left(\xi\right)\right|=\left(n+2\right)^{N}\cdot C_{\psi,N}\,,
\]
for some constant $C_{\psi,N}\in\left[0,\infty\right)$.

But because of $\supp\left(L_{n}\psi\right)\cap\supp\left(L_{i}\psi\right)=\emptyset$
for $i,n\in\Z$ with $i\neq n$, we have, for $m\geq n$, the identity
\[
\left\langle f_{m},g_{n}\right\rangle _{\Schwartz'}=\sum_{j=1}^{m}4^{j}\cdot\left\langle L_{j}\psi,L_{n}\psi\right\rangle _{\Schwartz'}=4^{n}\cdot\left\langle \psi,\psi\right\rangle _{\Schwartz'}\,\,,
\]
and thus
\begin{align*}
4^{n}\cdot\left\Vert \psi\right\Vert _{L^{2}}^{2} & =\lim_{m\rightarrow\infty}\left|\left\langle f_{m},g_{n}\right\rangle _{\Schwartz'}\right|=\left|\left\langle \smash{\widehat{f}}\,,g_{n}\right\rangle _{\Schwartz'}\right|\\
 & \leq C\cdot\sup_{\substack{\alpha\in\N_{0}\\
\alpha\leq N
}
}\:\sup_{\xi\in\R}\left(1+\left|\xi\right|\right)^{N}\cdot\left|\left(\partial^{\alpha}g_{n}\right)\left(\xi\right)\right|\leq CC_{\psi,N}\cdot\left(n+2\right)^{N}
\end{align*}
for all $n\in\N$, a contradiction.

Thus, there is no $f\in\CalD_{\Schwartz'}\left(\CalQ,L^{1},\ell_{u}^{1}\right)$
with $\left\Vert f-\Fourier^{-1}f_{n}\right\Vert _{\DecompSp{\CalQ}1{\ell_{u}^{1}}}\xrightarrow[n\rightarrow\infty]{}0$,
so that $\CalD_{\Schwartz'}\left(\CalQ,L^{1},\ell_{u}^{1}\right)$
is \emph{not} complete.
\end{example}

\section{Nested sequence spaces}

\label{sec:NestedSequenceSpaces}Recall from Section~\ref{subsec:IntroductionOverview}
that most of our embedding results will consist in showing that an
embedding between certain sequence spaces is sufficient (or necessary)
for the existence of an embedding
\[
\FourierDecompSp{\CalQ}{p_{1}}Y\hookrightarrow\FourierDecompSp{\CalP}{p_{2}}Z.
\]
Furthermore, as seen for instance in Theorems~\ref{thm:IntroductionFineIntoCoarse}
and \ref{thm:IntroductionCoarseIntoFine}, these sequence spaces are
often of a \emph{nested} nature, so that the (quasi)-norm is for example
given by
\[
\left\Vert \left(x_{i}\right)_{i\in I}\right\Vert _{X\left(\left[\ell_{u}^{q}\left(I_{k}\right)\right]_{k\in K}\right)}=\left\Vert \left(\left\Vert \left(u_{k,i}\cdot x_{i}\right)_{i\in I_{k}}\right\Vert _{\ell^{q}}\right)_{k\in K}\right\Vert _{X}.
\]
Thus, the present section is devoted to studying spaces of this type.
In particular, we will consider the existence of embeddings between
such sequence spaces. We begin with a rather general definition.
\begin{defn}
\label{def:NestedSequenceSpaces}Let $K$ be an index set and let
$X\subset\Compl^{K}$ be a solid sequence space. For each $k\in K$,
let $I_{k}$ be an index set, and let $X_{k}\subset\Compl^{I_{k}}$
be a solid sequence space. Then, setting $I:=\bigcup_{k\in K}I_{k}$,
we define the \textbf{nested sequence space}
\[
X\left(\left[X_{k}\right]_{k\in K}\right):=\left\{ x=\left(x_{i}\right)_{i\in I}\in\Compl^{I}\with\left(\forall\,k\in K:\,\left(x_{i}\right)_{i\in I_{k}}\in X_{k}\right)\text{ and }\left(\left\Vert \left(x_{i}\right)_{i\in I_{k}}\right\Vert _{X_{k}}\right)_{k\in K}\in X\right\} \,.
\]
If the context makes the intended meaning clear, we will also write
$X\left(X_{k}\right):=X\left(\left[X_{k}\right]_{k\in K}\right)$.

Finally, for $x=\left(x_{i}\right)_{i\in I}\in X\left(\left[X_{k}\right]_{k\in K}\right)$,
we define
\[
\left\Vert x\right\Vert _{X\left(\left[X_{k}\right]_{k\in K}\right)}:=\left\Vert \left[\left\Vert \left(x_{i}\right)_{i\in I_{k}}\right\Vert _{X_{k}}\right]_{k\in K}\right\Vert _{X}.\qedhere
\]
\end{defn}

\begin{rem*}
(1) We emphasize that the sets $\left(I_{k}\right)_{k\in K}$ need
\emph{not} be pairwise disjoint.

(2) We will often have $X_{k}=\ell_{w_{k}}^{q_{k}}\left(I_{k}\right)$
for a suitable weight $w_{k}=\left(w_{k,i}\right)_{i\in I_{k}}$.
In this case, we will write $X\left(\left[\ell_{w}^{q_{k}}\left(I_{k}\right)\right]_{k\in K}\right)$,
or even $X\left(\ell_{w}^{q_{k}}\left(I_{k}\right)\right)$, instead
of the (more correct) $X\left(\smash{\left[\ell_{w_{k}}^{q_{k}}\left(I_{k}\right)\right]_{k\in K}}\right)$.
\end{rem*}
The first question is whether $X\left(\left[X_{k}\right]_{k\in K}\right)$
is a solid sequence space over the index set $I=\bigcup_{k\in K}I_{k}$,
or even whether it is a vector space. Under very weak assumptions,
this is indeed the case:
\begin{lem}
\label{lem:NestedSequenceSpaceIsQuasiNormed}In the setting of Definition~\ref{def:NestedSequenceSpaces},
let $C_{k}\geq1$ be a triangle constant for $X_{k}$ for each $k\in K$.
Assume that the following expression (then a constant) is finite:
\[
C:=\sup_{k\in K}C_{k}\,.
\]
Then $X\left(\left[X_{k}\right]_{k\in K}\right)$ is a solid sequence
space over $I$, with triangle constant $C\cdot C_{0}$, where $C_{0}\geq1$
is a triangle constant for $X$.
\end{lem}

\begin{rem*}

\begin{enumerate}[leftmargin=0.7cm]
\item Under the assumptions from above, it is even true that $X\left(\left[X_{k}\right]_{k\in K}\right)$
is complete if $X$ and each $X_{k}$ are complete. Since we will
not need this fact, we omit the proof.\vspace{0.2cm}
\item If $X_{k}=\ell_{u^{\left(k\right)}}^{q_{k}}\left(I_{k}\right)$ for
a certain weight $u^{\left(k\right)}=\left(\smash{{u_{i}^{\left(k\right)}}}\right)_{i\in I_{k}}$,
then \cite[Exercise 1.1.5(c)]{GrafakosClassical} implies that a triangle
constant for $X_{k}$ is given by $C_{k}=\vphantom{\sum_{j}}\max\left\{ 1,\smash{2^{q_{k}^{-1}-1}}\right\} \vphantom{2^{q_{k}^{-1}-1}}$.
Thus, if $\inf_{k\in K}q_{k}>0$, then the lemma is applicable and
shows that $X\left(\vphantom{\ell_{u^{\left(k\right)}}^{q_{k}}}\smash{{\left[\ell_{u^{\left(k\right)}}^{q_{k}}\left(I_{k}\right)\right]_{k\in K}}}\right)$
is a solid sequence space over $I$.\qedhere
\end{enumerate}
\end{rem*}
\begin{proof}
For brevity, write $Y:=X\left(\left[X_{k}\right]_{k\in K}\right)$.
It is clear that $Y$ is closed under multiplication with complex
scalars and that $\left\Vert \mybullet\right\Vert _{Y}$ is homogeneous.

Furthermore, $\left\Vert \mybullet\right\Vert _{Y}$ is definite:
Assume that $x=\left(x_{i}\right)_{i\in I}\in Y$ satisfies $\left\Vert x\right\Vert _{Y}=0$.
For arbitrary $i_{0}\in I$, there is some $k_{0}\in K$ with $i_{0}\in I_{k_{0}}$.
But $\left\Vert x\right\Vert _{Y}=0$ implies $\left\Vert \smash{\left(x_{i}\right)_{i\in I_{k_{0}}}}\vphantom{\left(x_{i}\right)}\right\Vert _{X_{k_{0}}}=0$
and thus $x_{i_{0}}=0$. Since $i_{0}\in I$ was arbitrary, we see
$x=0$.

To see that $Y$ is solid, let $x=\left(x_{i}\right)_{i\in I}\in\Compl^{I}$
and $y=\left(y_{i}\right)_{i\in I}\in Y$ with $\left|x_{i}\right|\leq\left|y_{i}\right|$
for all $i\in I$. Since each $X_{k}$ is solid with $\left(y_{i}\right)_{i\in I_{k}}\in X_{k}$,
we get $\left(x_{i}\right)_{i\in I_{k}}\in X_{k}$ and $\varrho_{k}:=\left\Vert \left(x_{i}\right)_{i\in I_{k}}\right\Vert _{X_{k}}\leq\left\Vert \left(y_{i}\right)_{i\in I_{k}}\right\Vert _{X_{k}}=:\theta_{k}$.
But $y\in Y$ means $\theta=\left(\theta_{k}\right)_{k\in K}\in X$.
By solidity of $X$, we get $\varrho=\left(\varrho_{k}\right)_{k\in K}\in X$
and hence $x\in Y$, with $\left\Vert x\right\Vert _{Y}=\left\Vert \varrho\right\Vert _{X}\leq\left\Vert \theta\right\Vert _{X}=\left\Vert y\right\Vert _{Y}$,
as desired.

Finally, to show that $Y$ is a vector space and that $\left\Vert \mybullet\right\Vert _{Y}$
is a (quasi)-norm, let $x=\left(x_{i}\right)_{i\in I}\in Y$ and $y=\left(y_{i}\right)_{i\in I}\in Y$
and define $x^{\left(k\right)}:=\left\Vert \left(x_{i}\right)_{i\in I_{k}}\right\Vert _{X_{k}}$
and $y^{\left(k\right)}:=\left\Vert \left(y_{i}\right)_{i\in I_{k}}\right\Vert _{X_{k}}$
for $k\in K$. By definition of $Y$, we have $\left(x^{\left(k\right)}\right)_{k\in K},\left(y^{\left(k\right)}\right)_{k\in K}\in X$
and thus also $\left(C\cdot\left[x^{\left(k\right)}+y^{\left(k\right)}\right]\right)_{k\in K}\in X$.
But we also have
\[
z^{\left(k\right)}:=\left\Vert \left(x_{i}+y_{i}\right)_{i\in I_{k}}\right\Vert _{X_{k}}\leq C_{k}\cdot\left[\left\Vert \left(x_{i}\right)_{i\in I_{k}}\right\Vert _{X_{k}}+\left\Vert \left(y_{i}\right)_{i\in I_{k}}\right\Vert _{X_{k}}\right]\leq C\cdot\left[x^{\left(k\right)}+y^{\left(k\right)}\right]
\]
for all $k\in K$. By solidity of $X$, this yields $\left(z^{\left(k\right)}\right)_{k\in K}\in X$
and thus $x+y\in Y$ with
\begin{align*}
\left\Vert x+y\right\Vert _{Y}=\left\Vert \left(\smash{{z^{\left(k\right)}}}\right)_{k\in K}\right\Vert _{X} & \leq\left\Vert \left(C\cdot\left[\smash{x^{\left(k\right)}}+\smash{y^{\left(k\right)}}\right]\right)_{k\in K}\right\Vert _{X}\\
 & \leq C\cdot C_{0}\cdot\left(\left\Vert \left(\smash{{x^{\left(k\right)}}}\right)_{k\in K}\right\Vert _{X}+\left\Vert \left(\smash{{y^{\left(k\right)}}}\right)_{k\in K}\right\Vert _{X}\right)\\
 & =CC_{0}\cdot\left(\left\Vert x\right\Vert _{Y}+\left\Vert y\right\Vert _{Y}\right).
\end{align*}
Thus, $Y$ is a vector space and $CC_{0}$ is a triangle constant
for $\left\Vert \mybullet\right\Vert _{Y}$.
\end{proof}
An important property of (solid) sequence spaces is the Fatou property,
which is an abstraction of the conclusion of Fatou's lemma for the
spaces $\ell_{u}^{q}$ to general (solid) sequence spaces:
\begin{defn}
\label{def:FatouProperty}(cf.\@ \cite[Definition and Theorem 3 in §65]{ZaanenIntegration})

Let $I$ be an index set, and let $X\subset\Compl^{I}$ be a solid
sequence space. We say that $X$ has the \textbf{Fatou property},
if for each bounded sequence $\left(x^{\left(n\right)}\right)_{n\in\N}\in X^{\N}$
with $x^{\left(n\right)}=\left(\vphantom{x^{\left(n\right)}}\smash{{x_{i}^{\left(n\right)}}}\right)_{i\in I}$
and
\[
x_{i}:=\liminf_{n\to\infty}\bigl|x_{i}^{\left(n\right)}\bigr|\qquad\text{ for }i\in I,
\]
we have $x=\left(x_{i}\right)_{i\in I}\in X$ with $\left\Vert x\right\Vert _{X}\leq\liminf_{n\to\infty}\left\Vert x^{\left(n\right)}\right\Vert _{X}$.
\end{defn}

\begin{rem*}
By Lemma~\ref{lem:SolidSequenceSpaceEmbedsIntoWeightedLInfty}, we
have $X\hookrightarrow\ell_{u}^{\infty}\left(I\right)$ for a suitable
weight $u=\vphantom{\sum_{j}}\left(u_{i}\right)_{i\in I}$. Thus,
each bounded sequence in $X$ is also bounded component-wise, so that
$x_{i}=\liminf_{n\to\infty}\bigl|x_{i}^{\left(n\right)}\bigr|<\infty$
for all $i\in I$. 
\end{rem*}
For later use, we need to know that the Fatou property is inherited
by nested sequence spaces and also by weighted sequence spaces:
\begin{defn}
\label{def:WeightedSequenceSpaces}Let $I$ be an index set, let $X\subset\Compl^{I}$
be a sequence space, and let $u=\left(u_{i}\right)_{i\in I}\in\left(0,\infty\right)^{I}$
be a weight on $I$. The \textbf{weighted sequence space} $X_{u}$
is given by
\[
X_{u}:=\left\{ \left(x_{i}\right)_{i\in I}\in\Compl^{I}\with\left(u_{i}\cdot x_{i}\right)_{i\in I}\in X\right\} \qquad\text{ with }\qquad\left\Vert \left(x_{i}\right)_{i\in I}\right\Vert _{X_{u}}:=\left\Vert \left(u_{i}\cdot x_{i}\right)_{i\in I}\right\Vert _{X}.\qedhere
\]
\end{defn}

\begin{rem*}
It is easy to see that $X_{u}$ is a solid sequence space if $X$
is, even with the same triangle constant. Furthermore, $X$ is complete
iff $X_{u}$ is.
\end{rem*}
Now, we can state and prove the ``inheritance properties'' of the
Fatou property:
\begin{lem}
\label{lem:FatouPropertyIsInherited}

\begin{enumerate}[leftmargin=0.7cm]
\item In the setting of Lemma~\ref{lem:NestedSequenceSpaceIsQuasiNormed},
if $X$ and each $X_{k}$ satisfy the Fatou property, then so does
$X\left(\left[X_{k}\right]_{k\in K}\right)$.\vspace{0.2cm}
\item If $X\subset\Compl^{I}$ is a solid sequence space with the Fatou
property and $u=\left(u_{i}\right)_{i\in I}$ is a weight, then $X_{u}$
has the Fatou property.\qedhere
\end{enumerate}
\end{lem}

\begin{proof}
Ad (2): Let $\left(x^{\left(n\right)}\right)_{n\in\N}\in X_{u}^{\N}$
be a bounded sequence with $x^{\left(n\right)}=\left(\smash{{x_{i}^{\left(n\right)}}}\right)_{i\in I}$.
For 
\[
y^{\left(n\right)}:=\left(\smash{{y_{i}^{\left(n\right)}}}\right)_{i\in I}:=\left(\smash{{u_{i}\cdot x_{i}^{\left(n\right)}}}\right)_{i\in I}\,,
\]
this means that $\left(y^{\left(n\right)}\right)_{n\in\N}\in X^{\N}$
is a bounded sequence. Thus, for
\[
y_{i}:=\liminf_{n\to\infty}\bigl|y_{i}^{\left(n\right)}\bigr|=u_{i}\cdot\liminf_{n\to\infty}\bigl|x_{i}^{\left(n\right)}\bigr|\,,
\]
we have $y=\left(y_{i}\right)_{i\in I}\in X$ with $\left\Vert y\right\Vert _{X}\leq\liminf_{n\to\infty}\left\Vert \smash{y^{\left(n\right)}}\right\Vert _{X}=\liminf_{n\to\infty}\left\Vert \smash{x^{\left(n\right)}}\right\Vert _{X_{u}}$.
If we define $x_{i}:=\liminf_{n\to\infty}\left|\smash{x_{i}^{\left(n\right)}}\right|$,
this yields $x=\left(x_{i}\right)_{i\in I}\in X_{u}$ and $\left\Vert x\right\Vert _{X_{u}}=\left\Vert y\right\Vert _{X}\leq\liminf_{n\to\infty}\left\Vert x^{\left(n\right)}\right\Vert _{X_{u}}$,
as desired.

\medskip{}

Ad (1): For brevity, write $Y:=X\left(\left[X_{k}\right]_{k\in K}\right)$.
Let $\left(x^{\left(n\right)}\right)_{n\in\N}\in X^{\N}$ be a bounded
sequence, with $x^{\left(n\right)}=\left(\smash{{x_{i}^{\left(n\right)}}}\right)_{i\in I}$.
By solidity of $X$ and because of Lemma~\ref{lem:SolidSequenceSpaceEmbedsIntoWeightedLInfty},
we have $X\hookrightarrow\ell_{u}^{\infty}\left(I\right)$ for a certain
weight $u$ on $I$. In particular, we get
\[
z_{k}^{\left(n\right)}:=\left\Vert \left(\smash{x_{i}^{\left(n\right)}}\right)_{i\in I_{k}}\right\Vert _{X_{k}}\lesssim_{\,k}\,\left\Vert \left[\left\Vert \left(\smash{x_{i}^{\left(n\right)}}\right)_{i\in I_{\ell}}\right\Vert _{X_{\ell}}\right]_{\ell\in K}\right\Vert _{X}=\left\Vert \smash{x^{\left(n\right)}}\right\Vert _{Y}\leq C\qquad\forall\,k\in K.
\]
Therefore, for each $k\in K$, the sequence $\left(y^{\left(k,n\right)}\right)_{n\in\N}$,
given by $y^{\left(k,n\right)}:=\left(\smash{x_{i}^{\left(n\right)}}\right)_{i\in I_{k}}$
is bounded in $X_{k}$. Furthermore, we have $\left\Vert \left(\smash{z_{k}^{\left(n\right)}}\right)_{k\in K}\right\Vert _{X}=\left\Vert \smash{x^{\left(n\right)}}\right\Vert _{Y}\leq C$
for all $n\in\N$, so that the sequence ${\left(z^{\left(n\right)}\right)_{n\in\N}\in X^{\N}}$,
given by $z^{\left(n\right)}=\vphantom{z_{k}^{\left(n\right)}}\left(\smash{z_{k}^{\left(n\right)}}\right)_{k\in K}$,
is bounded.

Now, for $i\in I$, define $x_{i}:=\liminf_{n\to\infty}\bigl|x_{i}^{\left(n\right)}\bigr|$
and for $k\in K$, let $z_{k}:=\liminf_{n\to\infty}z_{k}^{\left(n\right)}$.
Using the Fatou property of $X_{k}$, we get $\left(x_{i}\right)_{i\in I_{k}}\in X_{k}$,
with
\[
0\leq\left\Vert \left(x_{i}\right)_{i\in I_{k}}\right\Vert _{X_{k}}\leq\liminf_{n\to\infty}\left\Vert \left(\smash{x_{i}^{\left(n\right)}}\right)_{i\in I_{k}}\right\Vert _{X_{k}}=\liminf_{n\to\infty}z_{k}^{\left(n\right)}=z_{k}\qquad\forall\,k\in K\,.
\]
But the Fatou property for $X$ yields $z=\left(z_{k}\right)_{k\in K}\in X$
with
\[
\left\Vert z\right\Vert _{X}\leq\liminf_{n\to\infty}\left\Vert \smash{z^{\left(n\right)}}\right\Vert _{X}=\liminf_{n\to\infty}\left\Vert \smash{x^{\left(n\right)}}\right\Vert _{Y}\,.
\]
By solidity of $X$, this finally yields $\left(\left\Vert \left(x_{i}\right)_{i\in I_{k}}\right\Vert _{X_{k}}\right)_{k\in K}\in X$
and thus $x=\left(x_{i}\right)_{i\in I}\in Y$ with
\[
\left\Vert x\right\Vert _{Y}=\left\Vert \left[\left\Vert \left(x_{i}\right)_{i\in I_{k}}\right\Vert _{X_{k}}\right]_{k\in K}\right\Vert _{X}\leq\left\Vert \left(z_{k}\right)_{k\in K}\right\Vert _{X}\leq\liminf_{n\to\infty}\left\Vert \smash{x^{\left(n\right)}}\right\Vert _{Y}.\qedhere
\]
\end{proof}
As noted above, we are interested in embeddings between sequence spaces.
Our next lemma shows that if the ``target space'' of the embedding
satisfies the Fatou property, it suffices to verify boundedness of
an embedding between \emph{solid} sequence spaces on the space of
finitely supported sequences. Although rather technical, this fact
will be very helpful for us in the remainder of the paper.
\begin{lem}
\label{lem:FinitelySupportedSequencesSufficeUnderFatouProperty}Let
$I$ be a \emph{countable(!)} set and let $X,Y\subset\Compl^{I}$
be solid sequence spaces. If $Y$ satisfies the Fatou property, then
$\iota:X\hookrightarrow Y$ is well-defined and bounded if and only
if 
\[
\iota_{0}:X\cap\ell_{0}\left(I\right)\hookrightarrow Y
\]
is well-defined and bounded. In this case, $\vertiii{\iota}=\vertiii{\iota_{0}}$.

Here, $\ell_{0}\left(I\right)=\left\langle \delta_{i}\with i\in I\right\rangle $
is the space of \textbf{finitely supported sequences} on $I$.
\end{lem}

\begin{proof}
The implication ``$\Rightarrow$'' and the estimate $\vertiii{\iota_{0}}\leq\vertiii{\iota}$
are trivial. Hence, we only prove ``$\Leftarrow$''. By assumption,
$I=\bigcup_{n\in\N}I^{\left(n\right)}$ for certain finite subsets
$I^{\left(n\right)}\subset I$ with $I^{\left(n\right)}\subset I^{\left(n+1\right)}$
for all $n\in\N$. Let $x=\left(x_{i}\right)_{i\in I}\in X$ be arbitrary
and define $x^{\left(n\right)}:=x\cdot\Indicator_{I^{\left(n\right)}}$.
Note (thanks to the solidity of $X$) that $x^{\left(n\right)}\in X\cap\ell_{0}\left(I\right)\subset Y$
and that
\[
\left|x_{i}\right|=\liminf_{n\to\infty}\left|\smash{x_{i}^{\left(n\right)}}\right|\qquad\forall\,i\in I.
\]
Furthermore, boundedness of $\iota_{0}$ and solidity of $X$ imply
$\left\Vert x^{\left(n\right)}\right\Vert _{Y}\leq\vertiii{\iota_{0}}\cdot\left\Vert x^{\left(n\right)}\right\Vert _{X}\leq\vertiii{\iota_{0}}\cdot\left\Vert x\right\Vert _{X}$
for all $n\in\N$. In particular, $\left(x^{\left(n\right)}\right)_{n\in\N}$
is a bounded sequence in $Y$.

Using the Fatou property and the solidity of $Y$, we conclude $x\in Y$,
with 
\[
\left\Vert x\right\Vert _{Y}=\left\Vert \left(\left|x_{i}\right|\right)_{i\in I}\right\Vert _{Y}\leq\liminf_{n\to\infty}\left\Vert \smash{x^{\left(n\right)}}\right\Vert _{Y}\leq\vertiii{\iota_{0}}\cdot\left\Vert x\right\Vert _{X}<\infty.
\]
Thus, $\iota$ is well-defined and bounded with $\vertiii{\iota}\leq\vertiii{\iota_{0}}$.
\end{proof}
Now, we analyze the existence of embeddings between a pair of \emph{nested}
sequence spaces. For this, we will \emph{first} assume that the ``inner
index sets'' $\left(I_{k}\right)_{k\in K}$ are disjoint. Furthermore,
we assume that the \emph{same} index sets are used on \emph{both}
sides of the embedding. In this case, the existence of the embedding
$X\left(\left[X_{k}\right]_{k\in K}\right)\hookrightarrow Y\left(\left[Y_{k}\right]_{k\in K}\right)$
is equivalent to the existence of the embedding $X\hookrightarrow Y_{\theta}$
for a suitable weight $\theta$:
\begin{lem}
\label{lem:NestedDisjointEmbedding}Let $K$ be an index set, and
let $X,Y\subset\Compl^{K}$ be two solid sequence spaces. For each
$k\in K$, let $I^{\left(k\right)}$ be an index set, and let $X_{k},Y_{k}\subset\Compl^{I^{\left(k\right)}}$
be solid sequence spaces. Finally, assume that the $\left(I^{\left(k\right)}\right)_{k\in K}$
are \emph{pairwise disjoint} and that the triangle constants for the
spaces $X_{k}$ and $Y_{k}$ are uniformly bounded.\vspace{0.1cm}

Set $I:=\biguplus_{k\in K}I^{\left(k\right)}$ and define the sequence
$\left(\theta_{k}\right)_{k\in K}\in\left[0,\infty\right]^{K}$ by\footnote{In case of $\delta_{k}\notin X$, one can choose $\theta_{k}\in\left(0,\infty\right)$
arbitrarily; the statement of the theorem still holds. We only make
the choice $\theta_{k}=1$ for definiteness.}
\[
\theta_{k}:=\begin{cases}
\vertiii{\iota_{k}} & \text{if }\delta_{k}\in X\text{ and }\iota_{k}:X_{k}\hookrightarrow Y_{k}\text{ is well-defined and bounded},\\
1, & \text{if }\delta_{k}\notin X,\\
\infty, & \text{otherwise}.
\end{cases}
\]
Then, the following are equivalent:

\begin{enumerate}
\item The embedding $\iota:X\left(\left[X_{k}\right]_{k\in K}\right)\hookrightarrow Y\left(\left[Y_{k}\right]_{k\in K}\right)$
is well-defined and bounded.
\item We have $\theta_{k}<\infty$ for all $k\in K$ and the map\footnote{The boundedness of $\gamma$ is essentially the same as the boundedness
of the embedding $X\hookrightarrow Y_{\theta}$. The only drawback
of this formulation is that $Y_{\theta}$ is in general not a well-defined
sequence space, if $\theta_{k}=0$ for some $k\in K$. This happens
if and only if $\delta_{k}\in X$ and $I^{\left(k\right)}=\emptyset$.
In particular, if $I^{\left(k\right)}\neq\emptyset$ for all $k\in K$,
the two formulations are equivalent.} $\gamma:X\to Y,\left(x_{k}\right)_{k\in K}\mapsto\left(\theta_{k}\cdot x_{k}\right)_{k\in K}$
is well-defined and bounded.
\end{enumerate}
In this case, we even have $\vertiii{\iota}=\vertiii{\gamma}$.

The implication $\left(2\right)\!\Rightarrow\!\left(1\right)$ (with
$\vertiii{\iota}\!\leq\!\vertiii{\gamma}$) even holds if the $\left(I^{\left(k\right)}\right)_{k\in K}$
are \emph{not} pairwise disjoint.
\end{lem}

\begin{proof}
``$\left(1\right)\Rightarrow\left(2\right)$'': Let $k_{0}\in K$.
In case of $\delta_{k_{0}}\notin X$, $\theta_{k_{0}}=1<\infty$ is
trivial. Thus, assume $\delta_{k_{0}}\in X$. Let $y=\left(y_{i}\right)_{i\in I^{\left(k_{0}\right)}}\in X_{k_{0}}$
be arbitrary and define
\[
x_{i}:=\begin{cases}
0, & \text{if }i\notin I^{\left(k_{0}\right)},\\
y_{i}, & \text{if }i\in I^{\left(k_{0}\right)}.
\end{cases}
\]
Since the sets $\left(I^{\left(k\right)}\right)_{k\in K}$ are pairwise
disjoint, we get
\[
\left\Vert \left(x_{i}\right)_{i\in I^{\left(k\right)}}\right\Vert _{X_{k}}=\begin{cases}
0=\left\Vert y\right\Vert _{X_{k_{0}}}\cdot\left(\delta_{k_{0}}\right)_{k}, & \text{if }k\neq k_{0},\\
\left\Vert y\right\Vert _{X_{k_{0}}}=\left\Vert y\right\Vert _{X_{k_{0}}}\cdot\left(\delta_{k_{0}}\right)_{k}, & \text{if }k=k_{0}.
\end{cases}
\]
Because of $\delta_{k_{0}}\in X$, this implies $x=\left(x_{i}\right)_{i\in I}\in X\left(\left[X_{k}\right]_{k\in K}\right)$
with 
\[
\left\Vert x\right\Vert _{X\left(\left[X_{k}\right]_{k\in K}\right)}=\left\Vert \left(\left\Vert y\right\Vert _{X_{k_{0}}}\cdot\left(\delta_{k_{0}}\right)_{k}\right)_{k\in K}\right\Vert _{X}=\left\Vert y\right\Vert _{X_{k_{0}}}\cdot\left\Vert \delta_{k_{0}}\right\Vert _{X}<\infty.
\]
Since $\iota$ is well-defined and bounded, we get $x\in Y\left(\left[Y_{k}\right]_{k\in K}\right)$.
Furthermore, Lemma~\ref{lem:SolidSequenceSpaceEmbedsIntoWeightedLInfty}
shows that each of the coordinate maps $Y\to\Compl,\left(z_{k}\right)_{k\in K}\mapsto z_{k_{0}}$
is bounded. Hence,
\begin{align*}
\left\Vert y\right\Vert _{Y_{k_{0}}}=\left\Vert \left(x_{i}\right)_{i\in I^{\left(k_{0}\right)}}\right\Vert _{Y_{k_{0}}} & \lesssim_{\,k_{0}}\left\Vert \left(\left\Vert \left(x_{i}\right)_{i\in I^{\left(k\right)}}\right\Vert _{Y_{k}}\right)_{k\in K}\right\Vert _{Y}=\left\Vert x\right\Vert _{Y\left(\left[Y_{k}\right]_{k\in K}\right)}\\
 & \leq\vertiii{\iota}\cdot\left\Vert x\right\Vert _{X\left(\left[X_{k}\right]_{k\in K}\right)}\leq\vertiii{\iota}\left\Vert \delta_{k_{0}}\right\Vert _{X}\cdot\left\Vert y\right\Vert _{X_{k_{0}}}.
\end{align*}
Therefore, $\iota_{k_{0}}:X_{k_{0}}\hookrightarrow Y_{k_{0}}$ is
well-defined and bounded with $\theta_{k_{0}}=\vertiii{\iota_{k_{0}}}\lesssim_{\,k_{0}}\vertiii{\iota}\cdot\left\Vert \delta_{k_{0}}\right\Vert _{X}<\infty$.

\medskip{}

Thus, it remains to show that $\gamma$ is well-defined and bounded,
with $\vertiii{\gamma}\leq\vertiii{\iota}$. To this end, let $\varepsilon\in\left(0,1\right)$
be arbitrary and let $x=\left(x_{k}\right)_{k\in K}\in X$. For each
$k\in K$ with $x_{k}\neq0$, we have $\delta_{k}\in X$ by solidity.
By what we just showed, this implies that $\iota_{k}$ is well-defined
and bounded. Thus (by definition of the operator norm), there is some
sequence $y^{\left(k\right)}=\left(\smash{y_{i}^{\left(k\right)}}\right)_{i\in I^{\left(k\right)}}$
with $\left\Vert y^{\left(k\right)}\right\Vert _{X_{k}}=1$ and such
that $\left\Vert y^{\left(k\right)}\right\Vert _{Y_{k}}\geq\left(1-\varepsilon\right)\vertiii{\iota_{k}}=\left(1-\varepsilon\right)\theta_{k}$.
For $k\in K$ with $x_{k}=0$, we set $y^{\left(k\right)}:=0\in X_{k}\cap Y_{k}\subset\Compl^{I^{\left(k\right)}}$.
Note that we always have
\begin{equation}
\left|x_{k}\right|\cdot\left\Vert \smash{y^{\left(k\right)}}\right\Vert _{Y_{k}}\geq\begin{cases}
\left|x_{k}\right|\cdot\left(1-\varepsilon\right)\theta_{k}\,, & \text{if }x_{k}\neq0,\\
0=\left|x_{k}\right|\cdot\left(1-\varepsilon\right)\theta_{k}\,, & \text{if }x_{k}=0.
\end{cases}\label{eq:NestedDisjointEmbeddingDoubleSequenceProperty1}
\end{equation}

Now, define a sequence $z=\left(z_{i}\right)_{i\in I}$ by
\[
z_{i}:=x_{k}\cdot y_{i}^{\left(k\right)}\text{ for the unique }k=k_{i}\in K\text{ with }i\in I^{\left(k\right)}.
\]
Note that $k$ is uniquely determined since $I=\biguplus_{k\in K}I^{\left(k\right)}$,
by assumption. Again by disjointness of the $\left(I^{\left(k\right)}\right)_{k\in K}$,
we have for $k\in K$ that
\[
\left\Vert \left(z_{i}\right)_{i\in I^{\left(k\right)}}\right\Vert _{X_{k}}=\left\Vert \left(x_{k}\cdot\smash{y_{i}^{\left(k\right)}}\right)_{i\in I^{\left(k\right)}}\right\Vert _{X_{k}}=\left|x_{k}\right|\cdot\left\Vert \smash{y^{\left(k\right)}}\right\Vert _{X_{k}}=\begin{cases}
0=\left|x_{k}\right|, & \text{if }x_{k}=0,\\
\left|x_{k}\right|, & \text{if }x_{k}\neq0.
\end{cases}
\]
By solidity of $X$ and since $x=\left(x_{k}\right)_{k\in K}\in X$,
we conclude $z\in X\left(\left[X_{k}\right]_{k\in K}\right)$.

Since $\iota$ is well-defined and bounded, we get $z\in Y\left(\left[Y_{k}\right]_{k\in K}\right)$,
with
\[
\left\Vert z\right\Vert _{Y\left(\left[Y_{k}\right]_{k\in K}\right)}\leq\vertiii{\iota}\cdot\left\Vert z\right\Vert _{X\left(\left[X_{k}\right]_{k\in K}\right)}=\vertiii{\iota}\cdot\left\Vert \left(\left|x_{k}\right|\right)_{k\in K}\right\Vert _{X}=\vertiii{\iota}\cdot\left\Vert x\right\Vert _{X}.
\]
But the disjointness of $\left(I^{\left(k\right)}\right)_{k\in K}$
finally implies
\[
\left\Vert \left(z_{i}\right)_{i\in I^{\left(k\right)}}\right\Vert _{Y_{k}}=\left\Vert \left(x_{k}\cdot\smash{y_{i}^{\left(k\right)}}\right)_{i\in I^{\left(k\right)}}\right\Vert _{Y_{k}}=\left|x_{k}\right|\cdot\left\Vert \smash{y^{\left(k\right)}}\right\Vert _{Y_{k}}\overset{\text{Eq. }\eqref{eq:NestedDisjointEmbeddingDoubleSequenceProperty1}}{\geq}\left(1-\varepsilon\right)\cdot\theta_{k}\cdot\left|x_{k}\right|\,,
\]
and thus—by solidity of $Y$—that $\left(1-\varepsilon\right)\cdot\left(\theta_{k}x_{k}\right)_{k\in K}\in Y$,
with
\begin{align*}
\left(1-\varepsilon\right)\cdot\left\Vert \gamma\left(x\right)\right\Vert _{Y} & =\left(1-\varepsilon\right)\cdot\left\Vert \left(\theta_{k}\left|x_{k}\right|\right)_{k\in K}\right\Vert _{Y}\\
 & \leq\left\Vert \left[\left\Vert \left(z_{i}\right)_{i\in I^{\left(k\right)}}\right\Vert _{Y_{k}}\right]_{k\in K}\right\Vert _{Y}=\left\Vert z\right\Vert _{Y\left(\left[Y_{k}\right]_{k\in K}\right)}\leq\vertiii{\iota}\cdot\left\Vert x\right\Vert _{X}<\infty.
\end{align*}
Since $\varepsilon\in\left(0,1\right)$ was arbitrary, we get $\gamma\left(x\right)\in Y$
and $\left\Vert \gamma\left(x\right)\right\Vert _{Y}\leq\vertiii{\iota}\cdot\left\Vert x\right\Vert _{X}$.
But $x\in X$ was arbitrary, so that we are done.

\medskip{}

``$\left(2\right)\!\Rightarrow\!\left(1\right)$'': Here, we do
\emph{not} assume $\left(I^{\left(k\right)}\right)_{k\in K}$ to be
pairwise disjoint. Let $x=\left(x_{i}\right)_{i\in I}\!\in\!X\!\left(\left[X_{k}\right]_{k\in K}\right)$
be arbitrary. By definition, $x^{\left(k\right)}:=\left(x_{i}\right)_{i\in I^{\left(k\right)}}\in X_{k}$
for each $k\in K$ and the sequence $y=\left(y_{k}\right)_{k\in K}$
defined by $y_{k}:=\left\Vert x^{\left(k\right)}\right\Vert _{X_{k}}$
satisfies $y\in X$ and $\left\Vert y\right\Vert _{X}=\left\Vert x\right\Vert _{X\left(\left[X_{k}\right]_{k\in K}\right)}$.
Hence—since $\gamma$ is well-defined and bounded—the sequence $z:=\gamma\left(y\right)=\left(\theta_{k}\cdot y_{k}\right)_{k\in K}$
satisfies $z\in Y$ and
\[
\left\Vert z\right\Vert _{Y}\leq\vertiii{\gamma}\cdot\left\Vert y\right\Vert _{X}=\vertiii{\gamma}\cdot\left\Vert x\right\Vert _{X\left(\left[X_{k}\right]_{k\in K}\right)}.
\]

Now, by assumption, we have $\theta_{k}<\infty$ for all $k\in K$.
Thus, for $k\in K$, there are two cases:

\begin{casenv}
\item $\delta_{k}\in X$. In this case, $\theta_{k}<\infty$ implies that
$\iota_{k}:X_{k_{0}}\hookrightarrow Y_{k_{0}}$ is well-defined and
bounded, so that we get $x^{\left(k\right)}\in Y_{k}$ with 
\[
\left\Vert \smash{x^{\left(k\right)}}\right\Vert _{Y_{k}}\leq\vertiii{\iota_{k}}\cdot\left\Vert \smash{x^{\left(k\right)}}\right\Vert _{X_{k}}=\theta_{k}\cdot y_{k}=z_{k}.
\]
\item $\delta_{k}\notin X$. In this case, solidity of $X$—together with
$y\in X$—yields $\left\Vert x^{\left(k\right)}\right\Vert _{X_{k}}=y_{k}=0$,
and hence $x^{\left(k\right)}=0$, so that trivially $x^{\left(k\right)}\in Y_{k}$
and
\[
\left\Vert \smash{x^{\left(k\right)}}\right\Vert _{Y_{k}}=0\leq z_{k}\,.
\]
\end{casenv}
Using $z\in Y$, the estimate from the case distinction, and the solidity
of $Y$, we get $\bigl(\left\Vert x^{\left(k\right)}\right\Vert _{Y_{k}}\bigr)_{k\in K}\in Y$
and
\[
\left\Vert x\right\Vert _{Y\left(\left[Y_{k}\right]_{k\in K}\right)}=\left\Vert \left(\left\Vert \smash{x^{\left(k\right)}}\right\Vert _{Y_{k}}\right)_{k\in K}\right\Vert _{Y}\leq\left\Vert z\right\Vert _{Y}\leq\vertiii{\gamma}\cdot\left\Vert x\right\Vert _{X\left(\left[X_{k}\right]_{k\in K}\right)}.
\]
Hence, $\iota$ is well-defined and bounded with $\vertiii{\iota}\leq\vertiii{\gamma}$.
\end{proof}
The most common type of solid sequence spaces that we will consider
are the weighted $\ell^{q}$ spaces $\ell_{w}^{q}$. Thus, the next
result is a useful companion to the preceding lemma.
\begin{lem}
\label{lem:EmbeddingBetweenWeightedSequenceSpaces}(cf.\@ \cite[Lemma 5.1]{DecompositionIntoSobolev})
Let $I$ be a set, let $p,q\in\left(0,\infty\right]$ and let $w=\left(w_{i}\right)_{i\in I}$
and $v=\left(v_{i}\right)_{i\in I}$ be two weights on $I$. Then,
$\iota:\ell_{w}^{q}\left(I\right)\hookrightarrow\ell_{v}^{p}\left(I\right)$
is well-defined and bounded if and only if the following expression
(then a constant) is finite:
\[
C:=\left\Vert \left(v_{i}/w_{i}\right)_{i\in I}\right\Vert _{\ell^{p\cdot\left(q/p\right)'}}\,\,.
\]
In this case, $\vertiii{\iota}=C$. Here, we use the convention
\[
p\cdot\left(q/p\right)'=\infty\text{ if }q\leq p.
\]
In the remaining case, where $p<q\leq\infty$, the expression $p\cdot\left(q/p\right)'$
is evaluated using the usual rules for computing the conjugate exponent.
\end{lem}

\begin{rem*}
We have
\begin{equation}
\frac{1}{p\cdot\left(q/p\right)'}=\left(\frac{1}{p}-\frac{1}{q}\right)_{+}.\label{eq:InverseOfSpecialExponent}
\end{equation}
In case of $p=\infty$, this is clear since the left-hand side is
$0$, while for the right-hand side, we have $p^{-1}-q^{-1}=-q^{-1}\leq0$
and hence $\left(p^{-1}-q^{-1}\right)_{+}=0$. Likewise, for $q\leq p<\infty$,
we have $p^{-1}\leq q^{-1}$, and thus $\left(p^{-1}-q^{-1}\right)_{+}=0$;
furthermore, the left-hand side of equation~(\ref{eq:InverseOfSpecialExponent})
vanishes as well.

Finally, if $p<q\leq\infty$, then
\[
\frac{1}{p\cdot\left(q/p\right)'}=\frac{1}{p}\cdot\left(1-\frac{1}{q/p}\right)=\frac{1}{p}\cdot\left(1-\frac{p}{q}\right)=\frac{1}{p}-\frac{1}{q}=\left(\frac{1}{p}-\frac{1}{q}\right)_{+}\:.
\]
In particular, equation~(\ref{eq:InverseOfSpecialExponent}) yields
\begin{equation}
p\cdot\left(q/p\right)'<\infty\:\Longleftrightarrow\:\frac{1}{p\cdot\left(q/p\right)'}>0\:\Longleftrightarrow\:\frac{1}{p}-\frac{1}{q}>0\:\Longleftrightarrow\:q>p.\qedhere\label{eq:SpecialExponentFiniteness}
\end{equation}
\end{rem*}
\begin{proof}[Proof of Lemma~\ref{lem:EmbeddingBetweenWeightedSequenceSpaces}]
``$\Leftarrow$'': Let $x=\left(x_{i}\right)_{i\in I}\in\ell_{w}^{q}\left(I\right)$
be arbitrary. Let us first assume $p<\infty$ and $p\leq q$. In this
case, we have $r:=q/p\in\left[1,\infty\right]$, so that we can apply
Hölder's inequality in the calculation
\begin{align*}
\left\Vert x\right\Vert _{\ell_{v}^{p}}^{p} & =\sum_{i\in I}\left(v_{i}\cdot\left|x_{i}\right|\right)^{p}=\sum_{i\in I}\left(\frac{v_{i}}{w_{i}}\right)^{p}\cdot\left(w_{i}\cdot\left|x_{i}\right|\right)^{p}\\
 & \leq\left\Vert \left(\left(v_{i}/w_{i}\right)^{p}\right)_{i\in I}\right\Vert _{\ell^{r'}}\cdot\left\Vert \left(\left(w_{i}\cdot\left|x_{i}\right|\right)^{p}\right)_{i\in I}\right\Vert _{\ell^{r}}\\
 & \overset{\left(\ast\right)}{=}\left\Vert \left(v_{i}/w_{i}\right)_{i\in I}\right\Vert _{\ell^{p\cdot\left(q/p\right)'}}^{p}\cdot\left\Vert \left(w_{i}\cdot x_{i}\right)_{i\in I}\right\Vert _{\ell^{q}}^{p}=C^{p}\cdot\left\Vert x\right\Vert _{\ell_{w}^{q}}^{p}<\infty.
\end{align*}
Here, the step marked with $\left(\ast\right)$ is a direct consequence
of the definition of the $\ell^{s}$-norm. A moment's thought shows
that this step is also valid in case of $\left(q/p\right)'=\infty$
or $q/p=\infty$.

Finally, taking $p$th roots of the above estimate completes the proof
for $p<\infty$ and $p\leq q$.

Next, assume $p=\infty$. In this case, we have $p\cdot\left(q/p\right)'=\infty$
and hence $v_{i}/w_{i}\leq C$ for all $i\in I$. But this implies
\[
\left|v_{i}\cdot x_{i}\right|=\frac{v_{i}}{w_{i}}\cdot\left|w_{i}\cdot x_{i}\right|\leq\frac{v_{i}}{w_{i}}\cdot\left\Vert x\right\Vert _{\ell_{w}^{q}}\leq C\cdot\left\Vert x\right\Vert _{\ell_{w}^{q}}<\infty
\]
for each $i\in I$, which easily yields the claim.

Finally, assume $q<p<\infty$. Again, $p\cdot\left(q/p\right)'=\infty$,
so that $v_{i}\leq C\cdot w_{i}$ for all $i\in I$. Note that we
have a (quasi)-norm decreasing embedding $\ell^{q}\hookrightarrow\ell^{p}$.
All in all, we get as desired that
\[
\left\Vert x\right\Vert _{\ell_{v}^{p}}\leq C\cdot\left\Vert x\right\Vert _{\ell_{w}^{p}}\leq C\cdot\left\Vert x\right\Vert _{\ell_{w}^{q}}<\infty.
\]

\medskip{}

``$\Rightarrow$'': First, let $i\in I$ be arbitrary. We have $v_{i}=\left\Vert \delta_{i}\right\Vert _{\ell_{v}^{p}}\leq\vertiii{\iota}\cdot\left\Vert \delta_{i}\right\Vert _{\ell_{w}^{q}}=\vertiii{\iota}\cdot w_{i}$,
and hence $v_{i}/w_{i}\leq\vertiii{\iota}$. This proves the claim
in all cases where $p\cdot\left(q/p\right)'=\infty$.

Thus, we can assume in the following that $p<q\leq\infty$. In this
case, we have $r:=q/p\in\left(1,\infty\right]$. Let $\theta_{i}:=\left(v_{i}/w_{i}\right)^{p}$
for $i\in I$ and note $p\cdot\left(q/p\right)'=p\cdot r'$, so that
we get
\[
C=\left\Vert \left(v_{i}/w_{i}\right)_{i\in I}\right\Vert _{\ell^{p\cdot\left(q/p\right)'}}=\left\Vert \theta\right\Vert _{\ell^{r'}}^{1/p}.
\]
It is thus sufficient to show $\left\Vert \theta\right\Vert _{\ell^{r'}}\leq\vertiii{\iota}^{p}$.
But since $r,r'\in\left[1,\infty\right]$, it is well-known that
\[
\left\Vert \theta\right\Vert _{\ell^{r'}}=\sup\left\{ \sum_{i\in I}\left|x_{i}\theta_{i}\right|\with x=\left(x_{i}\right)_{i\in I}\in\ell_{0}\left(I\right)\text{ and }\left\Vert x\right\Vert _{\ell^{r}}\leq1\right\} ,
\]
where $\ell_{0}\left(I\right)$ denotes the space of finitely supported
sequences over $I$.

Thus, let $x=\left(x_{i}\right)_{i\in I}\in\ell_{0}\left(I\right)$
with $\left\Vert x\right\Vert _{\ell^{r}}\leq1$ be arbitrary. Define
$y_{i}:=w_{i}^{-1}\cdot\left|x_{i}\right|^{1/p}$ for $i\in I$ and
note 
\[
\left\Vert y\right\Vert _{\ell_{w}^{q}}=\left\Vert \left(\vphantom{x^{i}}\smash{\left|x_{i}\right|^{1/p}}\right)_{i\in I}\right\Vert _{\ell^{q}}=\left\Vert x\right\Vert _{\ell^{q/p}}^{1/p}=\left\Vert x\right\Vert _{\ell^{r}}^{1/p}\leq1.
\]
Since the embedding $\iota$ is well-defined and bounded and because
of $p<\infty$, this implies 
\[
\vphantom{\sum_{i\in I}}\left[\,\smash{\sum_{i\in I}}\,\vphantom{\sum}\left|x_{i}\theta_{i}\right|\,\right]^{1/p}=\left[\,\smash{\sum_{i\in I}}\,\vphantom{\sum}\left(v_{i}/w_{i}\right)^{p}\cdot\left|x_{i}\right|\,\right]^{1/p}=\left\Vert \frac{v_{i}}{w_{i}}\cdot\left|x_{i}\right|^{1/p}\right\Vert _{\ell^{p}}=\left\Vert y\right\Vert _{\ell_{v}^{p}}\leq\vertiii{\iota}.
\]
All in all, we have shown $\sum_{i\in I}\left|x_{i}\theta_{i}\right|\leq\vertiii{\iota}^{p}$
for every finitely supported sequence $x=\left(x_{i}\right)_{i\in I}\in\ell_{0}\left(I\right)$
with $\left\Vert x\right\Vert _{\ell^{r}}\leq1$. As seen above, this
yields $\left\Vert \theta\right\Vert _{\ell^{r'}}\leq\vertiii{\iota}^{p}$,
which then implies the claim.
\end{proof}
The preceding results yield a satisfactory characterization of the
existence of an embedding between two nested sequence spaces if the
underlying sets $\left(I^{\left(k\right)}\right)_{k\in K}$ are disjoint
(and the same on both sides of the embedding). It is thus of interest
to have a criterion which allows to reduce matters to this case. Such
a criterion is given in the next lemma. Although the given result
might appear rather technical, we will see later that it is useful
and readily applicable in a number of situations; see Lemma~\ref{lem:SubordinatenessEnablesDisjointization}
and Corollaries \ref{cor:EmbeddingFineInCoarseSimplification} and
\ref{cor:EmbeddingCoarseIntoFineSimplification}.
\begin{lem}
\label{lem:NestedEmbeddingReductionToDisjointSets}Let $K$ be an
index set, and for each $k\in K$, let $I^{\left(k\right)}$ be a
set, and let $u_{k}=\left(u_{k,i}\right)_{i\in I^{\left(k\right)}}$
be a weight on $I^{\left(k\right)}$. Let $q\in\left(0,\infty\right]$.

Set $I:=\bigcup_{k\in K}I^{\left(k\right)}$ and define a relation
$\sim$ on $K$ by
\[
k\sim\ell\qquad:\Longleftrightarrow\qquad I^{\left(k\right)}\cap I^{\left(\ell\right)}\neq\emptyset.
\]
Let $X\subset\Compl^{K}$ be a solid sequence space for which the
\textbf{generalized clustering map}
\begin{equation}
\Theta:X\to X,x=\left(x_{k}\right)_{k\in K}\mapsto\left[\,\sum_{\ell\in K\text{ with }\ell\sim k}x_{\ell}\right]_{k\in K}\label{eq:ReductionToDisjointSetsGeneralizedClusteringMap}
\end{equation}
is well-defined and bounded. Furthermore, assume that we have a uniform
bound
\begin{equation}
N:=\sup_{k\in K}\left|\left[k\right]\right|<\infty\qquad\text{ where }\qquad\left[k\right]:=\left\{ \ell\in K\with\ell\sim k\right\} \,,\label{eq:ReductionToDisjointSetsCardinalityAssumption}
\end{equation}
and that there is some constant $C_{u}>0$ with $u_{k,i}\leq C_{u}\cdot u_{\ell,i}$
for all $k,\ell\in K$ and $i\in I^{\left(k\right)}\cap I^{\left(\ell\right)}$.
Finally, assume that for each $k\in K$, we are given a subset $I^{\left(k,\star\right)}\subset I^{\left(k\right)}$
with $I=\bigcup_{k\in K}I^{\left(k,\star\right)}$.

Then
\[
\left\Vert x\right\Vert _{X\left(\smash{\left[\ell_{u}^{q}\left(\smash{I^{\left(k\right)}}\right)\right]_{k\in K}}\right)}\geq\left\Vert x\right\Vert _{X\left(\smash{\left[\ell_{u}^{q}\left(\smash{I^{\left(k,\star\right)}}\right)\right]_{k\in K}}\right)}\geq C^{-1}\cdot\left\Vert x\right\Vert _{X\left(\smash{\left[\ell_{u}^{q}\left(\smash{I^{\left(k\right)}}\right)\right]_{k\in K}}\right)}
\]
for all sequences $x\in\Compl^{I}$ and some (fixed) constant $C=C\left(\vertiii{\Theta},C_{u},N,q\right)$.

In particular, for every $x\in\Compl^{I}$, we have $\vphantom{\ell^{I^{k}}}x\in X\bigl(\smash{\left[\ell_{u}^{q}\left(\smash{I^{\left(k\right)}}\right)\right]_{k\in K}}\bigr)\Longleftrightarrow x\in X\bigl(\smash{\left[\ell_{u}^{q}\left(\smash{I^{\left(k,\star\right)}}\right)\right]_{k\in K}}\bigr)$.
\end{lem}

\begin{proof}
Since $X$ is solid, the estimate $\left\Vert x\right\Vert _{X\left(\smash{\left[\ell_{u}^{q}\left(\smash{I^{\left(k\right)}}\right)\right]_{k\in K}}\right)}\geq\left\Vert x\right\Vert _{X\left(\smash{\left[\ell_{u}^{q}\left(\smash{I^{\left(k,\star\right)}}\right)\right]_{k\in K}}\right)}$
is a direct consequence of the estimate $\left\Vert y\right\Vert _{\ell_{u}^{q}\left(\smash{I^{\left(k,\star\right)}}\right)}\leq\left\Vert y\right\Vert _{\ell_{u}^{q}\left(\smash{I^{\left(k\right)}}\right)}$
for all $y\in\Compl^{I^{\left(k\right)}}$, which in turn follows
from $I^{\left(k,\ast\right)}\subset I^{\left(k\right)}$.

It is thus sufficient to show
\[
\left\Vert x\right\Vert _{X\left(\smash{\left[\ell_{u}^{q}\left(\smash{I^{\left(k\right)}}\right)\right]_{k\in K}}\right)}\leq C\cdot\left\Vert x\right\Vert _{X\left(\smash{\left[\ell_{u}^{q}\left(\smash{I^{\left(k,\star\right)}}\right)\right]_{k\in K}}\right)}\qquad\text{assuming}\qquad x\in X\bigl(\smash{\left[\ell_{u}^{q}\left(\smash{I^{\left(k,\star\right)}}\right)\right]_{k\in K}}\bigr).
\]
To this end, let $k\in K$ be arbitrary. For each $i\in I^{\left(k\right)}\subset I$,
there is—because of $I=\bigcup_{\ell\in K}I^{\left(\ell,\star\right)}$—some
$\ell_{i}\in K$ with $i\in I^{\left(\ell_{i},\star\right)}$. Note
that $i\in I^{\left(k\right)}\cap I^{\left(\ell_{i},\star\right)}\subset I^{\left(k\right)}\cap I^{\left(\ell_{i}\right)}$,
so that $\ell_{i}\in\left[k\right]$. Furthermore, $u_{k,i}\leq C_{u}\cdot u_{\ell,i}$
as long as $i\in I^{\left(k\right)}\cap I^{\left(\ell,\star\right)}\subset I^{\left(k\right)}\cap I^{\left(\ell\right)}$.
All in all, we get—in case of $q<\infty$—that
\begin{align*}
\left\Vert x\right\Vert _{\ell_{u}^{q}\left(\smash{I^{\left(k\right)}}\right)}^{q}=\sum_{i\in I^{\left(k\right)}}\left(u_{k,i}\cdot\left|x_{i}\right|\right)^{q} & \leq\sum_{i\in I^{\left(k\right)}}\left(\,u_{k,i}\vphantom{\sum}\smash{\cdot\sum_{\ell\in\left[k\right]}}\left[\Indicator_{I^{\left(\ell,\star\right)}}\left(i\right)\cdot\left|x_{i}\right|\right]\,\right)^{q}\\
 & \leq C_{u}^{q}\cdot\sum_{i\in I^{\left(k\right)}}\left(\,\vphantom{\sum}\smash{\sum_{\ell\in\left[k\right]}}\left[\Indicator_{I^{\left(\ell,\star\right)}}\left(i\right)\cdot u_{\ell,i}\cdot\left|x_{i}\right|\right]\,\right)^{q}.
\end{align*}

Now, note the general estimate
\begin{equation}
\left(\vphantom{\sum}\,\smash{\sum_{j\in J}}\alpha_{j}\,\right)^{q}\leq\Bigl(\left|J\right|\cdot\max\left\{ \alpha_{j}\with j\in J\right\} \Bigr)^{q}\leq\left|J\right|^{q}\cdot\sum_{j\in J}\alpha_{j}^{q}\label{eq:DragPowerThroughSumCheapEstimate}
\end{equation}
for finite sets $J$ and arbitrary sequences $\left(\alpha_{j}\right)_{j\in J}$
of nonnegative numbers. Applied to the above estimate, we conclude
\begin{align*}
\left\Vert x\right\Vert _{\ell_{u}^{q}\left(\smash{I^{\left(k\right)}}\right)}^{q} & \leq C_{u}^{q}\cdot\left|\left[k\right]\right|^{q}\cdot\sum_{i\in I^{\left(k\right)}}\;\sum_{\ell\in\left[k\right]}\bigl(\Indicator_{I^{\left(\ell,\star\right)}}\left(i\right)\cdot u_{\ell,i}\cdot\left|x_{i}\right|\bigr)^{q}\\
\left({\scriptstyle \text{since }\left|\left[k\right]\right|\leq N}\right) & \leq\left(C_{u}N\right)^{q}\cdot\sum_{\ell\in\left[k\right]}\;\sum_{i\in I^{\left(k\right)}}\bigl(\Indicator_{I^{\left(\ell,\star\right)}}\left(i\right)\cdot u_{\ell,i}\cdot\left|x_{i}\right|\bigr)^{q}\\
 & \leq\left(C_{u}N\right)^{q}\cdot\sum_{\ell\in\left[k\right]}\;\sum_{i\in I^{\left(\ell,\star\right)}}\left(u_{\ell,i}\cdot\left|x_{i}\right|\right)^{q}=\left(C_{u}N\right)^{q}\cdot\sum_{\ell\in\left[k\right]}\left\Vert x\right\Vert _{\ell_{u}^{q}\left(\smash{I^{\left(\ell,\star\right)}}\right)}^{q}\\
\left({\scriptstyle \text{since }\left|\left[k\right]\right|\leq N}\right) & \leq C_{u}^{q}N^{1+q}\cdot\left(\,\vphantom{\sum}\smash{\sum_{\ell\in\left[k\right]}}\left\Vert x\right\Vert _{\ell_{u}^{q}\left(I^{\left(\ell,\star\right)}\right)}\,\right)^{q}=C_{u}^{q}N^{1+q}\cdot\left[\left(\Theta\,y\right)_{k}\right]^{q},
\end{align*}
where we defined
\[
y_{\ell}:=\left\Vert x\right\Vert _{\ell_{u}^{q}\left(\smash{I^{\left(\ell,\star\right)}}\right)}\qquad\text{ for }\ell\in K.
\]

In case of $q=\infty$, we argue similarly: For arbitrary $i\in I^{\left(k\right)}$,
we saw above that $i\in I^{\left(\ell_{i},\star\right)}$ for some
$\ell_{i}\in\left[k\right]$. Furthermore, $u_{k,i}\leq C_{u}\cdot u_{\ell_{i},i}$
by our assumptions. Hence,
\begin{align*}
u_{k,i}\cdot\left|x_{i}\right|\leq C_{u}\cdot\Indicator_{I^{\left(\ell_{i},\star\right)}}\left(i\right)\cdot u_{\ell_{i},i}\cdot\left|x_{i}\right| & \leq C_{u}\cdot\sum_{\ell\in\left[k\right]}\Indicator_{I^{\left(\ell,\star\right)}}\left(i\right)\cdot u_{\ell,i}\cdot\left|x_{i}\right|\\
\left({\scriptstyle \text{since }q=\infty}\right) & \leq C_{u}\cdot\sum_{\ell\in\left[k\right]}\left\Vert x\right\Vert _{\ell_{u}^{q}\left(\smash{I^{\left(\ell,\star\right)}}\right)}=C_{u}\cdot\left(\Theta\,y\right)_{k}\,,
\end{align*}
with $y=\left(y_{\ell}\right)_{\ell\in K}$ as above. Since $i\in I^{\left(k\right)}$
was arbitrary, and since $q=\infty$, we conclude 
\[
\left\Vert x\right\Vert _{\ell_{u}^{q}\left(\smash{I^{\left(k\right)}}\right)}\leq C_{u}\cdot\left(\Theta\,y\right)_{k}\leq C_{u}N^{1+\frac{1}{q}}\cdot\left(\Theta\,y\right)_{k}.
\]

In summary, we have shown $0\leq\left\Vert x\right\Vert _{\ell_{u}^{q}\left(\smash{I^{\left(k\right)}}\right)}\leq C_{u}N^{1+\frac{1}{q}}\cdot\left(\Theta y\right)_{k}$
for all $k\in K$ and arbitrary $q\in\left(0,\infty\right]$. But
because of $x\in X\bigl(\left[\ell_{u}^{q}\left(\smash{I^{\left(k,\star\right)}}\right)\right]_{k\in K}\bigr)$,
we have $y\in X$ with $\left\Vert y\right\Vert _{X}=\left\Vert x\right\Vert _{X\left(\smash{\left[\ell_{u}^{q}\left(\smash{I^{\left(k,\star\right)}}\right)\right]_{k\in K}}\right)}$.
By solidity of $X$ and since we assume $\Theta$ to be well-defined
and bounded, we conclude $\vphantom{I^{I^{k}}}x\in X\bigl(\left[\ell_{u}^{q}\left(\smash{I^{\left(k\right)}}\right)\right]_{k\in K}\bigr)$
with
\begin{align*}
\left\Vert x\right\Vert _{X\left(\smash{\left[\ell_{u}^{q}\left(\smash{I^{\left(k\right)}}\right)\right]_{k\in K}}\right)}=\left\Vert \left(\left\Vert x\right\Vert _{\ell_{u}^{q}\left(\smash{I^{\left(k\right)}}\right)}\right)_{k\in K}\right\Vert _{X} & \leq C_{u}N^{1+\frac{1}{q}}\cdot\left\Vert \Theta\,y\right\Vert _{X}\leq C_{u}N^{1+\frac{1}{q}}\vertiii{\Theta}\cdot\left\Vert y\right\Vert _{X}\\
 & =C_{u}N^{1+\frac{1}{q}}\vertiii{\Theta}\cdot\left\Vert x\right\Vert _{X\left(\smash{\left[\ell_{u}^{q}\left(\smash{I^{\left(k,\star\right)}}\right)\right]_{k\in K}}\right)}<\infty.\qedhere
\end{align*}
\end{proof}
The main point of the preceding lemma is that it allows to switch
from the space $X\bigl(\left[\ell_{u}^{q}\left(\smash{I^{\left(k\right)}}\right)\right]_{k\in K}\bigr)$
to the slightly modified space $X\bigl(\left[\ell_{u}^{q}\left(\smash{I^{\left(k,\star\right)}}\right)\right]_{k\in K}\bigr)$.
In view of Lemma~\ref{lem:NestedDisjointEmbedding}, it is preferable
if the family $\left(\smash{I^{\left(k,\star\right)}}\right)_{k\in K}$
can be chosen to be disjoint. Thus, the following lemma is helpful,
since it shows that one can always choose such a disjoint family $\left(\smash{I^{\left(k,\star\right)}}\right)_{k\in K}$
which still satisfies $I=\bigcup_{k\in K}I^{\left(k,\star\right)}$.
\begin{lem}
\label{lem:NestedEmbeddingReductionExistenceOfDisjointifiedSets}Under
the assumptions of Lemma~\ref{lem:NestedEmbeddingReductionToDisjointSets},
the following are true:

\begin{enumerate}[leftmargin=0.7cm]
\item For each pairwise disjoint family of sets $\left(I^{\left(k,0\right)}\right)_{k\in K}$
with $I^{\left(k,0\right)}\subset I^{\left(k\right)}$, there is a
pairwise disjoint family $\left(I^{\left(k,\star\right)}\right)_{k\in K}$
satisfying $I^{\left(k,0\right)}\subset I^{\left(k,\star\right)}\subset I^{\left(k\right)}$
for all $k\in K$ and $I=\biguplus_{k\in K}I^{\left(k,\star\right)}$.
\item With $N$ as in equation~(\ref{eq:ReductionToDisjointSetsCardinalityAssumption}),
there is a family $\left(I^{\left(k,n\right)}\right)_{k\in K,n\in\underline{N}}$
with the following properties:

\begin{enumerate}
\item $I^{\left(k\right)}=\bigcup_{n=1}^{N}I^{\left(k,n\right)}$ for all
$k\in K$,
\item For each $n\in\underline{N}$, we have $I=\biguplus_{k\in K}I^{\left(k,n\right)}$.\qedhere
\end{enumerate}
\end{enumerate}
\end{lem}

\begin{proof}
Ad (1): Set
\[
\mathscr{A}:=\left\{ \left(\smash{I_{k}}\right)_{k\in K}\with\left(\smash{I_{k}}\right)_{k\in K}\text{ is pairwise disjoint and }I^{\left(k,0\right)}\subset I_{k}\subset I^{\left(k\right)}\text{ for all }k\in K\right\} .
\]
Note that $\left(I^{\left(k,0\right)}\right)_{k\in K}\in\mathscr{A}$,
so that $\mathscr{A}\neq\emptyset$. Now, define a partial order on
$\mathscr{A}$ by setting
\[
\left(\smash{I_{k}}\right)_{k\in K}\leq\left(\smash{J_{k}}\right)_{k\in K}\quad:\Longleftrightarrow\quad\forall\,k\in K:\:I_{k}\subset J_{k}.
\]

It is not hard to see that $\mathscr{A}$ is inductively ordered by
``$\leq$''; indeed, if $\left(\left(\smash{I_{k}^{\left(\alpha\right)}}\right)_{k\in K}\right)_{\alpha\in A}$
is a (nonempty) chain in $\mathscr{A}$, then $I_{k}:=\bigcup_{\alpha\in A}I_{k}^{\left(\alpha\right)}$
satisfies $I^{\left(k,0\right)}\subset I_{k}^{\left(\alpha\right)}\subset I_{k}\subset I^{\left(k\right)}$
for all $k\in K$ and arbitrary $\alpha\in A$. Furthermore, $\left(I_{k}\right)_{k\in K}$
is pairwise disjoint (and hence $\left(I_{k}\right)_{k\in K}\in\mathscr{A}$),
since otherwise there are $k,\ell\in K$ with $k\neq\ell$ and some
$x\in I_{k}\cap I_{\ell}$. By definition, there are thus $\alpha_{1},\alpha_{2}\in A$
with $x\in I_{k}^{\left(\alpha_{1}\right)}\cap I_{\ell}^{\left(\alpha_{2}\right)}$.
But since we are considering a chain, we have $I_{k}^{\left(\alpha_{1}\right)}\subset I_{k}^{\left(\alpha_{2}\right)}$
and hence $x\in I_{k}^{\left(\alpha_{2}\right)}\cap I_{\ell}^{\left(\alpha_{2}\right)}=\emptyset$
or $I_{\ell}^{\left(\alpha_{2}\right)}\subset I_{\ell}^{\left(\alpha_{1}\right)}$
and hence $x\in I_{k}^{\left(\alpha_{1}\right)}\cap I_{\ell}^{\left(\alpha_{1}\right)}=\emptyset$.
Both of these cases are absurd. Hence, $\left(I_{k}\right)_{k\in K}\in\mathscr{A}$
is an upper bound for the given chain.

By Zorn's lemma, $\mathscr{A}$ has a maximal element $\left(I^{\left(k,\star\right)}\right)_{k\in K}$.
It remains to show $I=\bigcup_{k\in K}I^{\left(k,\star\right)}$.
If this is false, there is some $i\in I=\bigcup_{k\in K}I^{\left(k\right)}$
satisfying $i\notin\bigcup_{k\in K}I^{\left(k,\star\right)}$. Choose
$k_{0}\in K$ with $i\in I^{\left(k_{0}\right)}$ and define
\[
J^{\left(k,\star\right)}:=\begin{cases}
I^{\left(k,\star\right)}, & \text{if }k\neq k_{0},\\
I^{\left(k,\star\right)}\cup\left\{ i\right\} , & \text{if }k=k_{0}.
\end{cases}
\]
It is then not hard to see $\left(J^{\left(k\right)}\right)_{k\in K}\in\mathscr{A}$
and $\left(I^{\left(k,\star\right)}\right)_{k\in K}\leq\left(J^{\left(k,\star\right)}\right)_{k\in K}\not\leq\left(I^{\left(k,\star\right)}\right)_{k\in K}$,
contradicting maximality of $\left(I^{\left(k,\star\right)}\right)_{k\in K}$.

\medskip{}

Ad (2): Let $K_{0}:=\left\{ k\in K\with I^{\left(k\right)}\neq\emptyset\right\} $.
The assumption of Lemma~\ref{lem:NestedEmbeddingReductionToDisjointSets}
implies that we have
\[
\left|\left[k\right]\right|\leq N\quad\text{and}\quad\left[k\right]=\left\{ \ell\in K_{0}\with I^{\left(k\right)}\cap I^{\left(\ell\right)}\neq\emptyset\right\} \qquad\forall\,k\in K_{0}.
\]
Furthermore, the relation $k\sim\ell\,:\Longleftrightarrow\,I^{\left(k\right)}\cap I^{\left(\ell\right)}\neq\emptyset$
is easily seen to be reflexive and symmetric on $K_{0}$. Thus, by
Lemma~\ref{lem:RelationDisjointification}, there is a finite partition
$K_{0}=\biguplus_{n=1}^{N}K_{0}^{\left(n\right)}$ such that $k\not\sim\ell$
for arbitrary $n\in\underline{N}$ and $k,\ell\in K_{0}^{\left(n\right)}$with
$k\neq\ell$.

Now, define
\[
I_{0}^{\left(k,n\right)}:=\begin{cases}
I^{\left(k\right)}, & \text{if }k\in K_{0}^{\left(n\right)},\\
\emptyset, & \text{if }k\notin K_{0}^{\left(n\right)}
\end{cases}\qquad\qquad\text{ for }k\in K\text{ and }n\in\underline{N}\,.
\]
We have $I_{0}^{\left(k,n\right)}\subset I^{\left(k\right)}$ for
all $k\in K$. Furthermore, for $k,\ell\in K$ with $k\neq\ell$,
we have $I_{0}^{\left(k,n\right)}\cap I_{0}^{\left(\ell,n\right)}=\emptyset$,
unless possibly if $k,\ell\in K_{0}^{\left(n\right)}$. But in the
latter case, we have $k\not\sim\ell$ by choice of $K_{0}^{\left(n\right)}$
and thus $I_{0}^{\left(k,n\right)}\cap I_{0}^{\left(\ell,n\right)}=I^{\left(k\right)}\cap I^{\left(\ell\right)}=\emptyset$.
Thus, for each $n\in\underline{N}$, the family $\left(\smash{I_{0}^{\left(k,n\right)}}\right)_{k\in K}$
is pairwise disjoint.

By the first part, there is thus for each $n\in\underline{N}$ a pairwise
disjoint family $\left(I^{\left(k,n\right)}\right)_{k\in K}$ satisfying
$I_{0}^{\left(k,n\right)}\subset I^{\left(k,n\right)}\subset I^{\left(k\right)}$
and $I=\biguplus_{k\in K}I^{\left(k,n\right)}$. It remains to check
$I^{\left(k\right)}\subset\bigcup_{n=1}^{N}I^{\left(k,n\right)}$
for each $k\in K$. For $k\in K\setminus K_{0}$, this is clear since
$I^{\left(k\right)}=\emptyset$. But for $k\in K_{0}$, we have $k\in K_{0}^{\left(n_{k}\right)}$
for some $n_{k}\in\underline{N}$ and hence $\bigcup_{n=1}^{N}I^{\left(k,n\right)}\supset I^{\left(k,n_{k}\right)}\supset I_{0}^{\left(k,n_{k}\right)}=I^{\left(k\right)}$,
as desired.
\end{proof}
Our next result indicates an important case in which the assumptions
of Lemma~\ref{lem:NestedEmbeddingReductionToDisjointSets} are satisfied:
\begin{lem}
\label{lem:SubordinatenessEnablesDisjointization}Let $\CalQ=\left(Q_{i}\right)_{i\in I}$
be a family of subsets of $\R^{\dimension}$ and assume that $\CalR=\left(R_{k}\right)_{k\in K}$
is an admissible covering of a set $X\subset\R^{\dimension}$.

Furthermore, assume that $\CalQ$ is almost subordinate to $\CalR$
and pick some $n\in\N_{0}$. For each $k\in K$, choose a set
\[
I^{\left(k\right)}\subset\left\{ i\in I\with Q_{i}\cap R_{k}^{n\ast}\neq\emptyset\right\} .
\]
If we define a relation $\sim$ on $K$ by setting $k\sim\ell\;:\Longleftrightarrow\;I^{\left(k\right)}\cap I^{\left(\ell\right)}\neq\emptyset$,
then $\left[k\right]:=\left\{ \ell\in K\with\ell\sim k\right\} $
satisfies
\[
\left[k\right]\subset k^{\left(2n+4k\left(\CalQ,\CalR\right)+5\right)\ast}\quad\text{and}\quad\left|\left[k\right]\right|\leq N_{\CalR}^{2n+4k\left(\CalQ,\CalR\right)+5}\qquad\qquad\forall\,k\in K.\qedhere
\]
\end{lem}

\begin{proof}
Fix $k\in K$ and let $\ell\in\left[k\right]$. Hence, there is some
$i\in I^{\left(k\right)}\cap I^{\left(\ell\right)}$, i.e.\@ with
$Q_{i}\cap R_{k}^{n\ast}\neq\emptyset\neq Q_{i}\cap R_{\ell}^{n\ast}$.
This yields $k_{0}\in k^{n\ast}$ and $\ell_{0}\in\ell^{n\ast}$ with
$Q_{i}\cap R_{k_{0}}\neq\emptyset\neq Q_{i}\cap R_{\ell_{0}}$. But
in view of Lemma~\ref{lem:SubordinatenessImpliesWeakSubordination},
this yields 
\[
\emptyset\neq Q_{i}\subset R_{k_{0}}^{\left(2k\left(\CalQ,\CalR\right)+2\right)\ast}\cap R_{\ell_{0}}^{\left(2k\left(\CalQ,\CalR\right)+2\right)\ast}.
\]
In particular, $\ell_{0}\in k_{0}^{\left(4k\left(\CalQ,\CalR\right)+5\right)\ast}\subset k^{\left(n+4k\left(\CalQ,\CalR\right)+5\right)\ast}$
and thus finally $\ell\in\ell_{0}^{n\ast}\subset k^{\left(2n+4k\left(\CalQ,\CalR\right)+5\right)\ast}$.

All in all, we have shown $\left[k\right]\subset k^{\left(2n+4k\left(\CalQ,\CalR\right)+5\right)\ast}$.
The estimate of $\left|\left[k\right]\right|$ is now a direct consequence
of Lemma~\ref{lem:SemiStructuredClusterInvariant}.
\end{proof}
The following is an important application of Lemmas \ref{lem:NestedDisjointEmbedding}–\ref{lem:NestedEmbeddingReductionExistenceOfDisjointifiedSets}.
It makes crucial use of the identity $\ell_{v}^{r}\left(I\right)=\ell^{r}\bigl(\left[\ell_{v}^{r}\left(\smash{I^{\left(k\right)}}\right)\right]_{k\in K}\bigr)$,
which is valid for \emph{every} partition $I=\biguplus_{k\in K}I^{\left(k\right)}$.
\begin{cor}
\label{cor:AlmostDisjointifiedEmbeddingLebesgue}Under the assumptions
of Lemma~\ref{lem:NestedEmbeddingReductionToDisjointSets}, choose
a family $\left(\smash{I^{\left(k,\natural\right)}}\right)_{k\in K}$
with $I^{\left(k,\natural\right)}\subset I^{\left(k\right)}$ for
all $k\in K$ and with $I=\bigcup_{k\in K}I^{\left(k,\natural\right)}$.
Then the following are true:

\begin{enumerate}[leftmargin=0.7cm]
\item \label{enu:AlmostDisjointifiedNestedIntoPure}For a given weight
$v=\left(v_{i}\right)_{i\in I}$ and $r\in\left(0,\infty\right]$,
the embedding
\[
\iota:X\left(\left[\ell_{u}^{q}\left(\smash{I^{\left(k\right)}}\right)\right]_{k\in K}\right)\hookrightarrow\ell_{v}^{r}\left(I\right)
\]
satisfies
\[
C^{-1}\cdot\vertiii{\gamma}\leq\vertiii{\iota}\leq\vertiii{\smash{\gamma^{\natural}}}\leq\vertiii{\gamma},
\]
for some constant $C=C\left(\vertiii{\Theta},N,C_{u},q,r\right)>0$,
where $\gamma,\gamma^{\natural}$ are given by
\begin{alignat*}{2}
\qquad\gamma: & X\to\ell^{r}\left(K\right), & \left(x_{k}\right)_{k\in K}\mapsto\left(\theta_{k}\cdot x_{k}\right)_{k\in k} & \text{ with }\theta_{k}:=\begin{cases}
\left\Vert \left(v_{i}/u_{k,i}\right)_{i\in I^{\left(k\right)}}\right\Vert _{\ell^{r\cdot\left(q/r\right)'}}, & \text{if }\delta_{k}\in X,\\
1, & \text{else},
\end{cases}\\
\qquad\gamma^{\natural}: & X\to\ell^{r}\left(K\right), & \left(x_{k}\right)_{k\in K}\mapsto\left(\smash{\theta_{k}^{\natural}}\cdot x_{k}\right)_{k\in k} & \text{ with }\theta_{k}^{\natural}:=\begin{cases}
\left\Vert \left(v_{i}/u_{k,i}\right)_{i\in I^{\left(k,\natural\right)}}\right\Vert _{\ell^{r\cdot\left(q/r\right)'}}, & \text{if }\delta_{k}\in X,\\
1, & \text{else}.
\end{cases}
\end{alignat*}
In particular, $\iota$ is well-defined and bounded if and only
if $\theta_{k}<\infty$ for all $k\in K$ and if $\gamma$ is well-defined
and bounded.
\item \label{enu:AlmostDisjointifiedPureIntoNested}Let $C_{X}\geq1$ be
a triangle constant for $X$. For a given weight $v=\left(v_{i}\right)_{i\in I}$
and $r\in\left(0,\infty\right]$, the embedding
\[
\iota:\ell_{v}^{r}\left(I\right)\hookrightarrow X\left(\left[\ell_{u}^{q}\left(\smash{I^{\left(k\right)}}\right)\right]_{k\in K}\right)
\]
satisfies
\[
C_{1}^{-1}\cdot\vertiii{\eta}\leq\vertiii{\iota}\leq C_{2}\cdot\vertiii{\eta^{\natural}}\leq C_{2}\cdot\vertiii{\eta}
\]
for certain constants $C_{1}=C_{1}\left(N,C_{X},q,r\right)>0$, $C_{2}=C_{2}\left(\vertiii{\Theta},C_{u},N,q\right)>0$
and
\begin{alignat*}{2}
\eta: & \ell^{r}\left(K\right)\to X, & \left(x_{k}\right)_{k\in K}\mapsto\left(\vartheta_{k}\cdot x_{k}\right)_{k\in K} & \quad\text{with}\quad\vartheta_{k}:=\left\Vert \left(u_{k,i}/v_{i}\right)_{i\in I^{\left(k\right)}}\right\Vert _{\ell^{q\cdot\left(r/q\right)'}}\\
\eta^{\natural}: & \ell^{r}\left(K\right)\to X, & \left(x_{k}\right)_{k\in K}\mapsto\left(\smash{\vartheta_{k}^{\natural}}\cdot x_{k}\right)_{k\in K} & \quad\text{with}\quad\vartheta_{k}^{\natural}:=\left\Vert \left(u_{k,i}/v_{i}\right)_{i\in I^{\left(k,\natural\right)}}\right\Vert _{\ell^{q\cdot\left(r/q\right)'}}\quad.
\end{alignat*}
In particular, $\iota$ is well-defined and bounded if and only if
$\vartheta_{k}<\infty$ for all $k\in K$ and if $\eta$ is well-defined
and bounded.
\item \label{enu:AlmostDisjointifiedLebesgueSimplified}In case of $X=\ell_{w}^{s}\left(K\right)$
for some weight $w=\left(w_{k}\right)_{k\in K}$ and some $s\in\left(0,\infty\right]$,
we have
\[
\vertiii{\gamma^{\natural}}=\left\Vert \left(w_{k}^{-1}\cdot\left\Vert \left(v_{i}/u_{k,i}\right)_{i\in I^{\left(k,\natural\right)}}\right\Vert _{\ell^{r\cdot\left(q/r\right)'}}\right)_{k\in K}\right\Vert _{\ell^{r\cdot\left(s/r\right)'}}
\]
and
\[
\vertiii{\eta^{\natural}}=\left\Vert \left(w_{k}\cdot\left\Vert \left(u_{k,i}/v_{i}\right)_{i\in I^{\left(k,\natural\right)}}\right\Vert _{\ell^{q\cdot\left(r/q\right)'}}\right)_{k\in K}\right\Vert _{\ell^{s\cdot\left(r/s\right)'}}\;.
\]
In particular, $\gamma^{\natural}$ or $\eta^{\natural}$ is well-defined
and bounded if and only if the respective right-hand side is finite.\qedhere
\end{enumerate}
\end{cor}

\begin{proof}
In case of $I=\emptyset$, all claims are trivially satisfied. Hence,
we can assume $I\neq\emptyset$, so that the constant $N$ defined
in equation~(\ref{eq:ReductionToDisjointSetsCardinalityAssumption})
satisfies $N>0$, i.e., $N\in\N$. We now consider each part of the
corollary separately.

Ad (1): ``$\Rightarrow$'': First assume that $\iota$ is well-defined
and bounded. By Lemma~\ref{lem:NestedEmbeddingReductionExistenceOfDisjointifiedSets},
there is for each $n\in\underline{N}$ a partition $I=\biguplus_{k\in K}I^{\left(k,n\right)}$
such that we have $I^{\left(k\right)}=\bigcup_{n=1}^{N}I^{\left(k,n\right)}$
for all $k\in K$. Now, Lemma~\ref{lem:NestedEmbeddingReductionToDisjointSets}
yields boundedness of $\iota^{\left(n\right)}:X\left(\left[\ell_{u}^{q}\left(\smash{I^{\left(k,n\right)}}\right)\right]_{k\in K}\right)\hookrightarrow X\left(\left[\ell_{u}^{q}\left(\smash{I^{\left(k\right)}}\right)\right]_{k\in K}\right)$
for each $n\in\underline{N}$, and furthermore $\vertiii{\iota^{\left(n\right)}}\leq C=C\left(\vertiii{\Theta},C_{u},N,q\right)$.

Next, note that $I=\biguplus_{k\in K}I^{\left(k,n\right)}$ implies
$\ell_{v}^{r}\left(I\right)=\ell^{r}\bigl(\left[\ell_{v}^{r}\left(\smash{I^{\left(k,n\right)}}\right)\right]_{k\in K}\bigr)$
with equal (quasi)-norms and $v_{k,i}=v_{i}$ for $k\in K$ and $i\in I^{\left(k,n\right)}$.
Thus, by composition, we get boundedness of
\[
\iota\circ\iota^{\left(n\right)}:X\bigl(\left[\ell_{u}^{q}\left(\smash{I^{\left(k,n\right)}}\right)\right]_{k\in K}\bigr)\hookrightarrow\ell_{v}^{r}\left(I\right)=\ell^{r}\bigl(\left[\ell_{v}^{r}\left(\smash{I^{\left(k,n\right)}}\right)\right]_{k\in K}\bigr).
\]
In view of Lemmas \ref{lem:NestedDisjointEmbedding} and \ref{lem:EmbeddingBetweenWeightedSequenceSpaces},
this yields $\theta_{k}^{\left(n\right)}<\infty$ for all $k\in K$,
where
\[
\theta_{k}^{\left(n\right)}=\begin{cases}
\vertiii{\ell_{u}^{q}\left(\smash{I^{\left(k,n\right)}}\right)\hookrightarrow\ell_{v}^{r}\left(\smash{I^{\left(k,n\right)}}\right)}=\left\Vert \left(v_{i}/u_{k,i}\right)_{i\in I^{\left(k,n\right)}}\right\Vert _{\ell^{r\cdot\left(q/r\right)'}}\;, & \text{if }\delta_{k}\in X,\\
1, & \text{if }\delta_{k}\notin X.
\end{cases}
\]
By the same lemmas, we also get the boundedness of
\[
\gamma^{\left(n\right)}:X\to\ell^{r}\left(K\right),\left(x_{k}\right)_{k\in K}\mapsto\left(\smash{\theta_{k}^{\left(n\right)}}\cdot x_{k}\right)_{k\in K}\qquad\text{with}\qquad\vertiii{\smash{\gamma^{\left(n\right)}}}=\vertiii{\iota\circ\smash{\iota^{\left(n\right)}}}\leq C\cdot\vertiii{\iota}.
\]

Finally, we recall $I^{\left(k\right)}=\bigcup_{n=1}^{N}I^{\left(k,n\right)}$.
Using the quasi-triangle inequality for $\ell^{r\cdot\left(q/r\right)'}$,
this implies
\[
\theta_{k}\leq C_{r,q,N}\cdot\sum_{n=1}^{N}\theta_{k}^{\left(n\right)}\qquad\forall\,k\in K\text{ for which }\delta_{k}\in X.
\]
If $\delta_{k}\notin X$, we trivially have $\theta_{k}=1\leq N=\sum_{n=1}^{N}\theta_{k}^{\left(n\right)}$.
All in all, these considerations show (using the (quasi)-triangle
inequality for $\ell^{r}$) that 
\begin{align*}
\left\Vert \gamma\left(x\right)\right\Vert _{\ell^{r}}=\left\Vert \left(\theta_{k}\cdot x_{k}\right)_{k\in K}\right\Vert _{\ell^{r}} & \leq C_{r,q,N}\cdot\left\Vert \left(\,\smash{\sum_{n=1}^{N}}\,\vphantom{\sum}\theta_{k}^{\left(n\right)}\cdot\left|x_{k}\right|\,\right)_{k\in K}\right\Vert _{\ell^{r}}\vphantom{\sum_{n=1}^{N}}\\
 & \leq C_{r,N}C_{r,q,N}\cdot\sum_{n=1}^{N}\left\Vert \left(\smash{\theta_{k}^{\left(n\right)}}\cdot\left|x_{k}\right|\right)_{k\in K}\right\Vert _{\ell^{r}}\\
 & \leq C_{r,N}C_{r,q,N}\cdot\left(\,\smash{\sum_{n=1}^{N}}\,\vphantom{\sum}\vertiii{\smash{\gamma^{\left(n\right)}}}\,\right)\vphantom{\sum_{n=1}^{N}}\cdot\left\Vert x\right\Vert _{X}<\infty
\end{align*}
for all $x\in X$. Thus, $\gamma$ is bounded with $\vertiii{\gamma}\leq C\cdot NC_{r,N}C_{r,q,N}\cdot\vertiii{\iota}$.

\medskip{}

``$\Leftarrow$'': Because of $I^{\left(k,\natural\right)}\subset I^{\left(k\right)}$
for all $k\in K$, it is clear that if $\gamma$ is bounded, then
so is $\gamma^{\natural}$, with $\vertiii{\gamma^{\natural}}\leq\vertiii{\gamma}$.
Thus, it suffices to assume that $\theta_{k}^{\natural}<\infty$ for
all $k\in K$ and that $\gamma^{\natural}$ is bounded. We need to
show that $\iota$ is well-defined and bounded with $\vertiii{\iota}\leq\vertiii{\gamma^{\natural}}$.

The first part of Lemma~\ref{lem:NestedEmbeddingReductionExistenceOfDisjointifiedSets}
(applied to $\bigl(I^{\left(k,\natural\right)}\bigr)_{k\in K}$ instead
of $\left(I^{\left(k\right)}\right)_{k\in K}$ and with $I^{\left(k,0\right)}:=\emptyset$
for all $k\in K$) yields a partition $I=\biguplus_{k\in K}I^{\left(k,\star\right)}$
with $I^{\left(k,\star\right)}\subset I^{\left(k,\natural\right)}\subset I^{\left(k\right)}$
for all $k\in K$. Since $X$ is solid, we have a (quasi)-norm decreasing
embedding $X\left(\left[\ell_{u}^{q}\left(\smash{I^{\left(k\right)}}\right)\right]_{k\in K}\right)\hookrightarrow X\left(\left[\ell_{u}^{q}\left(\smash{I^{\left(k,\star\right)}}\right)\right]_{k\in K}\right)$,
so that we have $\vertiii{\iota}\leq\vertiii{\iota_{1}}$ for
\[
\iota_{1}:X\left(\left[\ell_{u}^{q}\left(\smash{I^{\left(k,\star\right)}}\right)\right]_{k\in K}\right)\Xhookrightarrow !\ell_{v}^{r}\left(I\right)=\ell^{r}\left(\left[\ell_{v}^{r}\left(\smash{I^{\left(k,\star\right)}}\right)\right]_{k\in K}\right).
\]

Next, note that Lemma~\ref{lem:EmbeddingBetweenWeightedSequenceSpaces}
yields
\[
\widetilde{\theta_{k}}:=\vertiii{\ell_{u}^{q}\left(\smash{I^{\left(k,\star\right)}}\right)\hookrightarrow\ell_{v}^{r}\left(\smash{I^{\left(k,\star\right)}}\right)}=\left\Vert \left(v_{i}/u_{k,i}\right)_{i\in I^{\left(k,\star\right)}}\right\Vert _{\ell^{r\cdot\left(q/r\right)'}}\leq\left\Vert \left(v_{i}/u_{k,i}\right)_{i\in I^{\left(k,\natural\right)}}\right\Vert _{\ell^{r\cdot\left(q/r\right)'}}=\theta_{k}^{\natural}<\infty
\]
for all $k\in K$ with $\delta_{k}\in X$. Thus, if we set $\widetilde{\theta_{k}}:=1$
for those $k\in K$ with $\delta_{k}\notin X$, then Lemma~\ref{lem:NestedDisjointEmbedding}
shows that $\iota_{1}$ has the same norm as the map
\[
\gamma_{1}:X\to\ell^{r}\left(K\right),\left(x_{k}\right)_{k\in K}\mapsto\left(\smash{\widetilde{\theta_{k}}}\cdot x_{k}\right)_{k\in K}
\]
But using $\widetilde{\theta_{k}}\leq\theta_{k}^{\natural}$ for all
$k\in K$, we finally get $\vertiii{\iota}\leq\vertiii{\iota_{1}}=\vertiii{\gamma_{1}}\leq\vertiii{\gamma^{\natural}}<\infty$,
as desired.

\medskip{}

Ad (2): ``$\Rightarrow$'': Assume that $\iota$ is bounded. Using
Lemma~\ref{lem:NestedEmbeddingReductionExistenceOfDisjointifiedSets},
we get for each $n\in\underline{N}$ a partition $I=\biguplus_{k\in K}I^{\left(k,n\right)}$
such that we have $I^{\left(k\right)}=\bigcup_{n=1}^{N}I^{\left(k,n\right)}$
for all $k\in K$. But because of $I^{\left(k,n\right)}\subset I^{\left(k\right)}$
and by solidity of $X$, it is not hard to see that we have a (quasi)-norm
decreasing embedding
\[
\varrho_{n}:X\left(\left[\ell_{u}^{q}\left(\smash{I^{\left(k\right)}}\right)\right]_{k\in K}\right)\hookrightarrow X\left(\left[\ell_{u}^{q}\left(\smash{I^{\left(k,n\right)}}\right)\right]_{k\in K}\right)
\]
for each $n\in\underline{N}$. By composition, we get boundedness
of
\[
\iota^{\left(n\right)}:=\varrho_{n}\circ\iota:\ell_{v}^{r}\left(I\right)=\ell^{r}\left(\left[\ell_{v}^{r}\left(\smash{I^{\left(k,n\right)}}\right)\right]_{k\in K}\right)\hookrightarrow X\left(\left[\ell_{u}^{q}\left(\smash{I^{\left(k,n\right)}}\right)\right]_{k\in K}\right).
\]

In view of Lemmas \ref{lem:NestedDisjointEmbedding} and \ref{lem:EmbeddingBetweenWeightedSequenceSpaces}
(and since $\delta_{k}\in\ell^{r}\left(K\right)$ for all $k\in K$
and since $\left(I^{\left(k,n\right)}\right)_{k\in K}$ is pairwise
disjoint), this yields boundedness of
\[
\eta^{\left(n\right)}:\ell^{r}\left(K\right)\to X,\left(x_{k}\right)_{k\in K}\mapsto\left(\smash{\vartheta_{k}^{\left(n\right)}}\cdot x_{k}\right)_{k\in K}\qquad\text{with}\qquad\vertiii{\smash{\eta^{\left(n\right)}}}=\vertiii{\smash{\iota^{\left(n\right)}}}\leq\vertiii{\iota}\,,
\]
where the same lemmas also show that
\[
\vartheta_{k}^{\left(n\right)}:=\vertiii{\ell_{v}^{r}\left(\smash{I^{\left(k,n\right)}}\right)\hookrightarrow\ell_{u}^{q}\left(\smash{I^{\left(k,n\right)}}\right)}=\left\Vert \left(u_{k,i}/v_{i}\right)_{i\in I^{\left(k,n\right)}}\right\Vert _{\ell^{q\cdot\left(r/q\right)'}}<\infty\quad\text{ for all }k\in K.
\]

Now, using $I^{\left(k\right)}=\bigcup_{n=1}^{N}I^{\left(k,n\right)}$
for each $k\in K$, the (quasi)-triangle inequality for $\ell^{q\cdot\left(r/q\right)'}$
and Lemma~\ref{lem:EmbeddingBetweenWeightedSequenceSpaces} show
\[
\vartheta_{k}=\left\Vert \left(u_{k,i}/v_{i}\right)_{i\in I^{\left(k\right)}}\right\Vert _{\ell^{q\cdot\left(r/q\right)'}}\leq C_{N,q,r}\cdot\sum_{n=1}^{N}\left\Vert \left(u_{k,i}/v_{i}\right)_{i\in I^{\left(k,n\right)}}\right\Vert _{\ell^{q\cdot\left(r/q\right)'}}=C_{N,q,r}\cdot\sum_{n=1}^{N}\,\vartheta_{k}^{\left(n\right)}<\infty
\]
for all $k\in K$. All in all, using the solidity of $X$ and the
quasi-triangle inequality for $X$, we get
\begin{align*}
\left\Vert \left(\vartheta_{k}\cdot x_{k}\right)_{k\in K}\right\Vert _{X}\leq C_{N,q,r}\cdot\left\Vert \left(\,\smash{\sum_{n=1}^{N}}\,\vphantom{\sum}\vartheta_{k}^{\left(n\right)}\cdot\left|x_{k}\right|\,\right)_{k\in K}\right\Vert _{X}\vphantom{\sum_{n=1}^{N}} & \leq C_{N,C_{X}}C_{N,q,r}\cdot\sum_{n=1}^{N}\,\left\Vert \eta^{\left(n\right)}\left(x\right)\right\Vert _{X}\\
 & \leq C_{N,C_{X}}C_{N,q,r}\cdot\left(\,\smash{\sum_{n=1}^{N}}\,\,\vphantom{\sum}\vertiii{\smash{\eta^{\left(n\right)}}}\,\right)\vphantom{\sum_{n=1}^{N}}\cdot\left\Vert x\right\Vert _{\ell^{r}}\\
 & \leq C_{N,C_{X}}C_{N,q,r}\cdot N\vertiii{\iota}\cdot\left\Vert x\right\Vert _{\ell^{r}}<\infty
\end{align*}
for all $x=\left(x_{k}\right)_{k\in K}\in\ell^{r}\left(K\right)$,
so that $\eta$ is well-defined and bounded with $\vertiii{\eta}\leq C_{1}\cdot\vertiii{\iota}$
for some constant $C_{1}=C_{1}\left(N,C_{X},q,r\right)$.

\medskip{}

``$\Leftarrow$'': It is clear that if $\eta$ is bounded (and $\vartheta_{k}<\infty$
for all $k\in K$), then so is $\eta^{\natural}$, with $\vertiii{\eta^{\natural}}\leq\vertiii{\eta}$
(and with $\vartheta_{k}^{\natural}\leq\vartheta_{k}<\infty$ for
all $k\in K$). Thus, it suffices to show that $\iota$ is bounded
with $\vertiii{\iota}\leq C_{2}\cdot\vertiii{\eta^{\natural}}$, assuming
that $\eta^{\natural}$ is well-defined and bounded (and that $\vartheta_{k}^{\natural}<\infty$
for all $k\in K$).

The first part of Lemma~\ref{lem:NestedEmbeddingReductionExistenceOfDisjointifiedSets}
(applied to the family $\bigl(I^{\left(k,\natural\right)}\bigr)_{k\in K}$
instead of $\left(I^{\left(k\right)}\right)_{k\in K}$, and with $I^{\left(k,0\right)}=\emptyset$
for all $k\in K$) yields a partition $I=\biguplus_{k\in K}I^{\left(k,\star\right)}$
with $I^{\left(k,\star\right)}\subset I^{\left(k,\natural\right)}$
for all $k\in K$. In particular, we have
\[
\infty>\vartheta_{k}^{\natural}=\left\Vert \left(u_{k,i}/v_{i}\right)_{i\in I^{\left(k,\natural\right)}}\right\Vert _{\ell^{q\cdot\left(r/q\right)'}}\geq\left\Vert \left(u_{k,i}/v_{i}\right)_{i\in I^{\left(k,\star\right)}}\right\Vert _{\ell^{q\cdot\left(r/q\right)'}}\overset{\left(\ast\right)}{=}\vertiii{\ell_{v}^{r}\left(\smash{I^{\left(k,\star\right)}}\right)\hookrightarrow\ell_{u}^{q}\left(\smash{I^{\left(k,\star\right)}}\right)}=:\widetilde{\vartheta_{k}}
\]
for all $k\in K$, where the step marked with $\left(\ast\right)$
is justified by Lemma~\ref{lem:EmbeddingBetweenWeightedSequenceSpaces}.
Thus, by solidity of $X$, we see that the boundedness of $\eta^{\natural}$
implies boundedness of
\[
\widetilde{\eta}:\ell^{r}\left(K\right)\to X,\left(x_{k}\right)_{k\in K}\mapsto\left(\widetilde{\vartheta_{k}}\cdot x_{k}\right)_{k\in K}\,,\qquad\text{with}\qquad\vertiii{\widetilde{\eta}}\leq\vertiii{\eta^{\natural}}.
\]

But in view of Lemma~\ref{lem:NestedDisjointEmbedding}, this yields
boundedness of
\[
\widetilde{\iota}:\ell_{v}^{r}\left(I\right)\overset{\left(\blacklozenge\right)}{=}\ell^{r}\left(\left[\ell_{v}^{r}\left(\smash{I^{\left(k,\star\right)}}\right)\right]_{k\in K}\right)\hookrightarrow X\left(\left[\ell_{u}^{q}\left(\smash{I^{\left(k,\star\right)}}\right)\right]_{k\in K}\right),
\]
with $\vertiii{\widetilde{\iota}}\leq\vertiii{\widetilde{\eta}}\leq\vertiii{\eta^{\natural}}$,
where the (isometric(!\@)) identity marked with $\left(\blacklozenge\right)$
is a consequence of $I=\biguplus_{k\in K}I^{\left(k,\star\right)}$.

Finally, Lemma~\ref{lem:NestedEmbeddingReductionToDisjointSets}
shows that the identity map $\identity:X\left(\left[\ell_{u}^{q}\left(\smash{I^{\left(k,\star\right)}}\right)\right]_{k\in K}\right)\hookrightarrow X\left(\left[\ell_{u}^{q}\left(\smash{I^{\left(k\right)}}\right)\right]_{k\in K}\right)$
is bounded with $\vertiii{\identity}\leq C_{2}=C_{2}\left(\vertiii{\Theta},C_{u},N,q\right)$.
All in all, we see that $\iota=\identity\circ\widetilde{\iota}$ is
bounded with $\vertiii{\iota}\leq C_{2}\cdot\vertiii{\widetilde{\iota}}\leq C_{2}\cdot\vertiii{\eta^{\natural}}$,
as claimed.

\medskip{}

Ad (3): Note that we have $\delta_{k}\in X=\ell_{w}^{s}\left(K\right)$
for all $k\in K$, so that $\theta_{k}^{\natural}=\left\Vert \left(v_{i}/u_{k,i}\right)_{i\in I^{\left(k,\natural\right)}}\right\Vert _{\ell^{r\cdot\left(q/r\right)'}}$
for all $k\in K$. In particular, $\theta_{k}^{\natural}=0$ if and
only if $I^{\left(k,\natural\right)}=\emptyset$ and likewise $\vartheta_{k}^{\natural}=0$
if and only if $I^{\left(k,\natural\right)}=\emptyset$.

Thus, let $K_{0}:=\left\{ k\in K\with\smash{I^{\left(k,\natural\right)}}\neq\emptyset\right\} $.
It is then not hard to see that the embeddings
\[
\gamma^{\left(0\right)}:\ell_{w}^{s}\left(K_{0}\right)\hookrightarrow\ell_{\theta^{\natural}}^{r}\left(K_{0}\right)\qquad\text{ and }\qquad\eta^{\left(0\right)}:\ell^{r}\left(K_{0}\right)\hookrightarrow\ell_{w\cdot\vartheta^{\natural}}^{s}\left(K_{0}\right)
\]
satisfy $\vertiii{\gamma^{\natural}}=\vertiii{\gamma^{\left(0\right)}}$
and $\vertiii{\eta^{\natural}}=\vertiii{\eta^{\left(0\right)}}$.
But now, Lemma~\ref{lem:EmbeddingBetweenWeightedSequenceSpaces}
shows
\[
\vertiii{\smash{\gamma^{\left(0\right)}}}=\left\Vert \left(\smash{\theta_{k}^{\natural}}/w_{k}\right)_{k\in K_{0}}\right\Vert _{\ell^{r\cdot\left(s/r\right)'}}=\left\Vert \left(\smash{\theta_{k}^{\natural}}/w_{k}\right)_{k\in K}\right\Vert _{\ell^{r\cdot\left(s/r\right)'}}
\]
and
\[
\vertiii{\smash{\eta^{\left(0\right)}}}=\left\Vert \left(w_{k}\cdot\smash{\vartheta_{k}^{\natural}}\right)_{k\in K_{0}}\right\Vert _{\ell^{s\cdot\left(r/s\right)'}}=\left\Vert \left(w_{k}\cdot\smash{\vartheta_{k}^{\natural}}\right)_{k\in K}\right\Vert _{\ell^{s\cdot\left(r/s\right)'}}\:\,.\qedhere
\]
\end{proof}
We close this section by showing that $\CalQ$-regular sequence spaces
again yield $\CalQ$-regular sequence spaces when weighted with $\CalQ$-moderate
weights.
\begin{lem}
\label{lem:ModeratelyWeightedSpacesAreRegular}Let $\CalQ=\left(Q_{i}\right)_{i\in I}$
be a family of subsets of a set $\CalO$. If $X\subset\Compl^{I}$
is $\CalQ$-regular and if $u=\left(u_{i}\right)_{i\in I}$ is $\CalQ$-moderate,
then $X_{u}$ (as introduced in Definition~\ref{def:WeightedSequenceSpaces})
is $\CalQ$-regular, with
\[
\vertiii{\Gamma_{\CalQ}}_{X_{u}\to X_{u}}\leq\vertiii{\Gamma_{\CalQ}}_{X\to X}\cdot C_{u,\CalQ}.
\]
In particular, for any $q\in\left(0,\infty\right]$, the space $\ell_{u}^{q}\left(I\right)$
is $\CalQ$-regular with
\[
\vertiii{\Gamma_{\CalQ}}_{\ell_{u}^{q}\to\ell_{u}^{q}}\leq C_{u,\CalQ}\cdot N_{\CalQ}^{\max\left\{ 1,\smash{\frac{1}{q}}\right\} }\:.\qedhere
\]
\end{lem}

\begin{proof}
Let $x=\left(x_{i}\right)_{i\in I}\in X_{u}$ be arbitrary. We have
\[
u_{i}\cdot\left|x_{i}^{\ast}\right|\leq\sum_{\ell\in i^{\ast}}u_{i}\cdot\left|x_{\ell}\right|\leq C_{u,\CalQ}\cdot\sum_{\ell\in i^{\ast}}u_{\ell}\left|x_{\ell}\right|=C_{u,\CalQ}\cdot\left(u\cdot\left|x\right|\right)_{i}^{\ast}.
\]
By definition of $X_{u}$, we have $u\cdot x\in X$. By solidity of
$X$, this implies $u\cdot\left|x\right|=\left|u\cdot x\right|\in X$
and thus also $\left(u\cdot\left|x\right|\right)^{\ast}\in X$, since
$X$ is $\CalQ$-regular. Again by solidity of $X$, we get $u\cdot\left|x^{\ast}\right|\in X$,
and as desired
\begin{align*}
\left\Vert x^{\ast}\right\Vert _{X_{u}}=\left\Vert u\cdot x^{\ast}\right\Vert _{X}=\left\Vert u\cdot\left|x^{\ast}\right|\right\Vert _{X} & \leq\left\Vert C_{u,\CalQ}\cdot\left(u\cdot\left|x\right|\right)^{\ast}\right\Vert _{X}\\
 & \leq C_{u,\CalQ}\vertiii{\Gamma_{\CalQ}}_{X\to X}\cdot\left\Vert u\cdot\left|x\right|\right\Vert _{X}=C_{u,\CalQ}\vertiii{\Gamma_{\CalQ}}_{X\to X}\cdot\left\Vert x\right\Vert _{X_{u}}<\infty.
\end{align*}

For the second part, it suffices—in view of the first part—to show
that $\ell^{q}\left(I\right)$ is $\CalQ$-regular with $\vertiii{\Gamma_{\CalQ}}_{\ell^{q}\to\ell^{q}}\leq N_{\CalQ}^{\max\left\{ 1,\smash{\frac{1}{q}}\right\} }$.
Let us first consider the case $1\leq q<\infty$. Here, Hölder's inequality
shows because of $\left|i^{\ast}\right|\leq N_{\CalQ}$ that
\[
\left|\sum_{\ell\in i^{\ast}}x_{\ell}\right|^{q}\leq\left(\left\Vert \left(x_{\ell}\right)_{\ell\in i^{\ast}}\right\Vert _{\ell^{q}}\cdot\left\Vert \left(1\right)_{\ell\in i^{\ast}}\right\Vert _{\ell^{q'}}\right)^{q}\leq N_{\CalQ}^{q-1}\cdot\left\Vert \left(x_{\ell}\right)_{\ell\in i^{\ast}}\right\Vert _{\ell^{q}}^{q}=N_{\CalQ}^{\left(q-1\right)_{+}}\cdot\sum_{\ell\in i^{\ast}}\left|x_{\ell}\right|^{q}.
\]
Likewise, if $q\in\left(0,1\right)$, then $\left|\sum_{\ell\in i^{\ast}}x_{\ell}\right|^{q}\leq\sum_{\ell\in i^{\ast}}\left|x_{\ell}\right|^{q}=N_{\CalQ}^{\left(q-1\right)_{+}}\cdot\sum_{\ell\in i^{\ast}}\left|x_{\ell}\right|^{q}$.
Therefore, we see for arbitrary $q\in\left(0,\infty\right)$ that
\[
\left\Vert x^{\ast}\right\Vert _{\ell^{q}}^{q}\leq N_{\CalQ}^{\left(q-1\right)_{+}}\cdot\sum_{i\in I}\sum_{\ell\in i^{\ast}}\left|x_{\ell}\right|^{q}=N_{\CalQ}^{\left(q-1\right)_{+}}\cdot\sum_{\ell\in I}\left[\left|x_{\ell}\right|^{q}\cdot\sum_{i\in\ell^{\ast}}1\right]\leq N_{\CalQ}^{\left(q-1\right)_{+}+1}\cdot\left\Vert x\right\Vert _{\ell^{q}}^{q}<\infty\,,
\]
which proves the desired estimate, since $\frac{1}{q}\cdot\left[\left(q-1\right)_{+}+1\right]=\max\left\{ 1,q^{-1}\right\} $.

Finally, for $q=\infty$, simply note $\left|x_{i}^{\ast}\right|\leq\sum_{\ell\in i^{\ast}}\left|x_{\ell}\right|\leq\left|i^{\ast}\right|\cdot\left\Vert x\right\Vert _{\ell^{\infty}}\leq N_{\CalQ}\cdot\left\Vert x\right\Vert _{\ell^{\infty}}$,
which shows that $\left\Vert x^{\ast}\right\Vert _{\ell^{\infty}}\leq N_{\CalQ}\cdot\left\Vert x\right\Vert _{\ell^{\infty}}=N_{\CalQ}^{\max\left\{ 1,q^{-1}\right\} }\cdot\left\Vert x\right\Vert _{\ell^{q}}$,
as desired.
\end{proof}

\section{Sufficient conditions for embeddings}

\label{sec:SufficientConditions}In this section, we properly start
to investigate the existence of embeddings between decomposition spaces.
Precisely, we derive \emph{sufficient} conditions for the existence
of such embeddings. As noted in the introduction, these conditions
take the form of embeddings between certain (nested) sequence spaces.

All results in this section will be based on the following two lemmata
which investigate how the $L^{p}$-norm of a function $f$ is related
to the $L^{p}$-norm of its frequency localized parts $f_{i}=\Fourier^{-1}\left(\smash{\varphi_{i}\cdot\widehat{f}}\,\right)$.
Thanks to Plancherel's theorem and since the covering $\CalQ$ is
``almost disjoint''—because $\CalQ$ is admissible—this is easy
in case of $p=2$. Our more general results will be derived from this
basic case via interpolation.

Our first result shows how the $L^{p}$-norm of $f$ can be estimated
in terms of the $L^{p}$-norms of the $f_{i}$. This result is a simplified
form of \cite[Lemma 3.1]{DecompositionIntoSobolev} and similar to
\cite[Lemma 5.1.2]{VoigtlaenderPhDThesis}. Since it is crucial for
the remainder of the paper, we nevertheless provide a proof.
\begin{lem}
\label{lem:BasicEstimateFineIntoCoarse}Let $\CalQ=\left(Q_{i}\right)_{i\in I}$
be an $L^{1}$-decomposition covering of the open set $\emptyset\neq\CalO\subset\R^{\dimension}$.
Furthermore, let $I_{0}\subset I$ be \emph{finite}, let $p\in\left(0,\infty\right]$
and $k\in\N_{0}$ and assume that for each $i\in I_{0}$, we are given
some $f_{i}\in\Schwartz'\left(\R^{\dimension}\right)\cap L^{p}\left(\R^{\dimension}\right)$
with Fourier support $\supp\widehat{f_{i}}\subset\overline{Q_{i}^{k\ast}}$.

Then
\[
\vphantom{\sum_{i\in I_{0}}}\left\Vert \,\vphantom{\sum}\smash{\sum_{i\in I_{0}}}f_{i}\,\right\Vert _{L^{p}}\leq C\cdot\left\Vert \left(\left\Vert f_{i}\right\Vert _{L^{p}}\right)_{i\in I_{0}}\right\Vert _{\ell^{\LowerExpo p}}
\]
for $C=C_{\CalQ,\Phi,1}\cdot N_{\CalQ}^{2k+4}$, where $\Phi=\left(\varphi_{i}\right)_{i\in I}$
is some $L^{1}$-BAPU for $\CalQ$. Here, $\LowerExpo p=\min\left\{ p,p'\right\} $.
\end{lem}

\begin{proof}
We first handle the (easier) case $p\in\left(0,1\right]$. In this
case, we have $p'=\infty>p$ and hence $\LowerExpo p=p$\@. For $p\in\left(0,1\right]$,
it is well-known that $\left\Vert \mybullet\right\Vert _{L^{p}}$
is a $p$-norm, i.e.\@ we have $\left\Vert f+g\right\Vert _{L^{p}}^{p}\leq\left\Vert f\right\Vert _{L^{p}}^{p}+\left\Vert g\right\Vert _{L^{p}}^{p}$
for all measurable $f,g$. Indeed, this is an immediate consequence
of the estimate
\[
\left|a+b\right|^{p}\leq\left|a\right|^{p}+\left|b\right|^{p}
\]
which holds for all $a,b\in\Compl$, since $p\in\left(0,1\right]$.
Using a straightforward induction, we thus get (since $I_{0}$ is
finite) that
\[
\left\Vert \,\smash{\sum_{i\in I_{0}}}\,\vphantom{\sum}f_{i}\,\right\Vert _{L^{p}}^{p}\leq\sum_{i\in I_{0}}\left\Vert f_{i}\right\Vert _{L^{p}}^{p}=\left\Vert \left(\left\Vert f_{i}\right\Vert _{L^{p}}\right)_{i\in I_{0}}\right\Vert _{\ell^{\LowerExpo p}}^{p}\,.
\]
Taking $p$th roots completes the proof for $p\in\left(0,1\right]$.
In fact, we have shown that $C=1$ is a possible choice in this case\footnote{The choice $C=C_{\CalQ,\Phi,1}\cdot N_{\CalQ}^{2k+4}$ is also possible,
for the following reason: If $\Phi=\left(\varphi_{i}\right)_{i\in I}$
is an $L^{1}$-BAPU for $\CalQ$, then by Fourier inversion $\left|\varphi_{i}\left(\xi\right)\right|\leq\left\Vert \Fourier^{-1}\varphi_{i}\right\Vert _{L^{1}}\leq C_{\CalQ,\Phi,1}$
for all $\xi\in\R^{\dimension}$. But for any fixed $i_{0}\in I$
and $\xi\in Q_{i_{0}}$, we have $1=\left|\sum_{i\in I}\varphi_{i}\left(\xi\right)\right|\leq\sum_{i\in i_{0}^{\ast}}\left|\varphi_{i}\left(\xi\right)\right|\leq N_{\CalQ}\cdot C_{\CalQ,\Phi,1}\leq C$,
since $\varphi_{i}\left(\xi\right)=0$ for $i\in I\setminus i_{0}^{\ast}$.}. Thus, we can assume $p\geq1$ for the remainder of the proof.

Fix some $L^{1}$-BAPU $\Phi=\left(\varphi_{i}\right)_{i\in I}$ for
$\CalQ$. Our first aim is to show $\widehat{f_{i}}=\varphi_{i}^{\left(k+2\right)\ast}\widehat{f_{i}}$
for all $i\in I_{0}$. To this end, note that Lemma~\ref{lem:PartitionCoveringNecessary}
implies $\varphi_{i}^{\left(k+1\right)\ast}\equiv1$ on $Q_{i}^{k\ast}$
and hence also $\varphi_{i}^{\left(k+1\right)\ast}\equiv1$ on $\overline{Q_{i}^{k\ast}}\supset\supp\widehat{f_{i}}$.
Now, note that $\overline{Q_{i}^{k\ast}}\subset\left(Q_{i}^{\left(k+1\right)\ast}\right)^{\circ}=:U_{i}$,
since we just showed $\overline{Q_{i}^{k\ast}}\subset\left(\varphi_{i}^{\left(k+1\right)\ast}\right)^{-1}\left(\Compl^{\ast}\right)$,
which is an open(!) subset of $Q_{i}^{\left(k+1\right)\ast}$.

Finally, Lemma~\ref{lem:PartitionCoveringNecessary} shows $\varphi_{i}^{\left(k+2\right)\ast}\equiv1$
on $Q_{i}^{\left(k+1\right)\ast}\supset U_{i}\supset\supp\widehat{f_{i}}$.
Hence, for $\psi\in\TestFunctionSpace{\R^{\dimension}}$, we have
\[
\supp\!\left(\psi-\varphi_{i}^{\left(k+2\right)\ast}\psi\right)\!\subset U_{i}^{c}\subset\!\left(\supp\widehat{f_{i}}\right)^{c},\quad\!\!\!\text{and thus}\quad\!\!\!\left\langle \smash{\widehat{f_{i}}}\,,\psi\right\rangle _{\Schwartz'}=\left\langle \smash{\widehat{f_{i}}}\,,\,\varphi_{i}^{\left(k+2\right)\ast}\psi\right\rangle _{\Schwartz'}\!\!=\left\langle \varphi_{i}^{\left(k+2\right)\ast}\smash{\widehat{f_{i}}}\,,\,\psi\right\rangle _{\Schwartz'}\:.
\]
Since $\TestFunctionSpace{\R^{\dimension}}\subset\Schwartz\left(\R^{\dimension}\right)$
is dense, we have thus shown $\widehat{f_{i}}=\varphi_{i}^{\left(k+2\right)\ast}\widehat{f_{i}}$,
as desired. Hence,
\[
f_{i}=\Fourier^{-1}\left(\varphi_{i}^{\left(k+2\right)\ast}\cdot\widehat{f_{i}}\right)\qquad\text{ for all }i\in I.
\]

As a consequence, it suffices to show that the map
\[
\Phi_{p}:\ell^{\LowerExpo p}\!\!\left(I_{0};\,L^{p}\left(\R^{\dimension}\right)\right)\to L^{p}\left(\R^{\dimension}\right),\left(g_{i}\right)_{i\in I_{0}}\mapsto\sum_{i\in I_{0}}\Fourier^{-1}\left(\varphi_{i}^{\left(k+2\right)\ast}\cdot\widehat{g_{i}}\right)=\sum_{i\in I_{0}}\left[\left(\Fourier^{-1}\varphi_{i}^{\left(k+2\right)\ast}\right)\ast g_{i}\right]
\]
is bounded for all $p\in\left[1,\infty\right]$, with $\vertiii{\Phi_{p}}\leq C:=C_{\CalQ,\Phi,1}\cdot N_{\CalQ}^{2k+4}$,
since this will imply as desired that
\begin{align*}
C\cdot\left\Vert \left(\left\Vert f_{i}\right\Vert _{L^{p}}\right)_{i\in I_{0}}\right\Vert _{\ell^{\LowerExpo p}} & \geq\vertiii{\Phi_{p}}\cdot\left\Vert \left(f_{i}\right)_{i\in I_{0}}\right\Vert _{\ell^{\LowerExpo p}\left(I_{0};\,L^{p}\left(\R^{\dimension}\right)\right)}\\
 & \geq\left\Vert \Phi_{p}\left(\left(f_{i}\right)_{i\in I_{0}}\right)\right\Vert _{L^{p}}=\vphantom{\sum_{i\in I_{0}}}\left\Vert \,\smash{\sum_{i\in I_{0}}}\,\vphantom{\sum}\Fourier^{-1}\left(\varphi_{i}^{\left(k+2\right)\ast}\cdot\widehat{f_{i}}\right)\,\right\Vert _{L^{p}}=\left\Vert \,\smash{\sum_{i\in I_{0}}}\,\vphantom{\sum}f_{i}\,\right\Vert _{L^{p}}\:.
\end{align*}

Note that each $\Phi_{p}$ is well-defined, since $I_{0}$ is finite
and since $L^{p}\to L^{p},g\mapsto\left(\Fourier^{-1}\varphi_{i}^{\left(k+2\right)\ast}\right)\ast g$
is well-defined and bounded as a consequence of Young's convolution
inequality and since we have $\Fourier^{-1}\varphi_{i}^{\left(k+2\right)\ast}\in L^{1}\left(\R^{\dimension}\right)$
as a \emph{finite} sum of $L^{1}$-functions; in fact,
\begin{equation}
\left\Vert \Fourier^{-1}\varphi_{i}^{\left(k+2\right)\ast}\right\Vert _{L^{1}}\leq\sum_{\ell\in i^{\left(k+2\right)\ast}}\left\Vert \Fourier^{-1}\varphi_{\ell}\right\Vert _{L^{1}}\leq C_{\CalQ,\Phi,1}\cdot\left|i^{\left(k+2\right)\ast}\right|\leq C_{\CalQ,\Phi,1}\cdot N_{\CalQ}^{k+2},\label{eq:BasicEstimateFineIntoCoarseBAPUEstimate}
\end{equation}
where the last step is justified by Lemma~\ref{lem:SemiStructuredClusterInvariant}.

Let us first consider the case $p\in\left[1,2\right]$. Here, $\LowerExpo p=p$
and hence $\Phi_{p}:\ell^{p}\left(I_{0};\,L^{p}\left(\R^{\dimension}\right)\right)\to L^{p}\left(\R^{\dimension}\right)$.
It is thus sufficient to show $\vertiii{\Phi_{p}}\leq C$ for $p=1$
and $p=2$, since the general result for $p\in\left[1,2\right]$ then
follows by complex interpolation for vector-valued\footnote{In this case, one can even use the usual classical Riesz-Thorin interpolation
theorem (cf.\@ \cite[Theorem 6.27]{FollandRA}) by identifying $\ell^{p}\left(I_{0};\,L^{p}\left(\R^{\dimension}\right)\right)$
with $L^{p}\left(I_{0}\times\R^{\dimension},\mu\right)$, where $\mu$
is the product measure of the counting measure on $I_{0}$ and Lebesgue
measure on $\R^{\dimension}$. This argument is not applicable, however,
for $p\in\left[2,\infty\right]$, since the ``outer'' exponent $\LowerExpo p$
differs from the ``inner'' exponent $p$ for the space $\ell^{\LowerExpo p}\!\!\left(I_{0};\,L^{p}\left(\R^{\dimension}\right)\right)$.

Note furthermore that we have $\LowerExpo p\leq2<\infty$ for all
$p\in\left[1,\infty\right]$, so that \cite[Theorem 5.1.2]{BerghLoefstroemInterpolationSpaces}
is fully applicable.} $L^{p}$-spaces (cf.\@ \cite[Theorems 5.1.1 and 5.1.2]{BerghLoefstroemInterpolationSpaces}).

The case $p=1$ is easy: We simply use the triangle inequality for
$L^{1}\left(\R^{\dimension}\right)$ and Young's inequality for $L^{1}\left(\R^{\dimension}\right)$
to get
\begin{align*}
\left\Vert \Phi_{p}\left(\left(g_{i}\right)_{i\in I_{0}}\right)\right\Vert _{L^{1}}=\vphantom{\sum_{i\in I_{0}}}\left\Vert \,\smash{\sum_{i\in I_{0}}}\,\vphantom{\sum}\left(\Fourier^{-1}\varphi_{i}^{\left(k+2\right)\ast}\right)\ast g_{i}\,\right\Vert _{L^{1}} & \leq\sum_{i\in I_{0}}\left\Vert \Fourier^{-1}\varphi_{i}^{\left(k+2\right)\ast}\right\Vert _{L^{1}}\left\Vert g_{i}\right\Vert _{L^{1}}\\
\left({\scriptstyle \text{equation }\eqref{eq:BasicEstimateFineIntoCoarseBAPUEstimate}}\right) & \leq C_{\CalQ,\Phi,1}N_{\CalQ}^{k+2}\cdot\sum_{i\in I_{0}}\left\Vert g_{i}\right\Vert _{L^{1}}\\
 & =C_{\CalQ,\Phi,1}N_{\CalQ}^{k+2}\cdot\left\Vert \left(g_{i}\right)_{i\in I_{0}}\right\Vert _{\ell^{1}\left(I_{0};\,L^{1}\left(\R^{\dimension}\right)\right)}\\
 & \leq C_{\CalQ,\Phi,1}N_{\CalQ}^{2k+4}\cdot\left\Vert \left(g_{i}\right)_{i\in I_{0}}\right\Vert _{\ell^{1}\left(I_{0};\,L^{1}\left(\R^{\dimension}\right)\right)}\,.
\end{align*}

Now, let us consider the case $p=2$. Here, Plancherel's theorem implies
\begin{align*}
\left\Vert \Phi_{p}\left(\left(g_{i}\right)_{i\in I_{0}}\right)\right\Vert _{L^{2}}^{2}=\vphantom{\sum_{i\in I_{0}}}\left\Vert \Fourier^{-1}\left(\,\smash{\sum_{i\in I_{0}}}\,\vphantom{\sum}\varphi_{i}^{\left(k+2\right)\ast}\widehat{g_{i}}\,\right)\right\Vert _{L^{2}}^{2} & =\vphantom{\sum_{i\in I_{0}}}\left\Vert \,\smash{\sum_{i\in I_{0}}}\,\vphantom{\sum}\varphi_{i}^{\left(k+2\right)\ast}\cdot\widehat{g_{i}}\,\right\Vert _{L^{2}}^{2}\\
 & =\int_{\R^{\dimension}}\left|\,\smash{\sum_{i\in I_{0}}}\,\vphantom{\sum}\varphi_{i}^{\left(k+2\right)\ast}\left(\xi\right)\cdot\widehat{g_{i}}\left(\xi\right)\,\right|^{2}\,\d\xi\vphantom{\sum_{i\in I_{0}}}.
\end{align*}
Fix $\xi\in\R^{\dimension}$ for the moment and employ the Cauchy-Schwarz
inequality to derive
\[
\left|\smash{\sum_{i\in I_{0}}}\,\vphantom{\sum}\varphi_{i}^{\left(k+2\right)\ast}\left(\xi\right)\cdot\widehat{g_{i}}\left(\xi\right)\right|^{2}\leq\sum_{i\in I_{0}}\left|\varphi_{i}^{\left(k+2\right)\ast}\left(\xi\right)\right|^{2}\cdot\sum_{i\in I_{0}}\left|\widehat{g_{i}}\left(\xi\right)\right|^{2}\leq C_{\CalQ,\Phi,1}^{2}\cdot N_{\CalQ}^{3k+7}\cdot\sum_{i\in I_{0}}\left|\widehat{g_{i}}\left(\xi\right)\right|^{2}.
\]
The last estimate is justified as follows: In case of $\xi\in\R^{\dimension}\setminus\CalO$,
it is trivially true, since this entails $\varphi_{i}^{\left(k+2\right)\ast}\left(\xi\right)=0$
for all $i\in I\supset I_{0}$. In case of $\xi\in\CalO$, there is
some $i_{\xi}\in I$ satisfying $\xi\in Q_{i_{\xi}}$. Hence, if $\varphi_{i}^{\left(k+2\right)\ast}\left(\xi\right)\neq0$,
we have $\xi\in Q_{i}^{\left(k+2\right)\ast}\cap Q_{i_{\xi}}$ and
thus $i\in i_{\xi}^{\left(k+3\right)\ast}$. This implies
\[
\sum_{i\in I_{0}}\left|\varphi_{i}^{\left(k+2\right)\ast}\left(\xi\right)\right|^{2}\leq\sum_{i\in i_{\xi}^{\left(k+3\right)\ast}}\left|\varphi_{i}^{\left(k+2\right)\ast}\left(\xi\right)\right|^{2}\overset{\left(\ast\right)}{\leq}\left|i_{\xi}^{\left(k+3\right)\ast}\right|\cdot\left(C_{\CalQ,\Phi,1}\cdot N_{\CalQ}^{k+2}\right)^{2}\leq C_{\CalQ,\Phi,1}^{2}\cdot N_{\CalQ}^{3k+7}\,,
\]
where the last step used Lemma~\ref{lem:SemiStructuredClusterInvariant}.
Furthermore, the step marked with $\left(\ast\right)$ is justified
by recalling equation~(\ref{eq:BasicEstimateFineIntoCoarseBAPUEstimate})
and noting $\left\Vert \varphi_{i}^{\left(k+2\right)\ast}\right\Vert _{L^{\infty}}\leq\left\Vert \Fourier^{-1}\varphi_{i}^{\left(k+2\right)\ast}\right\Vert _{L^{1}}$,
by Fourier inversion, and because of $\left\Vert \widehat{g}\right\Vert _{L^{\infty}}\leq\left\Vert g\right\Vert _{L^{1}}$.\vspace{0.2cm}

All in all, we have shown
\begin{align*}
\left\Vert \Phi_{p}\left(\left(g_{i}\right)_{i\in I_{0}}\right)\right\Vert _{L^{2}}^{2} & \leq C_{\CalQ,\Phi,1}^{2}\cdot N_{\CalQ}^{3k+7}\cdot\int_{\R^{\dimension}}\sum_{i\in I_{0}}\left|\widehat{g_{i}}\left(\xi\right)\right|^{2}\,\d\xi\\
 & \leq\left(N_{\CalQ}^{2k+4}C_{\CalQ,\Phi,1}\right)^{2}\cdot\sum_{i\in I_{0}}\left\Vert \widehat{g_{i}}\right\Vert _{L^{2}}^{2}\\
\left({\scriptstyle \text{Plancherel}}\right) & =\left(N_{\CalQ}^{2k+4}C_{\CalQ,\Phi,1}\right)^{2}\cdot\sum_{i\in I_{0}}\left\Vert g_{i}\right\Vert _{L^{2}}^{2}\\
 & =\left[N_{\CalQ}^{2k+4}C_{\CalQ,\Phi,1}\cdot\left\Vert \left(g_{i}\right)_{i\in I_{0}}\right\Vert _{\ell^{2}\left(I_{0};\,L^{2}\left(\R^{\dimension}\right)\right)}\right]^{2}\,,
\end{align*}
as desired. This completes the proof for $p\in\left[1,2\right]$.

\medskip{}

Finally, we consider the case $p\in\left[2,\infty\right]$. Here,
$\LowerExpo p=p'$ and hence $\Phi_{p}:\ell^{p'}\left(I_{0};\,L^{p}\left(\R^{\dimension}\right)\right)\to L^{p}\left(\R^{\dimension}\right)$.
Note that every $p\in\left(2,\infty\right)$ satisfies $\frac{1}{p}=\frac{\theta}{2}+\frac{1-\theta}{\infty}=\frac{\theta}{p_{1}}+\frac{1-\theta}{p_{0}}$
for $\theta:=\frac{2}{p}\in\left(0,1\right)$ and $p_{0}=\infty$,
as well as $p_{1}=2$. Because of
\[
\frac{1}{p'}=1-\frac{1}{p}=\theta\left(1-\frac{1}{p_{1}}\right)+\left(1-\theta\right)\left(1-\frac{1}{p_{0}}\right)=\frac{\theta}{p_{1}'}+\frac{1-\theta}{p_{0}'}=\frac{\theta}{2}+\frac{1-\theta}{1},
\]
we get
\[
\ell^{p'}\left(I_{0};\,L^{p}\left(\R^{\dimension}\right)\right)=\left[\ell^{2}\left(I_{0};\,L^{2}\left(\R^{\dimension}\right)\right),\:\ell^{1}\left(I_{0};\,L^{\infty}\left(\R^{\dimension}\right)\right)\right]_{\theta}\,,
\]
where $\left[X,Y\right]_{\theta}$ denotes the space obtained by complex
interpolation between $X$ and $Y$ with parameter $\theta\in\left[0,1\right]$,
cf.\@ \cite[Theorems 5.1.1 and 5.1.2]{BerghLoefstroemInterpolationSpaces}.
Thus, it again suffices to prove $\vertiii{\Phi_{p}}\leq C$ for $p=2$
and $p=\infty$. But for $p=2$, we already established it above,
while proving the case $p\in\left[1,2\right]$.

Finally, for $p=\infty$, we can argue just as for $p=1$: The triangle
inequality for $L^{\infty}$ and Young's convolution inequality $L^{1}\ast L^{\infty}\hookrightarrow L^{\infty}$
yield
\begin{align*}
\left\Vert \Phi_{\infty}\left(\left(g_{i}\right)_{i\in I_{0}}\right)\right\Vert _{L^{\infty}} & =\vphantom{\sum_{i\in I_{0}}}\left\Vert \smash{\sum_{i\in I_{0}}}\,\vphantom{\sum}\left(\Fourier^{-1}\varphi_{i}^{\left(k+2\right)\ast}\right)\ast g_{i}\right\Vert _{L^{\infty}}\leq\sum_{i\in I_{0}}\left(\left\Vert \Fourier^{-1}\varphi_{i}^{\left(k+2\right)\ast}\right\Vert _{L^{1}}\cdot\left\Vert g_{i}\right\Vert _{L^{\infty}}\right)\\
\left({\scriptstyle \text{equation }\eqref{eq:BasicEstimateFineIntoCoarseBAPUEstimate}}\right) & \leq C_{\CalQ,\Phi,1}\cdot N_{\CalQ}^{k+2}\cdot\sum_{i\in I_{0}}\left\Vert g_{i}\right\Vert _{L^{\infty}}\leq C_{\CalQ,\Phi,1}N_{\CalQ}^{2k+4}\cdot\left\Vert \left(g_{i}\right)_{i\in I_{0}}\right\Vert _{\ell^{1}\left(I_{0};\,L^{\infty}\left(\R^{\dimension}\right)\right)}.
\end{align*}
Because of $\LowerExpo{\infty}=\infty'=1$, this completes the proof.
\end{proof}
The next lemma is concerned with a ``converse'' of Lemma~\ref{lem:BasicEstimateFineIntoCoarse},
i.e.\@ it shows how the $L^{p}$-norms of the individual ``pieces''
$f_{i}=\Fourier^{-1}\left(\smash{\varphi_{i}\cdot\widehat{f}}\,\right)$
can be estimated in terms of the $L^{p}$-norm of $f$ itself:
\begin{lem}
\label{lem:BasicEstimateCoarseInFine}Let $\CalQ=\left(Q_{i}\right)_{i\in I}$
be an $L^{1}$-decomposition covering of the open set $\emptyset\neq\CalO\subset\R^{\dimension}$
with $L^{1}$-BAPU $\Phi=\left(\varphi_{i}\right)_{i\in I}$. Then,
for each $p\in\left[1,\infty\right]$ and $k\in\N_{0}$, the map
\[
\Gamma_{p}^{\left(k\right)}:L^{p}\left(\R^{\dimension}\right)\to\ell^{\UpperExpo p}\!\!\left(I;\,L^{p}\left(\R^{\dimension}\right)\right),g\mapsto\left(\Fourier^{-1}\left(\varphi_{i}^{k\ast}\cdot\widehat{g}\right)\right)_{i\in I}
\]
is well-defined and bounded with
\[
\vertiii{\Gamma_{p}^{\left(k\right)}}\leq C_{\CalQ,\Phi,1}\cdot N_{\CalQ}^{2k+1}.
\]
Here, $\UpperExpo p=\max\left\{ p,p'\right\} $.
\end{lem}

\begin{proof}
As in the proof of Lemma~\ref{lem:BasicEstimateFineIntoCoarse},
it suffices by complex interpolation (precisely, by \cite[Theorems 5.1.1 and 5.6.3]{BerghLoefstroemInterpolationSpaces})
to establish the claim for $p\in\left\{ 1,2,\infty\right\} $. Furthermore,
Lemma~\ref{lem:SemiStructuredClusterInvariant} shows
\[
\left\Vert \Fourier^{-1}\varphi_{i}^{k\ast}\right\Vert _{L^{1}}\leq\sum_{\ell\in i^{k\ast}}\left\Vert \Fourier^{-1}\varphi_{\ell}\right\Vert _{L^{1}}\leq\left|i^{k\ast}\right|\cdot C_{\CalQ,\Phi,1}\leq N_{\CalQ}^{k}C_{\CalQ,\Phi,1}\,,
\]
so that Young's inequality implies for arbitrary $i\in I$, $p\in\left[1,\infty\right]$
and $g\in L^{p}\left(\R^{\dimension}\right)$ that
\[
\left\Vert \Fourier^{-1}\left(\varphi_{i}^{k\ast}\cdot\widehat{g}\right)\right\Vert _{L^{p}}\leq\left\Vert \Fourier^{-1}\varphi_{i}^{k\ast}\right\Vert _{L^{1}}\cdot\left\Vert g\right\Vert _{L^{p}}\leq N_{\CalQ}^{k}C_{\CalQ,\Phi,1}\cdot\left\Vert g\right\Vert _{L^{p}}\leq N_{\CalQ}^{2k+1}C_{\CalQ,\Phi,1}\cdot\left\Vert g\right\Vert _{L^{p}}.
\]
Since we have $\UpperExpo p=\infty$ for $p\in\left\{ 1,\infty\right\} $,
this yields the claim for $p=1$ and $p=\infty$.

\medskip{}

It remains to consider $p=2$. Here, we have $\UpperExpo p=2$. In
view of Plancherel's theorem, this implies for arbitrary finite sets
$I_{0}\subset I$ that
\[
\sum_{i\in I_{0}}\left\Vert \Fourier^{-1}\left(\varphi_{i}^{k\ast}\cdot\widehat{g}\right)\right\Vert _{L^{2}}^{2}=\sum_{i\in I_{0}}\int_{\R^{\dimension}}\left|\varphi_{i}^{k\ast}\left(\xi\right)\cdot\widehat{g}\left(\xi\right)\right|^{2}\,\d\xi=\int_{\R^{\dimension}}\left|\widehat{g}\left(\xi\right)\right|^{2}\cdot\sum_{i\in I_{0}}\left|\varphi_{i}^{k\ast}\left(\xi\right)\right|^{2}\,\d\xi.
\]
But for $\xi\in\R^{\dimension}\setminus\CalO$, we have $\sum_{i\in I_{0}}\left|\varphi_{i}^{k\ast}\left(\xi\right)\right|^{2}=0$.
If otherwise $\xi\in\CalO$, there is some $i_{\xi}\in I$ satisfying
$\xi\in Q_{i_{\xi}}$. For each $i\in I_{0}$ with $\varphi_{i}^{k\ast}\left(\xi\right)\neq0$,
this implies $\emptyset\neq Q_{i_{\xi}}\cap Q_{i}^{k\ast}$ and hence
$i\in i_{\xi}^{\left(k+1\right)\ast}$. But Lemma~\ref{lem:SemiStructuredClusterInvariant}
and the Hausdorff-Young inequality show
\[
\left\Vert \varphi_{i}^{k\ast}\right\Vert _{L^{\infty}}\leq\sum_{\ell\in i^{k\ast}}\left\Vert \varphi_{\ell}\right\Vert _{L^{\infty}}\leq\sum_{\ell\in i^{k\ast}}\left\Vert \Fourier^{-1}\varphi_{\ell}\right\Vert _{L^{1}}\leq\left|i^{k\ast}\right|C_{\CalQ,\Phi,1}\leq C_{\CalQ,\Phi,1}\cdot N_{\CalQ}^{k}\,,
\]
so that we get
\[
\sum_{i\in I_{0}}\left|\varphi_{i}^{k\ast}\left(\xi\right)\right|^{2}\leq\sum_{i\in i_{\xi}^{\left(k+1\right)\ast}}\left\Vert \varphi_{i}^{k\ast}\right\Vert _{L^{\infty}}^{2}\leq\left|i_{\xi}^{\left(k+1\right)\ast}\right|\cdot\left(C_{\CalQ,\Phi,1}N_{\CalQ}^{k}\right)^{2}\leq C_{\CalQ,\Phi,1}^{2}N_{\CalQ}^{3k+2}.
\]
All in all, we have shown
\begin{align*}
\sum_{i\in I_{0}}\left\Vert \Fourier^{-1}\left(\varphi_{i}^{k\ast}\cdot\widehat{g}\right)\right\Vert _{L^{2}}^{2} & \leq\int_{\R^{\dimension}}\left|\widehat{g}\left(\xi\right)\right|^{2}\cdot\sum_{i\in I_{0}}\left|\varphi_{i}^{k\ast}\left(\xi\right)\right|^{2}\,\d\xi\\
 & \leq\left(C_{\CalQ,\Phi,1}N_{\CalQ}^{2k+1}\cdot\left\Vert \widehat{g}\right\Vert _{L^{2}}\right)^{2}\\
\left({\scriptstyle \text{Plancherel}}\right) & =\left(C_{\CalQ,\Phi,1}N_{\CalQ}^{2k+1}\cdot\left\Vert g\right\Vert _{L^{2}}\right)^{2}.
\end{align*}

Since this holds for every finite subset $I_{0}\subset I$, we finally
arrive at
\[
\left\Vert \Gamma_{2}^{\left(k\right)}\left(g\right)\right\Vert _{\ell^{\UpperExpo 2}\left(I;\,\,L^{2}\left(\R^{\dimension}\right)\right)}=\left(\,\smash{\sum_{i\in I}}\,\vphantom{\sum}\left\Vert \Fourier^{-1}\left(\varphi_{i}^{k\ast}\cdot\widehat{g}\right)\right\Vert _{L^{2}}^{2}\,\right)^{1/2}\vphantom{\sum_{i\in I}}\leq C_{\CalQ,\Phi,1}N_{\CalQ}^{2k+1}\cdot\left\Vert g\right\Vert _{L^{2}}
\]
for all $g\in L^{2}\left(\R^{\dimension}\right)$. As seen above,
this completes the proof.
\end{proof}
We want to allow embeddings of the form $\FourierDecompSp{\CalQ}{p_{1}}Y\hookrightarrow\FourierDecompSp{\CalP}{p_{2}}Z$,
i.e.\@ the possibility $p_{1}\neq p_{2}$. To this end, the following
lemma is crucial.
\begin{lem}
\label{lem:LocalEmbeddingInHigherLpSpaces}Let $\CalQ=\left(Q_{i}\right)_{i\in I}=\left(T_{i}Q_{i}'+b_{i}\right)_{i\in I}$
be a semi-structured admissible covering of the open set $\emptyset\neq\CalO\subset\R^{\dimension}$.

Let $k\in\N_{0}$ and $p_{0},p_{1},p_{2}\in\left(0,\infty\right]$
with $p_{0}\leq p_{1}\leq p_{2}$. Then there is a constant $C=C\left(\dimension,k,p_{0},\CalQ\right)$
(independent of $p_{1},p_{2}$) such that
\[
\left\Vert \Fourier^{-1}\left(\gamma_{i}\,f\right)\right\Vert _{L^{p_{2}}}\leq C\cdot\left|\det T_{i}\right|^{\frac{1}{p_{1}}-\frac{1}{p_{2}}}\cdot\left\Vert \Fourier^{-1}\left(\gamma_{i}\,f\right)\right\Vert _{L^{p_{1}}}
\]
holds for each distribution $f\in\DistributionSpace{\CalO}$, all
$i\in I$ and each $\gamma_{i}\in\TestFunctionSpace{\CalO}$ with
$\gamma_{i}\equiv0$ on $\CalO\setminus Q_{i}^{k\ast}$.
\end{lem}

\begin{proof}
By Lemma~\ref{lem:SemiStructuredNormalizationNeighboring}, there
is some $R=R\left(R_{\CalQ},C_{\CalQ},k\right)>0$ such that $Q_{i}^{k\ast}\subset T_{i}\left(\overline{B_{R}}\left(0\right)\right)+b_{i}$
holds for all $i\in I$. By enlarging $R$, we can furthermore assume
that $L:=\lambda\left(B_{1}\left(0\right)\right)\cdot R^{d}\geq1$.
But note that this necessitates $R=R\left(R_{\CalQ},C_{\CalQ},k,\dimension\right)$.

Fix a function $\eta\in\TestFunctionSpace{\R^{\dimension}}$ with
$\eta|_{\overline{B_{R}}\left(0\right)}\equiv1$. Let $f\in\DistributionSpace{\CalO}$
and $i\in I$ with $\Fourier^{-1}\left(\gamma_{i}\,f\right)\in L^{p_{1}}\left(\R^{\dimension}\right)$
(otherwise, there is nothing to show).

Let us first consider the case $p_{1}\in\left(0,1\right)$. Here,
Corollary~\ref{cor:BandlimitedEmbedding} implies $\Fourier^{-1}\left(\gamma_{i}\,f\right)\in L^{p_{2}}\left(\R^{\dimension}\right)$
with
\begin{align*}
\left\Vert \Fourier^{-1}\left(\gamma_{i}\,f\right)\right\Vert _{L^{p_{2}}} & \leq\left[\lambda\left(T_{i}\left(\overline{B_{R}}\left(0\right)\right)+b_{i}\right)\right]^{\frac{1}{p_{1}}-\frac{1}{p_{2}}}\cdot\left\Vert \Fourier^{-1}\left(\gamma_{i}\,f\right)\right\Vert _{L^{p_{1}}}\\
 & =L^{p_{1}^{-1}-p_{2}^{-1}}\cdot\left|\det T_{i}\right|^{\frac{1}{p_{1}}-\frac{1}{p_{2}}}\cdot\left\Vert \Fourier^{-1}\left(\gamma_{i}\,f\right)\right\Vert _{L^{p_{1}}}\\
 & \leq L^{p_{0}^{-1}}\cdot\left|\det T_{i}\right|^{\frac{1}{p_{1}}-\frac{1}{p_{2}}}\cdot\left\Vert \Fourier^{-1}\left(\gamma_{i}\,f\right)\right\Vert _{L^{p_{1}}}.
\end{align*}
Here, we used 
\[
\supp\left(\widehat{\Fourier^{-1}\left(\gamma_{i}\,f\right)}\right)\subset\supp\gamma_{i}\subset\overline{Q_{i}^{k\ast}}\subset T_{i}\left(\overline{B_{R}}\left(0\right)\right)+b_{i}
\]
in the first line. In the last line, we used $p_{1}^{-1}-p_{2}^{-1}\leq p_{1}^{-1}\leq p_{0}^{-1}$
and $L\geq1$. This completes the proof in the case $p_{1}\in\left(0,1\right)$,
since $L$ only depends on $\dimension,k$ and $R_{\CalQ},C_{\CalQ}$.

\medskip{}

Now, let us assume $p_{1}\in\left[1,\infty\right]$, which also entails
$p_{2}\in\left[1,\infty\right]$, because of $p_{2}\geq p_{1}$. For
$i\in I$, define
\[
\eta_{i}:\R^{\dimension}\rightarrow\Compl,\xi\mapsto\eta\left(T_{i}^{-1}\left(\xi-b_{i}\right)\right).
\]
With this definition, we have $\eta_{i}\equiv1$ on $Q_{i}^{k\ast}$;
indeed, note that $\xi\in Q_{i}^{k\ast}\subset T_{i}\left(\overline{B_{R}}\left(0\right)\right)+b_{i}$
entails $T_{i}^{-1}\left(\xi-b_{i}\right)\in\overline{B_{R}}\left(0\right)$
and thus $\eta_{i}\left(\xi\right)=1$. Hence, $\eta_{i}\gamma_{i}=\gamma_{i}$
and furthermore
\[
\Fourier^{-1}\left(\gamma_{i}\,f\right)=\Fourier^{-1}\left(\eta_{i}\cdot\gamma_{i}\,f\right)=\left(\smash{\Fourier^{-1}}\eta_{i}\right)\ast\Fourier^{-1}\left(\gamma_{i}\,f\right).
\]

Because of $1\leq p_{1}\leq p_{2}$, we see $-1\leq-\frac{1}{p_{1}}\leq\frac{1}{p_{2}}-\frac{1}{p_{1}}\leq0$
and hence $1+\frac{1}{p_{2}}-\frac{1}{p_{1}}\in\left[0,1\right]$.
We can thus define $q\in\left[1,\infty\right]$ by $\frac{1}{q}=1+\frac{1}{p_{2}}-\frac{1}{p_{1}}$
with the understanding of $q=\infty$ in case of $1+\frac{1}{p_{2}}-\frac{1}{p_{1}}=0$.
By the general form of Young's inequality (cf.\@ \cite[Proposition 8.9(a)]{FollandRA}),
we derive
\[
\left\Vert \Fourier^{-1}\left(\gamma_{i}\,f\right)\right\Vert _{L^{p_{2}}}=\left\Vert \left(\Fourier^{-1}\eta_{i}\right)\ast\Fourier^{-1}\left(\gamma_{i}\,f\right)\right\Vert _{L^{p_{2}}}\leq\left\Vert \smash{\Fourier^{-1}}\eta_{i}\right\Vert _{L^{q}}\cdot\left\Vert \smash{\Fourier^{-1}}\left(\gamma_{i}\,f\right)\right\Vert _{L^{p_{1}}}<\infty\,,
\]
where we used $\eta_{i}\in\TestFunctionSpace{\R^{\dimension}}$ and
thus $\Fourier^{-1}\eta_{i}\in\Schwartz\left(\R^{\dimension}\right)\hookrightarrow L^{q}\left(\R^{\dimension}\right)$.

Since $\eta_{i}=L_{b_{i}}\left(\smash{\eta\circ T_{i}^{-1}}\right)$,
a straightforward calculation using elementary properties of the Fourier
transform (cf.\@ \cite[Theorem 8.22]{FollandRA}) yields
\[
\Fourier^{-1}\eta_{i}=\left|\det T_{i}\right|\cdot M_{b_{i}}\left(D_{T_{i}}\left[\smash{\Fourier^{-1}}\eta\right]\right),
\]
where we recall the notation $D_{h}f=f\circ h^{T}$. We conclude
\begin{align*}
\left\Vert \smash{\Fourier^{-1}}\eta_{i}\right\Vert _{L^{q}}=\left|\det T_{i}\right|\cdot\left\Vert D_{T_{i}}\left[\smash{\Fourier^{-1}}\eta\right]\right\Vert _{L^{q}} & =\left|\det T_{i}\right|^{1-\frac{1}{q}}\cdot\left\Vert \smash{\Fourier^{-1}}\eta\right\Vert _{L^{q}}\\
 & =\left|\det T_{i}\right|^{\frac{1}{p_{1}}-\frac{1}{p_{2}}}\cdot\left\Vert \smash{\Fourier^{-1}}\eta\right\Vert _{L^{q}}\\
 & \leq\max\left\{ \left\Vert \smash{\Fourier^{-1}}\eta\right\Vert _{L^{1}},\left\Vert \smash{\Fourier^{-1}}\eta\right\Vert _{L^{\infty}}\right\} \cdot\left|\det T_{i}\right|^{\frac{1}{p_{1}}-\frac{1}{p_{2}}}\,.
\end{align*}
Here, the last step used the estimate $\left\Vert f\right\Vert _{L^{q}}\leq\max\left\{ \left\Vert f\right\Vert _{L^{1}},\left\Vert f\right\Vert _{L^{\infty}}\right\} $
which is valid for all measurable $f:\R^{\dimension}\to\Compl$ and
$q\in\left[1,\infty\right]$, cf.\@ \cite[Proposition 6.10]{FollandRA}.
Further, recall that $\eta$ only depends on $\dimension$ and on
$R=R\left(\dimension,k,R_{\CalQ},C_{\CalQ}\right)$, so that $C=\max\left\{ \left\Vert \smash{\Fourier^{-1}}\eta\right\Vert _{L^{1}},\left\Vert \smash{\Fourier^{-1}}\eta\right\Vert _{L^{\infty}}\right\} $
is of the form $C=C\left(\dimension,k,p_{0},\CalQ\right)$, as desired.
Actually, $C$ does not depend on $p_{0}$ in this case.
\end{proof}
As a corollary, we see that every $L^{p}$-BAPU is also an $L^{q}$-BAPU
for all $q\geq p$.
\begin{cor}
\label{cor:LpBAPUsAreAlsoLqBAPUsForLargerq}Let $\CalQ=\left(Q_{i}\right)_{i\in I}=\left(T_{i}Q_{i}'+b_{i}\right)_{i\in I}$
be a semi-structured covering of $\emptyset\neq\CalO\subset\R^{\dimension}$,
let $p\in\left(0,\infty\right]$, and let $\Phi=\left(\varphi_{i}\right)_{i\in I}$
be an $L^{p}$-BAPU for $\CalQ$. Then, for every $q\in\left[p,\infty\right]$,
$\Phi$ is also an $L^{q}$-BAPU for $\CalQ$ with
\[
C_{\CalQ,\Phi,q}\leq C\cdot C_{\CalQ,\Phi,p}\,,\quad\text{for some constant }C=C\left(\dimension,\CalQ,p\right)\,.\qedhere
\]
\end{cor}

\begin{proof}
In case of $p\geq1$ (and $q\geq p$), the definition for an $L^{q}$-BAPU
is the same as that for an $L^{p}$-BAPU and $C_{\CalQ,\Phi,q}=C_{\CalQ,\Phi,p}$.
Thus, we can concentrate on the case $p\in\left(0,1\right)$.

Lemma~\ref{lem:LocalEmbeddingInHigherLpSpaces} (with $p_{0}=p=p_{1}$,
$p_{2}=r$, $k=0$, $f\equiv1$ and $\gamma_{i}=\varphi_{i}$) yields
a constant $C=C\left(\dimension,\CalQ,p\right)$ such that
\begin{align}
\left\Vert \Fourier^{-1}\varphi_{i}\right\Vert _{L^{r}}=\left\Vert \Fourier^{-1}\left(\varphi_{i}\cdot1\right)\right\Vert _{L^{r}} & \leq C\cdot\left|\det T_{i}\right|^{\frac{1}{p}-\frac{1}{r}}\cdot\left\Vert \Fourier^{-1}\left(\varphi_{i}\cdot1\right)\right\Vert _{L^{p}}\nonumber \\
 & \leq C\cdot\left|\det T_{i}\right|^{\frac{1}{p}-\frac{1}{r}}\cdot\left|\det T_{i}\right|^{1-\frac{1}{p}}\cdot C_{\CalQ,\Phi,p}\nonumber \\
 & =CC_{\CalQ,\Phi,p}\cdot\left|\det T_{i}\right|^{1-\frac{1}{r}}\label{eq:LpBAPUExponentChange}
\end{align}
holds for all $i\in I$ and $r\geq p$. For $q\in\left[p,1\right)$,
this implies (for $r=q\geq p$) that
\[
C_{\CalQ,\Phi,q}=\sup_{i\in I}\left|\det T_{i}\right|^{\frac{1}{q}-1}\left\Vert \Fourier^{-1}\varphi_{i}\right\Vert _{L^{q}}\leq C\cdot C_{\CalQ,\Phi,p}<\infty,
\]
as desired. Finally, for $q\in\left[1,\infty\right]$, we apply equation~(\ref{eq:LpBAPUExponentChange})
with $r=1\geq p$ to get
\[
C_{\CalQ,\Phi,q}=\sup_{i\in I}\left\Vert \Fourier^{-1}\varphi_{i}\right\Vert _{L^{1}}\leq C\cdot C_{\CalQ,\Phi,p}<\infty\,.\qedhere
\]
\end{proof}
Before we can state and prove our first sufficient criterion for embeddings
between decomposition spaces, we need one final technical lemma.
\begin{lem}
\label{lem:NestedNormOrderingInterchange}Let $I,J\neq\emptyset$
be two index-sets and let $p,q\in\left(0,\infty\right]$ with $p\leq q$.
Then, for an arbitrary sequence $\left(x_{i,j}\right)_{\left(i,j\right)\in I\times J}\in\Compl^{I\times J}$,
we have
\[
\left\Vert \left(\left\Vert \left(x_{i,j}\right)_{i\in I}\right\Vert _{\ell^{p}\left(I\right)}\right)_{j\in J}\right\Vert _{\ell^{q}\left(J\right)}\leq\left\Vert \left(\left\Vert \left(x_{i,j}\right)_{j\in J}\right\Vert _{\ell^{q}\left(J\right)}\right)_{i\in I}\right\Vert _{\ell^{p}\left(I\right)}\,.\qedhere
\]
\end{lem}

\begin{proof}
Let us first assume $q<\infty$, which also implies $p<\infty$, since
$p\leq q$. Then we get as desired
\[
\begin{alignedat}{2}\left\Vert \left(\left\Vert \left(x_{i,j}\right)_{i\in I}\right\Vert _{\ell^{p}}\right)_{j\in J}\right\Vert _{\ell^{q}} & =\vphantom{\sum_{j\in J}}\left[\,\smash{\sum_{j\in J}}\vphantom{\sum}\left(\,\smash{\sum_{i\in I}}\,\vphantom{\sum}\left|x_{i,j}\right|^{p}\,\right)^{q/p}\,\right]^{1/q} &  & =\left\Vert \left(\,\smash{\sum_{i\in I}}\,\vphantom{\sum}\left|x_{i,j}\right|^{p}\,\right)_{j\in J}\right\Vert _{\ell^{q/p}}^{1/p}\\
\left({\scriptstyle \text{triangle ineq. for }\ell^{q/p}\left(J\right),\text{ since }q/p\geq1}\right) & \leq\vphantom{\sum_{i\in I}}\left[\,\smash{\sum_{i\in I}}\,\vphantom{\sum}\left\Vert \left(\left|x_{i,j}\right|^{p}\right)_{j\in J}\right\Vert _{\ell^{q/p}}\,\right]^{1/p} &  & =\left[\,\smash{\sum_{i\in I}}\,\vphantom{\sum}\left(\,\vphantom{\sum}\smash{\sum_{j\in J}}\,\left|x_{i,j}\right|^{q}\,\right)^{p/q}\,\right]^{1/p}\vphantom{\sum_{i\in I}}\\
 & =\left\Vert \left(\left\Vert \left(x_{i,j}\right)_{j\in J}\right\Vert _{\ell^{q}}\right)_{i\in I}\right\Vert _{\ell^{p}}\:.
\end{alignedat}
\]

Now, we consider the case $q=\infty$. To this end, let $j\in J$
be arbitrary. Then, by solidity of $\ell^{p}$,
\[
\left\Vert \left(x_{i,j}\right)_{i\in I}\right\Vert _{\ell^{p}}=\left\Vert \left(\left|x_{i,j}\right|\right)_{i\in I}\right\Vert _{\ell^{p}}\leq\left\Vert \left(\:\smash{\sup_{j\in J}}\,\left|x_{i,j}\right|\vphantom{T^{i}}\right)_{i\in I}\right\Vert _{\ell^{p}}\quad\overset{q=\infty}{=}\quad\left\Vert \left(\left\Vert \left(x_{i,j}\right)_{j\in J}\right\Vert _{\ell^{q}}\right)_{i\in I}\right\Vert _{\ell^{p}}.
\]
Since this holds for arbitrary $j\in J$, we get the desired inequality.
\end{proof}
All of our sufficient conditions for the existence of embeddings between
decomposition spaces will be based on the following theorem. Its assumptions
are very general, but usually quite tedious to verify. Thus, in the
remainder of this section, we will derive several more specialized
consequences of this general theorem, whose assumptions can often
be verified more easily.
\begin{thm}
\label{thm:NoSubordinatenessWithConsiderationOfOverlaps}Let $\CalQ=\left(Q_{i}\right)_{i\in I}=\left(T_{i}Q_{i}'+b_{i}\right)_{i\in I}$
and $\CalP=\left(P_{j}\right)_{j\in J}=\left(S_{j}P_{j}'+c_{j}\right)_{j\in J}$
be two semi-structured admissible coverings of the open sets $\CalO\subset\R^{\dimension}$
and $\CalO'\subset\R^{\dimension}$, respectively. Let $Y\subset\Compl^{I}$
and $Z\subset\Compl^{J}$ be two solid sequence spaces which are $\CalQ$-regular
and $\CalP$-regular, respectively.

Let $p_{1},p_{2}\in\left(0,\infty\right]$ with $p_{1}\leq p_{2}$
and assume that $\CalQ$ admits an $L^{p_{1}}$-BAPU $\Phi=\left(\varphi_{i}\right)_{i\in I}$.

Fix $I_{0}\subset I$ and let $K\neq\emptyset$ be an index-set. For
each $k\in K$, let $I^{\left(k\right)}\subset I_{0}$ and $J^{\left(k\right)}\subset J$
be arbitrary, but assume that
\begin{equation}
\vphantom{\bigcup_{i\in I_{0}\setminus I^{\left(k\right)}}}\left(\,\smash{\bigcup_{i\in I_{0}\setminus I^{\left(k\right)}}}\vphantom{\bigcup}Q_{i}\,\right)\cap\left(\,\smash{\bigcup_{j\in J^{\left(k\right)}}}\vphantom{\bigcup}P_{j}\,\right)=\emptyset\qquad\forall\,k\in K\,.\label{eq:SpecialIndexSetConditionPSelected}
\end{equation}
Let $J_{0}\subset J_{00}:=\bigcup_{k\in K}J^{\left(k\right)}$. For
each $k\in K$, choose some $q_{k}\in\left[p_{1},p_{2}\right]$. Let
$q^{\left(0\right)}:=\inf_{k\in K}q_{k}$ and assume that $\CalP$
admits an $L^{q^{\left(0\right)}}$-BAPU $\Psi=\left(\psi_{j}\right)_{j\in J}$.

Let $X\subset\Compl^{K}$ be a solid sequence space on $K$ and define
two weights $v=\left(v_{k,i}\right)_{k\in K,i\in I^{\left(k\right)}}$
and $w=\left(w_{k,j}\right)_{k\in K,j\in J^{\left(k\right)}}$ by
\begin{equation}
w_{k,j}:=\begin{cases}
\left|\det S_{j}\right|^{p_{2}^{-1}-1}, & \text{if }q_{k}<1,\\
\vphantom{\rule{0.1cm}{0.55cm}}\left|\det S_{j}\right|^{p_{2}^{-1}-q_{k}^{-1}}, & \text{if }q_{k}\geq1
\end{cases}\label{eq:WeightForJDefinition}
\end{equation}
and
\begin{equation}
v_{k,i}:=\begin{cases}
\left|\det T_{i}\right|^{p_{1}^{-1}-q_{k}^{-1}}\cdot\left[\sup_{j\in J^{\left(k\right)}}\lambda\left(\,\overline{P_{j}}-\overline{Q_{i}}\,\right)\right]^{q_{k}^{-1}-1}, & \text{if }q_{k}<1,\\
\vphantom{\rule{0.1cm}{0.55cm}}\left|\det T_{i}\right|^{p_{1}^{-1}-q_{k}^{-1}}, & \text{if }q_{k}\geq1
\end{cases}\label{eq:WeightForIDefinition}
\end{equation}
where we implicitly assume (for $k\in K$ with $q_{k}<1$) that $\sup_{j\in J^{\left(k\right)}}\lambda\left(\,\overline{P_{j}}-\overline{Q_{i}}\,\right)<\infty$
for all $i\in I^{\left(k\right)}$.

If the maps
\begin{alignat}{5}
\eta_{1}\!:\, &  & X\!\!\left(\vphantom{\sum}\!\!\left[\vphantom{\ell_{w}^{M}}\smash{\ell_{w}^{\UpperExpo{q_{k}}}}\!\!\left(\smash{J^{\left(k\right)}}\right)\right]_{k\in K}\!\right) & \hookrightarrow Z, &  & \left(x_{j}\right)_{j\in J_{00}} & \mapsto & \left(x_{j}\right)_{j\in J} & \,\,\, & \text{with }x_{j}=0\text{ for }j\in J\!\setminus\!J_{00}\label{eq:AssumedDiscreteEmbedding1}\\
\eta_{2}\!:\, &  & Y & \to X\!\!\left(\vphantom{\sum}\!\left[\vphantom{\ell_{w}^{M}}\smash{\ell_{v}^{\LowerExpo{q_{k}}}}\!\!\!\left(\smash{I^{\left(k\right)}}\right)\right]_{k\in K}\right), &  & \left(x_{i}\right)_{i\in I} & \mapsto & \left(x_{i}\right)_{i\in L} & \,\,\, & \text{with }L:=\bigcup_{k\in K}\!I^{\left(k\right)}\subset I\label{eq:AssumedDiscreteEmbedding2}
\end{alignat}
are well-defined and bounded, the map
\[
\iota:\FourierDecompSp{\CalQ}{p_{1}}Y\to\FourierDecompSp{\CalP}{p_{2}}Z,f\mapsto\sum_{\left(i,j\right)\in I_{0}\times J_{0}}\psi_{j}\,\varphi_{i}\,f\:,
\]
is well-defined and bounded with $\vertiii{\iota}\leq C\cdot\vertiii{\eta_{1}}\cdot\vertiii{\eta_{2}}$
for a constant 
\[
C=C\left(\dimension,p_{1},p_{2},q^{\left(0\right)},\CalQ,\CalP,C_{\CalQ,\Phi,p_{1}},C_{\CalP,\Psi,q^{\left(0\right)}},\vertiii{\Gamma_{\CalP}}_{Z\to Z}\right).
\]

More precisely, we have (by definition)
\[
\left\langle \iota f,\,g\right\rangle _{\CalD'}=\sum_{\left(i,j\right)\in I_{0}\times J_{0}}\left\langle f,\,\varphi_{i}\,\psi_{j}\,g\right\rangle _{\CalD'}\qquad\forall\,g\in\TestFunctionSpace{\CalO'},
\]
with absolute convergence of the series for all $g\in\TestFunctionSpace{\CalO'}$.
\end{thm}

\begin{rem*}
A few remarks are in order:

\begin{enumerate}[leftmargin=0.7cm]
\item The most common case is $J_{0}=J$. In this case, note that since
$\left(\psi_{j}\right)_{j\in J}$ is a (locally finite) partition
of unity on $\CalO'$ (see Lemma~\ref{lem:PartitionCoveringNecessary}),
we have $g=\sum_{j\in J}\psi_{j}\,g=\sum_{j\in J_{0}}\psi_{j}\,g$
for arbitrary $g\in\TestFunctionSpace{\CalO'}$, where only finitely
many terms of the sum do not vanish identically. Hence,
\[
\left\langle \iota f,\,g\right\rangle _{\CalD'}=\sum_{\left(i,j\right)\in I_{0}\times J_{0}}\left\langle f,\,\varphi_{i}\,\psi_{j}\,g\right\rangle _{\CalD'}=\sum_{i\in I_{0}}\left\langle f,\,\varphi_{i}\,g\right\rangle _{\CalD'}.
\]
If furthermore $I_{0}=I$, a similar argument shows for $g\in\TestFunctionSpace{\CalO\cap\CalO'}$
that $\left\langle \iota f,\,g\right\rangle _{\CalD'}=\left\langle f,\,g\right\rangle _{\CalD'}$.

Thus, if $I_{0}=I$ and $J_{0}=J$ and if also $\CalO=\CalO'$, then
$\iota f=f$ for all $f\in\FourierDecompSp{\CalQ}{p_{1}}Y\subset\DistributionSpace{\CalO}$,
so that $\iota$ is indeed an (injective) embedding.
\item In the present case and in future theorems, we always use the BAPUs
$\Phi$ and $\Psi$ not only to define the map $\iota$, but also
to compute the (quasi)-norms of the spaces $\FourierDecompSp{\CalQ}{p_{1}}Y$
and $\FourierDecompSp{\CalP}{p_{2}}Z$, respectively. Of course, other
choices lead to equivalent norms (see Corollary~\ref{cor:DecompositionSpaceWellDefined});
but in this case, the constant $C$ from above would also depend on
the BAPUs which are used to calculate the (quasi)-norms on the two
decomposition spaces.
\item There is a slight variation of condition~(\ref{eq:SpecialIndexSetConditionPSelected})
which is sometimes useful: Assume that $I^{\left(k,0\right)}\subset I$
and $J^{\left(k\right)}\subset J$ satisfy
\[
\CalO\cap\bigcup_{j\in J^{\left(k\right)}}P_{j}\subset\bigcup_{i\in I^{\left(k,0\right)}}Q_{i}.
\]
Then the sets $J^{\left(k\right)}$ and $I^{\left(k\right)}:=\left(I^{\left(k,0\right)}\right)^{\ast}$
satisfy condition~(\ref{eq:SpecialIndexSetConditionPSelected}) with
$I_{0}:=I$.

Indeed, if this was false, there would be some $\ell\in I\setminus I^{\left(k\right)}$
and some $h\in J^{\left(k\right)}$ with $Q_{\ell}\cap P_{h}\neq\emptyset$.
Thus, there is some 
\[
\xi\in Q_{\ell}\cap P_{h}\subset\CalO\cap\bigcup_{j\in J^{\left(k\right)}}P_{j}\subset\bigcup_{i\in\smash{I^{\left(k,0\right)}}}Q_{i}\,,
\]
so that $\xi\in Q_{i}$ for some $i\in I^{\left(k,0\right)}$. But
then $\xi\in Q_{i}\cap Q_{\ell}\neq\emptyset$ and hence $\ell\in i^{\ast}\subset\left(I^{\left(k,0\right)}\right)^{\ast}=I^{\left(k\right)}$,
in contradiction to $\ell\in I\setminus I^{\left(k\right)}$.
\item Note that there is no symmetry between $\CalQ,\CalP$ in condition~(\ref{eq:SpecialIndexSetConditionPSelected}),
i.e.\@ the same condition is not ``automatically'' satisfied with
$\CalQ,\CalP$ interchanged. This is natural; as we will see in more
detail in the proof, if we want to estimate $\left\Vert \iota f\right\Vert _{\FourierDecompSp{\CalP}{p_{2}}Z}$,
we have to estimate the individual ``pieces''
\begin{align*}
\qquad\left\Vert \Fourier^{-1}\left(\psi_{j}\cdot\iota f\right)\right\Vert _{L^{p_{2}}}=\vphantom{\sum_{\left(i,\ell\right)\in I_{0}\times J_{0}}}\left\Vert \Fourier^{-1}\left(\psi_{j}\cdot\smash{\sum_{\left(i,\ell\right)\in I_{0}\times J_{0}}}\vphantom{\sum}\varphi_{i}\,\psi_{\ell}\,f\right)\right\Vert _{L^{p_{2}}} & =\vphantom{\sum_{\ell\in J_{0}\cap j^{\ast}}}\left\Vert \Fourier^{-1}\left(\,\vphantom{\sum}\smash{\sum_{\ell\in J_{0}\cap j^{\ast}}}\psi_{j}\,\psi_{\ell}\,\smash{\sum_{i\in I_{0}}}\varphi_{i}\,f\right)\right\Vert _{L^{p_{2}}}\\
 & \lesssim\sum_{\ell\in J_{0}\cap j^{\ast}}\left\Vert \Fourier^{-1}\left(\psi_{\ell}\cdot\smash{\sum_{i\in I_{0}}}\vphantom{\sum}\varphi_{i}\,f\right)\right\Vert _{L^{p_{2}}}.
\end{align*}
At this point, condition~(\ref{eq:SpecialIndexSetConditionPSelected})
is designed to ensure that we have $\psi_{\ell}\cdot\sum_{i\in I_{0}}\varphi_{i}f=\psi_{\ell}\cdot\sum_{i\in I^{\left(k\right)}}\varphi_{i}f$,
as long as $\ell\in J^{\left(k\right)}$ for some fixed $k\in K$.

In other words, given $\ell\in J^{\left(k\right)}$, condition~(\ref{eq:SpecialIndexSetConditionPSelected})
allows us to identify those $i\in I$ for which $\varphi_{i}f$ (potentially)
has an impact on $\psi_{\ell}\sum_{i\in I_{0}}\varphi_{i}f$. The
point here is that we want to estimate $\iota f$ when localized using
$\Psi=\left(\psi_{j}\right)_{j\in J}$, while we are given information
about $f$ localized using $\Phi=\left(\varphi_{i}\right)_{i\in I}$.
Hence, we need information on how the localization using $\Psi$ relates
to $\Phi$. This is exactly what condition~(\ref{eq:SpecialIndexSetConditionPSelected})
achieves.

Very briefly, there is no symmetry in (\ref{eq:SpecialIndexSetConditionPSelected})
regarding $\CalQ,\CalP$, since there is no symmetry in the embedding
$\FourierDecompSp{\CalQ}{p_{1}}Y\hookrightarrow\FourierDecompSp{\CalP}{p_{2}}Z$
regarding $\CalQ,\CalP$.\qedhere

\end{enumerate}
\end{rem*}
\begin{proof}[Proof of Theorem~\ref{thm:NoSubordinatenessWithConsiderationOfOverlaps}]
We divide the proof into $7$ steps.\vspace{0.15cm}

\textbf{Step 1}: Let $f\in\FourierDecompSp{\CalQ}{p_{1}}Y$, define
$c_{i}:=\left\Vert \Fourier^{-1}\left(\varphi_{i}f\right)\right\Vert _{L^{p_{1}}}$
for $i\in I$ and note that $c:=\left(c_{i}\right)_{i\in I}\in Y$
by definition of $\FourierDecompSp{\CalQ}{p_{1}}Y$. Since $\vphantom{\ell_{v}^{\LowerExpo{q_{k}}}}\eta_{2}:Y\to X\Bigl(\left[\smash{\ell_{v}^{\LowerExpo{q_{k}}}}\!\left(\smash{I^{\left(k\right)}}\right)\right]_{k\in K}\Bigr)$
is well-defined, we see in particular that $\left\Vert \left(v_{k,i}\cdot c_{i}\right)_{i\in I^{\left(k\right)}}\right\Vert _{\ell^{\LowerExpo{q_{k}}}}$
is finite for every $k\in K$. But because of $\LowerExpo{q_{k}}\leq2<\infty$,
this shows that there is a countable subset $I^{\left(k,0\right)}\subset I^{\left(k\right)}$
satisfying $\left\Vert \Fourier^{-1}\left(\varphi_{i}f\right)\right\Vert _{L^{p_{1}}}=c_{i}=0$,
and thus $\varphi_{i}f\equiv0$, for all $i\in I^{\left(k\right)}\setminus I^{\left(k,0\right)}$.

Now, an application of Lemma~\ref{lem:LocalEmbeddingInHigherLpSpaces}
shows (for $k\in K$ and $i\in I^{\left(k\right)}$, because of $q_{k}\geq p_{1}$)
that
\[
c_{i}^{\left(k\right)}:=\left\Vert \Fourier^{-1}\left(\varphi_{i}f\right)\right\Vert _{L^{q_{k}}}\leq C_{1}\cdot\left|\det T_{i}\right|^{p_{1}^{-1}-q_{k}^{-1}}\cdot\left\Vert \Fourier^{-1}\left(\varphi_{i}f\right)\right\Vert _{L^{p_{1}}}=C_{1}\cdot\left|\det T_{i}\right|^{p_{1}^{-1}-q_{k}^{-1}}\cdot c_{i}
\]
for some constant $C_{1}=C_{1}\left(\CalQ,\dimension,p_{1}\right)$.
Hence, we get for
\[
u_{k,i}:=\left|\det T_{i}\right|^{q_{k}^{-1}-p_{1}^{-1}}\cdot v_{k,i}=\begin{cases}
\sup_{j\in J^{\left(k\right)}}\left[\lambda\left(\,\overline{P_{j}}-\overline{Q_{i}}\,\right)\right]^{q_{k}^{-1}-1}, & \text{if }q_{k}<1,\\
\vphantom{\rule{0.1cm}{0.55cm}}1, & \text{if }q_{k}\geq1
\end{cases}
\]
that
\[
\left\Vert \left(\smash{c_{i}^{\left(k\right)}}\right)_{i\in I^{\left(k\right)}}\right\Vert _{\ell_{u}^{\LowerExpo{q_{k}}}}\leq C_{1}\cdot\left\Vert \left(v_{k,i}\cdot c_{i}\right)_{i\in I^{\left(k\right)}}\right\Vert _{\ell^{\LowerExpo{q_{k}}}}<\infty\qquad\forall\,k\in K\,.
\]

\medskip{}

\textbf{Step 2}: Fix $k\in K$. For $j\in J^{\left(k\right)}$ and
$i\in I^{\left(k\right)}$, define $\theta_{k,j,i}:=\left\Vert \Fourier^{-1}\left(\psi_{j}\varphi_{i}f\right)\right\Vert _{L^{q_{k}}}$
and
\begin{align}
d_{k,j,i} & :=w_{k,j}\cdot\left|\det S_{j}\right|^{q_{k}^{-1}-p_{2}^{-1}}\cdot\theta_{k,j,i}\nonumber \\
 & =w_{k,j}\cdot\left|\det S_{j}\right|^{q_{k}^{-1}-p_{2}^{-1}}\cdot\left\Vert \Fourier^{-1}\left(\psi_{j}\varphi_{i}f\right)\right\Vert _{L^{q_{k}}}\nonumber \\
 & =\begin{cases}
\left|\det S_{j}\right|^{q_{k}^{-1}-1}\cdot\left\Vert \Fourier^{-1}\left(\psi_{j}\varphi_{i}f\right)\right\Vert _{L^{q_{k}}}=\left|\det S_{j}\right|^{q_{k}^{-1}-1}\cdot\theta_{k,j,i}, & \text{if }q_{k}<1,\\
\vphantom{\rule{0.1cm}{0.55cm}}\left\Vert \Fourier^{-1}\left(\psi_{j}\varphi_{i}f\right)\right\Vert _{L^{q_{k}}}=\theta_{k,j,i}, & \text{if }q_{k}\geq1.
\end{cases}\label{eq:DSequenceDefinition}
\end{align}
As suggested by the definition of $d_{k,j,i}$ and of the weight $u=\left(u_{k,i}\right)_{k\in K,i\in I^{\left(k\right)}}$
from above, we now distinguish two cases.

\begin{casenv}
\item We have $q_{k}\in\left(0,1\right)$. Note that this implies $\UpperExpo{q_{k}}=\infty$.
Now, we note $\supp\varphi_{i}\subset\overline{Q_{i}}$ and $\supp\psi_{j}\subset\overline{P_{j}}$
and use Theorem~\ref{thm:QuasiBanachConvolution} to conclude for
each $j\in J^{\left(k\right)}$ and $i\in I^{\left(k\right)}$ that
\begin{align*}
\qquad d_{k,j,i} & =\left|\det S_{j}\right|^{q_{k}^{-1}-1}\cdot\left\Vert \Fourier^{-1}\left(\psi_{j}\,\varphi_{i}\,f\right)\right\Vert _{L^{q_{k}}}\\
 & \leq\left[\lambda\left(\,\overline{P_{j}}-\overline{Q_{i}}\,\right)\right]^{q_{k}^{-1}-1}\cdot\left|\det S_{j}\right|^{\frac{1}{q_{k}}-1}\cdot\left\Vert \Fourier^{-1}\psi_{j}\right\Vert _{L^{q_{k}}}\cdot\left\Vert \Fourier^{-1}\left(\varphi_{i}\,f\right)\right\Vert _{L^{q_{k}}}\\
 & \leq\left[\lambda\left(\,\overline{P_{j}}-\overline{Q_{i}}\,\right)\right]^{q_{k}^{-1}-1}\cdot C_{\CalP,\Psi,q_{k}}\cdot\left\Vert \Fourier^{-1}\left(\varphi_{i}\,f\right)\right\Vert _{L^{q_{k}}}\\
 & \leq C_{2}\cdot u_{k,i}\cdot\left\Vert \Fourier^{-1}\left(\varphi_{i}\,f\right)\right\Vert _{L^{q_{k}}}<\infty.
\end{align*}
Here, the last step used the definition of $u=\left(u_{k,i}\right)_{k\in K,\,i\in I^{\left(k\right)}}$
(and that $j\in J^{\left(k\right)}$). Furthermore, it was used that
$\Psi$ is an $L^{q^{\left(0\right)}}$-BAPU for $\CalP$ and that
$q_{k}\geq q^{\left(0\right)}$, so that Corollary~\ref{cor:LpBAPUsAreAlsoLqBAPUsForLargerq}
shows that $\Psi$ is also an $L^{q_{k}}$-BAPU for $\CalP$, with
$C_{\CalP,\Psi,q_{k}}\leq C_{2}=C_{2}\left(\dimension,\CalP,q^{\left(0\right)},C_{\CalP,\Psi,q^{\left(0\right)}}\right)$.

Because of $\UpperExpo{q_{k}}=\infty$ and since $j\in J^{\left(k\right)}$
was arbitrary, we finally get
\[
\left\Vert \left(d_{k,j,i}\right)_{j\in J^{\left(k\right)}}\right\Vert _{\ell^{\UpperExpo{q_{k}}}}\leq C_{2}\cdot u_{k,i}\cdot\left\Vert \Fourier^{-1}\left(\varphi_{i}f\right)\right\Vert _{L^{q_{k}}}<\infty\qquad\forall\,i\in I^{\left(k\right)}\,.
\]

\item We have $q_{k}\in\left[1,\infty\right]$. In this case, Lemma~\ref{lem:BasicEstimateCoarseInFine}
shows that for arbitrary $p\in\left[1,\infty\right]$, the map
\[
\Gamma_{p}^{\left(0\right)}:L^{p}\left(\R^{\dimension}\right)\to\ell^{\UpperExpo p}\!\!\left(J;\,L^{p}\left(\R^{\dimension}\right)\right),g\mapsto\left(\Fourier^{-1}\left(\psi_{j}\cdot\widehat{g}\right)\right)_{j\in J}
\]
is a well-defined, bounded, linear operator with $\vertiii{\smash{\Gamma_{p}^{\left(0\right)}}}\leq N_{\CalP}\cdot C_{\CalP,\Psi,1}=:C_{3}$.

We now use this with $p=q_{k}\in\left[1,\infty\right]$ for $g=\Fourier^{-1}\left(\varphi_{i}f\right)\in L^{q_{k}}\left(\R^{\dimension}\right)$,
to conclude for each $i\in I^{\left(k\right)}$ that
\begin{align*}
\qquad\qquad\left\Vert \left(d_{k,j,i}\right)_{j\in J^{\left(k\right)}}\right\Vert _{\ell^{\UpperExpo{q_{k}}}}=\left\Vert \left(\left\Vert \Fourier^{-1}\left(\psi_{j}\varphi_{i}f\right)\right\Vert _{L^{q_{k}}}\right)_{j\in J^{\left(k\right)}}\right\Vert _{\ell^{\UpperExpo{q_{k}}}} & \leq\left\Vert \left(\left\Vert \Fourier^{-1}\left(\psi_{j}\varphi_{i}f\right)\right\Vert _{L^{q_{k}}}\right)_{j\in J}\right\Vert _{\ell^{\UpperExpo{q_{k}}}}\\
 & \leq C_{3}\cdot\left\Vert \Fourier^{-1}\left(\varphi_{i}f\right)\right\Vert _{L^{q_{k}}}\\
 & =C_{3}\cdot u_{k,i}\cdot\left\Vert \Fourier^{-1}\left(\varphi_{i}f\right)\right\Vert _{L^{q_{k}}}<\infty.
\end{align*}

\end{casenv}
All in all, we have shown that
\begin{equation}
\left\Vert \left(d_{k,j,i}\right)_{j\in J^{\left(k\right)}}\right\Vert _{\ell^{\UpperExpo{q_{k}}}}\leq C_{4}\cdot u_{k,i}\cdot\left\Vert \Fourier^{-1}\left(\varphi_{i}f\right)\right\Vert _{L^{q_{k}}}<\infty\qquad\forall\,i\in I^{\left(k\right)}\label{eq:UpperConjugateEstimateComplete}
\end{equation}
holds in both cases, for a suitable constant $C_{4}=C_{4}\left(\dimension,\CalP,q^{\left(0\right)},C_{\CalP,\Psi,1},C_{\CalP,\Psi,q^{\left(0\right)}}\right)$.

\medskip{}

\textbf{Step 3}: Now, take the $\ell^{\LowerExpo{q_{k}}}\!\left(I^{\left(k\right)}\right)$-norm
of estimate~(\ref{eq:UpperConjugateEstimateComplete}) to derive
\begin{align*}
\left\Vert \left(\left\Vert \left(d_{k,j,i}\right)_{j\in J^{\left(k\right)}}\right\Vert _{\ell^{\UpperExpo{q_{k}}}}\right)_{i\in I^{\left(k\right)}}\right\Vert _{\ell^{\LowerExpo{q_{k}}}} & \leq C_{4}\cdot\left\Vert \left(u_{k,i}\cdot\left\Vert \Fourier^{-1}\left(\varphi_{i}f\right)\right\Vert _{L^{q_{k}}}\right)_{i\in I^{\left(k\right)}}\right\Vert _{\ell^{\LowerExpo{q_{k}}}}\\
 & =C_{4}\cdot\left\Vert \left(\smash{c_{i}^{\left(k\right)}}\right)_{i\in I^{\left(k\right)}}\right\Vert _{\ell_{u}^{\LowerExpo{q_{k}}}}\\
\left({\scriptstyle \text{Step }1}\right) & \leq C_{1}C_{4}\cdot\left\Vert \left(c_{i}\right)_{i\in I^{\left(k\right)}}\right\Vert _{\ell_{v}^{\LowerExpo{q_{k}}}}<\infty.
\end{align*}
Note that $\LowerExpo{q_{k}}=\min\left\{ q_{k},q_{k}'\right\} \leq\max\left\{ q_{k},q_{k}'\right\} =\UpperExpo{q_{k}}$.
Thus, an application of Lemma~\ref{lem:NestedNormOrderingInterchange}
finally implies
\begin{align}
\left\Vert \left(\left\Vert \left(d_{k,j,i}\right)_{i\in I^{\left(k\right)}}\right\Vert _{\ell^{\LowerExpo{q_{k}}}}\right)_{j\in J^{\left(k\right)}}\right\Vert _{\ell^{\UpperExpo{q_{k}}}} & \leq\left\Vert \left(\left\Vert \left(d_{k,j,i}\right)_{j\in J^{\left(k\right)}}\right\Vert _{\ell^{\UpperExpo{q_{k}}}}\right)_{i\in I^{\left(k\right)}}\right\Vert _{\ell^{\LowerExpo{q_{k}}}}\nonumber \\
 & \leq C_{1}C_{4}\cdot\left\Vert \left(c_{i}\right)_{i\in I^{\left(k\right)}}\right\Vert _{\ell_{v}^{\LowerExpo{q_{k}}}}<\infty\qquad\forall\,k\in K\,.\label{eq:DoubleNestedEstimateByQPieces}
\end{align}

\medskip{}

\textbf{Step 4}: The result from the previous step implies in particular
for each $k\in K$ and all $j\in J^{\left(k\right)}$ that
\begin{align}
\left\Vert \left(\theta_{k,j,i}\right)_{i\in I^{\left(k\right)}}\right\Vert _{\ell^{\LowerExpo{q_{k}}}} & =\left\Vert \left(\left\Vert \Fourier^{-1}\left(\psi_{j}\varphi_{i}f\right)\right\Vert _{L^{q_{k}}}\right)_{i\in I^{\left(k\right)}}\right\Vert _{\ell^{\LowerExpo{q_{k}}}}\nonumber \\
 & =\begin{cases}
\left|\det S_{j}\right|^{1-q_{k}^{-1}}\cdot\left\Vert \left(d_{k,j,i}\right)_{i\in I^{\left(k\right)}}\right\Vert _{\ell^{\LowerExpo{q_{k}}}}=:C_{k,j}<\infty, & \text{if }q_{k}<1,\\
\vphantom{\rule{0.1cm}{0.55cm}}\left\Vert \left(d_{k,j,i}\right)_{i\in I^{\left(k\right)}}\right\Vert _{\ell^{\LowerExpo{q_{k}}}}=:C_{k,j}<\infty, & \text{if }q_{k}\geq1.
\end{cases}\label{eq:AbsoluteConvergenceDiscreteNormFiniteness}
\end{align}
We will derive from this that the mapping
\begin{equation}
f_{k,j}:\Schwartz\left(\R^{\dimension}\right)\to\Compl,g\mapsto\sum_{i\in I^{\left(k\right)}}\left\langle \varphi_{i}\,f,\,\psi_{j}\,g\right\rangle _{\Schwartz'}\label{eq:DoubleLocalizedDistributionDefinition}
\end{equation}
is a well-defined, tempered distribution for all $k\in K$ and each
$j\in J^{\left(k\right)}$, with absolute convergence of the defining
series. Since this is only a qualitative statement, we will suppress
some unimportant constants in this step.

Fix $k\in K$ and $j\in J^{\left(k\right)}$. To prove the absolute
convergence, we first note $\supp\left(\varphi_{i}\,\psi_{j}\,f\right)\subset\overline{Q_{i}}$
and invoke Lemma~\ref{lem:BasicEstimateFineIntoCoarse} to derive
for arbitrary finite subsets $F^{\left(k\right)}\subset I^{\left(k\right)}$
and 
\[
f_{F^{\left(k\right)}}:=\sum_{i\in F^{\left(k\right)}}\varphi_{i}\,\psi_{j}\,f\in\Schwartz'\left(\R^{\dimension}\right)
\]
the estimate
\begin{equation}
\left\Vert \Fourier^{-1}f_{F^{\left(k\right)}}\right\Vert _{L^{q_{k}}}\leq C_{5}\cdot\left\Vert \left(\left\Vert \Fourier^{-1}\left(\varphi_{i}\psi_{j}f\right)\right\Vert _{L^{q_{k}}}\right)_{i\in F^{\left(k\right)}}\right\Vert _{\ell^{\LowerExpo{q_{k}}}}=C_{5}\cdot\left\Vert \left(\theta_{k,j,i}\right)_{i\in F^{\left(k\right)}}\right\Vert _{\ell^{\LowerExpo{q_{k}}}}\label{eq:SobolevLemmaEstimate}
\end{equation}
for some constant $C_{5}=C_{5}\left(N_{\CalQ},C_{\CalQ,\Phi,1}\right)$.
But Corollary~\ref{cor:LpBAPUsAreAlsoLqBAPUsForLargerq} yields $C_{\CalQ,\Phi,1}=C_{\CalQ,\Phi,\infty}\lesssim_{\,\CalQ,p_{1},\dimension}\,C_{\CalQ,\Phi,p_{1}}$,
so that we can choose $C_{5}=C_{5}\left(\CalQ,p_{1},\dimension,C_{\CalQ,\Phi,p_{1}}\right)$.

Now, we claim that there is some $r_{k}\in\left[1,\infty\right]$
such that we have
\begin{equation}
\left\Vert \Fourier^{-1}f_{F^{\left(k\right)}}\right\Vert _{L^{r_{k}}}\lesssim_{\,\CalQ,\Phi,k,j}\left\Vert \left(\theta_{k,j,i}\right)_{i\in F^{\left(k\right)}}\right\Vert _{\ell^{\LowerExpo{q_{k}}}}\label{eq:AbsoluteConvergenceMainEstimate}
\end{equation}
for all finite subsets $F^{\left(k\right)}\subset I^{\left(k\right)}$.
Here, the implied constant \emph{is allowed} to depend on $k\in K$
and on $j\in J^{\left(k\right)}$, but not on $F^{\left(k\right)}\subset I^{\left(k\right)}$.

In case of $q_{k}\in\left[1,\infty\right]$, we can simply take $r_{k}=q_{k}$,
so that we can assume $q_{k}\in\left(0,1\right)$. In this case, note
that we have $\supp f_{F^{\left(k\right)}}\subset\overline{P_{j}}$
(since we multiply with $\psi_{j}$), so that Corollary~\ref{cor:BandlimitedEmbedding}
implies
\begin{align*}
\left\Vert \Fourier^{-1}f_{F^{\left(k\right)}}\right\Vert _{L^{\infty}} & \leq\left[\lambda\left(\,\overline{P_{j}}\,\right)\right]^{1/q_{k}}\cdot\left\Vert \Fourier^{-1}f_{F^{\left(k\right)}}\right\Vert _{L^{q_{k}}}\\
\left({\scriptstyle \text{using eq. }\eqref{eq:SobolevLemmaEstimate}}\right) & \lesssim_{\,\CalQ,\Phi,k,j}\:\left\Vert \left(\theta_{k,j,i}\right)_{i\in F^{\left(k\right)}}\right\Vert _{\ell^{\LowerExpo{q_{k}}}}.
\end{align*}
Hence, we can choose $r_{k}=\infty$ in case of $q_{k}\in\left(0,1\right)$.

Now, we can prove the absolute convergence of the series defining
$\left\langle f_{k,j},\,g\right\rangle _{\Schwartz'}$ for $g\in\Schwartz\left(\R^{\dimension}\right)$.
Indeed, because of $\LowerExpo{q_{k}}\leq2<\infty$ and since $\left\Vert \left(\theta_{k,j,i}\right)_{i\in I^{\left(k\right)}}\right\Vert _{\ell^{\LowerExpo{q_{k}}}}<\infty$,
there is for arbitrary $\varepsilon>0$ some finite subset $F^{\left(k,\varepsilon\right)}\subset I^{\left(k\right)}$
(potentially depending on $j\in J^{\left(k\right)}$) satisfying
\[
\left\Vert \left(\theta_{k,j,i}\right)_{i\in F^{\left(k\right)}}\right\Vert _{\ell^{\LowerExpo{q_{k}}}}\leq\varepsilon\qquad\forall\:\,\,\text{finite subsets }F^{\left(k\right)}\subset I^{\left(k\right)}\setminus F^{\left(k,\varepsilon\right)}.
\]
But in view of estimate~(\ref{eq:AbsoluteConvergenceMainEstimate}),
this implies
\begin{align*}
\left|\vphantom{\sum}\smash{\sum_{i\in F^{\left(k\right)}}}\left\langle \varphi_{i}\,f,\,\psi_{j}\,g\right\rangle _{\Schwartz'}\right|=\left|\left\langle f_{F^{\left(k\right)}},\,g\right\rangle _{\Schwartz'}\right| & =\left|\left\langle \Fourier^{-1}f_{F^{\left(k\right)}},\,\widehat{g}\right\rangle _{\Schwartz'}\right|\\
 & \leq\left\Vert \Fourier^{-1}f_{F^{\left(k\right)}}\right\Vert _{L^{r_{k}}}\cdot\left\Vert \widehat{g}\right\Vert _{L^{r_{k}'}}\\
 & \lesssim_{\,\CalQ,\Phi,k,j}\:\left\Vert \left(\theta_{k,j,i}\right)_{i\in F^{\left(k\right)}}\right\Vert _{\ell^{\LowerExpo{q_{k}}}}\cdot\left\Vert \widehat{g}\right\Vert _{L^{r_{k}'}}\\
 & \leq\varepsilon\cdot\left\Vert \widehat{g}\right\Vert _{L^{r_{k}'}}
\end{align*}
for every finite subset $F^{\left(k\right)}\subset I^{\left(k\right)}\setminus F^{\left(k,\varepsilon\right)}$.
Since $\varepsilon>0$ was arbitrary, this proves the unconditional
(and hence absolute) convergence of the series $\sum_{i\in I}\left\langle \varphi_{i}\,f,\,\psi_{j}\,g\right\rangle _{\Schwartz'}=\left\langle f_{k,j},\,g\right\rangle _{\Schwartz'}$
for each $g\in\Schwartz\left(\R^{\dimension}\right)$.

Furthermore, as seen in Step $1$, there is a countable subset $I^{\left(k,0\right)}\subset I^{\left(k\right)}$
satisfying $\varphi_{i}f\equiv0$ for all $i\in I^{\left(k\right)}\setminus I^{\left(k,0\right)}$.
In particular, this implies $\left\langle \varphi_{i}f,\,\psi_{j}g\right\rangle _{\Schwartz'}=0$
for all $i\in I^{\left(k\right)}\setminus I^{\left(k,0\right)}$.
Thus, if $\left(I^{\left(k,N\right)}\right)_{N\in\N}$ is a nondecreasing
family of finite sets with $I^{\left(k,0\right)}=\bigcup_{N\in\N}I^{\left(k,N\right)}$,
then
\begin{align*}
\left|\left\langle f_{k,j},\,g\right\rangle _{\Schwartz'}\right|=\lim_{N\to\infty}\left|\vphantom{\sum}\smash{\sum_{i\in I^{\left(k,N\right)}}}\left\langle \varphi_{i}\,f,\,\psi_{j}\,g\right\rangle _{\Schwartz'}\right| & =\lim_{N\to\infty}\left|\left\langle f_{I^{\left(k,N\right)}},\,g\right\rangle _{\Schwartz'}\right|\\
 & =\lim_{N\to\infty}\left|\left\langle \Fourier^{-1}f_{I^{\left(k,N\right)}},\,\widehat{g}\right\rangle _{\Schwartz'}\right|\\
 & \leq\lim_{N\to\infty}\left\Vert \Fourier^{-1}f_{I^{\left(k,N\right)}}\right\Vert _{L^{r_{k}}}\cdot\left\Vert \widehat{g}\right\Vert _{L^{r_{k}'}}\\
\left({\scriptstyle \text{equation }\eqref{eq:AbsoluteConvergenceMainEstimate}}\right) & \lesssim_{\,\CalQ,\Phi,k,j}\:\lim_{N\to\infty}\left\Vert \left(\theta_{k,j,i}\right)_{i\in I^{\left(k,N\right)}}\right\Vert _{\ell^{\LowerExpo{q_{k}}}}\cdot\left\Vert \widehat{g}\right\Vert _{L^{r_{k}'}}\\
 & \leq\left\Vert \left(\theta_{k,j,i}\right)_{i\in I^{\left(k\right)}}\right\Vert _{\ell^{\LowerExpo{q_{k}}}}\cdot\left\Vert \widehat{g}\right\Vert _{L^{r_{k}'}}.
\end{align*}
Because of $\left\Vert \left(\theta_{k,j,i}\right)_{i\in I^{\left(k\right)}}\right\Vert _{\ell^{\LowerExpo{q_{k}}}}<\infty$
(see equation~(\ref{eq:AbsoluteConvergenceDiscreteNormFiniteness}))
and since the inclusion $\Schwartz\left(\R^{\dimension}\right)\hookrightarrow L^{r_{k}'}\!\left(\R^{\dimension}\right)$
and the Fourier transform $\Fourier:\Schwartz\left(\R^{\dimension}\right)\to\Schwartz\left(\R^{\dimension}\right)$
are continuous, we conclude $f_{k,j}\in\Schwartz'\!\left(\R^{\dimension}\right)$,
as claimed.

\medskip{}

\textbf{Step 5}: Here, we continue the previous step by showing
\begin{equation}
\left\Vert \Fourier^{-1}f_{k,j}\right\Vert _{L^{q_{k}}}\leq C_{5}\cdot\left\Vert \left(\theta_{k,j,i}\right)_{i\in I^{\left(k\right)}}\right\Vert _{\ell^{\LowerExpo{q_{k}}}}<\infty\qquad\forall\,\,k\in K\text{ and }j\in J^{\left(k\right)}\:.\label{eq:PsiPiecesLqkEstimate}
\end{equation}
To see this, note that (since we multiply by $\psi_{j}$, see equation~(\ref{eq:DoubleLocalizedDistributionDefinition})),
the tempered distribution $f_{k,j}$ has compact support $\supp f_{k,j}\subset\overline{P_{j}}$.
Thus, the Paley-Wiener theorem (cf.\@ \cite[Theorem 7.23]{RudinFA})
implies that $\Fourier^{-1}f_{k,j}$ is given by (integration against)
the smooth function
\begin{align*}
\left(\Fourier^{-1}f_{k,j}\right)\left(x\right)=\left\langle f_{k,j},\,e^{2\pi i\left\langle x,\mybullet\right\rangle }\right\rangle  & =\sum_{i\in I^{\left(k\right)}}\left\langle \varphi_{i}\,f,\,\psi_{j}\,e^{2\pi i\left\langle x,\mybullet\right\rangle }\right\rangle _{\Schwartz'}\\
\left({\scriptstyle \text{since }\varphi_{i}f=0\text{ for }i\in I^{\left(k\right)}\setminus I^{\left(k,0\right)}}\right) & =\sum_{i\in I^{\left(k,0\right)}}\left\langle \varphi_{i}\,f,\,\psi_{j}\,e^{2\pi i\left\langle x,\mybullet\right\rangle }\right\rangle _{\Schwartz'}\\
 & =\lim_{N\to\infty}\sum_{i\in I^{\left(k,N\right)}}\left\langle \varphi_{i}\,f,\,\psi_{j}\,e^{2\pi i\left\langle x,\mybullet\right\rangle }\right\rangle _{\Schwartz'}\\
 & =\lim_{N\to\infty}\left(\Fourier^{-1}f_{I^{\left(k,N\right)}}\right)\left(x\right),
\end{align*}
so that Fatou's lemma yields
\begin{align*}
\left\Vert \Fourier^{-1}f_{k,j}\right\Vert _{L^{q_{k}}} & \leq\liminf_{N\to\infty}\left\Vert \Fourier^{-1}f_{I^{\left(k,N\right)}}\right\Vert _{L^{q_{k}}}\\
\left({\scriptstyle \text{equation }\eqref{eq:SobolevLemmaEstimate}}\right) & \leq C_{5}\cdot\liminf_{N\to\infty}\left\Vert \left(\theta_{k,j,i}\right)_{i\in I^{\left(k,N\right)}}\right\Vert _{\ell^{\LowerExpo{q_{k}}}}\leq C_{5}\cdot\left\Vert \left(\theta_{k,j,i}\right)_{i\in I^{\left(k\right)}}\right\Vert _{\ell^{\LowerExpo{q_{k}}}}<\infty.
\end{align*}
Finiteness of the last term was shown in equation~(\ref{eq:AbsoluteConvergenceDiscreteNormFiniteness})
above.

\medskip{}

\textbf{Step 6}: Now, we show that $\iota f\in\DistributionSpace{\CalO'}$
is a well-defined distribution, with absolute convergence of the defining
series. To this end, let $k\in K$ be arbitrary and fix $j_{0}\in J^{\left(k\right)}$.

Now, for $\ell\in I_{0}\setminus I^{\left(k\right)}$, we have by
assumption (see equation~(\ref{eq:SpecialIndexSetConditionPSelected}))
that
\[
Q_{\ell}\cap P_{j_{0}}\subset\vphantom{\bigcup_{i\in I_{0}\setminus I^{\left(k\right)}}}\left(\,\vphantom{\bigcup}\smash{\bigcup_{i\in I_{0}\setminus I^{\left(k\right)}}}Q_{i}\,\right)\cap\left(\,\vphantom{\bigcup}\smash{\bigcup_{j\in J^{\left(k\right)}}}P_{j}\,\right)=\emptyset\:,
\]
and thus $\varphi_{\ell}\cdot\psi_{j_{0}}\equiv0$ for all $\ell\in I_{0}\setminus I^{\left(k\right)}$.

Thus, for arbitrary $g\in\Schwartz\left(\R^{\dimension}\right)$,
we have $\left\langle \varphi_{i}f,\,\psi_{j_{0}}g\right\rangle _{\Schwartz'}=0$
for all $i\in I_{0}\setminus I^{\left(k\right)}$, so that we get
(simply by dropping vanishing terms)
\[
\left\langle \,\smash{\sum_{i\in I_{0}}}\,\vphantom{\sum}\varphi_{i}\,\psi_{j_{0}}\,f,\,g\,\right\rangle _{\!\!\Schwartz'}=\sum_{i\in I_{0}}\left\langle \varphi_{i}\,f,\,\psi_{j_{0}}\,g\right\rangle _{\Schwartz'}=\sum_{i\in I^{\left(k\right)}}\left\langle \varphi_{i}\,f,\,\psi_{j_{0}}\,g\right\rangle _{\Schwartz'}=\left\langle f_{k,j_{0}}\:,\:g\right\rangle _{\Schwartz'}\;.
\]
In particular, by what we proved above (in Step $4$), the series
$\sum_{i\in I_{0}}\left\langle \varphi_{i}\,f,\,\psi_{j_{0}}\,g\right\rangle _{\Schwartz'}$
converges absolutely and we have
\begin{equation}
\sum_{i\in I_{0}}\varphi_{i}\,\psi_{j_{0}}\,f=f_{k,j_{0}}\qquad\forall\,k\in K\text{ and }j_{0}\in J^{\left(k\right)}.\label{eq:LocalizationIdentity}
\end{equation}
We call this the \textbf{localization identity}. As a crucial observation,
note that for fixed $j_{0}\in J_{00}$, the left-hand side of the
localization identity~(\ref{eq:LocalizationIdentity}) is \emph{independent}
of the choice $k\in K$ (satisfying $j_{0}\in J^{\left(k\right)}$).
Hence, we may (and will) in the following write
\[
f_{j_{0}}:=f_{k,j_{0}}\quad\text{for all }j_{0}\in J_{00}=\bigcup_{k\in K}J^{\left(k\right)}\text{ and arbitrary }k\in K\text{ with }j_{0}\in J^{\left(k\right)}\,.\tag{\ensuremath{\blacklozenge}}
\]

Now, let $M\subset\CalO'$ be an arbitrary compact set.  Since Lemma~\ref{lem:PartitionCoveringNecessary}
shows that $\left(\psi_{j}\right)_{j\in J}$ is a locally finite partition
of unity on $\CalO'$, the set
\[
J_{M}:=\left\{ j\in J_{0}\with M\cap\supp\psi_{j}\neq\emptyset\right\} 
\]
is finite. For $g\in\TestFunctionSpace{\CalO'}$ with $\supp g\subset M$,
this implies $\psi_{j}g\equiv0$ for $j\in J_{0}\setminus J_{M}$.
Next, note that because of $J_{M}\subset J_{0}\subset J_{00}=\bigcup_{k\in K}J^{\left(k\right)}$,
there is for each $j\in J_{M}$ some $k_{j}\in K$ with $j\in J^{\left(k_{j}\right)}$.
Hence, simply by dropping vanishing terms, we see
\[
\sum_{\left(i,j\right)\in I_{0}\times J_{0}}\left|\left\langle f,\,\varphi_{i}\,\psi_{j}\,g\right\rangle _{\CalD'}\right|=\sum_{\left(i,j\right)\in I_{0}\times J_{M}}\left|\left\langle f,\,\varphi_{i}\,\psi_{j}\,g\right\rangle _{\CalD'}\right|=\sum_{j\in J_{M}}\;\sum_{i\in I_{0}}\left|\left\langle f,\,\varphi_{i}\,\psi_{j}\,g\right\rangle _{\CalD'}\right|<\infty,
\]
since $J_{M}$ is finite and since each of the series $\sum_{i\in I_{0}}\left\langle f,\,\varphi_{i}\,\psi_{j}\,g\right\rangle _{\CalD'}=\left\langle f_{k_{j},j}\:,\:g\right\rangle _{\CalD'}$
for $j\in J_{M}$ converges absolutely, as seen above.

Entirely the same calculation, but without the absolute value, shows
\begin{align*}
\left\langle \iota f,\,g\right\rangle _{\CalD'}=\sum_{\left(i,j\right)\in I_{0}\times J_{0}}\left\langle f,\,\varphi_{i}\,\psi_{j}\,g\right\rangle _{\CalD'}=\sum_{j\in J_{M}}\;\sum_{i\in I_{0}}\left\langle f,\,\varphi_{i}\,\psi_{j}\,g\right\rangle _{\CalD'} & =\sum_{j\in J_{M}}\left\langle f_{k_{j},j}\:,\:g\right\rangle _{\Schwartz'}\\
 & =\left\langle \,\smash{\sum_{j\in J_{M}}}\,\vphantom{\sum}f_{k_{j},j}\:,\:g\,\right\rangle _{\Schwartz'}\vphantom{\sum_{j\in J_{M}}}.
\end{align*}
Here, the right-hand side $\sum_{j\in J_{M}}f_{k_{j},j}$ is a (tempered)
distribution, since it is a finite sum of tempered distributions.
In particular, we see that $g\mapsto\left\langle \iota f,\,g\right\rangle _{\CalD'}$
is continuous when restricted to
\[
C_{M}^{\infty}\left(\CalO'\right)=\left\{ g\in\TestFunctionSpace{\CalO'}\with\supp g\subset M\right\} ,
\]
for \emph{every} compact subset $M\subset\CalO'$. But by \cite[Theorem 6.6]{RudinFA},
this already implies $\iota f\in\DistributionSpace{\CalO'}$, as desired.

\medskip{}

\textbf{Step 7}: It remains to show $\iota f\in\FourierDecompSp{\CalP}{p_{2}}Z$,
together with an accompanying (quasi)-norm estimate. As we will see
below, with $f_{j}$ as in equation~($\blacklozenge$), the sequence
$\varrho=\left(\varrho_{j}\right)_{j\in J}$ given by
\[
\varrho_{j}:=\begin{cases}
\left\Vert \Fourier^{-1}f_{j}\right\Vert _{L^{p_{2}}}, & \text{if }j\in J_{0}\subset J_{00},\\
0, & \text{if }j\notin J_{0}
\end{cases}
\]
will be useful to us. First, we will derive a bound for $\varrho_{j}=\left\Vert \Fourier^{-1}f_{j}\right\Vert _{L^{p_{2}}}$,
for $j\in J^{\left(k\right)}\cap J_{0}$ (and arbitrary $k\in K$).
Note that we have $f_{j}=f_{k,j}$.

Note that $\psi_{j}^{\ast}\,\psi_{j}=\psi_{j}$, since $\psi_{j}^{\ast}\equiv1$
on $P_{j}$, see Lemma~\ref{lem:PartitionCoveringNecessary}. Hence,
\[
\left\langle f_{k,j},\,\psi_{j}^{\ast}\,g\right\rangle _{\Schwartz'}=\sum_{i\in I^{\left(k\right)}}\left\langle \varphi_{i}\,f,\,\psi_{j}\,\psi_{j}^{\ast}\,g\right\rangle _{\Schwartz'}=\sum_{i\in I^{\left(k\right)}}\left\langle \varphi_{i}\,f,\,\psi_{j}\,g\right\rangle _{\Schwartz'}=\left\langle f_{k,j},\,g\right\rangle _{\Schwartz'}\qquad\forall\,g\in\Schwartz\left(\R^{\dimension}\right),
\]
i.e.\@ $f_{k,j}=\psi_{j}^{\ast}\cdot f_{k,j}$. Hence, an application
of Lemma~\ref{lem:LocalEmbeddingInHigherLpSpaces} (with $k=1$
and $p_{0}=q^{\left(0\right)}$) yields—because of $p_{2}\geq q_{k}\geq q^{\left(0\right)}$—a
constant $C_{6}=C_{6}\left(\dimension,\CalP,q^{\left(0\right)}\right)$
which satisfies
\begin{align}
\varrho_{j} & =\left\Vert \Fourier^{-1}f_{k,j}\right\Vert _{L^{p_{2}}}=\left\Vert \Fourier^{-1}\left(\psi_{j}^{\ast}\cdot f_{k,j}\right)\right\Vert _{L^{p_{2}}}\nonumber \\
 & \leq C_{6}\cdot\left|\det S_{j}\right|^{q_{k}^{-1}-p_{2}^{-1}}\cdot\left\Vert \Fourier^{-1}\left(\psi_{j}^{\ast}\cdot f_{k,j}\right)\right\Vert _{L^{q_{k}}}\nonumber \\
 & =C_{6}\cdot\left|\det S_{j}\right|^{q_{k}^{-1}-p_{2}^{-1}}\cdot\left\Vert \Fourier^{-1}f_{k,j}\right\Vert _{L^{q_{k}}}\nonumber \\
\left({\scriptstyle \text{eq. }\eqref{eq:PsiPiecesLqkEstimate}}\right) & \leq C_{5}C_{6}\cdot\left|\det S_{j}\right|^{q_{k}^{-1}-p_{2}^{-1}}\cdot\left\Vert \left(\theta_{k,j,i}\right)_{i\in I^{\left(k\right)}}\right\Vert _{\ell^{\LowerExpo{q_{k}}}}<\infty\qquad\forall\,k\in K\text{ and }j\in J_{0}\cap J^{\left(k\right)}.\label{eq:VarrhoJEstimate}
\end{align}

Now, let $j\in J$ be arbitrary. Recall that $\psi_{j}\psi_{j_{0}}\not\equiv0$
implies $j_{0}\in j^{\ast}$. Thus, using the definition of $\iota f$
and the localization identity~(\ref{eq:LocalizationIdentity}), we
see
\[
\psi_{j}\cdot\iota f=\psi_{j}\cdot\sum_{\left(i,j_{0}\right)\in I_{0}\times J_{0}}\psi_{j_{0}}\varphi_{i}f=\sum_{j_{0}\in J_{0}\cap j^{\ast}}\left[\psi_{j}\cdot\sum_{i\in I_{0}}\psi_{j_{0}}\varphi_{i}f\right]=\sum_{j_{0}\in J_{0}\cap j^{\ast}}\psi_{j}f_{j_{0}},
\]
where we used the notation $f_{j_{0}}$ from Step $6$, equation~($\blacklozenge$).
We recall that $f_{j_{0}}$ satisfies $f_{j_{0}}=f_{k,j_{0}}$ for
every $k\in K$ with $j_{0}\in J^{\left(k\right)}$.

Now, we use the quasi-triangle inequality for $L^{p_{2}}$ and the
uniform bound $\left|j^{\ast}\right|\leq N_{\CalP}$ to obtain constants
$C_{7}=C_{7}\left(N_{\CalP},p_{2}\right)>0$ and $C_{8}=C_{8}\left(C_{\CalP,\Psi,1},C_{\CalP,\Psi,p_{2}},\CalP,\dimension,p_{2}\right)>0$
with
\begin{equation}
\begin{split}\left\Vert \Fourier^{-1}\left(\psi_{j}\cdot\iota f\right)\right\Vert _{L^{p_{2}}} & \leq C_{7}\cdot\sum_{j_{0}\in J_{0}\cap j^{\ast}}\left\Vert \Fourier^{-1}\left(\psi_{j}f_{j_{0}}\right)\right\Vert _{L^{p_{2}}}\\
 & \overset{\left(\ast\right)}{\leq}C_{7}C_{8}\cdot\sum_{j_{0}\in J_{0}\cap j^{\ast}}\left\Vert \Fourier^{-1}f_{j_{0}}\right\Vert _{L^{p_{2}}}\\
 & =C_{7}C_{8}\cdot\varrho_{j}^{\ast},
\end{split}
\label{eq:NormEstimatedByVarrhoJ}
\end{equation}
where $\varrho^{\ast}=\left(\varrho_{j}^{\ast}\right)_{j\in J}$ is
the ``clustered version'' of the sequence $\varrho$, given by $\varrho_{j}^{\ast}=\sum_{\ell\in j^{\ast}}\varrho_{\ell}$,
where the cluster $j^{\ast}$ is taken with respect to $\CalP$. To
justify the step marked with $\left(\ast\right)$, we distinguish
two cases:

\begin{casenv}
\item $p_{2}\in\left[1,\infty\right]$. In this case, we can simply use
Young's inequality and the uniform bound $\left\Vert \Fourier^{-1}\psi_{j}\right\Vert _{L^{1}}\leq C_{\CalP,\Psi,1}$
to conclude as desired that
\begin{align*}
\left\Vert \Fourier^{-1}\left(\psi_{j}f_{j_{0}}\right)\right\Vert _{L^{p_{2}}} & =\left\Vert \Fourier^{-1}\psi_{j}\ast\Fourier^{-1}f_{j_{0}}\right\Vert _{L^{p_{2}}}\\
 & \leq\left\Vert \Fourier^{-1}\psi_{j}\right\Vert _{L^{1}}\cdot\left\Vert \Fourier^{-1}f_{j_{0}}\right\Vert _{L^{p_{2}}}\\
 & \leq C_{\CalP,\Psi,1}\cdot\left\Vert \Fourier^{-1}f_{j_{0}}\right\Vert _{L^{p_{2}}}\:.
\end{align*}
\item $p_{2}\in\left(0,1\right)$. Here, by definition of an $L^{p_{2}}$-BAPU,
$\left\Vert \Fourier^{-1}\psi_{j}\right\Vert _{L^{p_{2}}}\leq C_{\CalP,\Psi,p_{2}}\cdot\left|\det S_{j}\right|^{1-p_{2}^{-1}}$
for all $j\in J$, so that Corollary~\ref{cor:QuasiBanachConvolutionSemiStructured}
implies because of $\supp f_{j_{0}}\subset\overline{P_{j_{0}}}$ (note
that we multiply with $\psi_{j_{0}}$ in the definition of $f_{j_{0}}=f_{k,j_{0}}$)
that
\begin{align*}
\left\Vert \Fourier^{-1}\left(\psi_{j}f_{j_{0}}\right)\right\Vert _{L^{p_{2}}} & \leq C\left(\CalP,\dimension,p_{2}\right)\cdot\left|\det S_{j}\right|^{p_{2}^{-1}-1}\cdot\left\Vert \Fourier^{-1}\psi_{j}\right\Vert _{L^{p_{2}}}\left\Vert \Fourier^{-1}f_{j_{0}}\right\Vert _{L^{p_{2}}}\\
 & \leq C\left(\CalP,\dimension,p_{2}\right)C_{\CalP,\Psi,p_{2}}\cdot\left\Vert \Fourier^{-1}f_{j_{0}}\right\Vert _{L^{p_{2}}}
\end{align*}
for all $j_{0}\in J_{0}\cap j^{\ast}$.
\end{casenv}
Now, we multiply equation~(\ref{eq:VarrhoJEstimate}) by $w_{k,j}$
and recall the definition of the sequence $\left(d_{k,j,i}\right)$
from equation~(\ref{eq:DSequenceDefinition}), to derive
\begin{equation}
w_{k,j}\cdot\varrho_{j}\leq C_{5}C_{6}\cdot\left|\det S_{j}\right|^{q_{k}^{-1}-p_{2}^{-1}}\cdot w_{k,j}\cdot\left\Vert \left(\theta_{k,j,i}\right)_{i\in I^{\left(k\right)}}\right\Vert _{\ell^{\LowerExpo{q_{k}}}}=C_{5}C_{6}\cdot\left\Vert \left(d_{k,j,i}\right)_{i\in I^{\left(k\right)}}\right\Vert _{\ell^{\LowerExpo{q_{k}}}}.\label{eq:VarrhoWithWeight}
\end{equation}

Finally, we use the $\CalP$-regularity of $Z$, the embeddings $\eta_{1},\eta_{2}$
and the solidity of $Z$ and $X$ to conclude
\begin{align*}
\left\Vert \iota f\right\Vert _{\FourierDecompSp{\CalP}{p_{2}}Z} & =\left\Vert \left(\left\Vert \Fourier^{-1}\left(\psi_{j}\cdot\iota f\right)\right\Vert _{L^{p_{2}}}\right)_{j\in J}\right\Vert _{Z}\\
\left({\scriptstyle \text{eq. }\eqref{eq:NormEstimatedByVarrhoJ}}\right) & \leq C_{7}C_{8}\cdot\left\Vert \left(\varrho_{j}^{\ast}\right)_{j\in J}\right\Vert _{Z}\\
\left({\scriptstyle Z\text{ is }\CalP\text{-regular}}\right) & \leq C_{7}C_{8}\cdot\vertiii{\Gamma_{\CalP}}_{Z\to Z}\cdot\left\Vert \varrho\right\Vert _{Z}\\
\left({\scriptstyle \vertiii{\eta_{1}}<\infty,\,\varrho_{j}=0\,\forall j\in J\setminus J_{0}}\right) & \leq C_{7}C_{8}\cdot\vertiii{\Gamma_{\CalP}}_{Z\to Z}\,\vertiii{\eta_{1}}\cdot\left\Vert \left(\left\Vert \left(\varrho_{j}\right)_{j\in J^{\left(k\right)}}\right\Vert _{\ell_{w}^{\UpperExpo{q_{k}}}}\right)_{k\in K}\right\Vert _{X}\\
 & =C_{7}C_{8}\cdot\vertiii{\Gamma_{\CalP}}_{Z\to Z}\,\vertiii{\eta_{1}}\cdot\left\Vert \left(\left\Vert \left(w_{k,j}\cdot\varrho_{j}\right)_{j\in J^{\left(k\right)}}\right\Vert _{\ell^{\UpperExpo{q_{k}}}}\right)_{k\in K}\right\Vert _{X}\\
\left({\scriptstyle \text{eq. }\eqref{eq:VarrhoWithWeight}}\right) & \leq C_{5}C_{6}C_{7}C_{8}\cdot\vertiii{\Gamma_{\CalP}}_{Z\to Z}\,\vertiii{\eta_{1}}\cdot\left\Vert \!\left(\left\Vert \left(\left\Vert \left(d_{k,j,i}\right)_{\!i\in I^{\left(k\right)}}\right\Vert _{\ell^{\LowerExpo{q_{k}}}}\!\!\right)_{\!\!j\in J^{\left(k\right)}}\right\Vert _{\ell^{\UpperExpo{q_{k}}}}\right)_{\!\!k\in K}\right\Vert _{X}\\
\left({\scriptstyle \text{eq. }\eqref{eq:DoubleNestedEstimateByQPieces}}\right) & \leq C_{1}C_{4}C_{5}C_{6}C_{7}C_{8}\cdot\vertiii{\Gamma_{\CalP}}_{Z\to Z}\,\vertiii{\eta_{1}}\cdot\left\Vert \left(\left\Vert \left(c_{i}\right)_{i\in I^{\left(k\right)}}\right\Vert _{\ell_{v}^{\LowerExpo{q_{k}}}}\right)_{k\in K}\right\Vert _{X}\\
\left({\scriptstyle \eta_{2}\text{ bounded}}\right) & \leq C_{1}C_{4}C_{5}C_{6}C_{7}C_{8}\cdot\vertiii{\Gamma_{\CalP}}_{Z\to Z}\,\vertiii{\eta_{1}}\cdot\vertiii{\eta_{2}}\cdot\left\Vert \left(c_{i}\right)_{i\in I}\right\Vert _{Y}\\
 & =C_{1}C_{4}C_{5}C_{6}C_{7}C_{8}\cdot\vertiii{\Gamma_{\CalP}}_{Z\to Z}\,\vertiii{\eta_{1}}\cdot\vertiii{\eta_{2}}\cdot\left\Vert f\right\Vert _{\FourierDecompSp{\CalQ}{p_{1}}Y}<\infty.
\end{align*}
Finally, Corollary~\ref{cor:LpBAPUsAreAlsoLqBAPUsForLargerq} yields
$C_{\CalP,\Psi,1}=C_{\CalP,\Psi,\infty}\lesssim_{\,\CalP,q^{\left(0\right)},\dimension}\:C_{\CalP,\Psi,q^{\left(0\right)}}$
and $C_{\CalP,\Psi,p_{2}}\lesssim_{\,\CalP,q^{\left(0\right)},\dimension}\:C_{\CalP,\Psi,q^{\left(0\right)}}$,
since $q^{\left(0\right)}\leq q_{k}\leq p_{2}\leq\infty$ for arbitrary
$k\in K$, so that the constant $C$ can be chosen as stated.
\end{proof}
Of course, verifying the assumptions of Theorem~\ref{thm:NoSubordinatenessWithConsiderationOfOverlaps}
is not easy in general. Thus, we formulate two important special cases
as corollaries. Our first result handles embeddings of a decomposition
space w.r.t. a ``fine'' covering into one w.r.t.\@ a ``coarse''
covering. This result is similar to \cite[Theorem 5.1.8]{VoigtlaenderPhDThesis}
from my PhD thesis, but slightly more general.
\begin{cor}
\label{cor:EmbeddingFineIntoCoarse}Let $\CalQ=\left(Q_{i}\right)_{i\in I}=\left(T_{i}Q_{i}'+b_{i}\right)_{i\in I}$
and $\CalP=\left(P_{j}\right)_{j\in J}=\left(S_{j}P_{j}'+c_{j}\right)_{j\in J}$
be semi-structured admissible coverings of the open sets $\CalO\subset\R^{\dimension}$
and $\CalO'\subset\R^{\dimension}$, respectively. Let $Y\subset\Compl^{I}$
and $Z\subset\Compl^{J}$ be two solid sequence spaces which are $\CalQ$-regular
and $\CalP$-regular, respectively.

Let $p_{1},p_{2}\in\left(0,\infty\right]$ with $p_{1}\leq p_{2}$
and assume that $\CalQ$ admits an $L^{p_{1}}$-BAPU $\Phi=\left(\varphi_{i}\right)_{i\in I}$
and that $\CalP$ admits an $L^{p_{2}}$-BAPU $\Psi=\left(\psi_{j}\right)_{j\in J}$.

Let $I_{0}\subset I$ be arbitrary with $\CalQ_{i}\subset\CalO'$
for all $i\in I_{0}$. In case of $p_{2}<1$, assume additionally
that the ``restricted'' family $\CalQ_{I_{0}}=\left(Q_{i}\right)_{i\in I_{0}}$
is almost subordinate to $\CalP$.

Finally, assume that the embedding\footnote{Recall the definition $Y|_{I_{0}}=\left\{ c=\left(c_{i}\right)_{i\in I_{0}}\in\Compl^{I_{0}}\with\tilde{c}\in Y\right\} $,
where $\tilde{c}=\left(c_{i}\right)_{i\in I}$, with $c_{i}=0$ for
$i\in I\setminus I_{0}$, see equation~(\ref{eq:DefinitionRestrictedSequenceSpace}).}
\[
\eta:Y|_{I_{0}}\hookrightarrow Z\left(\left[\,\vphantom{\sum}\smash{\ell_{\left|\det T_{i}\right|^{p_{1}^{-1}-p_{2}^{-1}}}^{\LowerExpo{p_{2}}}}\!\left(I_{0}\cap I_{j}\right)\,\right]_{j\in J}\right)
\]
is well-defined and bounded, with $I_{j}:=\left\{ i\in I\with Q_{i}\cap P_{j}\neq\emptyset\right\} $
for $j\in J$.

Then, the map
\[
\iota:\FourierDecompSp{\CalQ}{p_{1}}Y\to\FourierDecompSp{\CalP}{p_{2}}Z,f\mapsto\sum_{i\in I_{0}}\varphi_{i}f
\]
is well-defined and bounded with absolute convergence of the defining
series and with $\vertiii{\iota}\leq C\cdot\vertiii{\eta}$ for some
constant
\[
C=C\left(\dimension,p_{1},p_{2},k\left(\smash{\CalQ_{I_{0}}},\CalP\right),\CalQ,\CalP,C_{\CalQ,\Phi,p_{1}},C_{\CalP,\Psi,p_{2}},\vertiii{\Gamma_{\CalP}}_{Z\to Z}\right),
\]
where the dependence on $k\left(\smash{\CalQ_{I_{0}}},\CalP\right)$
can be dropped for $p_{2}\in\left[1,\infty\right]$.

Finally, if $I_{0}=I$, then the distribution $\iota f\in\DistributionSpace{\CalO'}$
is an extension of $f$ to $\TestFunctionSpace{\CalO'}\supset\TestFunctionSpace{\CalO}$,
for arbitrary $f\in\FourierDecompSp{\CalQ}{p_{1}}Y\subset\DistributionSpace{\CalO}$.
In particular, $\iota$ is injective if $I_{0}=I$.
\end{cor}

\begin{proof}
We want to apply Theorem~\ref{thm:NoSubordinatenessWithConsiderationOfOverlaps}
with $I_{0}$ as in the statement of the present corollary and the
following additional choices:
\[
K:=J,\qquad J^{\left(k\right)}:=\left\{ k\right\} ,\qquad I^{\left(k\right)}:=I_{k}\cap I_{0}=\left\{ i\in I_{0}\with Q_{i}\cap P_{k}\neq\emptyset\right\} \qquad\text{ and }\qquad q_{k}:=p_{2}
\]
for all $k\in K=J$. With these choices, it is clear that condition~(\ref{eq:SpecialIndexSetConditionPSelected})
is satisfied. Furthermore, $J_{00}=\bigcup_{k\in K}J^{\left(k\right)}=\bigcup_{k\in J}\left\{ k\right\} =J$,
so that we can choose $J_{0}:=J\subset J_{00}$. Finally, $q^{\left(0\right)}=\inf_{k\in K}q_{k}=p_{2}$.

Next, we estimate the weight $v$ from equation~(\ref{eq:WeightForIDefinition}).
To this end, let $k\in K$ and $i\in I^{\left(k\right)}$ be arbitrary.
For $p_{2}=q_{k}\geq1$, we simply have
\[
v_{k,i}=\left|\det T_{i}\right|^{p_{1}^{-1}-q_{k}^{-1}}=\left|\det T_{i}\right|^{p_{1}^{-1}-p_{2}^{-1}}.
\]
In case of $p_{2}=q_{k}<1$, our prerequisites include the assumption
that $\CalQ_{I_{0}}$ is almost subordinate to $\CalP$. Hence, set
$N:=k\left(\smash{\CalQ_{I_{0}}},\CalP\right)$. Now, since $J^{\left(k\right)}=\left\{ k\right\} $,
we have
\[
\sup_{j\in J^{\left(k\right)}}\lambda\left(\,\overline{P_{j}}-\overline{Q_{i}}\,\right)=\lambda\left(\,\overline{P_{k}}-\overline{Q_{i}}\,\right).
\]
Furthermore, for every $i\in I^{\left(k\right)}=I_{0}\cap I_{k}$,
there is some $j_{i}\in J$ with
\[
\emptyset\neq Q_{i}\cap P_{k}\subset P_{j_{i}}^{N\ast}\cap P_{k}
\]
and thus $k\in j_{i}^{\left(N+1\right)\ast}$, i.e.\@ $j_{i}\in k^{\left(N+1\right)\ast}$,
so that we get $Q_{i}\subset P_{j_{i}}^{N\ast}\subset P_{k}^{\left(2N+1\right)\ast}$.
Hence, Corollary~\ref{cor:SemiStructuredDifferenceSetsMeasureEstimate}
yields a constant $C_{1}=C_{1}\left(\CalP,N,\dimension\right)\geq1$
satisfying
\[
\smash{\sup_{j\in J^{\left(k\right)}}}\lambda\left(\,\overline{P_{j}}-\overline{Q_{i}}\,\right)=\lambda\left(\,\overline{P_{k}}-\overline{Q_{i}}\,\right)\leq\lambda\left(\,\overline{P_{k}^{\left(2N+1\right)\ast}}-\overline{P_{k}^{\left(2N+1\right)\ast}}\,\right)\leq C_{1}\cdot\left|\det S_{k}\right|.
\]
All in all, for $C_{2}:=C_{1}^{\frac{1}{p_{2}}-1}\geq1$ (in case
of $p_{2}<1$), we have shown
\[
v_{k,i}\leq\begin{cases}
C_{2}\cdot\left|\det T_{i}\right|^{p_{1}^{-1}-p_{2}^{-1}}\cdot\left|\det S_{k}\right|^{p_{2}^{-1}-1}, & \text{if }p_{2}<1,\\
\vphantom{\rule{0.1cm}{0.55cm}}\left|\det T_{i}\right|^{p_{1}^{-1}-p_{2}^{-1}}, & \text{if }p_{2}\geq1.
\end{cases}
\]

Now, our final choice for the application of Theorem~\ref{thm:NoSubordinatenessWithConsiderationOfOverlaps}
is to set $X:=Z_{1/w}\subset\Compl^{J}=\Compl^{K}$, where the weight
$w$ is chosen as in equation~(\ref{eq:WeightForJDefinition}), i.e.
\[
w_{j}:=w_{j,j}=w_{k,j}:=\begin{cases}
\left|\det S_{j}\right|^{p_{2}^{-1}-1}, & \text{if }p_{2}<1,\\
\vphantom{\rule{0.1cm}{0.55cm}}\left|\det S_{j}\right|^{p_{2}^{-1}-q_{k}^{-1}}=\left|\det S_{j}\right|^{p_{2}^{-1}-p_{2}^{-1}}=1, & \text{if }p_{2}\geq1
\end{cases}
\]
for $k\in K=J$ and $j\in J^{\left(k\right)}=\left\{ k\right\} $,
i.e. $j=k$. Using this choice and recalling $J^{\left(k\right)}=\left\{ k\right\} $
for $k\in K=J$, we see
\[
X\!\Bigl(\left[\vphantom{\sum}\smash{\ell_{w}^{\UpperExpo{q_{k}}}}\!\left(\smash{J^{\left(k\right)}}\right)\right]_{k\in K}\Bigr)=X_{w}=\left(Z_{1/w}\right)_{w}=Z\hookrightarrow Z,
\]
so that the map $\eta_{1}$ from condition~(\ref{eq:AssumedDiscreteEmbedding1})
of Theorem~\ref{thm:NoSubordinatenessWithConsiderationOfOverlaps}
is well-defined and bounded with $\vertiii{\eta_{1}}=1$.

To prove the boundedness of $\eta_{2}$ from condition~(\ref{eq:AssumedDiscreteEmbedding2})
of Theorem~\ref{thm:NoSubordinatenessWithConsiderationOfOverlaps},
first note that every $i\in I_{0}$ satisfies $\emptyset\neq Q_{i}\subset\CalO'$
and hence $Q_{i}\cap P_{j_{i}}\neq\emptyset$ for some $j_{i}\in J$.
In particular, there is some $k\in J=K$ satisfying $Q_{i}\cap P_{k}\neq\emptyset$,
i.e. $i\in I^{\left(k\right)}$, so that we have $L:=\bigcup_{k\in K}I^{\left(k\right)}=I_{0}$.
Furthermore, observe
\[
X\Bigl(\bigl[\ell_{v}^{\LowerExpo{q_{k}}}\!\left(\smash{I^{\left(k\right)}}\right)\bigr]\Bigr)_{k\in K}=Z\Bigl(\Bigl[\ell_{\left(v_{k,i}/w_{k}\right)_{i}}^{\LowerExpo{p_{2}}}\left(\smash{I^{\left(k\right)}}\right)\Bigr]_{k\in K}\Bigr)\text{ with identical (quasi)-norms},
\]
where
\begin{align}
\frac{v_{k,i}}{w_{k}} & \leq\begin{cases}
C_{2}\cdot\left|\det T_{i}\right|^{p_{1}^{-1}-p_{2}^{-1}}\cdot\left|\det S_{k}\right|^{p_{2}^{-1}-1}\cdot\left|\det S_{k}\right|^{1-p_{2}^{-1}}, & \text{if }p_{2}<1,\\
\vphantom{\rule{0.1cm}{0.55cm}}\left|\det T_{i}\right|^{p_{1}^{-1}-p_{2}^{-1}}, & \text{if }p_{2}\geq1
\end{cases}\nonumber \\
 & \leq C_{3}\cdot\left|\det T_{i}\right|^{p_{1}^{-1}-p_{2}^{-1}},\label{eq:FineInCoarseEmbeddingWeightQuotientEstimate}
\end{align}
with $C_{3}:=C_{2}$ if $p_{2}<1$ and $C_{3}:=1$ if $p_{2}\geq1$.

Thus, our assumption (boundedness of $\eta$) shows that we have
\begin{eqnarray*}
Y|_{I_{0}} & \Xhookrightarrow{\eta} & Z\Bigl(\Bigl[\ell_{\left|\det T_{i}\right|^{p_{1}^{-1}-p_{2}^{-1}}}^{\LowerExpo{p_{2}}}\left(I_{0}\cap I_{j}\right)\Bigr]_{j\in J}\Bigr)\\
 & \Xhookrightarrow{\text{eq. }\eqref{eq:FineInCoarseEmbeddingWeightQuotientEstimate}\text{ and }I^{\left(k\right)}=I_{0}\cap I_{k}} & Z\Bigl(\Bigl[\ell_{\left(v_{k,i}/w_{k}\right)_{i}}^{\LowerExpo{p_{2}}}\left(\smash{I^{\left(k\right)}}\right)\Bigr]_{k\in K}\Bigr)=X\Bigl(\bigl[\ell_{v}^{\LowerExpo{q_{k}}}\!\left(\smash{I^{\left(k\right)}}\right)\bigr]_{k\in K}\Bigr),
\end{eqnarray*}
which easily shows that the embedding $\eta_{2}$ from condition~(\ref{eq:AssumedDiscreteEmbedding2})
is well-defined and bounded with $\vertiii{\eta_{2}}\leq C_{3}\cdot\vertiii{\eta}$.
Hence, an application of Theorem~\ref{thm:NoSubordinatenessWithConsiderationOfOverlaps}
(together with the first part of the ensuing remark) implies that
$\iota$ is bounded with 
\[
\vertiii{\iota}\leq C_{4}\cdot\vertiii{\eta_{1}}\cdot\vertiii{\eta_{2}}\leq C_{3}C_{4}\cdot\vertiii{\eta}
\]
for some constant $C_{4}=C_{4}\left(\dimension,p_{1},p_{2},\CalQ,\CalP,C_{\CalQ,\Phi,p_{1}},C_{\CalP,\Psi,p_{2}},\vertiii{\Gamma_{\CalP}}_{Z\to Z}\right)$.
This completes the proof of the first part of the corollary.

Finally, in case of $I_{0}=I$, note that we have $Q_{i}\subset\CalO'$
for all $i\in I_{0}=I$ and hence $\CalO=\bigcup_{i\in I}Q_{i}\subset\CalO'$,
which implies $\TestFunctionSpace{\CalO}\subset\TestFunctionSpace{\CalO'}$.
Now, since $\left(\varphi_{i}\right)_{i\in I}$ is a (locally finite)
partition of unity on $\CalO$ (see Lemma~\ref{lem:PartitionCoveringNecessary}),
we have for every $g\in\TestFunctionSpace{\CalO}\subset\TestFunctionSpace{\CalO'}$
that $g=\sum_{i\in I}\varphi_{i}\,g$, where only finitely many terms
do not vanish identically. Hence,
\[
\left\langle \iota f,\,g\right\rangle _{\CalD'}=\sum_{i\in I}\left\langle f,\,\varphi_{i}\,g\right\rangle _{\CalD'}=\left\langle f,\,\smash{\sum_{i\in I}\,}\vphantom{\sum}\varphi_{i}\,g\right\rangle _{\CalD'}=\left\langle f,g\right\rangle _{\CalD'},
\]
so that $\iota f$ indeed extends $f$ if $I_{0}=I$.
\end{proof}
If one wants to apply the preceding corollary, one faces two major
challenges: First—at least for $p_{2}<1$—one has to verify that $\CalQ$
(or $\CalQ_{I_{0}}$) is almost subordinate to $\CalP$. Then, one
has to verify the boundedness of the embedding $\eta$; in fact, it
might also be desired to obtain a bound for $\vertiii{\eta}$, since
this yields a bound for $\vertiii{\iota}$. Luckily, using the results
from Section~\ref{sec:NestedSequenceSpaces}, one can greatly simplify
this second part, at least if $Y$ and $Z$ are both weighted $\ell^{q}$
spaces and if $\CalQ_{I_{0}}$ is almost subordinate to $\CalP$:
\begin{cor}
\label{cor:EmbeddingFineInCoarseSimplification}Let $\CalQ=\left(Q_{i}\right)_{i\in I}$
and $\CalP=\left(P_{j}\right)_{j\in J}$ be two admissible coverings
of the open sets $\emptyset\neq\CalO,\CalO'\subset\R^{\dimension}$,
respectively. Furthermore, let $w=\left(w_{i}\right)_{i\in I}$ and
$v=\left(v_{j}\right)_{j\in J}$ be $\CalQ$- and $\CalP$-moderate,
respectively. Finally, let $q_{1},q_{2},r\in\left(0,\infty\right]$
and let $u=\left(u_{i}\right)_{i\in I}$ be a further weight on $I$
(which is not necessarily $\CalQ$-moderate).

If $I_{0}\subset I$ is chosen such that $\CalQ_{I_{0}}=\left(Q_{i}\right)_{i\in I_{0}}$
is almost subordinate to $\CalP$ and if $J_{0}\subset J$ satisfies\footnote{The typical choice will simply be $J_{0}=J$. In this case, the conclusion
$I_{0}\subset\bigcup_{j\in J_{0}}I_{j}$ is always satisfied; indeed,
for $i\in I_{0}$, we have $\emptyset\neq Q_{i}\subset\CalO'$, since
$\CalQ_{I_{0}}$ is almost subordinate to $\CalP$. Since $\CalP$
covers $\CalO'$, this yields $Q_{i}\cap P_{j}\neq\emptyset$ for
some $j\in J$ and hence $i\in I_{j}$.} $I_{0}\subset\bigcup_{j\in J_{0}}I_{j}$, then we have
\begin{equation}
\vertiii{\eta}\asymp\left\Vert \left(v_{j}\cdot\left\Vert \left(u_{i}/w_{i}\right)_{i\in I_{0}\cap I_{j}}\right\Vert _{\ell^{r\cdot\left(q_{1}/r\right)'}}\right)_{j\in J_{0}}\right\Vert _{\ell^{q_{2}\cdot\left(q_{1}/q_{2}\right)'}}\label{eq:FineInCoarseSimplification}
\end{equation}
for
\[
\eta:\ell_{w}^{q_{1}}\left(I_{0}\right)\hookrightarrow\ell_{v}^{q_{2}}\left(\left[\ell_{u}^{r}\left(I_{0}\cap I_{j}\right)\right]_{j\in J_{0}}\right).
\]

Precisely, this means that $\eta$ is well-defined and bounded if
and only if the right-hand side of equation~(\ref{eq:FineInCoarseSimplification})
is finite. Further, there is a constant $C\geq1$ depending only on
$q_{1},q_{2},r,N_{\CalP},k\left(\smash{\CalQ_{I_{0}}},\CalP\right),C_{v,\CalP}$
which satisfies
\[
C^{-1}\cdot M\leq\vertiii{\eta}\leq C\cdot M\quad\text{ for }\quad M:=\left\Vert \left(v_{j}\cdot\left\Vert \left(u_{i}/w_{i}\right)_{i\in I_{0}\cap I_{j}}\right\Vert _{\ell^{r\cdot\left(q_{1}/r\right)'}}\right)_{j\in J_{0}}\right\Vert _{\ell^{q_{2}\cdot\left(q_{1}/q_{2}\right)'}}.
\]

\medskip{}

Finally, if

\begin{itemize}
\item $\CalQ=\left(T_{i}Q_{i}'+b_{i}\right)_{i\in I}$ and $\CalP=\left(S_{j}P_{j}'+c_{j}\right)_{j\in J}$
are \emph{tight} semi-structured coverings,
\item $\CalQ_{I_{0}}$ is relatively $\CalP$-moderate,
\item there is some $s\in\N_{0}$ and some $C_{0}>0$ such that
\[
\lambda\left(P_{j}\right)\leq C_{0}\cdot\vphantom{\bigcup_{i\in I_{0}\cap I_{j}}}\lambda\left(\,\vphantom{\bigcup}\smash{\bigcup_{i\in I_{0}\cap I_{j}}}Q_{i}^{s\ast}\,\right)\qquad\forall\,j\in J_{0}\text{ with }1\leq\left|I_{0}\cap I_{j}\right|<\infty,\,
\]
\item $u|_{I_{0}}$ and $w|_{I_{0}}$ are relatively $\CalP$-moderate and
\item for each $j\in J^{\left(0\right)}:=\left\{ j\in J_{0}\with I_{0}\cap I_{j}\neq\emptyset\right\} $,
some $i_{j}\in I_{0}\cap I_{j}$ is selected,
\end{itemize}
then
\[
\left\Vert \left(v_{j}\!\cdot\!\left\Vert \left(u_{i}/w_{i}\right)_{i\in I_{0}\cap I_{j}}\right\Vert _{\ell^{r\cdot\left(q_{1}/r\right)'}}\right)_{\!\!j\in J_{0}}\right\Vert _{\ell^{q_{2}\cdot\left(q_{1}/q_{2}\right)'}}\asymp\left\Vert \left(\,\frac{v_{j}\cdot u_{i_{j}}}{w_{i_{j}}}\!\cdot\!\smash{\left[\frac{\left|\det S_{j}\right|}{\left|\det T_{i_{j}}\right|}\right]^{\left(\frac{1}{r}-\frac{1}{q_{1}}\right)_{+}}}\vphantom{\left[\frac{\left|\det S_{j}\right|}{\left|\det T_{i_{j}}\right|}\right]}\,\right)_{\!\!j\in J^{\left(0\right)}}\right\Vert _{\ell^{q_{2}\cdot\left(q_{1}/q_{2}\right)'}}\;,
\]
where the implied constant only depends on 
\[
\dimension,s,r,q_{1},\CalQ,\CalP,\varepsilon_{\CalQ},\varepsilon_{\CalP},C_{0},k\left(\smash{\CalQ_{I_{0}}},\CalP\right),C_{u|_{I_{0}},\CalQ_{I_{0}},\CalP},C_{w|_{I_{0}},\CalQ_{I_{0}},\CalP},C_{{\rm mod}}\left(\smash{\CalQ_{I_{0}}},\CalP\right).\qedhere
\]
\end{cor}

\begin{rem}
\label{rem:SufficientFineInCoarseSimplification}In particular, if
we have $Y=\ell_{w}^{q_{1}}\left(I\right)$ and $Z=\ell_{v}^{q_{2}}\left(J\right)$
in Corollary~\ref{cor:EmbeddingFineIntoCoarse}, then the embedding
$\eta$ from that corollary satisfies
\[
\vertiii{\eta}\asymp\left\Vert \left(v_{j}\cdot\left\Vert \left(\left|\det T_{i}\right|^{p_{1}^{-1}-p_{2}^{-1}}/w_{i}\right)_{i\in I_{0}\cap I_{j}}\right\Vert _{\ell^{\LowerExpo{p_{2}}\cdot\left(q_{1}/\LowerExpo{p_{2}}\right)'}}\right)_{j\in J}\right\Vert _{\ell^{q_{2}\cdot\left(q_{1}/q_{2}\right)'}}\;,
\]
where the implied constant only depends on $q_{1},q_{2},p_{2},N_{\CalP},k\left(\smash{\CalQ_{I_{0}}},\CalP\right),C_{v,\CalP}$.

Furthermore, if the additional assumptions from the second part of
the corollary are satisfied, we have
\begin{align*}
\vertiii{\eta} & \asymp\left\Vert \left(\frac{v_{j}}{w_{i_{j}}}\cdot\left|\det T_{i_{j}}\right|^{\frac{1}{p_{1}}-\frac{1}{p_{2}}}\cdot\left[\left|\det S_{j}\right|/\left|\det T_{i_{j}}\right|\right]^{\left(\smash{\frac{1}{\LowerExpo{p_{2}}}}-\frac{1}{q_{1}}\right)_{+}}\right)_{\!\!j\in J^{\left(0\right)}}\right\Vert _{\ell^{q_{2}\cdot\left(q_{1}/q_{2}\right)'}}\\
 & =\left\Vert \left(\frac{v_{j}}{w_{i_{j}}}\cdot\left|\det T_{i_{j}}\right|^{\frac{1}{p_{1}}-\left(\smash{\frac{1}{\LowerExpo{p_{2}}}}-\frac{1}{q_{1}}\right)_{+}-\frac{1}{p_{2}}}\cdot\left|\det S_{j}\right|^{\left(\smash{\frac{1}{\LowerExpo{p_{2}}}}-\frac{1}{q_{1}}\right)_{+}}\right)_{\!\!j\in J^{\left(0\right)}}\right\Vert _{\ell^{q_{2}\cdot\left(q_{1}/q_{2}\right)'}}
\end{align*}
where the implied constant only depends on 
\[
\dimension,s,q_{1},q_{2},p_{1},p_{2},\CalQ,\CalP,\varepsilon_{\CalQ},\varepsilon_{\CalP},C_{0},k\left(\CalQ_{I_{0}},\CalP\right),C_{w|_{I_{0}},\CalQ_{I_{0}},\CalP},C_{v,\CalP},C_{{\rm mod}}\left(\CalQ_{I_{0}},\CalP\right).\qedhere
\]
\end{rem}

\begin{proof}
We begin with the first part of the corollary. For the proof, we want
to use the estimate provided by parts~(\ref{enu:AlmostDisjointifiedPureIntoNested})
and (\ref{enu:AlmostDisjointifiedLebesgueSimplified}) of Corollary~\ref{cor:AlmostDisjointifiedEmbeddingLebesgue}
with $K=J_{0}$, $I=I_{0}$, $X=\ell_{v}^{q_{2}}\left(K\right)=\ell_{v}^{q_{2}}\left(J_{0}\right)$
and 
\[
I^{\left(k,\natural\right)}:=I^{\left(k\right)}:=I_{0}\cap I_{k}=\left\{ i\in I_{0}\with Q_{i}\cap P_{k}\neq\emptyset\right\} 
\]
for $k\in K=J_{0}$. For this, we first have to verify the prerequisites
of Corollary~\ref{cor:AlmostDisjointifiedEmbeddingLebesgue}, which
are just the prerequisites of Lemma~\ref{lem:NestedEmbeddingReductionToDisjointSets}.

To this end, first note that we have $I_{0}=\bigcup_{k\in K}I^{\left(k\right)}$.
Indeed, ``$\supset$'' is trivial; and ``$\subset$'' follows
from our assumption $I_{0}\subset\bigcup_{j\in J_{0}}I_{j}$.

Now, define a relation $\sim$ on $J_{0}=K$ by $j\sim\ell\::\Longleftrightarrow\:I^{\left(j\right)}\cap I^{\left(\ell\right)}\neq\emptyset\Longleftrightarrow I_{0}\cap I_{j}\cap I_{\ell}\neq\emptyset$.
For brevity, let us write $k:=k\left(\CalQ_{I_{0}},\CalP\right)\in\N_{0}$,
which is well-defined, since $\CalQ_{I_{0}}$ is almost subordinate
to $\CalP$. Lemma~\ref{lem:SubordinatenessEnablesDisjointization}
(with $n=0$, with $\CalQ_{I_{0}}$ instead of $\CalQ$, and with
$\CalP$ instead of $\CalR$) implies\footnote{Strictly speaking, we need to extend the relation $\sim$ to all of
$J$ to apply Lemma~\ref{lem:SubordinatenessEnablesDisjointization}.
But this can be done simply by defining $I^{\left(k\right)}:=\emptyset$
for $k\in J\setminus J_{0}$.} that the classes $\left[j\right]$ of the given relation $\sim$
satisfy
\[
\left[j\right]\subset j^{\left(4k+5\right)\ast}\qquad\text{ as well as }\qquad\left|\left[j\right]\right|\leq N_{\CalP}^{4k+5}=:N\qquad\forall\,j\in J_{0}.
\]
In particular, condition~(\ref{eq:ReductionToDisjointSetsCardinalityAssumption})
of Lemma~\ref{lem:NestedEmbeddingReductionToDisjointSets} is satisfied.

Furthermore, with the $4k+5$-fold clustering map $\Theta_{4k+5}$
as in Lemma~\ref{lem:HigherOrderClusteringMap} (for $\CalP$ instead
of $\CalQ$), the generalized clustering map $\Theta$ from Lemma~\ref{lem:NestedEmbeddingReductionToDisjointSets}
satisfies 
\[
\left|\left(\Theta x\right)_{j}\right|\leq\sum_{\ell\in\left[j\right]}\left|x_{\ell}\right|\leq\sum_{\ell\in j^{\left(4k+5\right)\ast}}\Indicator_{J_{0}}\left(\ell\right)\cdot\left|x_{\ell}\right|=\left(\Theta_{4k+5}\left[\Indicator_{J_{0}}\cdot\left|x\right|\right]\right)_{j}
\]
for all $j\in J_{0}$ and $x=\left(x_{j}\right)_{j\in J_{0}}\in X$
(which we extended by $0$ to all of $J$ to obtain a sequence $\left(x_{j}\right)_{j\in J}$).
Hence,
\begin{align*}
\left\Vert \Theta x\right\Vert _{X}=\left\Vert \Theta x\right\Vert _{\ell_{v}^{q_{2}}\left(J_{0}\right)} & \leq\left\Vert \Theta_{4k+5}\left[\Indicator_{J_{0}}\cdot\left|x\right|\right]\right\Vert _{\ell_{v}^{q_{2}}\left(J_{0}\right)}\\
 & \leq\left\Vert \Theta_{4k+5}\left[\Indicator_{J_{0}}\cdot\left|x\right|\right]\right\Vert _{\ell_{v}^{q_{2}}\left(J\right)}\\
 & \leq\vertiii{\Theta_{4k+5}}_{\ell_{v}^{q_{2}}\left(J\right)\to\ell_{v}^{q_{2}}\left(J\right)}\cdot\left\Vert \Indicator_{J_{0}}\cdot\left|x\right|\right\Vert _{\ell_{v}^{q_{2}}\left(J\right)}\\
 & =\vertiii{\Theta_{4k+5}}_{\ell_{v}^{q_{2}}\left(J\right)\to\ell_{v}^{q_{2}}\left(J\right)}\cdot\left\Vert x\right\Vert _{\ell_{v}^{q_{2}}\left(J_{0}\right)}\\
\left({\scriptstyle \text{Lemmas }\ref{lem:HigherOrderClusteringMap}\text{ and }\ref{lem:ModeratelyWeightedSpacesAreRegular}}\right) & \leq\left(C_{v,\CalP}\cdot N_{\CalP}^{\max\left\{ 1,q_{2}^{-1}\right\} }\right)^{4k+5}\cdot\left\Vert x\right\Vert _{X}\:.
\end{align*}

For the final prerequisite of Lemma~\ref{lem:NestedEmbeddingReductionToDisjointSets}
(regarding the ``inner'' weight $u=\left(u_{j,i}\right)$), note
that we have in our case $u_{j,i}=u_{i}$ for all $j\in J_{0}$ and
$i\in I^{\left(j\right)}=I_{0}\cap I_{j}$, so that $u_{j,i}\leq u_{\ell,i}$
holds for all $j,\ell\in J_{0}$ and $i\in I^{\left(j\right)}\cap I^{\left(\ell\right)}$;
that is, we can choose $C_{u}=1$.

Now, we can finally apply Corollary~\ref{cor:AlmostDisjointifiedEmbeddingLebesgue}
(with slightly permuted roles of the weights $u,v,w$ and of the exponents
of the weighted $\ell^{q}$ spaces) to get
\[
\vertiii{\eta}\asymp\left\Vert \left(v_{j}\cdot\left\Vert \left(u_{i}/w_{i}\right)_{i\in I_{0}\cap I_{j}}\right\Vert _{\ell^{r\cdot\left(q_{1}/r\right)'}}\right)_{j\in J_{0}}\right\Vert _{\ell^{q_{2}\cdot\left(q_{1}/q_{2}\right)'}}\quad,
\]
where the implied constant only depends on $N_{\CalP},k,q_{1},q_{2},r,C_{v,\CalP}$,
since the triangle constant $C_{X}$ of $X=\ell_{v}^{q_{2}}\left(J_{0}\right)$
only depends on $q_{2}$.

\medskip{}

Now, we consider the second part of the corollary. Note that the present
assumptions include those of Lemma~\ref{lem:IntersectionCountForModerateCoverings}
(since $\CalP$-moderateness of $\CalQ_{I_{0}}$ implies $\CalP_{J_{0}}$-moderateness
of $\CalQ_{I_{0}}$; precisely, we have $C_{{\rm mod}}\left(\CalQ_{I_{0}},\CalP_{J_{0}}\right)\leq C_{{\rm mod}}\left(\CalQ_{I_{0}},\CalP\right)$).
Thus, that lemma yields
\[
\left|I_{0}\cap I_{j}\right|\asymp\left|\det S_{j}\right|\big/\left|\det T_{i_{j}}\right|
\]
for each $j\in J_{0}$ with $I_{0}\cap I_{j}\neq\emptyset$, i.e.\@
for each $j\in J^{\left(0\right)}$. Here, the implied constant only
depends on those quantities which are mentioned in the second part
of the present corollary.

But for $j\in J^{\left(0\right)}$ and $i\in I_{0}\cap I_{j}$, we
also have
\[
\left(C_{w|_{I_{0}},\CalQ_{I_{0}},\CalP}\cdot C_{u|_{I_{0}},\CalQ_{I_{0}},\CalP}\right)^{-1}\cdot\frac{u_{i_{j}}}{w_{i_{j}}}\leq\frac{u_{i}}{w_{i}}\leq C_{w|_{I_{0}},\CalQ_{I_{0}},\CalP}\cdot C_{u|_{I_{0}},\CalQ_{I_{0}},\CalP}\cdot\frac{u_{i_{j}}}{w_{i_{j}}}.
\]
All in all, this easily implies
\begin{align*}
\left\Vert \left(u_{i}/w_{i}\right)_{i\in I_{0}\cap I_{j}}\right\Vert _{\ell^{r\cdot\left(q_{1}/r\right)'}} & \asymp\frac{u_{i_{j}}}{w_{i_{j}}}\cdot\left|I_{0}\cap I_{j}\right|^{\left[r\cdot\left(q_{1}/r\right)'\right]^{-1}}\\
 & \asymp\frac{u_{i_{j}}}{w_{i_{j}}}\cdot\left[\left|\det S_{j}\right|\big/\left|\det T_{i_{j}}\right|\right]^{\left[r\cdot\left(q_{1}/r\right)'\right]^{-1}}\\
 & =\frac{u_{i_{j}}}{w_{i_{j}}}\cdot\left[\left|\det S_{j}\right|\big/\left|\det T_{i_{j}}\right|\right]^{\left(\frac{1}{r}-\frac{1}{q_{1}}\right)_{+}}\;,
\end{align*}
with implied constants as stated in the corollary. Here, the last
step used $\frac{1}{a\cdot\left(b/a\right)'}=\left(\frac{1}{a}-\frac{1}{b}\right)_{+}$,
see equation~(\ref{eq:InverseOfSpecialExponent}).

Since the left-hand side of the desired estimate is unaffected if
all terms $j\in J_{0}\setminus J^{\left(0\right)}$ are discarded,
this easily yields the claim.
\end{proof}
In Corollary~\ref{cor:EmbeddingFineIntoCoarse}, we assumed $\CalQ$
to be almost subordinate to $\CalP$. In our next result (which is
similar to \cite[Theorem 5.1.6]{VoigtlaenderPhDThesis} from my PhD
thesis), we make the ``reverse'' assumption. To simplify the proof,
we state the following technical lemma in advance.
\begin{lem}
\label{lem:EmbeddingInMixedSpaceWithSmallInnerSets}Let $\CalQ=\left(Q_{i}\right)_{i\in I}$
be an admissible covering of a set $\CalO$ and let $X,Y\subset\Compl^{I}$
be $\CalQ$-regular sequence spaces. Fix $n\in\N_{0}$ and assume
that for each $i\in I$, some $q_{i}\in\left(0,\infty\right]$ and
some subset $I_{i}\subset i^{n\ast}$ is given. Finally, assume that
$q:=\inf_{i\in I}q_{i}>0$ and let $v=\left(v_{i}\right)_{i\in I}$
be an arbitrary (positive) weight.

If $\iota:Y\hookrightarrow X_{v}$ is bounded, then so is $\theta:Y\hookrightarrow X\left(\left[\ell_{v}^{q_{i}}\left(I_{i}\right)\right]_{i\in I}\right)$,
with
\[
\vertiii{\theta}\leq N_{\CalQ}^{n/q}\cdot\vertiii{\Gamma_{\CalQ}}_{X\to X}^{n}\cdot\vertiii{\iota}\:.\qedhere
\]
\end{lem}

\begin{proof}
Let $x=\left(x_{i}\right)_{i\in I}\in Y$ be arbitrary. Let $i\in I$.
Using the uniform bound $\left|I_{i}\right|\leq\left|i^{n\ast}\right|\leq N_{\CalQ}^{n}$
(see Lemma~\ref{lem:SemiStructuredClusterInvariant}) and the norm-decreasing
embedding $\ell^{q}\hookrightarrow\ell^{q_{i}}$, we get (for $q<\infty$)
that
\begin{alignat*}{2}
\left\Vert \left(x_{\ell}\right)_{\ell\in I_{i}}\right\Vert _{\ell_{v}^{q_{i}}}\leq\left\Vert \left(x_{\ell}\right)_{\ell\in I_{i}}\right\Vert _{\ell_{v}^{q}} & =\left[\,\sum_{\ell\in I_{i}}\left(v_{\ell}\left|x_{\ell}\right|\right)^{q}\right]^{1/q} &  & \leq\left|I_{i}\right|^{1/q}\cdot\max_{\ell\in I_{i}}v_{\ell}\left|x_{\ell}\right|\\
 & \leq N_{\CalQ}^{n/q}\cdot\sum_{\ell\in i^{n\ast}}\left|v_{\ell}x_{\ell}\right| &  & =N_{\CalQ}^{n/q}\cdot\left(\Theta_{n}\left|v\cdot x\right|\right)_{i}\:,
\end{alignat*}
where $\Theta_{n}$ denotes the $n$-fold clustering map from Lemma~\ref{lem:HigherOrderClusteringMap}
and where $\left|v\cdot x\right|_{\ell}=v_{\ell}\cdot\left|x_{\ell}\right|$.
It is not hard to see that this estimate remains valid for $q=\infty$.

By solidity of $X$, we conclude
\begin{alignat*}{2}
\left\Vert x\right\Vert _{X\left(\left[\ell_{v}^{q_{i}}\left(I_{i}\right)\right]_{i\in I}\right)} & \leq N_{\CalQ}^{n/q}\cdot\left\Vert \Theta_{n}\left|v\cdot x\right|\right\Vert _{X} &  & \leq N_{\CalQ}^{n/q}\vertiii{\Theta_{n}}_{X\to X}\cdot\left\Vert \left|v\cdot x\right|\right\Vert _{X}\\
 & =N_{\CalQ}^{n/q}\vertiii{\Theta_{n}}_{X\to X}\cdot\left\Vert x\right\Vert _{X_{v}} &  & \leq N_{\CalQ}^{n/q}\vertiii{\Theta_{n}}_{X\to X}\vertiii{\iota}\cdot\left\Vert x\right\Vert _{Y}<\infty.
\end{alignat*}
Recalling the bound $\vertiii{\Theta_{n}}_{X\to X}\leq\vertiii{\Gamma_{\CalQ}}_{X\to X}^{n}$
from Lemma~\ref{lem:HigherOrderClusteringMap}, we get the desired
bound.
\end{proof}
\begin{cor}
\label{cor:EmbeddingCoarseIntoFine}Let $\CalQ=\left(Q_{i}\right)_{i\in I}=\left(T_{i}Q_{i}'+b_{i}\right)_{i\in I}$
and $\CalP=\left(P_{j}\right)_{j\in J}=\left(S_{j}P_{j}'+c_{j}\right)_{j\in J}$
be semi-structured admissible coverings of the open sets $\CalO\subset\R^{\dimension}$
and $\CalO'\subset\R^{\dimension}$, respectively. Let $Y\subset\Compl^{I}$
and $Z\subset\Compl^{J}$ be two solid sequence spaces which are $\CalQ$-regular
and $\CalP$-regular, respectively.

Let $p_{1},p_{2}\in\left(0,\infty\right]$ with $p_{1}\leq p_{2}$
and assume that $\CalQ$ and $\CalP$ admit $L^{p_{1}}$-BAPUs $\Phi=\left(\varphi_{i}\right)_{i\in I}$
and $\Psi=\left(\psi_{j}\right)_{j\in J}$, respectively.

Finally, let $J_{0}\subset J$ and assume that the ``restricted''
family $\CalP_{J_{0}}:=\left(P_{j}\right)_{j\in J_{0}}$ is almost
subordinate to $\CalQ$ and that
\[
\eta:Y\Bigl(\bigl[\ell_{u}^{\UpperExpo{p_{1}}}\!\left(J_{0}\cap J_{i}\right)\bigr]_{i\in I}\Bigr)\hookrightarrow Z|_{J_{0}}
\]
is well-defined and bounded, where $J_{i}:=\left\{ j\in J\with P_{j}\cap Q_{i}\neq\emptyset\right\} $
and
\begin{equation}
u_{i,j}:=\begin{cases}
\left|\det S_{j}\right|^{p_{2}^{-1}-1}\cdot\left|\det T_{i}\right|^{1-p_{1}^{-1}}, & \text{if }p_{1}<1,\\
\vphantom{\rule{0.1cm}{0.55cm}}\left|\det S_{j}\right|^{p_{2}^{-1}-p_{1}^{-1}}, & \text{if }p_{1}\geq1.
\end{cases}\label{eq:CoarseIntoFineInnerWeight}
\end{equation}
Then, the map
\[
\iota:\FourierDecompSp{\CalQ}{p_{1}}Y\to\FourierDecompSp{\CalP}{p_{2}}Z,f\mapsto\sum_{j\in J_{0}}\psi_{j}f
\]
is well-defined and bounded with absolute convergence of the defining
series and with $\vertiii{\iota}\leq C\cdot\vertiii{\eta}$ for some
constant 
\[
C=C\left(\dimension,p_{1},p_{2},k\left(\smash{\CalP_{J_{0}}},\CalQ\right),\CalQ,\CalP,C_{\CalQ,\Phi,p_{1}},C_{\CalP,\Psi,p_{1}},\vertiii{\Gamma_{\CalP}}_{Z\to Z},\vertiii{\Gamma_{\CalQ}}_{Y\to Y}\right).
\]

Finally, in case of $J_{0}=J$, we have $\CalO'\subset\CalO$ and
$\iota f=f|_{\TestFunctionSpace{\CalO'}}$ for all $f\in\FourierDecompSp{\CalQ}{p_{1}}Y$.
\end{cor}

\begin{proof}
By assumption, there is $N:=k\left(\smash{\CalP_{J_{0}}},\CalQ\right)\in\N_{0}$
such that for every $j\in J_{0}$, we have $P_{j}\subset Q_{i_{j}}^{N\ast}$
for some $i_{j}\in I$. Now, we want to apply Theorem~\ref{thm:NoSubordinatenessWithConsiderationOfOverlaps}
with the following choices:
\[
K:=I_{0}:=I,\qquad I^{\left(k\right)}:=k^{\left(2N+2\right)\ast},\qquad J^{\left(k\right)}:=J_{0}\cap J_{k}\qquad\text{ and }\qquad q_{k}:=p_{1}
\]
for all $k\in K=I$. After verifying the assumptions of Theorem~\ref{thm:NoSubordinatenessWithConsiderationOfOverlaps}
with these choices, we will show that the map $\iota$ from that theorem
coincides with the map $\iota$ as defined in the present corollary
(with the special case for $J_{0}=J$).

We first verify condition~(\ref{eq:SpecialIndexSetConditionPSelected}).
Thus, assume to the contrary that there is some $j\in J^{\left(k\right)}=J_{0}\cap J_{k}$
and some $i\in I\setminus I^{\left(k\right)}$ with $Q_{i}\cap P_{j}\neq\emptyset$.
Because of $j\in J_{0}$, we have $\emptyset\neq Q_{i}\cap P_{j}\subset Q_{i}\cap Q_{i_{j}}^{N\ast}$
and thus $i\in i_{j}^{\left(N+1\right)\ast}$. Furthermore, $j\in J_{0}\cap J_{k}$
implies $\emptyset\neq P_{j}\cap Q_{k}\subset Q_{i_{j}}^{N\ast}\cap Q_{k}$
and hence $k\in i_{j}^{\left(N+1\right)\ast}$, i.e.\@ $i_{j}\in k^{\left(N+1\right)\ast}$,
which finally yields $i\in i_{j}^{\left(N+1\right)\ast}\subset k^{\left(2N+2\right)\ast}=I^{\left(k\right)}$,
in contradiction to $i\in I\setminus I^{\left(k\right)}$. We have
thus verified condition~(\ref{eq:SpecialIndexSetConditionPSelected}).

Next, we want to verify conditions~(\ref{eq:AssumedDiscreteEmbedding1})
and (\ref{eq:AssumedDiscreteEmbedding2}) with $X:=Y_{1/\theta}\subset\Compl^{I}$,
where $\theta_{i}:=\left|\det T_{i}\right|^{\left(p_{1}^{-1}-1\right)_{+}}$.
Note that equation~(\ref{eq:DeterminantIsModerate}) shows that $\left(\left|\det T_{i}\right|\right)_{i\in I}$
is $\CalQ$-moderate with $C_{\left|\det T_{i}\right|,\CalQ}\leq C_{\CalQ}^{\dimension}$,
so that also $1/\theta$ is $\CalQ$-moderate with $C_{\theta,\CalQ}\leq C_{\CalQ}^{\dimension\left(p_{1}^{-1}-1\right)_{+}}$.
By Lemma~\ref{lem:ModeratelyWeightedSpacesAreRegular}, this implies
that $X$ is $\CalQ$-regular, with 
\[
\vertiii{\Gamma_{\CalQ}}_{X\to X}\leq C_{\CalQ}^{\dimension\left(p_{1}^{-1}-1\right)_{+}}\cdot\vertiii{\Gamma_{\CalQ}}_{Y\to Y}=:C_{1}.
\]

Note that we have $I^{\left(k\right)}\supset\left\{ k\right\} $ and
hence $L:=\bigcup_{k\in K}I^{\left(k\right)}=I$. Thus, condition~(\ref{eq:AssumedDiscreteEmbedding2})
requires boundedness of the embedding
\[
\eta_{2}:Y\overset{!}{\hookrightarrow}X\left(\left[\,\vphantom{\ell_{v}^{M}}\smash{\ell_{v}^{\LowerExpo{q_{k}}}}\!\left(\smash{I^{\left(k\right)}}\right)\,\right]_{k\in K}\right)=X\Bigl(\bigl[\,\ell_{v}^{\LowerExpo{p_{1}}}\!\left(\smash{k^{\left(2N+2\right)\ast}}\right)\,\bigr]_{k\in K}\Bigr),
\]
with
\begin{align*}
v_{k,i} & =\begin{cases}
\left|\det T_{i}\right|^{p_{1}^{-1}-q_{k}^{-1}}\cdot\left[\sup_{j\in J^{\left(k\right)}}\lambda\left(\,\overline{P_{j}}-\overline{Q_{i}}\,\right)\right]^{q_{k}^{-1}-1}, & \text{if }q_{k}<1,\\
\vphantom{\rule{0.1cm}{0.55cm}}\left|\det T_{i}\right|^{p_{1}^{-1}-q_{k}^{-1}}, & \text{if }q_{k}\geq1
\end{cases}\\
 & =\begin{cases}
\left[\sup_{j\in J^{\left(k\right)}}\lambda\left(\,\overline{P_{j}}-\overline{Q_{i}}\,\right)\right]^{p_{1}^{-1}-1}, & \text{if }p_{1}<1,\\
\vphantom{\rule{0.1cm}{0.55cm}}1, & \text{if }p_{1}\geq1
\end{cases}
\end{align*}
for $k\in K=I$ and $i\in I^{\left(k\right)}=k^{\left(2N+2\right)\ast}$.

We now estimate this weight further in case of $p_{1}<1$. For $j\in J^{\left(k\right)}=J_{0}\cap J_{k}$,
Lemma~\ref{lem:SubordinatenessImpliesWeakSubordination} yields $P_{j}\subset Q_{k}^{\left(2N+2\right)\ast}\subset Q_{i}^{\left(4N+4\right)\ast}$,
so that we get
\[
\lambda\left(\,\overline{P_{j}}-\overline{Q_{i}}\,\right)\leq\lambda\left(\,\overline{Q_{i}^{\left(4N+4\right)\ast}}-\overline{Q_{i}^{\left(4N+4\right)\ast}}\,\right)\leq C_{2}\cdot\left|\det T_{i}\right|
\]
for some constant $C_{2}=C_{2}\left(\CalQ,N,\dimension\right)\geq1$,
see Corollary~\ref{cor:SemiStructuredDifferenceSetsMeasureEstimate}.
All in all, we have shown that there is a constant $C_{3}=C_{3}\left(\CalQ,N,d,p_{1}\right)\geq1$
with
\[
v_{k,i}\leq\begin{cases}
C_{3}\cdot\left|\det T_{i}\right|^{p_{1}^{-1}-1}=C_{3}\cdot\theta_{i}, & \text{if }p_{1}<1,\\
\vphantom{\rule{0.1cm}{0.55cm}}1\leq C_{3}\cdot\theta_{i}, & \text{if }p_{1}\geq1.
\end{cases}
\]
Thus, Lemma~\ref{lem:EmbeddingInMixedSpaceWithSmallInnerSets} yields
(because of $Y=X_{\theta}$ and $K=I$) that
\begin{eqnarray*}
Y & \Xhookrightarrow{\vertiii{\cdot}\leq C_{4}=C_{4}\left(N,p_{1},N_{\CalQ},\vertiii{\Gamma_{\CalQ}}_{X\to X}\right)} & X\Bigl(\Bigl[\ell_{\theta}^{\LowerExpo{p_{1}}}\!\bigl(k^{\left(2N+2\right)\ast}\bigr)\Bigr]_{k\in K}\Bigr)\\
 & \Xhookrightarrow{\vertiii{\cdot}\leq C_{3}\text{ since }v_{k,i}\leq C_{3}\cdot\theta_{i}} & X\Bigl(\Bigl[\ell_{v}^{\LowerExpo{p_{1}}}\!\bigl(k^{\left(2N+2\right)\ast}\bigr)\Bigr]_{k\in K}\Bigr)=X\Bigl(\Bigl[\ell_{v}^{\LowerExpo{q_{k}}}\!\bigl(\smash{I^{\left(k\right)}}\bigr)\Bigr]_{k\in K}\Bigr).
\end{eqnarray*}
All in all, this shows that $\eta_{2}$ is bounded with $\vertiii{\eta_{2}}\leq C_{3}C_{4}$,
so that condition~(\ref{eq:AssumedDiscreteEmbedding2}) is satisfied.

Verification of condition~(\ref{eq:AssumedDiscreteEmbedding1}) is
easier: First note that each $j\in J_{0}$ satisfies $\emptyset\neq P_{j}\subset Q_{i_{j}}^{N\ast}$
for some $i_{j}\in I$. In particular, $P_{j}\cap Q_{k}\neq\emptyset$
for some $k\in K=I$, so that we get $j\in J_{0}\cap J_{k}=J^{\left(k\right)}$.
Hence—in the notation of Theorem~\ref{thm:NoSubordinatenessWithConsiderationOfOverlaps}—we
have $J_{00}=\bigcup_{k\in K}J^{\left(k\right)}=J_{0}$. Thus, (because
of $X=Y_{1/\theta}$), condition~(\ref{eq:AssumedDiscreteEmbedding1})
precisely requires boundedness of
\[
\eta_{1}:Y\!\Bigl(\Bigl[\ell_{\left(w_{i,j}/\theta_{i}\right)}^{\UpperExpo{p_{1}}}\bigl(J^{\left(i\right)}\bigr)\Bigr]_{i\in I}\Bigr)\!=\!X\!\Bigl(\Bigl[\ell_{w}^{\UpperExpo{q_{k}}}\!\bigl(J^{\left(k\right)}\bigr)\Bigr]_{k\in K}\Bigr)\overset{!}{\hookrightarrow}Z,\left(x_{j}\right)_{j\in J_{0}}\mapsto\left(x_{j}\right)_{j\in J}\text{ with }x_{j}=0\text{ for }j\in J\setminus J_{0},
\]
where (see equation~(\ref{eq:WeightForJDefinition}))
\[
\frac{w_{i,j}}{\theta_{i}}=\begin{cases}
\left|\det S_{j}\right|^{p_{2}^{-1}-1}\cdot\left|\det T_{i}\right|^{1-p_{1}^{-1}}=u_{i,j}\:, & \text{if }p_{1}<1,\\
\vphantom{\rule{0.1cm}{0.55cm}}\left|\det S_{j}\right|^{p_{2}^{-1}-p_{1}^{-1}}=u_{i,j}\:, & \text{if }p_{1}\geq1.
\end{cases}
\]
Thus, a moment's thought shows that we have $\vertiii{\eta_{1}}=\vertiii{\eta}$,
with $\eta$ as in the statement of the present corollary.

Now, Theorem~\ref{thm:NoSubordinatenessWithConsiderationOfOverlaps}
shows that
\[
\iota^{\left(0\right)}:\FourierDecompSp{\CalQ}{p_{1}}Y\to\FourierDecompSp{\CalP}{p_{2}}Z,f\mapsto\sum_{\left(i,j\right)\in I\times J_{0}}\varphi_{i}\,\psi_{j}\,f
\]
is well-defined and bounded with absolute convergence of the defining
series and with
\[
\vertiii{\smash{\iota^{\left(0\right)}}}\leq C_{5}\cdot\vertiii{\eta_{1}}\cdot\vertiii{\eta_{2}}\leq C_{3}C_{4}C_{5}\cdot\vertiii{\eta}
\]
for some constant $C_{5}=C_{5}\left(\dimension,p_{1},p_{2},\CalQ,\CalP,C_{\CalQ,\Phi,p_{1}},C_{\CalP,\Psi,p_{1}},\vertiii{\Gamma_{\CalP}}_{Z\to Z}\right)$.

\smallskip{}

But for $j\in J_{0}$, we have $P_{j}\subset Q_{i_{j}}^{N\ast}$ for
some $i_{j}\in I$. Because of $\overline{Q_{i_{j}}^{N\ast}}\subset\CalO$
(see Lemma~\ref{lem:PartitionCoveringNecessary}), this yields $\psi_{j}\in\TestFunctionSpace{\CalO}$,
so that $\psi_{j}\cdot f\in\Schwartz'\left(\R^{\dimension}\right)$
is well-defined for all $f\in\FourierDecompSp{\CalQ}{p_{1}}Y\subset\DistributionSpace{\CalO}$.
Finally, since $\left(\varphi_{i}\right)_{i\in I}$ is a (locally
finite) partition of unity on $\CalO$, we get $\psi_{j}=\sum_{i\in I}\varphi_{i}\psi_{j}$
and thus (using the absolute convergence of the series)
\[
\left\langle \,\smash{\iota^{\left(0\right)}}f,\,g\right\rangle _{\CalD'}=\sum_{j\in J_{0}}\;\sum_{i\in I}\left\langle f,\,\varphi_{i}\,\psi_{j}\,g\right\rangle _{\CalD'}=\sum_{j\in J_{0}}\left\langle f,\,\psi_{j}\,g\right\rangle _{\CalD'}=\left\langle \smash{\sum_{j\in J_{0}}}\vphantom{\sum}\psi_{j}\,f,\:g\right\rangle _{\CalD'}
\]
for all $g\in\TestFunctionSpace{\CalO'}$. Hence, $\iota^{\left(0\right)}f=\iota f$
with $\iota$ as in the statement of the present corollary. In particular,
this implies absolute convergence of the series defining $\left\langle \iota f,\,g\right\rangle _{\CalD'}$
for all $g\in\TestFunctionSpace{\CalO}$.

Finally, assume $J=J_{0}$. We just saw $P_{j}\subset\CalO$ for all
$j\in J_{0}$, so that $\CalO'=\bigcup_{j\in J}P_{j}\subset\CalO$.
Now, for $g\in\TestFunctionSpace{\CalO'}\subset\TestFunctionSpace{\CalO}$,
we have $g=\sum_{j\in J}\psi_{j}\,g$ with only finitely many terms
not vanishing, since $\left(\psi_{j}\right)_{j\in J}$ is a locally
finite partition of unity on $\CalO'$; see Lemma~\ref{lem:PartitionCoveringNecessary}.
Hence, we get as desired that
\[
\left\langle \iota f,\,g\right\rangle _{\CalD'}=\sum_{j\in J}\left\langle f,\:\psi_{j}\,g\right\rangle _{\CalD'}=\left\langle f,\,g\right\rangle _{\CalD'}=\left\langle f|_{\TestFunctionSpace{\CalO'}},\:g\right\rangle _{\CalD'}\qquad\forall\,g\in\TestFunctionSpace{\CalO'}\,.\qedhere
\]
\end{proof}
As for Corollary~\ref{cor:EmbeddingFineInCoarseSimplification} above,
one can greatly simplify the process of verifying boundedness of the
embedding $\eta$ from above, at least if $Y,Z$ are both weighted
$\ell^{q}$ spaces:
\begin{cor}
\label{cor:EmbeddingCoarseIntoFineSimplification}Let $\CalQ=\left(Q_{i}\right)_{i\in I}$
and $\CalP=\left(P_{j}\right)_{j\in J}$ be two admissible coverings
of the open sets $\emptyset\neq\CalO,\CalO'\subset\R^{\dimension}$,
respectively. Furthermore, let $w=\left(w_{i}\right)_{i\in I}$ and
$u^{\left(1\right)}=\left(\smash{u_{i}^{\left(1\right)}}\right)_{i\in I}$
be $\CalQ$-moderate and let $v=\left(v_{j}\right)_{j\in J}$ and
$u^{\left(2\right)}=\left(\smash{u_{j}^{\left(2\right)}}\right)_{j\in J}$
be weights on the index set $J$. Finally, let $q_{1},q_{2},r\in\left(0,\infty\right]$
and define $u:=\vphantom{u_{i}^{\left(1\right)}u_{j}^{\left(2\right)}}\left(\smash{u_{i}^{\left(1\right)}u_{j}^{\left(2\right)}}\right)_{i\in I,j\in J}$.

If $J_{0}\subset J$ is chosen so that $\CalP_{J_{0}}:=\left(P_{j}\right)_{j\in J_{0}}$
is almost subordinate to $\CalQ$, then we have
\begin{equation}
\vertiii{\eta}\asymp\left\Vert \left(w_{i}^{-1}\cdot\left\Vert \left(v_{j}/u_{i,j}\right)_{j\in J_{0}\cap J_{i}}\right\Vert _{\ell^{q_{2}\cdot\left(r/q_{2}\right)'}}\right)_{i\in I}\right\Vert _{\ell^{q_{2}\cdot\left(q_{1}/q_{2}\right)'}}\label{eq:CoarseInFineSimplification}
\end{equation}
for
\[
\eta:\ell_{w}^{q_{1}}\left(\left[\ell_{u}^{r}\left(J_{0}\cap J_{i}\right)\right]_{i\in I}\right)\hookrightarrow\ell_{v}^{q_{2}}\left(J_{0}\right).
\]

Precisely, this means that $\eta$ is well-defined and bounded if
and only if the right-hand side of equation~(\ref{eq:CoarseInFineSimplification})
is finite. Furthermore, there is a constant $C\geq1$ depending only
on $r,q_{1},q_{2}$ and on $C_{w,\CalQ},N_{\CalQ},k\left(\smash{\CalP_{J_{0}}},\CalQ\right),C_{u^{\left(1\right)},\CalQ}$
which satisfies
\[
C^{-1}\cdot M\leq\vertiii{\eta}\leq C\cdot M\quad\text{ for }\quad M:=\left\Vert \left(w_{i}^{-1}\cdot\left\Vert \left(v_{j}/u_{i,j}\right)_{j\in J_{0}\cap J_{i}}\right\Vert _{\ell^{q_{2}\cdot\left(r/q_{2}\right)'}}\right)_{i\in I}\right\Vert _{\ell^{q_{2}\cdot\left(q_{1}/q_{2}\right)'}}.
\]

\medskip{}

\noindent Finally, if

\begin{itemize}[leftmargin=0.7cm]
\item $\CalQ=\left(T_{i}Q_{i}'+b_{i}\right)_{i\in I}$ and $\CalP=\left(S_{j}P_{j}'+c_{j}\right)_{j\in J}$
are \emph{tight} semi-structured coverings,
\item $\CalP_{J_{0}}$ is relatively $\CalQ$-moderate,
\item there is some $s\in\N_{0}$ and some $C_{0}>0$ such that
\[
\lambda\left(Q_{i}\right)\leq C_{0}\vphantom{\bigcup_{j\in J_{0}\cap J_{i}}}\cdot\lambda\left(\,\smash{\bigcup_{j\in J_{0}\cap J_{i}}}\vphantom{\bigcup}P_{j}^{s\ast}\right)\qquad\forall\,i\in I\text{ with }1\leq\left|J_{0}\cap J_{i}\right|<\infty\,,
\]
\item $u^{\left(2\right)}|_{J_{0}}$ and $v|_{J_{0}}$ are relatively $\CalQ$-moderate,
and\vspace{0.1cm}
\item for each $i\in I^{\left(0\right)}:=\left\{ i\in I\with J_{0}\cap J_{i}\neq\emptyset\right\} $,
some $j_{i}\in J_{0}\cap J_{i}$ is selected,
\end{itemize}
then
\[
\left\Vert \!\left(w_{i}^{-1}\!\cdot\left\Vert \left(v_{j}/u_{i,j}\right)_{j\in J_{0}\cap J_{i}}\right\Vert _{\ell^{q_{2}\cdot\left(r/q_{2}\right)'}}\right)_{\!i\in I}\right\Vert _{\ell^{q_{2}\cdot\left(q_{1}/q_{2}\right)'}}\asymp\left\Vert \!\left(\frac{v_{j_{i}}}{w_{i}\cdot u_{i,j_{i}}}\!\cdot\!\smash{\left[\frac{\left|\det T_{i}\right|}{\left|\det S_{j_{i}}\right|}\right]^{\left(\frac{1}{q_{2}}-\frac{1}{r}\right)_{+}}}\vphantom{\left[\frac{\left|\det T_{i}\right|}{\left|\det S_{j_{i}}\right|}\right]}\right)_{\!\!i\in I^{\left(0\right)}}\right\Vert _{\ell^{q_{2}\cdot\left(q_{1}/q_{2}\right)'}}\;,
\]
where the implied constant only depends on 
\[
\dimension,s,r,q_{2},C_{0},k\left(\smash{\CalP_{J_{0}}},\CalQ\right),C_{{\rm mod}}\left(\smash{\CalP_{J_{0}}},\CalQ\right),\CalQ,\CalP,\varepsilon_{\CalQ},\varepsilon_{\CalP},C_{u^{\left(2\right)}|_{J_{0}},\CalP,\CalQ},C_{v|_{J_{0}},\CalP,\CalQ}\:.\qedhere
\]
\end{cor}

\begin{rem}
\label{rem:SufficientCoarseIntoFineSimplification}In particular,
if we have $Y=\ell_{w}^{q_{1}}\left(I\right)$ and $Z=\ell_{v}^{q_{2}}\left(J\right)$
in Corollary~\ref{cor:EmbeddingCoarseIntoFine}, then the embedding
$\eta$ from that corollary satisfies
\begin{align*}
\vertiii{\eta} & \asymp\left\Vert \left(w_{i}^{-1}\cdot\left\Vert \left(v_{j}/u_{i,j}\right)_{j\in J_{0}\cap J_{i}}\right\Vert _{\ell^{q_{2}\cdot\left(\UpperExpo{p_{1}}/q_{2}\right)'}}\right)_{i\in I}\right\Vert _{\ell^{q_{2}\cdot\left(q_{1}/q_{2}\right)'}}=:M\\
\text{with}\qquad u_{i,j}\; & =\begin{cases}
\left|\det S_{j}\right|^{p_{2}^{-1}-1}\cdot\left|\det T_{i}\right|^{1-p_{1}^{-1}}, & \text{if }p_{1}<1,\\
\vphantom{\rule{0.1cm}{0.55cm}}\left|\det S_{j}\right|^{p_{2}^{-1}-p_{1}^{-1}}, & \text{if }p_{1}\geq1,
\end{cases}
\end{align*}
where the implied constant only depends on $p_{1},q_{1},q_{2},C_{w,\CalQ},\CalQ,k\left(\smash{\CalP_{J_{0}}},\CalQ\right)$.
Here, we used the preceding corollary with
\[
u_{i}^{\left(1\right)}:=\begin{cases}
\left|\det T_{i}\right|^{1-p_{1}^{-1}}, & \text{if }p_{1}<1,\\
\vphantom{\rule{0.1cm}{0.55cm}}1, & \text{if }p_{1}\geq1
\end{cases}\qquad\text{ and }\qquad u_{j}^{\left(2\right)}:=\begin{cases}
\left|\det S_{j}\right|^{p_{2}^{-1}-1}, & \text{if }p_{1}<1,\\
\vphantom{\rule{0.1cm}{0.55cm}}\left|\det S_{j}\right|^{p_{2}^{-1}-p_{1}^{-1}}, & \text{if }p_{1}\geq1.
\end{cases}
\]

Finally, if the additional assumptions from the second part of the
corollary are satisfied (where the relative $\CalQ$-moderateness
of $u^{\left(2\right)}|_{J_{0}}$ is implied by the relative $\CalQ$-moderateness
of $\CalP_{J_{0}}$), we get
\begin{equation}
M\asymp\left\Vert \left(\frac{v_{j_{i}}}{w_{i}}\cdot\left|\det S_{j_{i}}\right|^{\frac{1}{p_{1}}-\left(\frac{1}{q_{2}}-\smash{\frac{1}{\SignedUpperExpo{p_{1}}}}\right)_{+}-\frac{1}{p_{2}}}\cdot\left|\det T_{i}\right|^{\left(\frac{1}{q_{2}}-\smash{\frac{1}{\SignedUpperExpo{p_{1}}}}\right)_{+}}\right)_{i\in I^{\left(0\right)}}\right\Vert _{\ell^{q_{2}\cdot\left(q_{1}/q_{2}\right)'}}\label{eq:SufficientCoarseIntoFineSimplificationRelativelyModerate}
\end{equation}
with $1/\SignedUpperExpo{p_{1}}=\min\left\{ p_{1}^{-1},1-p_{1}^{-1}\right\} $,
that is, with 
\[
\SignedUpperExpo{p_{1}}:=\begin{cases}
p_{1}, & \text{if }p_{1}\geq2,\\
\frac{p_{1}}{p_{1}-1}, & \text{if }0<p_{1}<2\text{ and }p_{1}\neq1,\\
\infty, & \text{if }p_{1}=1.
\end{cases}
\]

To see this, note that $\UpperExpo{p_{1}}=\SignedUpperExpo{p_{1}}$
for $p_{1}\in\left[1,\infty\right]$, while $\SignedUpperExpo{p_{1}}$
is negative for $0<p_{1}<1$. Furthermore, note that if we define
the \textbf{duality defect} $p_{d}$ of $p\in\left(0,\infty\right]$
by
\[
p_{d}:=\min\left\{ 0,1-p^{-1}\right\} =\begin{cases}
1-\frac{1}{p}, & \text{if }p<1,\\
0, & \text{if }p\geq1,
\end{cases}
\]
then $1/\SignedUpperExpo{p_{1}}=\left(p_{1}\right)_{d}+1/\UpperExpo{p_{1}}$,
as a simple case distinction shows. Furthermore, $u_{i}^{\left(1\right)}=\left|\det T_{i}\right|^{\left(p_{1}\right)_{d}}$
and $u_{j}^{\left(2\right)}=\left|\det S_{j}\right|^{\frac{1}{p_{2}}-\frac{1}{p_{1}}-\left(p_{1}\right)_{d}}$
for all $i\in I$ and $j\in J$, so that we get
\[
u_{i,j_{i}}^{-1}\cdot\left[\left|\det T_{i}\right|/\left|\det S_{j_{i}}\right|\right]^{\left(\frac{1}{q_{2}}-\smash{\frac{1}{\UpperExpo{p_{1}}}}\right)_{+}}=\left|\det T_{i}\right|^{\left(\frac{1}{q_{2}}-\smash{\frac{1}{\UpperExpo{p_{1}}}}\right)_{+}-\left(p_{1}\right)_{d}}\cdot\left|\det S_{j_{i}}\right|^{\frac{1}{p_{1}}-\frac{1}{p_{2}}-\left(\frac{1}{q_{2}}-\smash{\frac{1}{\UpperExpo{p_{1}}}}\right)_{+}+\left(p_{1}\right)_{d}}.
\]
Finally, note
\[
\left(\frac{1}{q_{2}}-\smash{\frac{1}{\UpperExpo{p_{1}}}}\right)_{+}-\left(p_{1}\right)_{d}=\left(\frac{1}{q_{2}}-\smash{\frac{1}{\SignedUpperExpo{p_{1}}}}\right)_{+}.
\]
In case of $p_{1}\in\left[1,\infty\right]$, this is clear. But for
$p_{1}\in\left(0,1\right)$, we have $\UpperExpo{p_{1}}=\infty$ and
hence
\[
\left(\frac{1}{q_{2}}-\smash{\frac{1}{\UpperExpo{p_{1}}}}\right)_{+}-\left(p_{1}\right)_{d}=\frac{1}{q_{2}}-\left(p_{1}\right)_{d}=\frac{1}{q_{2}}-1+\frac{1}{p_{1}}=\frac{1}{q_{2}}-\frac{1}{\SignedUpperExpo{p_{1}}}.
\]
Furthermore, since $\SignedUpperExpo{p_{1}}<0$ for $p_{1}\in\left(0,1\right)$,
we have $\left(\frac{1}{q_{2}}-\smash{\frac{1}{\SignedUpperExpo{p_{1}}}}\right)_{+}=\frac{1}{q_{2}}-\frac{1}{\SignedUpperExpo{p_{1}}}$.
All in all, these considerations establish equation~(\ref{eq:SufficientCoarseIntoFineSimplificationRelativelyModerate}),
where the implied constant only depends on 
\[
\dimension,s,p_{1},p_{2},q_{2},C_{0},k\left(\smash{\CalP_{J_{0}}},\CalQ\right),C_{{\rm mod}}\left(\smash{\CalP_{J_{0}}},\CalQ\right),\CalQ,\CalP,\varepsilon_{\CalQ},\varepsilon_{\CalP},C_{v|_{J_{0}},\CalP,\CalQ}\,.\qedhere
\]
\end{rem}

\begin{proof}[Proof of Corollary~\ref{cor:EmbeddingCoarseIntoFineSimplification}]
For the proof of the first part, we want to use the estimate provided
by Corollary~\ref{cor:AlmostDisjointifiedEmbeddingLebesgue}, with
$K=I$, $I=J_{0}$, $X=\ell_{w}^{q_{1}}\left(K\right)=\ell_{w}^{q_{1}}\left(I\right)$
and 
\[
I^{\left(k,\natural\right)}:=I^{\left(k\right)}:=J_{0}\cap J_{k}=\left\{ j\in J_{0}\with P_{j}\cap Q_{k}\neq\emptyset\right\} \qquad\text{ for }k\in K=I,
\]
to calculate the norm of this embedding. For this, we first have to
verify the assumptions of Corollary~\ref{cor:AlmostDisjointifiedEmbeddingLebesgue},
which are just the assumptions of Lemma~\ref{lem:NestedEmbeddingReductionToDisjointSets}.

To this end, first note that we have $J_{0}=\bigcup_{k\in K}I^{\left(k\right)}$.
Indeed, ``$\supset$'' is trivial; for the reverse inclusion, let
$j\in J_{0}$ be arbitrary. Since $\CalP_{J_{0}}$ is almost subordinate
to $\CalQ$, we have $P_{j}\subset\CalO$ and hence $P_{j}\cap Q_{i}\neq\emptyset$
for some $i\in I=K$. Hence, $j\in J_{0}\cap J_{i}=I^{\left(i\right)}$.

Now, define a relation $\sim$ on $K=I$ by $i\sim\ell\::\Longleftrightarrow\:I^{\left(i\right)}\cap I^{\left(\ell\right)}\neq\emptyset$.
For brevity, let us write $n:=k\left(\CalP_{J_{0}},\CalQ\right)\in\N_{0}$,
which is well-defined, since $\CalP_{J_{0}}$ is almost subordinate
to $\CalQ$. Lemma~\ref{lem:SubordinatenessEnablesDisjointization}
(with $n=0$ and $\CalP_{J_{0}},\CalQ$ for $\left(Q_{i}\right)_{i\in I},\CalR$)
implies that the the classes $\left[i\right]$ of this relation satisfy
\[
\left[i\right]\subset i^{\left(4n+5\right)\ast},\qquad\text{ as well as }\qquad\left|\left[i\right]\right|\leq N_{\CalQ}^{4n+5}=:N\qquad\forall\,i\in I.
\]
In particular, condition~(\ref{eq:ReductionToDisjointSetsCardinalityAssumption})
of Lemma~\ref{lem:NestedEmbeddingReductionToDisjointSets} is satisfied.
Furthermore, the generalized clustering map $\Theta$ from that lemma
satisfies 
\[
\left|\left(\Theta x\right)_{i}\right|\leq\sum_{\ell\in\left[i\right]}\left|x_{\ell}\right|\leq\sum_{\ell\in i^{\left(4n+5\right)\ast}}\left|x_{\ell}\right|=\left(\Theta_{4n+5}\left|x\right|\right)_{i}\qquad\forall\,x=\left(x_{i}\right)_{i\in I}\in X=\ell_{w}^{q_{1}}\left(I\right)
\]
for all $i\in I$, with $\Theta_{n}$ (for $n\in\N_{0}$) as in Lemma~\ref{lem:HigherOrderClusteringMap}.
Thus, by solidity of $X=\ell_{w}^{q_{1}}$, we see that $\Theta$
is well-defined and bounded with 
\[
\vertiii{\Theta}\leq\vertiii{\Theta_{4n+5}}\leq\vertiii{\Gamma_{\CalQ}}_{\ell_{w}^{q_{1}}\to\ell_{w}^{q_{1}}}^{4n+5}\leq\Bigl(C_{w,\CalQ}\cdot N_{\CalQ}^{\max\left\{ 1,q_{1}^{-1}\right\} }\Bigr)^{4n+5}\:,
\]
see Lemma~\ref{lem:ModeratelyWeightedSpacesAreRegular}.

For the final assumption of Lemma~\ref{lem:NestedEmbeddingReductionToDisjointSets}
(regarding the weight $u=\left(u_{i,j}\right)$) let $i,\ell\in K=I$
and $j\in I^{\left(i\right)}\cap I^{\left(\ell\right)}$. In particular,
this implies $i\in\left[\ell\right]\subset\ell^{\left(4n+5\right)\ast}$,
so that Lemma~\ref{lem:SemiStructuredClusterInvariant} yields
\[
u_{i,j}=u_{i}^{\left(1\right)}\cdot u_{j}^{\left(2\right)}\leq C_{u^{\left(1\right)},\CalQ}^{4n+5}\cdot u_{\ell}^{\left(1\right)}\cdot u_{j}^{\left(2\right)}=C_{u^{\left(1\right)},\CalQ}^{4n+5}\cdot u_{\ell,j}\:.
\]
Hence, we can choose (in the notation of Lemma~\ref{lem:NestedEmbeddingReductionToDisjointSets})
$C_{u}=C_{u^{\left(1\right)},\CalQ}^{4n+5}$.

Now, we can finally apply parts~(\ref{enu:AlmostDisjointifiedNestedIntoPure})
and (\ref{enu:AlmostDisjointifiedLebesgueSimplified}) of Corollary~\ref{cor:AlmostDisjointifiedEmbeddingLebesgue},
which yield
\[
\vertiii{\eta}\asymp\left\Vert \left(w_{i}^{-1}\cdot\left\Vert \left(v_{j}/u_{i,j}\right)_{j\in J_{0}\cap J_{i}}\right\Vert _{\ell^{q_{2}\cdot\left(r/q_{2}\right)'}}\right)_{i\in I}\right\Vert _{\ell^{q_{2}\cdot\left(q_{1}/q_{2}\right)'}}\quad,
\]
where the implied constant only depends on
\[
\vertiii{\Theta}\leq\Bigl(C_{w,\CalQ}\cdot N_{\CalQ}^{\max\left\{ 1,q_{1}^{-1}\right\} }\Bigr)^{4n+5},\;N=N_{\CalQ}^{4n+5},\;C_{u}=C_{u^{\left(1\right)},\CalQ}^{4n+5}\quad\text{and on}\quad r,q_{2}.
\]

\medskip{}

For the final part of the corollary, note that the present assumptions
include the prerequisites of Lemma~\ref{lem:IntersectionCountForModerateCoverings}
(with interchanged roles of $\CalQ,\CalP$, with $I_{0}=J_{0}$, and
with $J_{0}=I$), so that we get
\[
\left|J_{0}\cap J_{i}\right|\asymp\left|\det T_{i}\right|/\left|\det S_{j_{i}}\right|
\]
for all $i\in I$ with $J_{0}\cap J_{i}\neq\emptyset$, i.e.\@ for
all $i\in I^{\left(0\right)}$. Here, the implied constant only depends
on those quantities which are mentioned in the second part of the
present corollary.

Furthermore, relative $\CalQ$-moderateness of $u^{\left(2\right)}|_{J_{0}}$
and of $v|_{J_{0}}$ implies
\[
\frac{v_{j}}{u_{i,j}}=\frac{v_{j}}{u_{i}^{\left(1\right)}u_{j}^{\left(2\right)}}\leq C_{u^{\left(2\right)}|_{J_{0}},\CalP,\CalQ}\cdot C_{v|_{J_{0}},\CalP,\CalQ}\cdot\frac{v_{j_{i}}}{u_{i}^{\left(1\right)}u_{j_{i}}^{\left(2\right)}}=C_{u^{\left(2\right)}|_{J_{0}},\CalP,\CalQ}\cdot C_{v|_{J_{0}},\CalP,\CalQ}\cdot\frac{v_{j_{i}}}{u_{i,j_{i}}}
\]
and likewise
\[
\frac{v_{j_{i}}}{u_{i,j_{i}}}\leq C_{u^{\left(2\right)}|_{J_{0}},\CalP,\CalQ}\cdot C_{v|_{J_{0}},\CalP,\CalQ}\cdot\frac{v_{j}}{u_{i,j}}
\]
for all $i\in I^{\left(0\right)}$ and $j\in J_{0}\cap J_{i}$. All
in all, we derive
\[
\left\Vert \left(v_{j}/u_{i,j}\right)_{j\in J_{0}\cap J_{i}}\right\Vert _{\ell^{q_{2}\cdot\left(r/q_{2}\right)'}}\asymp\left|J_{0}\cap J_{i}\right|^{\left[q_{2}\cdot\left(r/q_{2}\right)'\right]^{-1}}\cdot\frac{v_{j_{i}}}{u_{i,j_{i}}}\asymp\left[\left|\det T_{i}\right|\big/\left|\det S_{j_{i}}\right|\right]^{\left(\frac{1}{q_{2}}-\frac{1}{r}\right)_{+}}\cdot\frac{v_{j_{i}}}{u_{i,j_{i}}}
\]
for all $i\in I^{\left(0\right)}$, with the implied constants as
in the statement of the corollary. Here, we used the formula $\frac{1}{q\cdot\left(r/q\right)'}=\left(\frac{1}{q}-\frac{1}{r}\right)_{+}$
(see equation~(\ref{eq:InverseOfSpecialExponent})) in the last step.

Since all indices $i\in I\setminus I^{\left(0\right)}$ can be neglected
in the left-hand side of the desired estimate, this completes the
proof.
\end{proof}
As a further corollary of Theorem~\ref{thm:NoSubordinatenessWithConsiderationOfOverlaps},
we derive a result which applies to coverings $\CalQ,\CalP$ that
exhibit a kind of ``mixed'' subordinateness. Roughly, we assume
that we can write $\CalO\cap\CalO'=A\cup B$, such that $\CalQ$ is
almost subordinate to $\CalP$ ``near $A$'' and vice versa ``near
$B$''. We remark that the following corollary is a generalized version
of \cite[Corollary 5.1.11]{VoigtlaenderPhDThesis} from my PhD thesis.
\begin{cor}
\label{cor:MixedSubordinateness}Let $\CalQ=\left(Q_{i}\right)_{i\in I}=\left(T_{i}Q_{i}'+b_{i}\right)_{i\in I}$
and $\CalP=\left(P_{j}\right)_{j\in J}=\left(S_{j}P_{j}'+c_{j}\right)_{j\in J}$
be two semi-structured coverings of the open sets $\emptyset\neq\CalO,\CalO'\subset\R^{\dimension}$,
respectively. Let $Y\subset\Compl^{I}$ and $Z\subset\Compl^{J}$
be $\CalQ$-regular and $\CalP$-regular sequence spaces, respectively
and let $p_{1},p_{2}\in\left(0,\infty\right]$ with $p_{1}\leq p_{2}$.
Assume that there are $L^{p_{1}}$-BAPUs $\Phi=\left(\varphi_{i}\right)_{i\in I}$
and $\Psi=\left(\psi_{j}\right)_{j\in J}$ for $\CalQ$ and $\CalP$,
respectively. Denote by $C_{Z}$ a triangle constant for the quasi-normed
space $Z\subset\Compl^{J}$.

Assume that there are subsets $A,B\subset\R^{\dimension}$ with the
following properties:

\begin{enumerate}
\item We have $\CalO\cap\CalO'=A\cup B$.
\item The family $\CalQ_{I_{A}}:=\left(Q_{i}\right)_{i\in I_{A}}$ is almost
subordinate to $\CalP$, where
\[
I_{A}:=\left\{ i\in I\with Q_{i}\cap A\neq\emptyset\right\} \,.
\]
\item The family $\CalP_{J_{B}}:=\left(P_{j}\right)_{j\in J_{B}}$ is almost
subordinate to $\CalQ$, where
\[
J_{B}:=\left\{ j\in J\with P_{j}\cap B\neq\emptyset\right\} \,.
\]
\end{enumerate}
Then, with $I_{j}=\left\{ i\in I\with Q_{i}\cap P_{j}\neq\emptyset\right\} $
and $J_{i}=\left\{ j\in J\with P_{j}\cap Q_{i}\neq\emptyset\right\} $
as usual, we have 
\[
L:=\bigcup_{j\in J\setminus J_{B}}I_{j}\subset I_{A}\:.
\]

Finally, define
\[
\theta_{i}:=\left|\det T_{i}\right|^{p_{1}^{-1}-p_{2}^{-1}}\qquad\text{ and }\qquad u_{i,j}:=\begin{cases}
\left|\det S_{j}\right|^{p_{2}^{-1}-1}\cdot\left|\det T_{i}\right|^{1-p_{1}^{-1}}, & \text{if }p_{1}<1,\\
\vphantom{\rule{0.1cm}{0.55cm}}\left|\det S_{j}\right|^{p_{2}^{-1}-p_{1}^{-1}}, & \text{if }p_{1}\geq1,
\end{cases}
\]
for $i\in I$ and $j\in J$, and assume additionally that the following
linear maps are bounded:
\begin{alignat*}{3}
\beta_{1}:\: & Y &  & \to Z|_{\!J\setminus J_{B}}\Bigl(\bigl[\ell_{\theta}^{\LowerExpo{p_{2}}}\!\bigl(I_{j}\bigr)\bigr]{}_{\!j\in J\setminus J_{B}}\Bigr), & \left(x_{i}\right)_{i\in I} & \mapsto\left(x_{i}\right)_{i\in L}\:,\\
\beta_{2}:\: & Y\!\Bigl(\bigl[\ell_{u}^{\UpperExpo{p_{1}}}\!\!\bigl(J_{i}\cap J_{B}\bigr)\bigr]{}_{i\in I}\Bigr) &  & \hookrightarrow Z|_{J_{B}}\:.
\end{alignat*}

Then, the map
\[
\iota:\FourierDecompSp{\CalQ}{p_{1}}Y\to\FourierDecompSp{\CalP}{p_{2}}Z,f\mapsto\sum_{i\in I}\varphi_{i}f
\]
is well-defined and bounded with $\vertiii{\iota}\leq C\cdot\left(\vertiii{\beta_{1}}+\vertiii{\beta_{2}}\right)$
for some constant
\[
C=C\left(\dimension,p_{1},p_{2},C_{Z},\CalQ,\CalP,k\left(\smash{\CalP_{J_{B}}},\CalQ\right),k\left(\smash{\CalQ_{I_{A}}},\CalP\right),C_{\CalQ,\Phi,p_{1}},C_{\CalP,\Psi,p_{1}},\vertiii{\Gamma_{\CalQ}}_{Y\to Y},\vertiii{\Gamma_{\CalP}}_{Z\to Z}\right)\:.
\]
The map $\iota$ has the following additional properties:

\begin{enumerate}
\item \label{enu:MixedSubordinatenessExtensionProperty}For $f\in\FourierDecompSp{\CalQ}{p_{1}}Y\subset\DistributionSpace{\CalO}$,
we have
\[
\left\langle \iota f,\,g\right\rangle _{\CalD'}=\left\langle f,\,g\right\rangle _{\CalD'}\qquad\forall\:g\in\TestFunctionSpace{\CalO\cap\CalO'}.
\]
In particular, if $\CalO=\CalO'$, then $\iota f=f$ for all $f\in\FourierDecompSp{\CalQ}{p_{1}}Y\subset\DistributionSpace{\CalO}$.
\item \label{enu:MixedSubordinatenessConsistencyForL1loc}If $f\in\FourierDecompSp{\CalQ}{p_{1}}Y$
is given by (integration against) a measurable function $f:\R^{\dimension}\to\Compl$
with $f\in L_{{\rm loc}}^{1}\left(\CalO\cup\CalO'\right)$ and with
$f=0$ almost everywhere on $\CalO'\setminus\CalO$, then $\iota f=f$
as elements of $\DistributionSpace{\CalO'}$.
\item \label{enu:MixedSubordinatenessConsistencyForTestFunctions}In particular,
$\iota f=f$ for $f\in\TestFunctionSpace{\CalO\cap\CalO'}\cap\FourierDecompSp{\CalQ}{p_{1}}Y$.
\item \label{enu:MixedSubordinatenessSupport}For $f\in\FourierDecompSp{\CalQ}{p_{1}}Y$,
we have $\supp\iota f\subset\CalO'\cap\overline{\supp f}$, where
the closure $\overline{\supp f}$ is taken in $\R^{\dimension}$.
\end{enumerate}
Finally, if $Y=\ell_{w}^{q_{1}}\left(I\right)$ and $Z=\ell_{v}^{q_{2}}\left(J\right)$
for certain $q_{1},q_{2}\in\left(0,\infty\right]$ and weights $w=\left(w_{i}\right)_{i\in I}$
and $v=\left(v_{j}\right)_{j\in J}$ which are $\CalQ$-moderate and
$\CalP$-moderate, respectively, then
\begin{align}
\vertiii{\beta_{1}} & \asymp\left\Vert \left(v_{j}\cdot\left\Vert \left(w_{i}^{-1}\cdot\left|\det T_{i}\right|^{p_{1}^{-1}-p_{2}^{-1}}\right)_{i\in I_{j}}\right\Vert _{\ell^{\LowerExpo{p_{2}}\cdot\left(q_{1}/\LowerExpo{p_{2}}\right)'}}\right)_{j\in J\setminus J_{B}}\right\Vert _{\ell^{q_{2}\cdot\left(q_{1}/q_{2}\right)'}}\label{eq:MixedSubordinatenessSimplificationBeta1Norm}\\
\vertiii{\beta_{2}} & \asymp\left\Vert \left(w_{i}^{-1}\cdot\left\Vert \left(v_{j}/u_{i,j}\right)_{j\in J_{i}\cap J_{B}}\right\Vert _{\ell^{q_{2}\cdot\left(\UpperExpo{p_{1}}/q_{2}\right)'}}\right)_{i\in I}\right\Vert _{\ell^{q_{2}\cdot\left(q_{1}/q_{2}\right)'}}\label{eq:MixedSubordinatenessSimplificationBeta2Norm}
\end{align}
where the implied constants only depend on 
\[
\dimension,p_{1},p_{2},q_{1},q_{2},\CalQ,\CalP,k\left(\smash{\CalQ_{I_{A}}},\CalP\right),k\left(\smash{\CalP_{J_{B}}},\CalQ\right),C_{w,\CalQ},C_{v,\CalP}.\qedhere
\]
\end{cor}

\begin{rem*}
Note that the corollary in particular yields $I_{j}\subset L\subset I_{A}$
for all $j\in J\setminus J_{B}$. Thus, under the assumptions of the
last part of the corollary, we have
\begin{align*}
\vertiii{\beta_{1}} & \asymp\left\Vert \left(v_{j}\cdot\left\Vert \left(w_{i}^{-1}\cdot\left|\det T_{i}\right|^{p_{1}^{-1}-p_{2}^{-1}}\right)_{i\in I_{j}}\right\Vert _{\ell^{\LowerExpo{p_{2}}\cdot\left(q_{1}/\LowerExpo{p_{2}}\right)'}}\right)_{j\in J\setminus J_{B}}\right\Vert _{\ell^{q_{2}\cdot\left(q_{1}/q_{2}\right)'}}\\
 & \leq\left\Vert \left(v_{j}\cdot\left\Vert \left(w_{i}^{-1}\cdot\left|\det T_{i}\right|^{p_{1}^{-1}-p_{2}^{-1}}\right)_{i\in I_{j}\cap I_{A}}\right\Vert _{\ell^{\LowerExpo{p_{2}}\cdot\left(q_{1}/\LowerExpo{p_{2}}\right)'}}\right)_{j\in J}\right\Vert _{\ell^{q_{2}\cdot\left(q_{1}/q_{2}\right)'}},
\end{align*}
which is more similar to the estimate given for $\vertiii{\beta_{2}}$.
Note, though, that this only yields an \emph{upper} bound for $\vertiii{\beta_{1}}$.
\end{rem*}
\begin{proof}
To begin with, we show $L\subset I_{A}$. Indeed, if $i\in I_{j}$
for some $j\in J\setminus J_{B}$, there is some $\xi\in Q_{i}\cap P_{j}\subset\CalO\cap\CalO'=A\cup B$.
In case of $\xi\in B$, we would have $\xi\in P_{j}\cap B$ and hence
$j\in J_{B}$, a contradiction. Hence, $\xi\in A$, i.e.\@ $\xi\in Q_{i}\cap A\neq\emptyset$
and therefore $i\in I_{A}$.

\medskip{}

Now, let us begin with the actual proof. Below, we will use Theorem~\ref{thm:NoSubordinatenessWithConsiderationOfOverlaps}
to show that the following two maps are well-defined and bounded:
\begin{align*}
\iota_{1}:\FourierDecompSp{\CalQ}{p_{1}}Y\to\FourierDecompSp{\CalP}{p_{2}}Z, & f\mapsto\sum_{\left(i,j\right)\in I_{A}\times\left(J\setminus J_{B}\right)}\varphi_{i}\,\psi_{j}\,f,\\
\iota_{2}:\FourierDecompSp{\CalQ}{p_{1}}Y\to\FourierDecompSp{\CalP}{p_{2}}Z, & f\mapsto\sum_{\left(i,j\right)\in I\times J_{B}}\varphi_{i}\,\psi_{j}\,f,
\end{align*}
both with absolute convergence of the series defining $\left\langle \iota_{1}f,\,g\right\rangle _{\CalD'}$
and $\left\langle \iota_{2}f,\,g\right\rangle _{\CalD'}$ for $g\in\TestFunctionSpace{\CalO'}$.
Let us assume for the moment that this holds.

We first show that $\iota_{1}f=\sum_{\left(i,j\right)\in I\times\left(J\setminus J_{B}\right)}\varphi_{i}\,\psi_{j}\,f$
for $f\in\FourierDecompSp{\CalQ}{p_{1}}Y$. To this end, it suffices
to show $\varphi_{i}\,\psi_{j}\equiv0$ for $i\in I\setminus I_{A}$
and $j\in J\setminus J_{B}$. But $\varphi_{i}\,\psi_{j}\not\equiv0$
would imply $i\in I_{j}\subset L\subset I_{A}$, since $j\in J\setminus J_{B}$.
Since this is impossible for $i\in I\setminus I_{A}$, we have shown
$\varphi_{i}\,\psi_{j}\equiv0$ for $i\in I\setminus I_{A}$ and $j\in J\setminus J_{B}$,
as needed. Therefore, we see for any $f\in\FourierDecompSp{\CalQ}{p_{1}}Y$
that
\[
\iota_{1}f+\iota_{2}f=\sum_{\left(i,j\right)\in I\times\left(J\setminus J_{B}\right)}\varphi_{i}\,\psi_{j}\,f+\sum_{\left(i,j\right)\in I\times J_{B}}\varphi_{i}\,\psi_{j}\,f=\sum_{\left(i,j\right)\in I\times J}\varphi_{i}\,\psi_{j}\,f,
\]
again with absolute convergence of the defining series, since all
we did was to add some vanishing terms and to add two absolutely convergent
(defining) series.

But for $g\in\TestFunctionSpace{\CalO'}$, we have $g=\sum_{j\in J}\psi_{j}g$,
with only finitely many non-vanishing terms, since $\left(\psi_{j}\right)_{j\in J}$
is a locally finite partition of unity on $\CalO'$; see Lemma~\ref{lem:PartitionCoveringNecessary}.
This implies
\[
\left\langle \iota_{1}f+\iota_{2}f,\,g\right\rangle _{\CalD'}=\sum_{i\in I}\,\sum_{j\in J}\left\langle \varphi_{i}\,f,\,\psi_{j}\,g\right\rangle _{\Schwartz'}=\sum_{i\in I}\left\langle \varphi_{i}\,f,\,g\right\rangle _{\Schwartz'}\;,
\]
with absolute convergence of the series. Hence, we see that the map
$\iota$ as defined in the present corollary is well-defined and bounded,
with absolute convergence of the defining series.

\medskip{}

It remains to establish the additional properties of $\iota$ (and
the boundedness of $\iota_{1},\iota_{2}$). Validity of property (\ref{enu:MixedSubordinatenessExtensionProperty})
is shown just as above: For $g\in\TestFunctionSpace{\CalO\cap\CalO'}$,
we have $g=\sum_{i\in I}\varphi_{i}\,g$ (with only finitely many
terms not vanishing) and this easily yields $\left\langle \iota f,\,g\right\rangle _{\CalD'}=\left\langle f,\,g\right\rangle _{\CalD'}$.

\medskip{}

Property (\ref{enu:MixedSubordinatenessConsistencyForTestFunctions})
is an immediate consequence of property (\ref{enu:MixedSubordinatenessConsistencyForL1loc}),
which we prove next: For each $\xi\in\CalO$, there is some $i_{\xi}\in I$
with $\varphi_{i_{\xi}}\left(\xi\right)\neq0$ and hence $\xi\in Q_{i_{\xi}}^{\circ}$.
Since $\CalO$ is second countable, there is thus a sequence $\left(i_{\xi_{n}}\right)_{n\in\N}$
satisfying $\CalO=\bigcup_{n\in\N}Q_{i_{\xi_{n}}}^{\circ}$. But this
implies $\varphi_{i}\equiv0$ for all $i\in I\setminus I_{0}$, with
the \emph{countable} set $I_{0}:=\bigcup_{n\in\N}i_{\xi_{n}}^{\ast}$.
Indeed, if $\varphi_{i}\not\equiv0$, then $\varphi_{i}\left(\xi\right)\neq0$
for some $\xi\in Q_{i}\subset\CalO$. By what we just saw, $\xi\in Q_{i_{\xi_{n}}}$
for some $n\in\N$, so that we get $Q_{i}\cap Q_{i_{\xi_{n}}}\neq\emptyset$
and hence $i\in i_{\xi_{n}}^{\ast}\subset I_{0}$.

Now, recall that Fourier inversion implies
\[
\left|\varphi_{i}\left(\xi\right)\right|=\left|\widehat{\Fourier^{-1}\varphi_{i}}\left(\xi\right)\right|\leq\left\Vert \Fourier^{-1}\varphi_{i}\right\Vert _{L^{1}}\leq C_{\Phi}<\infty
\]
for all $\xi\in\R^{\dimension}$ and $i\in I$. As we saw in the previous
paragraph, for each $\xi\in\CalO$, there is some $n=n_{\xi}\in\N$
with $\varphi_{i}\left(\xi\right)=0$ for $i\in I\setminus i_{\xi_{n}}^{\ast}$.
This implies
\[
\sum_{i\in I_{0}}\left|\varphi_{i}\left(\xi\right)\right|=\sum_{i\in i_{\xi_{n}}^{\ast}}\left|\varphi_{i}\left(\xi\right)\right|\leq C_{\Phi}\left|i_{\xi_{n}}^{\ast}\right|\leq C_{\Phi}N_{\CalQ}<\infty
\]
for all $\xi\in\CalO$. For $\xi\in\R^{\dimension}\setminus\CalO$,
the left-hand side vanishes, so that the estimate is true for all
$\xi\in\R^{\dimension}$.

Now, let $g\in\TestFunctionSpace{\CalO'}$ be arbitrary and set $K:=\supp g$.
Note that $K\subset\CalO\cup\CalO'$ is compact, so that $\Indicator_{K}\cdot f\in L^{1}\left(\R^{\dimension}\right)$,
since $f\in L_{{\rm loc}}^{1}\left(\CalO\cup\CalO'\right)$. Furthermore,
\[
\sum_{i\in I_{0}}\left|\varphi_{i}\left(\xi\right)f\left(\xi\right)g\left(\xi\right)\right|\leq C_{\Phi}N_{\CalQ}\left\Vert g\right\Vert _{\sup}\cdot\left(\Indicator_{K}\cdot\left|f\right|\right)\left(\xi\right)\in L^{1}\left(\R^{\dimension}\right).
\]
Thus, the dominated convergence theorem implies as desired that
\begin{align*}
\left\langle \iota f,\,g\right\rangle _{\CalD'} & =\sum_{i\in I}\left\langle \varphi_{i}\,f,\,g\right\rangle _{\Schwartz'}\\
\left({\scriptstyle \varphi_{i}\equiv0\text{ for }i\in I\setminus I_{0}}\right) & =\sum_{i\in I_{0}}\int_{\R^{\dimension}}\varphi_{i}\left(\xi\right)f\left(\xi\right)g\left(\xi\right)\,\d\xi\\
\left({\scriptstyle I_{0}\text{ countable, dominated conv.}}\right) & =\int_{\R^{\dimension}}\left[\,\smash{\sum_{i\in I_{0}}}\,\vphantom{\sum}\varphi_{i}\left(\xi\right)\,\right]\vphantom{\sum_{i\in I_{0}}}\cdot f\left(\xi\right)g\left(\xi\right)\,\d\xi\\
\left({\scriptstyle \varphi_{i}\equiv0\text{ for }i\in I\setminus I_{0}}\right) & =\int_{\R^{\dimension}}\left[\,\smash{\sum_{i\in I}}\,\vphantom{\sum}\varphi_{i}\left(\xi\right)\,\right]\vphantom{\sum_{i\in I}}\cdot f\left(\xi\right)g\left(\xi\right)\,\d\xi\\
\left({\scriptstyle \sum_{i\in I}\varphi_{i}\equiv1\text{ on }\CalO,\:\varphi_{i}\equiv0\text{ on }\R^{\dimension}\setminus\CalO}\right) & =\int_{\CalO}f\left(\xi\right)g\left(\xi\right)\,\d\xi\\
\left({\scriptstyle f=0\text{ a.e. on }\CalO'\setminus\CalO\text{ and }g=0\text{ on }\R^{\dimension}\setminus\CalO'}\right) & =\int_{\R^{\dimension}}f\left(\xi\right)g\left(\xi\right)\,\d\xi=\left\langle f,\,g\right\rangle _{\CalD'}\,\,.
\end{align*}

\medskip{}

Lastly, we establish property (\ref{enu:MixedSubordinatenessSupport}):
Let $f\in\FourierDecompSp{\CalQ}{p_{1}}Y$ be arbitrary and set $S:=\overline{\supp f}$,
where the closure is taken in $\R^{\dimension}$. Since $S\subset\R^{\dimension}$
is closed, $U:=\CalO'\setminus S$ is open. Now, let $g\in\TestFunctionSpace{\CalO'}$
with $\supp g\subset U$. Because of $\supp\varphi_{i}\subset\CalO$
for all $i\in I$, we have
\[
\supp\left(\varphi_{i}\,g\right)\subset\CalO\cap U=\CalO\cap\left(\CalO'\setminus S\right)\subset\CalO\setminus S=\CalO\setminus\overline{\supp f}\subset\CalO\setminus\supp f,
\]
where we note that $\CalO\setminus\supp f$ is open, since $\supp f$
is closed in $\CalO$, because of $f\in\DistributionSpace{\CalO}$.
Since $f$ vanishes on the open set $\CalO\setminus\supp f$, we get
$\left\langle f,\,\varphi_{i}\,g\right\rangle _{\CalD'}=0$ for all
$i\in I$, and thus
\[
\left\langle \iota f,\,g\right\rangle _{\CalD'}=\sum_{i\in I}\left\langle \varphi_{i}\,f,\,g\right\rangle _{\CalD'}=0.
\]

All in all, $\iota f$ vanishes on the open(!) set $U$, so that we
get as claimed that
\[
\supp\left(\iota f\right)\subset\CalO'\setminus U=\CalO'\setminus\left[\CalO'\setminus S\right]\subset\CalO'\cap S=\CalO'\cap\overline{\supp f}\,.
\]

\medskip{}

To complete the proof, we finally prove the boundedness of $\iota_{1}$
and $\iota_{2}$.

\textbf{Boundedness of $\iota_{1}$}: Here, we apply Theorem~\ref{thm:NoSubordinatenessWithConsiderationOfOverlaps}
with $I_{0}:=I_{A}\subset I$ and $K:=J_{0}:=J\setminus J_{B}$ and
with
\[
J^{\left(k\right)}:=\left\{ k\right\} ,\qquad\text{ as well as }\qquad I^{\left(k\right)}:=I_{k}=\left\{ i\in I\with Q_{i}\cap P_{k}\neq\emptyset\right\} \qquad\text{ for }k\in K.
\]
Recall from the beginning of the proof that we have $I^{\left(k\right)}=I_{k}\subset I_{A}=I_{0}$
for all $k\in K=J\setminus J_{B}$, as required in Theorem~\ref{thm:NoSubordinatenessWithConsiderationOfOverlaps}.
Furthermore, note that we have (in the notation of Theorem~\ref{thm:NoSubordinatenessWithConsiderationOfOverlaps})
\[
J_{00}=\bigcup_{k\in K}J^{\left(k\right)}=\bigcup_{k\in K}\left\{ k\right\} =K=J_{0}
\]
and hence $J_{0}\subset J_{00}$, as required.

Still in the notation of Theorem~\ref{thm:NoSubordinatenessWithConsiderationOfOverlaps},
we select $q_{k}:=p_{2}\in\left[p_{1},p_{2}\right]$ for all $k\in K$.
This entails $q^{\left(0\right)}=\inf_{k\in K}q_{k}=p_{2}$, so that
$\Psi=\left(\psi_{j}\right)_{j\in J}$ is indeed an $L^{q^{\left(0\right)}}$-BAPU
for $\CalP$, as required in Theorem~\ref{thm:NoSubordinatenessWithConsiderationOfOverlaps}.
Here, we used Corollary~\ref{cor:LpBAPUsAreAlsoLqBAPUsForLargerq}
to conclude that the $L^{p_{1}}$-BAPU $\Psi$ is also an $L^{p_{2}}$-BAPU,
since $p_{1}\leq p_{2}$.

For the application of Theorem~\ref{thm:NoSubordinatenessWithConsiderationOfOverlaps},
we first verify condition~(\ref{eq:SpecialIndexSetConditionPSelected}).
Thus, assume towards a contradiction that for some $k\in K=J\setminus J_{B}$,
there is some
\[
\xi\in\vphantom{\bigcup_{i\in I_{0}\setminus I^{\left(k\right)}}}\left(\,\smash{\bigcup_{i\in I_{0}\setminus I^{\left(k\right)}}}\vphantom{\bigcup}Q_{i}\,\right)\cap\left(\,\smash{\bigcup_{j\in J^{\left(k\right)}}}\vphantom{\bigcup}P_{j}\,\right)=\bigcup_{i\in I_{A}\setminus I_{k}}\left(Q_{i}\cap P_{k}\right).
\]
This yields some $i\in I_{A}\setminus I_{k}$ with $\xi\in Q_{i}\cap P_{k}\neq\emptyset$,
in contradiction to $i\notin I_{k}$.

Next, we estimate the weights from equations~(\ref{eq:WeightForJDefinition})
and (\ref{eq:WeightForIDefinition}). For the weight $w$, we have—because
of $q_{k}=p_{2}$—for $k\in K=J\setminus J_{B}$ and $j\in J^{\left(k\right)}=\left\{ k\right\} $
that
\begin{equation}
w_{k,j}=w_{k}^{\left(0\right)}:=\begin{cases}
\left|\det S_{k}\right|^{p_{2}^{-1}-1}=\left|\det S_{j}\right|^{p_{2}^{-1}-1}, & \text{if }p_{2}<1,\\
1, & \text{if }p_{2}\geq1.
\end{cases}\label{eq:MixedSubordinatenessWChoice}
\end{equation}

Estimating $v$ is more involved. Here, we first estimate for $k\in K=J\setminus J_{B}$
and $i\in I^{\left(k\right)}=I_{k}$ the quantity
\[
\sup_{j\in J^{\left(k\right)}}\lambda\left(\,\overline{P_{j}}-\overline{Q_{i}}\,\right)=\lambda\left(\,\overline{P_{k}}-\overline{Q_{i}}\,\right).
\]
As noted above, we have $i\in I_{k}\subset I_{A}$. By assumption,
this means $Q_{i}\subset P_{j_{i}}^{n\ast}$ for some $j_{i}\in J$,
where $n:=k\left(\smash{\CalQ_{I_{A}}},\CalP\right)$. But because
of $i\in I_{k}$, we also have $\emptyset\neq P_{k}\cap Q_{i}\subset P_{k}\cap P_{j_{i}}^{n\ast}$
and hence $j_{i}\in k^{\left(n+1\right)\ast}$. In view of Corollary~\ref{cor:SemiStructuredDifferenceSetsMeasureEstimate},
this yields
\[
\lambda\left(\,\overline{P_{k}}-\overline{Q_{i}}\,\right)\leq\lambda\left(\,\overline{P_{k}}-\overline{P_{j_{i}}^{n\ast}}\,\right)\leq\lambda\left(\,\overline{P_{k}^{\left(n+1\right)\ast}}-\overline{P_{j_{i}}^{\left(n+1\right)\ast}}\,\right)\leq C_{1}\cdot\left|\det S_{k}\right|,
\]
for some constant $C_{1}=C_{1}\left(\dimension,n,\CalP\right)$.

All in all, this shows $v_{k,i}\leq C_{2}\cdot v_{k,i}^{\left(0\right)}$
for all $k\in K$ and $i\in I^{\left(k\right)}=I_{k}$ for a suitable
constant $C_{2}=C_{2}\left(p_{2},\dimension,n,\CalP\right)$ and
\[
v_{k,i}^{\left(0\right)}:=\begin{cases}
\left|\det T_{i}\right|^{p_{1}^{-1}-p_{2}^{-1}}\cdot\left|\det S_{k}\right|^{p_{2}^{-1}-1}, & \text{if }p_{2}<1,\\
\vphantom{\rule{0.1cm}{0.55cm}}\left|\det T_{i}\right|^{p_{1}^{-1}-p_{2}^{-1}}, & \text{if }p_{2}\geq1.
\end{cases}
\]

It remains to establish boundedness of the embeddings $\eta_{1},\eta_{2}$
from equations (\ref{eq:AssumedDiscreteEmbedding1}) and (\ref{eq:AssumedDiscreteEmbedding2}),
for a suitable solid sequence space $X\subset\Compl^{K}$. Here, we
choose $X:=\left(Z|_{K}\right)_{1/w^{\left(0\right)}}$, where we
recall $K=J\setminus J_{B}$ and the definition of $w^{\left(0\right)}=\left(\smash{w_{j}^{\left(0\right)}}\right)_{j\in J}$
from equation~(\ref{eq:MixedSubordinatenessWChoice}).

It is not hard to see that the embedding from equation~(\ref{eq:AssumedDiscreteEmbedding1})
is just $X\left(\left[\vphantom{\ell_{w}^{q}}\smash{\ell_{w}^{\UpperExpo{q_{k}}}}\!\left(\smash{J^{\left(k\right)}}\right)\right]_{k\in K}\right)\hookrightarrow Z|_{J_{00}}$,
where in our case $J_{00}=J_{0}=K=J\setminus J_{B}$. But since we
have $J^{\left(k\right)}=\left\{ k\right\} $ for $k\in K$ and because
of $w_{k,j}=w_{k}^{\left(0\right)}$, we easily see 
\[
X\Bigl(\bigl[\ell_{w}^{\UpperExpo{q_{k}}}\!\bigl(J^{\left(k\right)}\bigr)\bigr]_{k\in K}\Bigr)=\left(Z|_{K}\right)_{w^{\left(0\right)}/w^{\left(0\right)}}\Bigl(\bigl[\ell^{\UpperExpo{q_{k}}}\!\left(\left\{ k\right\} \right)\bigr]_{k\in K}\Bigr)=Z|_{K},
\]
so that the embedding $\eta_{1}$ from equation~(\ref{eq:AssumedDiscreteEmbedding1})
is trivially bounded, with $\vertiii{\eta_{1}}\leq1$.

For the embedding $\eta_{2}$ from equation~(\ref{eq:AssumedDiscreteEmbedding2}),
first note that we have
\begin{align*}
v_{k,i}^{\left(0\right)}/w_{k}^{\left(0\right)} & =\begin{cases}
\left|\det T_{i}\right|^{p_{1}^{-1}-p_{2}^{-1}}\cdot\left|\det S_{k}\right|^{p_{2}^{-1}-1}/\left|\det S_{k}\right|^{p_{2}^{-1}-1}, & \text{if }p_{2}<1,\\
\vphantom{\rule{0.1cm}{0.55cm}}\left|\det T_{i}\right|^{p_{1}^{-1}-p_{2}^{-1}}/\,1, & \text{if }p_{2}\geq1
\end{cases}\\
 & =\left|\det T_{i}\right|^{p_{1}^{-1}-p_{2}^{-1}}=\theta_{i}
\end{align*}
for all $k\in K=J\setminus J_{B}$ and all $i\in I^{\left(k\right)}=I_{k}$,
with $\theta_{i}$ as in the statement of the present corollary. Thus,
our estimate $v_{k,i}\leq C_{2}\cdot v_{k,i}^{\left(0\right)}$ from
above, together with the boundedness of the map $\beta_{1}$ from
the statement of the present corollary, shows that we have
\begin{eqnarray*}
Y & \Xhookrightarrow{\beta_{1}} & Z|_{J\setminus J_{B}}\Bigl(\bigl[\ell_{\theta}^{\LowerExpo{p_{2}}}\!\bigl(I_{j}\bigr)\bigr]_{j\in J\setminus J_{B}}\Bigr)\\
 & = & X_{w^{\left(0\right)}}\biggl(\Bigl[\ell_{\left(v_{k,i}^{\left(0\right)}/w_{k}^{\left(0\right)}\right)_{i}}^{\LowerExpo{p_{2}}}\bigl(I^{\left(k\right)}\bigr)\Bigr]_{k\in J\setminus J_{B}}\biggr)\\
 & = & X\Bigl(\Bigl[\ell_{v^{\left(0\right)}}^{\LowerExpo{q_{k}}}\!\bigl(I^{\left(k\right)}\bigr)\Bigr]_{k\in K}\Bigr)\\
\left({\scriptstyle \text{since }v\leq C_{2}\cdot v^{\left(0\right)}}\right) & \hookrightarrow & X\Bigl(\Bigl[\bigl[\ell_{v}^{\LowerExpo{q_{k}}}\!\bigl(I^{\left(k\right)}\bigr)\bigr]{}_{k\in K}\Bigr),
\end{eqnarray*}
as required. The quantitative version of these considerations yields
$\vertiii{\eta_{2}}\leq C_{2}\cdot\vertiii{\beta_{1}}$.

All in all, Theorem~\ref{thm:NoSubordinatenessWithConsiderationOfOverlaps}
shows that $\iota_{1}$ is well-defined and bounded and satisfies
\[
\vertiii{\iota_{1}}\leq C_{3}\cdot\vertiii{\eta_{1}}\cdot\vertiii{\eta_{2}}\leq C_{2}C_{3}\cdot\vertiii{\beta_{1}}
\]
for some constant
\[
C_{3}=C_{3}\left(\dimension,p_{1},p_{2},\CalQ,\CalP,C_{\CalQ,\Phi,p_{1}},C_{\CalP,\Psi,p_{2}},\vertiii{\Gamma_{\CalP}}_{Z\to Z}\right).
\]
Note that $C_{\CalP,\Psi,p_{2}}$ can be estimated in terms of $C_{\CalP,\Psi,p_{1}}$
and $\dimension,p_{1},\CalP$, thanks to Corollary~\ref{cor:LpBAPUsAreAlsoLqBAPUsForLargerq}.

\medskip{}

\textbf{Boundedness of $\iota_{2}$}: Here, we apply Theorem~\ref{thm:NoSubordinatenessWithConsiderationOfOverlaps}
with $K:=I_{0}:=I$ and $J_{0}:=J_{B}$ and with
\[
J^{\left(k\right)}:=J_{B}\cap J_{k}=\left\{ j\in J_{B}\with P_{j}\cap Q_{k}\neq\emptyset\right\} ,\qquad\text{ as well as }\qquad I^{\left(k\right)}:=k^{\left(2m+2\right)\ast}\qquad\text{ for }k\in K,
\]
where $m:=k\left(\smash{\CalP_{J_{B}}},\CalQ\right)$. Furthermore,
we select $q_{k}:=p_{1}\in\left[p_{1},p_{2}\right]$ for all $k\in K$.
This yields $q^{\left(0\right)}=\inf_{k\in K}q_{k}=p_{1}$. Note that
$\Psi$ is an $L^{p_{1}}$-BAPU for $\CalP$, as needed for applying
Theorem~\ref{thm:NoSubordinatenessWithConsiderationOfOverlaps}.

Observe $J_{0}=J_{B}\subset\bigcup_{k\in K}J^{\left(k\right)}=J_{00}$,
as required in Theorem~\ref{thm:NoSubordinatenessWithConsiderationOfOverlaps};
indeed, for $j\in J_{B}$, there is some $\xi\in P_{j}\cap B\subset A\cup B=\CalO\cap\CalO'$.
But since $\CalQ=\left(Q_{i}\right)_{i\in I}$ covers $\CalO$, there
is some $k\in I$ satisfying $\xi\in Q_{k}$. Hence, $\xi\in P_{j}\cap Q_{k}\neq\emptyset$,
which yields $j\in J_{B}\cap J_{k}=J^{\left(k\right)}$.

Now, let us show that condition~(\ref{eq:SpecialIndexSetConditionPSelected})
holds with our choices from above. Hence, assume towards a contradiction
that there is some
\[
\xi\in\vphantom{\bigcup_{i\in I_{0}\setminus I^{\left(k\right)}}}\left(\,\vphantom{\bigcup}\smash{\bigcup_{i\in I_{0}\setminus I^{\left(k\right)}}}Q_{i}\,\right)\cap\left(\,\smash{\bigcup_{j\in J^{\left(k\right)}}}\vphantom{\bigcup}P_{j}\,\right)=\left(\,\smash{\bigcup_{i\in I\setminus k^{\left(2m+2\right)\ast}}}Q_{i}\,\vphantom{\bigcup}\right)\cap\left(\,\vphantom{\bigcup}\smash{\bigcup_{j\in J_{B}\cap J_{k}}}P_{j}\,\right).
\]
Thus, there are $i\in I\setminus k^{\left(2m+2\right)\ast}$ and $j\in J_{k}\cap J_{B}$
with $\xi\in Q_{i}\cap P_{j}$. Because of $j\in J_{B}$, our assumptions
imply $\xi\in P_{j}\subset Q_{i_{j}}^{m\ast}$ for some $i_{j}\in I$.
But because of $\xi\in Q_{i}$, this implies $Q_{i}\cap Q_{i_{j}}^{m\ast}\neq\emptyset$
and thus $i_{j}\in i^{\left(m+1\right)\ast}$. Finally, since $j\in J_{k}$,
we get $\emptyset\neq P_{j}\cap Q_{k}\subset Q_{i_{j}}^{m\ast}\cap Q_{k}$
and hence $k\in i_{j}^{\left(m+1\right)\ast}\subset i^{\left(2m+2\right)\ast}$,
i.e., $i\in k^{\left(2m+2\right)\ast}$, in contradiction to $i\in I\setminus k^{\left(2m+2\right)\ast}$.

Next, we estimate the weights from equations~(\ref{eq:WeightForJDefinition})
and (\ref{eq:WeightForIDefinition}). For the weight $w$ from equation~(\ref{eq:WeightForJDefinition}),
we note for $k\in K$ and $j\in J^{\left(k\right)}=J_{k}\cap J_{B}$
(because of $q_{k}=p_{1}$) that
\[
w_{k,j}=w_{j}^{\left(0\right)}:=\begin{cases}
\left|\det S_{j}\right|^{p_{2}^{-1}-1}, & \text{if }p_{1}<1,\\
\left|\det S_{j}\right|^{p_{2}^{-1}-p_{1}^{-1}}, & \text{if }p_{1}\geq1.
\end{cases}
\]

For the weight $v$ from equation~(\ref{eq:WeightForIDefinition}),
we have to work slightly harder. First note for $k\in K$ and $i\in I^{\left(k\right)}=k^{\left(2m+2\right)\ast}$,
as well as $j\in J^{\left(k\right)}=J_{k}\cap J_{B}$ that (by our
assumptions) $P_{j}\subset Q_{i_{j}}^{m\ast}$ for some $i_{j}\in I$.
Furthermore, there is some $\xi\in P_{j}\cap Q_{k}\subset Q_{i_{j}}^{m\ast}\cap Q_{k}$,
so that $i_{j}\in k^{\left(m+1\right)\ast}\subset i^{\left(3m+3\right)\ast}$.
Hence, $P_{j}\subset Q_{i_{j}}^{m\ast}\subset Q_{i}^{\left(4m+3\right)\ast}$.
All in all, we get
\[
\lambda\left(\,\overline{P_{j}}-\overline{Q_{i}}\,\right)\leq\lambda\left(\,\overline{Q_{i}^{\left(4m+3\right)\ast}}-\overline{Q_{i}^{\left(4m+3\right)\ast}}\,\right)\leq C_{4}\cdot\left|\det T_{i}\right|,
\]
for some constant $C_{4}=C_{4}\left(\dimension,m,\CalQ\right)\geq1$
which is provided by Corollary~\ref{cor:SemiStructuredDifferenceSetsMeasureEstimate}.
All in all, we derive (because of $q_{k}=p_{1}$) that
\[
v_{k,i}=\begin{cases}
\left[\sup_{j\in J^{\left(k\right)}}\lambda\left(\,\overline{P_{j}}-\overline{Q_{i}}\,\right)\right]^{p_{1}^{-1}-1}\leq C_{4}^{p_{1}^{-1}-1}\cdot\left|\det T_{i}\right|^{p_{1}^{-1}-1}, & \text{if }p_{1}<1,\\
1, & \text{if }p_{1}\geq1.
\end{cases}
\]
We denote the right-hand side of this estimate by
\[
v_{i}^{\left(0\right)}:=\begin{cases}
\left|\det T_{i}\right|^{\frac{1}{p_{1}}-1}, & \text{if }p_{1}<1,\\
1, & \text{if }p_{1}\geq1,
\end{cases}
\]
so that we have shown $v_{k,i}\leq C_{5}\cdot v_{i}^{\left(0\right)}$
for all $k\in K$, $i\in I^{\left(k\right)}$ and some constant $C_{5}=C_{5}\left(\dimension,p_{1},m,\CalQ\right)$.
Equation~(\ref{eq:DeterminantIsModerate}) shows that $v^{\left(0\right)}$
is $\CalQ$-moderate, with $C_{v^{\left(0\right)},\CalQ}\leq C_{6}=C_{6}\left(\dimension,p_{1},\CalQ\right)$.
Hence, Lemma~\ref{lem:ModeratelyWeightedSpacesAreRegular} shows
that $X:=Y_{1/v^{\left(0\right)}}\subset\Compl^{I}=\Compl^{K}$ is
$\CalQ$-regular, with 
\[
\vertiii{\Gamma_{\CalQ}}_{X\to X}\leq C_{v^{\left(0\right)},\CalQ}\cdot\vertiii{\Gamma_{\CalQ}}_{Y\to Y}\leq C_{6}\cdot\vertiii{\Gamma_{\CalQ}}_{Y\to Y}=:C_{7}.
\]

It remains to verify boundedness of the embeddings $\eta_{1},\eta_{2}$
from equations~(\ref{eq:AssumedDiscreteEmbedding1}) and (\ref{eq:AssumedDiscreteEmbedding2})
of Theorem~\ref{thm:NoSubordinatenessWithConsiderationOfOverlaps},
with the preceding choice of $X=Y_{1/v^{\left(0\right)}}$. To this
end, note that $Y\hookrightarrow X_{v^{\left(0\right)}}$ and that
$X=Y_{1/v^{\left(0\right)}}$ is $\CalQ$-moderate, as we just saw.
Hence, Lemma~\ref{lem:EmbeddingInMixedSpaceWithSmallInnerSets} yields
\[
Y\hookrightarrow X\Bigl(\Bigl[\ell_{v^{\left(0\right)}}^{\LowerExpo{p_{1}}}\bigl(k^{\left(2m+2\right)\ast}\bigr)\Bigr]_{k\in I}\Bigr)\overset{\left(\ast\right)}{\hookrightarrow}X\Bigl(\bigl[\ell_{v}^{\LowerExpo{q_{k}}}\!\bigl(I^{\left(k\right)}\bigr)\bigr]_{k\in K}\Bigr).
\]
Here, the step marked with $\left(\ast\right)$ used $v_{k,i}\leq C_{5}\cdot v_{i}^{\left(0\right)}$,
as well as $I^{\left(k\right)}=k^{\left(2m+2\right)\ast}$ and $q_{k}=p_{1}$
for all $k\in K=I$. But this embedding precisely yields boundedness
of $\eta_{2}$ from equation~(\ref{eq:AssumedDiscreteEmbedding2}).
One can easily make the preceding arguments quantitative (see Lemma~\ref{lem:EmbeddingInMixedSpaceWithSmallInnerSets}),
to derive $\vertiii{\eta_{2}}\leq C_{8}$ for some constant $C_{8}=C_{8}\left(\dimension,p_{1},m,\CalQ,\vertiii{\Gamma_{\CalQ}}_{Y\to Y}\right)$
.

Finally, for the embedding $\eta_{1}$ from equation~(\ref{eq:AssumedDiscreteEmbedding1}),
we note for $k\in K$ and $j\in J^{\left(k\right)}=J_{B}\cap J_{k}$
that
\[
\frac{w_{k,j}}{v_{k}^{\left(0\right)}}=\frac{w_{j}^{\left(0\right)}}{v_{k}^{\left(0\right)}}=\begin{cases}
\left|\det S_{j}\right|^{p_{2}^{-1}-1}/\left|\det T_{k}\right|^{p_{1}^{-1}-1}=u_{k,j}, & \text{if }p_{1}<1,\\
\vphantom{\rule{0.1cm}{0.55cm}}\left|\det S_{j}\right|^{p_{2}^{-1}-p_{1}^{-1}}=u_{k,j}, & \text{if }p_{1}\geq1.
\end{cases}
\]
Thus, noting $J^{\left(k\right)}=J_{B}\cap J_{k}$ and recalling the
definition of the (assumed) embedding $\beta_{2}$ from the statement
of the present corollary, we derive
\[
X\Bigl(\bigl[\ell_{w}^{\UpperExpo{q_{k}}}\!\bigl(J^{\left(k\right)}\bigr)\bigr]_{k\in K}\Bigr)=Y_{1/v^{\left(0\right)}}\Bigl(\bigl[\ell_{w}^{\UpperExpo{p_{1}}}\!\bigl(J^{\left(k\right)}\bigr)\bigr]_{k\in K}\Bigr)=Y\Bigl(\bigl[\ell_{u}^{\UpperExpo{p_{1}}}\!\bigl(J_{i}\cap J_{B}\bigr)\bigr]_{i\in I}\Bigr)\Xhookrightarrow{\beta_{2}}Z|_{J_{B}},
\]
which precisely yields boundedness of $\eta_{1}$ from equation~(\ref{eq:AssumedDiscreteEmbedding1}),
since we have $J_{B}=J_{00}$, as seen above. Quantitatively, we get
$\vertiii{\eta_{1}}\leq\vertiii{\beta_{2}}$.

All in all, Theorem~\ref{thm:NoSubordinatenessWithConsiderationOfOverlaps}
shows that $\iota_{2}$ as defined above is bounded; more precisely,
Theorem~\ref{thm:NoSubordinatenessWithConsiderationOfOverlaps} yields
$\vertiii{\iota_{2}}\leq C_{9}\cdot\vertiii{\eta_{1}}\cdot\vertiii{\eta_{2}}\leq C_{8}C_{9}\cdot\vertiii{\beta_{2}}$
for some constant
\[
C_{9}=C_{9}\left(\dimension,p_{1},p_{2},\CalQ,\CalP,C_{\CalQ,\Phi,p_{1}},C_{\CalP,\Psi,p_{1}},\vertiii{\Gamma_{\CalP}}_{Z\to Z}\right).
\]
Altogether, this yields the desired estimate 
\[
\vertiii{\iota}\leq C_{10}\cdot\left[\vertiii{\iota_{1}}+\vertiii{\iota_{2}}\right]\lesssim\vertiii{\beta_{1}}+\vertiii{\beta_{2}},
\]
with an implied constant as in the statement of the corollary. Here,
$C_{10}$ is the triangle constant for $\FourierDecompSp{\CalP}{p_{2}}Z$,
which can be bounded only in terms of $p_{2}$ and $C_{Z}$, see Theorem~\ref{thm:DecompositionSpaceComplete}.

\medskip{}

All that remains is to establish the final claim of the corollary,
i.e., to estimate $\vertiii{\beta_{1}}$ and $\vertiii{\beta_{2}}$
in case of $Y=\ell_{w}^{q_{1}}\left(I\right)$ and $Z=\ell_{v}^{q_{2}}\left(J\right)$.

For $\beta_{2}$, which in this case reads
\[
\beta_{2}:\ell_{w}^{q_{1}}\Bigl(\bigl[\ell_{u}^{\UpperExpo{p_{1}}}\!\bigl(J_{i}\cap J_{B}\bigr)\bigr]_{i\in I}\Bigr)\hookrightarrow\ell_{v}^{q_{2}}\left(J_{B}\right),
\]
we can simply use Corollary~\ref{cor:EmbeddingCoarseIntoFineSimplification},
with $r=\UpperExpo{p_{1}}$, $J_{0}=J_{B}$ (since $\CalP_{J_{B}}$
is almost subordinate to $\CalQ$) and with
\[
u_{i}^{\left(1\right)}:=\begin{cases}
\left|\det T_{i}\right|^{1-p_{1}^{-1}}, & \text{if }p_{1}<1,\\
1 & \text{if }p_{1}\geq1,
\end{cases}\qquad\text{ and }\qquad u_{j}^{\left(2\right)}:=\begin{cases}
\left|\det S_{j}\right|^{p_{2}^{-1}-1}, & \text{if }p_{1}<1,\\
\left|\det S_{j}\right|^{p_{2}^{-1}-p_{1}^{-1}}, & \text{if }p_{1}\geq1.
\end{cases}
\]
Note that equation~(\ref{eq:DeterminantIsModerate}) shows that $u^{\left(1\right)}$
is $\CalQ$-moderate, with $C_{u^{\left(1\right)},\CalQ}\leq C\left(\dimension,p_{1},\CalQ\right)$.

Hence, Corollary~\ref{cor:EmbeddingCoarseIntoFineSimplification}
yields
\[
\vertiii{\beta_{2}}\asymp\left\Vert \left(w_{i}^{-1}\cdot\left\Vert \left(v_{j}/u_{i,j}\right)_{j\in J_{B}\cap J_{i}}\right\Vert _{\ell^{q_{2}\cdot\left(\UpperExpo{p_{1}}/q_{2}\right)'}}\right)_{i\in I}\right\Vert _{\ell^{q_{2}\cdot\left(q_{1}/q_{2}\right)'}},
\]
where the implied constant only depends on $\dimension,p_{1},q_{1},q_{2},\CalQ,\CalP,k\left(\smash{\CalP_{J_{B}}},\CalQ\right),C_{w,\CalQ}$,
as desired.

Finally, for $\beta_{1}$, it is easy to see that the boundedness
of $\beta_{1}$ is equivalent to that of 
\[
\widetilde{\beta_{1}}:\:\ell_{w}^{q_{1}}\left(L\right)\hookrightarrow\ell_{v}^{q_{2}}\Bigl(\bigl[\ell_{\theta}^{\LowerExpo{p_{2}}}\!\bigl(I_{j}\bigr)\bigr]_{\!j\in J\setminus J_{B}}\Bigr),
\]
with $\vertiii{\beta_{1}}=\vertiii{\smash{\widetilde{\beta_{1}}}}$.
Thus, we use Corollary~\ref{cor:EmbeddingFineInCoarseSimplification},
with $r=\LowerExpo{p_{2}}$, $I_{0}=L$ and $J_{0}=J\setminus J_{B}$,
as well as $u=\theta$ to obtain (because of $I_{j}\subset L$ and
hence $I_{0}\cap I_{j}=L\cap I_{j}=I_{j}$ for $j\in J\setminus J_{B}$)
that 
\[
\vertiii{\widetilde{\beta_{1}}}\asymp\left\Vert \left(v_{j}\cdot\left\Vert \left(w_{i}^{-1}\cdot\theta_{i}\right)_{i\in I_{j}}\right\Vert _{\ell^{\LowerExpo{p_{2}}\cdot\left(q_{1}/\LowerExpo{p_{2}}\right)'}}\right)_{j\in J\setminus J_{B}}\right\Vert _{\ell^{q_{2}\cdot\left(q_{1}/q_{2}\right)'}}\;,
\]
where the implied constant only depends on $p_{2},q_{1},q_{2},\CalP,k\left(\smash{\CalQ_{I_{A}}},\CalP\right),C_{v,\CalP}$.
Note that we indeed have $I_{0}=L=\bigcup_{j\in J_{0}}I_{j}$, as
needed for the application of Corollary~\ref{cor:EmbeddingFineInCoarseSimplification}.
Furthermore, as seen above, we have $I_{0}=L\subset I_{A}$, so that
$\CalQ_{L}=\CalQ_{I_{0}}$ is indeed almost subordinate to $\CalP$,
with $k\left(\CalQ_{L},\CalP\right)\leq k\left(\smash{\CalQ_{I_{A}}},\CalP\right)$.
\end{proof}

\section{Necessary conditions for embeddings}

\label{sec:NecessaryConditions}In this section, we study the sharpness
of the sufficient criteria for the existence of embeddings between
decomposition spaces that we developed before.

Roughly speaking, the setting considered in this section is as follows:
We assume that an embedding of the form
\begin{equation}
\FourierDecompSp{\CalQ}{p_{1}}Y\hookrightarrow\FourierDecompSp{\CalP}{p_{2}}Z\label{eq:DecompositionSpaceEmbeddingSimple}
\end{equation}
is true, where the two decomposition spaces are built with respect
to the ((tight) semi-structured) coverings $\CalQ=\left(T_{i}Q_{i}'+b_{i}\right)_{i\in I}$
of $\CalO$ and $\CalP=\left(\smash{S_{j}P_{j}'+c_{j}}\right)_{j\in J}$
of $\CalO'$. We will then show that the existence of such an embedding
necessarily implies $p_{1}\leq p_{2}$ and also the existence of certain
embeddings between discrete sequence spaces; these embeddings for
the sequence spaces will be similar (and often identical) to the sufficient
conditions which we derived in the previous section.

The section consists of four subsections. In the first subsection,
we impose no special restrictions on the relation between $\CalQ$
and $\CalP$. Nevertheless, we will show that the condition $p_{1}\leq p_{2}$—which
was required in \emph{all} sufficient conditions from the previous
section—is a necessary consequence of the boundedness of the embedding~(\ref{eq:DecompositionSpaceEmbeddingSimple}).
Furthermore, in case of $p_{1}=p_{2}$, we will see that necessarily
$\left\Vert \delta_{j}\right\Vert _{Z}\lesssim\left\Vert \delta_{i}\right\Vert _{Y}$
for all $i\in I$ and $j\in J$ for which $Q_{i}^{\circ}\cap P_{j}^{\circ}\neq\emptyset$.

In Subsection \ref{subsec:CoincidenceOfDecompositionSpaces}, we employ
these elementary necessary criteria to show that an \emph{equality}
\[
\FourierDecompSp{\CalQ}{p_{1}}Y=\FourierDecompSp{\CalP}{p_{2}}Z
\]
of two decomposition spaces is only possible if we have $p_{1}=p_{2}$
and if the coverings $\CalQ,\CalP$ are \emph{weakly equivalent}.
Note though that in general this only holds for $p_{1}\neq2$. For
the case of weighted $\ell^{q}$ spaces $Y=\ell_{w}^{q_{1}}\left(I\right)$
and $Z=\ell_{v}^{q_{2}}\left(J\right)$, however, we will be able
to extend this result: In this case, we necessarily have $q_{1}=q_{2}$
and $\CalQ$ and $\CalP$ are also equivalent for $p_{1}=2$, as long
as $q_{1}\neq2$. Finally, we remark that in addition to the elementary
results from the previous subsection, the proof of the equivalence
between $\CalQ,\CalP$ also uses simplified forms of the arguments
used in Subsection \ref{subsec:ImprovedNecessaryConditions}. Thus,
Subsection \ref{subsec:CoincidenceOfDecompositionSpaces} serves as
a gentle introduction to the remainder of the section.

From Subsection \ref{subsec:ImprovedNecessaryConditions} on, we always
assume (essentially) that $\CalQ$ is almost subordinate to $\CalP$
(or vice versa). Under this assumption, we will show that the boundedness
of the embedding~(\ref{eq:DecompositionSpaceEmbeddingSimple}) (essentially)
implies boundedness of the embedding
\[
Y\hookrightarrow Z\left(\vphantom{\ell_{\left|\det T_{i}\right|^{p_{1}^{-1}-p_{2}^{-1}}}^{p_{2}}}\left[\,\vphantom{\sum}\smash{\ell_{\left|\det T_{i}\right|^{p_{1}^{-1}-p_{2}^{-1}}}^{p_{2}}}\left(I_{j}\right)\,\right]_{j\in J}\right).
\]
This necessary condition almost coincides with the sufficient condition
from Corollary~\ref{cor:EmbeddingFineIntoCoarse}; the only difference
is that the ``inner norm'' for the necessary condition is $\ell^{p_{2}}$,
while it is $\ell^{\LowerExpo{p_{2}}}$ for the sufficient condition.
At least for $p_{2}\in\left(0,2\right]$, we have $\LowerExpo{p_{2}}=p_{2}$,
so that both conditions coincide. Thus, for a certain range of exponents,
we achieve a \emph{complete characterization}. A similar statement
also holds if $\CalP$ is almost subordinate to $\CalQ$.

Finally, in Subsection \ref{subsec:RelativelyModerateCase}, in addition
to $\CalQ$ being almost subordinate to $\CalP$, we require $\CalQ$
to be \emph{relatively moderate} with respect to $\CalP$ and we only
consider weighted $\ell^{q}$ spaces as the ``global components''
for the decomposition spaces. Under these more restrictive assumptions,
we will show that the sufficient conditions from Corollary~\ref{cor:EmbeddingFineIntoCoarse}
are also necessary, so that we obtain a \emph{complete characterization}
for all possible exponents $p_{2}$. With suitable changes, the same
holds for the conditions from Corollary~\ref{cor:EmbeddingCoarseIntoFine},
i.e.\@ if $\CalP$ is almost subordinate to $\CalQ$.

\medskip{}

Before we properly begin our investigation of necessary criteria for
the existence of embeddings, we briefly indicate our proof strategy.
Assuming an embedding $\iota:\FourierDecompSp{\CalQ}{p_{1}}Y\hookrightarrow\FourierDecompSp{\CalP}{p_{2}}Z$,
we will ``test'' this embedding using suitably crafted functions:

For simplicity, we will assume $p_{1}=p_{2}=p\in\left[1,\infty\right]$.
Given any (finitely supported) sequence $\left(c_{j}\right)_{j\in J}$,
where $\CalP=\left(P_{j}\right)_{j\in J}$, we consider functions
of the form
\[
f=\sum_{j\in J}c_{j}\gamma_{j}\qquad\text{ where }\qquad\gamma_{j}\in\TestFunctionSpace{P_{j}}.
\]
If $\Phi=\left(\varphi_{i}\right)_{i\in I}$ is an $L^{p}$-BAPU for
$\CalQ$, then $\varphi_{i}\equiv0$ outside of $Q_{i}$. In particular,
$\varphi_{i}\gamma_{j}\equiv0$ for $j\notin J_{i}=\left\{ \ell\in J\with P_{\ell}\cap Q_{i}\neq\emptyset\right\} $,
and thus
\begin{equation}
\vphantom{\sum_{j\in J_{i}}}\left\Vert \Fourier^{-1}\left(\varphi_{i}\cdot f\right)\right\Vert _{L^{p}}=\left\Vert \Fourier^{-1}\left(\varphi_{i}\cdot\vphantom{\sum}\smash{\sum_{j\in J_{i}}}c_{j}\gamma_{j}\right)\right\Vert _{L^{p}}\lesssim\left\Vert \Fourier^{-1}\left(\vphantom{\sum}\smash{\sum_{j\in J_{i}}}c_{j}\gamma_{j}\right)\right\Vert _{L^{p}}=:d_{i}\qquad\forall\,i\in I.\label{eq:IntroductionArbitrageExplanationDIDefinition}
\end{equation}
Here, the last step used that we have $p\in\left[1,\infty\right]$
and $\left\Vert \Fourier^{-1}\varphi_{i}\right\Vert _{L^{1}}\lesssim1$
for all $i\in I$, as part of the definition of an $L^{p}$-BAPU.
Consequently, using the boundedness of $\iota$, we get
\[
\left\Vert f\right\Vert _{\FourierDecompSp{\CalP}pZ}\lesssim\left\Vert f\right\Vert _{\FourierDecompSp{\CalQ}pY}\lesssim\left\Vert \left(d_{i}\right)_{i\in I}\right\Vert _{Y}<\infty,
\]
at least if we assume $\left\Vert \left(d_{i}\right)_{i\in I}\right\Vert _{Y}<\infty$.
Thus, our next goal is to obtain a \emph{lower} bound on $\left\Vert f\right\Vert _{\FourierDecompSp{\CalP}pZ}$.

To this end, we use that admissibility of the covering $\CalP$ implies
(see the ``\emph{disjointization lemma}'', Lemma~\ref{lem:DisjointizationPrinciple})
that there is a finite partition $J=\biguplus_{r=1}^{r_{0}}J^{\left(r\right)}$
such that no two different sets $P_{j}$ from the same index set $J^{\left(r\right)}$
are neighbors (of degree $3$), i.e., such that $P_{j}^{\ast}\cap P_{\ell}^{\ast}=\emptyset$
for all $j,\ell\in J^{\left(r\right)}$ with $j\neq\ell$. Thus, if
we fix $r\in\underline{r_{0}}$ and assume $c_{j}=0$ for all $j\in J\setminus J^{\left(r\right)}$,
we get 
\[
f=\sum_{j\in J^{\left(r\right)}}c_{j}\gamma_{j}.
\]
The advantage of this additional assumption on the coefficients $\left(c_{j}\right)_{j\in J}$
is—among other things—that the different summands have disjoint support,
which prevents cancellations in the sum. As a consequence of this
and the fact that the family $\left(\gamma_{j}\right)_{j\in J}$ is
adapted to the covering $\CalP$, one can show (see Lemma~\ref{lem:GeneralizedEasyNormEquivalenceFineLowerBound})
that
\[
\left\Vert \left(d_{i}\right)_{i\in I}\right\Vert _{Y}\gtrsim\left\Vert f\right\Vert _{\FourierDecompSp{\CalP}pZ}\gtrsim\left\Vert \left(\left\Vert \Fourier^{-1}\left(c_{j}\gamma_{j}\right)\right\Vert _{L^{p}}\right)_{j\in J}\right\Vert _{Z}=:\left\Vert \left(\left|c_{j}\right|\cdot e_{j}\right)_{j\in J}\right\Vert _{Z}.
\]

Now, we have reached a crucial point: The sequence $\left(e_{j}\right)_{j\in J}$
given by $e_{j}=\left\Vert \Fourier^{-1}\gamma_{j}\right\Vert _{L^{p}}$
is \emph{highly invariant} under certain transformations, for example
under switching from $\gamma_{j}$ to the modulated version $\widetilde{\gamma_{j}}=M_{z_{j}}\gamma_{j}$,
with arbitrary $z_{j}\in\R^{\dimension}$. Note that this transformation
retains the property $\widetilde{\gamma_{j}}\in\TestFunctionSpace{P_{j}}$.
In contrast to $e_{j}$, however, the sequence $\left(d_{i}\right)_{i\in I}$
defined in equation~(\ref{eq:IntroductionArbitrageExplanationDIDefinition})
is (in general) \emph{not} invariant under these transformations.
More precisely, even though it is hard (or even impossible) to evaluate
$d_{i}$ in general, we will see (see Corollary~\ref{cor:AsymptoticModulationBehaviour})
that a suitable choice of the modulations $z_{j}$ entails
\[
d_{i}=\left\Vert \Fourier^{-1}\left(\,\vphantom{\sum_{i}}\smash{\sum_{j\in J_{i}}}c_{j}\cdot M_{z_{j}}\gamma_{j}\right)\right\Vert _{L^{p}}\asymp\left\Vert \left(\left|c_{j}\right|\cdot\left\Vert \Fourier^{-1}\gamma_{j}\right\Vert _{L^{p}}\right)_{j\in J_{i}}\right\Vert _{\ell^{p}}=\left\Vert \left(\left|c_{j}\right|\cdot e_{j}\right)_{j\in J_{i}}\right\Vert _{\ell^{p}}
\]
Thus, in the spirit of \cite[Section 1.9]{TaoStructureRandomness},
we found a possibility to \emph{arbitrage} the embedding, i.e., to
exploit a differing degree of symmetry between both sides of the embedding.

All in all, we arrive at
\[
\left\Vert \left(\left|c_{j}\right|\cdot e_{j}\right)_{j\in J}\right\Vert _{Z}\lesssim\left\Vert \left(d_{i}\right)_{i\in I}\right\Vert _{Y}\asymp\left\Vert \left(\left\Vert \left(\left|c_{j}\right|\cdot e_{j}\right)_{j\in J_{i}}\right\Vert _{\ell^{p}}\right)_{i\in I}\right\Vert _{Y},
\]
for (almost) arbitrary sequences $\left(c_{j}\right)_{j\in J}$ supported
in $J^{\left(r\right)}$. Since $r\in\underline{r_{0}}$ was arbitrary,
this restriction on the support can easily be removed. Careful readers
will notice that the argument sketched here is used (in a slightly
more precise, general and elaborate form) in the proof of Theorem~\ref{thm:BurnerNecessaryConditionCoarseInFine}.

\medskip{}

Now that we have illustrated the general idea, we are ready to properly
start our investigations.

\subsection{Elementary necessary conditions}

\label{subsec:ElementaryNecessaryConditions}In the following, we
will actually consider a slightly more general setting than in equation~(\ref{eq:DecompositionSpaceEmbeddingSimple}).
Precisely, we will assume that there is an embedding
\begin{equation}
\left(\CalD_{K}\:,\left\Vert \mybullet\right\Vert _{\FourierDecompSp{\CalQ}{p_{1}}Y}\right)\hookrightarrow\FourierDecompSp{\CalP}{p_{2}}Z,\label{eq:GeneralEmbeddingRequirement}
\end{equation}
for some subset $K\subset\CalO\cap\CalO'$, where 
\[
\CalD_{K}:=\CalD_{K}^{\CalQ,p_{1},Y}:=\FourierDecompSp{\CalQ}{p_{1}}Y\cap\left\{ f\in\TestFunctionSpace{\R^{\dimension}}\with\supp f\subset K\right\} 
\]
are all those test functions with support in $K$ which also belong
to the Fourier-side decomposition space $\FourierDecompSp{\CalQ}{p_{1}}Y$.
As long as there is no chance for confusion, we will simply write
$\CalD_{K}$ instead of $\CalD_{K}^{\CalQ,p_{1},Y}$. Finally, if
$\delta_{i}\in Y$ for all $i\in I$, it is not hard to see $\TestFunctionSpace{\CalO}\subset\FourierDecompSp{\CalQ}{p_{1}}Y$,
which implies that $\CalD_{K}^{\CalQ,p_{1},Y}=\left\{ f\in\TestFunctionSpace{\R^{\dimension}}\with\supp f\subset K\right\} $
is independent of $\CalQ,p_{1},Y$ (as long as $Y$ contains all finitely
supported sequences). A similar argument shows that $\CalD_{K}^{\CalQ,p_{1},Y}$
is independent of the choice of $p_{1}$ (even if $Y$ does not contain
all finitely supported sequences), but we will not use this in the
following.

Our first necessary condition is rather simple: We observe that in
Theorem~\ref{thm:NoSubordinatenessWithConsiderationOfOverlaps} and
in Corollaries \ref{cor:EmbeddingFineIntoCoarse}, \ref{cor:EmbeddingCoarseIntoFine}
and \ref{cor:MixedSubordinateness}, the inequality $p_{1}\leq p_{2}$
was always part of the assumptions. Our first necessary condition
will show that this is inevitable. It is worth noting that we do
not need to impose any additional assumptions (like subordinateness)
on the coverings $\CalQ,\CalP$.
\begin{lem}
\label{lem:SimpleNecessaryCondition}Let $\emptyset\neq\CalO,\CalO'\subset\R^{\dimension}$
be open, let $p_{1},p_{2}\in\left(0,\infty\right]$ and let $\CalQ=\left(Q_{i}\right)_{i\in I}$
be an $L^{p_{1}}$-decomposition covering of $\CalO$, and let $\CalP=\left(P_{j}\right)_{j\in J}$
be an $L^{p_{2}}$-decomposition covering of $\CalO'$. Let $Y\subset\Compl^{I}$
and $Z\subset\Compl^{J}$ be $\CalQ$-regular and $\CalP$-regular,
respectively.

Finally, let $K\subset\R^{\dimension}$ and assume that there are
$i\in I$ and $j\in J$ with $K^{\circ}\cap Q_{i}^{\circ}\cap P_{j}^{\circ}\neq\emptyset$
and with $\delta_{i}:=\Indicator_{\left\{ i\right\} }\in Y$. If the
identity map
\[
\iota:\left(\CalD_{K}^{\CalQ,p_{1},Y},\left\Vert \mybullet\right\Vert _{\FourierDecompSp{\CalQ}{p_{1}}Y}\right)\rightarrow\FourierDecompSp{\CalP}{p_{2}}Z,f\mapsto f
\]
is well-defined and bounded, then we have $p_{1}\leq p_{2}$.

In case of $p_{1}=p_{2}$, there exists a constant $C>0$ depending
only on 
\[
\dimension,p_{1},p_{2},\CalQ,\CalP,C_{\CalQ,\Phi,p_{1}},C_{\CalP,\Psi,p_{2}},\vertiii{\Gamma_{\CalQ}}_{Y\to Y},\vertiii{\Gamma_{\CalP}}_{Z\to Z}
\]
such that $\delta_{j}\in Z$ and $\left\Vert \delta_{j}\right\Vert _{Z}\leq C\vertiii{\iota}\cdot\left\Vert \delta_{i}\right\Vert _{Y}$
holds for all $\left(i,j\right)\in I\times J$ with $K^{\circ}\cap Q_{i}^{\circ}\cap P_{j}^{\circ}\neq\emptyset$
and with $\delta_{i}\in Y$.

Here, the $L^{p_{1}}$-BAPU $\Phi$ and the $L^{p_{2}}$-BAPU $\Psi$
have to be used to compute the (quasi)-norms on $\FourierDecompSp{\CalQ}{p_{1}}Y$
and $\FourierDecompSp{\CalP}{p_{2}}Z$, respectively.
\end{lem}

\begin{proof}
Let $\Phi=\left(\varphi_{i}\right)_{i\in I}$ be an $L^{p_{1}}$-BAPU
for $\CalQ$ and let $\Psi=\left(\psi_{j}\right)_{j\in J}$ be an
$L^{p_{2}}$-BAPU for $\CalP$. Fix a nontrivial test function $\theta\in\TestFunctionSpace{B_{1}\left(0\right)}$.
We will prove both parts of the lemma (more or less) simultaneously.
The assumptions yield $\xi_{0}\in\R^{\dimension}$ and $\varepsilon>0$
with $B_{\varepsilon}\left(\xi_{0}\right)\subset K^{\circ}\cap Q_{i}^{\circ}\cap P_{j}^{\circ}$.

For $n\in\N$, define 
\[
\theta_{n}:\R^{\dimension}\rightarrow\Compl,\xi\mapsto\theta\left(\frac{n}{\varepsilon}\left(\xi-\xi_{0}\right)\right),
\]
i.e.\@ $\theta_{n}=L_{\xi_{0}}\left[D_{n/\varepsilon}\theta\right]$.
We then have 
\[
\theta_{n}\in\TestFunctionSpace{B_{\varepsilon/n}\left(\xi_{0}\right)}\subset\TestFunctionSpace{B_{\varepsilon}\left(\xi_{0}\right)}\subset\TestFunctionSpace{K^{\circ}\cap P_{j}^{\circ}\cap Q_{i}^{\circ}}.
\]
But Lemma~\ref{lem:PartitionCoveringNecessary} implies $\psi_{j}^{\ast}\equiv1$
on $P_{j}\supset\supp\theta_{n}$, and hence $\psi_{j}^{\ast}\theta_{n}=\theta_{n}$.

Now note that if $\varphi_{\ell}\cdot\theta_{n}\not\equiv0$, then
there is some $\xi\in\supp\theta_{n}\subset Q_{i}$ with $\varphi_{\ell}\left(\xi\right)\neq0$.
Hence, $\xi\in Q_{\ell}\cap Q_{i}\neq\emptyset$, i.e.\@ $\ell\in i^{\ast}$.
But for $\ell\in i^{\ast}$, there are two cases: For $p_{1}\in\left[1,\infty\right]$,
Young's inequality ($L^{1}\ast L^{p_{1}}\hookrightarrow L^{p_{1}}$)
yields
\[
\left\Vert \Fourier^{-1}\!\left(\varphi_{\ell}\,\theta_{n}\right)\right\Vert _{L^{p_{1}}}=\left\Vert \left(\Fourier^{-1}\varphi_{\ell}\right)\ast\left(\Fourier^{-1}\theta_{n}\right)\right\Vert _{L^{p_{1}}}\leq\left\Vert \Fourier^{-1}\varphi_{\ell}\right\Vert _{L^{1}}\cdot\left\Vert \Fourier^{-1}\theta_{n}\right\Vert _{L^{p_{1}}}\leq C_{\CalQ,\Phi,p_{1}}\cdot\left\Vert \Fourier^{-1}\theta_{n}\right\Vert _{L^{p_{1}}}\!<\!\infty.
\]
Otherwise, in case of $p_{1}\in\left(0,1\right)$, we know (by definition
of an $L^{p_{1}}$-decomposition covering) that $\CalQ=\left(Q_{i}\right)_{i\in I}=\left(T_{i}Q_{i}'+b_{i}\right)_{i\in I}$
is a semi-structured covering. Furthermore, we note $\supp\varphi_{\ell}\subset\overline{Q_{\ell}}\subset\overline{Q_{\ell}^{\ast}}$
and $\supp\theta_{n}\subset Q_{i}^{\circ}\subset\overline{Q_{\ell}^{\ast}}$,
so that Corollary~\ref{cor:QuasiBanachConvolutionSemiStructured}
yields a constant $C_{1}=C_{1}\left(\CalQ,\dimension,p_{1}\right)>0$
with
\[
\left\Vert \Fourier^{-1}\left(\varphi_{\ell}\,\theta_{n}\right)\right\Vert _{L^{p_{1}}}\leq C_{1}\cdot\left|\det T_{\ell}\right|^{\frac{1}{p_{1}}-1}\cdot\left\Vert \Fourier^{-1}\varphi_{\ell}\right\Vert _{L^{p_{1}}}\left\Vert \Fourier^{-1}\theta_{n}\right\Vert _{L^{p_{1}}}\leq C_{1}C_{\CalQ,\Phi,p_{1}}\cdot\left\Vert \Fourier^{-1}\theta_{n}\right\Vert _{L^{p_{1}}}\qquad\forall\,\ell\in i^{\ast}\,.
\]

If we set $C_{1}:=1$ for $p_{1}\in\left[1,\infty\right]$, then a
combination of the two cases from above easily yields
\[
\left\Vert \Fourier^{-1}\left(\varphi_{\ell}\,\theta_{n}\right)\right\Vert _{L^{p_{1}}}\leq C_{1}C_{\CalQ,\Phi,p_{1}}\cdot\left\Vert \Fourier^{-1}\theta_{n}\right\Vert _{L^{p_{1}}}\cdot\Indicator_{i^{\ast}}\left(\ell\right)\qquad\forall\,\ell\in I.
\]
But since $Y$ is a $\CalQ$-regular—in particular solid—sequence
space with $\delta_{i}\in Y$ and since $\Indicator_{i^{\ast}}=\Gamma_{\CalQ}\,\delta_{i}$,
we get $\theta_{n}\in\FourierDecompSp{\CalQ}{p_{1}}Y$ (and hence
$\theta_{n}\in\CalD_{K}^{\CalQ,p_{1},Y}$) with
\begin{align*}
\left\Vert \theta_{n}\right\Vert _{\FourierDecompSp{\CalQ}{p_{1}}Y}=\left\Vert \left(\left\Vert \Fourier^{-1}\left(\varphi_{\ell}\,\theta_{n}\right)\right\Vert _{L^{p_{1}}}\right)_{\ell\in I}\right\Vert _{Y} & \leq C_{1}C_{\CalQ,\Phi,p_{1}}\cdot\left\Vert \Fourier^{-1}\theta_{n}\right\Vert _{L^{p_{1}}}\cdot\left\Vert \Indicator_{i^{\ast}}\right\Vert _{Y}\\
 & \leq C_{1}C_{\CalQ,\Phi,p_{1}}\vertiii{\Gamma_{\CalQ}}_{Y\to Y}\cdot\left\Vert \delta_{i}\right\Vert _{Y}\cdot\left\Vert \Fourier^{-1}\theta_{n}\right\Vert _{L^{p_{1}}}\\
 & =:C_{2}\cdot\left\Vert \delta_{i}\right\Vert _{Y}\cdot\left\Vert \Fourier^{-1}\theta_{n}\right\Vert _{L^{p_{1}}}<\infty\qquad\forall\,n\in\N\,.
\end{align*}

Since we assume $\iota$ to be well-defined and bounded, we get $\theta_{n}\in\FourierDecompSp{\CalP}{p_{2}}Z$
and
\[
\left\Vert \theta_{n}\right\Vert _{\FourierDecompSp{\CalP}{p_{2}}Z}\leq\vertiii{\iota}\cdot\left\Vert \theta_{n}\right\Vert _{\FourierDecompSp{\CalQ}{p_{1}}Y}\leq C_{2}\cdot\vertiii{\iota}\cdot\left\Vert \delta_{i}\right\Vert _{Y}\cdot\left\Vert \Fourier^{-1}\theta_{n}\right\Vert _{L^{p_{1}}}.
\]
But by Remark~\ref{rem:ClusteredBAPUYIeldsBoundedControlSystem},
the family $\Gamma:=\left(\gamma_{\ell}\right)_{\ell\in J}:=\left(\psi_{\ell}^{\ast}\right)_{\ell\in J}$
is an $L^{p_{2}}$-bounded family for $\CalP$, with $C_{\CalP,\Gamma,p_{2}}\leq C\left(\CalP,C_{\CalP,\Psi,p_{2}},\dimension,p_{2}\right)$.
Thus, Theorem~\ref{thm:BoundedControlSystemEquivalentQuasiNorm}
yields
\begin{align*}
\left\Vert \theta_{n}\right\Vert _{\FourierDecompSp{\CalP}{p_{2}}Z} & \geq C_{3}^{-1}\cdot\left\Vert \left(\left\Vert \Fourier^{-1}\left(\psi_{\ell}^{\ast}\theta_{n}\right)\right\Vert _{L^{p_{2}}}\right)_{\ell\in J}\right\Vert _{Z}\\
\left({\scriptstyle \text{since }Z\text{ is solid and }\psi_{j}^{\ast}\theta_{n}=\theta_{n}}\right) & \geq C_{3}^{-1}\cdot\left\Vert \Fourier^{-1}\theta_{n}\right\Vert _{L^{p_{2}}}\cdot\left\Vert \delta_{j}\right\Vert _{Z}\,,
\end{align*}
with a constant $C_{3}=C_{3}\left(\CalP,p_{2},\dimension,C_{\CalP,\Gamma,p_{2}},\vertiii{\Gamma_{\CalP}}_{Z\to Z}\right)>0$.
In particular, $\delta_{j}\in Z$.

Altogether, we finally arrive at
\begin{equation}
\left\Vert \Fourier^{-1}\theta_{n}\right\Vert _{L^{p_{2}}}\cdot\left\Vert \delta_{j}\right\Vert _{Z}\leq C_{2}C_{3}\cdot\vertiii{\iota}\cdot\left\Vert \delta_{i}\right\Vert _{Y}\cdot\left\Vert \Fourier^{-1}\theta_{n}\right\Vert _{L^{p_{1}}}\qquad\forall\,n\in\N\,.\label{eq:SimpleNecessaryConditionMainEstimate}
\end{equation}
Furthermore, we can estimate 
\[
C_{2}C_{3}\leq C_{4}=C_{4}\left(\dimension,p_{1},p_{2},\CalQ,\CalP,C_{\CalQ,\Phi,p_{1}},C_{\CalP,\Psi,p_{2}},\vertiii{\Gamma_{\CalQ}}_{Y\to Y},\vertiii{\Gamma_{\CalP}}_{Z\to Z}\right).
\]

Now, note $\Fourier^{-1}\theta_{n}=\left(\varepsilon/n\right)^{\dimension}\cdot e^{2\pi i\left\langle \mybullet,\xi_{0}\right\rangle }\cdot\left(\Fourier^{-1}\theta\right)\left(\frac{\varepsilon}{n}\mybullet\right)$,
and hence
\[
\left\Vert \mathcal{F}^{-1}\theta_{n}\right\Vert _{L^{p}}=\left(\varepsilon/n\right)^{\dimension}\cdot\left\Vert D_{\frac{\varepsilon}{n}}\left(\Fourier^{-1}\theta\right)\right\Vert _{L^{p}}=\left(\varepsilon/n\right)^{\dimension\left(1-\frac{1}{p}\right)}\cdot\left\Vert \Fourier^{-1}\theta\right\Vert _{L^{p}}
\]
for all $p\in\left(0,\infty\right]$. Thus, estimate~(\ref{eq:SimpleNecessaryConditionMainEstimate})
yields
\[
\left\Vert \Fourier^{-1}\theta\right\Vert _{L^{p_{2}}}\cdot\left(\varepsilon/n\right)^{\dimension\left(1-\frac{1}{p_{2}}\right)}\cdot\left\Vert \delta_{j}\right\Vert _{Z}\leq C_{4}\cdot\vertiii{\iota}\cdot\left\Vert \delta_{i}\right\Vert _{Y}\cdot\left(\varepsilon/n\right)^{\dimension\left(1-\frac{1}{p_{1}}\right)}\cdot\left\Vert \Fourier^{-1}\theta\right\Vert _{L^{p_{1}}}
\]
and hence
\[
n^{\frac{1}{p_{2}}-\frac{1}{p_{1}}}\leq C_{4}\vertiii{\iota}\varepsilon^{\dimension\left(\frac{1}{p_{2}}-\frac{1}{p_{1}}\right)}\cdot\frac{\left\Vert \delta_{i}\right\Vert _{Y}}{\left\Vert \delta_{j}\right\Vert _{Z}}\cdot\frac{\left\Vert \Fourier^{-1}\theta\right\Vert _{L^{p_{1}}}}{\left\Vert \Fourier^{-1}\theta\right\Vert _{L^{p_{2}}}}\qquad\forall\,n\in\N,
\]
which can only hold if $p_{2}^{-1}-p_{1}^{-1}\leq0$, i.e.\@ if $p_{1}\leq p_{2}$.
This yields the first claim.

Finally, in case of $p_{1}=p_{2}$, most terms cancel, so that we
get $\left\Vert \delta_{j}\right\Vert _{Z}\leq C_{4}\vertiii{\iota}\cdot\left\Vert \delta_{i}\right\Vert _{Y}$
with $C_{4}$ independent of $i,j$.
\end{proof}

\subsection{Coincidence of decomposition spaces}

\label{subsec:CoincidenceOfDecompositionSpaces}In this subsection,
we assume that we are given two coverings $\CalQ,\CalP$ of a common
set $\emptyset\neq\CalO\subset\R^{\dimension}$ such that we have
\begin{equation}
\FourierDecompSp{\CalQ}{p_{1}}{\ell_{w}^{q_{1}}}=\FourierDecompSp{\CalP}{p_{2}}{\ell_{v}^{q_{2}}}\label{eq:EqualityOfDecompositionSpaces}
\end{equation}
for certain $p_{1},p_{2},q_{1},q_{2}\in\left(0,\infty\right]$ and
certain weights $w,v$ which are $\CalQ$-moderate and $\CalP$-moderate,
respectively. Given these assumptions, we will be able to show that
the ``ingredients'' for the two decomposition spaces are essentially
the same. Precisely, we will show $p_{1}=p_{2}$, $q_{1}=q_{2}$ and
$w_{i}\asymp v_{j}$ given $Q_{i}\cap P_{j}\neq\emptyset$. Finally,
in case of $\left(p_{1},q_{1}\right)\neq\left(2,2\right)$, we will
also show that the two coverings $\CalQ,\CalP$ are \textbf{weakly
equivalent}. Recall from Definition~\ref{def:RelativeIndexClustersSubordinateCoveringsModerateCoverings}
that this means that the sets
\[
I_{j}:=\left\{ i\in I\with Q_{i}\cap P_{j}\neq\emptyset\right\} \qquad\text{ and }\qquad J_{i}:=\left\{ j\in J\with P_{j}\cap Q_{i}\neq\emptyset\right\} 
\]
satisfy $\left|I_{j}\right|\leq C$ and $\left|J_{i}\right|\leq C$
for all $i\in I$ and $j\in J$. As the final result of this section,
we show that these properties together are sufficient for the equality
(\ref{eq:EqualityOfDecompositionSpaces}) to hold, at least for $p_{1}\in\left[1,\infty\right]$.
In the quasi-Banach regime $p_{1}\in\left(0,1\right)$, we will have
to impose certain additional technical assumptions.

We thus obtain an (almost) complete characterization of the equality
of two decomposition spaces in terms of their ``ingredients''. Note
though that in order to get the full characterization as explained
above, we need to restrict ourselves to global components of the form
$Y=\ell_{w}^{q_{1}}$ instead of general ($\CalQ$)-regular sequence
spaces. For general ($\CalQ$)-regular sequence spaces, we will nevertheless
be able to show that $\FourierDecompSp{\CalQ}{p_{1}}Y=\FourierDecompSp{\CalP}{p_{2}}Z$
can only hold if $p_{1}=p_{2}$; in case of $p_{1}\neq2$, a further
necessary condition is that $\CalQ$ and $\CalP$ are weakly equivalent.

Finally, the methods of proof in this subsection provide a nice preparation
for the following subsections, which are slightly more involved. One
general technique which we will use again and again is to ``test''
an embedding between (or an equality of) two Fourier-side decomposition
spaces by considering functions of the form
\begin{equation}
f_{z,c,\varepsilon}=\sum_{i\in I}M_{z_{i}}\left(\varepsilon_{i}c_{i}\cdot\gamma_{i}\right)\label{eq:UsualTestFunctionFine}
\end{equation}
for certain functions $\gamma_{i}$ which are supported (essentially)
in $Q_{i}$ and for (more or less) arbitrary sequences $\varepsilon=\left(\varepsilon_{i}\right)_{i\in I}\in\left\{ \pm1\right\} ^{I}$
and $c=\left(c_{i}\right)_{i\in I}\in\ell_{0}\left(I\right)$. With
this choice, the (quasi)-norm $\left\Vert f_{z,c,\varepsilon}\right\Vert _{\FourierDecompSp{\CalQ}{p_{1}}Y}$
is (essentially) independent of the choice of the modulations $z_{i}\in\R^{\dimension}$
and of the signs $\varepsilon_{i}$, since modulation does not change
the support (on the Fourier side) of the individual summands and since
the individual summands are (essentially) disjointly supported. Note
that we are considering the quasi-norm with respect to $\CalQ$, i.e.
with respect to the covering to which the $\gamma_{i}$ are adapted.

In general, this independence from the choice of $z=\left(z_{i}\right)_{i\in I}$
and $\varepsilon=\left(\varepsilon_{i}\right)_{i\in I}$ does \emph{not}
hold for the (quasi)-norm $\left\Vert f_{z,c,\varepsilon}\right\Vert _{\FourierDecompSp{\CalP}{p_{2}}Z}$.
As we will see (see Corollary~\ref{cor:AsymptoticModulationBehaviour}),
a suitable choice of $z$ often makes it possible to derive a certain
embedding between sequence spaces as a result of the embedding for
decomposition spaces. A similar statement holds for a suitable choice
of $\varepsilon$. In the spirit of \cite[Section 1.9]{TaoStructureRandomness},
we have thus found a possibility to \emph{arbitrage} the embedding,
see also the introduction to this section.

Our next lemma provides one of two estimates for the (essential) independence
of $\left\Vert f_{z,c,\varepsilon}\right\Vert _{\FourierDecompSp{\CalQ}{p_{1}}Y}$
from the choice of $z$.
\begin{lem}
\label{lem:EasyNormEquivalenceFineCovering}Let $\emptyset\neq\CalO\subset\R^{\dimension}$
be open, let $p\in\left(0,\infty\right]$ and let $\CalQ=\left(Q_{i}\right)_{i\in I}$
be an $L^{p}$-decomposition covering of $\CalO$ with $L^{p}$-BAPU
$\Phi=\left(\varphi_{i}\right)_{i\in I}$. Finally, let $Y\subset\Compl^{I}$
be $\CalQ$-regular.

For each $k\in\N_{0}$, there is a constant
\[
C=C\left(k,\dimension,p,\CalQ,C_{\CalQ,\Phi,p},\vertiii{\Gamma_{\CalQ}}_{Y\to Y}\right)>0
\]
such that for arbitrary coefficients $\left(c_{i}\right)_{i\in I}\in\ell_{0}\left(I\right)$,
functions $\left(\gamma_{i}\right)_{i\in I}$ with $\gamma_{i}\in\TestFunctionSpace{\CalO}$
and $\gamma_{i}\equiv0$ on $\left(Q_{i}^{k\ast}\right)^{c}$, signs
$\varepsilon=\left(\varepsilon_{i}\right)_{i\in I}\in\left\{ \pm1\right\} ^{I}$
and modulations $\left(z_{i}\right)_{i\in I}\in\left(\R^{\dimension}\right)^{I}$,
the estimate
\begin{equation}
\left\Vert \sum_{i\in I}M_{z_{i}}\left(\varepsilon_{i}c_{i}\cdot\gamma_{i}\right)\right\Vert _{\BAPUFourierDecompSp{\CalQ}pY{\Phi}}\leq C\cdot\left\Vert \left(c_{i}\cdot\left\Vert \Fourier^{-1}\gamma_{i}\right\Vert _{L^{p}}\right)_{i\in I}\right\Vert _{Y}\label{eq:StandardEasyNormEstimate}
\end{equation}
holds, provided that $\left(c_{i}\cdot\left\Vert \Fourier^{-1}\gamma_{i}\right\Vert _{L^{p}}\right)_{i\in I}\in Y$.
In particular, if we have $\left(c_{i}\cdot\left\Vert \Fourier^{-1}\gamma_{i}\right\Vert _{L^{p}}\right)_{i\in I}\in Y$,
then $\sum_{i\in I}M_{z_{i}}\left(\varepsilon_{i}c_{i}\cdot\gamma_{i}\right)\in\FourierDecompSp{\CalQ}pY$.
\end{lem}

\begin{rem*}
Since the sequence $\left(c_{i}\right)_{i\in I}\in\ell_{0}\left(I\right)$
is finitely supported and since $\gamma_{i}\in\TestFunctionSpace{\CalO}\subset\Schwartz\left(\R^{\dimension}\right)$,
we have $\left(c_{i}\cdot\left\Vert \Fourier^{-1}\gamma_{i}\right\Vert _{L^{p}}\right)_{i\in I}\in Y$
as soon as $\delta_{i}\in Y$ for all $i\in I$, which holds for essentially
all reasonable ``global components'' $Y$, in particular for $Y=\ell_{w}^{q}\left(I\right)$.
\end{rem*}
\begin{proof}
Let $\ell\in I$ be arbitrary. For $i\in I$, we note that
\[
0\not\equiv\varphi_{\ell}\cdot M_{z_{i}}\left(\varepsilon_{i}c_{i}\cdot\gamma_{i}\right)=\varepsilon_{i}c_{i}\cdot M_{z_{i}}\left(\varphi_{\ell}\cdot\gamma_{i}\right)
\]
is only possible if $\emptyset\neq Q_{\ell}\cap Q_{i}^{k\ast}$ and
thus $i\in\ell^{\left(k+1\right)\ast}$.

As usual, we distinguish the two cases $p\in\left[1,\infty\right]$
and $p\in\left(0,1\right)$. For $p\in\left[1,\infty\right]$, Young's
inequality ($L^{1}\ast L^{p}\hookrightarrow L^{p}$), together with
the identity $\Fourier^{-1}\left[M_{z}g\right]=L_{-z}\left[\Fourier^{-1}g\right]$
and with the translation invariance of $\left\Vert \mybullet\right\Vert _{L^{p}}$,
implies
\[
\left\Vert \Fourier^{-1}\left[M_{z_{i}}\left(\varphi_{\ell}\cdot\gamma_{i}\right)\right]\right\Vert _{L^{p}}=\left\Vert \Fourier^{-1}\left[\varphi_{\ell}\cdot\gamma_{i}\right]\right\Vert _{L^{p}}\leq\left\Vert \Fourier^{-1}\varphi_{\ell}\right\Vert _{L^{1}}\cdot\left\Vert \Fourier^{-1}\gamma_{i}\right\Vert _{L^{p}}\leq C_{\CalQ,\Phi,p}\cdot\left\Vert \Fourier^{-1}\gamma_{i}\right\Vert _{L^{p}}.
\]
In case of $p\in\left(0,1\right)$, the definition of an $L^{p}$-decomposition
covering shows that $\CalQ=\left(T_{i}Q_{i}'+b_{i}\right)_{i\in I}$
is semi-structured. Next, we observe $\supp\varphi_{\ell}\subset\overline{Q_{\ell}}\subset\overline{Q_{\ell}^{\left(2k+1\right)\ast}}$
and $\supp\gamma_{i}\subset\overline{Q_{i}^{k\ast}}\subset\overline{Q_{\ell}^{\left(2k+1\right)\ast}}$,
because of $i\in\ell^{\left(k+1\right)\ast}$. Hence, Corollary~\ref{cor:QuasiBanachConvolutionSemiStructured}
yields a constant $C_{1}=C_{1}\left(\CalQ,k,\dimension,p\right)$
with
\begin{align*}
\left\Vert \Fourier^{-1}\left[M_{z_{i}}\left(\varphi_{\ell}\cdot\gamma_{i}\right)\right]\right\Vert _{L^{p}}=\left\Vert \Fourier^{-1}\left[\varphi_{\ell}\cdot\gamma_{i}\right]\right\Vert _{L^{p}} & \leq C_{1}\cdot\left|\det T_{\ell}\right|^{\frac{1}{p}-1}\left\Vert \Fourier^{-1}\varphi_{\ell}\right\Vert _{L^{p}}\cdot\left\Vert \Fourier^{-1}\gamma_{i}\right\Vert _{L^{p}}\\
 & \leq C_{1}C_{\CalQ,\Phi,p}\cdot\left\Vert \smash{\Fourier^{-1}}\gamma_{i}\right\Vert _{L^{p}}.
\end{align*}
Thus, if we set $C_{1}:=1$ for $p\in\left[1,\infty\right]$, this
estimate is valid for all $p\in\left(0,\infty\right]$.

Since $\left\Vert \mybullet\right\Vert _{L^{p}}$ is a quasi-norm
(with triangle constant only depending on $p$) and because of the
uniform bound $\left|\smash{\ell^{\left(k+1\right)\ast}}\right|\leq N_{\CalQ}^{k+1}$,
which was established in Lemma~\ref{lem:SemiStructuredClusterInvariant},
there is a constant $C_{2}=C_{2}\left(p,k,\CalQ\right)$ satisfying
\begin{align*}
\left\Vert \Fourier^{-1}\!\biggl[\varphi_{\ell}\cdot\sum_{i\in I}M_{z_{i}}\left(\varepsilon_{i}c_{i}\cdot\gamma_{i}\right)\biggr]\!\right\Vert _{L^{p}}\! & \overset{\left(\dagger\right)}{\leq}C_{2}\sum_{i\in\ell^{\left(k+1\right)\ast}}\left[\left|c_{i}\right|\cdot\left\Vert \Fourier^{-1}\left[M_{z_{i}}\left(\varphi_{\ell}\cdot\gamma_{i}\right)\right]\right\Vert _{L^{p}}\right]\\
 & \leq C_{1}C_{2}C_{\CalQ,\Phi,p}\cdot\sum_{i\in\ell^{\left(k+1\right)\ast}}\left[\left|c_{i}\right|\cdot\left\Vert \Fourier^{-1}\gamma_{i}\right\Vert _{L^{p}}\right]\\
 & =C_{1}C_{2}C_{\CalQ,\Phi,p}\cdot\left(\Theta_{k+1}\,\zeta\right)_{\ell}\qquad\forall\,\ell\in I.
\end{align*}
In the last step, we introduced the sequence $\zeta=\left(\zeta_{i}\right)_{i\in I}$
defined by $\zeta_{i}:=\left|c_{i}\right|\cdot\left\Vert \Fourier^{-1}\gamma_{i}\right\Vert _{L^{p}}$
for $i\in I$. Furthermore, we used the \textbf{$\left(k+1\right)$-fold
clustering map} $\Theta_{k+1}:Y\to Y$ as defined in Lemma~\ref{lem:HigherOrderClusteringMap}.
At $\left(\dagger\right)$, we used that $\varphi_{\ell}\cdot M_{z_{i}}\left(\varepsilon_{i}c_{i}\cdot\gamma_{i}\right)\not\equiv0$
can only hold for $i\in\ell^{\left(k+1\right)\ast}$, as observed
at the beginning of the proof.

Let $f:=\sum_{i\in I}M_{z_{i}}\left(\varepsilon_{i}c_{i}\cdot\gamma_{i}\right)$.
Note that our assumptions and the solidity of $Y$ imply $\zeta\in Y$
and thus $\Theta_{k+1}\,\zeta\in Y$. But the estimate above showed
$\left\Vert \Fourier^{-1}\left(\varphi_{\ell}\cdot f\right)\right\Vert _{L^{p}}\lesssim\left(\Theta_{k+1}\,\zeta\right)_{\ell}$
for all $\ell\in I$. Hence, the solidity of $Y$ yields $\left(\left\Vert \Fourier^{-1}\left(\varphi_{\ell}\cdot f\right)\right\Vert _{L^{p}}\right)_{\ell\in I}\in Y$
and thus $f\in\FourierDecompSp{\CalQ}pY$. To obtain a quantitative
estimate, note that Lemma~\ref{lem:HigherOrderClusteringMap} yields
$\vertiii{\Theta_{k+1}}\leq C_{3}$ for some constant $C_{3}=C_{3}\left(k,\vertiii{\Gamma_{\CalQ}}_{Y\to Y}\right)$.
Finally, define $C_{4}:=C_{1}C_{2}C_{3}C_{\CalQ,\Phi,p}$ and conclude
(using the solidity of $Y$) that
\begin{align*}
\vphantom{\sum_{i\in I}}\left\Vert \,\smash{\sum_{i\in I}}\,\vphantom{\sum}M_{z_{i}}\left(\varepsilon_{i}c_{i}\cdot\gamma_{i}\right)\,\right\Vert _{\BAPUFourierDecompSp{\CalQ}pY{\Phi}} & \leq C_{1}C_{2}C_{\CalQ,\Phi,p}\cdot\left\Vert \Theta_{k+1}\,\zeta\right\Vert _{Y}\\
 & \leq C_{4}\cdot\left\Vert \left(\left|c_{i}\right|\cdot\left\Vert \smash{\Fourier^{-1}}\gamma_{i}\right\Vert _{L^{p}}\right)_{i\in I}\right\Vert _{Y}\\
 & =C_{4}\cdot\left\Vert \left(c_{i}\cdot\left\Vert \smash{\Fourier^{-1}}\gamma_{i}\right\Vert _{L^{p}}\right)_{i\in I}\right\Vert _{Y}.\qedhere
\end{align*}
\end{proof}
The reverse of the previous estimate is not true in general, since
there could be cancellations between neighboring indices, for example
we could have $\gamma_{i}+\gamma_{\ell}=0$ for certain $i\neq\ell$.
We will see in the next lemma, however, that this can be excluded
by introducing suitable assumptions:
\begin{lem}
\label{lem:GeneralizedEasyNormEquivalenceFineLowerBound}Let $\emptyset\neq\CalO\subset\R^{\dimension}$
be open, let $p\in\left(0,\infty\right]$, and let $\CalQ=\left(Q_{i}\right)_{i\in I}$
be an $L^{p}$-decomposition covering of $\CalO$. Finally, let $Y\subset\Compl^{I}$
be $\CalQ$-regular.

For each $k\in\N_{0}$, there is a constant
\[
C=C\left(\CalQ,p,\dimension,k,\vertiii{\Gamma_{\CalQ}}_{Y\to Y}\right)>0
\]
such that the following holds: If

\begin{itemize}[leftmargin=0.7cm]
\item $I_{0}\subset I$ with $Q_{i}^{k\ast}\cap Q_{\ell}^{k\ast}=\emptyset$
for all $i,\ell\in I_{0}$ with $i\neq\ell$,
\item $\left(c_{i}\right)_{i\in I_{0}}\in\ell_{0}\left(I_{0}\right)$,
\item $\left(\phi_{i}\right)_{i\in I}$ is an $L^{p}$-bounded system for
$\CalQ$ with $\phi_{i}\equiv0$ on $\CalO\setminus Q_{i}^{k\ast}$
for all $i\in I$,
\item $\left(\gamma_{i}\right)_{i\in I_{0}}$ satisfies $\gamma_{i}\in\TestFunctionSpace{\CalO}$
and $\phi_{i}\gamma_{i}=\gamma_{i}$ for all $i\in I_{0}$,
\item $\varepsilon=\left(\varepsilon_{i}\right)_{i\in I_{0}}\in\left\{ \pm1\right\} ^{I_{0}}$
and $\left(z_{i}\right)_{i\in I_{0}}\in\left(\R^{\dimension}\right)^{I_{0}}$,
\item $f\in\FourierDecompSp{\CalQ}pY$ with $\phi_{i}f=\phi_{i}g$ for all
$i\in I_{0}$, where $g:=\sum_{i\in I_{0}}M_{z_{i}}\left(\varepsilon_{i}c_{i}\cdot\gamma_{i}\right)$,
\end{itemize}
then
\begin{equation}
\left\Vert \left(c_{i}\cdot\left\Vert \Fourier^{-1}\gamma_{i}\right\Vert _{L^{p}}\right)_{i\in I_{0}}\right\Vert _{Y|_{I_{0}}}\leq C\cdot C_{\CalQ,\left(\phi_{i}\right)_{i},p}\cdot\left\Vert f\right\Vert _{\BAPUFourierDecompSp{\CalQ}pY{\Phi}}\label{eq:GeneralizedEasyNormEstimateLowerBound}
\end{equation}
for each $L^{p}$-BAPU $\Phi$ for $\CalQ$. In particular, $\left(c_{i}\cdot\left\Vert \Fourier^{-1}\gamma_{i}\right\Vert _{L^{p}}\right)_{i\in I_{0}}\in Y|_{I_{0}}$.
\end{lem}

\begin{proof}
Our assumption on the $L^{p}$-bounded system $\left(\phi_{i}\right)_{i\in I}$
shows that we can choose $\ell_{\left(\phi_{i}\right)_{i},\CalQ}=k$.
Therefore, Theorem~\ref{thm:BoundedControlSystemEquivalentQuasiNorm}
(with the partition $\left(I^{\left(r\right)}\right)_{r=1}$ of $I$
consisting of the single element $I^{\left(1\right)}:=I$ and with
the $L^{p}$-bounded system $\left(\phi_{i}\right)_{i\in I}$) yields
a constant $C_{1}=C_{1}\left(\CalQ,p,\dimension,k,\vertiii{\Gamma_{\CalQ}}_{Y\to Y}\right)>0$
with
\begin{align*}
C_{1}C_{\CalQ,\left(\phi_{i}\right)_{i},p}\cdot\left\Vert f\right\Vert _{\CalD_{\Fourier,\Phi}\left(\CalQ,L^{p},Y\right)} & \geq\left\Vert f\right\Vert _{\left(\phi_{i}\right)_{i},\left(I^{\left(r\right)}\right)_{r},L^{p},Y}\\
 & =\left\Vert \left(\left\Vert \Fourier^{-1}\left(\phi_{i}\cdot f\right)\right\Vert _{L^{p}}\right)_{i\in I}\right\Vert _{Y}=\left\Vert \zeta^{\left(f\right)}\right\Vert _{Y}\qquad\forall\,f\in\FourierDecompSp{\CalQ}pY\,,
\end{align*}
where we defined $\zeta^{\left(f\right)}$ by $\zeta_{i}^{\left(f\right)}:=\left\Vert \Fourier^{-1}\left(\phi_{i}\cdot f\right)\right\Vert _{L^{p}}$.
In particular, we have thus shown $\zeta^{\left(f\right)}\in Y$ for
all $f\in\FourierDecompSp{\CalQ}pY$.

Now, we will apply this to the given function $f$ from the assumption.
To this end, let $\ell\in I_{0}$ be arbitrary. For $i\in I_{0}$,
there are two cases:

\begin{casenv}
\item We have $i\neq\ell$. By our assumption on $I_{0}$, this implies
$Q_{i}^{k\ast}\cap Q_{\ell}^{k\ast}=\emptyset$. But because of $\phi_{i}\equiv0$
on $\CalO\setminus Q_{i}^{k\ast}$, this yields $\phi_{i}\,\phi_{\ell}\equiv0$.
Using $\gamma_{i}=\phi_{i}\,\gamma_{i}$, we finally get $\phi_{\ell}\,\gamma_{i}=\phi_{\ell}\,\phi_{i}\,\gamma_{i}\equiv0$.
\item We have $i=\ell$. In this case, our assumptions yield $\phi_{\ell}\,\gamma_{i}=\phi_{i}\,\gamma_{i}=\gamma_{i}$.
\end{casenv}
All in all, these considerations yield—together with our assumption
$\phi_{\ell}\,f=\phi_{\ell}\,g$ for all $\ell\in I_{0}$—that
\[
\phi_{\ell}\,f=\phi_{\ell}\,g=\phi_{\ell}\cdot\sum_{i\in I_{0}}M_{z_{i}}\left(\varepsilon_{i}c_{i}\cdot\gamma_{i}\right)=\sum_{i\in I_{0}}M_{z_{i}}\left(\varepsilon_{i}c_{i}\cdot\phi_{\ell}\gamma_{i}\right)=M_{z_{\ell}}\left(\varepsilon_{\ell}\,c_{\ell}\cdot\gamma_{\ell}\right)\qquad\forall\,\ell\in I_{0}\,.
\]
But this means 
\begin{align*}
\zeta_{\ell}^{\left(f\right)} & =\left\Vert \Fourier^{-1}\left(\phi_{\ell}\cdot f\right)\right\Vert _{L^{p}}=\left\Vert \Fourier^{-1}\left[M_{z_{\ell}}\left(\varepsilon_{\ell}c_{\ell}\cdot\gamma_{\ell}\right)\right]\right\Vert _{L^{p}}\\
 & =\left|c_{\ell}\right|\cdot\left\Vert L_{-z_{\ell}}\left[\Fourier^{-1}\gamma_{\ell}\right]\right\Vert _{L^{p}}=\left|c_{\ell}\right|\cdot\left\Vert \Fourier^{-1}\gamma_{\ell}\right\Vert _{L^{p}}=:\omega_{\ell}\qquad\forall\,\ell\in I_{0}\,.
\end{align*}
If we set $\omega_{\ell}:=0$ for $\ell\in I\setminus I_{0}$, we
have thus shown $0\leq\omega_{\ell}\leq\zeta_{\ell}^{\left(f\right)}$
for all $\ell\in I$. Using the solidity of $Y$, this implies $\omega:=\left(\omega_{\ell}\right)_{\ell\in I}\in Y$,
and also the desired estimate
\[
\left\Vert \left(c_{i}\cdot\left\Vert \Fourier^{-1}\gamma_{i}\right\Vert _{L^{p}}\right)_{i\in I_{0}}\right\Vert _{Y|_{I_{0}}}=\left\Vert \omega\right\Vert _{Y}\leq\left\Vert \zeta^{\left(f\right)}\right\Vert _{Y}\leq C_{1}C_{\CalQ,\left(\phi_{i}\right)_{i},p}\cdot\left\Vert f\right\Vert _{\BAPUFourierDecompSp{\CalQ}pY{\Phi}}\:.\qedhere
\]
\end{proof}
Occasionally, the following specialized version of the preceding lemma
is more convenient to use.
\begin{cor}
\label{cor:EasyNormEquivalenceFineLowerBound}Let $\emptyset\neq\CalO\subset\R^{\dimension}$
be open, let $p\in\left(0,\infty\right]$, and let $\CalQ=\left(Q_{i}\right)_{i\in I}$
be an $L^{p}$-decomposition covering of $\CalO$ with $L^{p}$-BAPU
$\Phi=\left(\varphi_{i}\right)_{i\in I}$. Finally, let $Y\subset\Compl^{I}$
be $\CalQ$-regular.

For each $k\in\N_{0}$, there is a constant
\[
C=C\left(\CalQ,p,\dimension,k,C_{\CalQ,\Phi,p},\vertiii{\Gamma_{\CalQ}}_{Y\to Y}\right)>0
\]
such that the following holds: If

\begin{itemize}[leftmargin=0.7cm]
\item $I_{0}\subset I$ with $Q_{i}^{\left(k+1\right)\ast}\cap Q_{\ell}^{\left(k+1\right)\ast}=\emptyset$
for all $i,\ell\in I_{0}$ with $i\neq\ell$,
\item $\left(c_{i}\right)_{i\in I_{0}}\in\ell_{0}\left(I_{0}\right)$,
\item $\left(\gamma_{i}\right)_{i\in I_{0}}$ satisfies $\gamma_{i}\in\TestFunctionSpace{\CalO}$
and $\gamma_{i}\equiv0$ on $\CalO\setminus Q_{i}^{k\ast}$ for all
$i\in I_{0}$,
\item $\varepsilon=\left(\varepsilon_{i}\right)_{i\in I_{0}}\in\left\{ \pm1\right\} ^{I_{0}}$
and $\left(z_{i}\right)_{i\in I_{0}}\in\left(\R^{\dimension}\right)^{I_{0}}$,
\end{itemize}
then
\begin{equation}
\vphantom{\sum_{i\in I_{0}}}\left\Vert \,\smash{\sum_{i\in I_{0}}}\,\vphantom{\sum}M_{z_{i}}\left(\varepsilon_{i}c_{i}\cdot\gamma_{i}\right)\,\right\Vert _{\BAPUFourierDecompSp{\CalQ}pY{\Phi}}\geq C^{-1}\cdot\left\Vert \left(c_{i}\cdot\left\Vert \Fourier^{-1}\gamma_{i}\right\Vert _{L^{p}}\right)_{i\in I_{0}}\right\Vert _{Y|_{I_{0}}},\label{eq:StandardEasyNormEstimateLowerBound}
\end{equation}
provided that $\sum_{i\in I_{0}}M_{z_{i}}\left(\varepsilon_{i}c_{i}\cdot\gamma_{i}\right)\in\FourierDecompSp{\CalQ}pY$.
Therefore, if $\sum_{i\in I_{0}}M_{z_{i}}\left(\varepsilon_{i}c_{i}\cdot\gamma_{i}\right)\in\FourierDecompSp{\CalQ}pY$,
then $\left(c_{i}\cdot\left\Vert \Fourier^{-1}\gamma_{i}\right\Vert _{L^{p}}\right)_{i\in I_{0}}\in Y|_{I_{0}}$.
\end{cor}

\begin{rem*}
For brevity, let $f:=\sum_{i\in I_{0}}M_{z_{i}}\left(\varepsilon_{i}c_{i}\cdot\gamma_{i}\right)$.
The corollary shows for $\zeta:=\left(c_{i}\cdot\left\Vert \Fourier^{-1}\gamma_{i}\right\Vert _{L^{p}}\right)_{i\in I_{0}}$
that $\zeta\in Y|_{I_{0}}$ is a necessary condition for $f\in\FourierDecompSp{\CalQ}pY$.
Furthermore, it yields a corresponding (quasi)-norm estimate.

Conversely, an application of Lemma~\ref{lem:EasyNormEquivalenceFineCovering}
(to the extended sequences $\left(c_{i}\right)_{i\in I}$ and $\left(\gamma_{i}\right)_{i\in I}$
with $c_{i}=0$ and $\gamma_{i}\equiv0$ for $i\in I\setminus I_{0}$)
shows that we have $f\in\FourierDecompSp{\CalQ}pY$ as soon as $\left(c_{i}\cdot\left\Vert \Fourier^{-1}\gamma_{i}\right\Vert _{L^{p}}\right)_{i\in I}\in Y$,
i.e. as soon as $\zeta\in Y|_{I_{0}}$. Furthermore, in this case,
the same lemma shows $\left\Vert f\right\Vert _{\BAPUFourierDecompSp{\CalQ}pY{\Phi}}\lesssim\left\Vert \zeta\right\Vert _{Y|_{I_{0}}}$,
where the constant depends on the same quantities as in Corollary~\ref{cor:EasyNormEquivalenceFineLowerBound}.

In summary, a combination of Corollary~\ref{cor:EasyNormEquivalenceFineLowerBound}
and Lemma~\ref{lem:EasyNormEquivalenceFineCovering} shows $f\in\FourierDecompSp{\CalQ}pY\Longleftrightarrow\zeta\in Y|_{I_{0}}$
and also
\begin{equation}
\left\Vert f\right\Vert _{\BAPUFourierDecompSp{\CalQ}pY{\Phi}}\asymp\left\Vert \left(c_{i}\cdot\left\Vert \Fourier^{-1}\gamma_{i}\right\Vert _{L^{p}}\right)_{i\in I_{0}}\right\Vert _{Y|_{I_{0}}}\;,\label{eq:EasyNormEquivalence}
\end{equation}
where the implied constant only depends on $\CalQ,p,\dimension,k,C_{\CalQ,\Phi,p},\vertiii{\Gamma_{\CalQ}}_{Y\to Y}$.
\end{rem*}
\begin{proof}
We apply Lemma~\ref{lem:GeneralizedEasyNormEquivalenceFineLowerBound}
with $k+1$ instead of $k$ and with $\phi_{i}:=\varphi_{i}^{\left(k+1\right)\ast}$.
To verify the prerequisites of Lemma~\ref{lem:GeneralizedEasyNormEquivalenceFineLowerBound},
note that we indeed have $\phi_{i}\equiv0$ on $\CalO\setminus Q_{i}^{\left(k+1\right)\ast}$
and that Remark~\ref{rem:ClusteredBAPUYIeldsBoundedControlSystem}
shows that the family $\Lambda:=\left(\phi_{i}\right)_{i\in I}:=\left(\smash{\varphi_{i}^{\left(k+1\right)\ast}}\right)_{i\in I}$
is indeed an $L^{p}$-bounded control system for $\CalQ$ with 
\[
C_{\CalQ,\Lambda,p}\leq C_{1}=C_{1}\left(\CalQ,C_{\CalQ,\Phi,p},\dimension,p,k\right).
\]

Furthermore, since $\gamma_{i}\equiv0$ on $\CalO\setminus Q_{i}^{k\ast}$
and since $\phi_{i}=\varphi_{i}^{\left(k+1\right)\ast}\equiv1$ on
$Q_{i}^{k\ast}$ (see Lemma~\ref{lem:PartitionCoveringNecessary}),
we also get $\phi_{i}\,\gamma_{i}=\gamma_{i}$ for all $i\in I_{0}$,
as required in Lemma~\ref{lem:GeneralizedEasyNormEquivalenceFineLowerBound}.
The claim now easily follows from that lemma.
\end{proof}
Above, it was claimed that—using a suitable choice of $z=\left(z_{i}\right)_{i\in I}$—one
can achieve suitable values of $\left\Vert f_{z,c,\varepsilon}\right\Vert _{\FourierDecompSp{\CalP}{p_{2}}Z}$,
with $f_{z,c,\varepsilon}$ as in equation~(\ref{eq:UsualTestFunctionFine}).
The goal of the next few lemmata is to show how this can be done.
The main point is that
\[
\vphantom{\sum_{i=1}^{n}}\left\Vert \Fourier^{-1}\left[\,\smash{\sum_{i=1}^{n}}\,\vphantom{\sum}M_{z_{i}}f_{i}\,\right]\right\Vert _{L^{p}}=\left\Vert \,\smash{\sum_{i=1}^{n}}\,\vphantom{\sum}L_{-z_{i}}\left[\Fourier^{-1}f_{i}\right]\,\right\Vert _{L^{p}}\xrightarrow[\min\limits _{i\neq j}\left|z_{i}-z_{j}\right|\to\infty]{}\left\Vert \left(\left\Vert \Fourier^{-1}f_{i}\right\Vert _{L^{p}}\right)_{i\in\underline{n}}\right\Vert _{\ell^{p}}\;.
\]
Validity of the stated limit follows (at least if $\Fourier^{-1}f_{i}\in L^{p}\cap C_{0}$
for all $i\in\underline{n}$) from the following lemma:
\begin{lem}
\label{lem:AsymptoticTranslationNormEstimate}Let $n\in\N$, $p\in\left(0,\infty\right]$
and $f_{1},\dots,f_{n}\in L^{p}\left(\R^{\dimension}\right)$. For
$p=\infty$, assume additionally that
\[
f_{i}\in\overline{\left\{ f\in L^{\infty}\left(\R^{\dimension}\right)\with\supp f\text{ compact}\right\} }\qquad\text{ for all }i\in\underline{n},
\]
where the closure is taken in $L^{\infty}\left(\R^{\dimension}\right)$.

Then we have
\[
\vphantom{\sum_{i=1}^{n}}\left\Vert \,\smash{\sum_{i=1}^{n}}\,\vphantom{\sum}L_{x_{i}}f_{i}\,\right\Vert _{L^{p}}\xrightarrow[\min\limits _{i\neq j}\left|x_{i}-x_{j}\right|\to\infty]{}\left\Vert \left(\left\Vert f_{i}\right\Vert _{L^{p}}\right)_{i\in\underline{n}}\right\Vert _{\ell^{p}}\,.
\]
In particular, there exists $R=R\left(p,\smash{\left(f_{i}\right)_{i\in\underline{n}}}\right)>0$
with
\[
\frac{1}{2}\cdot\left\Vert \left(\left\Vert f_{i}\right\Vert _{L^{p}}\right)_{i\in\underline{n}}\right\Vert _{\ell^{p}}\leq\vphantom{\sum_{i=1}^{n}}\left\Vert \,\smash{\sum_{i=1}^{n}}\,\vphantom{\sum}L_{x_{i}}f_{i}\,\right\Vert _{L^{p}}\leq2\cdot\left\Vert \left(\left\Vert f_{i}\right\Vert _{L^{p}}\right)_{i\in\underline{n}}\right\Vert _{\ell^{p}}
\]
for all $x_{1},\dots,x_{n}\in\R^{\dimension}$ with $\left|x_{i}-x_{j}\right|\geq R$
for all $i,j\in\underline{n}$ with $i\neq j$.
\end{lem}

\begin{rem*}
The additional assumption regarding the $f_{i}$ is certainly satisfied
for $f_{1},\dots,f_{n}\in C_{0}\left(\R^{\dimension}\right)$, where
$C_{0}\left(\R^{\dimension}\right)$ is the space of continuous functions
vanishing at infinity, i.e.\@ with $\lim_{\left|x\right|\to\infty}f\left(x\right)=0$.
\end{rem*}
\begin{proof}
A complete proof of this statement is given in \cite[Lemma 4.1]{DecompositionIntoSobolev}.
Here, we only sketch the idea.

If $f_{1},\dots,f_{n}$ are compactly supported, say $\supp f_{i}\subset B_{R}\left(0\right)$
for all $i\in\underline{n}$ and if $\left|x_{i}-x_{j}\right|>2R$
for all $i\neq j$, then the supports of $\left(L_{x_{i}}f_{i}\right)_{i\in\underline{n}}$
are pairwise disjoint. Hence, for $p<\infty$, we have
\[
\vphantom{\sum_{i=1}^{n}}\left|\,\smash{\sum_{i=1}^{n}}\,\vphantom{\sum}L_{x_{i}}f_{i}\left(x\right)\,\right|^{p}=\sum_{i=1}^{n}\left|L_{x_{i}}f_{i}\left(x\right)\right|^{p},
\]
since for each $x\in\R^{\dimension}$, at most one summand is nonzero.
This easily implies the claim (with equality for $\min_{i\neq j}\left|x_{i}-x_{j}\right|>2R$)
for compactly supported $f_{1},\dots,f_{n}$ and $p<\infty$. The
general case follows by approximating $f_{1},\dots,f_{n}$ by compactly
supported $g_{1},\dots,g_{n}\in L^{p}\left(\R^{\dimension}\right)$.
Finally, the case $p=\infty$ can be handled using similar arguments.
\end{proof}
The following version of the preceding lemma is more adapted to the
setting that we will consider.
\begin{cor}
\label{cor:AsymptoticModulationBehaviour}Let $p\in\left(0,\infty\right]$
and let $M\neq\emptyset$ be a \emph{finite} index set such that for
each $i\in M$ a Schwartz-function $f_{i}\in\Schwartz\left(\R^{\dimension}\right)$
is given.

For $S\subset M$ and any family $z=\left(z_{i}\right)_{i\in M}$
in $\R^{\dimension}$, define
\[
f_{S}^{\left(z\right)}:=\sum_{i\in S}M_{z_{i}}f_{i}\in\Schwartz\left(\R^{\dimension}\right).
\]

Then there is a constant $R=R\left(\left(f_{i}\right)_{i\in M},M,p\right)>0$
such that for every family $\left(z_{i}\right)_{i\in M}$ in $\R^{\dimension}$
which satisfies $\left|z_{i}-z_{j}\right|\geq R$ for all $i,j\in M$
with $i\neq j$, the estimate
\begin{equation}
\frac{1}{2}\cdot\left\Vert \left(\left\Vert \smash{\Fourier^{-1}}f_{i}\right\Vert _{L^{p}}\right)_{i\in S}\vphantom{\sum}\right\Vert _{\ell^{p}}\leq\left\Vert \Fourier^{-1}f_{S}^{\left(z\right)}\right\Vert _{L^{p}}\leq2\cdot\left\Vert \left(\left\Vert \smash{\Fourier^{-1}}f_{i}\right\Vert _{L^{p}}\right)_{i\in S}\vphantom{\sum}\right\Vert _{\ell^{p}}\label{eq:AsymptoticModulationBehaviour}
\end{equation}
is true for every $S\subset M$. Furthermore, such a family $\left(z_{i}\right)_{i\in M}$
always exists.
\end{cor}

\begin{proof}
We first show the claim for the case $S=M$. In case of $f_{i}\equiv0$
for all $i\in M$, choose $R=1$ and note that equation~(\ref{eq:AsymptoticModulationBehaviour})
holds trivially. Otherwise, there is some $i\in M$ with $f_{i}\not\equiv0$.
In particular, $M\neq\emptyset$, so that we can assume $M=\left\{ i_{1},\dots,i_{n}\right\} $
with $n:=\left|M\right|$ and (necessarily) pairwise distinct $i_{1},\dots,i_{n}\in M$.

Because of $\Fourier^{-1}f_{i}\in\Schwartz\left(\R^{\dimension}\right)\subset L^{p}\left(\R^{\dimension}\right)$,
the following expression (then a constant) is finite and positive:
\[
\varepsilon:=\frac{1}{2}\cdot\left\Vert \left(\left\Vert \smash{\Fourier^{-1}}f_{i}\right\Vert _{L^{p}}\right)_{i\in M}\right\Vert _{\ell^{p}}>0\,.
\]

For $\ell\in\underline{n}$, let $g_{\ell}:=\Fourier^{-1}f_{i_{\ell}}\in\Schwartz\left(\R^{\dimension}\right)\subset L^{p}\left(\R^{\dimension}\right)$
and note that we even have $g_{\ell}\in\Schwartz\left(\R^{\dimension}\right)\subset C_{0}\left(\R^{\dimension}\right)$,
so that Lemma~\ref{lem:AsymptoticTranslationNormEstimate} yields
some $R>0$ such that
\begin{align*}
\vphantom{\sum_{\ell=1}^{n}}\left|\:\left\Vert \Fourier^{-1}\left[\,\smash{\sum_{\ell=1}^{n}}\,\vphantom{\sum}M_{-y_{\ell}}f_{i_{\ell}}\,\right]\right\Vert _{L^{p}}-2\varepsilon\:\right| & =\vphantom{\sum_{\ell=1}^{n}}\left|\:\left\Vert \,\smash{\sum_{\ell=1}^{n}}\,\vphantom{\sum}L_{y_{\ell}}\left(\smash{\Fourier^{-1}}f_{i_{\ell}}\right)\,\right\Vert _{L^{p}}-\left\Vert \left(\left\Vert \smash{\Fourier^{-1}}f_{i}\right\Vert _{L^{p}}\right)_{i\in M}\vphantom{\left(\left\Vert g_{\ell}\right\Vert _{L^{p}}\right)_{\ell\in\underline{n}}}\right\Vert _{\ell^{p}}\:\right|\\
 & =\vphantom{\sum_{\ell=1}^{n}}\left|\:\left\Vert \,\smash{\sum_{\ell=1}^{n}}\,\vphantom{\sum}L_{y_{\ell}}g_{\ell}\,\right\Vert _{L^{p}}-\left\Vert \left(\left\Vert g_{\ell}\right\Vert _{L^{p}}\right)_{\ell\in\underline{n}}\right\Vert _{\ell^{p}}\:\right|<\varepsilon
\end{align*}
holds for every family $\left(y_{\ell}\right)_{\ell\in\underline{n}}$
which satisfies $\left|y_{\ell}-y_{k}\right|\geq R$ for all $\ell,k\in\underline{n}$
with $\ell\neq k$. Here, we used the identity $\Fourier^{-1}\left[M_{-z}f\right]=L_{z}\left[\smash{\Fourier^{-1}}f\right]$.

For $y_{\ell}:=-z_{i_{\ell}}$, this implies
\[
\left\Vert \Fourier^{-1}f_{M}^{\left(z\right)}\right\Vert _{L^{p}}=\vphantom{\sum_{i}}\left\Vert \Fourier^{-1}\left[\,\vphantom{\sum}\smash{\sum_{\ell=1}^{n}}\,M_{-y_{\ell}}f_{i_{\ell}}\,\right]\right\Vert _{L^{p}}\in\left(2\varepsilon-\varepsilon,2\varepsilon+\varepsilon\right)\,,
\]
as soon as we have $\left|y_{k}-y_{\ell}\right|=\left|z_{i_{\ell}}-z_{i_{k}}\right|\geq R$
for all $\ell\neq k$. Using the definition of $\varepsilon$, this
easily implies the claim for $S=M$.

Now the above (applied to $\left(f_{i}\right)_{i\in S}$ for $S\subset M$)
yields some $R_{S}>0$ for each $S\subset M$ so that equation~(\ref{eq:AsymptoticModulationBehaviour})
is satisfied for every family $\left(z_{i}\right)_{i\in S}$ satisfying
$\left|z_{i}-z_{j}\right|\geq R_{S}$ for all $i,j\in S$ with $i\neq j$.

Since $M$ is finite, the same is true of the power set $\mathcal{P}\left(M\right)$,
so that $R:=\max_{S\subset M}R_{S}$ is finite. It is now easy to
see that the claim holds for this choice of $R$.

For the existence of a family $\left(z_{i}\right)_{i\in M}$ as in
the statement of the corollary, note that we can take $z_{i_{\ell}}:=\ell\cdot R\cdot\left(1,0,\dots,0\right)\in\R^{\dimension}$,
since this implies $\left|z_{i_{\ell}}-z_{i_{k}}\right|=\left|\ell-k\right|\cdot R\geq R$
for $\ell,k\in\underline{n}$ with $\ell\neq k$.
\end{proof}
A convenient tool for choosing suitable signs $\varepsilon=\left(\varepsilon_{i}\right)_{i\in I}$
for the function $f_{z,c,\varepsilon}$ (see equation~(\ref{eq:UsualTestFunctionFine}))
is given by \textbf{Khintchine's inequality}, which we state next.
A proof can be found for example in \cite[Proposition 4.5]{WolffHarmonicAnalysis}
or in \cite[Section 10.3, Theorem 1]{ChowProbabilityTheory}.
\begin{thm}
\label{thm:KhintchineInequality}Let $p\in\left(0,\infty\right)$.
Then there is a constant $C_{p}>0$ such that
\[
C_{p}^{-1}\cdot\left(\,\smash{\sum_{n=1}^{N}}\,\vphantom{\sum}\left|a_{n}\right|^{2}\,\right)^{p/2}\vphantom{\sum_{n=1}^{N}}\leq\mathbb{E}_{\omega}\,\left|\,\smash{\sum_{n=1}^{N}}\,\vphantom{\sum}\omega_{n}a_{n}\,\right|^{p}\leq C_{p}\cdot\left(\,\smash{\sum_{n=1}^{N}}\,\vphantom{\sum}\left|a_{n}\right|^{2}\,\right)^{p/2}
\]
holds for all $N\in\N$ and all $a_{1},\dots,a_{N}\in\Compl$.

Here, the expectation $\mathbb{E}_{\omega}$ is taken with respect
to the random variable $\omega\in\left\{ \pm1\right\} ^{N}$ which
is assumed to be uniformly distributed in that set, i.e.\@ $\mathbb{E}_{\omega}\left[f\left(\omega\right)\right]=\frac{1}{2^{N}}\sum_{\omega\in\left\{ \pm1\right\} ^{N}}f\left(\omega\right)$.
\end{thm}

Using all of these ingredients, we can now prove the following technical
lemma which we will use to prove our desired necessary criteria for
the coincidence of two decomposition spaces.
\begin{lem}
\label{lem:CoincidencePreparationLemma}Let $\emptyset\neq\CalO\subset\R^{\dimension}$
be open, let $p\in\left(0,\infty\right]$ and let $\CalQ=\left(Q_{i}\right)_{i\in I}$
and $\CalP=\left(P_{j}\right)_{j\in J}$ be two \emph{open} $L^{p}$-decomposition
coverings of $\CalO$, with associated $L^{p}$-BAPUs $\Phi=\left(\varphi_{i}\right)_{i\in I}$
and $\Psi=\left(\psi_{j}\right)_{j\in J}$. Finally, let $Y\subset\Compl^{I}$
and $Z\subset\Compl^{J}$ be $\CalQ$-regular and $\CalP$-regular,
respectively, with $\ell_{0}\left(I\right)\subset Y$ and with $\ell_{0}\left(J\right)\subset Z$.

Assume that
\begin{equation}
\left\Vert f\right\Vert _{\FourierDecompSp{\CalQ}pY}\asymp\left\Vert f\right\Vert _{\FourierDecompSp{\CalP}pZ}\qquad\forall\,f\in\TestFunctionSpace{\CalO}.\label{eq:CoincidencePreparationLemmaAssumption}
\end{equation}
Then, for arbitrary $i_{0}\in I$ and pairwise distinct $j_{1},\dots,j_{N}\in J_{i_{0}}=\left\{ j\in J\with P_{j}\cap Q_{i_{0}}\neq\emptyset\right\} $,
we have
\begin{equation}
\left\Vert \Indicator_{\left\{ j_{1},\dots,j_{N}\right\} }\right\Vert _{Z}\asymp\left\Vert \delta_{i_{0}}\right\Vert _{Y}\cdot N^{1/p},\label{eq:CoincidencePreparationLemmaP}
\end{equation}
where the implied constant only depends on the implied constants in
equation~(\ref{eq:CoincidencePreparationLemmaAssumption}) and on
\[
\dimension,p,\CalQ,\CalP,\vertiii{\Gamma_{\CalQ}}_{Y\to Y},\vertiii{\Gamma_{\CalP}}_{Z\to Z},C_{\CalQ,\Phi,p},C_{\CalP,\Psi,p},C_{Z},
\]
where $C_{Z}\geq1$ is a triangle constant for $Z$. In particular,
the implied constant is independent of $i_{0},j_{1},\dots,j_{N}$
and of $N$.

Furthermore, in case of $p\in\left(0,\infty\right)$, we also have
\begin{equation}
\left\Vert \Indicator_{\left\{ j_{1},\dots,j_{N}\right\} }\right\Vert _{Z}\asymp\left\Vert \delta_{i_{0}}\right\Vert _{Y}\cdot N^{1/2},\label{eq:CoincidencePreparationLemma2}
\end{equation}
where the implied constant depends on the same quantities as above.

Finally, for $p=\infty$, we have
\begin{equation}
\left\Vert \Indicator_{\left\{ j_{1},\dots,j_{N}\right\} }\right\Vert _{Z}\asymp\left\Vert \delta_{i_{0}}\right\Vert _{Y}\cdot N,\label{eq:CoincidencePreparationLemma1}
\end{equation}
where the implied constant depends on the same quantities as above.
\end{lem}

\begin{proof}
Since we assume $\ell_{0}\left(I\right)\leq Y$ and $\ell_{0}\left(J\right)\leq Z$,
it is not hard to see that we have 
\[
\TestFunctionSpace{\CalO}\subset\FourierDecompSp{\CalQ}pY\cap\FourierDecompSp{\CalP}pZ,
\]
which we will use without comment in the remainder of the proof.

Fix a nontrivial, \emph{nonnegative} function $\gamma\in\TestFunctionSpace{B_{1}\left(0\right)}\setminus\left\{ 0\right\} $
for the remainder of the proof. Furthermore, set $r_{0}:=N_{\CalP}^{3}=N_{\CalP}^{2\cdot1+1}$,
so that the disjointization lemma (Lemma~\ref{lem:DisjointizationPrinciple})
yields a finite partition $J=\biguplus_{r=1}^{r_{0}}J^{\left(r\right)}$
satisfying $P_{j}^{\ast}\cap P_{\ell}^{\ast}=\emptyset$ for all $j,\ell\in J^{\left(r\right)}$
with $j\neq\ell$ and arbitrary $r\in\underline{r_{0}}$. 

Since we assume $\CalQ,\CalP$ to be open coverings and because of
$j_{1},\dots,j_{N}\in J_{i_{0}}$, there is some $\varepsilon>0$
and $\xi_{1},\dots,\xi_{N}\in\R^{\dimension}$ satisfying $B_{\varepsilon}\left(\xi_{\ell}\right)\subset Q_{i_{0}}\cap P_{j_{\ell}}$
for all $\ell\in\underline{N}$. Note that $\varepsilon>0$ and $\xi_{1},\dots,\xi_{N}$
depend heavily on $j_{1},\dots,j_{N}$ and on $i_{0}$, but that all
occurrences of $\varepsilon$ in our estimates will cancel in the
end.

Define $\xi_{j}$ for $j\in J_{0}:=\left\{ j_{1},\dots,j_{N}\right\} $
by $\xi_{j_{\ell}}:=\xi_{\ell}$ for $\ell\in\underline{N}$ and set
$\gamma_{j}^{\left(\varepsilon\right)}:=L_{\xi_{j}}\left[\gamma\left(\varepsilon^{-1}\mybullet\right)\right]$
for $j\in J_{0}$. Note that $\gamma_{j}^{\left(\varepsilon\right)}\in\TestFunctionSpace{B_{\varepsilon}\left(\xi_{j}\right)}\subset\TestFunctionSpace{P_{j}\cap Q_{i_{0}}}$
for all $j\in J_{0}$. Now, for $r\in\underline{r_{0}}$, set
\[
f_{z,\theta}^{\left(r\right)}:=\!\!\sum_{j\in J_{0}\cap J^{\left(r\right)}}\!\!\!\theta_{j}\cdot M_{z_{j}}\gamma_{j}^{\left(\varepsilon\right)}\!=\!\!\sum_{j\in J_{0}\cap J^{\left(r\right)}}\!\!\!M_{z_{j}}\!\left[\theta_{j}\gamma_{j}^{\left(\varepsilon\right)}\right]\quad\text{for }z=\left(z_{j}\right)_{j\in J_{0}}\in\left(\R^{\dimension}\right)^{J_{0}}\!\text{ and }\theta=\left(\theta_{j}\right)_{j\in J_{0}}\!\in\left\{ \pm1\right\} ^{J_{0}}\!\!.
\]

Now, Lemma~\ref{lem:EasyNormEquivalenceFineCovering} (with $\CalP$
instead of $\CalQ$, with $k=0$ and with $c_{j}=\Indicator_{J_{0}\cap J^{\left(r\right)}}$)
yields a constant $C_{1}=C_{1}\left(\dimension,p,\CalP,C_{\CalP,\Psi,p},\vertiii{\Gamma_{\CalP}}_{Z\to Z}\right)>0$
satisfying
\begin{align*}
\left\Vert f_{z,\theta}^{\left(r\right)}\right\Vert _{\FourierDecompSp{\CalP}pZ} & \leq C_{1}\cdot\left\Vert \left(\Indicator_{J_{0}\cap J^{\left(r\right)}}\left(j\right)\cdot\left\Vert \Fourier^{-1}\gamma_{j}^{\left(\varepsilon\right)}\right\Vert _{L^{p}}\right)_{j\in J}\right\Vert _{Z}\\
 & =C_{1}\cdot\left\Vert \Fourier^{-1}\left[\gamma\left(\varepsilon^{-1}\mybullet\right)\right]\right\Vert _{L^{p}}\cdot\left\Vert \Indicator_{J_{0}\cap J^{\left(r\right)}}\right\Vert _{Z}\\
\left({\scriptstyle Z\text{ is solid}}\right) & \leq C_{1}\cdot\varepsilon^{\dimension\left(1-\frac{1}{p}\right)}\left\Vert \Fourier^{-1}\gamma\right\Vert _{L^{p}}\cdot\left\Vert \Indicator_{J_{0}}\right\Vert _{Z},
\end{align*}
for arbitrary $r\in\underline{r_{0}}$. Here, we also used
\[
\left\Vert \Fourier^{-1}\gamma_{j}^{\left(\varepsilon\right)}\right\Vert _{L^{p}}=\left\Vert \Fourier^{-1}\left(L_{\xi_{j}}\left[\gamma\left(\varepsilon^{-1}\mybullet\right)\right]\right)\right\Vert _{L^{p}}=\left\Vert M_{\xi_{j}}\left(\Fourier^{-1}\left[\gamma\left(\varepsilon^{-1}\mybullet\right)\right]\right)\right\Vert _{L^{p}}=\left\Vert \Fourier^{-1}\left[\gamma\left(\varepsilon^{-1}\mybullet\right)\right]\right\Vert _{L^{p}},
\]
by standard properties of the Fourier transform. Similar identities
will be used in the remainder of the proof without comment.

Conversely, Corollary~\ref{cor:EasyNormEquivalenceFineLowerBound}
(again with $\CalP$ instead of $\CalQ$, with $k=0$ and with $I_{0}=J_{0}\cap J^{\left(r\right)}$,
as well as $c_{j}=1$ for all $j\in J_{0}\cap J^{\left(r\right)}=I_{0}$)
yields a constant $C_{2}=C_{2}\left(\dimension,p,\CalP,C_{\CalP,\Psi,p},\vertiii{\Gamma_{\CalP}}_{Z\to Z}\right)>0$
satisfying
\begin{align*}
\left\Vert f_{z,\theta}^{\left(r\right)}\right\Vert _{\FourierDecompSp{\CalP}pZ}\geq C_{2}^{-1}\cdot\left\Vert \left(\left\Vert \Fourier^{-1}\gamma_{j}^{\left(\varepsilon\right)}\right\Vert _{L^{p}}\right)_{j\in I_{0}}\right\Vert _{Z|_{I_{0}}} & =C_{2}^{-1}\cdot\left\Vert \Fourier^{-1}\left[\gamma\left(\varepsilon^{-1}\mybullet\right)\right]\right\Vert _{L^{p}}\cdot\left\Vert \Indicator_{J_{0}\cap J^{\left(r\right)}}\right\Vert _{Z}\\
 & =C_{2}^{-1}\cdot\varepsilon^{\dimension\left(1-\frac{1}{p}\right)}\left\Vert \Fourier^{-1}\gamma\right\Vert _{L^{p}}\cdot\left\Vert \Indicator_{J_{0}\cap J^{\left(r\right)}}\right\Vert _{Z}
\end{align*}
for arbitrary $r\in\underline{r_{0}}$. Note that Corollary~\ref{cor:EasyNormEquivalenceFineLowerBound}
requires $P_{j}^{\ast}\cap P_{\ell}^{\ast}=P_{j}^{\left(k+1\right)\ast}\cap P_{\ell}^{\left(k+1\right)\ast}\overset{!}{=}\emptyset$
for arbitrary $j,\ell\in I_{0}=J_{0}\cap J^{\left(r\right)}$ with
$j\neq\ell$. This holds for every $r\in\underline{r_{0}}$ by choice
of the partition $J=\biguplus_{r=1}^{r_{0}}J^{\left(r\right)}$.

Now, using the identity $J_{0}=\biguplus_{r=1}^{r_{0}}\left(J_{0}\cap J^{\left(r\right)}\right)$
and the (quasi)-triangle inequality for $Z$, we arrive at $\left\Vert \Indicator_{J_{0}}\right\Vert _{Z}\leq C_{3}\cdot\sum_{r=1}^{r_{0}}\left\Vert \Indicator_{J_{0}\cap J^{\left(r\right)}}\right\Vert _{Z}$
for some constant $C_{3}=C_{3}\left(r_{0},C_{Z}\right)=C_{3}\left(\CalP,C_{Z}\right)$.
Hence,
\begin{align}
\varepsilon^{\dimension\left(1-\frac{1}{p}\right)}\left\Vert \Fourier^{-1}\gamma\right\Vert _{L^{p}}\cdot\left\Vert \Indicator_{J_{0}}\right\Vert _{Z} & \leq C_{3}\cdot\sum_{r=1}^{r_{0}}\varepsilon^{\dimension\left(1-\frac{1}{p}\right)}\left\Vert \Fourier^{-1}\gamma\right\Vert _{L^{p}}\cdot\left\Vert \Indicator_{J_{0}\cap J^{\left(r\right)}}\right\Vert _{Z}\nonumber \\
 & \leq C_{2}C_{3}\cdot\sum_{r=1}^{r_{0}}\left\Vert f_{z^{\left(r\right)},\theta^{\left(r\right)}}^{\left(r\right)}\right\Vert _{\FourierDecompSp{\CalP}pZ}\nonumber \\
 & \leq C_{1}C_{2}C_{3}r_{0}\cdot\varepsilon^{\dimension\left(1-\frac{1}{p}\right)}\left\Vert \Fourier^{-1}\gamma\right\Vert _{L^{p}}\cdot\left\Vert \Indicator_{J_{0}}\right\Vert _{Z},\label{eq:DecompositionCoincidenceLemmaFirstMainAsymptotic}
\end{align}
which is our first main estimate. Note that this estimate holds for
\emph{all} choices of $z^{\left(r\right)}=\left(\smash{z_{j}^{\left(r\right)}}\right)_{j\in J_{0}}$
and of $\theta^{\left(r\right)}=\left(\smash{\theta_{j}^{\left(r\right)}}\right)_{j\in J_{0}}$,
where for each $r\in\underline{r_{0}}$ a different choice is possible.

\medskip{}

Now, recall that—by construction—$\supp\gamma_{j}^{\left(\varepsilon\right)}\subset Q_{i_{0}}$
for all $j\in J_{0}$, so that also $\supp f_{z,\theta}^{\left(r\right)}\subset Q_{i_{0}}$.
Hence, we can apply Lemma~\ref{lem:EasyNormEquivalenceFineCovering}
with $k=0$, $\gamma_{i_{0}}=f_{z,\theta}^{\left(r\right)}$ and $\gamma_{i}\equiv0$
for all $i\in I\setminus\left\{ i_{0}\right\} $, as well as $c_{i}=\delta_{i_{0}}\left(i\right)$,
$z_{i}=0$ and $\varepsilon_{i}=1$ for all $i\in I$ to conclude
\[
\left\Vert f_{z,\theta}^{\left(r\right)}\right\Vert _{\FourierDecompSp{\CalQ}pY}\leq C_{4}\cdot\left\Vert \Fourier^{-1}f_{z,\theta}^{\left(r\right)}\right\Vert _{L^{p}}\cdot\left\Vert \delta_{i_{0}}\right\Vert _{Y}
\]
for some constant $C_{4}=C_{4}\left(\dimension,p,\CalQ,C_{\CalQ,\Phi,p},\vertiii{\Gamma_{\CalQ}}_{Y\to Y}\right)$
and arbitrary choice of $z,\theta,r$. Note that Lemma~\ref{lem:EasyNormEquivalenceFineCovering}
is indeed applicable with the above choices, since we have $\delta_{i_{0}}\in\ell_{0}\left(I\right)\leq Y$.

Conversely, we can also apply Corollary~\ref{cor:EasyNormEquivalenceFineLowerBound}
with $I_{0}=\left\{ i_{0}\right\} $, $\gamma_{i_{0}}=f_{z,\theta}^{\left(r\right)}$
and $c_{i_{0}}=\varepsilon_{i_{0}}=1$ and $z_{i_{0}}=0$ to conclude
\[
\left\Vert f_{z,\theta}^{\left(r\right)}\right\Vert _{\FourierDecompSp{\CalQ}pY}\geq C_{5}^{-1}\cdot\left\Vert \Fourier^{-1}f_{z,\theta}^{\left(r\right)}\right\Vert _{L^{p}}\cdot\left\Vert \delta_{i_{0}}\right\Vert _{Y}
\]
for a further constant $C_{5}=C_{5}\left(p,\dimension,\CalQ,C_{\CalQ,\Phi,p},\vertiii{\Gamma_{\CalQ}}_{Y\to Y}\right)$
and arbitrary choices of $z,\theta,r$.

In summary, with equation~(\ref{eq:DecompositionCoincidenceLemmaFirstMainAsymptotic})
and with our assumption $\left\Vert \mybullet\right\Vert _{\FourierDecompSp{\CalQ}pY}\asymp\left\Vert \mybullet\right\Vert _{\FourierDecompSp{\CalP}pZ}$
on $\TestFunctionSpace{\CalO}$, we have thus shown
\begin{align}
\varepsilon^{\dimension\left(1-\frac{1}{p}\right)}\left\Vert \Fourier^{-1}\gamma\right\Vert _{L^{p}}\cdot\left\Vert \Indicator_{J_{0}}\right\Vert _{Z}\asymp\sum_{r=1}^{r_{0}}\left\Vert f_{z^{\left(r\right)},\theta^{\left(r\right)}}^{\left(r\right)}\right\Vert _{\FourierDecompSp{\CalP}pZ} & \asymp\sum_{r=1}^{r_{0}}\left\Vert f_{z^{\left(r\right)},\theta^{\left(r\right)}}^{\left(r\right)}\right\Vert _{\FourierDecompSp{\CalQ}pY}\nonumber \\
 & \asymp\left\Vert \delta_{i_{0}}\right\Vert _{Y}\cdot\sum_{r=1}^{r_{0}}\left\Vert \Fourier^{-1}f_{z^{\left(r\right)},\theta^{\left(r\right)}}^{\left(r\right)}\right\Vert _{L^{p}},\label{eq:CoincidenceMainArbitrageOpportunity}
\end{align}
where the implied constants only depend on $p,\dimension,C_{Z},\CalQ,\CalP,C_{\CalQ,\Phi,p},C_{\CalP,\Psi,p},\vertiii{\Gamma_{\CalQ}}_{Y\to Y},\vertiii{\Gamma_{\CalP}}_{Z\to Z}$
and on the implied constants in $\left\Vert \mybullet\right\Vert _{\FourierDecompSp{\CalQ}pY}\asymp\left\Vert \mybullet\right\Vert _{\FourierDecompSp{\CalP}pZ}$
(see equation~(\ref{eq:CoincidencePreparationLemmaAssumption})).

\medskip{}

With estimate~(\ref{eq:CoincidenceMainArbitrageOpportunity}) we
have found an opportunity for \emph{arbitrage}: The left-hand side
is independent of $z^{\left(r\right)},\theta^{\left(r\right)}$, while
the right-hand side is not, as we will see now: If we choose $\theta_{j}^{\left(r\right)}=1$
for all $j\in J_{0}$, we get using Corollary~\ref{cor:AsymptoticModulationBehaviour}
(which is applicable since $\gamma_{j}^{\left(\varepsilon\right)}\in\Schwartz\left(\R^{\dimension}\right)$
for all $j\in J_{0}$) that
\begin{align*}
\left\Vert \Fourier^{-1}f_{z^{\left(r\right)},\theta^{\left(r\right)}}^{\left(r\right)}\right\Vert _{L^{p}} & =\vphantom{\sum_{j\in J_{0}\cap J^{\left(r\right)}}}\left\Vert \Fourier^{-1}\left[\,\vphantom{\sum}\smash{\sum_{j\in J_{0}\cap J^{\left(r\right)}}}M_{z_{j}^{\left(r\right)}}\gamma_{j}^{\left(\varepsilon\right)}\,\right]\right\Vert _{L^{p}}\\
\left({\scriptstyle \text{for suitable }\left(\smash{z_{j}^{\left(r\right)}}\right)_{j\in J_{0}}\in\left(\R^{\dimension}\right)^{J_{0}}}\right) & \asymp\left\Vert \left(\left\Vert \Fourier^{-1}\gamma_{j}^{\left(\varepsilon\right)}\right\Vert _{L^{p}}\right)_{j\in J_{0}\cap J^{\left(r\right)}}\right\Vert _{\ell^{p}}=\varepsilon^{\dimension\left(1-\frac{1}{p}\right)}\left\Vert \Fourier^{-1}\gamma\right\Vert _{L^{p}}\cdot\left|J_{0}\cap J^{\left(r\right)}\right|^{1/p}.
\end{align*}
Now, note that we have $\left|J_{0}\cap J^{\left(r\right)}\right|\leq\left|J_{0}\right|=N$
for all $r\in\underline{r_{0}}$ and that $J_{0}=\biguplus_{r=1}^{r_{0}}\left(J_{0}\cap J^{\left(r\right)}\right)$,
so that there is some $r\in\underline{r_{0}}$ with $\left|J_{0}\cap J^{\left(r\right)}\right|\geq\left|J_{0}\right|/r_{0}=N/r_{0}$.
Thus, using equation~(\ref{eq:CoincidenceMainArbitrageOpportunity}),
we get (for suitable choices of the coefficients $z^{\left(r\right)}$
and $\theta^{\left(r\right)}$) that
\begin{align}
\varepsilon^{\dimension\left(1-\frac{1}{p}\right)}\left\Vert \Fourier^{-1}\gamma\right\Vert _{L^{p}}\cdot\left\Vert \Indicator_{J_{0}}\right\Vert _{Z} & \asymp\left\Vert \delta_{i_{0}}\right\Vert _{Y}\cdot\sum_{r=1}^{r_{0}}\left[\varepsilon^{\dimension\left(1-\frac{1}{p}\right)}\left\Vert \Fourier^{-1}\gamma\right\Vert _{L^{p}}\cdot\left|J_{0}\cap J^{\left(r\right)}\right|^{1/p}\right]\nonumber \\
\left({\scriptstyle \text{implied constant depending on }\CalP,p}\right) & \asymp\varepsilon^{\dimension\left(1-\frac{1}{p}\right)}\left\Vert \Fourier^{-1}\gamma\right\Vert _{L^{p}}\cdot\left\Vert \delta_{i_{0}}\right\Vert _{Y}\cdot N^{1/p}.\label{eq:CoincidenceFirstAsymptotic}
\end{align}
Now, we can cancel the common factor $\varepsilon^{\dimension\left(1-\frac{1}{p}\right)}\left\Vert \Fourier^{-1}\gamma\right\Vert _{L^{p}}$
on both sides to obtain equation~(\ref{eq:CoincidencePreparationLemmaP}).
Note in particular that all occurrences of $\varepsilon$ canceled.

\medskip{}

Next, let us consider the case $p\in\left(0,\infty\right)$. In this
case, choose $z_{j}^{\left(r\right)}=0$ for all $j\in J_{0}$ and
let $\omega=\left(\omega_{j}\right)_{j\in J_{0}}$ be a random variable
which is uniformly distributed in $\left\{ \pm1\right\} ^{J_{0}}$.
The expected value in the following calculation is taken with respect
to $\omega$. Elementary properties of the Fourier transform and Khintchine's
inequality (Theorem~\ref{thm:KhintchineInequality}) yield
\begin{align}
\mathbb{E}\left\Vert \Fourier^{-1}f_{z^{\left(r\right)},\omega}^{\left(r\right)}\right\Vert _{L^{p}}^{p} & =\mathbb{E}\int_{\R^{\dimension}}\vphantom{\sum_{j\in J_{0}\cap J^{\left(r\right)}}}\left|\,\vphantom{\sum}\smash{\sum_{j\in J_{0}\cap J^{\left(r\right)}}}\omega_{j}\cdot\left(\Fourier^{-1}\gamma_{j}^{\left(\varepsilon\right)}\right)\left(x\right)\,\right|^{p}\,\d x\nonumber \\
 & =\mathbb{E}\int_{\R^{\dimension}}\vphantom{\sum_{j\in J_{0}\cap J^{\left(r\right)}}}\left|\,\vphantom{\sum}\smash{\sum_{j\in J_{0}\cap J^{\left(r\right)}}}\omega_{j}\cdot e^{2\pi i\left\langle \xi_{j},x\right\rangle }\cdot\varepsilon^{\dimension}\cdot\left(\Fourier^{-1}\gamma\right)\left(\varepsilon x\right)\,\right|^{p}\,\d x\nonumber \\
 & =\varepsilon^{\dimension p}\cdot\int_{\R^{\dimension}}\left|\left(\Fourier^{-1}\gamma\right)\left(\varepsilon x\right)\right|^{p}\cdot\mathbb{E}\vphantom{\sum_{j\in J_{0}\cap J^{\left(r\right)}}}\left|\,\smash{\sum_{j\in J_{0}\cap J^{\left(r\right)}}}\vphantom{\sum}\omega_{j}\cdot e^{2\pi i\left\langle \xi_{j},x\right\rangle }\,\right|^{p}\,\d x\nonumber \\
\left({\scriptstyle \text{Theorem }\ref{thm:KhintchineInequality}}\right) & \asymp\varepsilon^{\dimension p}\cdot\int_{\R^{\dimension}}\left|\left(\Fourier^{-1}\gamma\right)\left(\varepsilon x\right)\right|^{p}\cdot\vphantom{\sum_{j\in J_{0}\cap J^{\left(r\right)}}}\left(\,\vphantom{\sum}\smash{\sum_{j\in J_{0}\cap J^{\left(r\right)}}}\left|e^{2\pi i\left\langle \xi_{j},x\right\rangle }\right|^{2}\,\right)^{p/2}\,\d x\nonumber \\
 & =\varepsilon^{\dimension p}\cdot\left|J_{0}\cap J^{\left(r\right)}\right|^{p/2}\cdot\int_{\R^{\dimension}}\left|\left(\Fourier^{-1}\gamma\right)\left(\varepsilon x\right)\right|^{p}\,\d x\nonumber \\
 & =\left(\varepsilon^{\dimension\left(1-\frac{1}{p}\right)}\cdot\left|J_{0}\cap J^{\left(r\right)}\right|^{1/2}\cdot\left\Vert \Fourier^{-1}\gamma\right\Vert _{L^{p}}\right)^{p},\label{eq:CoincidencePreparationKhintchine}
\end{align}
where the implied constant only depends on $p\in\left(0,\infty\right)$.
Note that the interchange of the expectation and the integral in the
preceding calculation is justified, since the expectation is just
a finite sum. Alternatively, we could have used Fubini's theorem.

In particular, estimate~(\ref{eq:CoincidencePreparationKhintchine})
yields deterministic realizations $\theta^{\left(U,r\right)}=\vphantom{\theta_{j}^{\left(U,r\right)}}\left(\smash{\theta_{j}^{\left(U,r\right)}}\right)_{j\in J_{0}}\in\left\{ \pm1\right\} ^{J_{0}}$
and $\theta^{\left(L,r\right)}=\vphantom{\theta_{j}^{\left(L,r\right)}}\left(\smash{\theta_{j}^{\left(L,r\right)}}\right)_{j\in J_{0}}\in\left\{ \pm1\right\} ^{J_{0}}$
satisfying
\[
\left\Vert \Fourier^{-1}f_{z^{\left(r\right)},\theta^{\left(L,r\right)}}^{\left(r\right)}\right\Vert _{L^{p}}\lesssim\varepsilon^{\dimension\left(1-\frac{1}{p}\right)}\left\Vert \Fourier^{-1}\gamma\right\Vert _{L^{p}}\cdot\left|J_{0}\cap J^{\left(r\right)}\right|^{1/2}\lesssim\left\Vert \Fourier^{-1}f_{z^{\left(r\right)},\theta^{\left(U,r\right)}}^{\left(r\right)}\right\Vert _{L^{p}},
\]
where the implied constants again only depend on $p\in\left(0,\infty\right)$.
Using essentially the same arguments as before equation~(\ref{eq:CoincidenceFirstAsymptotic}),
we thus conclude
\begin{align*}
\varepsilon^{\dimension\left(1-\frac{1}{p}\right)}\left\Vert \Fourier^{-1}\gamma\right\Vert _{L^{p}}\cdot\left\Vert \delta_{i_{0}}\right\Vert _{Y}\cdot N^{1/2} & \lesssim\left\Vert \delta_{i_{0}}\right\Vert _{Y}\cdot\sum_{r=1}^{r_{0}}\left\Vert \Fourier^{-1}f_{z^{\left(r\right)},\theta^{\left(L,r\right)}}^{\left(r\right)}\right\Vert _{L^{p}}\\
\left({\scriptstyle \text{equation }\eqref{eq:CoincidenceMainArbitrageOpportunity}}\right) & \asymp\varepsilon^{\dimension\left(1-\frac{1}{p}\right)}\left\Vert \Fourier^{-1}\gamma\right\Vert _{L^{p}}\cdot\left\Vert \Indicator_{J_{0}}\right\Vert _{Z}\\
\left({\scriptstyle \text{equation }\eqref{eq:CoincidenceMainArbitrageOpportunity}}\right) & \asymp\left\Vert \delta_{i_{0}}\right\Vert _{Y}\cdot\sum_{r=1}^{r_{0}}\left\Vert \Fourier^{-1}f_{z^{\left(r\right)},\theta^{\left(U,r\right)}}^{\left(r\right)}\right\Vert _{L^{p}}\\
 & \lesssim\varepsilon^{\dimension\left(1-\frac{1}{p}\right)}\left\Vert \Fourier^{-1}\gamma\right\Vert _{L^{p}}\cdot\left\Vert \delta_{i_{0}}\right\Vert _{Y}\cdot N^{1/2},
\end{align*}
which establishes equation~(\ref{eq:CoincidencePreparationLemma2}).

\medskip{}

Finally, in case of $p=\infty$, we choose $z_{j}^{\left(r\right)}=0$
and $\theta_{j}^{\left(r\right)}=1$ for all $j\in J_{0}$ and $r\in\underline{r_{0}}$
to arrive at
\begin{align*}
\left\Vert \Fourier^{-1}f_{z^{\left(r\right)},\theta^{\left(r\right)}}^{\left(r\right)}\right\Vert _{L^{p}}=\vphantom{\sum_{j\in J_{0}\cap J^{\left(r\right)}}}\left\Vert \smash{\sum_{j\in J_{0}\cap J^{\left(r\right)}}}\vphantom{\sum}\Fourier^{-1}\gamma_{j}^{\left(\varepsilon\right)}\right\Vert _{L^{\infty}} & \geq\vphantom{\sum_{j\in J_{0}\cap J^{\left(r\right)}}}\left|\smash{\sum_{j\in J_{0}\cap J^{\left(r\right)}}}\vphantom{\sum}\left(\Fourier^{-1}\gamma_{j}^{\left(\varepsilon\right)}\right)\left(0\right)\right|\\
\left({\scriptstyle \text{since }\gamma_{j}^{\left(\varepsilon\right)}\geq0}\right) & =\sum_{j\in J_{0}\cap J^{\left(r\right)}}\left\Vert \gamma_{j}^{\left(\varepsilon\right)}\right\Vert _{L^{1}}\\
\left({\scriptstyle \text{Riemann-Lebesgue}}\right) & \geq\sum_{j\in J_{0}\cap J^{\left(r\right)}}\left\Vert \Fourier^{-1}\gamma_{j}^{\left(\varepsilon\right)}\right\Vert _{L^{\infty}}.
\end{align*}
The reverse estimate $\left\Vert \Fourier^{-1}f_{z^{\left(r\right)},\theta^{\left(r\right)}}^{\left(r\right)}\right\Vert _{L^{p}}\leq\sum_{j\in J_{0}\cap J^{\left(r\right)}}\left\Vert \Fourier^{-1}\gamma_{j}^{\left(\varepsilon\right)}\right\Vert _{L^{\infty}}$
is a direct consequence of the triangle inequality. Hence, recalling
$p=\infty$, we get
\[
\left\Vert \Fourier^{-1}f_{z^{\left(r\right)},\theta^{\left(r\right)}}^{\left(r\right)}\right\Vert _{L^{p}}=\sum_{j\in J_{0}\cap J^{\left(r\right)}}\left\Vert \Fourier^{-1}\gamma_{j}^{\left(\varepsilon\right)}\right\Vert _{L^{\infty}}=\varepsilon^{\dimension\left(1-\frac{1}{p}\right)}\cdot\left\Vert \Fourier^{-1}\gamma\right\Vert _{L^{p}}\cdot\left|J_{0}\cap J^{\left(r\right)}\right|.
\]
Exactly as in equation~(\ref{eq:CoincidenceFirstAsymptotic}), this
yields
\[
\varepsilon^{\dimension\left(1-\frac{1}{p}\right)}\left\Vert \Fourier^{-1}\gamma\right\Vert _{L^{p}}\cdot\left\Vert \Indicator_{J_{0}}\right\Vert _{Z}\asymp\left\Vert \delta_{i_{0}}\right\Vert _{Y}\cdot\sum_{r=1}^{r_{0}}\left\Vert \Fourier^{-1}f_{z^{\left(r\right)},\theta^{\left(r\right)}}^{\left(r\right)}\right\Vert _{L^{p}}\asymp\varepsilon^{\dimension\left(1-\frac{1}{p}\right)}\left\Vert \Fourier^{-1}\gamma\right\Vert _{L^{p}}\cdot\left\Vert \delta_{i_{0}}\right\Vert _{Y}\cdot N,
\]
which establishes equation~(\ref{eq:CoincidencePreparationLemma1}).
\end{proof}
Now, we can finally prove the announced necessary criteria for the
equality of two decomposition spaces:
\begin{thm}
\label{thm:NecessaryCriterionForCoincidenceOfDecompositionSpaces}Let
$\emptyset\neq\CalO\subset\R^{\dimension}$ be open, let $p_{1},p_{2}\in\left(0,\infty\right]$,
let $\CalQ=\left(Q_{i}\right)_{i\in I}$ be an \emph{open} $L^{p_{1}}$-decomposition
covering of $\CalO$, and let $\CalP=\left(P_{j}\right)_{j\in J}$
be an \emph{open} $L^{p_{2}}$-decomposition covering of $\CalO$.
Let $Y\subset\Compl^{I}$ and $Z\subset\Compl^{J}$ be $\CalQ$-regular
and $\CalP$-regular, respectively and assume $\ell_{0}\left(I\right)\subset Y$
and $\ell_{0}\left(J\right)\subset Z$.

Finally, assume that
\[
\FourierDecompSp{\CalQ}{p_{1}}Y=\FourierDecompSp{\CalP}{p_{2}}Z
\]
holds, with equivalent (quasi)-norms. Then the following hold:

\begin{enumerate}
\item We have $p_{1}=p_{2}$.
\item \label{enu:CoincidenceSequenceSpaceSimilar}There is a constant $C_{1}>0$
with 
\[
C_{1}^{-1}\cdot\left\Vert \delta_{i}\right\Vert _{Y}\leq\left\Vert \delta_{j}\right\Vert _{Z}\leq C_{1}\cdot\left\Vert \delta_{i}\right\Vert _{Y}\text{ for all }i\in I\text{ and }j\in J\text{ with }Q_{i}\cap P_{j}\neq\emptyset.
\]
\item In case of $p_{1}\neq2$, $\CalQ$ and $\CalP$ are weakly equivalent
(see Definition~\ref{def:RelativeIndexClustersSubordinateCoveringsModerateCoverings}).
\item In case of $Y=\ell_{w}^{q_{1}}\left(I\right)$ and $Z=\ell_{v}^{q_{2}}\left(J\right)$
for a $\CalQ$-moderate weight $w=\left(w_{i}\right)_{i\in I}$ and
a $\CalP$-moderate weight $v=\left(v_{j}\right)_{j\in J}$ and certain
$q_{1},q_{2}\in\left(0,\infty\right]$, the following hold:

\begin{enumerate}
\item We have $q_{1}=q_{2}$.
\item \label{enu:CoincidenceWeightEquivalence}There is a constant $C_{1}>0$
with 
\[
C_{1}^{-1}\cdot w_{i}\leq v_{j}\leq C_{1}\cdot w_{i}\text{ for all }i\in I\text{ and }j\in J\text{ with }Q_{i}\cap P_{j}\neq\emptyset.
\]
\item If $\left(p_{1},q_{1}\right)\neq\left(2,2\right)$, the coverings
$\CalQ,\CalP$ are weakly equivalent.\qedhere
\end{enumerate}
\end{enumerate}
\end{thm}

\begin{rem*}
(1) Note that the theorem only yields \emph{weak} equivalence of $\CalQ,\CalP$.
Nevertheless, if all sets of the coverings $\CalQ,\CalP$ are connected
(which is usually the case), then these sets are even path-connected,
since both notions of connectedness are equivalent for open subsets
of $\R^{\dimension}$ and since we assume $\CalQ,\CalP$ to be open
coverings. But then, Corollary~\ref{cor:WeakSubordinationImpliesSubordinationIfConnected}
shows that $\CalQ,\CalP$ are equivalent, not only weakly equivalent.

\medskip{}

(2) It is not necessary to assume that the (quasi)-norms on both spaces
are equivalent. Indeed, using the completeness of decomposition spaces,
a form of the closed graph theorem (see \cite[Theorem 2.15]{RudinFA})
and the continuous embeddings $\FourierDecompSp{\CalQ}{p_{1}}Y\hookrightarrow\DistributionSpace{\CalO}$
and $\FourierDecompSp{\CalP}{p_{2}}Z\hookrightarrow\DistributionSpace{\CalO}$
(see Theorem~\ref{thm:DecompositionSpaceComplete}), it is not hard
to see that an equality of the two decomposition spaces as sets already
implies that the (quasi)-norms are equivalent.

\medskip{}

(3) Instead of assuming $\FourierDecompSp{\CalQ}{p_{1}}Y=\FourierDecompSp{\CalP}{p_{2}}Z$,
the proof shows that it suffices to assume $\left\Vert f\right\Vert _{\FourierDecompSp{\CalQ}{p_{1}}Y}\asymp\left\Vert f\right\Vert _{\FourierDecompSp{\CalP}{p_{2}}Z}$
for all $f\in\TestFunctionSpace{\CalO}$.

\medskip{}

(4) The requirement that $\CalQ,\CalP$ are \emph{open} coverings
is made mostly for convenience. It is not hard to see that we could
always replace $\CalQ$ and $\CalP$ by the open coverings $\CalQ^{\circ}$
and $\CalP^{\circ}$ and get the same decomposition spaces. The main
reason for this is that for any subordinate partition of unity $\left(\varphi_{i}\right)_{i\in I}$,
we already have $\varphi_{i}\equiv0$ outside of $Q_{i}^{\circ}$.
Note though that $Q_{i}\cap P_{j}\neq\emptyset$ is not necessarily
equivalent to $Q_{i}^{\circ}\cap P_{j}^{\circ}\neq\emptyset$. To
avoid distinctions of this type, we directly assume $\CalQ$ and $\CalP$
to be open.

\medskip{}

(5) In contrast to our other theorems, the conditions in the theorem
above are all stated \emph{qualitatively}, that is, no (explicit)
bounds on the quantities $N\left(\CalQ,\CalP\right)$ and $N\left(\CalP,\CalQ\right)$
are provided. As the proof shows, it is still true that we have $N\left(\CalQ,\CalP\right)\leq C$
and $N\left(\CalP,\CalQ\right)\leq C$, where $C>0$ depends on the
bounds of the norm equivalence $\left\Vert \mybullet\right\Vert _{\FourierDecompSp{\CalQ}{p_{1}}Y}\asymp\left\Vert \mybullet\right\Vert _{\FourierDecompSp{\CalP}{p_{2}}Z}$
and on the usual quantities like $p_{1},q_{1},\CalQ,\CalP$, etc.
In this case, we prefer the compact form of the theorem, since we
will not need to know the precise dependency of the constants in any
relevant application.
\end{rem*}
\begin{proof}[Proof of Theorem~\ref{thm:NecessaryCriterionForCoincidenceOfDecompositionSpaces}]
A straightforward application of Lemma~\ref{lem:SimpleNecessaryCondition}
(with $K=\CalO$) shows that we have $p_{1}\leq p_{2}$ and $p_{2}\leq p_{1}$
and thus $p_{1}=p_{2}$. Given this identity, the same lemma also
yields existence of $C_{1}>0$ as required in (\ref{enu:CoincidenceSequenceSpaceSimilar})
and (\ref{enu:CoincidenceWeightEquivalence}). This uses that (by
assumption) $\delta_{i}\in Y$ and $\delta_{j}\in Z$ for arbitrary
$i\in I$ and $j\in J$ and that $\CalQ,\CalP$ are open coverings
of $\CalO$, so that $Q_{i}\cap P_{j}\neq\emptyset$ is equivalent
to $K^{\circ}\cap Q_{i}^{\circ}\cap P_{j}^{\circ}\neq\emptyset$.
Furthermore, we used $\left\Vert \delta_{i}\right\Vert _{\ell_{w}^{q_{1}}}=w_{i}$
and $\left\Vert \delta_{j}\right\Vert _{\ell_{v}^{q_{2}}}=v_{j}$
for all $i\in I$ and $j\in J$.

\medskip{}

Let us first show that $\CalQ,\CalP$ are weakly equivalent if $p_{1}\neq2$.
Note that (because of $p_{1}=p_{2}$), our assumptions are symmetric
in $\CalQ,\CalP$, so that it suffices to show that $\CalQ$ is weakly
subordinate to $\CalP$, i.e.\@ that $\sup_{i\in I}\left|J_{i}\right|<\infty$.
To this end, first assume $p_{1}\in\left(0,\infty\right)$ and let
$i_{0}\in I$ be arbitrary. For pairwise distinct $j_{1},\dots,j_{N}\in J_{i_{0}}$,
Lemma~\ref{lem:CoincidencePreparationLemma} (precisely, equations
(\ref{eq:CoincidencePreparationLemmaP}) and (\ref{eq:CoincidencePreparationLemma2}))
yields
\[
\left\Vert \delta_{i_{0}}\right\Vert _{Y}\cdot N^{1/p_{1}}\asymp\left\Vert \Indicator_{\left\{ j_{1},\dots,j_{N}\right\} }\right\Vert _{Z}\asymp\left\Vert \delta_{i_{0}}\right\Vert _{Y}\cdot N^{1/2}\:,
\]
and hence $N^{1/p_{1}}\asymp N^{1/2}$, where the implied constants
are independent of $i_{0}\in I$ and of $N$. Since we assume $p_{1}\neq2$,
it is easy to see that this can only hold if $N\leq C$ for some constant
$C>0$ (independent of $i_{0}\in I$). This yields $\sup_{i\in I}\left|J_{i}\right|<\infty$
as desired.

Finally, in case of $p_{1}=\infty$, we again use Lemma~\ref{lem:CoincidencePreparationLemma}
(this time equations (\ref{eq:CoincidencePreparationLemmaP}) and
(\ref{eq:CoincidencePreparationLemma1})) to deduce $1=N^{1/p_{1}}\asymp N$,
which again yields $\sup_{i\in I}\left|J_{i}\right|<\infty$, as desired.

\medskip{}

We have now completed the proof of the first three claims of the theorem,
so that we can assume $Y=\ell_{w}^{q_{1}}\left(I\right)$ and $Z=\ell_{v}^{q_{2}}\left(J\right)$
for the remainder of the proof.

For brevity, set $p:=p_{1}=p_{2}$. Let $\Phi=\left(\varphi_{i}\right)_{i\in I}$
and $\Psi=\left(\psi_{j}\right)_{j\in J}$ be $L^{p}$-BAPUs for $\CalQ$
and $\CalP$, respectively. Furthermore, fix a nontrivial, \emph{nonnegative}
function $\gamma\in\TestFunctionSpace{B_{1}\left(0\right)}$ for the
rest of the proof. For $\varepsilon>0$, define $\gamma^{\left(\varepsilon\right)}:=\gamma\left(\varepsilon^{-1}\mybullet\right)\in\TestFunctionSpace{B_{\varepsilon}\left(0\right)}$.
Finally, fix $C\geq1$ with
\[
C^{-1}\cdot\left\Vert f\right\Vert _{\FourierDecompSp{\CalQ}p{\ell_{w}^{q_{1}}}}\leq\left\Vert f\right\Vert _{\FourierDecompSp{\CalP}p{\ell_{v}^{q_{2}}}}\leq C\cdot\left\Vert f\right\Vert _{\FourierDecompSp{\CalQ}p{\ell_{w}^{q_{1}}}}\qquad\forall\,f\in\TestFunctionSpace{\CalO}.
\]

Let us first show $q_{1}=q_{2}$. By symmetry, it suffices to show
$q_{1}\leq q_{2}$. To this end, let $N\in\N$ be arbitrary and assume
that there are sequences $i_{1},\dots,i_{N}\in I$ and $j_{1},\dots,j_{N}\in J$
with $Q_{i_{\ell}}\cap P_{j_{\ell}}\neq\emptyset$ for all $\ell\in\underline{N}$
and with $Q_{i_{\ell}}^{\ast}\cap Q_{i_{k}}^{\ast}=\emptyset=P_{j_{\ell}}^{\ast}\cap P_{j_{k}}^{\ast}$
for all $\ell,k\in\underline{N}$ with $\ell\neq k$. We will see
below that such a sequence indeed exists. Since the family $\left(Q_{i_{\ell}}\cap P_{j_{\ell}}\right)_{\ell\in\underline{N}}$
is a finite family of nonempty open(!) subsets of $\CalO$, there
is some $\varepsilon>0$ and a sequence $\xi_{1},\dots,\xi_{N}\in\CalO$
with $B_{\varepsilon}\left(\xi_{\ell}\right)\subset Q_{i_{\ell}}\cap P_{j_{\ell}}$
for all $\ell\in\underline{N}$. Note that $\varepsilon>0$ may depend
heavily on $N\in\N$ and on $i_{1},\dots,i_{N}$ and $j_{1},\dots,j_{N}$,
but—as we will see—all dependencies on $\varepsilon$ that are relevant
for us will cancel in the end.

As an auxiliary result, note that $\ell\neq k$ implies $i_{\ell}\neq i_{k}$;
indeed, we have $Q_{i_{\ell}}^{\ast}\cap Q_{i_{k}}^{\ast}=\emptyset$
since $\ell\neq k$. But $i_{\ell}=i_{k}$ would imply $Q_{i_{\ell}}^{\ast}\cap Q_{i_{k}}^{\ast}=Q_{i_{\ell}}^{\ast}\supset Q_{i_{\ell}}\supset Q_{i_{\ell}}\cap P_{j_{\ell}}\neq\emptyset$.
The same argument also shows that $\ell\neq k$ implies $j_{\ell}\neq j_{k}$.

Now, set $\gamma_{\ell}^{\left(\varepsilon\right)}:=L_{\xi_{\ell}}\gamma_{\ell}^{\left(\varepsilon\right)}\in\TestFunctionSpace{B_{\varepsilon}\left(\xi_{\ell}\right)}$
for $\ell\in\underline{N}$ and define
\[
f_{N}:=\sum_{\ell=1}^{N}w_{i_{\ell}}^{-1}\cdot\gamma_{\ell}^{\left(\varepsilon\right)}\in\TestFunctionSpace{\CalO}\subset\FourierDecompSp{\CalQ}p{\ell_{w}^{q_{1}}}\cap\FourierDecompSp{\CalP}p{\ell_{v}^{q_{2}}}.
\]
By assumption, we have
\begin{equation}
\left\Vert f_{N}\right\Vert _{\FourierDecompSp{\CalP}p{\ell_{v}^{q_{2}}}}\leq C\cdot\left\Vert f_{N}\right\Vert _{\FourierDecompSp{\CalQ}p{\ell_{w}^{q_{1}}}}\qquad\forall\,N\in\N.\label{eq:NecessaryForCoincidenceQsCoincideMainEstimate}
\end{equation}
To exploit this estimate, we will now obtain an upper bound on $\left\Vert f_{N}\right\Vert _{\FourierDecompSp{\CalQ}p{\ell_{w}^{q_{1}}}}$
and a lower bound on $\left\Vert f_{N}\right\Vert _{\FourierDecompSp{\CalP}p{\ell_{v}^{q_{2}}}}$.
Since each $\gamma_{\ell}^{\left(\varepsilon\right)}$ vanishes outside
of $B_{\varepsilon}\left(\xi_{\ell}\right)\subset Q_{i_{\ell}}$ (and
since we have $i_{\ell}\neq i_{m}$ for $\ell\neq m$), Lemma~\ref{lem:EasyNormEquivalenceFineCovering}
yields a constant $L_{1}>0$ (which does \emph{not} depend on $N$)
with
\begin{align}
\left\Vert f_{N}\right\Vert _{\FourierDecompSp{\CalQ}p{\ell_{w}^{q_{1}}}} & \leq L_{1}\cdot\left\Vert \left(w_{i_{\ell}}\cdot\left\Vert \Fourier^{-1}\left[w_{i_{\ell}}^{-1}\cdot\gamma_{\ell}^{\left(\varepsilon\right)}\right]\right\Vert _{L^{p}}\right)_{\ell\in\underline{N}}\right\Vert _{\ell^{q_{1}}}\nonumber \\
 & =L_{1}\cdot\left\Vert \left(\varepsilon^{\dimension\left(1-\frac{1}{p}\right)}\cdot\left\Vert \Fourier^{-1}\gamma\right\Vert _{L^{p}}\right)_{\ell\in\underline{N}}\right\Vert _{\ell^{q_{1}}}\nonumber \\
 & =L_{1}\cdot\varepsilon^{\dimension\left(1-\frac{1}{p}\right)}\cdot\left\Vert \Fourier^{-1}\gamma\right\Vert _{L^{p}}\cdot N^{1/q_{1}}.\label{eq:NecessaryForCoincidenceQsCoincideUpperEstimate}
\end{align}

Likewise, Corollary~\ref{cor:EasyNormEquivalenceFineLowerBound}
(with $\CalP$ instead of $\CalQ$ and with $I_{0}=\left\{ j_{1},\dots,j_{N}\right\} $
and $k=0$) yields a constant $L_{2}>0$ (independent of $N\in\N$)
with
\begin{align}
\left\Vert f_{N}\right\Vert _{\FourierDecompSp{\CalP}p{\ell_{v}^{q_{2}}}} & \geq L_{2}^{-1}\cdot\left\Vert \left(v_{j_{\ell}}\cdot\left\Vert \Fourier^{-1}\left(w_{i_{\ell}}^{-1}\cdot\gamma_{\ell}^{\left(\varepsilon\right)}\right)\right\Vert _{L^{p}}\right)_{\ell\in\underline{N}}\right\Vert _{\ell^{q_{2}}}\nonumber \\
\left({\scriptstyle v_{j_{\ell}}\geq C_{1}^{-1}\cdot w_{i_{\ell}}\text{ since }Q_{i_{\ell}}\cap P_{j_{\ell}}\neq\emptyset}\right) & \geq\left(C_{1}L_{2}\right)^{-1}\cdot\varepsilon^{\dimension\left(1-\frac{1}{p}\right)}\cdot\left\Vert \Fourier^{-1}\gamma\right\Vert _{L^{p}}\cdot N^{1/q_{2}}.\label{eq:NecessaryForCoincidenceQsCoincideLowerEstimate}
\end{align}
Here, the prerequisites of Corollary~\ref{cor:EasyNormEquivalenceFineLowerBound}
are indeed satisfied, since we have $\supp\gamma_{\ell}^{\left(\varepsilon\right)}\subset P_{j_{\ell}}=P_{j_{\ell}}^{k\ast}$
for all $\ell\in\underline{N}$ and because of $P_{j_{\ell}}^{\left(k+1\right)\ast}\cap P_{j_{n}}^{\left(k+1\right)\ast}=P_{j_{\ell}}^{\ast}\cap P_{j_{n}}^{\ast}=\emptyset$
as soon as $j_{\ell}\neq j_{n}$ (which implies $\ell\neq n$). Finally,
we also used $j_{\ell}\neq j_{n}$ for $\ell\neq n$—which was shown
above—since this implies that the map $I_{0}=\left\{ j_{1},\dots,j_{N}\right\} \to\TestFunctionSpace{\CalO},j_{\ell}\mapsto w_{i_{\ell}}^{-1}\cdot\gamma_{\ell}^{\left(\varepsilon\right)}$
is well-defined.

Finally, a combination of inequalities (\ref{eq:NecessaryForCoincidenceQsCoincideMainEstimate}),
(\ref{eq:NecessaryForCoincidenceQsCoincideUpperEstimate}) and (\ref{eq:NecessaryForCoincidenceQsCoincideLowerEstimate})
yields—after canceling the common factor $\varepsilon^{\dimension\left(1-\frac{1}{p}\right)}\cdot\left\Vert \Fourier^{-1}\gamma\right\Vert _{L^{p}}>0$—the
estimate
\[
N^{1/q_{2}}\leq CC_{1}L_{1}L_{2}\cdot N^{1/q_{1}}.
\]
Since the constant $CC_{1}L_{1}L_{2}$ is independent of $N$ and
since this estimate holds for all $N\in\N$, we conclude $q_{2}^{-1}\leq q_{1}^{-1}$
and thus $q_{1}\leq q_{2}$ as desired.

\medskip{}

To complete the proof of $q_{1}=q_{2}$, it remains to construct sequences
of indices $i_{1},\dots,i_{N}$ and $j_{1},\dots,j_{N}$ with $Q_{i_{\ell}}\cap P_{j_{\ell}}\neq\emptyset$
for all $\ell\in\underline{N}$ and so that $\left(Q_{i_{\ell}}^{\ast}\right)_{\ell\in\underline{N}}$
and $\left(P_{j_{\ell}}^{\ast}\right)_{\ell\in\underline{N}}$ are
two sequences of pairwise disjoint sets. This is done by induction
on $N\in\N$: If $i_{1},\dots,i_{N-1}$ and $j_{1},\dots,j_{N-1}$
are already constructed, note that
\[
M:=\bigcup_{\ell=1}^{N-1}\overline{Q_{i_{\ell}}^{3\ast}\cup P_{j_{\ell}}^{3\ast}}
\]
is a compact subset of $\CalO$ by Lemma~\ref{lem:PartitionCoveringNecessary}.
Since $\CalO\subset\R^{\dimension}$ is a nonempty open set, we conclude\footnote{Otherwise, $M=\CalO$ would be compact and open, so that connectedness
of $\R^{\dimension}$ yields $M=\CalO\in\left\{ \emptyset,\R^{\dimension}\right\} $,
which is impossible, since $\CalO$ is nonempty and $M$ is compact.} $M\subsetneq\CalO$. Thus, there is some $\xi\in\CalO\setminus M$.
Since $\CalQ,\CalP$ both cover $\CalO$, there are $i_{N}\in I$
and $j_{N}\in J$ with $\xi\in Q_{i_{N}}\cap P_{j_{N}}\neq\emptyset$.
It remains to verify that $\left(Q_{i_{\ell}}^{\ast}\right)_{\ell\in\underline{N}}$
and $\left(P_{j_{\ell}}^{\ast}\right)_{\ell\in\underline{N}}$ are
both pairwise disjoint. By symmetry, it suffices to consider the first
sequence. Also, since $\left(Q_{i_{\ell}}^{\ast}\right)_{\ell\in\underline{N-1}}$
is pairwise disjoint, we only need to show $Q_{i_{N}}^{\ast}\cap Q_{i_{\ell}}^{\ast}=\emptyset$
for all $\ell\in\underline{N-1}$. If this was false, we would have
$Q_{i}\cap Q_{i_{\ell}}^{\ast}\neq\emptyset$ for some $i\in i_{N}^{\ast}$.
Hence, $i\in i_{\ell}^{2\ast}$, which implies $i_{N}\in i_{\ell}^{3\ast}$
and thus $\xi\in Q_{i_{N}}\subset Q_{i_{\ell}}^{3\ast}\subset M$,
in contradiction to $\xi\in\CalO\setminus M$. All in all, we have
thus shown $q_{1}=q_{2}$. For brevity, let us set $q:=q_{1}$ for
the remainder of the proof.

\medskip{}

It remains to show that $\CalQ$ and $\CalP$ are weakly equivalent
if $\left(p_{1},q_{1}\right)\neq\left(2,2\right)$. In case of $p_{1}\neq2$,
this follows from our general considerations above. Hence, we can
assume $p_{1}=2$ and thus $q_{1}=q_{2}\neq2$, since $\left(p_{1},q_{1}\right)\neq\left(2,2\right)$.
As above, our assumptions are symmetric in $\CalQ,\CalP$, so that
it suffices to show that $\CalQ$ is weakly subordinate to $\CalP$,
i.e.\@ that $\sup_{i\in I}\left|J_{i}\right|<\infty$. To this end,
let $i_{0}\in I$ be arbitrary and let $j_{1},\dots,j_{N}\in J_{i_{0}}$
be pairwise distinct.

Note that we have $p_{1}=2<\infty$. Thus, Lemma~\ref{lem:CoincidencePreparationLemma}
(precisely, equation~(\ref{eq:CoincidencePreparationLemma2})) shows
that
\begin{align*}
w_{i_{0}}\cdot N^{1/2}=\left\Vert \delta_{i_{0}}\right\Vert _{Y}\cdot N^{1/2} & \asymp\left\Vert \Indicator_{\left\{ j_{1},\dots,j_{N}\right\} }\right\Vert _{Z}\\
 & =\left\Vert \left(v_{j}\cdot\Indicator_{\left\{ j_{1},\dots,j_{N}\right\} }\left(j\right)\right)_{j\in J}\right\Vert _{\ell^{q_{2}}}\\
\left({\scriptstyle \text{since }v_{j}\asymp w_{i_{0}}\text{ for }j\in\left\{ j_{1},\dots,j_{N}\right\} \subset J_{i_{0}}}\right) & \asymp w_{i_{0}}\cdot\left\Vert \Indicator_{\left\{ j_{1},\dots,j_{N}\right\} }\right\Vert _{\ell^{q_{2}}}\\
\left({\scriptstyle \text{since }q_{2}=q_{1}}\right) & =w_{i_{0}}\cdot N^{1/q_{1}},
\end{align*}
where the implied constants are independent of $i_{0}\in I$ and of
$N$. But since we assume $q_{1}\neq2$, this is only possible if
$N\leq C$, for some constant $C>0$ (independent of $i_{0}\in I$).
This shows $\sup_{i\in I}\left|J_{i}\right|<\infty$, as desired.
\end{proof}
In the theorem above, we assumed that both $\CalQ$ and $\CalP$ are
coverings of the \emph{same} set $\CalO\subset\R^{\dimension}$. To
complement this result, our next theorem shows that the decomposition
spaces $\DecompSp{\CalQ}{p_{1}}{\ell_{w}^{q_{1}}}$ and $\DecompSp{\CalP}{p_{2}}{\ell_{v}^{q_{2}}}$
can—except in trivial cases—\emph{never} coincide if $\CalQ$ and
$\CalP$ cover two \emph{different} sets $\CalO,\CalO'$.

\begin{customthm}{\thethm\relax$\frac{\boldsymbol{1}}{\boldsymbol{2}}$}\manuallabel{\thethm\relax$\frac{1}{2}$}\label{thm:NoncoincidenceForDifferentOrbits}Let
$\emptyset\neq\CalO,\CalO'\subset\R^{\dimension}$ be open, let $p_{1},p_{2}\in\left(0,\infty\right]$,
let $\CalQ=\left(Q_{i}\right)_{i\in I}$ be an \emph{open} $L^{p_{1}}$-decomposition
covering of $\CalO$, and let $\CalP=\left(P_{j}\right)_{j\in J}$
be an \emph{open} $L^{p_{2}}$-decomposition covering of $\CalO'$.
Furthermore, let $Y\subset\Compl^{I}$ and $Z\subset\Compl^{J}$
be $\CalQ$-regular and $\CalP$-regular, respectively, and assume
that $\ell_{0}\left(I\right)\leq Y$ and that $\ell_{0}\left(J\right)\leq Z$.

Assume that $\CalO'\cap\partial\CalO\neq\emptyset$ and that we have
\[
\left\Vert f\right\Vert _{\FourierDecompSp{\CalQ}{p_{1}}Y}\asymp\left\Vert f\right\Vert _{\FourierDecompSp{\CalP}{p_{2}}Z}\qquad\forall\,f\in\TestFunctionSpace{\CalO\cap\CalO'},
\]
where the implied constant is independent of $f$. Then the following
hold:
\begin{enumerate}
\item We have $p_{1}=p_{2}=2$.
\item We have $\left\Vert \delta_{i}\right\Vert _{Y}\asymp\left\Vert \delta_{j}\right\Vert _{Z}$
for all $i\in I$ and $j\in J$ with $Q_{i}\cap P_{j}\neq\emptyset$,
with the implied constant independent of $i,j$.
\item If $Y=\ell_{w}^{q_{1}}\left(I\right)$ with a $\CalQ$-moderate weight
$w=\left(w_{i}\right)_{i\in I}$ and some $q_{1}\in\left(0,\infty\right]$,
then $q_{1}=2$.
\item Assume that

\begin{enumerate}
\item $\CalO\cap\CalO'$ is unbounded, or $\overline{\CalO\cap\CalO'}\nsubseteq\CalO\cup\CalO'$.
\item $Y=\ell_{w}^{q_{1}}\left(I\right)$ and $Z=\ell_{v}^{q_{2}}\left(J\right)$
for certain $q_{1},q_{2}\in\left(0,\infty\right]$ and weights $w=\left(w_{i}\right)_{i\in I}$
and $v=\left(v_{j}\right)_{j\in J}$ which are $\CalQ$-moderate and
$\CalP$-moderate, respectively.
\end{enumerate}
\noindent Then $q_{1}=q_{2}=2$ and $w_{i}\asymp v_{j}$ for all $i\in I$
and $j\in J$ with $Q_{i}\cap P_{j}\neq\emptyset$.\qedhere
\end{enumerate}
\end{customthm}
\begin{rem*}
The reason for the unusual theorem numbering is that I wanted to avoid
changing theorem numbers with respect to the earlier version of the
manuscript.
\end{rem*}
\begin{proof}
By assumption, there is some $\xi_{0}\in\CalO'\cap\partial\CalO\subset\CalO'\cap\overline{\CalO}$.
Since $\CalO'$ is open, this easily implies that $K:=\CalO\cap\CalO'$
is nonempty (and open). With this choice of $K$, Lemma~\ref{lem:SimpleNecessaryCondition}
yields $p_{1}=p_{2}$. For brevity, we will write $p:=p_{1}=p_{2}$
in the following. Having $p_{1}=p_{2}$, another application of Lemma~\ref{lem:SimpleNecessaryCondition}
shows $\left\Vert \delta_{i}\right\Vert _{Y}\asymp\left\Vert \delta_{j}\right\Vert _{Z}$
in case of $Q_{i}\cap P_{j}\neq\emptyset$.

Next, let us fix a \emph{nonnegative}, nonzero function $\gamma\in\TestFunctionSpace{B_{1}\left(0\right)}$.
By the disjointization lemma (Lemma~\ref{lem:DisjointizationPrinciple}),
we get finite partitions $I=\biguplus_{r=1}^{R_{1}}I^{\left(r\right)}$
and $J=\biguplus_{r=1}^{R_{2}}J^{\left(r\right)}$ such that $Q_{i}^{\ast}\cap Q_{\ell}^{\ast}=\emptyset$
for all $i,\ell\in I^{\left(r\right)}$ (for any $r\in\underline{R_{1}}$)
with $i\neq\ell$, and analogously for the covering $\CalP$. Furthermore,
since $\CalO'\cap\partial\CalO\neq\emptyset$, Lemma~\ref{lem:DifferentOrbitsPreventModerateness}
shows that there is some $j_{0}\in J$ such that $I_{j_{0}}$ is infinite;
thus, there is some $r\in\underline{R_{1}}$ such that $I_{j_{0}}\cap I^{\left(r\right)}$
is infinite.

Now, let $N\in\N$ be arbitrary, and choose pairwise distinct indices
$i_{1},\dots,i_{N}\in I_{j_{0}}\cap I^{\left(r\right)}$. Since each
of the sets $P_{j_{0}}\cap Q_{i_{\ell}}$ is open and nonempty, there
is some $\varepsilon>0$ and for each $\ell\in\underline{N}$ some
$\omega_{\ell}\in\R^{\dimension}$ with $B_{\varepsilon}\left(\omega_{\ell}\right)\subset P_{j_{0}}\cap Q_{i_{\ell}}\subset\CalO\cap\CalO'$.
For $n\in\underline{N}$, let $\gamma_{n}:=\gamma_{i_{n}}:=L_{\omega_{n}}\left(\gamma\circ\varepsilon^{-1}\identity\right)$
and note $\gamma_{n}\in\TestFunctionSpace{B_{\varepsilon}\left(\omega_{n}\right)}\subset\TestFunctionSpace{P_{j_{0}}\cap Q_{i_{n}}}\subset\TestFunctionSpace{\CalO\cap\CalO'}$.
Now, for vectors $c=\left(c_{1},\dots,c_{N}\right)\in\Compl^{N}$,
$\sigma=\left(\sigma_{1},\dots,\sigma_{N}\right)\in\left\{ \pm1\right\} ^{N}$
and $\Xi=\left(\xi_{1},\dots,\xi_{N}\right)\in\left(\R^{\dimension}\right)^{N}$,
let 
\[
f_{N}^{\left(c,\sigma,\Xi\right)}:=\sum_{n=1}^{N}M_{\xi_{n}}\left(\sigma_{n}c_{n}\cdot\gamma_{n}\right)\in\TestFunctionSpace{\CalO\cap\CalO'},
\]
and set $\xi_{i_{n}}:=\xi_{n}$, $c_{i_{n}}:=c_{n}$ and $\sigma_{i_{n}}:=\sigma_{n}$
for $n\in\underline{N}$. With these definitions, a combination of
Lemma~\ref{lem:EasyNormEquivalenceFineCovering} and Corollary~\ref{cor:EasyNormEquivalenceFineLowerBound}
(with $I_{0}=\left\{ i_{1},\dots,i_{N}\right\} $) yields
\[
\begin{split}\left\Vert f_{N}^{\left(c,\sigma,\Xi\right)}\right\Vert _{\FourierDecompSp{\CalQ}{p_{1}}Y}=\left\Vert \smash{\sum_{i\in I_{0}}}\vphantom{\sum}M_{\xi_{i}}\left(\sigma_{i}c_{i}\cdot\gamma_{i}\right)\right\Vert _{\FourierDecompSp{\CalQ}pY} & \asymp\left\Vert \left(c_{i}\cdot\left\Vert \Fourier^{-1}\gamma_{i}\right\Vert _{L^{p}}\right)_{i\in I_{0}}\right\Vert _{Y|_{I_{0}}}\\
 & =\left\Vert \Fourier^{-1}\gamma\right\Vert _{L^{p}}\cdot\varepsilon^{\dimension\left(1-\frac{1}{p}\right)}\cdot\left\Vert \left(c_{i}\right)_{i\in I_{0}}\right\Vert _{Y|_{I_{0}}}\,\,,
\end{split}
\]
where the implied constant is independent of the choice of $N\in\N$
and $\varepsilon>0$ and of $c,\sigma,\Xi$. Note that $\varepsilon$
depends crucially on the choice of $N\in\N$; but fortunately, all
occurrences of $\varepsilon$ will cancel in the end.

Next, another combination of Lemma~\ref{lem:EasyNormEquivalenceFineCovering}
and Corollary~\ref{cor:EasyNormEquivalenceFineLowerBound} (this
time with $\CalP$ instead of $\CalQ$ and with $I_{0}=\left\{ j_{0}\right\} $)
shows because of $f_{N}^{\left(c,\sigma,\Xi\right)}\in\TestFunctionSpace{P_{j_{0}}}$
that
\begin{align*}
\left\Vert \smash{f_{N}^{\left(c,\sigma,\Xi\right)}}\right\Vert _{\FourierDecompSp{\CalP}{p_{2}}Z} & \asymp\left\Vert \delta_{j_{0}}\right\Vert _{Z}\cdot\left\Vert \Fourier^{-1}\left(\smash{f_{N}^{\left(c,\sigma,\Xi\right)}}\right)\right\Vert _{L^{p}}\\
\left({\scriptstyle \text{implied constant may depend on }j_{0},\text{ but not on }N,c,\sigma,\Xi,\varepsilon}\right) & \asymp\left\Vert \,\smash{\sum_{n=1}^{N}}\,\vphantom{\sum}L_{-\xi_{n}}\left(\sigma_{n}c_{n}\cdot\Fourier^{-1}\gamma_{n}\right)\,\right\Vert _{L^{p}}\vphantom{\sum_{n=1}^{N}}=:E_{N}^{\left(c,\sigma,\Xi\right)}.
\end{align*}
But our assumptions ensure $\vphantom{f^{\left(c,\sigma,\Xi\right)}}\left\Vert \smash{f_{N}^{\left(c,\sigma,\Xi\right)}}\right\Vert _{\FourierDecompSp{\CalQ}{p_{1}}Y}\asymp\left\Vert \smash{f_{N}^{\left(c,\sigma,\Xi\right)}}\right\Vert _{\FourierDecompSp{\CalP}{p_{2}}Z}$,
so that we get
\[
\left\Vert \Fourier^{-1}\gamma\right\Vert _{L^{p}}\cdot\varepsilon^{\dimension\left(1-\frac{1}{p}\right)}\cdot\left\Vert \left(c_{i}\right)_{i\in I_{0}}\right\Vert _{Y|_{I_{0}}}\asymp\left\Vert f_{N}^{\left(c,\sigma,\Xi\right)}\right\Vert _{\FourierDecompSp{\CalQ}{p_{1}}Y}\asymp\left\Vert f_{N}^{\left(c,\sigma,\Xi\right)}\right\Vert _{\FourierDecompSp{\CalP}{p_{2}}Z}\asymp E_{N}^{\left(c,\sigma,\Xi\right)},\tag{\ensuremath{\ast}}
\]
with the implied constants independent of the choice of $N\in\N$,
$i_{1},\dots,i_{N}\in I_{j_{0}}\cap I^{\left(r\right)}$, $\varepsilon>0$
and $c,\sigma,\Xi$. We will now ``arbitrage'' this estimate to
obtain the desired claims.

\medskip{}

First, we note that Lemma~\ref{lem:AsymptoticTranslationNormEstimate}
shows
\begin{align*}
E_{N}^{\left(c,\sigma,\Xi\right)}=\vphantom{\sum^{N}}\left\Vert \smash{\sum_{n=1}^{N}}\vphantom{\sum}L_{-\xi_{n}}\left(\sigma_{n}c_{n}\cdot\Fourier^{-1}\gamma_{n}\right)\right\Vert _{L^{p}} & \xrightarrow[\min_{i\neq j}\left|\xi_{i}-\xi_{j}\right|\to\infty]{}\left\Vert \left(\left\Vert \sigma_{n}c_{n}\cdot\Fourier^{-1}\gamma_{n}\right\Vert _{L^{p}}\right)_{n\in\underline{N}}\right\Vert _{\ell^{p}}\\
 & =\left\Vert \Fourier^{-1}\gamma\right\Vert _{L^{p}}\cdot\varepsilon^{\dimension\left(1-\frac{1}{p}\right)}\cdot\left\Vert \left(c_{n}\right)_{n\in\underline{N}}\right\Vert _{\ell^{p}}\;.
\end{align*}
In connection with equation~$\left(\ast\right)$, we thus see (with
``p.d.\@'' abbreviating ``pairwise distinct'') that
\[
\left\Vert \left(c_{n}\right)_{n\in\underline{N}}\right\Vert _{\ell^{p}}\!\!\asymp\!\left\Vert \left(c_{i}\right)_{i\in\left\{ i_{1},\dots,i_{N}\right\} }\right\Vert _{Y|_{\left\{ i_{1},\dots,i_{N}\right\} }}\forall\,N\!\in\N,\,c_{1},\dots,c_{N}\!\in\Compl\text{ and p.d. }i_{1},\dots,i_{N}\!\in I_{j_{0}}\cap I^{\left(r\right)}.\tag{\ensuremath{\lozenge}}
\]

Now, our first ``real'' goal is to show $p<\infty$. To see this,
we assume towards a contradiction that $p=\infty$. Then, we choose
$c_{1}=\dots=c_{N}=1$, $\Xi=0$ and $\sigma_{1}=\dots=\sigma_{N}=1$.
Furthermore, we recall that if $f\in L^{1}\left(\R^{d}\right)$ is
nonnegative, then $\widehat{f}\left(0\right)=\left\Vert f\right\Vert _{L^{1}}$;
by continuity of $\widehat{f}$, this entails $\left\Vert \smash{\widehat{f}}\,\right\Vert _{L^{\infty}}\geq\left\Vert f\right\Vert _{L^{1}}=\int f\d x$.
Since we have $\sum_{n=1}^{N}\gamma_{n}\geq0$ and $p=\infty$, this
implies
\[
E_{N}^{\left(c,\sigma,\Xi\right)}=\left\Vert \Fourier^{-1}\left(\,\smash{\sum_{n=1}^{N}}\vphantom{\sum}\gamma_{n}\,\right)\right\Vert _{L^{p}}\geq\int_{\R^{\dimension}}\sum_{n=1}^{N}\gamma_{n}\left(x\right)\d x=N\cdot\varepsilon^{\dimension}\cdot\left\Vert \gamma\right\Vert _{L^{1}}=\left\Vert \gamma\right\Vert _{L^{1}}\cdot\varepsilon^{\dimension\left(1-\frac{1}{p}\right)}\cdot N.
\]
In view of equations $\left(\ast\right)$ and $\left(\lozenge\right)$
and because of $I_{0}=\left\{ i_{1},\dots,i_{N}\right\} $ and $c_{1}=\dots=c_{N}=1$,
this yields
\[
\left\Vert \gamma\right\Vert _{L^{1}}\cdot N\leq\varepsilon^{\dimension\left(\frac{1}{p}-1\right)}\cdot E_{N}^{\left(c,\sigma,\Xi\right)}\asymp\left\Vert \Fourier^{-1}\gamma\right\Vert _{L^{p}}\cdot\left\Vert \left(1\right)_{i\in I_{0}}\right\Vert _{Y|_{I_{0}}}\!\!\smash{\overset{\text{eq. }\left(\lozenge\right)}{\asymp}}\:\left\Vert \Fourier^{-1}\gamma\right\Vert _{L^{p}}\cdot\left\Vert \left(1\right)_{n\in\underline{N}}\right\Vert _{\ell^{p}}=\left\Vert \Fourier^{-1}\gamma\right\Vert _{L^{\infty}},
\]
where the last step used again that $p=\infty$. Since the implied
constants in the preceding equation are independent of $N\in\N$,
we obtain the desired contradiction. Therefore, $p<\infty$.

Now, we want to show $p=2$. To this end, we choose $\Xi=0$, and
we choose $\sigma=\left(\sigma_{1},\dots,\sigma_{N}\right)$ to be
uniformly distributed in $\left\{ \pm1\right\} ^{N}$. The expectation
with respect to $\sigma$ is simply denoted by $\mathbb{E}$. Using
Khintchine's inequality (Theorem~\ref{thm:KhintchineInequality}),
we then get because of $\Xi=0$ that
\begin{align*}
\mathbb{E}\left(E_{N}^{\left(c,\sigma,\Xi\right)}\right)^{p} & =\vphantom{\sum_{n}^{N}}\mathbb{E}\left\Vert \smash{\sum_{n=1}^{N}}\vphantom{\sum}\sigma_{n}c_{n}\cdot\Fourier^{-1}\gamma_{n}\right\Vert _{L^{p}}^{p}=\mathbb{E}\left[\varepsilon^{\dimension p}\cdot\left\Vert \smash{\sum_{n=1}^{N}}\vphantom{\sum}\sigma_{n}c_{n}\cdot M_{\omega_{n}}\left(\left[\Fourier^{-1}\gamma\right]\circ\varepsilon\identity\right)\right\Vert _{L^{p}}^{p}\right]\\
 & =\vphantom{\sum_{n}^{N}}\varepsilon^{\dimension p}\cdot\int_{\R^{d}}\mathbb{E}\left|\smash{\sum_{n=1}^{N}}\vphantom{\sum}\sigma_{n}\cdot c_{n}\cdot e^{2\pi i\left\langle \omega_{n},x\right\rangle }\cdot\left(\Fourier^{-1}\gamma\right)\left(\varepsilon x\right)\right|^{p}\d x\\
\left({\scriptstyle \text{Khintchine}}\right) & \asymp\varepsilon^{\dimension p}\cdot\vphantom{\sum_{n}^{N}}\int_{\R^{d}}\left(\,\smash{\sum_{n=1}^{N}}\vphantom{\sum}\left|c_{n}\cdot e^{2\pi i\left\langle \omega_{n},x\right\rangle }\cdot\left(\Fourier^{-1}\gamma\right)\left(\varepsilon x\right)\right|^{2}\,\right)^{p/2}\d x\\
 & =\varepsilon^{\dimension p}\cdot\vphantom{\sum^{N}}\left(\smash{\sum_{n=1}^{N}}\vphantom{\sum}\left|c_{n}\right|^{2}\right)^{p/2}\cdot\int_{\R^{d}}\left|\left(\Fourier^{-1}\gamma\right)\left(\varepsilon x\right)\right|^{p}\d x=\left[\varepsilon^{\dimension\left(1-\frac{1}{p}\right)}\cdot\left\Vert \left(c_{n}\right)_{n\in\underline{N}}\right\Vert _{\ell^{2}}\cdot\left\Vert \Fourier^{-1}\gamma\right\Vert _{L^{p}}\right]^{p},
\end{align*}
where the implied constant only depends on $p\in\left(0,\infty\right)$.
In combination with equations $\left(\ast\right)$ and $\left(\lozenge\right)$,
we thus see $\left\Vert \mybullet\right\Vert _{\ell^{2}\left(\underline{N}\right)}\asymp\left\Vert \mybullet\right\Vert _{\ell^{p}\left(\underline{N}\right)}$,
where the implied constant is independent of the choice of $N\in\N$.
Hence, necessarily $p_{1}=p_{2}=p=2$. We have thus established the
first two claims.

\smallskip{}

Now, let us additionally assume that $Y=\ell_{w}^{q_{1}}\left(I\right)$
with a $\CalQ$-moderate weight $w=\left(w_{i}\right)_{i\in I}$.
From the second claim of the current theorem, we get $w_{i}=\left\Vert \delta_{i}\right\Vert _{Y}\asymp\left\Vert \delta_{j_{0}}\right\Vert _{Z}$
for all $i\in I$ with $Q_{i}\cap P_{j_{0}}\neq\emptyset$, and thus
in particular for every $i=i_{\ell}$ with $i_{\ell}\in I_{j_{0}}\cap I^{\left(r\right)}$.
In combination with equation~$\left(\lozenge\right)$ (choosing $c_{1}=\dots=c_{N}=1$),
we thus get
\begin{align*}
\left\Vert \delta_{j_{0}}\right\Vert _{Z}\cdot N^{1/q_{1}} & =\left\Vert \delta_{j_{0}}\right\Vert _{Z}\cdot\left\Vert \left(1\right)_{n\in\underline{N}}\right\Vert _{\ell^{q_{1}}}\\
\left({\scriptstyle \text{since }Y=\ell_{w}^{q_{1}}\left(I\right)}\right) & \asymp\vphantom{\left\Vert \left(c_{i}\right)_{i\in\left\{ i_{1},\dots,i_{N}\right\} }\right\Vert _{Y|_{\left\{ i_{1},\dots,i_{N}\right\} }}}\smash{\left\Vert \left(c_{i}\right)_{i\in\left\{ i_{1},\dots,i_{N}\right\} }\right\Vert _{Y|_{\left\{ i_{1},\dots,i_{N}\right\} }}}\\
\left({\scriptstyle \text{by equation }\left(\lozenge\right)}\right) & \asymp\left\Vert \left(c_{n}\right)_{n\in\underline{N}}\right\Vert _{\ell^{p}}=N^{1/p},
\end{align*}
where the implied constant is independent of $N\in\N$. Hence, $q_{1}=p=2$.

\smallskip{}

It remains to prove the last claim. Hence, assume $Y=\ell_{w}^{q_{1}}\left(I\right)$
and $Z=\ell_{v}^{q_{2}}\left(J\right)$. From the second claim of
the current theorem, we get $w_{i}=\left\Vert \delta_{i}\right\Vert _{Y}\asymp\left\Vert \delta_{j}\right\Vert _{Z}=v_{j}$
in case of $Q_{i}\cap P_{j}\neq\emptyset$.

Let us first consider the case where $\overline{\CalO\cap\CalO'}\nsubseteq\CalO\cup\CalO'$,
i.e., there is some $\eta_{0}\in\overline{\CalO\cap\CalO'}$ with
$\eta_{0}\notin\CalO\cup\CalO'$. Hence, there is a sequence $\left(\eta_{n}\right)_{n\in\N}$
in $\CalO\cap\CalO'$ with $\eta_{n}\to\eta_{0}$. For each $n\in\N$,
choose $i_{n}\in I$ and $j_{n}\in J$ with $\eta_{n}\in Q_{i_{n}}\cap P_{j_{n}}$.
If we had $i_{n}=i$ for some fixed $i\in I$ and infinitely many
$n\in\N$, we would get $\eta_{0}\in\overline{Q_{i}}\subset\CalO$,
see Lemma~\ref{lem:PartitionCoveringNecessary}. But since $\eta_{0}\notin\CalO$,
we see that $i_{n}=i$ can only hold for finitely many $n\in\N$,
for each $i\in I$. Thus, by passing to a subsequence, we can assume
$i_{n}\neq i_{m}$ for $n\neq m$. Completely similarly, we can also
assume $j_{n}\neq j_{m}$ for $n\neq m$.

If instead of $\overline{\CalO\cap\CalO'}\nsubseteq\CalO\cup\CalO'$
we have that $\CalO\cap\CalO'$ is unbounded, we can choose a sequence
$\left(\eta_{n}\right)_{n\in\N}$ in $\CalO\cap\CalO'$ with $\left|\eta_{n}\right|\to\infty$.
As above, let $\eta_{n}\in Q_{i_{n}}\cap P_{j_{n}}$ for each $n\in\N$.
Now, since each of the sets $Q_{i}$ and $P_{j}$ is bounded, we get
as above that $i=i_{n}$ can only hold for finitely many $n\in\N$,
for each $i\in I$. Thus, as above we can assume $i_{n}\neq i_{m}$
and $j_{n}\neq j_{m}$ for $n\neq m$, by passing to a subsequence.

\medskip{}

Next, there are $\left(r,s\right)\in\underline{R_{1}}\times\underline{R_{2}}$
such that $\left(i_{n},j_{n}\right)\in I^{\left(r\right)}\times J^{\left(s\right)}$
for infinitely many $n\in\N$. By passing to a subsequence, we can
(and will) assume that this holds for all $n\in\N$. In particular,
we thus get $Q_{i_{n}}^{\ast}\cap Q_{i_{m}}^{\ast}=\emptyset=P_{j_{n}}^{\ast}\cap P_{j_{\ell}}^{\ast}$
if $n\neq\ell$.

Let $N\in\N$ be arbitrary and choose $\varepsilon=\varepsilon\left(N\right)>0$
with $B_{\varepsilon}\left(\eta_{n}\right)\subset Q_{i_{n}}\cap P_{j_{n}}$
for all $n\in\underline{N}$, and set $\gamma_{i_{n}}:=\gamma_{j_{n}}:=\gamma_{n}:=L_{\eta_{n}}\left(\gamma\circ\varepsilon^{-1}\identity\right)\in\TestFunctionSpace{Q_{i_{n}}\cap P_{j_{n}}}$.
Furthermore, let $u_{n}:=u_{j_{n}}:=u_{i_{n}}:=w_{i_{n}}$ for $n\in\underline{N}$
and note because of $\eta_{n}\in Q_{i_{n}}\cap P_{j_{n}}\neq\emptyset$
that $u_{j_{n}}=u_{n}=w_{i_{n}}\asymp v_{j_{n}}$. By recalling $p_{1}=p_{2}=q_{1}=2$
and by combining Lemma~\ref{lem:EasyNormEquivalenceFineCovering}
and Corollary~\ref{cor:EasyNormEquivalenceFineLowerBound} (one time
with $I_{0}=\left\{ i_{1},\dots,i_{N}\right\} $ and one time with
$\CalP$ instead of $\CalQ$ and with $I_{0}=\left\{ j_{1},\dots,j_{N}\right\} $),
we thus get
\begin{align*}
\varepsilon^{\dimension/2}\cdot\left\Vert \Fourier^{-1}\gamma\right\Vert _{L^{2}}\cdot N^{1/2}=\left\Vert \left(u_{i}^{-1}\cdot\left\Vert \Fourier^{-1}\gamma_{i}\right\Vert _{L^{2}}\right)_{i\in\left\{ i_{1},\dots,i_{N}\right\} }\right\Vert _{\ell_{w}^{q_{1}}} & \asymp\vphantom{\sum_{n=1}^{N}}\left\Vert \smash{\sum_{n=1}^{N}}\vphantom{\sum}u_{n}^{-1}\cdot\gamma_{n}\right\Vert _{\FourierDecompSp{\CalQ}{p_{1}}{\ell_{w}^{q_{1}}}}\\
\left({\scriptstyle \text{by assumption of the current theorem}}\right) & \asymp\vphantom{\sum^{N}}\left\Vert \smash{\sum_{n=1}^{N}}\vphantom{\sum}u_{n}^{-1}\cdot\gamma_{n}\right\Vert _{\FourierDecompSp{\CalP}{p_{2}}{\ell_{v}^{q_{2}}}}\\
 & \asymp\left\Vert \left(u_{j}^{-1}\cdot\left\Vert \Fourier^{-1}\gamma_{j}\right\Vert _{L^{2}}\right)_{j\in\left\{ j_{1},\dots,j_{N}\right\} }\right\Vert _{\ell_{v}^{q_{2}}}\\
\left({\scriptstyle \text{since }u_{j_{n}}=w_{i_{n}}\asymp v_{j_{n}}}\right) & \asymp\varepsilon^{\dimension/2}\cdot\left\Vert \Fourier^{-1}\gamma\right\Vert _{L^{2}}\cdot N^{1/q_{2}},
\end{align*}
where the implied constants are independent of $N\in\N$ and of $\varepsilon>0$.
Hence, $q_{2}=2$, as desired.
\end{proof}
We close this subsection by showing that the necessary criteria from
Theorem~\ref{thm:NecessaryCriterionForCoincidenceOfDecompositionSpaces}
for the equality of two decomposition spaces—with weighted $\ell^{q}$
spaces as global components—are sufficient for the coincidence of
the two decomposition spaces. At least for $p\in\left[1,\infty\right]$,
this is true as stated. In the quasi-Banach regime $p\in\left(0,1\right)$,
we will need to impose certain additional restrictions.

As our starting point, note that Theorem~\ref{thm:NecessaryCriterionForCoincidenceOfDecompositionSpaces}
only shows that $\CalQ,\CalP$ are weakly equivalent as long as $\left(p_{1},q_{1}\right)\neq\left(2,2\right)$.
This is natural, since all (Fourier side) decomposition spaces $\FourierDecompSp{\CalQ}p{\ell_{w}^{q}}$
with $p=q=2$ are just certain weighted $L^{2}$-spaces:
\begin{lem}
\label{lem:DecompositionSpacesHilbertCase}Let $\emptyset\neq\CalO\subset\R^{\dimension}$
be open, and let $\CalQ=\left(Q_{i}\right)_{i\in I}$ be an $L^{2}$-decomposition
covering of $\CalO$. Finally, let $w=\left(w_{i}\right)_{i\in I}$
be $\CalQ$-moderate. Then there is a measurable weight 
\[
w_{0}:\CalO\to\left(0,\infty\right)\qquad\text{ with }\qquad w_{0}\left(\xi\right)\asymp w_{i}\text{ for all }\xi\in Q_{i}\text{ and arbitrary }i\in I,
\]
where the implied constant is independent of $i,\xi$. For any such
weight $w_{0}$, both $w_{0}$ and $w_{0}^{-1}$ are locally bounded
on $\CalO$.

Furthermore, we have
\begin{equation}
\FourierDecompSp{\CalQ}2{\ell_{w}^{2}}=L_{w_{0}}^{2}\left(\CalO\right)\qquad\text{with equivalent norms},\label{eq:DecompositionSpaceHilbertCase}
\end{equation}
for the weighted $L^{2}$ space $L_{w_{0}}^{2}\left(\CalO\right):=\left\{ f:\CalO\to\Compl\with w_{0}\cdot f\in L^{2}\left(\CalO\right)\right\} $,
with norm $\left\Vert f\right\Vert _{L_{w_{0}}^{2}}:=\left\Vert w_{0}\cdot f\right\Vert _{L^{2}}$.

In particular, if $\CalP=\left(P_{j}\right)_{j\in J}$ is another
$L^{2}$-decomposition covering of $\CalO$ and if $v=\left(v_{j}\right)_{j\in J}$
is $\CalP$-moderate with $w_{i}\asymp v_{j}$ in case of $Q_{i}\cap P_{j}\neq\emptyset$,
then
\[
\FourierDecompSp{\CalQ}2{\ell_{w}^{2}}=\FourierDecompSp{\CalP}2{\ell_{v}^{2}}\quad\text{with equivalent norms}.\qedhere
\]
\end{lem}

\begin{proof}
Let $\Phi=\left(\varphi_{i}\right)_{i\in I}$ be an $L^{2}$-BAPU
for $\CalO$. We start with two technical results: First, we show
that $I$ is necessarily countable. Then, we show that one can choose
a weight $w_{0}$ as in the statement of the lemma, but with the additional
property that $w_{0}$ is smooth.

To see that $I$ is countable, note as a consequence of Lemma~\ref{lem:PartitionCoveringNecessary}
that the family of interiors $\left(Q_{i}^{\circ}\right)_{i\in I}$
covers $\CalO$. Since $\CalO\subset\R^{\dimension}$ is second-countable,
there is thus a sequence $\left(i_{n}\right)_{n\in\N}$ with $\CalO=\bigcup_{n\in\N}Q_{i_{n}}^{\circ}$.
But this implies that $I=\bigcup_{n\in\N}i_{n}^{\ast}$ is countable
as a union of countably many finite sets. To see that the last identity
is valid, let $i\in I$ be arbitrary. Since $\CalQ$ is an admissible
covering, we have $Q_{i}\neq\emptyset$; see Definition~\ref{defn:AdmissibleCoveringModerateWeight}.
Choosing any $\xi\in Q_{i}\subset\CalO$, we see that there is some
$n\in\N$ with $\xi\in Q_{i_{n}}^{\circ}$, so that $i\in i_{n}^{\ast}$,
as claimed.

Next, for constructing a smooth weight $w_{0}$ as in the statement
of the lemma, we first construct a smooth, \emph{nonnegative} partition
of unity $\left(\gamma_{i}\right)_{i\in I}$ subordinate to $\CalQ$.
We cannot simply take $\gamma_{i}=\varphi_{i}$, since the elements
of a BAPU are not required to be nonnegative. To construct $\left(\gamma_{i}\right)_{i\in I}$,
set $\psi_{i}:=\left(\mathrm{Re}\,\varphi_{i}\right)^{2}$ for $i\in I$,
and note $\supp\psi_{i}\subset\supp\varphi_{i}$. Lemma~\ref{lem:PartitionCoveringNecessary}
shows that the family $\left(\supp\varphi_{i}\right)_{i\in I}$ is
locally finite in $\CalO$. Thus, the function $\Psi:\CalO\to\left\lceil 0,\infty\right),\xi\mapsto\sum_{i\in I}\psi_{i}$
is well-defined and smooth as a locally finite sum of smooth, nonnegative
functions. In fact, $\Psi$ is positive on $\CalO$: For any $\xi\in\CalO$,
we have $\sum_{i\in I}\varphi_{i}\left(\xi\right)=1$, so that there
is some $i\in I$ with $\mathrm{Re}\,\varphi_{i}\left(\xi\right)\neq0$,
and thus $\psi_{i}\left(\xi\right)>0$, whence $\Psi\left(\xi\right)>0$.
Thus, if we set $\gamma_{i}:=\psi_{i}/\Psi$, then it is not hard
to see that $\gamma_{i}\in C_{c}^{\infty}\left(\CalO\right)$ is nonnegative
with $\gamma_{i}\left(\xi\right)=0$ for $\xi\in\R^{\dimension}\setminus Q_{i}$,
and with $\sum_{i\in I}\gamma_{i}\equiv1$ on $\CalO$. Furthermore,
$\supp\gamma_{i}\subset\supp\varphi_{i}$, so that $\left(\gamma_{i}\right)_{i\in I}$
is a smooth locally finite partition of unity on $\CalO$.

Therefore, the function
\begin{equation}
w_{0}:\CalO\to\left[0,\infty\right),\xi\mapsto\sum_{i\in I}w_{i}\cdot\gamma_{i}\left(\xi\right)\label{eq:ContinuousWeightForHilbertSpaceCase}
\end{equation}
is smooth. But for $i,\ell\in I$ and $\xi\in Q_{\ell}$, we have
$\gamma_{i}\left(\xi\right)=0$ unless $i\in\ell^{\ast}$. But in
case of $i\in\ell^{\ast}$, we have $C_{w,\CalQ}^{-1}\cdot w_{\ell}\leq w_{i}\leq C_{w,\CalQ}\cdot w_{\ell}$.
Thus,
\[
0<C_{w,\CalQ}^{-1}\cdot w_{\ell}=C_{w,\CalQ}^{-1}\cdot w_{\ell}\cdot\sum_{i\in I}\gamma_{i}\left(\xi\right)\leq w_{0}\left(\xi\right)\leq C_{w,\CalQ}\cdot w_{\ell}\cdot\sum_{i\in I}\gamma_{i}\left(\xi\right)=C_{w,\CalQ}\cdot w_{\ell}.
\]
We have thus shown $w_{0}\left(\xi\right)\asymp w_{\ell}$ for $\xi\in Q_{\ell}$
and $\ell\in I$, where the implied constant only depends on $C_{w,\CalQ}$.
In particular, $w_{0}:\CalO\to\left(0,\infty\right)$, so that $w_{0}$
is as desired.

Next, we show that any weight $w_{0}$ as in the statement of the
lemma is locally bounded on $\CalO$. Indeed, let $w_{0}:\CalO\to\left(0,\infty\right)$
be measurable with $C^{-1}\cdot w_{i}\leq w_{0}\left(\xi\right)\leq C\cdot w_{i}$
for all $i\in I$ and all $\xi\in Q_{i}$, with a fixed constant $C\geq1$.
Then, for any $\xi\in\CalO$, Lemma~\ref{lem:PartitionCoveringNecessary}
shows that there is some $i\in I$ with $\xi\in Q_{i}^{\circ}$. But
on this set, we have $w_{0}\left(\eta\right)\leq C\cdot w_{i}$ and
$\left[w_{0}\left(\eta\right)\right]^{-1}\leq C\cdot w_{i}^{-1}$,
so that $w_{0}$ and $w_{0}^{-1}$ are bounded on the neighborhood
$Q_{i}^{\circ}$ of $\xi$. Hence, $w_{0}$ and $w_{0}^{-1}$ are
locally bounded on $\CalO$. Thus, it remains to prove equation~(\ref{eq:DecompositionSpaceHilbertCase}).

\medskip{}

For the proof of ``$\hookleftarrow$'' in equation~(\ref{eq:DecompositionSpaceHilbertCase}),
let $f\in L_{w_{0}}^{2}\left(\CalO\right)$ be arbitrary. Since $w_{0}^{-1}$
is locally bounded, $w_{0}$ is locally bounded from below, so that
we get $L_{w_{0}}^{2}\left(\CalO\right)\hookrightarrow L_{{\rm loc}}^{2}\left(\CalO\right)\hookrightarrow L_{{\rm loc}}^{1}\left(\CalO\right)\hookrightarrow\DistributionSpace{\CalO}$.

Now, by continuity of each of the $\varphi_{i}$, we have
\begin{equation}
\left|\varphi_{i}\left(\xi\right)\right|\leq\left\Vert \varphi_{i}\right\Vert _{L^{\infty}}=\left\Vert \Fourier\Fourier^{-1}\varphi_{i}\right\Vert _{L^{\infty}}\leq\left\Vert \Fourier^{-1}\varphi\right\Vert _{L^{1}}\leq C_{\CalQ,\Phi,2}\qquad\forall\,i\in I\text{ and }\xi\in\R^{\dimension}\,.\label{eq:BAPUIsPointwiseBounded}
\end{equation}
Furthermore, Plancherel's theorem implies for any $i\in I$ that
\begin{align*}
\left\Vert \Fourier^{-1}\left(\varphi_{i}\cdot f\right)\right\Vert _{L^{2}}^{2} & =\left\Vert \varphi_{i}\cdot f\right\Vert _{L^{2}}^{2}\\
\left({\scriptstyle \text{since }\supp\varphi_{i}\subset\CalO}\right) & =w_{i}^{-2}\cdot\int_{\CalO}\left|w_{i}\cdot\varphi_{i}\left(\xi\right)\cdot f\left(\xi\right)\right|^{2}\,\d\xi\\
\left({\scriptstyle \text{since }\varphi_{i}\equiv0\text{ on }Q_{i}^{c}\text{ and }w_{i}\leq C\cdot w_{0}\left(\xi\right)\text{ on }Q_{i}}\right) & \leq C^{2}\cdot w_{i}^{-2}\cdot\int_{\CalO}\left|\varphi_{i}\left(\xi\right)\cdot\left(w_{0}\cdot f\right)\left(\xi\right)\right|^{2}\,\d\xi.
\end{align*}
Multiplying with $w_{i}$ and summing over $i\in I$ yields
\begin{align*}
\left\Vert f\right\Vert _{\FourierDecompSp{\CalQ}2{\ell_{w}^{2}}}^{2}=\sum_{i\in I}\left(w_{i}\cdot\left\Vert \Fourier^{-1}\left(\varphi_{i}\cdot f\right)\right\Vert _{L^{2}}\right)^{2} & \leq C^{2}\cdot\sum_{i\in I}\int_{\CalO}\left|\varphi_{i}\left(\xi\right)\cdot\left(w_{0}\cdot f\right)\left(\xi\right)\right|^{2}\,\d\xi\\
\left({\scriptstyle \text{monotone convergence, }I\text{ countable}}\right) & =C^{2}\cdot\int_{\CalO}\left(\,\smash{\sum_{i\in I}}\,\vphantom{\sum}\left|\varphi_{i}\left(\xi\right)\right|^{2}\,\right)\vphantom{\sum_{i\in I}}\cdot\left|\left(w_{0}\cdot f\right)\left(\xi\right)\right|^{2}\,\d\xi.
\end{align*}
But equation~(\ref{eq:BAPUIsPointwiseBounded}) yields for $\xi\in\CalO$
that
\[
\sum_{i\in I}\left|\varphi_{i}\left(\xi\right)\right|^{2}\leq C_{\CalQ,\Phi,2}^{2}\cdot\sum_{i\in I}\Indicator_{Q_{i}}\left(\xi\right)\leq N_{\CalQ}C_{\CalQ,\Phi,2}^{2}\,\,,
\]
where the last step used that $\CalQ$ is admissible. All in all,
we have thus shown
\begin{align*}
\left\Vert f\right\Vert _{\FourierDecompSp{\CalQ}2{\ell_{w}^{2}}}^{2} & \leq C^{2}\cdot\int_{\CalO}\left(\,\smash{\sum_{i\in I}}\,\vphantom{\sum}\left|\varphi_{i}\left(\xi\right)\right|^{2}\,\right)\vphantom{\sum_{i\in I}}\cdot\left|\left(w_{0}\cdot f\right)\left(\xi\right)\right|^{2}\,\d\xi\\
 & \leq N_{\CalQ}\cdot\left(CC_{\CalQ,\Phi,2}\right)^{2}\cdot\int_{\CalO}\left|\left(w_{0}\cdot f\right)\left(\xi\right)\right|^{2}\,\d\xi=N_{\CalQ}\cdot\left(CC_{\CalQ,\Phi,2}\right)^{2}\cdot\left\Vert f\right\Vert _{L_{w_{0}}^{2}\left(\CalO\right)}^{2}<\infty,
\end{align*}
that is, $L_{w_{0}}^{2}\left(\CalO\right)\hookrightarrow\FourierDecompSp{\CalQ}2{\ell_{w}^{2}}$.

\medskip{}

For the reverse inclusion, let $f\in\FourierDecompSp{\CalQ}2{\ell_{w}^{2}}\subset\DistributionSpace{\CalO}$.
Note that $L_{w_{0}}^{2}\left(\CalO\right)=L_{w_{1}}^{2}\left(\CalO\right)$
if $w_{0},w_{1}$ are measurable and equivalent, i.e., if $w_{0}\asymp w_{1}$.
Furthermore, any two weights $w_{0},w_{1}$ satisfying $w_{\ell}\left(\xi\right)\asymp w_{i}$
for $\xi\in Q_{i}$, $\ell\in\left\{ 0,1\right\} $ and $i\in I$
are equivalent. Hence, we can assume without loss of generality that
$w_{0}$ is as in equation~(\ref{eq:ContinuousWeightForHilbertSpaceCase}).
As seen above, this ensures $w_{0}\in C^{\infty}\left(\CalO\right)$.

Thus, for $g\in\TestFunctionSpace{\CalO}$, we have $w_{0}\cdot g\in\TestFunctionSpace{\CalO}$
as well. Now, note
\begin{align*}
\left|\left\langle f,\,w_{0}\cdot g\right\rangle _{\CalD'}\right| & =\left|\,\smash{\sum_{i\in I}}\,\vphantom{\sum}\left\langle \varphi_{i}\cdot f,\,w_{0}\cdot g\right\rangle _{\Schwartz'}\,\right|\vphantom{\sum_{i\in I}}\\
\left({\scriptstyle \text{since }\varphi_{i}^{\ast}\equiv1\text{ on }Q_{i}\text{ and }\varphi_{i}\equiv0\text{ on }Q_{i}^{c}}\right) & =\left|\,\smash{\sum_{i\in I}}\,\vphantom{\sum}\left\langle \varphi_{i}\cdot f,\,\varphi_{i}^{\ast}\cdot w_{0}\cdot g\right\rangle _{\Schwartz'}\,\right|\vphantom{\sum_{i\in I}}\\
 & =\left|\,\smash{\sum_{i\in I}}\,\vphantom{\sum}\left\langle \Fourier^{-1}\left(\varphi_{i}\cdot f\right),\,\Fourier\left(\varphi_{i}^{\ast}\cdot w_{0}\cdot g\right)\right\rangle _{\Schwartz'}\,\right|\vphantom{\sum_{i\in I}}\\
 & \leq\sum_{i\in I}\left[w_{i}\cdot\left\Vert \Fourier^{-1}\left(\varphi_{i}\cdot f\right)\right\Vert _{L^{2}}\cdot w_{i}^{-1}\cdot\left\Vert \Fourier\left(\varphi_{i}^{\ast}\cdot w_{0}\cdot g\right)\right\Vert _{L^{2}}\right]\\
\left({\scriptstyle \text{Cauchy-Schwarz and Plancherel}}\right) & \leq\left[\,\smash{\sum_{i\in I}}\,\vphantom{\sum}\left(w_{i}\cdot\left\Vert \Fourier^{-1}\left(\varphi_{i}\cdot f\right)\right\Vert _{L^{2}}\right)^{2}\,\right]^{1/2}\!\!\cdot\left[\,\smash{\sum_{i\in I}}\,\vphantom{\sum}\left(w_{i}^{-1}\cdot\left\Vert \varphi_{i}^{\ast}\cdot w_{0}\cdot g\right\Vert _{L^{2}}\right)^{2}\,\right]^{1/2}\vphantom{\sum_{i\in I}}\!\!.
\end{align*}
The first factor in the last line is exactly $\left\Vert f\right\Vert _{\FourierDecompSp{\CalQ}2{\ell_{w}^{2}}}$.

Now, let us consider the second factor: First, note
\begin{align*}
\sum_{i\in I}\left(w_{i}^{-1}\cdot\left\Vert \varphi_{i}^{\ast}\cdot w_{0}\cdot g\right\Vert _{L^{2}}\right)^{2} & =\sum_{i\in I}\left(w_{i}^{-2}\cdot\int_{\R^{\dimension}}\left|\varphi_{i}^{\ast}\left(\xi\right)\right|^{2}\cdot\left|\left(w_{0}\cdot g\right)\left(\xi\right)\right|^{2}\,\d\xi\right)\\
\left({\scriptstyle \text{since }g\in\TestFunctionSpace{\CalO},\text{ and by mon. conv., since }I\text{ countable}}\right) & =\int_{\CalO}\left|g\left(\xi\right)\right|^{2}\cdot\sum_{i\in I}\left|w_{i}^{-1}\cdot\varphi_{i}^{\ast}\left(\xi\right)\cdot w_{0}\left(\xi\right)\right|^{2}\,\d\xi.
\end{align*}
But as a consequence of equation~(\ref{eq:BAPUIsPointwiseBounded}),
we have
\[
\left|\varphi_{i}^{\ast}\left(\xi\right)\right|\leq\sum_{\ell\in i^{\ast}}\left|\varphi_{\ell}\left(\xi\right)\right|\leq C_{\CalQ,\Phi,2}\cdot\left|i^{\ast}\right|\cdot\Indicator_{Q_{i}^{\ast}}\left(\xi\right)\leq C_{\CalQ,\Phi,2}\cdot N_{\CalQ}\cdot\Indicator_{Q_{i}^{\ast}}\left(\xi\right)\qquad\forall\,\,i\in I\text{ and }\xi\in\R^{\dimension}\,.
\]
Furthermore, for any $i\in I$ with $\varphi_{i}^{\ast}\left(\xi\right)\neq0$,
we have $\xi\in Q_{\ell}$ for some $\ell\in i^{\ast}$; by our assumptions
on $w_{0}$, this implies $w_{0}\left(\xi\right)\leq C\cdot w_{\ell}\leq CC_{w,\CalQ}\cdot w_{i}$.
Thus, $w_{i}^{-1}\cdot w_{0}\left(\xi\right)\leq CC_{w,\CalQ}$ in
case of $\varphi_{i}^{\ast}\left(\xi\right)\neq0$. Hence,
\begin{align*}
\sum_{i\in I}\left|w_{i}^{-1}\cdot\varphi_{i}^{\ast}\left(\xi\right)\cdot w_{0}\left(\xi\right)\right|^{2}\leq\left(CC_{w,\CalQ}\right)^{2}\cdot\sum_{i\in I}\left|\varphi_{i}^{\ast}\left(\xi\right)\right|^{2} & \leq\left(CC_{w,\CalQ}C_{\CalQ,\Phi,2}N_{\CalQ}\right)^{2}\cdot\sum_{i\in I}\Indicator_{Q_{i}^{\ast}}\left(\xi\right)\\
 & \leq\left(CC_{w,\CalQ}C_{\CalQ,\Phi,2}N_{\CalQ}\right)^{2}\cdot N_{\CalQ^{\ast}}\\
\left({\scriptstyle \text{see Lemma }\ref{lem:SemiStructuredClusterInvariant}}\right) & \leq\left(CC_{w,\CalQ}C_{\CalQ,\Phi,2}N_{\CalQ}\right)^{2}\cdot N_{\CalQ}^{3}.
\end{align*}

Putting everything together, we see
\[
\left|\left\langle f,\,w_{0}\cdot g\right\rangle _{\CalD'}\right|\leq\left\Vert f\right\Vert _{\FourierDecompSp{\CalQ}2{\ell_{w}^{2}}}\cdot CC_{w,\CalQ}C_{\CalQ,\Phi,2}N_{\CalQ}^{5/2}\cdot\sqrt{\int_{\CalO}\left|g\left(\xi\right)\right|^{2}\,\d\xi}=:C_{1}\cdot\left\Vert f\right\Vert _{\FourierDecompSp{\CalQ}2{\ell_{w}^{2}}}\cdot\left\Vert g\right\Vert _{L^{2}\left(\CalO\right)}.
\]
Since this holds for all $g\in\TestFunctionSpace{\CalO}$ and since
$\TestFunctionSpace{\CalO}\leq L^{2}\left(\CalO\right)$ is dense,
the Riesz representation theorem for Hilbert spaces yields some $h\in L^{2}\left(\CalO\right)\hookrightarrow\DistributionSpace{\CalO}$
with $\left\Vert h\right\Vert _{L^{2}}\leq C_{1}\cdot\left\Vert f\right\Vert _{\FourierDecompSp{\CalQ}2{\ell_{w}^{2}}}$
and with
\[
\left\langle f,\,w_{0}\cdot g\right\rangle _{\CalD'}=\left\langle h,\,g\right\rangle _{\CalD'}=\left\langle w_{0}^{-1}\cdot h,\,w_{0}\cdot g\right\rangle _{\CalD'}\qquad\forall\,\,g\in\TestFunctionSpace{\CalO}.
\]
Since $\TestFunctionSpace{\CalO}\to\TestFunctionSpace{\CalO},g\mapsto w_{0}\cdot g$
is a bijection (because $w_{0}:\CalO\to\left(0,\infty\right)$ is
smooth), this means $f=w_{0}^{-1}\cdot h$ as elements of $\DistributionSpace{\CalO}$.
Hence, $f=w_{0}^{-1}\cdot h\in L_{w_{0}}^{2}\left(\CalO\right)$ with
\[
\left\Vert f\right\Vert _{L_{w_{0}}^{2}}=\left\Vert w_{0}^{-1}h\right\Vert _{L_{w_{0}}^{2}}=\left\Vert h\right\Vert _{L^{2}}\leq C_{1}\cdot\left\Vert f\right\Vert _{\FourierDecompSp{\CalQ}2{\ell_{w}^{2}}}<\infty.
\]

\medskip{}

For the final statement of the lemma, simply note that if $w_{0},v_{0}:\CalO\to\left(0,\infty\right)$
are the associated continuous weights for the discrete weights $w,v$
and if $\xi\in\CalO$, then $\xi\in Q_{i}\cap P_{j}$ for suitable
$i\in I$ and $j\in J$. Hence, $w_{0}\left(\xi\right)\asymp w_{i}\asymp v_{j}\asymp v_{0}\left(\xi\right)$,
so that we get $w_{0}\asymp v_{0}$, and thus
\[
\FourierDecompSp{\CalQ}2{\ell_{w}^{2}}=L_{w_{0}}^{2}\left(\CalO\right)=L_{v_{0}}^{2}\left(\CalO\right)=\FourierDecompSp{\CalP}2{\ell_{v}^{2}}\quad\text{with equivalent norms.}\qedhere
\]
\end{proof}
Finally, we consider the case in which $\CalQ,\CalP$ are weakly equivalent.
Again, we only consider the case of weighted $\ell^{q}$ spaces as
the global component. The main reason for this is that it is difficult
to formulate what it means for two generic (solid) sequence spaces
$Y\subset\Compl^{I}$ and $Z\subset\Compl^{J}$—which live on different
sets $I\neq J$—to be ``equivalent''. In the present case of weighted
$\ell^{q}$ spaces, we simply require the exponents $q_{1},q_{2}$
of $\ell_{w}^{q_{1}}\left(I\right)$ and $\ell_{v}^{q_{2}}\left(J\right)$
to coincide, $q_{1}=q_{2}$, and furthermore require that $w_{i}\asymp v_{j}$
if $Q_{i}\cap P_{j}\neq\emptyset$.
\begin{lem}
\label{lem:EquivalentCoveringsYieldSameDecompositionSpaces}Let $p,q\in\left(0,\infty\right]$
and let $\emptyset\neq\CalO\subset\R^{\dimension}$ be open. Let $\CalQ=\left(Q_{i}\right)_{i\in I}$
and $\CalP=\left(P_{j}\right)_{j\in J}$ be two $L^{p}$-decomposition
coverings of $\CalO$. Let $w=\left(w_{i}\right)_{i\in I}$ be $\CalQ$-moderate
and let $v=\left(v_{j}\right)_{j\in J}$ be $\CalP$-moderate. Assume
that

\begin{itemize}[leftmargin=0.6cm]
\item There is a constant $C_{0}>0$ with
\begin{equation}
C_{0}^{-1}\cdot w_{i}\leq v_{j}\leq C_{0}\cdot w_{i}\qquad\forall\,i\in I\text{ and }j\in J\text{ with }Q_{i}\cap P_{j}\neq\emptyset\,.\label{eq:SameDecompositionSpacesWeightEquivalence}
\end{equation}
\item $\CalQ,\CalP$ are weakly equivalent.
\end{itemize}
Then we have $\FourierDecompSp{\CalQ}p{\ell_{w}^{q}}=\FourierDecompSp{\CalP}p{\ell_{v}^{q}}$
with equivalent quasi-norms, as long as $p\in\left[1,\infty\right]$.

\medskip{}

For $p\in\left(0,1\right)$, the same holds under the following additional
assumption: Because of $p\in\left(0,1\right)$, the coverings $\CalQ=\left(T_{i}Q_{i}'+b_{i}\right)_{i\in I}$
and $\CalP=\left(\smash{S_{j}P_{j}'+c_{j}}\right)_{j\in J}$ are semi-structured.
We assume that there is a constant $C_{1}>0$ such that one of the
following two conditions holds:

\begin{enumerate}[leftmargin=0.8cm]
\item For arbitrary $i\in I$ and $j\in J$ with $Q_{i}\cap P_{j}\neq\emptyset$,
we have
\[
\left\Vert \smash{T_{i}^{-1}}S_{j}\right\Vert \leq C_{1}\;\text{ and }\;\left\Vert \smash{S_{j}^{-1}}T_{i}\right\Vert \leq C_{1}.
\]
\item $\CalP$ is almost subordinate to $\CalQ$ and we have
\[
\left|\det\left(\smash{S_{j}^{-1}}T_{i}\right)\right|\leq C_{1}\qquad\forall\,i\in I\text{ and }j\in J\text{ with }Q_{i}\cap P_{j}\neq\emptyset\,.\qedhere
\]
\end{enumerate}
\end{lem}

\begin{rem*}

\begin{enumerate}[leftmargin=0.8cm]
\item The two kinds of assumptions for $p\in\left(0,1\right)$ seem to
be somewhat artificial. But there is a more natural condition which
implies the second condition. Indeed, assume that

\begin{enumerate}
\item $\CalQ,\CalP$ are (open) semi-structured coverings of $\CalO$,
\item $\CalQ$ is tight,
\item $\CalQ$ and $\CalP$ are weakly equivalent,
\item all sets in $\CalQ$ and $\CalP$ are connected.
\end{enumerate}
Then the sets in $\CalQ$ and $\CalP$ are in fact path-connected
(as connected, open subsets of $\R^{\dimension}$). Since the sets
in $\CalQ$ and $\CalP$ are open, Corollary~\ref{cor:WeakSubordinationImpliesSubordinationIfConnected}
shows that $\CalQ$ and $\CalP$ are equivalent, not just weakly equivalent.
In particular, $\CalP$ is almost subordinate to $\CalQ$, as required.
Furthermore, there is $k\in\N_{0}$ and for each $i\in I$ some $j_{i}\in J$
with $Q_{i}\subset P_{j_{i}}^{k\ast}$. Hence, for $i\in I$ and $j\in J$
with $Q_{i}\cap P_{j}\neq\emptyset$, Lemma~\ref{lem:SubordinatenessImpliesWeakSubordination}
yields $Q_{i}\subset P_{j}^{\left(2k+2\right)\ast}$. By tightness
of $\CalQ$, this implies
\begin{align*}
\left|\det T_{i}\right|\lesssim\lambda\left(Q_{i}\right)\leq\lambda\left(\smash{P_{j}^{\left(2k+2\right)\ast}}\right) & \leq\sum_{\ell\in j^{\left(2k+2\right)\ast}}\lambda\left(P_{\ell}\right)\lesssim\sum_{\ell\in j^{\left(2k+2\right)\ast}}\left|\det S_{\ell}\right|\\
\left({\scriptstyle \left|\smash{j^{\left(2k+2\right)\ast}}\right|\leq N_{\CalP}^{2k+2}\text{ and eq. }\eqref{eq:DeterminantIsModerate}}\right) & \lesssim\left|\det S_{j}\right|,
\end{align*}
where all implied constants are independent of $i,j$. Thus, $\left|\det\left(S_{j}^{-1}T_{i}\right)\right|=\left|\det T_{i}\right|/\left|\det S_{j}\right|\lesssim1$,
so that the second additional assumption from the theorem is indeed
satisfied.
\item The proof shows that an alternative assumption for the case $p\in\left(0,1\right)$
would be to assume 
\[
\lambda\left(\,\overline{Q_{i}}-\overline{P_{j}}\,\right)\lesssim\min\left\{ \left|\det T_{i}\right|,\left|\det S_{j}\right|\right\} \qquad\forall\,i\in I\text{ and }j\in J\text{ with }Q_{i}\cap P_{j}\neq\emptyset.
\]
But this estimate is usually very hard to check without establishing
one of the additional assumptions from the lemma above.\qedhere
\end{enumerate}
\end{rem*}
\begin{proof}
Let $\Phi=\left(\varphi_{i}\right)_{i\in I}$ and $\Psi=\left(\psi_{j}\right)_{j\in J}$
be $L^{p}$-BAPUs for $\CalQ$ and $\CalP$, respectively. As usual,
let $J_{i}=\left\{ j\in J\with P_{j}\cap Q_{i}\neq\emptyset\right\} $,
and define $I_{j}$ analogously for $i\in I$ and $j\in J$. We first
show $\varphi_{i}\,\psi_{J_{i}}=\varphi_{i}$. For $\xi\in\R^{\dimension}$
with $\varphi_{i}\left(\xi\right)=0$, this is clear. If otherwise
$\varphi_{i}\left(\xi\right)\neq0$, then $\xi\in Q_{i}\subset\CalO$,
and each $j\in J$ with $\psi_{j}\left(\xi\right)\neq0$ hence satisfies
$\xi\in Q_{i}\cap P_{j}\neq\emptyset$, i.e.\@ $j\in J_{i}$. Thus,
\[
1=\sum_{j\in J}\psi_{j}\left(\xi\right)=\sum_{j\in J_{i}}\psi_{j}\left(\xi\right)=\psi_{J_{i}}\left(\xi\right),
\]
which implies $\varphi_{i}\left(\xi\right)\,\psi_{J_{i}}\left(\xi\right)=\varphi_{i}\left(\xi\right)$
also in this case. Note that the identity $\sum_{j\in J}\psi_{j}\left(\xi\right)=1$
crucially used that $\CalQ,\CalP$ both cover the \emph{same} set
$\CalO$. Analogously, one can show $\psi_{j}\,\varphi_{I_{j}}=\psi_{j}$.

\medskip{}

We first prove the claim for $p\in\left[1,\infty\right]$. Here, it
suffices to show $\left\Vert \mybullet\right\Vert _{\FourierDecompSp{\CalQ}p{\ell_{w}^{q}}}\lesssim\left\Vert \mybullet\right\Vert _{\FourierDecompSp{\CalP}p{\ell_{v}^{q}}}$,
since the assumptions are symmetric in $\CalQ,\CalP$. Hence, let
$f\in\FourierDecompSp{\CalP}p{\ell_{v}^{q}}$ be arbitrary. As seen
above, we have $\varphi_{i}=\varphi_{i}\,\psi_{J_{i}}$. Together
with Young's inequality $L^{1}\ast L^{p}\hookrightarrow L^{p}$ and
with the triangle inequality for $L^{p}$, this yields
\begin{align}
w_{i}\cdot\left\Vert \Fourier^{-1}\left(\varphi_{i}\,f\right)\right\Vert _{L^{p}} & =w_{i}\cdot\left\Vert \Fourier^{-1}\left(\varphi_{i}\,\psi_{J_{i}}\,f\right)\right\Vert _{L^{p}}\nonumber \\
 & \leq\sum_{j\in J_{i}}\left[w_{i}\cdot\left\Vert \Fourier^{-1}\left(\varphi_{i}\,\psi_{j}\,f\right)\right\Vert _{L^{p}}\right]\nonumber \\
 & \leq C_{0}\cdot\sum_{j\in J_{i}}\left[v_{j}\cdot\left\Vert \Fourier^{-1}\varphi_{i}\right\Vert _{L^{1}}\cdot\left\Vert \Fourier^{-1}\left(\psi_{j}\,f\right)\right\Vert _{L^{p}}\right]\nonumber \\
 & \leq C_{0}C_{\CalQ,\Phi,p}\cdot\sum_{j\in J_{i}}\left[v_{j}\cdot\left\Vert \Fourier^{-1}\left(\psi_{j}\,f\right)\right\Vert _{L^{p}}\right]\qquad\forall\,i\in I\,.\label{eq:EquivalentCoveringsYieldSameDecompositionSpacesFundamentalBanach}
\end{align}
In case of $q=\infty$, this implies
\begin{align*}
\left\Vert f\right\Vert _{\FourierDecompSp{\CalQ}p{\ell_{w}^{q}}}=\sup_{i\in I}w_{i}\cdot\left\Vert \Fourier^{-1}\left(\varphi_{i}\,f\right)\right\Vert _{L^{p}} & \leq\sup_{i\in I}C_{0}C_{\CalQ,\Phi,p}\cdot\left|J_{i}\right|\cdot\left\Vert f\right\Vert _{\FourierDecompSp{\CalP}p{\ell_{v}^{q}}}\\
 & \leq C_{0}C_{\CalQ,\Phi,p}\cdot N\left(\CalQ,\CalP\right)\cdot\left\Vert f\right\Vert _{\FourierDecompSp{\CalP}p{\ell_{v}^{q}}}<\infty.
\end{align*}
For $q\in\left(0,\infty\right)$, we use the uniform bound $\left|J_{i}\right|\leq N\left(\CalQ,\CalP\right)=:N$
to obtain $C_{2}=C_{2}\left(N,q\right)>0$ with
\[
\vphantom{\sum_{j\in J_{i}}}\left[\,\smash{\sum_{j\in J_{i}}}\vphantom{\sum}v_{j}\cdot\left\Vert \Fourier^{-1}\left(\psi_{j}\,f\right)\right\Vert _{L^{p}}\,\right]^{q}\leq C_{2}\cdot\sum_{j\in J_{i}}\left[v_{j}\cdot\left\Vert \Fourier^{-1}\left(\psi_{j}\,f\right)\right\Vert _{L^{p}}\right]^{q}
\]
for all $i\in I$. This implies
\begin{align*}
\left\Vert f\right\Vert _{\FourierDecompSp{\CalQ}p{\ell_{w}^{q}}}^{q}=\sum_{i\in I}\left[w_{i}\cdot\left\Vert \Fourier^{-1}\left(\varphi_{i}\,f\right)\right\Vert _{L^{p}}\right]^{q} & \leq\left(C_{0}C_{\CalQ,\Phi,p}\right)^{q}C_{2}\cdot\sum_{i\in I}\:\sum_{j\in J_{i}}\left[v_{j}\cdot\left\Vert \Fourier^{-1}\left(\psi_{j}\,f\right)\right\Vert _{L^{p}}\right]^{q}\\
 & =\left(C_{0}C_{\CalQ,\Phi,p}\right)^{q}C_{2}\cdot\sum_{j\in J}\:\sum_{i\in I_{j}}\left[v_{j}\cdot\left\Vert \Fourier^{-1}\left(\psi_{j}\,f\right)\right\Vert _{L^{p}}\right]^{q}\\
 & \leq\left(C_{0}C_{\CalQ,\Phi,p}\right)^{q}C_{2}\cdot N\left(\CalP,\CalQ\right)\cdot\sum_{j\in J}\left[v_{j}\cdot\left\Vert \Fourier^{-1}\left(\psi_{j}\,f\right)\right\Vert _{L^{p}}\right]^{q}\\
 & =\left(C_{0}C_{\CalQ,\Phi,p}\right)^{q}C_{2}\cdot N\left(\CalP,\CalQ\right)\cdot\left\Vert f\right\Vert _{\FourierDecompSp{\CalP}p{\ell_{v}^{q}}}^{q}<\infty.
\end{align*}
Here, we used the equivalence $i\in I_{j}\Leftrightarrow j\in J_{i}$,
as well as the estimate $\left|I_{j}\right|\leq N\left(\CalP,\CalQ\right)$
for all $j\in J$.

\medskip{}

Now, we consider the case $p\in\left(0,1\right)$. Here, we first
show that both of the additional assumptions imply that there are
constants $C_{3},R>0$ which satisfy
\begin{equation}
\left|\det\left(\smash{S_{j}^{-1}}T_{i}\right)\right|\leq C_{3}\qquad\forall\,i\in I\text{ and }j\in J\text{ with }Q_{i}\cap P_{j}\neq\emptyset\,,\label{eq:EquivalentCoveringsSameSpacesDeterminantEstimate}
\end{equation}
and
\begin{equation}
P_{j}\subset T_{i}\left[\,\overline{B_{R}}\left(0\right)\right]+b_{i}\qquad\forall\,i\in I\text{ and }j\in J\text{ with }Q_{i}\cap P_{j}\neq\emptyset\,.\label{eq:EquivalentCoveringsSameSpacesSpecialCoveringInclusion}
\end{equation}

The inequality concerning the determinant is clear (with $C_{3}=C_{1}$)
in case of the second assumption. In case of the first assumption,
we use Hadamard's inequality $\left|\det A\right|\leq\left\Vert A\right\Vert ^{\dimension}$
for arbitrary matrices $A\in\R^{\dimension\times\dimension}$ (see
for instance \cite[Section 75]{RieszFunctionalAnalysis}) to deduce
\[
\left|\det\left(\smash{S_{j}^{-1}}T_{i}\right)\right|\leq\left\Vert \smash{S_{j}^{-1}}T_{i}\right\Vert ^{\dimension}\leq C_{1}^{\dimension}\qquad\forall\,i\in I\text{ and }j\in J\text{ with }Q_{i}\cap P_{j}\neq\emptyset\,.
\]

For the inclusion $P_{j}\subset T_{i}\left(\,\overline{B_{R}}\left(0\right)\right)+b_{i}$,
let us first consider the case of the second assumption. Since $\CalP$
is almost subordinate to $\CalQ$, the constant $k:=k\left(\CalP,\CalQ\right)\in\N_{0}$
is well-defined and we have $P_{j}\subset Q_{i_{j}}^{k\ast}$ for
all $j\in J$ and suitable $i_{j}\in I$. In case of $Q_{i}\cap P_{j}\neq\emptyset$,
this implies $P_{j}\subset Q_{i}^{\left(2k+2\right)\ast}$ by Lemma~\ref{lem:SubordinatenessImpliesWeakSubordination}.
But finally, Lemma~\ref{lem:SemiStructuredNormalizationNeighboring}
yields
\[
P_{j}\subset Q_{i}^{\left(2k+2\right)\ast}\subset T_{i}\left[\overline{B_{\left(2C_{\CalQ}+1\right)^{2k+2}R_{\CalQ}}}\left(0\right)\right]+b_{i}\qquad\forall\,i\in I\text{ and }j\in J\text{ with }Q_{i}\cap P_{j}\neq\emptyset\,.
\]

In case of the first assumption, define $R_{0}:=\max\left\{ R_{\CalQ},R_{\CalP}\right\} $
and let $\xi\in Q_{i}\cap P_{j}\neq\emptyset$. This yields $\omega\in Q_{i}'\subset\overline{B_{R_{0}}}\left(0\right)$
and $\eta\in P_{j}'\subset\overline{B_{R_{0}}}\left(0\right)$ with
$\xi=T_{i}\,\omega+b_{i}=S_{j}\,\eta+c_{j}$. Hence, $b_{i}=S_{j}\,\eta+c_{j}-T_{i}\,\omega$.
Now, let $x\in\overline{B_{R_{0}}}\left(0\right)$ be arbitrary. We
have
\begin{align*}
\left|T_{i}^{-1}\left(S_{j}\,x+c_{j}-b_{i}\right)\right| & =\left|T_{i}^{-1}\left(S_{j}\,x+c_{j}-\left(S_{j}\,\eta+c_{j}-T_{i}\,\omega\right)\right)\right|\\
 & \leq\left|T_{i}^{-1}S_{j}\,x\right|+\left|T_{i}^{-1}S_{j}\,\eta\right|+\left|\omega\right|\\
 & \leq C_{1}\left|x\right|+C_{1}\left|\eta\right|+\left|\omega\right|\\
 & \leq R_{0}\left(1+2C_{1}\right)=:R\,,
\end{align*}
and hence $S_{j}\,x+c_{j}\in T_{i}\left[\,\overline{B_{R}}\left(0\right)\right]+b_{i}$.
Since $x\in\overline{B_{R_{0}}}\left(0\right)$ was arbitrary, this
implies 
\[
P_{j}=S_{j}P_{j}'+c_{j}\subset S_{j}\left[\,\overline{B_{R_{0}}}\left(0\right)\right]+c_{j}\subset T_{i}\left[\,\overline{B_{R}}\left(0\right)\right]+b_{i}\qquad\forall\,i\in I\text{ and }j\in J\text{ with }Q_{i}\cap P_{j}\neq\emptyset\,.
\]

With this preparation, we can prove the actual claim: The argument
is almost the same as for ${p\in\left[1,\infty\right]}$; only the
application of Young's inequality and of the triangle inequality for
$L^{p}$ in equation~(\ref{eq:EquivalentCoveringsYieldSameDecompositionSpacesFundamentalBanach})
need to be replaced. In short, we only need to establish equation~(\ref{eq:EquivalentCoveringsYieldSameDecompositionSpacesFundamentalBanach})
also for $p\in\left(0,1\right)$; the rest of the proof then proceeds
as for $p\in\left[1,\infty\right]$. The only caveat is that since
the assumptions are \emph{not} symmetric in $\CalQ,\CalP$ anymore,
we also need to show the ``reverse'' version of equation~(\ref{eq:EquivalentCoveringsYieldSameDecompositionSpacesFundamentalBanach}).

Since $\left\Vert \mybullet\right\Vert _{L^{p}}$ is a quasi-norm
and because of the uniform bound $\left|J_{i}\right|\leq N\left(\CalQ,\CalP\right)=N$,
there is a constant $C_{4}=C_{4}\left(N,p\right)>0$ with
\begin{align}
w_{i}\cdot\left\Vert \Fourier^{-1}\left(\varphi_{i}\,f\right)\right\Vert _{L^{p}} & =w_{i}\cdot\left\Vert \Fourier^{-1}\left(\varphi_{i}\,\psi_{J_{i}}\,f\right)\right\Vert _{L^{p}}\nonumber \\
 & \leq C_{4}\cdot w_{i}\cdot\sum_{j\in J_{i}}\left\Vert \Fourier^{-1}\left(\varphi_{i}\,\psi_{j}\,f\right)\right\Vert _{L^{p}}\nonumber \\
 & \overset{\left(\dagger\right)}{\leq}C_{4}C_{5}\cdot w_{i}\cdot\sum_{j\in J_{i}}\left|\det T_{i}\right|^{\frac{1}{p}-1}\cdot\left\Vert \Fourier^{-1}\varphi_{i}\right\Vert _{L^{p}}\cdot\left\Vert \Fourier^{-1}\left(\psi_{j}\,f\right)\right\Vert _{L^{p}}\nonumber \\
 & \leq C_{0}C_{4}C_{5}C_{\CalQ,\Phi,p}\cdot\sum_{j\in J_{i}}v_{j}\cdot\left\Vert \Fourier^{-1}\left(\psi_{j}\,f\right)\right\Vert _{L^{p}}\qquad\forall\,i\in I\,.\label{eq:EquivalentCoveringsYieldSameSpacesQuasiBanach1}
\end{align}
Here, the step marked with $\left(\dagger\right)$ is justified as
follows: By possibly enlarging the constant $R$ from equation~(\ref{eq:EquivalentCoveringsSameSpacesSpecialCoveringInclusion}),
we can assume $R\geq R_{0}=\max\left\{ R_{\CalQ},R_{\CalP}\right\} $.
This yields the inclusions
\[
\supp\varphi_{i}\subset\overline{Q_{i}}\subset T_{i}\left[\,\overline{B_{R}}\left(0\right)\right]+b_{i}=:M_{i}\quad\text{and}\quad\supp\psi_{j}\subset\overline{P_{j}}\subset T_{i}\left[\,\overline{B_{R}}\left(0\right)\right]+b_{i}=M_{i}\qquad\forall\,j\in J_{i}\,,
\]
see equation~(\ref{eq:EquivalentCoveringsSameSpacesSpecialCoveringInclusion}).
Thus, Theorem~\ref{thm:QuasiBanachConvolution} shows 
\begin{align}
\left\Vert \Fourier^{-1}\left(\varphi_{i}\,\psi_{j}\,f\right)\right\Vert _{L^{p}} & \leq\left[\lambda\left(M_{i}-M_{i}\right)\right]^{\frac{1}{p}-1}\cdot\left\Vert \Fourier^{-1}\varphi_{i}\right\Vert _{L^{p}}\cdot\left\Vert \Fourier^{-1}\left(\psi_{j}\,f\right)\right\Vert _{L^{p}}\nonumber \\
 & =\left[\lambda\left(T_{i}\left[\overline{B_{R}}\left(0\right)-\overline{B_{R}}\left(0\right)\right]\right)\right]^{\frac{1}{p}-1}\cdot\left\Vert \Fourier^{-1}\varphi_{i}\right\Vert _{L^{p}}\cdot\left\Vert \Fourier^{-1}\left(\psi_{j}\,f\right)\right\Vert _{L^{p}}\nonumber \\
 & \leq\left[\lambda\left(\overline{B_{2R}}\left(0\right)\right)\right]^{\frac{1}{p}-1}\cdot\left|\det T_{i}\right|^{\frac{1}{p}-1}\cdot\left\Vert \Fourier^{-1}\varphi_{i}\right\Vert _{L^{p}}\cdot\left\Vert \Fourier^{-1}\left(\psi_{j}\,f\right)\right\Vert _{L^{p}}\nonumber \\
 & =:C_{5}\cdot\left|\det T_{i}\right|^{\frac{1}{p}-1}\cdot\left\Vert \Fourier^{-1}\varphi_{i}\right\Vert _{L^{p}}\cdot\left\Vert \Fourier^{-1}\left(\psi_{j}\,f\right)\right\Vert _{L^{p}},\label{eq:EquivalentCoveringsYieldSameSpacesQuasiBanach2}
\end{align}
where the constant $C_{5}$ only depends on $R>0$, on $\dimension\in\N$
and on $p\in\left(0,1\right)$.

For the ``reverse'' version of equation~(\ref{eq:EquivalentCoveringsYieldSameDecompositionSpacesFundamentalBanach}),
a similar argument applies: Because of the uniform bound $\left|I_{j}\right|\leq N\left(\CalP,\CalQ\right)=:N_{2}$,
there is a constant $C_{6}=C_{6}\left(N_{2},p\right)>0$ with
\begin{align*}
v_{j}\cdot\left\Vert \Fourier^{-1}\left(\psi_{j}\,f\right)\right\Vert _{L^{p}} & =v_{j}\cdot\left\Vert \Fourier^{-1}\left(\psi_{j}\,\varphi_{I_{j}}\,f\right)\right\Vert _{L^{p}}\\
 & \leq C_{6}\cdot v_{j}\cdot\sum_{i\in I_{j}}\left\Vert \Fourier^{-1}\left(\psi_{j}\,\varphi_{i}\,f\right)\right\Vert _{L^{p}}\\
 & \overset{\left(\dagger\right)}{\leq}C_{5}C_{6}\cdot v_{j}\cdot\sum_{i\in I_{j}}\left|\det T_{i}\right|^{\frac{1}{p}-1}\cdot\left\Vert \Fourier^{-1}\psi_{j}\right\Vert _{L^{p}}\cdot\left\Vert \Fourier^{-1}\left(\varphi_{i}\,f\right)\right\Vert _{L^{p}}\\
 & \overset{\left(\ddagger\right)}{\leq}C_{0}C_{3}^{\frac{1}{p}-1}C_{5}C_{6}\cdot\sum_{i\in I_{j}}\left[\left|\det S_{j}\right|^{\frac{1}{p}-1}\cdot\left\Vert \Fourier^{-1}\psi_{j}\right\Vert _{L^{p}}\cdot w_{i}\cdot\left\Vert \Fourier^{-1}\left(\varphi_{i}\,f\right)\right\Vert _{L^{p}}\right]\\
 & \leq C_{0}C_{3}^{\frac{1}{p}-1}C_{5}C_{6}C_{\CalP,\Psi,p}\cdot\sum_{i\in I_{j}}\left[w_{i}\cdot\left\Vert \Fourier^{-1}\left(\varphi_{i}\,f\right)\right\Vert _{L^{p}}\right].
\end{align*}
Here, we used $\frac{1}{p}-1>0$, together with estimate~(\ref{eq:EquivalentCoveringsSameSpacesDeterminantEstimate})
at $\left(\ddagger\right)$. The justification for $\left(\dagger\right)$
is exactly as in equation~(\ref{eq:EquivalentCoveringsYieldSameSpacesQuasiBanach2})
above.

Now, using the preceding estimate and equation~(\ref{eq:EquivalentCoveringsYieldSameSpacesQuasiBanach1})
instead of equation~(\ref{eq:EquivalentCoveringsYieldSameDecompositionSpacesFundamentalBanach}),
the proof for $p\in\left(0,1\right)$ can be completed just as for
$p\in\left[1,\infty\right]$.
\end{proof}

\subsection{Improved necessary conditions}

\label{subsec:ImprovedNecessaryConditions}In this subsection, we
use many of the techniques from the preceding subsection to improve
upon the ``elementary'' necessary conditions that we developed in
Subsection \ref{subsec:ElementaryNecessaryConditions}.

As in the previous subsection, we will use functions of the form
\[
f=\sum_{i\in I_{0}}M_{z_{i}}\left(\varepsilon_{i}c_{i}\cdot\gamma_{i}\right),\quad\text{with}\quad\gamma_{i}\in\TestFunctionSpace{Q_{i}},\quad z_{i}\in\R^{\dimension},\quad\varepsilon_{i}\in\left\{ \pm1\right\} ,\quad\text{and}\quad c_{i}\in\Compl,
\]
to ``test'' the embedding $\FourierDecompSp{\CalQ}{p_{1}}Y\hookrightarrow\FourierDecompSp{\CalP}{p_{2}}Z$.
Our next lemma indicates how we will choose the functions $\gamma_{i}$.
\begin{lem}
\label{lem:NormOfClusteredBAPUAndTestFunctionBuildingBlocks}Let $\emptyset\neq\CalQ\subset\R^{\dimension}$
be open, and let $\CalQ=\left(Q_{i}\right)_{i\in I}=\left(T_{i}Q_{i}'+b_{i}\right)_{i\in I}$
be a semi-structured covering of $\CalO$ with $L^{p}$-BAPU $\Phi=\left(\varphi_{i}\right)_{i\in I}$
for some $p\in\left(0,\infty\right]$.

Let $q\in\left[\min\left\{ 1,p\right\} ,\infty\right]$ and $\ell\in\N_{0}$.
Then there is a constant $C_{1}=C_{1}\left(\dimension,p,\ell,\CalQ,C_{\CalQ,\Phi,p}\right)>0$
such that
\begin{equation}
\left\Vert \smash{\Fourier^{-1}}\varphi_{M_{i}}\right\Vert _{L^{q}}\leq C_{1}\cdot\left|\det T_{i}\right|^{1-\frac{1}{q}}\qquad\forall\,i\in I\text{ and all sets }M_{i}\subset i^{\ell\ast}\,.\label{eq:UpperBAPUNormEstimate}
\end{equation}

\medskip{}

If $\CalQ$ is tight, there are two families of functions $\left(\gamma_{i}\right)_{i\in I},\left(\phi_{i}\right)_{i\in I}$
with the following properties:

\begin{enumerate}[leftmargin=0.9cm]
\item We have $\gamma_{i},\phi_{i}\in\TestFunctionSpace{\CalO}$ and $\supp\gamma_{i},\supp\phi_{i}\subset Q_{i}$,
as well as $\gamma_{i},\phi_{i}\geq0$ for all $i\in I$.
\item We have $\phi_{i}\gamma_{i}=\gamma_{i}$ for all $i\in I$.
\item For each $s\in\left(0,\infty\right]$, there are constants $C_{2}=C_{2}\left(\dimension,s\right)>0$
and $C_{3}=C_{3}\left(\dimension,s\right)>0$ with
\[
\quad\quad\quad\left\Vert \smash{\Fourier^{-1}}\gamma_{i}\right\Vert _{L^{s}}=C_{2}\cdot\varepsilon_{\CalQ}^{\dimension\left(1-\frac{1}{s}\right)}\cdot\left|\det T_{i}\right|^{1-\frac{1}{s}}\quad\!\text{and}\!\quad\left\Vert \smash{\Fourier^{-1}}\phi_{i}\right\Vert _{L^{s}}=C_{3}\cdot\varepsilon_{\CalQ}^{\dimension\left(1-\frac{1}{s}\right)}\cdot\left|\det T_{i}\right|^{1-\frac{1}{s}}\!\quad\!\forall\,i\in I.\qedhere
\]
\end{enumerate}
\end{lem}

\begin{proof}
Let $r:=\min\left\{ 1,p\right\} $. By definition of an $L^{p}$-BAPU,
we know that 
\[
C_{\CalQ,\Phi,p}=C_{\CalQ,\Phi,r}=\sup_{i\in I}\left[\left|\det T_{i}\right|^{\frac{1}{r}-1}\cdot\left\Vert \smash{\Fourier^{-1}}\varphi_{i}\right\Vert _{L^{r}}\right]
\]
is finite. By applying Lemma~\ref{lem:LocalEmbeddingInHigherLpSpaces}
(with $k=0$, $p_{0}=p_{1}=r$ and $p_{2}=q\geq r=p_{1}$, as well
as $f\equiv1$), we get a constant $K_{0}=K_{0}\left(\dimension,r,\CalQ\right)=K_{0}\left(\dimension,p,\CalQ\right)>0$
with
\[
\left\Vert \smash{\Fourier^{-1}}\varphi_{i}\right\Vert _{L^{q}}\leq K_{0}\cdot\left|\det T_{i}\right|^{\frac{1}{r}-\frac{1}{q}}\cdot\left\Vert \smash{\Fourier^{-1}}\varphi_{i}\right\Vert _{L^{r}}\leq K_{0}C_{\CalQ,\Phi,p}\left|\det T_{i}\right|^{\frac{1}{r}-\frac{1}{q}}\cdot\left|\det T_{i}\right|^{1-\frac{1}{r}}=K_{0}C_{\CalQ,\Phi,p}\cdot\left|\det T_{i}\right|^{1-\frac{1}{q}}
\]
for all $i\in I$, because of $\varphi_{i}\equiv0$ on $\CalO\setminus Q_{i}$.

Since $\left\Vert \mybullet\right\Vert _{L^{q}}$ is a quasi-norm
with triangle constant $2^{\max\left\{ 0,q^{-1}-1\right\} }\leq2^{r^{-1}-1}$
(see \cite[equation (1.1.4)]{GrafakosClassical}), and because of
the uniform bound $\left|M_{i}\right|\leq\left|i^{\ell\ast}\right|\leq N_{\CalQ}^{\ell}$
(which was shown in Lemma~\ref{lem:SemiStructuredClusterInvariant}),
we obtain a constant $K_{1}=K_{1}\left(N_{\CalQ},\ell,r\right)=K_{1}\left(N_{\CalQ},\ell,p\right)>0$
with
\[
\left\Vert \Fourier^{-1}\varphi_{M_{i}}\right\Vert _{L^{q}}\leq K_{1}\cdot\sum_{j\in M_{i}}\left\Vert \Fourier^{-1}\varphi_{j}\right\Vert _{L^{q}}\leq K_{0}K_{1}C_{\CalQ,\Phi,p}\cdot\sum_{j\in M_{i}}\left|\det T_{j}\right|^{1-\frac{1}{q}}\qquad\forall\,i\in I\text{ and all }M_{i}\subset i^{\ell\ast}\,.
\]

Next, recall from equation~(\ref{eq:DeterminantIsModerate}) that
the weight $\left(\left|\det T_{i}\right|\right)_{i\in I}$ is $\CalQ$-moderate
with ``moderateness constant'' $C_{\CalQ}^{\dimension}$. In combination
with Lemma~\ref{lem:SemiStructuredClusterInvariant}, we thus see
$\left|\det T_{j}\right|\leq C_{\CalQ}^{\dimension\ell}\cdot\left|\det T_{i}\right|$
for all $i\in I$ and $j\in i^{\ell\ast}$. There are now two cases:
In case of $q\geq1$, we have $0\leq1-q^{-1}\leq1$, and thus $\left|\det T_{j}\right|^{1-q^{-1}}\leq C_{\CalQ}^{\dimension\ell\left(1-q^{-1}\right)}\cdot\left|\det T_{i}\right|^{1-q^{-1}}\leq C_{\CalQ}^{\dimension\ell}\cdot\left|\det T_{i}\right|^{1-q^{-1}}$,
since $C_{\CalQ}\geq1$. Otherwise, if $q\in\left(0,1\right)$, then
$0<p\leq q<1$ and thus $q^{-1}-1\leq p^{-1}-1\leq p^{-1}$. Therefore,
by applying $\left|\det T_{j}\right|\leq C_{\CalQ}^{\dimension\ell}\cdot\left|\det T_{i}\right|$
with interchanged roles of $i,j$, we get
\[
\left|\det T_{j}\right|^{1-q^{-1}}=\left(1/\left|\det T_{j}\right|\right)^{q^{-1}-1}\leq C_{\CalQ}^{\dimension\ell\left(q^{-1}-1\right)}\cdot\left(1/\left|\det T_{i}\right|\right)^{q^{-1}-1}\leq C_{\CalQ}^{\dimension\ell/p}\cdot\left|\det T_{i}\right|^{1-q^{-1}}\,.
\]
Because of the equivalence $j\in i^{\ell\ast}\Longleftrightarrow i\in j^{\ell\ast}$,
we thus see that in each case, there is a constant $K_{2}=K_{2}\left(C_{\CalQ},\dimension,\ell,p\right)\geq1$
with
\[
K_{2}^{-1}\cdot\left|\det T_{i}\right|^{1-\frac{1}{q}}\leq\left|\det T_{j}\right|^{1-\frac{1}{q}}\leq K_{2}\cdot\left|\det T_{i}\right|^{1-\frac{1}{q}}\qquad\forall\,i\in I\text{ and }j\in i^{\ell\ast}\supset M_{i}\,.
\]

All in all, we conclude
\begin{align*}
\left\Vert \Fourier^{-1}\varphi_{M_{i}}\right\Vert _{L^{q}} & \leq K_{0}K_{1}K_{2}C_{\CalQ,\Phi,p}\cdot\sum_{j\in M_{i}}\left|\det T_{i}\right|^{1-\frac{1}{q}}\\
\left({\scriptstyle \text{since }\left|M_{i}\right|\leq\left|i^{\ell\ast}\right|\leq N_{\CalQ}^{\ell}}\right) & \leq K_{0}K_{1}K_{2}N_{\CalQ}^{\ell}C_{\CalQ,\Phi,p}\cdot\left|\det T_{i}\right|^{1-\frac{1}{q}}\qquad\forall\,i\in I\text{ and all }M_{i}\subset i^{\ell\ast}\,.
\end{align*}
We have thus established equation~(\ref{eq:UpperBAPUNormEstimate}),
for $C_{1}:=K_{0}K_{1}K_{2}N_{\CalQ}^{\ell}C_{\CalQ,\Phi,p}$.

\medskip{}

For the construction of $\left(\gamma_{i}\right)_{i\in I},\left(\phi_{i}\right)_{i\in I}$,
fix a non-trivial nonnegative function $\gamma\in\TestFunctionSpace{B_{1/2}\left(0\right)}$
and a nonnegative function $\phi\in\TestFunctionSpace{B_{1}\left(0\right)}$
with $\phi\equiv1$ on $B_{1/2}\left(0\right)$ and let $\varepsilon:=\varepsilon_{\CalQ}$.
For each $i\in I$, choose $c_{i}\in\R^{\dimension}$ with $B_{\varepsilon}\left(c_{i}\right)\subset Q_{i}'$
(see Remark~\ref{rem:CoveringTypesRemark} for the existence of $c_{i}$)
and define
\[
\gamma_{i}:=L_{b_{i}}\left[\left(L_{c_{i}}\left[\gamma\circ\smash{\varepsilon^{-1}}{\rm id}\right]\right)\circ T_{i}^{-1}\right]\qquad\text{ and }\qquad\phi_{i}:=L_{b_{i}}\left[\left(L_{c_{i}}\left[\phi\circ\smash{\varepsilon^{-1}}{\rm id}\right]\right)\circ T_{i}^{-1}\right]\,.
\]
Clearly, $\gamma_{i}\geq0$ and $\phi_{i}\geq0$. Furthermore, $\supp\gamma\subset\supp\phi$
and hence
\begin{align*}
\supp\gamma_{i}\subset\supp\phi_{i} & =\supp\left(\left(L_{c_{i}}\left[\phi\circ\smash{\varepsilon^{-1}}{\rm id}\right]\right)\circ T_{i}^{-1}\right)+b_{i}\\
 & =T_{i}\left(\supp\left(L_{c_{i}}\left[\phi\circ\smash{\varepsilon^{-1}}{\rm id}\right]\right)\right)+b_{i}\\
 & =T_{i}\left(\supp\left(\phi\circ\smash{\varepsilon^{-1}}{\rm id}\right)+c_{i}\right)+b_{i}\\
 & \subset T_{i}\left[B_{\varepsilon}\left(c_{i}\right)\right]+b_{i}\subset T_{i}Q_{i}'+b_{i}=Q_{i}\subset\CalO,
\end{align*}
and thus $\gamma_{i},\phi_{i}\in\TestFunctionSpace{\CalO}$. Finally,
since $\phi\equiv1$ on $B_{1/2}\left(0\right)\supset\supp\gamma$,
it is easy to see $\phi\gamma=\gamma$, which then implies $\phi_{i}\gamma_{i}=\gamma_{i}$,
as desired.

To complete the proof, note that a calculation using standard properties
of the Fourier transform (see \cite[Theorem 8.22]{FollandRA}) shows
\[
\Fourier^{-1}\gamma_{i}=\varepsilon^{\dimension}\cdot\left|\det T_{i}\right|\cdot M_{b_{i}}\left[\left(M_{c_{i}}\left[\left(\smash{\Fourier^{-1}}\gamma\right)\circ\varepsilon{\rm id}\right]\right)\circ T_{i}^{T}\right]\,,
\]
and thus
\[
\left\Vert \Fourier^{-1}\gamma_{i}\right\Vert _{L^{s}}=\varepsilon^{\dimension}\cdot\left|\det T_{i}\right|\cdot\left|\det T_{i}\right|^{-\frac{1}{s}}\cdot\left\Vert \left(\Fourier^{-1}\gamma\right)\circ\varepsilon{\rm id}\right\Vert _{L^{s}}=\left\Vert \Fourier^{-1}\gamma\right\Vert _{L^{s}}\cdot\left|\det T_{i}\right|^{1-\frac{1}{s}}\cdot\varepsilon^{\dimension\left(1-\frac{1}{s}\right)}
\]
for all $i\in I$ and $s\in\left(0,\infty\right]$. Thus, since $\gamma$
only depends on the dimension $\dimension$, the family $\left(\gamma_{i}\right)_{i\in I}$
has the desired properties, with $C_{2}\left(\dimension,s\right)=\left\Vert \Fourier^{-1}\gamma\right\Vert _{L^{s}}$.
With $C_{3}\left(\dimension,s\right)=\left\Vert \Fourier^{-1}\phi\right\Vert _{L^{s}}$,
exactly the same calculation shows 
\[
\left\Vert \Fourier^{-1}\phi_{i}\right\Vert _{L^{s}}=C_{3}\left(\dimension,s\right)\cdot\varepsilon^{\dimension\left(1-\frac{1}{s}\right)}\cdot\left|\det T_{i}\right|^{1-\frac{1}{s}}\qquad\forall\,i\in I\text{ and }s\in\left(0,\infty\right]\,.\qedhere
\]
\end{proof}
Given these preparations, we can now state and prove the first principal
result of this subsection. We remark that the following theorem is
an improved version of (the second part of) \cite[Theorem 5.3.6]{VoigtlaenderPhDThesis}
from my PhD thesis.
\begin{thm}
\label{thm:BurnerNecessaryConditionCoarseInFine}Let $\emptyset\neq\CalO,\CalO'\subset\R^{\dimension}$
be open, let $p_{1},p_{2}\in\left(0,\infty\right]$ and let $\CalP=\left(P_{j}\right)_{j\in J}=\left(S_{j}P_{j}'+c_{j}\right)_{j\in J}$
be a \emph{tight} semi-structured $L^{p_{2}}$-decomposition covering
of $\CalO'$. Furthermore, let $\CalQ=\left(Q_{i}\right)_{i\in I}$
be an $L^{p_{1}}$-decomposition covering of $\CalO$. Finally, let
$Y\subset\Compl^{I}$ and $Z\subset\Compl^{J}$ be $\CalQ$-regular
and $\CalP$-regular with triangle constants $C_{Y},C_{Z}\geq1$,
respectively.

Let $J_{0}\subset J$ be arbitrary with $P_{j}\subset\CalO$ for all
$j\in J_{0}$. Define
\[
K:=\bigcup_{j\in J_{0}}P_{j}\subset\CalO\cap\CalO'
\]
and—with $\CalD_{K}^{\CalQ,p_{1},Y}$ as after equation~(\ref{eq:GeneralEmbeddingRequirement})—assume
that there is a bounded linear map
\begin{equation}
\iota:\left(\CalD_{K}^{\CalQ,p_{1},Y},\left\Vert \mybullet\right\Vert _{\FourierDecompSp{\CalQ}{p_{1}}Y}\right)\to\FourierDecompSp{\CalP}{p_{2}}Z\label{eq:EmbeddingCoarseInFine}
\end{equation}
satisfying $\left\langle \iota f,\,\varphi\right\rangle _{\CalD'}=\left\langle f,\,\varphi\right\rangle _{\CalD'}$
for all $\varphi\in\TestFunctionSpace{\CalO\cap\CalO'}$ and all $f\in\CalD_{K}^{\CalQ,p_{1},Y}$.
In case of $p_{1}\in\left(0,1\right)$, assume additionally that $\CalP_{J_{0}}:=\left(P_{j}\right)_{j\in J_{0}}$
is almost subordinate to $\CalQ$.

Then the embedding
\begin{alignat}{2}
\eta: & \; & \ell_{0}\left(J_{0}\right)\cap Y\left(\left[\ell^{p_{1}}\left(J_{i}\cap J_{0}\right)\right]_{i\in I}\right) & \hookrightarrow Z_{\left|\det S_{j}\right|^{p_{1}^{-1}-p_{2}^{-1}}}\;,\label{eq:BurnerNecessaryCoarseInFineDiscreteEmbedding}\\
 &  & \left(x_{j}\right)_{j\in J_{0}} & \mapsto\left(x_{j}\right)_{j\in J}\text{ with }x_{j}=0\text{ for }j\in J\setminus J_{0}\nonumber 
\end{alignat}
is well-defined and bounded, with
\[
\vertiii{\eta}\leq C\cdot\vertiii{\iota}
\]
for some constant
\[
C=C\left(\dimension,p_{1},p_{2},k\left(\CalP_{J_{0}},\CalQ\right),\CalQ,\CalP,\varepsilon_{\CalP},C_{\CalQ,\Phi,p_{1}},C_{Z},\vertiii{\Gamma_{\CalP}}_{Z\to Z}\right),
\]
where the dependence on $k\left(\CalP_{J_{0}},\CalQ\right)$ can be
dropped for $p_{1}\in\left[1,\infty\right]$. Here, as usual, the
$L^{p_{1}}$-BAPU $\Phi=\left(\varphi_{i}\right)_{i\in I}$ has to
be used to calculate the (quasi)-norm on the decomposition space $\FourierDecompSp{\CalQ}{p_{1}}Y$
when calculating $\vertiii{\iota}$.

\medskip{}

If $Z$ satisfies the Fatou property (see Definition~\ref{def:FatouProperty}),
the embedding $\eta$ is bounded (with the same estimate for the norm)
even \emph{without} restricting to $\ell_{0}\left(J_{0}\right)$.
\end{thm}

\begin{rem}
\label{rem:BurnerNecessaryConditionCoarseInFine}Two remarks are in
order:\medskip{}

(1) A remarkable property of the preceding result is—at least for
$p_{1}\in\left[1,\infty\right]$—that we need to impose virtually
no restrictions upon the relation between the coverings $\CalQ$ and
$\CalP$—only that $P_{j}\subset\CalO$ for all $j\in J_{0}$. Only
for $p_{1}\in\left(0,1\right)$, we assume that $\CalP_{J_{0}}$ is
almost subordinate to $\CalQ$.

However, under these minimal assumptions, the theorem is only of limited
value since it is difficult to (dis)prove boundedness of the embedding~(\ref{eq:BurnerNecessaryCoarseInFineDiscreteEmbedding})
in this general setting, i.e., without assuming subordinateness. In
contrast, if $\CalP_{J_{0}}$ is almost subordinate to $\CalQ$ and
if the ``global'' components $Y,Z$ are in fact weighted $\ell^{q}$
spaces, one can use Corollary~\ref{cor:EmbeddingCoarseIntoFineSimplification}
to greatly simplify the task of checking whether the embedding~(\ref{eq:BurnerNecessaryCoarseInFineDiscreteEmbedding})
is bounded; see also Theorem~\ref{thm:SummaryCoarseIntoFine} below.

\medskip{}

(2) Despite the pessimistic view of the preceding point, one can still
derive nontrivial conditions from the boundedness of the embedding~(\ref{eq:BurnerNecessaryCoarseInFineDiscreteEmbedding}):
Indeed, let $j_{0}\in J_{0}$ and assume $\Indicator_{I_{j_{0}}}\in Y$.
Now, for $\ell\in I$ there are two cases:

\begin{casenv}
\item We have $j_{0}\in J_{\ell}\cap J_{0}$. This implies $\ell\in I_{j_{0}}$
and thus $\left\Vert \delta_{j_{0}}\right\Vert _{\ell^{p_{1}}\left(J_{\ell}\cap J_{0}\right)}=1\leq\Indicator_{I_{j_{0}}}\left(\ell\right)$.
\item We have $j_{0}\notin J_{\ell}\cap J_{0}$. This implies $\left\Vert \delta_{j_{0}}\right\Vert _{\ell^{p_{1}}\left(J_{\ell}\cap J_{0}\right)}=0\leq\Indicator_{I_{j_{0}}}\left(\ell\right)$.
\end{casenv}
All in all, we get $\delta_{j_{0}}\in\ell_{0}\left(J_{0}\right)\cap Y\left(\left[\ell^{p_{1}}\left(J_{\ell}\cap J_{0}\right)\right]_{\ell\in I}\right)$
with
\begin{align*}
\left|\det S_{j_{0}}\right|^{p_{1}^{-1}-p_{2}^{-1}}\cdot\left\Vert \delta_{j_{0}}\right\Vert _{Z}=\left\Vert \eta\left(\delta_{j_{0}}\right)\right\Vert _{Z_{\left|\det S_{j}\right|^{p_{1}^{-1}-p_{2}^{-1}}}} & \leq\vertiii{\eta}\cdot\left\Vert \delta_{j_{0}}\right\Vert _{Y\left(\left[\ell^{p_{1}}\left(J_{\ell}\cap J_{0}\right)\right]_{\ell\in I}\right)}\\
 & =\vertiii{\eta}\cdot\left\Vert \left(\left\Vert \delta_{j_{0}}\right\Vert _{\ell^{p_{1}}\left(J_{\ell}\cap J_{0}\right)}\right)_{\ell\in I}\right\Vert _{Y}\leq\vertiii{\eta}\cdot\left\Vert \Indicator_{I_{j_{0}}}\right\Vert _{Y}.
\end{align*}

Of course, in general, it is not easy to guarantee $\Indicator_{I_{j_{0}}}\in Y$,
or to estimate $\left\Vert \Indicator_{I_{j_{0}}}\right\Vert _{Y}$.
But if $\CalP_{J_{0}}$ is almost subordinate to $\CalQ$, the above
estimate can be improved: Indeed, simply assume that $\delta_{i}\in Y$
for some $i\in I$ with $Q_{i}\cap P_{j_{0}}\neq\emptyset$. Setting
$k:=k\left(\CalP_{J_{0}},\CalQ\right)$, Lemma~\ref{lem:SubordinatenessImpliesWeakSubordination}
shows $I_{j_{0}}\subset i^{\left(2k+2\right)\ast}$ and hence $\Indicator_{I_{j_{0}}}\leq\Theta_{2k+2}\,\delta_{i}$
with the $\ell$-fold clustering map $\Theta_{\ell}$ from Lemma~\ref{lem:HigherOrderClusteringMap}.
But in view of that lemma and since $Y$ is solid, this implies $\Indicator_{I_{j_{0}}}\in Y$
with $\left\Vert \Indicator_{I_{j_{0}}}\right\Vert _{Y}\leq\vertiii{\Gamma_{\CalQ}}_{Y\to Y}^{2k+2}\cdot\left\Vert \delta_{i}\right\Vert _{Y}$.
Hence, if $\CalP_{J_{0}}$ is almost subordinate to $\CalQ$, then
\begin{equation}
\left|\det S_{j_{0}}\right|^{p_{1}^{-1}-p_{2}^{-1}}\cdot\left\Vert \delta_{j_{0}}\right\Vert _{Z}\lesssim\left\Vert \delta_{i}\right\Vert _{Y}\qquad\forall\,j_{0}\in J_{0}\text{ and }i\in I\text{ with }\delta_{i}\in Y\text{ and }Q_{i}\cap P_{j_{0}}\neq\emptyset.\label{eq:ElementaryNecessaryConditionCoarseInFine}
\end{equation}
This estimate can be seen as a generalization to the case $p_{1}\neq p_{2}$
of the estimate $\left\Vert \delta_{j}\right\Vert _{Z}\lesssim\left\Vert \delta_{i}\right\Vert _{Y}$
from Lemma~\ref{lem:SimpleNecessaryCondition}, which required $Q_{i}\cap P_{j}\neq\emptyset$
and $p_{1}=p_{2}$. In many cases, this—very simple—estimate already
suffices to show that a certain embedding between decomposition spaces
can \emph{not} exist.
\end{rem}

\begin{proof}[Proof of Theorem~\ref{thm:BurnerNecessaryConditionCoarseInFine}]
For $j\in J_{0}$, we have $P_{j}\subset\CalO$ and $P_{j}\neq\emptyset$
(by tightness of $\CalP$) and hence $P_{j}\cap Q_{i}\neq\emptyset$
for some $i\in I$, since $\CalQ$ covers $\CalO$. Hence, $j\in J_{i}\cap J_{0}$
for some $i\in I$. All in all, this shows $J_{0}=\bigcup_{i\in I}\left(J_{i}\cap J_{0}\right)$,
so that
\[
V:=Y\left(\left[\ell^{p_{1}}\left(J_{i}\cap J_{0}\right)\right]_{i\in I}\right)\subset\Compl^{J_{0}}\qquad\text{ and }\qquad V_{0}:=\ell_{0}\left(J_{0}\right)\cap V
\]
are solid sequence spaces on $J_{0}$.

We begin by setting up some quantities which we will need for constructing
suitable functions to ``test'' the embedding $\iota$. For $r_{0}:=N_{\CalP}^{3}$,
the disjointization lemma (Lemma~\ref{lem:DisjointizationPrinciple})
yields a partition $J=\biguplus_{r=1}^{r_{0}}J^{\left(r\right)}$
satisfying $P_{j}^{\ast}\cap P_{\ell}^{\ast}=\emptyset$ for all $j,\ell\in J^{\left(r\right)}$
with $j\neq\ell$ and arbitrary $r\in\underline{r_{0}}$. Let $\varepsilon:=\varepsilon_{\CalP}$.
With this notation, Lemma~\ref{lem:NormOfClusteredBAPUAndTestFunctionBuildingBlocks}
(applied to the tight semi-structured covering $\CalP$) yields two
families $\left(\gamma_{j}\right)_{j\in J},\left(\phi_{j}\right)_{j\in J}$
of nonnegative functions $\gamma_{j},\phi_{j}\in\TestFunctionSpace{\CalO'}$
with $\phi_{j}\gamma_{j}=\gamma_{j}$ and $\supp\gamma_{j},\supp\phi_{j}\subset P_{j}$,
as well as 
\[
\left\Vert \Fourier^{-1}\gamma_{j}\right\Vert _{L^{p}}=C_{1}^{\left(p\right)}\cdot\varepsilon^{\dimension\left(1-\frac{1}{p}\right)}\cdot\left|\det S_{j}\right|^{1-\frac{1}{p}}\qquad\text{ and }\qquad\left\Vert \Fourier^{-1}\phi_{j}\right\Vert _{L^{p}}=C_{2}^{\left(p\right)}\cdot\varepsilon^{\dimension\left(1-\frac{1}{p}\right)}\cdot\left|\det S_{j}\right|^{1-\frac{1}{p}}
\]
for all $j\in J$ and $p\in\left(0,\infty\right]$, with $C_{1}^{\left(p\right)}=C_{1}^{\left(p\right)}\!\left(\dimension\right)$
and $C_{2}^{\left(p\right)}=C_{2}^{\left(p\right)}\!\left(\dimension\right)$.
Note that 
\[
\supp\gamma_{j},\supp\phi_{j}\subset P_{j}\subset K\subset\CalO\qquad\forall\,j\in J_{0}.
\]
Let $\Psi=\left(\psi_{j}\right)_{j\in J}$ be an $L^{p_{2}}$-BAPU
for $\CalP$.

\medskip{}

Now, let $c=\left(c_{j}\right)_{j\in J_{0}}\in V_{0}$ be arbitrary
and let $M:=\supp c$, which is a finite subset of $J_{0}$. Define
$\zeta_{j}:=\left|\det S_{j}\right|^{p_{1}^{-1}-1}\cdot c_{j}$ for
$j\in J_{0}$ and note $\supp\zeta=\supp c=M$ for $\zeta=\left(\zeta_{j}\right)_{j\in J_{0}}$.
Now, for $z=\left(z_{j}\right)_{j\in J_{0}}\in\left(\R^{\dimension}\right)^{J_{0}}$
and $r\in\underline{r_{0}}$, define $J_{0}^{\left(r\right)}:=J_{0}\cap J^{\left(r\right)}$,
as well as $M^{\left(r\right)}:=M\cap J^{\left(r\right)}$ and
\begin{equation}
f_{z}^{\left(r\right)}:=\sum_{j\in J_{0}^{\left(r\right)}}M_{z_{j}}\left(\zeta_{j}\cdot\gamma_{j}\right)=\sum_{j\in M^{\left(r\right)}}M_{z_{j}}\left(\zeta_{j}\cdot\gamma_{j}\right).\label{eq:BurnerNecessaryCoarseInFineTestFunction}
\end{equation}
Note that we have $\supp\gamma_{j}\subset K\subset\CalO$ for all
$j\in J_{0}$, so that we have $f_{z}^{\left(r\right)}\in\CalD_{K}^{\CalQ,p_{1},Y}$
if we can show $f_{z}^{\left(r\right)}\in\FourierDecompSp{\CalQ}{p_{1}}Y$,
at least for suitable values of $z$.

\medskip{}

We will now show that this is indeed true. To this end, first recall
for $i\in I$ that $\varphi_{i}\gamma_{j}\not\equiv0$ can only hold
for $Q_{i}\cap P_{j}\neq\emptyset$, that is, if $j\in J_{i}$. Therefore,
\begin{align*}
\varphi_{i}\cdot f_{z}^{\left(r\right)}=\sum_{j\in J_{0}^{\left(r\right)}}\left[\zeta_{j}\cdot M_{z_{j}}\left(\varphi_{i}\cdot\gamma_{j}\right)\right] & =\sum_{j\in J_{i}\cap J_{0}^{\left(r\right)}}\left[\zeta_{j}\cdot M_{z_{j}}\left(\varphi_{i}\cdot\gamma_{j}\right)\right]\\
 & =\varphi_{i}\cdot\sum_{j\in J_{i}\cap J_{0}^{\left(r\right)}}\left[M_{z_{j}}\left(\zeta_{j}\cdot\gamma_{j}\right)\right]\\
 & =:\varphi_{i}\cdot F_{i}^{\left(r,z\right)}.
\end{align*}

Now, we invoke Corollary~\ref{cor:AsymptoticModulationBehaviour}
(with $f_{j}=\zeta_{j}\gamma_{j}\in\Schwartz\left(\R^{\dimension}\right)$
for $j\in M$) to obtain a family of modulations $\left(z_{j}\right)_{j\in M}$
(depending on the chosen $f_{j}$ and thus on $\zeta_{j}$, as well
as on $p_{1}$) such that we have for every subset $S\subset M$ the
estimate
\begin{align*}
\left\Vert \Fourier^{-1}\left[\,\smash{\sum_{j\in S}}\vphantom{\sum}M_{z_{j}}\left(\zeta_{j}\,\gamma_{j}\right)\,\right]\right\Vert _{L^{p_{1}}} & \leq2\left\Vert \left(\left\Vert \Fourier^{-1}\left(\zeta_{j}\,\gamma_{j}\right)\right\Vert _{L^{p_{1}}}\right)_{j\in S}\right\Vert _{\ell^{p_{1}}}\\
 & =2C_{1}^{\left(p_{1}\right)}\varepsilon^{\dimension\left(1-p_{1}^{-1}\right)}\cdot\left\Vert \left(\zeta_{j}\cdot\left|\det S_{j}\right|^{1-p_{1}^{-1}}\right)_{j\in S}\right\Vert _{\ell^{p_{1}}}\\
 & =2C_{1}^{\left(p_{1}\right)}\varepsilon^{\dimension\left(1-p_{1}^{-1}\right)}\cdot\left\Vert \left(c_{j}\right)_{j\in S}\right\Vert _{\ell^{p_{1}}}\:.
\end{align*}
In particular, if we set $z_{j}=0$ for $j\in J\setminus M$, we get
\begin{align*}
\left\Vert \Fourier^{-1}F_{i}^{\left(r,z\right)}\right\Vert _{L^{p_{1}}} & =\left\Vert \Fourier^{-1}\left[\,\smash{\sum_{j\in J_{i}\cap M^{\left(r\right)}}}\vphantom{\sum}M_{z_{j}}\left(\zeta_{j}\,\gamma_{j}\right)\,\right]\right\Vert _{L^{p_{1}}}\vphantom{\sum_{j\in J_{i}\cap M^{\left(r\right)}}}\\
 & \leq C_{2}\cdot\left\Vert \left(c_{j}\right)_{j\in J_{i}\cap M^{\left(r\right)}}\right\Vert _{\ell^{p_{1}}}\leq C_{2}\cdot\left\Vert \left(c_{j}\right)_{j\in J_{i}\cap J_{0}}\right\Vert _{\ell^{p_{1}}}=C_{2}\cdot\mu_{i}
\end{align*}
for all $i\in I$, with $C_{2}:=2C_{1}^{\left(p_{1}\right)}\varepsilon^{\dimension\left(1-p_{1}^{-1}\right)}$
and $\mu_{i}:=\left\Vert \left(c_{j}\right)_{j\in J_{i}\cap J_{0}}\right\Vert _{\ell^{p_{1}}}$.
Observe that $\mu:=\left(\mu_{i}\right)_{i\in I}\in Y$ with $\left\Vert \mu\right\Vert _{Y}=\left\Vert c\right\Vert _{V}$,
because of $c\in V=Y\left(\left[\ell^{p_{1}}\left(J_{i}\cap J_{0}\right)\right]_{i\in I}\right)$.

Now, as usual, there are two cases for $p_{1}$:

\begin{casenv}
\item $p_{1}\in\left[1,\infty\right]$. In this case, we set $C_{3}:=1$
and invoke Young's inequality to derive
\begin{align*}
\qquad\left\Vert \Fourier^{-1}\left(\varphi_{i}\cdot\smash{f_{z}^{\left(r\right)}}\,\right)\right\Vert _{L^{p_{1}}} & =\left\Vert \Fourier^{-1}\left(\varphi_{i}\cdot\smash{F_{i}^{\left(r,z\right)}}\,\right)\right\Vert _{L^{p_{1}}}\\
 & \leq\left\Vert \Fourier^{-1}\varphi_{i}\right\Vert _{L^{1}}\cdot\left\Vert \Fourier^{-1}F_{i}^{\left(r,z\right)}\right\Vert _{L^{p_{1}}}\\
 & \leq C_{\CalQ,\Phi,p_{1}}\cdot\left\Vert \Fourier^{-1}F_{i}^{\left(r,z\right)}\right\Vert _{L^{p_{1}}}\leq C_{2}C_{3}C_{\CalQ,\Phi,p_{1}}\cdot\mu_{i}\,.
\end{align*}
\item $p_{1}\in\left(0,1\right)$. In this case, $\CalQ=\left(T_{i}Q_{i}'+b_{i}\right)_{i\in I}$
is semi-structured and $\CalP_{J_{0}}$ is almost subordinate to $\CalQ$.
In fact, this is the only part of the proof where this is used. Because
of this subordinateness, if we set $k:=k\left(\CalP_{J_{0}},\CalQ\right)$,
then Lemma~\ref{lem:SubordinatenessImpliesWeakSubordination} implies
$\supp\gamma_{j}\subset P_{j}\subset Q_{i}^{\left(2k+2\right)\ast}$
for all $j\in J_{i}\cap J_{0}$, and thus $\supp F_{i}^{\left(r,z\right)}\subset Q_{i}^{\left(2k+2\right)\ast}$.
Since also $\supp\varphi_{i}\subset\overline{Q_{i}}\subset\overline{Q_{i}^{\left(2k+2\right)\ast}}$,
we can invoke Corollary~\ref{cor:QuasiBanachConvolutionSemiStructured}
to derive
\begin{align*}
\qquad\qquad\left\Vert \Fourier^{-1}\left(\varphi_{i}\cdot\smash{f_{z}^{\left(r\right)}}\,\right)\right\Vert _{L^{p_{1}}} & =\left\Vert \Fourier^{-1}\left(\varphi_{i}\cdot\smash{F_{i}^{\left(r,z\right)}}\,\right)\right\Vert _{L^{p_{1}}}\\
 & \leq C_{3}\cdot\left|\det T_{i}\right|^{\frac{1}{p_{1}}-1}\cdot\left\Vert \Fourier^{-1}\varphi_{i}\right\Vert _{L^{p_{1}}}\cdot\left\Vert \Fourier^{-1}F_{i}^{\left(r,z\right)}\right\Vert _{L^{p_{1}}}\\
 & \leq C_{3}C_{\CalQ,\Phi,p_{1}}\cdot\left\Vert \Fourier^{-1}F_{i}^{\left(r,z\right)}\right\Vert _{L^{p_{1}}}\leq C_{2}C_{3}C_{\CalQ,\Phi,p_{1}}\cdot\mu_{i}\,.
\end{align*}
for some constant $C_{3}=C_{3}\left(\CalQ,k,\dimension,p_{1}\right)$.
\end{casenv}
Taken together, both cases show that we have
\[
0\leq\varrho_{i}:=\left\Vert \Fourier^{-1}\left(\varphi_{i}\cdot\smash{f_{z}^{\left(r\right)}}\right)\right\Vert _{L^{p_{1}}}\leq C_{4}\cdot\mu_{i}\qquad\forall\,i\in I,
\]
with $C_{4}:=C_{2}C_{3}C_{\CalQ,\Phi,p_{1}}$. By solidity of $Y$
and because of $\mu\in Y$, this yields $\varrho=\left(\varrho_{i}\right)_{i\in I}\in Y$
and hence $f_{z}^{\left(r\right)}\in\FourierDecompSp{\CalQ}{p_{1}}Y$
with 
\[
\left\Vert \smash{f_{z}^{\left(r\right)}}\right\Vert _{\BAPUFourierDecompSp{\CalQ}{p_{1}}Y{\Phi}}=\left\Vert \varrho\right\Vert _{Y}\leq C_{4}\cdot\left\Vert \mu\right\Vert _{Y}=C_{4}\cdot\left\Vert c\right\Vert _{V}<\infty.
\]

\medskip{}

Now, we are almost done. As seen above, what we have shown implies
$f_{z}^{\left(r\right)}\in\CalD_{K}^{\CalQ,p_{1},Y}$ with the specific
choice of $z$ from above. Since $\iota$ is well-defined and bounded,
this implies $\iota f_{z}^{\left(r\right)}\in\FourierDecompSp{\CalP}{p_{2}}Z$.
Now, note for $j\in J_{0}$ that we have $\supp\phi_{j}\subset P_{j}\subset\CalO\cap\CalO'$,
so that our assumptions on $\iota$ imply for arbitrary $g\in\Schwartz\left(\R^{\dimension}\right)$
and $f\in\CalD_{K}^{\CalQ,p_{1},Y}$ that
\[
\left\langle \phi_{j}\cdot\iota f,\,g\right\rangle _{\Schwartz'}=\left\langle \iota f,\,\phi_{j}\,g\right\rangle _{\CalD'}=\left\langle f,\,\phi_{j}\,g\right\rangle _{\CalD'}=\left\langle \phi_{j}f,\,g\right\rangle _{\Schwartz'}\:,
\]
i.e.\@ $\phi_{j}\cdot\iota f=\phi_{j}\,f$. In particular, $\phi_{j}\cdot\iota f_{z}^{\left(r\right)}=\phi_{j}\,f_{z}^{\left(r\right)}$.
Furthermore, $\left(\phi_{j}\right)_{j\in J}$ is an $L^{p_{2}}$-bounded
system for $\CalP$, since we have $\supp\phi_{j}\subset P_{j}$ for
all $j\in J$ and $\sup_{j\in J}\left\Vert \Fourier^{-1}\phi_{j}\right\Vert _{L^{1}}=C_{2}^{\left(1\right)}$,
as well as
\[
\sup_{j\in J}\left(\left|\det S_{j}\right|^{\frac{1}{q}-1}\cdot\left\Vert \Fourier^{-1}\phi_{j}\right\Vert _{L^{q}}\right)=C_{2}^{\left(q\right)}\cdot\varepsilon^{\dimension\left(1-\frac{1}{q}\right)}\qquad\text{ for arbitrary }q\in\left(0,1\right).
\]
Hence, $C_{5}:=C_{\CalP,\left(\phi_{i}\right)_{i},p_{2}}=C\bigl(\varepsilon_{\CalP},C_{2}^{\left(\min\left\{ 1,p_{2}\right\} \right)},p_{2}\bigr)=C\left(\dimension,p_{2},\varepsilon_{\CalP}\right)$.

Next, recall from above that $f_{z}^{\left(r\right)}=\sum_{j\in M^{\left(r\right)}}M_{z_{j}}\left(\zeta_{j}\cdot\gamma_{j}\right)$,
with $\supp\gamma_{j}\subset P_{j}$ for all $j\in M^{\left(r\right)}$
and with $P_{j}^{\ast}\cap P_{\ell}^{\ast}=\emptyset$ for $j,\ell\in M^{\left(r\right)}\subset J^{\left(r\right)}$
with $j\neq\ell$. Finally, the ``coefficient sequence'' $\zeta$
is finitely supported. Thus, Lemma~\ref{lem:GeneralizedEasyNormEquivalenceFineLowerBound}
(applied to $\CalP$ instead of $\CalQ$, with $f=\iota f_{z}^{\left(r\right)}$
and $g=f_{z}^{\left(r\right)}$ as well as $I_{0}=M^{\left(r\right)}$)
yields a constant $C_{6}=C_{6}\left(\dimension,p_{2},\CalP,\vertiii{\Gamma_{\CalP}}_{Z\to Z}\right)>0$
satisfying
\begin{align*}
C_{5}C_{6}\cdot\left\Vert \iota\smash{f_{z}^{\left(r\right)}}\right\Vert _{\BAPUFourierDecompSp{\CalP}{p_{2}}Z{\Psi}} & \geq\left\Vert \left(\zeta_{j}\cdot\left\Vert \Fourier^{-1}\gamma_{j}\right\Vert _{L^{p_{2}}}\right)_{j\in M^{\left(r\right)}}\right\Vert _{Z|_{M^{\left(r\right)}}}\\
\left({\scriptstyle \text{since }\zeta_{j}=0\text{ for }j\in J_{0}^{\left(r\right)}\setminus M^{\left(r\right)}}\right) & =C_{1}^{\left(p_{2}\right)}\cdot\varepsilon^{\dimension\left(1-p_{2}^{-1}\right)}\cdot\left\Vert \left(\zeta_{j}\cdot\left|\det S_{j}\right|^{1-p_{2}^{-1}}\right)_{j\in J_{0}^{\left(r\right)}}\right\Vert _{Z|_{J_{0}^{\left(r\right)}}}\\
 & =C_{1}^{\left(p_{2}\right)}\cdot\varepsilon^{\dimension\left(1-p_{2}^{-1}\right)}\cdot\left\Vert \left(c_{j}\cdot\left|\det S_{j}\right|^{p_{1}^{-1}-p_{2}^{-1}}\right)_{j\in J_{0}^{\left(r\right)}}\right\Vert _{Z|_{J_{0}^{\left(r\right)}}}.
\end{align*}
In particular, part of the statement of Lemma~\ref{lem:GeneralizedEasyNormEquivalenceFineLowerBound}
is that the sequence $\left(c_{j}\cdot\left|\det S_{j}\right|^{p_{1}^{-1}-p_{2}^{-1}}\right)_{j\in J_{0}^{\left(r\right)}}$
is a member of $Z|_{J_{0}^{\left(r\right)}}$.

But since we have $J_{0}=\biguplus_{r=1}^{r_{0}}J_{0}^{\left(r\right)}$
and in view of the quasi-triangle inequality for $Z$, we get a constant
$C_{7}=C_{7}\left(C_{Z},r_{0}\right)=C_{7}\left(C_{Z},N_{\CalP}\right)$
with
\begin{align*}
\left\Vert \eta\left(c\right)\right\Vert _{Z_{\left|\det S_{j}\right|^{p_{1}^{-1}-p_{2}^{-1}}}} & \leq C_{7}\cdot\sum_{r=1}^{r_{0}}\left\Vert \left(\left|\det S_{j}\right|^{p_{1}^{-1}-p_{2}^{-1}}\cdot c_{j}\right)_{j\in J_{0}^{\left(r\right)}}\right\Vert _{Z|_{J_{0}^{\left(r\right)}}}\\
 & \leq\frac{C_{5}C_{6}C_{7}}{C_{1}^{\left(p_{2}\right)}\varepsilon^{\dimension\left(1-p_{2}^{-1}\right)}}\cdot\sum_{r=1}^{r_{0}}\left\Vert \iota\smash{f_{z}^{\left(r\right)}}\right\Vert _{\BAPUFourierDecompSp{\CalP}{p_{2}}Z{\Psi}}\\
 & \leq\vertiii{\iota}\frac{C_{5}C_{6}C_{7}}{C_{1}^{\left(p_{2}\right)}\varepsilon^{\dimension\left(1-p_{2}^{-1}\right)}}\cdot\sum_{r=1}^{r_{0}}\left\Vert \smash{f_{z}^{\left(r\right)}}\right\Vert _{\BAPUFourierDecompSp{\CalQ}{p_{1}}Y{\Phi}}\\
 & \leq\vertiii{\iota}\frac{r_{0}C_{4}C_{5}C_{6}C_{7}}{C_{1}^{\left(p_{2}\right)}\varepsilon^{\dimension\left(1-p_{2}^{-1}\right)}}\cdot\left\Vert c\right\Vert _{V}<\infty.
\end{align*}
This is precisely the desired estimate.

\medskip{}

The final statement of the theorem is an immediate consequence of
Lemmas~\ref{lem:FatouPropertyIsInherited} and \ref{lem:FinitelySupportedSequencesSufficeUnderFatouProperty}.
The required countability of $J$ for Lemma~\ref{lem:FinitelySupportedSequencesSufficeUnderFatouProperty}
holds, since $\CalP$ is a locally finite covering of the second countable
space $\CalO'$ (see Lemma~\ref{lem:PartitionCoveringNecessary})
with the additional property that $P_{j}\neq\emptyset$ for all $j\in J$,
by definition of an admissible covering.
\end{proof}
Above, we simply assumed $P_{j}\subset\CalO=\bigcup_{i\in I}Q_{i}$
for all $j\in J_{0}$, or—for $p_{1}\in\left(0,1\right)$—that $\CalP_{J_{0}}$
is almost subordinate to $\CalQ$. In our next theorem, we make the
``reverse'' assumption, i.e., we assume $\CalQ_{I_{0}}$ to be almost
subordinate to $\CalP$. As we will see, the techniques used in the
proof are similar, but the disjointization argument becomes slightly
more complex: We will apply the disjointization principle to $\CalP$
and then use this disjointization to construct a suitable disjointization
of the index set $I$ of the covering $\CalQ$. We remark that the
following theorem is a slightly improved version of the first part
of \cite[Theorem 5.3.6]{VoigtlaenderPhDThesis} from my PhD thesis.
\begin{thm}
\label{thm:BurnerNecessaryConditionFineInCoarse}Let $\emptyset\neq\CalO,\CalO'\subset\R^{\dimension}$
be open, let $p_{1},p_{2}\in\left(0,\infty\right]$ and let $\CalQ=\left(Q_{i}\right)_{i\in I}=\left(T_{i}Q_{i}'+b_{i}\right)_{i\in I}$
be a \emph{tight} semi-structured $L^{p_{1}}$-decomposition covering
of $\CalO$. Furthermore, let $\CalP=\left(P_{j}\right)_{j\in J}$
be an $L^{p_{2}}$-decomposition covering of $\CalO'$. Finally, let
$Y\subset\Compl^{I}$ and $Z\subset\Compl^{J}$ be $\CalQ$-regular
and $\CalP$-regular with triangle constants $C_{Y},C_{Z}\geq1$,
respectively.

Let $I_{0}\subset I$ and assume that $\CalQ_{I_{0}}:=\left(Q_{i}\right)_{i\in I_{0}}$
is almost subordinate to $\CalP$. Define
\[
K:=\bigcup_{i\in I_{0}}Q_{i}\subset\CalO\cap\CalO'
\]
and—with $\CalD_{K}^{\CalQ,p_{1},Y}$ as after equation~(\ref{eq:GeneralEmbeddingRequirement})—assume
that the identity map
\begin{equation}
\iota:\left(\CalD_{K}^{\CalQ,p_{1},Y},\left\Vert \mybullet\right\Vert _{\FourierDecompSp{\CalQ}{p_{1}}Y}\right)\to\FourierDecompSp{\CalP}{p_{2}}Z,f\mapsto f\label{eq:EmbeddingFineInCoarse}
\end{equation}
is well-defined and bounded.

Then the embedding
\[
\eta:\ell_{0}\left(I_{0}\right)\cap\left(Y|_{I_{0}}\right)_{\left|\det T_{i}\right|^{p_{2}^{-1}-p_{1}^{-1}}}\hookrightarrow Z\left(\left[\ell^{p_{2}}\left(I_{0}\cap I_{j}\right)\right]_{j\in J}\right)
\]
is well-defined and bounded with
\[
\vertiii{\iota}\leq C\cdot\vertiii{\iota}
\]
for some constant
\[
C=C\left(\dimension,p_{1},p_{2},k\left(\CalQ_{I_{0}},\CalP\right),C_{Z},\varepsilon_{\CalQ},\CalQ,C_{\CalQ,\Phi,p_{1}},\CalP,C_{\CalP,\Psi,p_{2}},\vertiii{\Gamma_{\CalQ}}_{Y\to Y},\vertiii{\Gamma_{\CalP}}_{Z\to Z}\right).
\]
Here, the $L^{p_{1}}$-BAPU $\Phi=\left(\varphi_{i}\right)_{i\in I}$
and the $L^{p_{2}}$-BAPU $\Psi=\left(\psi_{j}\right)_{j\in J}$ are
those which are used to calculate the (quasi)-norms on the two decomposition
spaces when calculating $\vertiii{\iota}$.

\medskip{}

Furthermore, the following hold:

\begin{itemize}[leftmargin=0.8cm]
\item If $Z$ satisfies the Fatou property (see Definition~\ref{def:FatouProperty}),
then $\eta$ remains well-defined and bounded (with the same estimate
for $\vertiii{\eta}$), even without intersecting with $\ell_{0}\left(I_{0}\right)$.\vspace{0.1cm}
\item If $p_{2}=\infty$, then $\eta$ remains well-defined and bounded
(with the same estimate for $\vertiii{\eta}$) even if the ``inner
norm'' $\ell^{p_{2}}\left(I_{0}\cap I_{j}\right)=\ell^{\infty}\left(I_{0}\cap I_{j}\right)$
is changed to $\ell^{\LowerExpo{p_{2}}}\!\left(I_{0}\cap I_{j}\right)=\ell^{1}\left(I_{0}\cap I_{j}\right)$.\qedhere
\end{itemize}
\end{thm}

\begin{rem}
\label{rem:ElementaryBurnerNecessaryConditionFineInCoarse}As for
Theorem~\ref{thm:BurnerNecessaryConditionCoarseInFine}, one can
obtain a generalization to the case $p_{1}\neq p_{2}$ of the estimate
$\left\Vert \delta_{j}\right\Vert _{Z}\lesssim\left\Vert \delta_{i}\right\Vert _{Y}$
for $Q_{i}\cap P_{j}$ that was shown in Lemma~\ref{lem:SimpleNecessaryCondition}.
Indeed, if $i_{0}\in I_{0}$ with $\delta_{i_{0}}\in Y$ and $j_{0}\in J$
with $Q_{i_{0}}\cap P_{j_{0}}\neq\emptyset$, then $i_{0}\in I_{0}\cap I_{j_{0}}$
and hence
\begin{align}
\left\Vert \delta_{j_{0}}\right\Vert _{Z}=\left\Vert \delta_{j_{0}}\right\Vert _{Z}\cdot\left\Vert \left(\delta_{i_{0}}\left(i\right)\right)_{i\in I_{0}\cap I_{j_{0}}}\right\Vert _{\ell^{p_{2}}} & \leq\left\Vert \delta_{i_{0}}\right\Vert _{Z\left(\left[\ell^{p_{2}}\left(I_{0}\cap I_{j}\right)\right]_{j\in J}\right)}\nonumber \\
 & \leq\vertiii{\eta}\cdot\left\Vert \delta_{i_{0}}\right\Vert _{\left(Y|_{I_{0}}\right)_{\left|\det T_{i}\right|^{p_{2}^{-1}-p_{1}^{-1}}}}\nonumber \\
 & =\vertiii{\eta}\cdot\left|\det T_{i_{0}}\right|^{p_{2}^{-1}-p_{1}^{-1}}\cdot\left\Vert \delta_{i_{0}}\right\Vert _{Y}.\label{eq:ElementaryBurnerNecessaryConditionFineInCoarse}
\end{align}
In case of $p_{1}=p_{2}$, this again yields the estimate $\left\Vert \delta_{j_{0}}\right\Vert _{Z}\lesssim\left\Vert \delta_{i_{0}}\right\Vert _{Y}$
from Lemma~\ref{lem:SimpleNecessaryCondition}.
\end{rem}

\begin{proof}[Proof of Theorem~\ref{thm:BurnerNecessaryConditionFineInCoarse}]
We begin by constructing several auxiliary objects and by deriving
a few properties of these objects; these properties and objects will
be helpful later in the proof.

For brevity, we write $V:=\left(Y|_{I_{0}}\right)_{\left|\det T_{i}\right|^{p_{2}^{-1}-p_{1}^{-1}}}\subset\Compl^{I_{0}}$
and $V_{0}:=\ell_{0}\left(I_{0}\right)\cap V$, as well as $k:=k\left(\CalQ_{I_{0}},\CalP\right)$.
Now, for each $i\in I_{0}$, we have $Q_{i}\subset P_{j_{i}}^{k\ast}\subset\CalO'$
for a suitable $j_{i}\in J$. Furthermore, $Q_{i}\neq\emptyset$ since
$\CalQ$ is tight. But this implies $Q_{i}\cap P_{j}\neq\emptyset$
for some $j\in J$, and in particular $i\in I_{0}\cap I_{j}$. All
in all, we have shown $I_{0}=\bigcup_{j\in J}\left(I_{0}\cap I_{j}\right)$,
so that $W:=Z\bigl(\left[\ell^{p_{2}}\left(I_{0}\cap I_{j}\right)\right]_{j\in J}\bigr)$
is a solid sequence space on $I_{0}$.

Now, choose $\left(\gamma_{i}\right)_{i\in I}$ as in the second part
of Lemma~\ref{lem:NormOfClusteredBAPUAndTestFunctionBuildingBlocks}
and set $\varepsilon:=\varepsilon_{\CalQ}$. Note that Lemma~\ref{lem:NormOfClusteredBAPUAndTestFunctionBuildingBlocks}
is applicable, since $\CalQ$ is tight semi-structured. By choice
of $\left(\gamma_{i}\right)_{i\in I}$, there are certain constants
${C_{1}=C_{1}\left(\dimension,p_{2},\varepsilon\right)>0}$ and $C_{2}=C_{2}\left(\dimension,p_{1},\varepsilon\right)>0$
with
\begin{equation}
\left\Vert \Fourier^{-1}\gamma_{i}\right\Vert _{L^{p_{2}}}=C_{1}^{-1}\cdot\left|\det T_{i}\right|^{1-p_{2}^{-1}}\quad\text{and}\quad\left\Vert \Fourier^{-1}\gamma_{i}\right\Vert _{L^{p_{1}}}=C_{2}\cdot\left|\det T_{i}\right|^{1-p_{1}^{-1}}\qquad\forall\,i\in I\,.\label{eq:NecessaryConditionFineInCoarseLpEstimatesOfBAPU}
\end{equation}

\medskip{}

Let $r_{0}:=N_{\CalP}^{2\left(2k+3\right)+1}=N_{\CalP}^{4k+7}$. With
this choice, the disjointization lemma (Lemma~\ref{lem:DisjointizationPrinciple})
yields a partition $J=\biguplus_{r=1}^{r_{0}}J^{\left(r\right)}$
such that we have $P_{j}^{\left(2k+3\right)\ast}\cap P_{\ell}^{\left(2k+3\right)\ast}=\emptyset$
for all $j,\ell\in J^{\left(r\right)}$ with $j\neq\ell$, for arbitrary
$r\in\underline{r_{0}}$.

A crucial property of the partition $J=\biguplus_{r=1}^{r_{0}}J^{\left(r\right)}$
is the following: We have 
\begin{equation}
\forall\,j,\ell\in J^{\left(r\right)}\text{ with }j\neq\ell:\qquad\left(I_{j}\cap I_{0}\right)\cap\left(I_{\ell}\cap I_{0}\right)=\emptyset.\label{eq:NecessaryFineInCoarseIntersectionIndexSetsDisjoint}
\end{equation}
Indeed, in case of $i\in I_{j}\cap I_{\ell}\cap I_{0}$, we would
have $Q_{i}\cap P_{j}\neq\emptyset\neq Q_{i}\cap P_{\ell}$—in particular,
$Q_{i}\neq\emptyset$—so that Lemma~\ref{lem:SubordinatenessImpliesWeakSubordination}
would yield $\emptyset\neq Q_{i}\subset P_{j}^{\left(2k+2\right)\ast}\cap P_{\ell}^{\left(2k+2\right)\ast}\subset P_{j}^{\left(2k+3\right)\ast}\cap P_{\ell}^{\left(2k+3\right)\ast}$,
in contradiction to $j,\ell\in J^{\left(r\right)}$ with $j\neq\ell$.

These considerations show that 
\[
I^{\left(r\right)}:=\biguplus_{j\in J^{\left(r\right)}}\left(I_{j}\cap I_{0}\right)\subset I_{0}
\]
is well-defined. As a simple consequence of this definition, we observe
\begin{equation}
\forall\,j\in J^{\left(r\right)}:\qquad I^{\left(r\right)}\cap I_{j}=I_{0}\cap I_{j}.\label{eq:ClusterDisjointizationIdentity}
\end{equation}

Next, for $i\in I_{0}$ and $j\in J$ with $Q_{i}\cap P_{j}\neq\emptyset$,
Lemma~\ref{lem:SubordinatenessImpliesWeakSubordination} yields $\supp\gamma_{i}\subset Q_{i}\subset P_{j}^{\left(2k+2\right)\ast}$,
and Lemma~\ref{lem:PartitionCoveringNecessary} implies $\psi_{j}^{\left(2k+3\right)\ast}\equiv1$
on $P_{j}^{\left(2k+2\right)\ast}$. Taken together, this establishes
\begin{equation}
\forall\,i\in I_{0}\,\forall\,j\in J\text{ with }Q_{i}\cap P_{j}\neq\emptyset:\qquad\gamma_{i}\cdot\psi_{j}^{\left(2k+3\right)\ast}=\gamma_{i}.\label{eq:NecessaryFineInCoarseBAPUMultiplication}
\end{equation}

\medskip{}

Now, we properly start the proof. Let $c=\left(c_{i}\right)_{i\in I_{0}}\in V_{0}=\ell_{0}\left(I_{0}\right)\cap\left(Y|_{I_{0}}\right)_{\left|\det T_{i}\right|^{p_{2}^{-1}-p_{1}^{-1}}}$
be arbitrary and extend this sequence to all of $I$ by setting $c_{i}:=0$
for $i\in I\setminus I_{0}$. Next, define $\zeta_{i}:=\left|\det T_{i}\right|^{p_{2}^{-1}-1}\cdot\left|c_{i}\right|$
for $i\in I$ and set $M:=\supp c=\supp\zeta\subset I_{0}$. With
these choices, define 
\begin{equation}
\zeta^{\left(r\right)}:=\zeta\cdot\Indicator_{I^{\left(r\right)}}\qquad\text{ and }\qquad g_{z}^{\left(r\right)}:=\sum_{i\in I^{\left(r\right)}}M_{z_{i}}\left(\zeta_{i}\cdot\gamma_{i}\right)=\sum_{i\in I_{0}}M_{z_{i}}\left(\smash{\zeta_{i}^{\left(r\right)}}\cdot\gamma_{i}\right)\label{eq:NecessaryFineInCoarseTestFunctionConstruction}
\end{equation}
for $r\in\underline{r_{0}}$ and any family $\left(z_{i}\right)_{i\in I_{0}}\in\left(\R^{\dimension}\right)^{I_{0}}$
of modulations. Observe that we have $\gamma_{i}\in\TestFunctionSpace{\CalO}$
with $\supp\gamma_{i}\subset Q_{i}\subset K$ for all $i\in I_{0}$.
Because of $\zeta_{i}^{\left(r\right)}=0$ for all but finitely many
$i\in I_{0}$, we thus have $g_{z}^{\left(r\right)}\in\CalD_{K}^{\CalQ,p_{1},Y}$,
as soon as we have shown $g_{z}^{\left(r\right)}\in\FourierDecompSp{\CalQ}{p_{1}}Y$.

But this is ensured by Lemma~\ref{lem:EasyNormEquivalenceFineCovering};
indeed, that lemma shows $g_{z}^{\left(r\right)}\in\FourierDecompSp{\CalQ}{p_{1}}Y$
and yields
\begin{align*}
\left\Vert \smash{g_{z}^{\left(r\right)}}\right\Vert _{\BAPUFourierDecompSp{\CalQ}{p_{1}}Y{\Phi}}=\left\Vert \sum_{i\in I}M_{z_{i}}\left(\smash{\zeta_{i}^{\left(r\right)}}\cdot\gamma_{i}\right)\right\Vert _{\BAPUFourierDecompSp{\CalQ}{p_{1}}Y{\Phi}} & \leq C_{3}\cdot\left\Vert \left(\zeta_{i}^{\left(r\right)}\cdot\left\Vert \Fourier^{-1}\gamma_{i}\right\Vert _{L^{p_{1}}}\right)_{i\in I}\right\Vert _{Y}\\
 & =C_{2}C_{3}\cdot\left\Vert \left(\zeta_{i}^{\left(r\right)}\cdot\left|\det T_{i}\right|^{1-p_{1}^{-1}}\right)_{i\in I}\right\Vert _{Y}\\
\left({\scriptstyle \text{since }\zeta_{i}=\left|c_{i}\right|\cdot\left|\det T_{i}\right|^{p_{2}^{-1}-1}\text{ and }\zeta_{i}=0\text{ for }i\notin I_{0}}\right) & \leq\vphantom{\sum^{T}}C_{2}C_{3}\cdot\left\Vert c\right\Vert _{\left(Y|_{I_{0}}\right)_{\left|\det T_{i}\right|^{p_{2}^{-1}-p_{1}^{-1}}}}<\infty
\end{align*}
for some constant $C_{3}=C_{3}\left(\dimension,p_{1},\CalQ,C_{\CalQ,\Phi,p_{1}},\vertiii{\Gamma_{\CalQ}}_{Y\to Y}\right)$.

\medskip{}

Since by assumption $\iota$ is well-defined and bounded, we get $g_{z}^{\left(r\right)}\in\FourierDecompSp{\CalP}{p_{2}}Z$
with
\begin{equation}
\left\Vert g_{z}^{\left(r\right)}\right\Vert _{\BAPUFourierDecompSp{\CalP}{p_{2}}Z{\Psi}}\leq C_{2}C_{3}\cdot\vertiii{\iota}\cdot\left\Vert c\right\Vert _{V}\,,\label{eq:NecessaryFineInCoarseUpperEstimate}
\end{equation}
for every family $z$ of modulations. Our next goal is to obtain a
lower bound on $\left\Vert \smash{g_{z}^{\left(r\right)}}\right\Vert _{\BAPUFourierDecompSp{\CalP}{p_{2}}Z{\Psi}}$.
The idea we use is similar to the proof of Lemma~\ref{lem:GeneralizedEasyNormEquivalenceFineLowerBound}.
But here, it is slightly easier to give a direct proof than to reduce
to the setting of that lemma.

Thus, let $r\in\underline{r_{0}}$, $j\in J^{\left(r\right)}$ and
$i\in I^{\left(r\right)}\subset I_{0}$ be arbitrary. There are two
cases:

\begin{casenv}
\item We have $i\in I_{j}$. In this case, equation~(\ref{eq:NecessaryFineInCoarseBAPUMultiplication})
implies
\begin{equation}
\psi_{j}^{\left(2k+3\right)\ast}\cdot M_{z_{i}}\left(\zeta_{i}\cdot\gamma_{i}\right)=M_{z_{i}}\left(\zeta_{i}\cdot\smash{\psi_{j}^{\left(2k+3\right)\ast}}\cdot\gamma_{i}\right)=M_{z_{i}}\left(\zeta_{i}\cdot\gamma_{i}\right).\label{eq:NecessaryConditionFineInCoarseBAPUReproduction}
\end{equation}
\item We have $i\notin I_{j}$. Because of $i\in I^{\left(r\right)}=\biguplus_{\ell\in J^{\left(r\right)}}\left(I_{0}\cap I_{\ell}\right)$,
this implies $i\in I_{\ell}$ for some $\ell\in J^{\left(r\right)}\setminus\left\{ j\right\} $.
But since $j,\ell\in J^{\left(r\right)}$ with $j\neq\ell$, we get
$P_{j}^{\left(2k+3\right)\ast}\cap P_{\ell}^{\left(2k+3\ast\right)}=\emptyset$
which entails $\psi_{j}^{\left(2k+3\right)\ast}\psi_{\ell}^{\left(2k+3\right)\ast}\equiv0$
and thus (using equation~(\ref{eq:NecessaryConditionFineInCoarseBAPUReproduction})
for $\ell$ instead of $j$) that
\[
\psi_{j}^{\left(2k+3\right)\ast}\cdot M_{z_{i}}\left(\zeta_{i}\cdot\gamma_{i}\right)=\psi_{j}^{\left(2k+3\right)\ast}\cdot\psi_{\ell}^{\left(2k+3\right)\ast}\cdot M_{z_{i}}\left(\zeta_{i}\cdot\gamma_{i}\right)\equiv0.
\]
\end{casenv}
In summary, these considerations imply
\begin{equation}
\psi_{j}^{\left(2k+3\right)\ast}\cdot M_{z_{i}}\left(\zeta_{i}\cdot\gamma_{i}\right)=\begin{cases}
M_{z_{i}}\left(\zeta_{i}\cdot\gamma_{i}\right), & \text{if }i\in I_{j},\\
0, & \text{if }i\notin I_{j}
\end{cases}\qquad\forall\,r\in\underline{r_{0}}\text{ and }j\in J^{\left(r\right)},\;i\in I^{\left(r\right)}\,.\label{eq:HyperIdenticalLocalizationIdentity}
\end{equation}

Next, combining Theorem~\ref{thm:BoundedControlSystemEquivalentQuasiNorm}
and Remark~\ref{rem:ClusteredBAPUYIeldsBoundedControlSystem} yields
$C_{4}=C_{4}\left(\CalP,p_{2},k,\dimension,C_{\CalP,\Psi,p_{2}},\vertiii{\Gamma_{\CalP}}_{Z\to Z}\right)>0$
such that
\[
\left\Vert \left(\left\Vert \Fourier^{-1}\left(\smash{\psi_{j}^{\left(2k+3\right)\ast}}\cdot g\right)\right\Vert _{L^{p_{2}}}\right)_{j\in J}\right\Vert _{Z}\leq C_{4}\cdot\left\Vert g\right\Vert _{\BAPUFourierDecompSp{\CalP}{p_{2}}Z{\Psi}}\qquad\forall\,g\in\FourierDecompSp{\CalP}{p_{2}}Z\,.
\]
In particular, we have $\left(\left\Vert \Fourier^{-1}\left(\smash{\psi_{j}^{\left(2k+3\right)\ast}}\cdot g\right)\right\Vert _{L^{p_{2}}}\right)_{j\in J}\in Z$
for each such $g$. By applying this to $g=g_{z}^{\left(r\right)}$,
we derive
\begin{align}
\left\Vert g_{z}^{\left(r\right)}\right\Vert _{\BAPUFourierDecompSp{\CalP}{p_{2}}Z{\Psi}} & \geq C_{4}^{-1}\cdot\left\Vert \left(\left\Vert \Fourier^{-1}\left(\psi_{j}^{\left(2k+3\right)\ast}\cdot g_{z}^{\left(r\right)}\right)\right\Vert _{L^{p_{2}}}\right)_{j\in J}\right\Vert _{Z}\nonumber \\
 & \geq C_{4}^{-1}\cdot\left\Vert \left(\left\Vert \Fourier^{-1}\left(\psi_{j}^{\left(2k+3\right)\ast}\cdot g_{z}^{\left(r\right)}\right)\right\Vert _{L^{p_{2}}}\right)_{j\in J^{\left(r\right)}}\right\Vert _{Z|_{J^{\left(r\right)}}}\nonumber \\
 & =C_{4}^{-1}\cdot\vphantom{\sum_{i\in I^{\left(r\right)}}}\left\Vert \left(\left\Vert \Fourier^{-1}\left[\,\smash{\sum_{i\in I^{\left(r\right)}}}\vphantom{\sum}\psi_{j}^{\left(2k+3\right)\ast}\cdot M_{z_{i}}\left(\zeta_{i}\cdot\gamma_{i}\right)\,\right]\right\Vert _{L^{p_{2}}}\right)_{\!j\in J^{\left(r\right)}}\right\Vert _{Z|_{J^{\left(r\right)}}}\nonumber \\
\left({\scriptstyle \text{eq. }\eqref{eq:HyperIdenticalLocalizationIdentity}}\right) & =C_{4}^{-1}\cdot\vphantom{\sum_{i\in I^{\left(r\right)}\cap I_{j}}}\left\Vert \left(\left\Vert \Fourier^{-1}\left[\,\smash{\sum_{i\in I^{\left(r\right)}\cap I_{j}}}\vphantom{\sum}M_{z_{i}}\left(\zeta_{i}\cdot\gamma_{i}\right)\,\right]\right\Vert _{L^{p_{2}}}\right)_{\!j\in J^{\left(r\right)}}\right\Vert _{Z|_{J^{\left(r\right)}}}.\label{eq:NecessaryFineInCoarseLowerBound}
\end{align}
At this point, we need to separate two cases, namely $p_{2}=\infty$
and $p_{2}<\infty$.

\medskip{}

We begin with the case $p_{2}=\infty$, in which we will even show
the (stronger) estimate with $\ell^{p_{2}}\left(I_{0}\cap I_{j}\right)$
replaced by $\ell^{1}\left(I_{0}\cap I_{j}\right)$, see the last
part of the theorem. Recall that we obtained the $\gamma_{i}$ from
Lemma~\ref{lem:NormOfClusteredBAPUAndTestFunctionBuildingBlocks}.
In particular, $\gamma_{i}\geq0$ for all $i\in I$, which yields
\[
\left(\Fourier^{-1}\gamma_{i}\right)\left(0\right)=\int_{\R^{\dimension}}\gamma_{i}\left(\xi\right)\,\d\xi=\left\Vert \gamma_{i}\right\Vert _{L^{1}}\geq\left\Vert \Fourier^{-1}\gamma_{i}\right\Vert _{L^{\infty}}=\left\Vert \Fourier^{-1}\gamma_{i}\right\Vert _{L^{p_{2}}}=C_{1}^{-1}\cdot\left|\det T_{i}\right|^{1-p_{2}^{-1}},
\]
as a consequence of the Riemann-Lebesgue lemma and of equation~(\ref{eq:NecessaryConditionFineInCoarseLpEstimatesOfBAPU}).
Thus, if we simply set $z_{i}=0$ for all $i\in I$, we get for each
$j\in J^{\left(r\right)}$ that
\begin{align*}
\vphantom{\sum_{i\in I^{\left(r\right)}\cap I_{j}}}\left\Vert \Fourier^{-1}\left[\,\smash{\sum_{i\in I^{\left(r\right)}\cap I_{j}}}\vphantom{\sum}M_{z_{i}}\left(\zeta_{i}\cdot\gamma_{i}\right)\,\right]\right\Vert _{L^{p_{2}}} & =\vphantom{\sum_{i\in I^{\left(r\right)}\cap I_{j}}}\left\Vert \,\smash{\sum_{i\in I^{\left(r\right)}\cap I_{j}}}\vphantom{\sum}\zeta_{i}\cdot\Fourier^{-1}\gamma_{i}\,\right\Vert _{L^{\infty}}\\
\left({\scriptstyle \text{function is continuous}}\right) & \geq\vphantom{\sum_{i\in I^{\left(r\right)}\cap I_{j}}}\left|\,\smash{\sum_{i\in I^{\left(r\right)}\cap I_{j}}}\vphantom{\sum}\zeta_{i}\cdot\left(\Fourier^{-1}\gamma_{i}\right)\left(0\right)\,\right|\\
\left({\scriptstyle \zeta_{i}\geq0}\right) & \geq C_{1}^{-1}\cdot\sum_{i\in I^{\left(r\right)}\cap I_{j}}\left[\zeta_{i}\cdot\left|\det T_{i}\right|^{1-p_{2}^{-1}}\right]\\
 & =C_{1}^{-1}\cdot\sum_{i\in I^{\left(r\right)}\cap I_{j}}\left[\left|c_{i}\right|\cdot\left|\det T_{i}\right|^{p_{2}^{-1}-1}\cdot\left|\det T_{i}\right|^{1-p_{2}^{-1}}\right]\\
 & =C_{1}^{-1}\cdot\left\Vert c\right\Vert _{\ell^{1}\left(I_{j}\cap I^{\left(r\right)}\right)}\\
\left({\scriptstyle \text{eq. }\eqref{eq:ClusterDisjointizationIdentity}\text{ and }j\in J^{\left(r\right)}}\right) & =C_{1}^{-1}\cdot\left\Vert c\right\Vert _{\ell^{1}\left(I_{j}\cap I_{0}\right)}.
\end{align*}

By solidity of $Z$ and in view of equation~(\ref{eq:NecessaryFineInCoarseLowerBound}),
we have thus shown $\left(\left\Vert c\right\Vert _{\ell^{1}\left(I_{j}\cap I_{0}\right)}\cdot\Indicator_{J^{\left(r\right)}}\left(j\right)\right)_{j\in J}\in Z$
with
\[
\left\Vert g_{z}^{\left(r\right)}\right\Vert _{\BAPUFourierDecompSp{\CalP}{p_{2}}Z{\Psi}}\geq\frac{1}{C_{1}C_{4}}\cdot\left\Vert \left(\left\Vert c\right\Vert _{\ell^{1}\left(I_{j}\cap I_{0}\right)}\right)_{j\in J^{\left(r\right)}}\right\Vert _{Z|_{J^{\left(r\right)}}}.
\]
Now, using the quasi-triangle inequality for $Z$ and recalling the
identity $J=\biguplus_{r=1}^{r_{0}}J^{\left(r\right)}$, we finally
get $c\in Z\bigl(\left[\ell^{1}\!\left(I_{j}\cap I_{0}\right)\right]_{j\in J}\bigr)$,
with
\begin{align*}
\left\Vert c\right\Vert _{Z\left(\left[\ell^{1}\left(I_{j}\cap I_{0}\right)\right]_{j\in J}\right)}=\left\Vert \left(\left\Vert c\right\Vert _{\ell^{1}\left(I_{j}\cap I_{0}\right)}\right)_{j\in J}\right\Vert _{Z} & \leq C_{5}\cdot\sum_{r=1}^{r_{0}}\left\Vert \left(\left\Vert c\right\Vert _{\ell^{1}\left(I_{j}\cap I_{0}\right)}\right)_{j\in J^{\left(r\right)}}\right\Vert _{Z|_{J^{\left(r\right)}}}\\
 & \leq C_{1}C_{4}C_{5}\cdot\sum_{r=1}^{r_{0}}\left\Vert g_{z}^{\left(r\right)}\right\Vert _{\BAPUFourierDecompSp{\CalP}{p_{2}}Z{\Psi}}\\
\left({\scriptstyle \text{eq. }\eqref{eq:NecessaryFineInCoarseUpperEstimate}}\right) & \leq C_{1}C_{2}C_{3}C_{4}C_{5}r_{0}\cdot\vertiii{\iota}\cdot\left\Vert c\right\Vert _{V}<\infty
\end{align*}
for some constant $C_{5}=C_{5}\left(r_{0},C_{Z}\right)=C_{5}\left(k,N_{\CalP},C_{Z}\right)$.
This establishes the boundedness of $\eta$ (with the ``inner norm''
$\ell^{p_{2}}\left(I_{0}\cap I_{j}\right)$ replaced by $\ell^{1}\left(I_{0}\cap I_{j}\right)$),
as desired. In particular, this establishes the last part of the theorem.

\medskip{}

Now, let us consider the general case. The proof given here is also
applicable for $p_{2}=\infty$, but yields a weaker conclusion (using
$\ell^{p_{2}}\left(I_{0}\cap I_{j}\right)$ instead of $\ell^{1}\left(I_{0}\cap I_{j}\right)$).
Here, we cannot simply choose all modulations $z_{i}=0$. Instead,
we apply Corollary~\ref{cor:AsymptoticModulationBehaviour} to the
family $\left(h_{i}\right)_{i\in M}$ defined by $h_{i}:=\zeta_{i}^{\left(r\right)}\cdot\gamma_{i}$
for $i\in M=\supp\zeta$. This yields a family $z^{\left(r\right)}=\left(\smash{z_{i}^{\left(r\right)}}\right)_{i\in M}$
of modulations so that
\begin{align}
\vphantom{\sum_{i\in S}}\left\Vert \Fourier^{-1}\left(\,\smash{\sum_{i\in S}}\,\vphantom{\sum}M_{z_{i}^{\left(r\right)}}\left(\smash{\zeta_{i}^{\left(r\right)}}\cdot\gamma_{i}\right)\,\right)\right\Vert _{L^{p_{2}}} & \geq\frac{1}{2}\cdot\left\Vert \left(\left\Vert \Fourier^{-1}\left(\smash{\zeta_{i}^{\left(r\right)}}\cdot\gamma_{i}\right)\right\Vert _{L^{p_{2}}}\right)_{i\in S}\right\Vert _{\ell^{p_{2}}}\nonumber \\
\left({\scriptstyle \text{eq. }\eqref{eq:NecessaryConditionFineInCoarseLpEstimatesOfBAPU}\text{ and }\zeta_{i}^{\left(r\right)}=0\text{ for }i\notin I^{\left(r\right)}}\right) & =\frac{1}{2C_{1}}\cdot\left\Vert \left(\left|c_{i}\right|\cdot\left|\det T_{i}\right|^{p_{2}^{-1}-1}\cdot\left|\det T_{i}\right|^{1-p_{2}^{-1}}\right)_{i\in S\cap I^{\left(r\right)}}\right\Vert _{\ell^{p_{2}}}\nonumber \\
 & =\frac{1}{2C_{1}}\cdot\left\Vert \left(c_{i}\right)_{i\in S\cap I^{\left(r\right)}}\right\Vert _{\ell^{p_{2}}}\qquad\forall\,S\subset M\,.\label{eq:NecessaryConditionAsymptoticModulationApplication}
\end{align}
For simplicity, set $z_{i}^{\left(r\right)}:=0$ for $i\in I\setminus M$.

With this choice of $z^{\left(r\right)}$, we have
\begin{align*}
\vphantom{\sum_{i\in I^{\left(r\right)}\cap I_{j}}}\left\Vert \Fourier^{-1}\left[\,\smash{\sum_{i\in I^{\left(r\right)}\cap I_{j}}}\vphantom{\sum}M_{z_{i}^{\left(r\right)}}\left(\zeta_{i}\cdot\gamma_{i}\right)\,\right]\right\Vert _{L^{p_{2}}} & =\vphantom{\sum_{i\in M\cap I_{j}}}\left\Vert \Fourier^{-1}\left[\,\smash{\sum_{i\in M\cap I_{j}}}\vphantom{\sum}M_{z_{i}^{\left(r\right)}}\left(\zeta_{i}^{\left(r\right)}\cdot\gamma_{i}\right)\,\right]\right\Vert _{L^{p_{2}}}\\
 & \geq\frac{1}{2C_{1}}\cdot\left\Vert \left(c_{i}\right)_{i\in M\cap I_{j}\cap I^{\left(r\right)}}\right\Vert _{\ell^{p_{2}}}\\
\left({\scriptstyle \text{since }c_{i}=0\text{ for }i\in I\setminus M}\right) & =\frac{1}{2C_{1}}\cdot\left\Vert \left(c_{i}\right)_{i\in I_{j}\cap I^{\left(r\right)}}\right\Vert _{\ell^{p_{2}}}\\
\left({\scriptstyle \text{eq. }\eqref{eq:ClusterDisjointizationIdentity}}\right) & =\frac{1}{2C_{1}}\cdot\left\Vert \left(c_{i}\right)_{i\in I_{j}\cap I_{0}}\right\Vert _{\ell^{p_{2}}}\qquad\forall\,j\in J^{\left(r\right)}\,.
\end{align*}
By solidity of $Z$ and using eq.\@ (\ref{eq:NecessaryFineInCoarseLowerBound}),
we thus get $\left(\left\Vert c\right\Vert _{\ell^{p_{2}}\left(I_{j}\cap I_{0}\right)}\cdot\Indicator_{J^{\left(r\right)}}\left(j\right)\right)_{j\in J}\in Z$
with
\[
\left\Vert g_{z^{\left(r\right)}}^{\left(r\right)}\right\Vert _{\BAPUFourierDecompSp{\CalP}{p_{2}}Z{\Psi}}\geq\frac{1}{2C_{1}C_{4}}\cdot\left\Vert \left(\left\Vert c\right\Vert _{\ell^{p_{2}}\left(I_{j}\cap I_{0}\right)}\right)_{j\in J^{\left(r\right)}}\right\Vert _{Z|_{J^{\left(r\right)}}}.
\]
Now, using the quasi-triangle inequality for $Z$ and $J=\biguplus_{r=1}^{r_{0}}J^{\left(r\right)}$,
we finally get $c\in Z\bigl(\left[\ell^{p_{2}}\!\left(I_{j}\cap I_{0}\right)\right]_{j\in J}\bigr)$,
with
\begin{align*}
\left\Vert c\right\Vert _{Z\left(\left[\ell^{p_{2}}\left(I_{j}\cap I_{0}\right)\right]_{j\in J}\right)}=\left\Vert \left(\left\Vert c\right\Vert _{\ell^{p_{2}}\left(I_{j}\cap I_{0}\right)}\right)_{j\in J}\right\Vert _{Z} & \leq C_{5}\cdot\sum_{r=1}^{r_{0}}\left\Vert \left(\left\Vert c\right\Vert _{\ell^{p_{2}}\left(I_{j}\cap I_{0}\right)}\right)_{j\in J^{\left(r\right)}}\right\Vert _{Z|_{J^{\left(r\right)}}}\\
 & \leq2C_{1}C_{4}C_{5}\cdot\sum_{r=1}^{r_{0}}\left\Vert g_{z^{\left(r\right)}}^{\left(r\right)}\right\Vert _{\BAPUFourierDecompSp{\CalP}{p_{2}}Z{\Psi}}\\
\left({\scriptstyle \text{eq. }\eqref{eq:NecessaryFineInCoarseUpperEstimate}}\right) & \leq2C_{1}C_{2}C_{3}C_{4}C_{5}r_{0}\cdot\vertiii{\iota}\cdot\left\Vert c\right\Vert _{V}<\infty
\end{align*}
for the same constant $C_{5}=C_{5}\left(k,N_{\CalP},C_{Z}\right)$
as for $p_{2}=\infty$. This completes the proof of the main statement
of the theorem.

\medskip{}

Finally, if $Z$ satisfies the Fatou property, it is an immediate
consequence of Lemmas \ref{lem:FatouPropertyIsInherited} and \ref{lem:FinitelySupportedSequencesSufficeUnderFatouProperty}
that 
\[
\widetilde{\eta}:\left(Y|_{I_{0}}\right)_{\left|\det T_{i}\right|^{p_{2}^{-1}-p_{1}^{-1}}}\hookrightarrow Z\left(\left[\ell^{p_{2}}\left(I_{j}\cap I_{0}\right)\right]_{j\in J}\right)
\]
is bounded if and only if $\eta$ is, with $\vertiii{\widetilde{\eta}}=\vertiii{\eta}$.
The required countability of $I_{0}\subset I$ for Lemma~\ref{lem:FinitelySupportedSequencesSufficeUnderFatouProperty}
holds, since $\CalQ$ is a locally finite covering of the second countable
space $\CalO$ (see Lemma~\ref{lem:PartitionCoveringNecessary})
with the additional property that $Q_{i}\neq\emptyset$ for all $i\in I$,
by definition of an admissible covering.
\end{proof}

\subsection{Further necessary conditions in case of \texorpdfstring{$p_{1}=p_{2}$}{p₁=p₂}}

\label{subsec:NecessaryForP1EqualP2}In the previous subsections,
we obtained necessary conditions for the existence (of slight generalizations)
of the embedding
\begin{equation}
\FourierDecompSp{\CalQ}{p_{1}}Y\hookrightarrow\FourierDecompSp{\CalP}{p_{2}}Z,\label{eq:KhinchinConditionIntroductionEmbedding}
\end{equation}
either assuming that (a subfamily of) $\CalQ=\left(T_{i}Q_{i}'+b_{i}\right)_{i\in I}$
is almost subordinate to $\CalP=\left(S_{j}P_{j}'+c_{j}\right)_{j\in J}$,
or vice versa. Slightly simplified, the derived necessary conditions
were $p_{1}\leq p_{2}$ and
\begin{equation}
Y_{\left|\det T_{i}\right|^{p_{2}^{-1}-p_{1}^{-1}}}\hookrightarrow Z\bigl(\left[\ell^{p_{2}}\left(I_{j}\right)\right]_{j\in J}\bigr)\label{eq:KhinchinConditionIntroductionQSubordinateToP}
\end{equation}
in case $\CalQ$ is almost subordinate to $\CalP$, and
\begin{equation}
Y\bigl(\left[\ell^{p_{1}}\left(J_{i}\right)\right]_{i\in I}\bigr)\hookrightarrow Z_{\left|\det S_{j}\right|^{p_{1}^{-1}-p_{2}^{-1}}}\label{eq:KhinchinConditionIntroductionPSubordinateToQ}
\end{equation}
in case $\CalP$ is almost subordinate to $\CalQ$. Indeed, at least
if $Z$ satisfies the Fatou property and if $\CalQ,\CalP$ are tight
semi-structured, these embeddings are consequences of Theorem~\ref{thm:BurnerNecessaryConditionFineInCoarse}
(with $I_{0}=I$) or of Theorem~\ref{thm:BurnerNecessaryConditionCoarseInFine}
(with $J_{0}=J$), respectively, while $p_{1}\leq p_{2}$ is a consequence
of Lemma~\ref{lem:SimpleNecessaryCondition}.

Conversely, if $\CalQ$ is almost subordinate to $\CalP$, then Corollary~\ref{cor:EmbeddingFineIntoCoarse}
shows that a \emph{sufficient} condition for the existence of the
embedding~(\ref{eq:KhinchinConditionIntroductionEmbedding}) is $p_{1}\leq p_{2}$
and
\[
Y_{\left|\det T_{i}\right|^{p_{2}^{-1}-p_{1}^{-1}}}\hookrightarrow Z\Bigl(\bigl[\ell^{\LowerExpo{p_{2}}}\!\bigl(I_{j}\bigr)\bigr]_{j\in J}\Bigr).
\]
If instead $\CalP$ is almost subordinate to $\CalQ$, then Corollary~\ref{cor:EmbeddingCoarseIntoFine}
shows—at least for $p_{1}\geq1$—that $p_{1}\leq p_{2}$ and
\[
Y\Bigl(\bigl[\ell^{\UpperExpo{p_{1}}}\!\bigl(J_{i}\bigr)\bigr]_{i\in I}\Bigr)\hookrightarrow Z_{\left|\det S_{j}\right|^{p_{1}^{-1}-p_{2}^{-1}}}
\]
are sufficient for the existence of the embedding~(\ref{eq:KhinchinConditionIntroductionEmbedding}).

In summary, we thus achieve a \emph{complete characterization} of
the existence of the embedding~(\ref{eq:KhinchinConditionIntroductionEmbedding})
in terms of embeddings for discrete (nested) sequence spaces in the
following cases:
\begin{itemize}[leftmargin=0.8cm]
\item If $\CalQ$ is almost subordinate to $\CalP$ and if $p_{2}\in\left(0,2\right]\cup\left\{ \infty\right\} $.
Indeed, for $p_{2}\in\left(0,2\right]$, we have $\LowerExpo{p_{2}}=p_{2}$,
so that both conditions coincide. Finally, in case of $p_{2}=\infty$,
Theorem~\ref{thm:BurnerNecessaryConditionFineInCoarse} shows that
we can replace $\ell^{p_{2}}\left(I_{j}\right)$ by $\ell^{1}\left(I_{j}\right)=\ell^{\LowerExpo{\infty}}\!\left(I_{j}\right)$
in the necessary condition~(\ref{eq:KhinchinConditionIntroductionQSubordinateToP}).
\item If $\CalP$ is almost subordinate to $\CalQ$ and if $p_{1}\in\left[2,\infty\right]$,
since this implies $p_{1}\geq1$, so that the sufficient condition
from above is applicable and since we have $\UpperExpo{p_{1}}=p_{1}$
in this range.
\end{itemize}
In the remaining cases, there is a gap between the necessary and sufficient
conditions. Under additional hypotheses, this gap can be closed, as
the next subsection shows.

In this subsection, our aim is slightly different: Note that we always
have $\UpperExpo{p_{1}}\geq2$ and $\LowerExpo{p_{1}}\leq2$. Thus,
there seems to be something special about the ``critical'' exponent
$2$. In this subsection, we will show that this is indeed the case.
Assuming $p_{1}=p_{2}$, we will show that
\[
Y\hookrightarrow Z\bigl(\left[\ell^{2}\left(I_{j}\right)\right]_{j\in J}\bigr)\qquad\text{ or }\qquad Y\bigl(\left[\ell^{2}\left(J_{i}\right)\right]_{i\in I}\bigr)\hookrightarrow Z
\]
are always necessary conditions for the existence of the embedding~(\ref{eq:KhinchinConditionIntroductionEmbedding}),
respectively, if $\CalQ$ is almost subordinate to $\CalP$ or vice
versa. This strengthens the necessary conditions (\ref{eq:KhinchinConditionIntroductionQSubordinateToP})
and (\ref{eq:KhinchinConditionIntroductionPSubordinateToQ}), while
still not fully reaching the sufficient conditions. Thus, we \emph{tighten}
the gap, but fail to close it completely. As remarked above, closing
the gap completely—but only under additional assumptions—is the goal
of the next subsection.

We begin with the following theorem which is analogous to Theorem~\ref{thm:BurnerNecessaryConditionCoarseInFine}
from the previous subsection. A curious feature of the present result
is that we do not need to assume $\CalQ$ or $\CalP$ to be tight—not
even semi-structured. Furthermore, we do not need to assume $\CalP$
to be almost subordinate to $\CalQ$, not even for $p\in\left(0,1\right)$,
in contrast to Theorem~\ref{thm:BurnerNecessaryConditionCoarseInFine}.
\begin{thm}
\label{thm:KhinchinNecessaryCoarseInFine}Let $\emptyset\neq\CalO,\CalO'\subset\R^{\dimension}$
be open, let $p\in\left(0,\infty\right]$ and let $\CalQ=\left(Q_{i}\right)_{i\in I}$
and $\CalP=\left(P_{j}\right)_{j\in J}$ be two $L^{p}$-decomposition
coverings of $\CalO$ and $\CalO'$, respectively. Finally, let $Y\subset\Compl^{I}$
be $\CalQ$-regular and $Z\subset\Compl^{J}$ be $\CalP$-regular,
with triangle constants $C_{Y}\geq1$ and $C_{Z}\geq1$, respectively.

Let $J_{0}\subset J$ and assume that $P_{j}$ is open and $P_{j}\subset\CalO$
for each $j\in J_{0}$. Define $K:=\bigcup_{j\in J_{0}}P_{j}\subset\CalO\cap\CalO'$.

If—with $\CalD_{K}^{\CalQ,p,Y}$ as after equation~(\ref{eq:GeneralEmbeddingRequirement})—the
identity map
\begin{equation}
\iota:\left(\CalD_{K}^{\CalQ,p,Y},\left\Vert \mybullet\right\Vert _{\FourierDecompSp{\CalQ}pY}\right)\rightarrow\FourierDecompSp{\CalP}pZ,f\mapsto f\label{eq:KhinchinEmbeddingCoarseInFine}
\end{equation}
is well-defined and bounded, then so is the embedding
\[
\eta:\ell_{0}\left(J_{0}\right)\cap Y\Bigl(\bigl[\ell^{\max\left\{ 2,p\right\} }\left(J_{0}\cap J_{i}\right)\bigr]_{i\in I}\Bigr)\hookrightarrow Z|_{J_{0}}\:,
\]
with $\vertiii{\eta}\leq C\cdot\vertiii{\iota}$ for some constant
\[
C=C\left(\dimension,p,C_{Z},\CalQ,\CalP,C_{\CalQ,\Phi,p},C_{\CalP,\Psi,p},\vertiii{\Gamma_{\CalQ}}_{Y\to Y},\vertiii{\Gamma_{\CalP}}_{Z\to Z}\right)>0.
\]
Here, the $L^{p}$-BAPUs $\Phi=\left(\varphi_{i}\right)_{i\in I}$
and $\Psi=\left(\psi_{j}\right)_{j\in J}$ are those which are used
to compute the (quasi)-norms on the respective decomposition spaces
when computing the norm $\vertiii{\iota}$.

Finally, if $Z$ satisfies the Fatou property, the same statement
as above also holds for the embedding $\eta:Y\bigl(\left[\ell^{\max\left\{ 2,p\right\} }\left(J_{0}\cap J_{i}\right)\right]_{i\in I}\bigr)\hookrightarrow Z$;
that is, without restricting to $\ell_{0}\left(J_{0}\right)$. 
\end{thm}

\begin{rem*}
Note that we always have norm-decreasing embeddings $\ell^{2}\left(J_{0}\cap J_{i}\right)\hookrightarrow\ell^{\max\left\{ 2,p\right\} }\left(J_{0}\cap J_{i}\right)$
and $\ell^{p}\left(J_{0}\cap J_{i}\right)\hookrightarrow\ell^{\max\left\{ 2,p\right\} }\left(J_{0}\cap J_{i}\right)$.
Thus, it is easy to see that the theorem implies in particular that
\[
\eta_{1}:\ell_{0}\left(J_{0}\right)\cap Y\left(\left[\ell^{p}\left(J_{0}\cap J_{i}\right)\right]_{i\in I}\right)\hookrightarrow Z|_{J_{0}}\quad\text{and}\quad\eta_{2}:\ell_{0}\left(J_{0}\right)\cap Y\bigl(\left[\ell^{2}\left(J_{0}\cap J_{i}\right)\right]_{i\in I}\bigr)\hookrightarrow Z|_{J_{0}}
\]
are well-defined and bounded with $\vertiii{\smash{\eta_{\ell}}}\leq\vertiii{\eta}$
for $\ell\in\left\{ 1,2\right\} $.
\end{rem*}
\begin{proof}
For the sake of brevity, set $s:=\max\left\{ 2,p\right\} \in\left[2,\infty\right]$
and $V:=Y\left(\left[\ell^{s}\left(J_{0}\cap J_{i}\right)\right]_{i\in I}\right)$,
as well as $V_{0}:=\ell_{0}\left(J_{0}\right)\cap V$. Furthermore,
fix a nonzero function $\gamma\in\TestFunctionSpace{B_{1}\left(0\right)}$
for the rest of the proof. Finally, let $r_{0}:=N_{\CalP}^{3}$.
With this choice, Lemma~\ref{lem:DisjointizationPrinciple} yields
a partition $J=\biguplus_{r=1}^{r_{0}}J_{\left(r\right)}$ such that
$P_{j}^{\ast}\cap P_{\ell}^{\ast}=\emptyset$ holds for all $j,\ell\in J_{\left(r\right)}$
with $j\neq\ell$ and arbitrary $r\in\underline{r_{0}}$.

Note that we have $\emptyset\neq P_{j}\subset\CalO$ for all $j\in J_{0}$
(by definition of an admissible covering, all sets of $\CalP$ are
nonempty), so that we get $Q_{i}\cap P_{j}\neq\emptyset$ for some
$i\in I$, since $\CalQ$ covers $I$. Hence, $J_{0}=\bigcup_{i\in I}\left(J_{0}\cap J_{i}\right)$,
so that $V$ is a solid sequence space over $J_{0}$. 

Let $c=\left(c_{j}\right)_{j\in J_{0}}\in V_{0}$ be arbitrary and
let $M:=\supp c$, which is a finite subset of $J_{0}$. Since $M\subset J_{0}$
is finite and since each set $P_{j}$ with $j\in J_{0}$ is open,
there is some $\delta=\delta\left(c\right)>0$ and for each $j\in M$
some $\xi_{j}\in\R^{\dimension}$ (possibly depending on $c$) with
$B_{\delta}\left(\xi_{j}\right)\subset P_{j}\subset\CalO$. Now, Lemma~\ref{lem:PartitionCoveringNecessary}
shows that $\CalQ^{\circ}=\left(Q_{i}^{\circ}\right)_{i\in I}$ covers
$\CalO$, so that for each $j\in M$ there is some $\ell_{j}\in I$
with $\xi_{j}\in Q_{\ell_{j}}^{\circ}$. This yields some $\varepsilon_{j}\in\left(0,\delta\right)$
satisfying $B_{\varepsilon_{j}}\left(\xi_{j}\right)\subset Q_{\ell_{j}}$.
Again by finiteness of $M$, we see that $\varepsilon:=\min_{j\in M}\varepsilon_{j}\in\left(0,\delta\right)$
is well-defined with $B_{\varepsilon}\left(\xi_{j}\right)\subset P_{j}\cap Q_{\ell_{j}}$
for all $j\in M$.

It is important to observe that the quantities $\left(\xi_{j}\right)_{j\in M}$
and $\varepsilon>0$ depend highly on the sequence $c=\left(c_{j}\right)_{j\in J_{0}}$.
Nevertheless, all constants $C_{1},C_{2},\dots$ chosen in the remainder
of the proof will be independent of the sequence $c$. Finally, we
will see that the precise choice of $\xi_{j}$ will be immaterial
for the resulting estimates and that all occurrences of $\varepsilon$
will cancel in the end.

For $j\in M$, let $\gamma_{j}:=L_{\xi_{j}}\left(\gamma\left(\varepsilon^{-1}\mybullet\right)\right)$
and note $\supp\gamma_{j}\subset B_{\varepsilon}\left(\xi_{j}\right)\subset P_{j}\cap Q_{\ell_{j}}$,
as well as
\begin{equation}
\Fourier^{-1}\gamma_{j}=\varepsilon^{\dimension}\cdot M_{\xi_{j}}\left[\left(\Fourier^{-1}\gamma\right)\left(\varepsilon\mybullet\right)\right]\,,\label{eq:KhinchinNecessaryConditionCoarseInFineFourierTransformOfPieces}
\end{equation}
which implies
\begin{equation}
\left\Vert \Fourier^{-1}\gamma_{j}\right\Vert _{L^{p}}=C_{1}\cdot\varepsilon^{\dimension\left(1-\frac{1}{p}\right)}\qquad\text{ for }\qquad C_{1}=C_{1}\left(\dimension,p\right):=\left\Vert \Fourier^{-1}\gamma\right\Vert _{L^{p}}.\label{eq:KhinchinNecessaryConditionCoarseInFineNormOfPieces}
\end{equation}

\medskip{}

Now, we make some technical observations which will become important
in the remainder of the proof. First, for $\ell\in I$, we define
\[
M^{\left(\ell\right)}:=\left\{ j\in M\with\ell_{j}=\ell\right\} 
\]
and note $M=\biguplus_{\ell\in I}M^{\left(\ell\right)}$, as well
as $j\in M^{\left(\ell_{j}\right)}$ for all $j\in M$. Next, suppose
that we have
\[
f=\sum_{j\in M}M_{z_{j}}\left(d_{j}\cdot\gamma_{j}\right)\quad\text{for certain sequences }\left(d_{j}\right)_{j\in M}\in\Compl^{M}\text{ and }\left(z_{j}\right)_{j\in M}\in\left(\R^{\dimension}\right)^{M}\,.
\]
We are interested in a simplified form of $\varphi_{i}\cdot f$, for
arbitrary $i\in I$.

To this end, note for $j\in M$ that we have $\varphi_{i}\,\gamma_{j}\equiv0$
unless $\emptyset\neq Q_{i}\cap\supp\gamma_{j}\subset Q_{i}\cap Q_{\ell_{j}}$.
But in this case, we get $\ell_{j}\in i^{\ast}$ and thus $j\in M^{\left(\ell_{j}\right)}\subset\bigcup_{\ell\in i^{\ast}}M^{\left(\ell\right)}$.
Thus, since the sequence $\left(M^{\left(\ell\right)}\right)_{\ell\in I}$
is pairwise disjoint, we finally get
\begin{align}
\varphi_{i}\cdot f=\varphi_{i}\cdot\sum_{j\in M}M_{z_{j}}\left(d_{j}\cdot\gamma_{j}\right)=\sum_{j\in M}d_{j}\cdot M_{z_{j}}\left(\varphi_{i}\cdot\gamma_{j}\right) & =\sum_{j\in\bigcup_{\ell\in i^{\ast}}M^{\left(\ell\right)}}d_{j}\cdot M_{z_{j}}\left(\varphi_{i}\cdot\gamma_{j}\right)\nonumber \\
 & =\sum_{\ell\in i^{\ast}}\,\sum_{j\in M^{\left(\ell\right)}}d_{j}\cdot M_{z_{j}}\left(\varphi_{i}\cdot\gamma_{j}\right)\nonumber \\
\left({\scriptstyle \text{with }f^{\left(\ell\right)}:=\sum_{j\in M^{\left(\ell\right)}}M_{z_{j}}\left(d_{j}\cdot\gamma_{j}\right)}\right) & =\sum_{\ell\in i^{\ast}}\left[\varphi_{i}\cdot\smash{f^{\left(\ell\right)}}\right].\label{eq:KhinchinNecessaryConditionCoarseInFineLocalizationIdentity}
\end{align}

Now, we want to derive an estimate for $\left\Vert \Fourier^{-1}\left(\varphi_{i}\cdot f\right)\right\Vert _{L^{p}}$.
In case of $p\in\left[1,\infty\right]$, this is straightforward:
Using the triangle inequality for $L^{p}$ and Young's inequality
$L^{1}\ast L^{p}\hookrightarrow L^{p}$, we get
\begin{align*}
\vphantom{\sum_{j\in M}}\left\Vert \Fourier^{-1}\left(\,\varphi_{i}\cdot\smash{\sum_{j\in M}}\vphantom{\sum}M_{z_{j}}\left(d_{j}\cdot\gamma_{j}\right)\,\right)\right\Vert _{L^{p}} & \leq\sum_{\ell\in i^{\ast}}\left\Vert \Fourier^{-1}\left[\,\varphi_{i}\cdot\smash{\sum_{j\in M^{\left(\ell\right)}}}\vphantom{\sum}M_{z_{j}}\left(d_{j}\cdot\gamma_{j}\right)\,\right]\right\Vert _{L^{p}}\vphantom{\sum}\\
 & \leq\sum_{\ell\in i^{\ast}}\left\Vert \Fourier^{-1}\varphi_{i}\right\Vert _{L^{p}}\cdot\left\Vert \Fourier^{-1}\left[\,\smash{\sum_{j\in M^{\left(\ell\right)}}}\vphantom{\sum}M_{z_{j}}\left(d_{j}\cdot\gamma_{j}\right)\,\right]\right\Vert _{L^{p}}\vphantom{\sum_{j\in M^{\left(\ell\right)}}^{T}}\\
 & \leq C_{\CalQ,\Phi,p}\cdot\sum_{\ell\in i^{\ast}}\left\Vert \Fourier^{-1}\left[\,\smash{\sum_{j\in M^{\left(\ell\right)}}}\vphantom{\sum}M_{z_{j}}\left(d_{j}\cdot\gamma_{j}\right)\,\right]\right\Vert _{L^{p}}\vphantom{\sum_{j\in M^{\left(\ell\right)}}}.
\end{align*}
In case of $p\in\left(0,1\right)$, the argument is slightly more
involved: For $\ell\in i^{\ast}$ and $j\in M^{\left(\ell\right)}$,
we have $\ell_{j}=\ell\in i^{\ast}$, and hence 
\[
\supp\gamma_{j}\subset P_{j}\cap Q_{\ell_{j}}\subset Q_{\ell}\subset Q_{i}^{\ast}\subset\overline{Q_{i}^{\ast}},
\]
which yields $\supp f^{\left(\ell\right)}\subset\overline{Q_{i}^{\ast}}$,
since $j\in M^{\left(\ell\right)}$ was arbitrary. Furthermore, $\supp\varphi_{i}\subset\overline{Q_{i}}\subset\overline{Q_{i}^{\ast}}$.
Finally, because $\CalQ$ is an $L^{p}$-decomposition covering of
$\CalO$ and since $p\in\left(0,1\right)$, $\CalQ=\left(T_{i}Q_{i}'+b_{i}\right)_{i\in I}$
is semi-structured. All in all, we can apply Corollary~\ref{cor:QuasiBanachConvolutionSemiStructured},
which yields a constant $C_{2}=C_{2}\left(\CalQ,\dimension,p\right)$
satisfying
\[
\left\Vert \Fourier^{-1}\left(\varphi_{i}\cdot\smash{f^{\left(\ell\right)}}\right)\right\Vert _{L^{p}}\leq C_{2}\cdot\left|\det T_{i}\right|^{\frac{1}{p}-1}\cdot\left\Vert \Fourier^{-1}\varphi_{i}\right\Vert _{L^{p}}\cdot\left\Vert \Fourier^{-1}\smash{f^{\left(\ell\right)}}\right\Vert _{L^{p}}\leq C_{2}C_{\CalQ,\Phi,p}\cdot\left\Vert \Fourier^{-1}\smash{f^{\left(\ell\right)}}\right\Vert _{L^{p}}.
\]
Finally, using the uniform estimate $\left|i^{\ast}\right|\leq N_{\CalQ}$
and the quasi-triangle inequality for $L^{p}$, we obtain a constant
$C_{3}=C_{3}\left(N_{\CalQ},p\right)$ with
\begin{align}
\left\Vert \Fourier^{-1}\left(\,\vphantom{\sum}\varphi_{i}\cdot\smash{\sum_{j\in M}}M_{z_{j}}\left(d_{j}\cdot\gamma_{j}\right)\,\right)\right\Vert _{L^{p}} & \overset{\text{eq. }\eqref{eq:KhinchinNecessaryConditionCoarseInFineLocalizationIdentity}}{=}\left\Vert \Fourier^{-1}\left(\,\smash{\sum_{\ell\in i^{\ast}}}\vphantom{\sum}\varphi_{i}\cdot f^{\left(\ell\right)}\,\right)\right\Vert _{L^{p}}\vphantom{\sum_{\ell\in i^{\ast}}}\nonumber \\
 & \leq C_{3}\cdot\sum_{\ell\in i^{\ast}}\left\Vert \Fourier^{-1}\left(\varphi_{i}\cdot\smash{f^{\left(\ell\right)}}\right)\right\Vert _{L^{p}}\nonumber \\
 & \leq C_{2}C_{3}C_{\CalQ,\Phi,p}\cdot\sum_{\ell\in i^{\ast}}\left\Vert \Fourier^{-1}\smash{f^{\left(\ell\right)}}\right\Vert _{L^{p}}\nonumber \\
 & =C_{2}C_{3}C_{\CalQ,\Phi,p}\cdot\sum_{\ell\in i^{\ast}}\left\Vert \Fourier^{-1}\left[\,\smash{\sum_{j\in M^{\left(\ell\right)}}}\vphantom{\sum}M_{z_{j}}\left(d_{j}\cdot\gamma_{j}\right)\,\right]\right\Vert _{L^{p}}\vphantom{\sum_{j\in M^{\left(\ell\right)}}}.\label{eq:KhinchinNecessaryConditionGenericUpperEstimate}
\end{align}
Thus, if we set $C_{2}:=C_{3}:=1$ for $p\in\left[1,\infty\right]$,
we see that the previous estimate remains true for all $p\in\left(0,\infty\right]$.

\medskip{}

Now, we divide the proof into two parts. For the first part, we assume
$p\leq2<\infty$, so that we get $s=\max\left\{ 2,p\right\} =2$.
In the second part, we will consider the case $p>2$, where $s=p$.

\medskip{}

\textbf{Part~1}: Here, we have $p\leq2<\infty$. Let us fix some
$r\in\underline{r_{0}}$. Now, for arbitrary $\ell\in I$, we consider
the random variable $\omega^{\left(\ell\right)}=\omega^{\left(\ell,r\right)}=\left(\smash{\omega_{j}^{\left(\ell\right)}}\right)_{j\in M^{\left(\ell\right)}\cap J_{\left(r\right)}}\in\left\{ \pm1\right\} ^{M^{\left(\ell\right)}\cap J_{\left(r\right)}}$,
which we take to be uniformly distributed in $\left\{ \pm1\right\} ^{M^{\left(\ell\right)}\cap J_{\left(r\right)}}$.
The expectation in the following computation is to be understood with
respect to $\omega^{\left(\ell\right)}$. Let $C_{4}=C_{4}\left(p\right)$
be the constant supplied by Khintchine's inequality (Theorem~\ref{thm:KhintchineInequality}).
In the following calculation, we can interchange the expectation with
the integral either by Fubini's theorem, or by noting that the expectation
is just a finite sum. Therefore,
\begin{align}
\vphantom{\sum_{j\in M^{\left(\ell\right)}\cap J_{\left(r\right)}}}\mathbb{E}\left\Vert \Fourier^{-1}\!\left[\,\smash{\sum_{j\in M^{\left(\ell\right)}\cap J_{\left(r\right)}}}\vphantom{\sum}\!\!\!c_{j}\cdot\omega_{j}^{\left(\ell\right)}\cdot\gamma_{j}\,\right]\right\Vert _{L^{p}}^{p} & =\mathbb{E}\int_{\R^{\dimension}}\left|\,\smash{\sum_{j\in M^{\left(\ell\right)}\cap J_{\left(r\right)}}}\vphantom{\sum}c_{j}\cdot\omega_{j}^{\left(\ell\right)}\cdot\left(\Fourier^{-1}\gamma_{j}\right)\left(x\right)\,\right|^{p}\,\d x\vphantom{\sum_{j\in M^{\left(\ell\right)}\cap J_{\left(r\right)}}}\nonumber \\
\left({\scriptstyle \text{eq. }\eqref{eq:KhinchinNecessaryConditionCoarseInFineFourierTransformOfPieces}}\right) & =\varepsilon^{\dimension p}\cdot\int_{\R^{\dimension}}\mathbb{E}\left|\,\smash{\sum_{j\in M^{\left(\ell\right)}\cap J_{\left(r\right)}}}\vphantom{\sum}\omega_{j}^{\left(\ell\right)}\cdot c_{j}\cdot e^{2\pi i\left\langle \xi_{j},x\right\rangle }\cdot\left(\Fourier^{-1}\gamma\right)\left(\varepsilon x\right)\,\right|^{p}\,\d x\vphantom{\sum_{j\in M^{\left(\ell\right)}\cap J_{\left(r\right)}}}\nonumber \\
\left({\scriptstyle \text{by Khintchine's ineq.}}\right) & \leq C_{4}\varepsilon^{\dimension p}\cdot\int_{\R^{\dimension}}\!\left(\,\smash{\sum_{j\in M^{\left(\ell\right)}\cap J_{\left(r\right)}}}\vphantom{\sum}\!\!\left|c_{j}\cdot e^{2\pi i\left\langle \xi_{j},x\right\rangle }\cdot\left(\Fourier^{-1}\gamma\right)\left(\varepsilon x\right)\right|^{2}\,\right)^{p/2}\!\!\d x\vphantom{\sum_{j\in M^{\left(\ell\right)}\cap J_{\left(r\right)}}}\nonumber \\
 & =C_{4}\varepsilon^{\dimension p}\cdot\left[\,\smash{\sum_{j\in M^{\left(\ell\right)}\cap J_{\left(r\right)}}}\vphantom{\sum}\left|c_{j}\right|^{2}\,\right]^{p/2}\vphantom{\sum_{j\in M^{\left(\ell\right)}\cap J_{\left(r\right)}}}\cdot\int_{\R^{\dimension}}\left|\left(\Fourier^{-1}\gamma\right)\left(\varepsilon x\right)\right|^{p}\,\d x\nonumber \\
\left({\scriptstyle \text{since }C_{1}=\left\Vert \Fourier^{-1}\gamma\right\Vert _{L^{p}}}\right) & =C_{1}^{p}C_{4}\cdot\varepsilon^{\dimension\left(p-1\right)}\cdot\left\Vert \left(c_{j}\right)_{j\in M^{\left(\ell\right)}\cap J_{\left(r\right)}}\right\Vert _{\ell^{2}}^{p}\nonumber \\
 & \leq\left[C_{1}C_{4}^{1/p}\cdot\varepsilon^{\dimension\left(1-\frac{1}{p}\right)}\cdot\left\Vert \left(c_{j}\right)_{j\in J_{0}\cap J_{\ell}}\right\Vert _{\ell^{2}}\right]^{p}.\label{eq:KhinchinNecessaryCoarseInFineExpectationEstimate}
\end{align}
Here, the last step used that we have $\xi_{j}\in Q_{\ell_{j}}\cap P_{j}=Q_{\ell}\cap P_{j}$
and hence $j\in M\cap J_{\ell}\subset J_{0}\cap J_{\ell}$ for $j\in M^{\left(\ell\right)}$.

Estimate (\ref{eq:KhinchinNecessaryCoarseInFineExpectationEstimate})
yields a (deterministic) realization $\theta^{\left(\ell\right)}=\theta^{\left(\ell,r\right)}=\left(\smash{\theta_{j}^{\left(\ell\right)}}\right)_{j\in M^{\left(\ell\right)}\cap J_{\left(r\right)}}\in\left\{ \pm1\right\} ^{M^{\left(\ell\right)}\cap J_{\left(r\right)}}$
with
\begin{equation}
\left\Vert \Fourier^{-1}\left[\,\smash{\sum_{j\in M^{\left(\ell\right)}\cap J_{\left(r\right)}}}\vphantom{\sum}c_{j}\cdot\theta_{j}^{\left(\ell\right)}\cdot\gamma_{j}\,\right]\right\Vert _{L^{p}}\vphantom{\sum_{j\in M^{\left(\ell\right)}\cap J_{\left(r\right)}}}\leq C_{1}C_{4}^{1/p}\cdot\varepsilon^{\dimension\left(1-\frac{1}{p}\right)}\cdot\left\Vert \left(c_{j}\right)_{j\in J_{0}\cap J_{\ell}}\right\Vert _{\ell^{2}}.\label{eq:KhinchinNecessaryCoarseInFineNiceRealization}
\end{equation}
Since we can choose such a realization for every $\ell\in I$ and
since the sets $\left(M^{\left(\ell\right)}\right)_{\ell\in I}$ are
pairwise disjoint with $M=\biguplus_{\ell\in I}M^{\left(\ell\right)}$,
we obtain a well-defined ``global'' realization $\theta=\theta^{\left(r\right)}\in\left\{ \pm1\right\} ^{M\cap J_{\left(r\right)}}$
defined by $\theta_{j}:=\theta_{j}^{\left(\ell_{j}\right)}$ for every
$j\in M\cap J_{\left(r\right)}$. Note that we have $\theta_{j}=\theta_{j}^{\left(\ell\right)}=\theta_{j}^{\left(\ell,r\right)}$
for every $j\in M^{\left(\ell\right)}\cap J_{\left(r\right)}$.

Now, define
\begin{equation}
g^{\left(r\right)}=g_{\theta}^{\left(r\right)}:=\sum_{j\in M\cap J_{\left(r\right)}}c_{j}\cdot\theta_{j}\cdot\gamma_{j}.\label{eq:KhinchinNecessaryConditionCoarseInFineTestFunctionDefinition}
\end{equation}
Observe $\supp\gamma_{j}\subset P_{j}\subset K$ for each $j\in M\subset J_{0}$
and hence $\supp g^{\left(r\right)}\subset K$, since the sum is finite.
Thus, $g^{\left(r\right)}\in\CalD_{K}^{\CalQ,p,Y}$, as soon as we
have shown $g^{\left(r\right)}\in\FourierDecompSp{\CalQ}pY$.

To see $g^{\left(r\right)}\in\FourierDecompSp{\CalQ}pY$, define $\varrho=\left(\varrho_{\ell}\right)_{\ell\in I}$
by $\varrho_{\ell}:=\left\Vert \left(c_{j}\right)_{j\in J_{0}\cap J_{\ell}}\right\Vert _{\ell^{2}}$.
Note that we have $\varrho\in Y$ with $\left\Vert \varrho\right\Vert _{Y}=\left\Vert c\right\Vert _{V}$,
since $c\in V=Y\left(\left[\ell^{s}\left(J_{0}\cap J_{i}\right)\right]_{i\in I}\right)$
and because of $s=2$, since $p\leq2$. Since $Y$ is $\CalQ$-regular,
this also implies $\Gamma_{\CalQ}\,\varrho\in Y$, which will soon
become helpful.

Indeed, let $i\in I$. Using estimate~(\ref{eq:KhinchinNecessaryConditionGenericUpperEstimate})
(with $z_{j}=0$ for all $j\in M$ and $d_{j}=c_{j}\theta_{j}\cdot\Indicator_{M\cap J_{\left(r\right)}}\left(j\right)$),
we derive
\begin{align*}
\left\Vert \Fourier^{-1}\left(\varphi_{i}\cdot\smash{g^{\left(r\right)}}\,\right)\right\Vert _{L^{p}} & \leq C_{2}C_{3}C_{\CalQ,\Phi,p}\cdot\sum_{\ell\in i^{\ast}}\left\Vert \Fourier^{-1}\left[\,\smash{\sum_{j\in M^{\left(\ell\right)}\cap J_{\left(r\right)}}}\vphantom{\sum}c_{j}\cdot\theta_{j}\cdot\gamma_{j}\,\right]\right\Vert _{L^{p}}\vphantom{\sum}\\
\left({\scriptstyle \text{eq. }\eqref{eq:KhinchinNecessaryCoarseInFineNiceRealization}\text{ and }\theta_{j}=\theta_{j}^{\left(\ell\right)}\text{ for }j\in M^{\left(\ell\right)}\cap J_{\left(r\right)}}\right) & \leq C_{1}C_{2}C_{3}C_{4}^{1/p}\cdot C_{\CalQ,\Phi,p}\cdot\varepsilon^{\dimension\left(1-\frac{1}{p}\right)}\cdot\vphantom{\sum_{\ell\in i^{\ast}}^{T}}\sum_{\ell\in i^{\ast}}\left\Vert \left(c_{j}\right)_{j\in J_{0}\cap J_{\ell}}\right\Vert _{\ell^{2}}\\
 & =C_{5}\cdot\varepsilon^{\dimension\left(1-\frac{1}{p}\right)}\cdot\left(\Gamma_{\CalQ}\,\varrho\right)_{i},
\end{align*}
with $C_{5}:=C_{1}C_{2}C_{3}C_{4}^{1/p}\cdot C_{\CalQ,\Phi,p}$. Since
the right-hand side is an element of $Y$, the solidity of $Y$ implies
$g^{\left(r\right)}\in\FourierDecompSp{\CalQ}pY$, with
\begin{align}
\left\Vert \smash{g^{\left(r\right)}}\right\Vert _{\BAPUFourierDecompSp{\CalQ}pY{\Phi}}=\left\Vert \left(\left\Vert \Fourier^{-1}\left(\varphi_{i}\cdot\smash{g^{\left(r\right)}}\right)\right\Vert _{L^{p}}\right)_{i\in I}\right\Vert _{Y} & \leq C_{5}\cdot\varepsilon^{\dimension\left(1-\frac{1}{p}\right)}\cdot\left\Vert \Gamma_{\CalQ}\,\varrho\right\Vert _{Y}\nonumber \\
 & \leq C_{5}\cdot\vertiii{\Gamma_{\CalQ}}_{Y\to Y}\cdot\varepsilon^{\dimension\left(1-\frac{1}{p}\right)}\cdot\left\Vert \varrho\right\Vert _{Y}\nonumber \\
 & =C_{6}\cdot\varepsilon^{\dimension\left(1-\frac{1}{p}\right)}\cdot\left\Vert c\right\Vert _{V}<\infty,\label{eq:KhinchinNecessaryConditionUpperEstimateComplete}
\end{align}
where $C_{6}:=C_{5}\cdot\vertiii{\Gamma_{\CalQ}}_{Y\to Y}$.

As noted above, we have now shown $g^{\left(r\right)}\in\CalD_{K}^{\CalQ,p,Y}$.
Since $\iota$ is well-defined and bounded, we conclude $g^{\left(r\right)}\in\FourierDecompSp{\CalP}pZ$,
with
\[
\left\Vert \smash{g^{\left(r\right)}}\right\Vert _{\BAPUFourierDecompSp{\CalP}pZ{\Psi}}\leq\vertiii{\iota}\cdot\left\Vert \smash{g^{\left(r\right)}}\right\Vert _{\BAPUFourierDecompSp{\CalQ}pY{\Phi}}\leq C_{6}\cdot\vertiii{\iota}\cdot\varepsilon^{\dimension\left(1-\frac{1}{p}\right)}\cdot\left\Vert c\right\Vert _{V}.
\]
It remains to obtain a lower bound for $\left\Vert g^{\left(r\right)}\right\Vert _{\BAPUFourierDecompSp{\CalP}pZ{\Psi}}$.
But this is a direct consequence of Corollary~\ref{cor:EasyNormEquivalenceFineLowerBound}
(with $\CalP$ instead of $\CalQ$ and with $I_{0}=M\cap J_{\left(r\right)}$):
By choice of $J_{\left(r\right)}$, we have for $j,\ell\in M\cap J_{\left(r\right)}$
with $j\neq\ell$ that $P_{j}^{\ast}\cap P_{\ell}^{\ast}=\emptyset$.
Furthermore, we have $\left(c_{j}\right)_{j\in M\cap J_{\left(r\right)}}\in\ell_{0}\left(M\cap J_{\left(r\right)}\right)$,
simply because $M$ is finite. Finally, $\gamma_{j}\in\TestFunctionSpace{\CalO'}$
with $\supp\gamma_{j}\subset P_{j}$ for all $j\in M\supset M\cap J_{\left(r\right)}$.
Thus, Corollary~\ref{cor:EasyNormEquivalenceFineLowerBound} shows
\[
C_{1}\cdot\varepsilon^{\dimension\left(1-\frac{1}{p}\right)}\cdot\left(c_{j}\right)_{j\in M\cap J_{\left(r\right)}}\overset{\text{eq. }\eqref{eq:KhinchinNecessaryConditionCoarseInFineNormOfPieces}}{=}\,\,\left(c_{j}\cdot\left\Vert \Fourier^{-1}\gamma_{j}\right\Vert _{L^{p}}\right)_{j\in M\cap J_{\left(r\right)}}\in Z|_{M\cap J_{\left(r\right)}}
\]
and yields a constant $C_{7}=C_{7}\left(\CalP,p,\dimension,C_{\CalP,\Psi,p},\vertiii{\Gamma_{\CalP}}_{Z\to Z}\right)>0$
satisfying
\begin{align*}
\left\Vert \smash{g^{\left(r\right)}}\right\Vert _{\BAPUFourierDecompSp{\CalP}pZ{\Psi}} & \geq C_{7}^{-1}\cdot\left\Vert \left(c_{j}\cdot\left\Vert \Fourier^{-1}\gamma_{j}\right\Vert _{L^{p}}\right)_{j\in M\cap J_{\left(r\right)}}\right\Vert _{Z|_{M\cap J^{\left(r\right)}}}\\
\left({\scriptstyle \text{eq. }\eqref{eq:KhinchinNecessaryConditionCoarseInFineNormOfPieces}\text{ and }c_{j}=0\text{ for }j\in J_{0}\setminus M}\right) & =C_{1}\cdot\varepsilon^{\dimension\left(1-\frac{1}{p}\right)}\cdot C_{7}^{-1}\cdot\left\Vert \left(c_{j}\cdot\Indicator_{J_{0}\cap J_{\left(r\right)}}\left(j\right)\right)_{j\in J_{0}}\right\Vert _{Z|_{J_{0}}}.
\end{align*}
Finally, using $J_{0}=\biguplus_{r=1}^{r_{0}}\left(J_{0}\cap J_{\left(r\right)}\right)$
and the triangle inequality for $Z$, we arrive at $c\in Z|_{J_{0}}$
with
\begin{align*}
\left\Vert c\right\Vert _{Z|_{J_{0}}}\leq C_{8}\cdot\sum_{r=1}^{r_{0}}\left\Vert c\cdot\Indicator_{J_{0}\cap J_{\left(r\right)}}\right\Vert _{Z|_{J_{0}}} & \leq C_{8}\cdot\sum_{r=1}^{r_{0}}\frac{C_{7}}{C_{1}\varepsilon^{\dimension\left(1-\frac{1}{p}\right)}}\left\Vert \smash{g^{\left(r\right)}}\right\Vert _{\BAPUFourierDecompSp{\CalP}pZ{\Psi}}\\
 & \leq C_{8}\cdot\sum_{r=1}^{r_{0}}\frac{C_{7}}{C_{1}\varepsilon^{\dimension\left(1-\frac{1}{p}\right)}}C_{6}\vertiii{\iota}\cdot\varepsilon^{\dimension\left(1-\frac{1}{p}\right)}\cdot\left\Vert c\right\Vert _{V}\\
 & =\frac{C_{6}C_{7}C_{8}r_{0}}{C_{1}}\cdot\vertiii{\iota}\cdot\left\Vert c\right\Vert _{V}\:,
\end{align*}
for a suitable constant $C_{8}=C_{8}\left(C_{Z},r_{0}\right)=C_{8}\left(C_{Z},N_{\CalP}\right)$.
This is precisely the desired embedding, since $c\in V_{0}$ was arbitrary.

\medskip{}

\textbf{Part~2}: Here, we have $p\geq2$ and $s=\max\left\{ 2,p\right\} =p$.
The proof in this case is a mixture of that of Theorem~\ref{thm:BurnerNecessaryConditionCoarseInFine}
and of that of the first part. We only provide it here for the sake
of completeness. 

Corollary~\ref{cor:AsymptoticModulationBehaviour} (with $M\cap J_{\left(r\right)}$
in place of $M$ and with $f_{j}=c_{j}\cdot\gamma_{j}$ for $j\in M\cap J_{\left(r\right)}$)
yields a family $z=\left(z_{j}\right)_{j\in M\cap J_{\left(r\right)}}\in\left(\R^{\dimension}\right)^{M\cap J_{\left(r\right)}}$
of modulations (possibly depending on $c$) satisfying
\begin{align}
\vphantom{\sum_{j\in S}}\left\Vert \Fourier^{-1}\left[\,\smash{\sum_{j\in S}}\vphantom{\sum}M_{z_{j}}\left(c_{j}\cdot\gamma_{j}\right)\,\right]\right\Vert _{L^{p}} & \leq2\cdot\left\Vert \left(\left\Vert \Fourier^{-1}\left(c_{j}\cdot\gamma_{j}\right)\right\Vert _{L^{p}}\right)_{j\in S}\right\Vert _{\ell^{p}}\nonumber \\
\left({\scriptstyle \text{eq. }\eqref{eq:KhinchinNecessaryConditionCoarseInFineNormOfPieces}}\right) & =2C_{1}\cdot\varepsilon^{\dimension\left(1-\frac{1}{p}\right)}\cdot\left\Vert \left(c_{j}\right)_{j\in S}\right\Vert _{\ell^{p}}\qquad\forall\,S\subset M\cap J_{\left(r\right)}\:.\label{eq:KhinchinNecessaryCoarseInFineNiceRealization2}
\end{align}
With this choice of $z$, define
\[
g^{\left(r\right)}:=\sum_{j\in M\cap J_{\left(r\right)}}M_{z_{j}}\left(c_{j}\cdot\gamma_{j}\right).
\]
As above, we have $\supp g^{\left(r\right)}\subset K$ and hence $g^{\left(r\right)}\in\CalD_{K}^{\CalQ,p,Y}$,
once we have shown $g^{\left(r\right)}\in\FourierDecompSp{\CalQ}pY$.

To see $g^{\left(r\right)}\in\FourierDecompSp{\CalQ}pY$, define $\varrho=\left(\varrho_{\ell}\right)_{\ell\in I}$
by $\varrho_{\ell}:=\left\Vert \left(c_{j}\right)_{j\in J_{0}\cap J_{\ell}}\right\Vert _{\ell^{p}}$.
Note that we have $\varrho\in Y$ with $\left\Vert \varrho\right\Vert _{Y}=\left\Vert c\right\Vert _{V}$,
since $c\in V=Y\left(\left[\ell^{s}\left(J_{0}\cap J_{i}\right)\right]_{i\in I}\right)$
and because of $s=p$, since $p\geq2$. Since $Y$ is $\CalQ$-regular,
this also implies $\Gamma_{\CalQ}\,\varrho\in Y$, which will soon
become helpful.

Indeed, let $i\in I$ be arbitrary. Using estimate~(\ref{eq:KhinchinNecessaryConditionGenericUpperEstimate})
(with $d_{j}=c_{j}\cdot\Indicator_{M\cap J_{\left(r\right)}}\left(j\right)$),
we derive
\begin{align}
\left\Vert \Fourier^{-1}\left(\varphi_{i}\cdot\smash{g^{\left(r\right)}}\,\right)\right\Vert _{L^{p}} & \leq C_{2}C_{3}C_{\CalQ,\Phi,p}\cdot\sum_{\ell\in i^{\ast}}\left\Vert \Fourier^{-1}\left[\,\smash{\sum_{j\in M^{\left(\ell\right)}\cap J_{\left(r\right)}}}\vphantom{\sum}c_{j}\cdot\gamma_{j}\,\right]\right\Vert _{L^{p}}\vphantom{\sum_{j\in M^{\left(\ell\right)}\cap J_{\left(r\right)}}}\nonumber \\
\left({\scriptstyle \text{eq. }\eqref{eq:KhinchinNecessaryCoarseInFineNiceRealization2}}\right) & \leq2C_{1}C_{2}C_{3}\cdot C_{\CalQ,\Phi,p}\cdot\varepsilon^{\dimension\left(1-\frac{1}{p}\right)}\cdot\vphantom{\sum_{\ell\in i^{\ast}}^{T}}\sum_{\ell\in i^{\ast}}\left\Vert \left(c_{j}\right)_{j\in M^{\left(\ell\right)}\cap J_{\left(r\right)}}\right\Vert _{\ell^{p}}\nonumber \\
 & \overset{\left(\ast\right)}{\leq}2C_{1}C_{2}C_{3}\cdot C_{\CalQ,\Phi,p}\cdot\varepsilon^{\dimension\left(1-\frac{1}{p}\right)}\cdot\sum_{\ell\in i^{\ast}}\left\Vert \left(c_{j}\right)_{j\in J_{0}\cap J_{\ell}}\right\Vert _{\ell^{p}}\nonumber \\
 & =C_{9}\cdot\varepsilon^{\dimension\left(1-\frac{1}{p}\right)}\cdot\left(\Gamma_{\CalQ}\,\varrho\right)_{i},\label{eq:KhinchinNecessaryCoarseInFinePart2UpperCoefficientEstimate}
\end{align}
with $C_{9}:=2C_{1}C_{2}C_{3}\cdot C_{\CalQ,\Phi,p}$. Here, the step
marked with $\left(\ast\right)$ used that we have $\xi_{j}\in P_{j}\cap Q_{\ell_{j}}=P_{j}\cap Q_{\ell}$
and hence $j\in J_{0}\cap J_{\ell}$ for $j\in M^{\left(\ell\right)}\subset M\subset J_{0}$.

Since the right-hand side of estimate~(\ref{eq:KhinchinNecessaryCoarseInFinePart2UpperCoefficientEstimate})
lies in $Y$, the solidity of $Y$ implies $g^{\left(r\right)}\in\FourierDecompSp{\CalQ}pY$,
with
\begin{align*}
\left\Vert \smash{g^{\left(r\right)}}\right\Vert _{\BAPUFourierDecompSp{\CalQ}pY{\Phi}}=\left\Vert \left(\left\Vert \Fourier^{-1}\left(\varphi_{i}\cdot\smash{g^{\left(r\right)}}\right)\right\Vert _{L^{p}}\right)_{i\in I}\right\Vert _{Y} & \leq C_{9}\cdot\varepsilon^{\dimension\left(1-\frac{1}{p}\right)}\cdot\left\Vert \Gamma_{\CalQ}\,\varrho\right\Vert _{Y}\\
 & \leq C_{9}\cdot\vertiii{\Gamma_{\CalQ}}_{Y\to Y}\cdot\varepsilon^{\dimension\left(1-\frac{1}{p}\right)}\cdot\left\Vert \varrho\right\Vert _{Y}\\
 & =C_{10}\cdot\varepsilon^{\dimension\left(1-\frac{1}{p}\right)}\cdot\left\Vert c\right\Vert _{V}<\infty,
\end{align*}
where $C_{10}:=C_{9}\cdot\vertiii{\Gamma_{\CalQ}}_{Y\to Y}$. This
is the exact analog of equation~(\ref{eq:KhinchinNecessaryConditionUpperEstimateComplete})
in the first part of the proof. The remainder of the proof is now
essentially identical to that of the first part and hence omitted.
Simply note that Corollary~\ref{cor:EasyNormEquivalenceFineLowerBound}
allows for an arbitrary modulation of the individual summands.

\medskip{}

The final claim of the theorem—in case $Z$ satisfies the Fatou property—is
an easy consequence of Lemmas \ref{lem:FatouPropertyIsInherited}
and \ref{lem:FinitelySupportedSequencesSufficeUnderFatouProperty}.
For more details, see the end of the proof of Theorem~\ref{thm:BurnerNecessaryConditionFineInCoarse}.
\end{proof}
Now, in our final necessary criterion in this subsection, we will
assume $\CalQ$ (or a subfamily of $\CalQ$) to be almost subordinate
to $\CalP$. Note though that no additional assumptions—like tightness—are
necessary.
\begin{thm}
\label{thm:KhinchinNecessaryFineInCoarse}Let $\emptyset\neq\CalO,\CalO'\subset\R^{\dimension}$
be open, let $p\in\left(0,\infty\right]$ and let $\CalQ=\left(Q_{i}\right)_{i\in I}$
and $\CalP=\left(P_{j}\right)_{j\in J}$ be two $L^{p}$-decomposition
coverings of $\CalO$ and $\CalO'$, respectively. Finally, let $Y\subset\Compl^{I}$
be $\CalQ$-regular and let $Z\subset\Compl^{J}$ be $\CalP$-regular,
with triangle constants $C_{Y}\geq1$ and $C_{Z}\geq1$, respectively.

Choose a subset $I_{0}\subset I$ such that we have $Q_{i}^{\circ}\neq\emptyset$
for all $i\in I_{0}$ and such that the restricted family $\CalQ_{I_{0}}:=\left(Q_{i}\right)_{i\in I_{0}}$
is almost subordinate to $\CalP$. Define 
\[
K:=\bigcup_{i\in I_{0}}Q_{i}\subset\CalO\cap\CalO'
\]
and assume—with $\CalD_{K}^{\CalQ,p,Y}$ as after equation~(\ref{eq:GeneralEmbeddingRequirement})—that
the identity map
\begin{equation}
\iota:\left(\CalD_{K}^{\CalQ,p,Y},\left\Vert \mybullet\right\Vert _{\FourierDecompSp{\CalQ}pY}\right)\rightarrow\FourierDecompSp{\CalP}pZ,f\mapsto f\label{eq:KhinchinEmbeddingFineInCoarse}
\end{equation}
is well-defined and bounded.

Let
\begin{equation}
s:=\begin{cases}
1, & \text{if }p=\infty,\\
\min\left\{ 2,p\right\} , & \text{if }p<\infty.
\end{cases}\label{eq:KhinchinFineInCoarseExponentDefinition}
\end{equation}
Then the embedding
\[
\eta:\ell_{0}\left(I_{0}\right)\cap Y|_{I_{0}}\hookrightarrow Z\bigl(\left[\ell^{s}\left(I_{0}\cap I_{j}\right)\right]_{j\in J}\bigr)
\]
is well-defined and bounded, with $\vertiii{\eta}\leq C\cdot\vertiii{\iota}$
for some constant
\[
C=C\left(\dimension,p,k\left(\CalQ_{I_{0}},\CalP\right),C_{Z},\CalQ,\CalP,C_{\CalQ,\Phi,p},C_{\CalP,\Psi,p},\vertiii{\Gamma_{\CalQ}}_{Y\to Y},\vertiii{\Gamma_{\CalP}}_{Z\to Z}\right)>0.
\]
Here, the $L^{p}$-BAPUs $\Phi=\left(\varphi_{i}\right)_{i\in I}$
and $\Psi=\left(\psi_{j}\right)_{j\in J}$ are those which are used
to compute the (quasi)-norms on the respective decomposition spaces
when computing the operator norm $\vertiii{\iota}$.

Finally, if $Z$ satisfies the Fatou property, the same statement
as above also holds for the embedding $\eta:Y|_{I_{0}}\hookrightarrow Z\bigl(\left[\ell^{s}\left(I_{0}\cap I_{j}\right)\right]_{j\in J}\bigr)$,
that is, without restricting to $\ell_{0}\left(I_{0}\right)$.
\end{thm}

\begin{rem*}
Note that we always have $s\leq2$ and $s\leq p$. Thus, we have norm-decreasing
embeddings $\ell^{s}\left(I_{0}\cap I_{j}\right)\hookrightarrow\ell^{2}\left(I_{0}\cap I_{j}\right)$
and $\ell^{s}\left(I_{0}\cap I_{j}\right)\hookrightarrow\ell^{p}\left(I_{0}\cap I_{j}\right)$.
Consequently,  it is easy to see that the theorem implies boundedness
of the embeddings
\[
\eta_{1}:\ell_{0}\left(I_{0}\right)\cap Y|_{I_{0}}\hookrightarrow Z\bigl([\ell^{2}\left(I_{0}\cap I_{j}\right)]_{j\in J}\bigr)\quad\text{and}\quad\eta_{2}:\ell_{0}\left(I_{0}\right)\cap Y|_{I_{0}}\hookrightarrow Z\bigl(\left[\,\ell^{p}\left(I_{0}\cap I_{j}\right)\right]_{j\in J}\bigr),
\]
with $\vertiii{\eta_{\ell}}\leq\vertiii{\eta}$ for $\ell\in\underline{2}$.
If $Z$ satisfies the Fatou property, then this holds even without
restricting to $\ell_{0}\left(I_{0}\right)$.
\end{rem*}
We remark that parts of the following proof are very similar to parts
of the proofs of Theorems \ref{thm:BurnerNecessaryConditionFineInCoarse}
and \ref{thm:KhinchinNecessaryCoarseInFine}. For the sake of completeness—and
to make the proof self-contained—we still provide almost all details.
\begin{proof}[Proof of Theorem~\ref{thm:KhinchinNecessaryFineInCoarse}]
As usual, we begin by setting up a few auxiliary objects. For brevity,
set $V:=Z\bigl(\left[\ell^{s}\left(I_{0}\cap I_{j}\right)\right]_{j\in J}\bigr)$
and $k:=k\left(\CalQ_{I_{0}},\CalP\right)$. For arbitrary $i\in I_{0}$,
there is some $j_{i}\in J$ with $Q_{i}\subset P_{j_{i}}^{k\ast}\subset\CalO'$.
Furthermore, $Q_{i}^{\circ}\neq\emptyset$ by assumption so that we
get $Q_{i}\cap P_{j}\neq\emptyset$ for some $j\in j_{i}^{k\ast}\subset J$.
In particular, $i\in I_{j}\cap I_{0}$ and thus $I_{0}=\bigcup_{j\in J}\left(I_{j}\cap I_{0}\right)$.
All in all, this shows that $V\subset\Compl^{I_{0}}$ is a solid sequence
space over $I_{0}$.

Fix a nonzero, nonnegative $\gamma\in\TestFunctionSpace{B_{1}\left(0\right)}$
for the rest of the proof. Set $r_{0}:=N_{\CalP}^{2\left(2k+3\right)+1}=N_{\CalP}^{4k+7}$
and note that Lemma~\ref{lem:DisjointizationPrinciple} yields a
partition $J=\biguplus_{r=1}^{r_{0}}J^{\left(r\right)}$ such that
$P_{j}^{\left(2k+3\right)\ast}\cap P_{\ell}^{\left(2k+3\right)\ast}=\emptyset$
holds for all $j,\ell\in J^{\left(r\right)}$ with $j\neq\ell$ and
all $r\in\underline{r_{0}}$.

As a crucial consequence of this, we have 
\begin{equation}
\forall\,j,\ell\in J^{\left(r\right)}\text{ with }j\neq\ell:\qquad\left(I_{j}\cap I_{0}\right)\cap\left(I_{\ell}\cap I_{0}\right)=\emptyset,\label{eq:KhinchinNecessaryFineInCoarseIntersectionIndexSetsDisjoint}
\end{equation}
because for $i\in I_{j}\cap I_{\ell}\cap I_{0}$, we would have $Q_{i}\cap P_{j}\neq\emptyset\neq Q_{i}\cap P_{\ell}$,
so that Lemma~\ref{lem:SubordinatenessImpliesWeakSubordination}
would yield
\[
\emptyset\neq Q_{i}\subset P_{j}^{\left(2k+2\right)\ast}\cap P_{\ell}^{\left(2k+2\right)\ast}\subset P_{j}^{\left(2k+3\right)\ast}\cap P_{\ell}^{\left(2k+3\right)\ast}
\]
in contradiction to $j,\ell\in J^{\left(r\right)}$ with $j\neq\ell$.

These considerations show that 
\[
I^{\left(r\right)}:=\biguplus_{j\in J^{\left(r\right)}}\left(I_{j}\cap I_{0}\right)\subset I_{0}
\]
is well-defined. As a simple consequence of this definition, we observe
\begin{equation}
\forall\,j\in J^{\left(r\right)}:\qquad I^{\left(r\right)}\cap I_{j}=I_{0}\cap I_{j}.\label{eq:KhinchinClusterDisjointizationIdentity}
\end{equation}

\medskip{}

Now, we properly start the proof. Let $c=\left(c_{i}\right)_{i\in I_{0}}\in\ell_{0}\left(I_{0}\right)\cap Y|_{I_{0}}$
and set $M:=\supp c\subset I_{0}$. Since each $Q_{i}$ with $i\in M\subset I_{0}$
has nonempty interior and because $M$ is finite, there is some $\varepsilon=\varepsilon\left(c\right)>0$
and for each $i\in M$ some $\xi_{i}\in\R^{\dimension}$ (possibly
depending on $c$) with $B_{\varepsilon}\left(\xi_{i}\right)\subset Q_{i}$.
Note that although $\varepsilon$ and the family $\left(\xi_{i}\right)_{i\in M}$
may depend heavily on the specific sequence $c$—more precisely on
its support—none of the constants $C_{1},C_{2},\dots$ in this proof
will depend on $c$; in particular, all occurrences of $\varepsilon$
will cancel in the end.

For $i\in M$, define $\gamma_{i}:=L_{\xi_{i}}\left[\gamma\left(\varepsilon^{-1}\mybullet\right)\right]$
as in the proof of Theorem~\ref{thm:KhinchinNecessaryCoarseInFine}
and note that equations (\ref{eq:KhinchinNecessaryConditionCoarseInFineFourierTransformOfPieces})
and (\ref{eq:KhinchinNecessaryConditionCoarseInFineNormOfPieces})
(with $C_{1}:=\left\Vert \Fourier^{-1}\gamma\right\Vert _{L^{p}}$,
i.e.\@ $C_{1}=C_{1}\left(p,\dimension\right)$) still apply, with
$j$ replaced by $i$. We extend the sequence $\left(c_{i}\right)_{i\in I_{0}}$
to all of $I$ by setting $c_{i}:=0$ for $i\in I\setminus I_{0}$.
Let us fix some $r\in\underline{r_{0}}$ (almost) until the end of
the proof.

In the following, we will ``test'' the embedding $\iota$ with a
function of the form
\[
g_{z,\theta}^{\left(r\right)}:=\sum_{i\in M\cap I^{\left(r\right)}}M_{z_{i}}\left(\theta_{i}\cdot\left|c_{i}\right|\cdot\gamma_{i}\right)
\]
for suitable modulations $z=\left(z_{i}\right)_{i\in M\cap I^{\left(r\right)}}\in\left(\R^{\dimension}\right)^{M\cap I^{\left(r\right)}}$
and signs $\theta=\left(\theta_{i}\right)_{i\in M\cap I^{\left(r\right)}}\in\left\{ \pm1\right\} ^{M\cap I^{\left(r\right)}}$.
To this end, first note that equation~(\ref{eq:KhinchinNecessaryConditionCoarseInFineNormOfPieces})
and the solidity of $Y$ yield 
\[
\left(\left|c_{i}\right|\cdot\left\Vert \Fourier^{-1}\gamma_{i}\right\Vert _{L^{p}}\right)_{i\in I}=\left(C_{1}\cdot\varepsilon^{\dimension\left(1-\frac{1}{p}\right)}\cdot\left|c_{i}\right|\right)_{i\in I}\in Y,
\]
so that Lemma~\ref{lem:EasyNormEquivalenceFineCovering} (with $k=0$
and $c_{i}=\gamma_{i}=0$ for $i\in I\setminus M$) implies $g_{z,\theta}^{\left(r\right)}\in\FourierDecompSp{\CalQ}pY$
with
\begin{align*}
\left\Vert g_{z,\theta}^{\left(r\right)}\right\Vert _{\BAPUFourierDecompSp{\CalQ}pY{\Phi}} & \leq C_{2}\cdot\left\Vert \left(\left|c_{i}\right|\cdot\left\Vert \Fourier^{-1}\gamma_{i}\right\Vert _{L^{p}}\right)_{i\in M\cap I^{\left(r\right)}}\right\Vert _{Y|_{M\cap I^{\left(r\right)}}}\\
 & =C_{1}C_{2}\cdot\varepsilon^{\dimension\left(1-\frac{1}{p}\right)}\cdot\left\Vert \left(\left|c_{i}\right|\right)_{i\in M\cap I^{\left(r\right)}}\right\Vert _{Y|_{M\cap I^{\left(r\right)}}}\\
\left({\scriptstyle \text{since }Y\text{ is solid}}\right) & \leq C_{1}C_{2}\cdot\varepsilon^{\dimension\left(1-\frac{1}{p}\right)}\cdot\left\Vert c\right\Vert _{Y|_{I_{0}}}<\infty
\end{align*}
for some constant $C_{2}=C_{2}\left(\dimension,p,\CalQ,C_{\CalQ,\Phi,p},\vertiii{\Gamma_{\CalQ}}_{Y\to Y}\right)$.
Here, we used $\supp\gamma_{i}\subset B_{\varepsilon}\left(\xi_{i}\right)\subset Q_{i}$
for all $i\in M\cap I^{\left(r\right)}$.

Since we have $\gamma_{i}\in\TestFunctionSpace{\CalO}$ with $\supp\gamma_{i}\subset Q_{i}\subset K$
for all $i\in I_{0}$ and since the sum defining $g_{z,\theta}^{\left(r\right)}$
is finite, we get $g_{z,\theta}^{\left(r\right)}\in\TestFunctionSpace{\CalO}$
with $\supp g_{z,\theta}^{\left(r\right)}\subset K$ and thus $g_{z,\theta}^{\left(r\right)}\in\CalD_{K}^{\CalQ,p,Y}$.
Since the identity map $\iota$ from the statement of the theorem
is well-defined and bounded, this implies $g_{z,\theta}^{\left(r\right)}\in\FourierDecompSp{\CalP}pZ$,
with
\begin{equation}
\left\Vert g_{z,\theta}^{\left(r\right)}\right\Vert _{\BAPUFourierDecompSp{\CalP}pZ{\Psi}}\leq\vertiii{\iota}\cdot\left\Vert g_{z,\theta}^{\left(r\right)}\right\Vert _{\BAPUFourierDecompSp{\CalQ}pY{\Phi}}\leq C_{1}C_{2}\cdot\varepsilon^{\dimension\left(1-\frac{1}{p}\right)}\cdot\vertiii{\iota}\cdot\left\Vert c\right\Vert _{Y|_{I_{0}}}<\infty.\label{eq:KhinchinFineInCoarseUpperEstimateComplete}
\end{equation}
In the remainder of the proof, we will obtain lower bounds on $\left\Vert g_{z,\theta}^{\left(r\right)}\right\Vert _{\BAPUFourierDecompSp{\CalP}pZ{\Psi}}$,
for suitable values of $z$ and $\theta$.

\medskip{}

To this end, we will now prove the following observation: For arbitrary
$j\in J^{\left(r\right)}$ and $i\in M\cap I^{\left(r\right)}$, we
have
\begin{equation}
\psi_{j}^{\left(2k+3\right)\ast}\cdot\gamma_{i}=\begin{cases}
\gamma_{i}, & \text{if }i\in I_{j},\\
0, & \text{if }i\notin I_{j}.
\end{cases}\label{eq:KhinchinHyperIdenticalLocalizationIdentity}
\end{equation}
Indeed, for $i\in M\subset I_{0}$ and $j\in J$ with $Q_{i}\cap P_{j}\neq\emptyset$,
Lemma~\ref{lem:SubordinatenessImpliesWeakSubordination} yields $\supp\gamma_{i}\subset Q_{i}\subset P_{j}^{\left(2k+2\right)\ast}$,
and Lemma~\ref{lem:PartitionCoveringNecessary} implies $\psi_{j}^{\left(2k+3\right)\ast}\equiv1$
on $P_{j}^{\left(2k+2\right)\ast}$. We have thus established the
following:
\begin{equation}
\forall\,i\in M\quad\forall\,j\in J\text{ with }Q_{i}\cap P_{j}\neq\emptyset:\qquad\psi_{j}^{\left(2k+3\right)\ast}\cdot\gamma_{i}=\gamma_{i}.\label{eq:KhinchinNecessaryFineInCoarseBAPUMultiplication}
\end{equation}
Now, let $j\in J^{\left(r\right)}$ and $i\in M\cap I^{\left(r\right)}\subset I_{0}$
be arbitrary. In case of $i\in I_{j}$, equation~(\ref{eq:KhinchinNecessaryFineInCoarseBAPUMultiplication})
yields $\psi_{j}^{\left(2k+3\right)\ast}\cdot\gamma_{i}=\gamma_{i}$.
Otherwise, if $i\notin I_{j}$, then the definition of $I^{\left(r\right)}$
yields some $\ell\in J^{\left(r\right)}\setminus\left\{ j\right\} $
with $i\in I_{\ell}\cap I_{0}$. Note that this entails $\psi_{\ell}^{\left(2k+3\right)\ast}\cdot\gamma_{i}=\gamma_{i}$
by equation~(\ref{eq:KhinchinNecessaryFineInCoarseBAPUMultiplication}).
But because of $j,\ell\in J^{\left(r\right)}$ with $j\neq\ell$,
we get $P_{j}^{\left(2k+3\right)\ast}\cap P_{\ell}^{\left(2k+3\ast\right)}=\emptyset$,
whence $\psi_{j}^{\left(2k+3\right)\ast}\cdot\psi_{\ell}^{\left(2k+3\right)\ast}\equiv0$
and thus
\[
\psi_{j}^{\left(2k+3\right)\ast}\cdot\gamma_{i}=\psi_{j}^{\left(2k+3\right)\ast}\cdot\psi_{\ell}^{\left(2k+3\right)\ast}\cdot\gamma_{i}\equiv0.
\]
All in all, we have established identity~(\ref{eq:KhinchinHyperIdenticalLocalizationIdentity})
from above in both cases.

Next, Remark~\ref{rem:ClusteredBAPUYIeldsBoundedControlSystem} shows
that the family $\Lambda:=\left(\psi_{j}^{\left(2k+3\right)\ast}\right)_{j\in J}$
is an $L^{p}$-bounded family for $\CalP$, with $C_{\CalP,\Lambda,p}\leq C_{3}=C_{3}\left(\CalP,C_{\CalP,\Psi,p},\dimension,p,k\right)$
and $\ell_{\Lambda,\CalP}=2k+3$. Consequently, Theorem~\ref{thm:BoundedControlSystemEquivalentQuasiNorm}
yields a further constant $C_{4}=C_{4}\left(\CalP,p,\dimension,k,\vertiii{\Gamma_{\CalP}}_{Z\to Z}\right)>0$
satisfying
\begin{align*}
\left\Vert \left(\left\Vert \Fourier^{-1}\left(\psi_{j}^{\left(2k+3\right)\ast}\cdot g_{z,\theta}^{\left(r\right)}\right)\right\Vert _{L^{p}}\right)_{j\in J^{\left(r\right)}}\right\Vert _{Z|_{J^{\left(r\right)}}} & \leq\left\Vert \left(\left\Vert \Fourier^{-1}\left(\psi_{j}^{\left(2k+3\right)\ast}\cdot g_{z,\theta}^{\left(r\right)}\right)\right\Vert _{L^{p}}\right)_{j\in J}\right\Vert _{Z}\\
 & \leq C_{4}\cdot C_{\CalP,\Lambda,p}\cdot\left\Vert g_{z,\theta}^{\left(r\right)}\right\Vert _{\BAPUFourierDecompSp{\CalP}pZ{\Psi}}\\
\left({\scriptstyle \text{eq. }\eqref{eq:KhinchinFineInCoarseUpperEstimateComplete}}\right) & \leq C_{1}C_{2}C_{3}C_{4}\cdot\varepsilon^{\dimension\left(1-\frac{1}{p}\right)}\cdot\vertiii{\iota}\cdot\left\Vert c\right\Vert _{Y|_{I_{0}}}.
\end{align*}

But equation~(\ref{eq:KhinchinHyperIdenticalLocalizationIdentity})
and the special form of $g_{z,\theta}^{\left(r\right)}$ show  for
$j\in J^{\left(r\right)}$ that
\[
\psi_{j}^{\left(2k+3\right)\ast}\cdot g_{z,\theta}^{\left(r\right)}=\sum_{i\in M\cap I^{\left(r\right)}}\left[\theta_{i}\cdot\left|c_{i}\right|\cdot M_{z_{i}}\left(\psi_{j}^{\left(2k+3\right)\ast}\cdot\gamma_{i}\right)\right]=\sum_{i\in M\cap I^{\left(r\right)}\cap I_{j}}M_{z_{i}}\left(\theta_{i}\cdot\left|c_{i}\right|\cdot\gamma_{i}\right),
\]
so that we get—for $C_{5}:=C_{1}C_{2}C_{3}C_{4}$—the estimate
\begin{equation}
\vphantom{\sum_{i\in M\cap I^{\left(r\right)}\cap I_{j}}}\left\Vert \!\left(\left\Vert \vphantom{\sum_{i}}\Fourier^{-1}\left[\,\smash{\sum_{i\in M\cap I^{\left(r\right)}\cap I_{j}}}\vphantom{\sum}\!\!\!\!\!M_{z_{i}}\!\left(\theta_{i}\cdot\left|c_{i}\right|\cdot\gamma_{i}\right)\,\right]\right\Vert _{L^{p}}\right)_{\!\!\!j\in J^{\left(r\right)}}\right\Vert _{Z|_{J^{\left(r\right)}}}\!\!\leq C_{5}\cdot\varepsilon^{\dimension\left(1-\frac{1}{p}\right)}\vertiii{\iota}\cdot\left\Vert c\right\Vert _{Y|_{I_{0}}}\label{eq:KhinchinNecessaryFineInCoarseFundamentalEstimate}
\end{equation}
for all $r\in\underline{r_{0}}$ and $z\in\left(\R^{\dimension}\right)^{M\cap I^{\left(r\right)}}$,
as well as $\theta\in\left\{ \pm1\right\} ^{M\cap I^{\left(r\right)}}$.

In the remainder of the proof, we will show for each $r\in\underline{r_{0}}$
that one can choose the parameters $z^{\left(r\right)}=\left(\smash{z_{i}^{\left(r\right)}}\right)_{i\in M\cap I^{\left(r\right)}}\in\left(\R^{\dimension}\right)^{M\cap I^{\left(r\right)}}$
and $\theta^{\left(r\right)}=\left(\smash{\theta_{i}^{\left(r\right)}}\right)_{i\in M\cap I^{\left(r\right)}}\in\left\{ \pm1\right\} ^{M\cap I^{\left(r\right)}}$
in such a way that we get
\begin{equation}
\varepsilon^{\dimension\left(1-\frac{1}{p}\right)}\cdot\left\Vert \left(c_{i}\right)_{i\in I_{0}\cap I_{j}}\right\Vert _{\ell^{s}}\leq C_{6}\cdot\left\Vert \Fourier^{-1}\left[\,\smash{\sum_{i\in M\cap I_{j}\cap I^{\left(r\right)}}}\vphantom{\sum}M_{z_{i}^{\left(r\right)}}\left(\,\!\smash{\theta_{i}^{\left(r\right)}}\cdot\left|c_{i}\right|\cdot\gamma_{i}\right)\,\right]\right\Vert _{L^{p}}\vphantom{\sum_{i\in M\cap I_{j}\cap I^{\left(r\right)}}^{T}}\qquad\forall\,j\in J^{\left(r\right)}\label{eq:KhinchinNecessaryFineInCoarseRemainingEstimate}
\end{equation}
for a suitable constant $C_{6}=C_{6}\left(\dimension,p\right)$. Once
this is done, the quasi-triangle inequality for $Z$ and the identity
$J=\biguplus_{r=1}^{r_{0}}J^{\left(r\right)}$ yield a constant $C_{7}=C_{7}\left(C_{Z},r_{0}\right)=C_{7}\left(C_{Z},k,N_{\CalP}\right)$
satisfying
\begin{align*}
\varepsilon^{\dimension\left(1-\frac{1}{p}\right)}\cdot\left\Vert c\right\Vert _{V} & =\left\Vert \left(\varepsilon^{\dimension\left(1-\frac{1}{p}\right)}\cdot\left\Vert \left(c_{i}\right)_{i\in I_{0}\cap I_{j}}\right\Vert _{\ell^{s}}\right)_{j\in J}\right\Vert _{Z}\\
 & \leq C_{7}\cdot\sum_{r=1}^{r_{0}}\left\Vert \left(\varepsilon^{\dimension\left(1-\frac{1}{p}\right)}\cdot\left\Vert \left(c_{i}\right)_{i\in I_{0}\cap I_{j}}\right\Vert _{\ell^{s}}\right)_{j\in J^{\left(r\right)}}\right\Vert _{Z|_{J^{\left(r\right)}}}\\
\left({\scriptstyle \text{equation }\eqref{eq:KhinchinNecessaryFineInCoarseRemainingEstimate}}\right) & \leq C_{6}C_{7}\cdot\sum_{r=1}^{r_{0}}\left\Vert \left(\left\Vert \vphantom{\sum_{i}}\Fourier^{-1}\left(\smash{\sum_{i\in M\cap I_{j}\cap I^{\left(r\right)}}}\vphantom{\sum}M_{z_{i}^{\left(r\right)}}\left(\smash{\theta_{i}^{\left(r\right)}}\cdot\left|c_{i}\right|\cdot\gamma_{i}\right)\right)\right\Vert _{L^{p}}\right)_{j\in J^{\left(r\right)}}\right\Vert _{Z|_{J^{\left(r\right)}}}\vphantom{\sum_{i\in M\cap I_{j}\cap I^{\left(r\right)}}}\\
\left({\scriptstyle \text{equation }\eqref{eq:KhinchinNecessaryFineInCoarseFundamentalEstimate}}\right) & \leq r_{0}C_{5}C_{6}C_{7}\cdot\varepsilon^{\dimension\left(1-\frac{1}{p}\right)}\vertiii{\iota}\cdot\left\Vert c\right\Vert _{Y|_{I_{0}}},
\end{align*}
so that canceling the common factor $\varepsilon^{\dimension\left(1-\frac{1}{p}\right)}$
yields the claim. Note that the quantitative arguments from above
(and solidity of $Z$) in particular yield the qualitative conclusion
$c\in V$.

It remains to show that we can choose $z^{\left(r\right)},\theta^{\left(r\right)}$
in such a way that estimate~(\ref{eq:KhinchinNecessaryFineInCoarseRemainingEstimate})
is fulfilled. To this end we will distinguish the three cases which
are indicated by the definition of $s$; see equation~(\ref{eq:KhinchinFineInCoarseExponentDefinition}).

\medskip{}

\textbf{Case 1}: $p=\infty$. In this case, choose $z_{i}^{\left(r\right)}=0$
and $\theta_{i}^{\left(r\right)}=1$ for all $i\in M\cap I^{\left(r\right)}$.
Then the continuity of $\Fourier^{-1}\left(\sum_{i\in M\cap I_{j}}\left|c_{i}\right|\cdot\gamma_{i}\right)\in\Schwartz\left(\R^{\dimension}\right)$
implies because of $M\cap I_{j}=M\cap I_{0}\cap I_{j}=M\cap I_{j}\cap I^{\left(r\right)}$
(see equation~(\ref{eq:KhinchinClusterDisjointizationIdentity}))
that
\begin{align*}
\vphantom{\sum_{i\in M\cap I_{j}\cap I^{\left(r\right)}}}\left\Vert \Fourier^{-1}\left[\,\smash{\sum_{i\in M\cap I_{j}\cap I^{\left(r\right)}}}\vphantom{\sum}M_{z_{i}^{\left(r\right)}}\left(\,\!\smash{\theta_{i}^{\left(r\right)}}\cdot\left|c_{i}\right|\cdot\gamma_{i}\right)\,\right]\right\Vert _{L^{p}} & =\left\Vert \Fourier^{-1}\left(\,\smash{\sum_{i\in M\cap I_{j}}}\vphantom{\sum}\left|c_{i}\right|\cdot\gamma_{i}\,\right)\right\Vert _{L^{\infty}}\vphantom{\sum_{i\in M\cap I_{j}}}\\
 & \geq\left|\smash{\sum_{i\in M\cap I_{j}}}\vphantom{\sum}\left|c_{i}\right|\cdot\left(\Fourier^{-1}\gamma_{i}\right)\left(0\right)\right|=\left|\smash{\sum_{i\in M\cap I_{j}}}\vphantom{\sum}\left|c_{i}\right|\cdot\int_{\R^{\dimension}}\gamma_{i}\left(\xi\right)\,\d\xi\right|\vphantom{\sum_{i\in M\cap I_{j}}}\\
\left({\scriptstyle \text{since }\gamma_{i}\geq0}\right) & =\sum_{i\in M\cap I_{j}}\left|c_{i}\right|\cdot\left\Vert \gamma_{i}\right\Vert _{L^{1}}\\
\left({\scriptstyle \text{Riemann-Lebesgue}}\right) & \geq\sum_{i\in M\cap I_{j}}\left|c_{i}\right|\cdot\left\Vert \Fourier^{-1}\gamma_{i}\right\Vert _{L^{\infty}}\\
\left({\scriptstyle \text{eq. }\eqref{eq:KhinchinNecessaryConditionCoarseInFineNormOfPieces}\text{ and }c_{i}=0\text{ if }i\in I_{0}\setminus M}\right) & =C_{1}\varepsilon^{\dimension\left(1-\frac{1}{p}\right)}\cdot\left\Vert \left(c_{i}\right)_{i\in I_{0}\cap I_{j}}\right\Vert _{\ell^{1}}\qquad\forall\,j\in J^{\left(r\right)}\,.
\end{align*}
Since we have $s=1$ for $p=\infty$, this is nothing but estimate~(\ref{eq:KhinchinNecessaryFineInCoarseRemainingEstimate}),
with $C_{6}=C_{1}^{-1}$. \medskip{}

\textbf{Case 2}: $2<p<\infty$. For $j\in J^{\left(r\right)}$ consider
the random variable $\omega^{\left(j\right)}=\left(\smash{\omega_{i}^{\left(j\right)}}\right)_{i\in M\cap I_{j}}\in\left\{ \pm1\right\} ^{M\cap I_{j}}$
which we take to be uniformly distributed in $\left\{ \pm1\right\} ^{M\cap I_{j}}$.
Let $C_{8}=C_{8}\left(p\right)>0$ be the constant provided by Khintchine's
inequality (Theorem~\ref{thm:KhintchineInequality}). By applying
that inequality, as well as equation~(\ref{eq:KhinchinNecessaryConditionCoarseInFineFourierTransformOfPieces}),
i.e.\@ $\Fourier^{-1}\gamma_{i}=\varepsilon^{\dimension}\cdot M_{\xi_{i}}\left[\left(\Fourier^{-1}\gamma\right)\left(\varepsilon\mybullet\right)\right]$,
we conclude 
\begin{align*}
\vphantom{\sum_{i\in I_{j}\cap M}}\mathbb{E}\left\Vert \Fourier^{-1}\left[\,\smash{\sum_{i\in I_{j}\cap M}}\vphantom{\sum}\left|c_{i}\right|\omega_{i}^{\left(j\right)}\gamma_{i}\,\right]\right\Vert _{L^{p}}^{p} & =\mathbb{E}\int_{\R^{\dimension}}\left|\,\smash{\sum_{i\in I_{j}\cap M}}\vphantom{\sum}\left|c_{i}\right|\omega_{i}^{\left(j\right)}\cdot\left(\Fourier^{-1}\gamma_{i}\right)\left(x\right)\,\right|^{p}\,\d x\vphantom{\sum_{i\in I_{j}\cap M}}\\
\left({\scriptstyle \text{equation }\eqref{eq:KhinchinNecessaryConditionCoarseInFineFourierTransformOfPieces}}\right) & =\varepsilon^{\dimension p}\cdot\int_{\R^{\dimension}}\mathbb{E}\left|\,\smash{\sum_{i\in I_{j}\cap M}}\vphantom{\sum}\left|c_{i}\right|\omega_{i}^{\left(j\right)}\cdot e^{2\pi i\left\langle \xi_{i},x\right\rangle }\cdot\left(\Fourier^{-1}\gamma\right)\left(\varepsilon x\right)\,\right|^{p}\,\d x\vphantom{\sum_{i\in I_{j}\cap M}}\\
\left({\scriptstyle \text{Khintchine's inequality}}\right) & \geq C_{8}^{-1}\varepsilon^{\dimension p}\cdot\int_{\R^{\dimension}}\left(\,\smash{\sum_{i\in I_{j}\cap M}}\vphantom{\sum}\left|\,\left|c_{i}\right|\cdot e^{2\pi i\left\langle \xi_{i},x\right\rangle }\cdot\left(\Fourier^{-1}\gamma\right)\left(\varepsilon x\right)\right|^{2}\,\right)^{p/2}\,\d x\vphantom{\sum_{i\in I_{j}\cap M}}\\
 & =C_{8}^{-1}\varepsilon^{\dimension p}\cdot\left(\,\smash{\sum_{i\in I_{j}\cap M}}\vphantom{\sum}\left|c_{i}\right|^{2}\,\right)^{p/2}\cdot\int_{\R^{\dimension}}\left|\left(\Fourier^{-1}\gamma\right)\left(\varepsilon x\right)\right|^{p}\,\d x\vphantom{\sum_{i\in I_{j}\cap M}}\\
 & =\frac{\left\Vert \Fourier^{-1}\gamma\right\Vert _{L^{p}}^{p}}{C_{8}}\cdot\varepsilon^{\dimension\left(p-1\right)}\cdot\left\Vert \left(c_{i}\right)_{i\in I_{j}\cap M}\right\Vert _{\ell^{2}}^{p}\\
\left({\scriptstyle \text{since }c_{i}=0\text{ for }i\in I_{0}\setminus M}\right) & =\left[C_{8}^{-1/p}\cdot\left\Vert \Fourier^{-1}\gamma\right\Vert _{L^{p}}\cdot\varepsilon^{\dimension\left(1-\frac{1}{p}\right)}\cdot\left\Vert \left(c_{i}\right)_{i\in I_{0}\cap I_{j}}\right\Vert _{\ell^{2}}\right]^{p}.
\end{align*}

In particular, the estimate above yields a \emph{deterministic} realization
$\omega^{\left(j\right)}=\left(\smash{\omega_{i}^{\left(j\right)}}\right)_{i\in I_{j}\cap M}\in\left\{ \pm1\right\} ^{I_{j}\cap M}$
with
\begin{equation}
\vphantom{\sum_{i\in I_{j}\cap M}}\left\Vert \Fourier^{-1}\left[\,\smash{\sum_{i\in I_{j}\cap M}}\vphantom{\sum}\left|c_{i}\right|\omega_{i}^{\left(j\right)}\gamma_{i}\,\right]\right\Vert _{L^{p}}\geq\frac{\left\Vert \Fourier^{-1}\gamma\right\Vert _{L^{p}}}{C_{8}^{1/p}}\cdot\varepsilon^{\dimension\left(1-\frac{1}{p}\right)}\cdot\left\Vert \left(c_{i}\right)_{i\in I_{0}\cap I_{j}}\right\Vert _{\ell^{2}}.\label{eq:KhinchinNecessaryFineInCoarseNiceRealization}
\end{equation}
Recall from above that the sets $\left(I_{0}\cap I_{j}\right)_{j\in J^{\left(r\right)}}$
form a partition of the set $I^{\left(r\right)}\subset I_{0}$. Since
$M\subset I_{0}$, we thus see that $\left(I_{j}\cap M\right)_{j\in J^{\left(r\right)}}$
partitions $M\cap I^{\left(r\right)}$. Thus, we can combine the individual
``sign realizations'' $\omega^{\left(j\right)}$ for $j\in J^{\left(r\right)}$
into the ``global'' sequence $\theta^{\left(r\right)}=\left(\smash{\theta_{i}^{\left(r\right)}}\right)_{i\in M\cap I^{\left(r\right)}}$
given by $\theta_{i}^{\left(r\right)}=\omega_{i}^{\left(j\right)}$,
where—for given $i\in M\cap I^{\left(r\right)}$—the index $j=j\left(i\right)$
is the unique $j\in J^{\left(r\right)}$ with $i\in M\cap I_{j}$.

Note that we have $\theta_{i}^{\left(r\right)}=\omega_{i}^{\left(j\right)}$
for all $i\in M\cap I_{j}$ and arbitrary $j\in J^{\left(r\right)}$.
Furthermore, recall from equation~(\ref{eq:KhinchinClusterDisjointizationIdentity})
that $M\cap I_{j}\cap I^{\left(r\right)}=M\cap I_{0}\cap I_{j}=M\cap I_{j}$.
Thus, given this ``joint sign choice'' $\theta^{\left(r\right)}$,
our previous considerations easily imply that equation~(\ref{eq:KhinchinNecessaryFineInCoarseRemainingEstimate})
is indeed satisfied, with $z_{i}^{\left(r\right)}=0$ for all $i\in M\cap I^{\left(r\right)}$
and with $C_{6}=C_{8}^{1/p}\cdot\left\Vert \Fourier^{-1}\gamma\right\Vert _{L^{p}}^{-1}$,
since we have $s=2$ in the present case $2<p<\infty$. \medskip{}

\textbf{Case 3}: We have $0<p\leq2$. In this case, we choose $\theta_{i}^{\left(r\right)}=1$
for all $i\in M\cap I^{\left(r\right)}$. Next, Corollary~\ref{cor:AsymptoticModulationBehaviour},
applied to the family $\left(f_{i}\right)_{i\in M\cap I^{\left(r\right)}}$
with $f_{i}:=\left|c_{i}\right|\cdot\gamma_{i}\in\Schwartz\left(\R^{\dimension}\right)$,
yields a family $z^{\left(r\right)}=\left(\smash{z_{i}^{\left(r\right)}}\right)_{i\in M\cap I^{\left(r\right)}}$
of modulations satisfying\vspace{0.15cm}
\begin{align*}
\vphantom{\sum_{i\in S}}\left\Vert \Fourier^{-1}\left[\,\smash{\sum_{i\in S}}\vphantom{\sum}M_{z_{i}^{\left(r\right)}}\left(\smash{\theta_{i}^{\left(r\right)}}\cdot\left|c_{i}\right|\cdot\gamma_{i}\right)\,\right]\right\Vert _{L^{p}} & =\left\Vert \Fourier^{-1}\left[\,\smash{\sum_{i\in S}}\vphantom{\sum}M_{z_{i}^{\left(r\right)}}f_{i}\,\right]\right\Vert _{L^{p}}\vphantom{\sum_{i\in S}}\\
 & \geq\frac{1}{2}\cdot\left\Vert \left(\left\Vert \Fourier^{-1}f_{i}\right\Vert _{L^{p}}\right)_{i\in S}\right\Vert _{\ell^{p}}=\frac{1}{2}\cdot\left\Vert \left(\left|c_{i}\right|\cdot\left\Vert \Fourier^{-1}\gamma_{i}\right\Vert _{L^{p}}\right)_{i\in S}\right\Vert _{\ell^{p}}\\
\left({\scriptstyle \text{eq. }\eqref{eq:KhinchinNecessaryConditionCoarseInFineNormOfPieces}}\right) & =\frac{C_{1}}{2}\cdot\varepsilon^{\dimension\left(1-\frac{1}{p}\right)}\cdot\left\Vert \left(c_{i}\right)_{i\in S}\right\Vert _{\ell^{p}}\qquad\forall\,S\subset M\cap I^{\left(r\right)}\,.
\end{align*}
If we apply this with $S=M\cap I_{j}=M\cap I_{0}\cap I_{j}=M\cap I^{\left(r\right)}\cap I_{j}$
(for $j\in J^{\left(r\right)}$) and recall $s=p$ (because of $0<p\leq2$),
we get
\begin{align*}
\vphantom{\sum_{i\in M\cap I_{j}\cap I^{\left(r\right)}}}\left\Vert \Fourier^{-1}\left[\,\smash{\sum_{i\in M\cap I_{j}\cap I^{\left(r\right)}}}\vphantom{\sum}M_{z_{i}^{\left(r\right)}}\left(\smash{\theta_{i}^{\left(r\right)}}\cdot\left|c_{i}\right|\cdot\gamma_{i}\right)\,\right]\right\Vert _{L^{p}} & \geq\frac{C_{1}}{2}\cdot\varepsilon^{\dimension\left(1-\frac{1}{p}\right)}\cdot\left\Vert \left(c_{i}\right)_{i\in M\cap I_{j}}\right\Vert _{\ell^{p}}\\
\left({\scriptstyle s=p\text{ and }c_{i}=0\text{ for }i\in I_{0}\setminus M}\right) & =\frac{C_{1}}{2}\cdot\varepsilon^{\dimension\left(1-\frac{1}{p}\right)}\cdot\left\Vert \left(c_{i}\right)_{i\in I_{0}\cap I_{j}}\right\Vert _{\ell^{s}},
\end{align*}
which is nothing but equation~(\ref{eq:KhinchinNecessaryFineInCoarseRemainingEstimate}),
with $C_{6}=2C_{1}^{-1}$.

\medskip{}

This completes the proof, with the exception of the additional claim
in case $Z$ satisfies the Fatou property. But this claim is an easy
consequence of Lemmas \ref{lem:FatouPropertyIsInherited} and \ref{lem:FinitelySupportedSequencesSufficeUnderFatouProperty}.
For more details, see the end of the proof of Theorem~\ref{thm:BurnerNecessaryConditionFineInCoarse}.
\end{proof}

\subsection{Complete characterizations for relatively moderate coverings}

\label{subsec:RelativelyModerateCase}In this subsection, we develop
further necessary conditions for embeddings between decomposition
spaces. These will show that in general, one has to use the exponents
$\LowerExpo{p_{2}}$ and $\UpperExpo{p_{1}}$ (or $\SignedUpperExpo{p_{1}}$)
to calculate the ``inner norm'', as in the sufficient conditions
from Remarks~\ref{rem:SufficientFineInCoarseSimplification} and
\ref{rem:SufficientCoarseIntoFineSimplification}. Note that up to
now, all our necessary criteria were only able to show that these
conditions are necessary for the existence of the embedding if one
replaces $\LowerExpo{p_{2}}$ by $p_{2}$ and $\UpperExpo{p_{1}}$
by $p_{1}$, see Theorems~\ref{thm:BurnerNecessaryConditionFineInCoarse}
and \ref{thm:BurnerNecessaryConditionCoarseInFine}; only in case
of $p_{1}=p_{2}$, we were able to obtain a slight improvement, see
Theorems~\ref{thm:KhinchinNecessaryCoarseInFine} and \ref{thm:KhinchinNecessaryFineInCoarse}.

For the proof of our stronger necessary conditions, we place more
severe restrictions on the relation between the two coverings and
on the ``global'' components of the two decomposition spaces. Precisely,
we will only consider weighted $\ell^{q}$ spaces as our global components.
Furthermore, we will assume that one covering is almost subordinate
as well as \emph{relatively moderate} with respect to the other and
that the weight of the ``subordinate'' covering is also moderate
with respect to the ``coarse'' covering.

Note that in the preceding subsections, the functions which we used
to ``test'' the embedding were always adapted to the \emph{finer}
of the two coverings. In contrast, in this subsection, our ``test
functions'' will be adapted to the \emph{coarser} of the two coverings.
We remark that the basic idea of this construction is taken from Han
and Wang\cite{HanWangAlphaModulationEmbeddings}, who use a similar
construction for the special case of embeddings between $\alpha$-modulation
spaces.

As a preparation, we first establish a necessary condition which does
\emph{not} assume relative moderateness of the two coverings. Afterwards,
we specialize this to the relatively moderate case. Occasionally,
the following lemma—especially the ensuing remark—will also be useful
for coverings which are not relatively moderate to each other.
\begin{lem}
\label{lem:NecessaryConjugateCoarseInFineWithoutModerateness}Let
$\emptyset\neq\CalO,\CalO'\subset\R^{\dimension}$ be open, let $\CalQ=\left(Q_{i}\right)_{i\in I}=\left(T_{i}Q_{i}'+b_{i}\right)_{i\in I}$
be a semi-structured $L^{p_{1}}$-decomposition covering of $\CalO$
and let $\CalP=\left(P_{j}\right)_{j\in J}=\left(\smash{S_{j}P_{j}'+c_{j}}\right)_{j\in J}$
be a \emph{tight} semi-structured $L^{p_{2}}$-decomposition covering
of $\CalO'$ for certain $p_{1},p_{2}\in\left(0,\infty\right]$.
Let $Y\subset\Compl^{I}$ and $Z\subset\Compl^{J}$ be $\CalQ$-regular
and $\CalP$-regular, with triangle constants $C_{Y}\geq1$ and $C_{Z}\geq1$,
respectively.

Assume that $J_{0}\subset J$ is such that $\CalP_{J_{0}}:=\left(P_{j}\right)_{j\in J_{0}}$
is almost subordinate to $\CalQ$, and define 
\[
I_{0}:=\left\{ i\in I\with J_{0}\cap J_{i}\neq\emptyset\right\} \,.
\]
Set\footnote{The inclusion $K\subset\CalO$ is a direct consequence of Lemma~\ref{lem:PartitionCoveringNecessary}.}
\[
k:=k\left(\CalP_{J_{0}},\CalQ\right)\quad\text{and}\quad K:=\bigcup_{i\in I_{0}}\overline{Q_{i}^{\left(2k+3\right)\ast}}\subset\CalO,
\]
and assume—with $\CalD_{K}^{\CalQ,p_{1},Y}$ as after equation~(\ref{eq:GeneralEmbeddingRequirement})—that
there is a bounded linear map
\[
\iota:\left(\CalD_{K}^{\CalQ,p_{1},Y},\left\Vert \mybullet\right\Vert _{\FourierDecompSp{\CalQ}{p_{1}}Y}\right)\to\FourierDecompSp{\CalP}{p_{2}}Z
\]
which satisfies $\left\langle \iota f,\,\varphi\right\rangle _{\CalD'}=\left\langle f,\varphi\right\rangle _{\CalD'}$
for all $\varphi\in\TestFunctionSpace{\CalO\cap\CalO'}$ and all $f\in\CalD_{K}^{\CalQ,p_{1},Y}$.

\medskip{}

Then the map
\[
\eta:\ell_{0}\left(I_{0}\right)\cap Y|_{I_{0}}\to Z,\left(x_{i}\right)_{i\in I_{0}}\mapsto\sum_{i\in I_{0}}\left[x_{i}\cdot\left|\det T_{i}\right|^{p_{1}^{-1}-1}\cdot\left(\left|\det S_{j}\right|^{1-p_{2}^{-1}}\cdot\Indicator_{J_{0}\cap J_{i}}\left(j\right)\right)_{j\in J}\right]
\]
is well-defined and bounded, with $\vertiii{\eta}\leq C\cdot\vertiii{\iota}$,
for a constant $C>0$ of the form
\[
C=C\left(\dimension,C_{Z},p_{1},p_{2},k\left(\CalP_{J_{0}},\CalQ\right),\CalQ,\CalP,\varepsilon_{\CalP},C_{\CalQ,\Phi,p_{1}},\vertiii{\Gamma_{\CalQ}}_{Y\to Y},\vertiii{\Gamma_{\CalP}}_{Z\to Z}\right).
\]
Here, as usual, the $L^{p_{1}}$-BAPU $\Phi=\left(\varphi_{i}\right)_{i\in I}$
has to be used to calculate the (quasi)-norm on the decomposition
space $\FourierDecompSp{\CalQ}{p_{1}}Y$ for computing $\vertiii{\iota}$.

\medskip{}

Finally, in case of $Z=\ell_{v}^{q_{2}}\left(J\right)$ for a $\CalP$-moderate
weight $v=\left(v_{j}\right)_{j\in J}$ and some $q_{2}\in\left(0,\infty\right]$,
we get that the embedding 
\[
\gamma:Y|_{I_{0}}\hookrightarrow\ell_{u}^{q_{2}}\left(I_{0}\right)\qquad\text{ with }\qquad u_{i}:=\left|\det T_{i}\right|^{p_{1}^{-1}-1}\cdot\left\Vert \left(v_{j}\cdot\left|\det S_{j}\right|^{1-p_{2}^{-1}}\right)_{j\in J_{0}\cap J_{i}}\right\Vert _{\ell^{q_{2}}}
\]
is well-defined and bounded, with $\vertiii{\gamma}\leq C'\cdot\vertiii{\iota}$
for some constant 
\[
C'=C'\left(\dimension,q_{2},p_{1},p_{2},k\left(\CalP_{J_{0}},\CalQ\right),\CalQ,\CalP,\varepsilon_{\CalP},C_{v,\CalP},C_{\CalQ,\Phi,p_{1}},\vertiii{\Gamma_{\CalQ}}_{Y\to Y}\right).
\]
In particular, $u_{i}<\infty$ for all $i\in I_{0}$ with $\delta_{i}\in Y$.
\end{lem}

\begin{rem}
\label{rem:NecessaryConjugateCoarseInFineElementaryConsequence}As
we have seen previously, even applying the derived embeddings for
sequence spaces to sequences which are supported on a single point
can be useful. Indeed, in the present case, let $j_{0}\in J_{0}$
be arbitrary and let $i_{0}\in I$ with $Q_{i_{0}}\cap P_{j_{0}}\neq\emptyset$.
This yields $j_{0}\in J_{0}\cap J_{i_{0}}$ and hence $i_{0}\in I_{0}$.
Now, if $\delta_{i_{0}}\in Y$, we can apply the lemma to conclude
\[
\left|\det T_{i_{0}}\right|^{p_{1}^{-1}-1}\cdot\left|\det S_{j_{0}}\right|^{1-p_{2}^{-1}}\cdot\left\Vert \delta_{j_{0}}\right\Vert _{Z}\leq\left|\det T_{i_{0}}\right|^{p_{1}^{-1}-1}\cdot\left\Vert \left(\left|\det S_{j}\right|^{1-p_{2}^{-1}}\right)_{j\in J_{0}\cap J_{i_{0}}}\right\Vert _{Z|_{J_{0}\cap J_{i_{0}}}}\leq\vertiii{\eta}\cdot\left\Vert \delta_{i_{0}}\right\Vert _{Y}.
\]
In several cases, this estimate yields sharper conditions than estimate~(\ref{eq:ElementaryNecessaryConditionCoarseInFine})
from Remark~\ref{rem:BurnerNecessaryConditionCoarseInFine}; see
the proof of Theorem~\ref{thm:InhomogeneousIntoHomogeneousBesov}
for an example where this occurs.
\end{rem}

\begin{proof}
Choose $\left(\gamma_{j}\right)_{j\in J}$ as in Lemma~\ref{lem:NormOfClusteredBAPUAndTestFunctionBuildingBlocks}
(applied to $\CalP$ instead of $\CalQ$). Let $r_{0}:=N_{\CalQ}^{2\left(2k+3\right)+1}=N_{\CalQ}^{4k+7}$,
so that Lemma~\ref{lem:DisjointizationPrinciple} yields a partition
$I=\biguplus_{r=1}^{r_{0}}I^{\left(r\right)}$ for which $Q_{i}^{\left(2k+3\right)\ast}\cap Q_{\ell}^{\left(2k+3\right)\ast}=\emptyset$
holds for all $r\in\underline{r_{0}}$ and all $i,\ell\in I^{\left(r\right)}$
with $i\neq\ell$.

We first note that the assumption on $\iota$ implies 
\begin{equation}
\gamma_{j}\cdot\iota f=\gamma_{j}\cdot f\text{ for all }f\in\CalD_{K}^{\CalQ,p_{1},Y}\text{ and }j\in J_{0},\label{eq:NecessaryConjugateCoseInFineWithoutModeratenessConsistency}
\end{equation}
where the equality has to be understood in the sense of tempered distributions.
To see this, let $g\in\Schwartz\left(\R^{\dimension}\right)$ be arbitrary.
By definition, we have
\[
\left\langle \gamma_{j}\cdot\iota f,\,g\right\rangle _{\Schwartz'}=\left\langle \iota f,\,\gamma_{j}\,g\right\rangle _{\CalD'}\overset{\left(\dagger\right)}{=}\left\langle f,\,\gamma_{j}\,g\right\rangle _{\CalD'}=\left\langle \gamma_{j}\cdot f,\,g\right\rangle _{\Schwartz'}\:,
\]
where we used at $\left(\dagger\right)$ that $\supp\gamma_{j}\subset P_{j}\subset Q_{i_{j}}^{k\ast}\subset\CalO$
for some $i_{j}\in I$ and also $\supp\gamma_{j}\subset P_{j}\subset\CalO'$.
Taken together, we see $\gamma_{j}\,g\in\TestFunctionSpace{\CalO\cap\CalO'}$,
so that $\left\langle \iota f,\,\gamma_{j}\,g\right\rangle _{\CalD'}=\left\langle f,\,\gamma_{j}\,g\right\rangle _{\CalD'}$
holds by assumption on $\iota$.

\medskip{}

For $r\in\underline{r_{0}}$, the family $\left(J_{0}\cap J_{i}\right)_{i\in I^{\left(r\right)}}$
is pairwise disjoint, because for $i,\ell\in I^{\left(r\right)}$
with $j\in J_{0}\cap J_{i}\cap J_{\ell}$, Lemma~\ref{lem:SubordinatenessImpliesWeakSubordination}
would imply
\[
\emptyset\neq P_{j}\subset Q_{i}^{\left(2k+2\right)\ast}\cap Q_{\ell}^{\left(2k+2\right)\ast}\subset Q_{i}^{\left(2k+3\right)\ast}\cap Q_{\ell}^{\left(2k+3\right)\ast}
\]
which in turn yields $i=\ell$ by choice of $\left(I^{\left(r\right)}\right)_{r\in\underline{r_{0}}}$.
Thus, we can define 
\[
J^{\left(r\right)}:=\biguplus_{i\in I^{\left(r\right)}}\left(J_{0}\cap J_{i}\right)\subset J_{0}.
\]
Note that this definition yields
\[
J^{\left(r\right)}\cap J_{i}=J_{0}\cap J_{i}\qquad\forall\,i\in I^{\left(r\right)}.
\]

As a next step, we note for $i\in I^{\left(r\right)}$ and $j\in J^{\left(r\right)}$
that
\begin{equation}
\gamma_{j}\cdot\varphi_{i}^{\left(2k+3\right)\ast}=\begin{cases}
\gamma_{j}, & \text{if }j\in J_{0}\cap J_{i},\\
0, & \text{otherwise}.
\end{cases}\label{eq:NecessaryModerateMutualLocalization}
\end{equation}
To see this, first assume $j\in J_{0}\cap J_{i}$. Lemma~\ref{lem:SubordinatenessImpliesWeakSubordination}
yields $P_{j}\subset Q_{i}^{\left(2k+2\right)\ast}$ and Lemma~\ref{lem:PartitionCoveringNecessary}
implies $\varphi_{i}^{\left(2k+3\right)\ast}\equiv1$ on $Q_{i}^{\left(2k+2\right)\ast}\supset P_{j}$.
Since $\gamma_{j}$ vanishes outside of $P_{j}$, this establishes
the first case. For the case $j\notin J_{0}\cap J_{i}$, observe that
we have $j\in J^{\left(r\right)}$ and thus $j\in J_{0}\cap J_{\ell}$
for some $\ell\in I^{\left(r\right)}$. As we just saw, this implies
$\gamma_{j}=\gamma_{j}\cdot\varphi_{\ell}^{\left(2k+3\right)\ast}$.
But $j\notin J_{0}\cap J_{i}$ yields $i\neq\ell$ and thus $Q_{i}^{\left(2k+3\right)\ast}\cap Q_{\ell}^{\left(2k+3\right)\ast}=\emptyset$
because of $i,\ell\in I^{\left(r\right)}$. This leads to $\gamma_{j}\cdot\varphi_{i}^{\left(2k+3\right)\ast}=\gamma_{j}\cdot\varphi_{\ell}^{\left(2k+3\right)\ast}\cdot\varphi_{j}^{\left(2k+3\right)\ast}\equiv0$,
so that equation~(\ref{eq:NecessaryModerateMutualLocalization})
is established.

\medskip{}

Now, let $d=\left(d_{i}\right)_{i\in I_{0}}\in\ell_{0}\left(I_{0}\right)\cap Y|_{I_{0}}$
be arbitrary and consider $d$ as an element of $\ell_{0}\left(I\right)$
by extending it trivially. Set 
\[
\zeta_{i}:=\left|\det T_{i}\right|^{p_{1}^{-1}-1}\cdot\left|d_{i}\right|\quad\text{for }i\in I,\qquad\text{and}\qquad g^{\left(r\right)}:=\sum_{i\in I_{0}\cap I^{\left(r\right)}}\zeta_{i}\cdot\varphi_{i}^{\left(2k+3\right)\ast}\quad\text{for }r\in\underline{r_{0}}\,.
\]
 Note that we have $g^{\left(r\right)}\in\TestFunctionSpace{\CalO}$
with $\supp g^{\left(r\right)}\subset K$, as a \emph{finite} sum
of functions which satisfy the same properties.

Now, recall from Lemma~\ref{lem:NormOfClusteredBAPUAndTestFunctionBuildingBlocks}
that $\left\Vert \Fourier^{-1}\varphi_{i}^{\left(2k+3\right)\ast}\right\Vert _{L^{p_{1}}}\leq C_{1}\cdot\left|\det T_{i}\right|^{1-p_{1}^{-1}}$
for all $i\in I$, with a constant $C_{1}=C_{1}\left(\dimension,p_{1},k,\CalQ,C_{\CalQ,\Phi,p_{1}}\right)$.
Thus, Lemma~\ref{lem:EasyNormEquivalenceFineCovering}, combined
with the solidity of $Y$, implies $g^{\left(r\right)}\in\FourierDecompSp{\CalQ}{p_{1}}Y$,
with
\begin{equation}
\left\Vert \smash{g^{\left(r\right)}}\right\Vert _{\FourierDecompSp{\CalQ}{p_{1}}Y}\leq\!C_{2}\cdot\left\Vert \!\left(\zeta_{i}\cdot\left\Vert \Fourier^{-1}\varphi_{i}^{\left(2k+3\right)\ast}\right\Vert _{L^{p_{1}}}\right)_{\!i\in I}\right\Vert _{Y}\!\leq\!C_{1}C_{2}\cdot\left\Vert d\right\Vert _{Y}\!=C_{1}C_{2}\cdot\left\Vert d\right\Vert _{Y|_{I_{0}}}\!<\infty\,,\label{eq:NecessaryConjugateCoarseInFineQNormUpperEstimate}
\end{equation}
for some constant $C_{2}=C_{2}\left(k,\dimension,p_{1},\CalQ,C_{\CalQ,\Phi,p_{1}},\vertiii{\Gamma_{\CalQ}}_{Y\to Y}\right)$.
All in all, we have shown $g^{\left(r\right)}\in\CalD_{K}^{\CalQ,p_{1},Y}$,
so that $\iota\,g^{\left(r\right)}\in\FourierDecompSp{\CalP}{p_{2}}Z$
is well-defined.

\medskip{}

Now, we want to establish a lower bound for $\left\Vert \iota\,g^{\left(r\right)}\right\Vert _{\FourierDecompSp{\CalP}{p_{2}}Z}$.
To this end, note that equations (\ref{eq:NecessaryModerateMutualLocalization})
and (\ref{eq:NecessaryConjugateCoseInFineWithoutModeratenessConsistency})
and the pairwise disjointness of the family $\left(J_{0}\cap J_{\ell}\right)_{\ell\in I^{\left(r\right)}}$
yield
\begin{equation}
\left\Vert \Fourier^{-1}\left(\gamma_{j}\cdot\smash{\iota g^{\left(r\right)}}\right)\right\Vert _{L^{p_{2}}}\!=\zeta_{i}\cdot\left\Vert \Fourier^{-1}\gamma_{j}\right\Vert _{L^{p_{2}}}\!=C_{3}\cdot\zeta_{i}\cdot\left|\det S_{j}\right|^{1-p_{2}^{-1}}\quad\forall\,i\in I_{0}\cap I^{\left(r\right)}\text{ and }j\in J_{0}\cap J_{i}\,,\label{eq:NecessaryConjugateCoarseInFinePLocalizedLowerEstimate}
\end{equation}
where the constant $C_{3}=C_{3}\left(\dimension,p_{2},\varepsilon_{\CalP}\right)$
is taken from Lemma~\ref{lem:NormOfClusteredBAPUAndTestFunctionBuildingBlocks}.

Similarly, Lemma~\ref{lem:NormOfClusteredBAPUAndTestFunctionBuildingBlocks}
shows that $\left(\gamma_{j}\right)_{j\in J}$ satisfies for $s:=\min\left\{ 1,p_{2}\right\} $
the estimate 
\[
\left|\det S_{j}\right|^{s^{-1}-1}\cdot\left\Vert \Fourier^{-1}\gamma_{j}\right\Vert _{L^{s}}\leq C=C\left(\dimension,p_{2},\varepsilon_{\CalP}\right)\qquad\forall\,j\in J\,.
\]
Since we also have $\supp\gamma_{j}\subset P_{j}$, we see that $\Gamma=\left(\gamma_{j}\right)_{j\in J}$
is an $L^{p_{2}}$-bounded system for $\CalP$, with $C_{\CalP,\Gamma,p_{2}}\leq C\left(\dimension,p_{2},\varepsilon_{\CalP}\right)$
and $\ell_{\Gamma,\CalP}=0$. By virtue of Theorem~\ref{thm:BoundedControlSystemEquivalentQuasiNorm},
this yields a positive constant $C_{4}=C_{4}\left(\dimension,p_{2},\CalP,\varepsilon_{\CalP},\vertiii{\Gamma_{\CalP}}_{Z\to Z}\right)$
with
\begin{align}
\left\Vert \iota\,\smash{g^{\left(r\right)}}\right\Vert _{\FourierDecompSp{\CalP}{p_{2}}Z} & \geq C_{4}^{-1}\cdot\left\Vert \left(\left\Vert \Fourier^{-1}\left(\gamma_{j}\cdot\iota\,\smash{g^{\left(r\right)}}\right)\right\Vert _{L^{p_{2}}}\right)_{j\in J}\right\Vert _{Z}\nonumber \\
 & \geq C_{4}^{-1}\cdot\left\Vert \left(\left\Vert \Fourier^{-1}\left(\gamma_{j}\cdot\iota\,\smash{g^{\left(r\right)}}\right)\right\Vert _{L^{p_{2}}}\cdot\Indicator_{J^{\left(r\right)}}\left(j\right)\right)_{j\in J}\right\Vert _{Z}.\label{eq:NecessaryConjugateCoarseInFineLowerEstimate}
\end{align}
In particular, Theorem~\ref{thm:BoundedControlSystemEquivalentQuasiNorm}
shows that the sequence on the right-hand side of equation~(\ref{eq:NecessaryConjugateCoarseInFineLowerEstimate})
is indeed an element of $Z$.

But we have $J^{\left(r\right)}=\biguplus_{i\in I^{\left(r\right)}}\left(J_{0}\cap J_{i}\right)$
and $J_{0}\cap J_{i}=\emptyset$ in case of $i\notin I_{0}$, i.e.
$J^{\left(r\right)}=\biguplus_{i\in I_{0}\cap I^{\left(r\right)}}\left(J_{0}\cap J_{i}\right)$.
Thus, for any $j\in J$, we have
\begin{align*}
\left\Vert \Fourier^{-1}\left(\gamma_{j}\cdot\iota\,\smash{g^{\left(r\right)}}\right)\right\Vert _{L^{p_{2}}}\cdot\Indicator_{J^{\left(r\right)}}\left(j\right) & =\sum_{i\in I^{\left(r\right)}\cap I_{0}}\left\Vert \Fourier^{-1}\left(\gamma_{j}\cdot\iota\,\smash{g^{\left(r\right)}}\right)\right\Vert _{L^{p_{2}}}\cdot\Indicator_{J_{0}\cap J_{i}}\left(j\right)\\
\left({\scriptstyle \text{eq. }\eqref{eq:NecessaryConjugateCoarseInFinePLocalizedLowerEstimate}}\right) & =C_{3}\cdot\sum_{i\in I^{\left(r\right)}\cap I_{0}}\zeta_{i}\cdot\left(\left|\det S_{j}\right|^{1-p_{2}^{-1}}\cdot\Indicator_{J_{0}\cap J_{i}}\left(j\right)\right)\\
\left({\scriptstyle \text{since }d_{i}=0\text{ for }i\in I\setminus I_{0}}\right) & =C_{3}\cdot\sum_{i\in I^{\left(r\right)}}\left[\left|d_{i}\right|\cdot\left|\det T_{i}\right|^{p_{1}^{-1}-1}\cdot\left(\left|\det S_{j}\right|^{1-p_{2}^{-1}}\cdot\Indicator_{J_{0}\cap J_{i}}\left(j\right)\right)\right].
\end{align*}
Now, summing over $r\in\underline{r_{0}}$, using $I=\biguplus_{r=1}^{r_{0}}I^{\left(r\right)}$,
and recalling that $Z$ is solid, we see that the map $\eta$ defined
in the current lemma is well-defined. Furthermore, the (quasi)-triangle
inequality for $Z$ ensures existence of a constant $C_{5}=C_{5}\left(r_{0},C_{Z}\right)=C_{5}\left(N_{\CalQ},k,C_{Z}\right)$
satisfying
\begin{align*}
\left\Vert \eta\left(d\right)\right\Vert _{Z}\leq\left\Vert \eta\left(\left|d\right|\right)\right\Vert _{Z} & \overset{\left(\ast\right)}{=}\left\Vert \sum_{i\in I}\left|d_{i}\right|\cdot\left|\det T_{i}\right|^{p_{1}^{-1}-1}\cdot\left(\left|\det S_{j}\right|^{1-p_{2}^{-1}}\cdot\Indicator_{J_{0}\cap J_{i}}\left(j\right)\right)_{j\in J}\right\Vert _{Z}\\
 & \leq C_{5}\cdot\sum_{r=1}^{r_{0}}\left\Vert \sum_{i\in I^{\left(r\right)}}\left|d_{i}\right|\cdot\left|\det T_{i}\right|^{p_{1}^{-1}-1}\cdot\left(\left|\det S_{j}\right|^{1-p_{2}^{-1}}\cdot\Indicator_{J_{0}\cap J_{i}}\left(j\right)\right)_{j\in J}\right\Vert _{Z}\\
 & =\frac{C_{5}}{C_{3}}\cdot\sum_{r=1}^{r_{0}}\left\Vert \left(\Indicator_{J^{\left(r\right)}}\left(j\right)\cdot\left\Vert \Fourier^{-1}\left(\gamma_{j}\cdot\iota\,\smash{g^{\left(r\right)}}\right)\right\Vert _{L^{p_{2}}}\right)_{j\in J}\right\Vert _{Z}\\
\left({\scriptstyle \text{eq. }\eqref{eq:NecessaryConjugateCoarseInFineLowerEstimate}}\right) & \leq\frac{C_{4}C_{5}}{C_{3}}\cdot\sum_{r=1}^{r_{0}}\left\Vert \iota\,\smash{g^{\left(r\right)}}\right\Vert _{\FourierDecompSp{\CalP}{p_{2}}Z}\leq\frac{C_{4}C_{5}}{C_{3}}\vertiii{\iota}\cdot\sum_{r=1}^{r_{0}}\left\Vert \smash{g^{\left(r\right)}}\right\Vert _{\FourierDecompSp{\CalQ}{p_{1}}Y}\\
\left({\scriptstyle \text{eq. }\eqref{eq:NecessaryConjugateCoarseInFineQNormUpperEstimate}}\right) & \leq r_{0}\frac{C_{1}C_{2}C_{4}C_{5}}{C_{3}}\vertiii{\iota}\cdot\left\Vert d\right\Vert _{Y|_{I_{0}}}.
\end{align*}
Here, the first step (marked with $\left(\ast\right)$) used that
$d_{i}=0$ for $i\in I\setminus I_{0}$. All in all, we have shown
that $\eta$ is well-defined and bounded, so that it remains to consider
the second part of the lemma.

\medskip{}

To this end, note that $J^{\left(r\right)}=\biguplus_{i\in I^{\left(r\right)}}\left(J_{0}\cap J_{i}\right)$
implies
\[
\left\Vert \left(c_{j}\right)_{j\in J^{\left(r\right)}}\right\Vert _{\ell^{s}}=\left\Vert \left(\left\Vert \left(c_{j}\right)_{j\in J_{0}\cap J_{i}}\right\Vert _{\ell^{s}}\right)_{i\in I^{\left(r\right)}}\right\Vert _{\ell^{s}}
\]
for arbitrary sequences $\left(c_{j}\right)_{j\in J^{\left(r\right)}}$
and any $s\in\left(0,\infty\right]$. If we apply this in equation~(\ref{eq:NecessaryConjugateCoarseInFineLowerEstimate}),
we get
\begin{align*}
\infty>\left\Vert \iota\,\smash{g^{\left(r\right)}}\right\Vert _{\FourierDecompSp{\CalP}{p_{2}}Z} & \geq C_{4}^{-1}\cdot\left\Vert \left(v_{j}\cdot\left\Vert \Fourier^{-1}\left(\gamma_{j}\cdot\iota\,\smash{g^{\left(r\right)}}\right)\right\Vert _{L^{p_{2}}}\right)_{j\in J^{\left(r\right)}}\right\Vert _{\ell^{q_{2}}}\\
 & =C_{4}^{-1}\cdot\left\Vert \left(\left\Vert \left(v_{j}\cdot\left\Vert \Fourier^{-1}\left(\gamma_{j}\cdot\iota\,\smash{g^{\left(r\right)}}\right)\right\Vert _{L^{p_{2}}}\right)_{j\in J_{0}\cap J_{i}}\right\Vert _{\ell^{q_{2}}}\right)_{i\in I^{\left(r\right)}}\right\Vert _{\ell^{q_{2}}}\\
 & \geq C_{4}^{-1}\cdot\left\Vert \left(\left\Vert \left(v_{j}\cdot\left\Vert \Fourier^{-1}\left(\gamma_{j}\cdot\iota\,\smash{g^{\left(r\right)}}\right)\right\Vert _{L^{p_{2}}}\right)_{j\in J_{0}\cap J_{i}}\right\Vert _{\ell^{q_{2}}}\right)_{i\in I^{\left(r\right)}\cap I_{0}}\right\Vert _{\ell^{q_{2}}}\\
\left({\scriptstyle \text{eq. }\eqref{eq:NecessaryConjugateCoarseInFinePLocalizedLowerEstimate}}\right) & =\frac{C_{3}}{C_{4}}\cdot\left\Vert \left(\left\Vert \left(v_{j}\cdot\zeta_{i}\cdot\left|\det S_{j}\right|^{1-p_{2}^{-1}}\right)_{j\in J_{0}\cap J_{i}}\right\Vert _{\ell^{q_{2}}}\right)_{i\in I^{\left(r\right)}\cap I_{0}}\right\Vert _{\ell^{q_{2}}}\\
 & =\frac{C_{3}}{C_{4}}\cdot\left\Vert \left(\left|d_{i}\right|\cdot\left|\det T_{i}\right|^{p_{1}^{-1}-1}\cdot\left\Vert \left(v_{j}\cdot\left|\det S_{j}\right|^{1-p_{2}^{-1}}\right)_{j\in J_{0}\cap J_{i}}\right\Vert _{\ell^{q_{2}}}\right)_{i\in I^{\left(r\right)}\cap I_{0}}\right\Vert _{\ell^{q_{2}}}\\
 & =\frac{C_{3}}{C_{4}}\cdot\left\Vert \left(d_{i}\right)_{i\in I^{\left(r\right)}\cap I_{0}}\right\Vert _{\ell_{u}^{q_{2}}\left(I_{0}\cap I^{\left(r\right)}\right)}.
\end{align*}
Now, as above, we sum over $r\in\underline{r_{0}}$ and use the (quasi)-triangle
inequality for $\ell^{q_{2}}$ to obtain
\begin{align*}
\left\Vert d\right\Vert _{\ell_{u}^{q_{2}}\left(I_{0}\right)} & \leq C_{6}\cdot\sum_{r=1}^{r_{0}}\left\Vert \left(d_{i}\cdot\Indicator_{I^{\left(r\right)}\cap I_{0}}\left(i\right)\right)_{i\in I_{0}}\right\Vert _{\ell_{u}^{q_{2}}\left(I_{0}\right)}\\
 & \leq\frac{C_{4}C_{6}}{C_{3}}\cdot\sum_{r=1}^{r_{0}}\left\Vert \iota\,\smash{g^{\left(r\right)}}\right\Vert _{\FourierDecompSp{\CalP}{p_{2}}Z}\\
\left({\scriptstyle \text{eq. }\eqref{eq:NecessaryConjugateCoarseInFineQNormUpperEstimate}}\right) & \leq r_{0}\frac{C_{1}C_{2}C_{4}C_{6}}{C_{3}}\vertiii{\iota}\cdot\left\Vert d\right\Vert _{Y|_{I_{0}}}<\infty
\end{align*}
for some constant $C_{6}=C_{6}\left(r_{0},q_{2}\right)=C_{6}\left(N_{\CalQ},k,q_{2}\right)$.
For $C':=r_{0}\frac{C_{1}C_{2}C_{4}C_{6}}{C_{3}}$, this shows $\vertiii{\gamma_{0}}\leq C'\cdot\vertiii{\iota}$
for $\gamma_{0}:\ell_{0}\left(I_{0}\right)\cap Y|_{I_{0}}\hookrightarrow\ell_{u}^{q_{2}}\left(I_{0}\right)$.

Finally, since $\ell_{u}^{q_{2}}\left(I_{0}\right)$ satisfies the
Fatou property, Lemma~\ref{lem:FinitelySupportedSequencesSufficeUnderFatouProperty}
yields boundedness of $\gamma$, as well as $\vertiii{\gamma}\leq C'\cdot\vertiii{\iota}$,
as desired. Regarding the dependencies of the constants, note that
Lemma~\ref{lem:ModeratelyWeightedSpacesAreRegular} yields $\vertiii{\Gamma_{\CalP}}_{\ell_{v}^{q_{2}}\to\ell_{v}^{q_{2}}}\leq C\left(q_{2},N_{\CalP},C_{v,\CalP}\right)$.
\end{proof}
Now, we specialize the above lemma to the case in which $\CalP_{J_{0}}$
is relatively moderate with respect to $\CalQ$. It is worth noting
that the imposed requirements are rather strict. In particular, for
$J_{0}=J$, they can never be fulfilled in case of $\CalO\cap\partial\CalO'\neq\emptyset$,
as a consequence of Lemma~\ref{lem:DifferentOrbitsPreventModerateness}
(with interchanged roles of $\CalQ,\CalP$). Furthermore, we remark
that the following theorem is a slightly improved version of \cite[Theorem 5.3.12]{VoigtlaenderPhDThesis}
from my PhD thesis.
\begin{thm}
\label{thm:NecessaryConditionForModerateCoveringCoarseInFine}Under
the assumptions of Lemma~\ref{lem:NecessaryConjugateCoarseInFineWithoutModerateness},
assume additionally that $Y=\ell_{w}^{q_{1}}\left(I\right)$ and $Z=\ell_{v}^{q_{2}}\left(J\right)$
for certain $q_{1},q_{2}\in\left(0,\infty\right]$ and for certain
weights $w=\left(w_{i}\right)_{i\in I}$ and $v=\left(v_{j}\right)_{j\in J}$
which are $\CalQ$-moderate and $\CalP$-moderate, respectively.

Furthermore, assume that $\CalQ$ is tight and that

\begin{enumerate}
\item \label{enu:NecessaryConditionModerateCoarseInFineCoveringModerate}$\CalP_{J_{0}}$
is relatively $\CalQ$-moderate.
\item \label{enu:NecessaryConditionModerateCoarseInFineWeightModerate}The
weight $v|_{J_{0}}$ is relatively $\CalQ$-moderate\footnote{We recall that relative $\CalQ$-moderateness of $v|_{J_{0}}$ means
that there is some $L=C_{v|_{J_{0}},\CalP,\CalQ}>0$ such that $v_{j}\leq L\cdot v_{\ell}$
holds for all $j,\ell\in J_{0}\cap J_{i}$ for arbitrary $i\in I$.}.
\item \label{enu:NecessaryConditionModerateCoarseInFineReverseSubordinateness}There
is some $r\in\N_{0}$ and some $C_{0}>0$ such that we have
\[
\lambda\left(Q_{i}\right)\leq C_{0}\cdot\lambda\left(\,\smash{\bigcup_{j\in J_{0}\cap J_{i}}}\,\vphantom{\bigcup}P_{j}^{r\ast}\,\right)\vphantom{\bigcup_{j\in J_{0}\cap J_{i}}}\qquad\forall\,i\in I_{0}:=\left\{ i\in I\with J_{0}\cap J_{i}\neq\emptyset\right\} \,.
\]
\end{enumerate}
Let\footnote{Recall from Remark~\ref{rem:SufficientCoarseIntoFineSimplification}
that $\SignedUpperExpo p\in\R\cup\left\{ \infty\right\} $ is defined
by $1/\SignedUpperExpo p=\min\left\{ p^{-1},1-p^{-1}\right\} $.} $s:=\!\left(\frac{1}{q_{2}}-\smash{\frac{1}{\SignedUpperExpo{p_{1}}}}\right)_{\!+}$
and choose for each $i\in I_{0}$ some $j_{i}\in J_{0}\cap J_{i}$.
Then
\begin{equation}
\left\Vert \!\left(\frac{v_{j_{i}}}{w_{i}}\cdot\left|\det T_{i}\right|^{s}\cdot\left|\det S_{j_{i}}\right|^{p_{1}^{-1}-p_{2}^{-1}-s}\right)_{\!\!i\in I_{0}}\right\Vert _{\ell^{q_{2}\cdot\left(q_{1}/q_{2}\right)'}}\leq C\cdot\vertiii{\iota},\label{eq:NecessarySimplifiedConjugateExponentCoarseInFine}
\end{equation}
for some constant
\[
C=C\left(\dimension,r,C_{0},p_{1},p_{2},q_{1},q_{2},\CalQ,\varepsilon_{\CalQ},\CalP,\varepsilon_{\CalP},k\left(\smash{\CalP_{J_{0}}},\CalQ\right),C_{{\rm mod}}\left(\smash{\CalP_{J_{0}}},\CalQ\right),C_{w,\CalQ},C_{v,\CalP},C_{v|_{J_{0}},\CalP,\CalQ},C_{\CalQ,\Phi,p_{1}}\right).
\]
Here, the $L^{p_{1}}$-BAPU $\Phi=\left(\varphi_{i}\right)_{i\in I}$
has to be used to calculate the norm $\vertiii{\iota}$.
\end{thm}

\begin{rem*}
As in the remark after Lemma~\ref{lem:IntersectionCountForModerateCoverings},
we observe that assumption (\ref{enu:NecessaryConditionModerateCoarseInFineReverseSubordinateness})
is automatically satisfied (with $r=0$ and $C_{0}=1$) if we have
$\CalO=\CalO'$ and $J_{0}=J$.
\end{rem*}
\begin{proof}
Lemma~\ref{lem:NecessaryConjugateCoarseInFineWithoutModerateness}
yields a constant $C_{1}>0$ which depends only on quantities mentioned
in the current theorem\footnote{This uses Lemma~\ref{lem:ModeratelyWeightedSpacesAreRegular} to
estimate $\vertiii{\Gamma_{\CalQ}}_{Y\to Y}$ and $\vertiii{\Gamma_{\CalP}}_{Z\to Z}$
for $Y=\ell_{w}^{q_{1}}\left(I\right)$ and $Z=\ell_{v}^{q_{2}}\left(J\right)$
in terms of $C_{w,\CalQ},C_{v,\CalP}$, $N_{\CalQ},N_{\CalP}$ and
$q_{1},q_{2}$. Furthermore, it is used that $C_{Z}$ only depends
on $q_{2}$.} and which satisfies $\vertiii{\gamma}\leq C_{1}\cdot\vertiii{\iota}$,
where
\[
\gamma:\ell_{w}^{q_{1}}\left(I_{0}\right)\hookrightarrow\ell_{u}^{q_{2}}\left(I_{0}\right)\qquad\text{ with }\qquad u_{i}:=\left|\det T_{i}\right|^{p_{1}^{-1}-1}\cdot\left\Vert \left(v_{j}\cdot\left|\det S_{j}\right|^{1-p_{2}^{-1}}\right)_{j\in J_{0}\cap J_{i}}\right\Vert _{\ell^{q_{2}}}\,.
\]
In view of Lemma~\ref{lem:EmbeddingBetweenWeightedSequenceSpaces},
this yields
\begin{equation}
\begin{split} & \left\Vert \left(\frac{\left|\det T_{i}\right|^{p_{1}^{-1}-1}}{w_{i}}\cdot\left\Vert \left(v_{j}\cdot\left|\det S_{j}\right|^{1-p_{2}^{-1}}\right)_{\!j\in J_{0}\cap J_{i}}\right\Vert _{\ell^{q_{2}}}\right)_{\!\!i\in I_{0}}\right\Vert _{\ell^{q_{2}\cdot\left(q_{1}/q_{2}\right)'}}\\
 & =\left\Vert \left(u_{i}/w_{i}\right)_{i\in I_{0}}\right\Vert _{\ell^{q_{2}\cdot\left(q_{1}/q_{2}\right)'}}=\vertiii{\gamma}\leq C_{1}\cdot\vertiii{\iota}\,.
\end{split}
\label{eq:NecessarySimplifiedConjugateExponentLemmaApplication}
\end{equation}

Next, set $C_{2}:=C_{{\rm mod}}\left(\smash{\CalP_{J_{0}}},\CalQ\right)$
and $L:=C_{v|_{J_{0}},\CalP,\CalQ}>0$, so that $v_{j}\leq L\cdot v_{\ell}$
holds for all $j,\ell\in J_{0}\cap J_{i}$ for arbitrary $i\in I$.
In particular, this implies $v_{j}\geq L^{-1}\cdot v_{j_{i}}$ for
all $j\in J_{0}\cap J_{i}$ and arbitrary $i\in I_{0}$. In addition,
by choice of $C_{2}$, we have
\begin{equation}
C_{2}^{-1}\cdot\left|\det S_{j}\right|\leq\left|\det S_{j_{i}}\right|\leq C_{2}\cdot\left|\det S_{j}\right|\qquad\forall\,i\in I_{0}\text{ and }j\in J_{0}\cap J_{i},\label{eq:NecessaryConjugateModerateModeratenessAssumption}
\end{equation}
since $j_{i}\in J_{0}\cap J_{i}$ as well. But this yields a constant
$C_{3}=C_{3}\left(C_{{\rm mod}}\left(\smash{\CalP_{J_{0}}},\CalQ\right),p_{2}\right)\geq1$
with
\begin{equation}
\left|\det S_{j}\right|^{1-p_{2}^{-1}}\geq C_{3}^{-1}\cdot\left|\det S_{j_{i}}\right|^{1-p_{2}^{-1}}\qquad\forall\,i\in I_{0}\text{ and }j\in J_{0}\cap J_{i}.\label{eq:NecessaryConjugateModeratenessAssumptionPowed}
\end{equation}
Note that this remains true also in case of $1-p_{2}^{-1}<0$.

Next, note that the assumptions of the present theorem include the
prerequisites of Lemma~\ref{lem:IntersectionCountForModerateCoverings}
(with interchanged roles of $\CalQ$ and $\CalP$), so that we get
\begin{equation}
\left|\det T_{i}\right|\big/\left|\det S_{j_{i}}\right|\leq C_{4}\cdot\left|J_{0}\cap J_{i}\right|\qquad\forall\,i\in I_{0}\,,\label{eq:NecessaryConjugateExponentCoarseInFineIntersectionCount}
\end{equation}
for a suitable constant $C_{4}=C_{4}\left(C_{0},\dimension,r,\CalQ,\CalP,\varepsilon_{\CalQ},C_{{\rm mod}}\left(\smash{\CalP_{J_{0}}},\CalQ\right)\right)>0$.
Combined with equation~(\ref{eq:NecessaryConjugateModerateModeratenessAssumption}),
this implies
\begin{align*}
\left\Vert \left(v_{j}\cdot\left|\det S_{j}\right|^{1-p_{2}^{-1}}\right)_{j\in J_{0}\cap J_{i}}\right\Vert _{\ell^{q_{2}}} & \geq C_{3}^{-1}L^{-1}\cdot\left|\det S_{j_{i}}\right|^{1-p_{2}^{-1}}v_{j_{i}}\cdot\left|J_{0}\cap J_{i}\right|^{q_{2}^{-1}}\\
 & \geq\left(C_{3}C_{4}^{q_{2}^{-1}}L\right)^{-1}\cdot\left|\det S_{j_{i}}\right|^{1-p_{2}^{-1}-q_{2}^{-1}}\cdot\left|\det T_{i}\right|^{q_{2}^{-1}}\cdot v_{j_{i}}\qquad\forall\,i\in I_{0}\,.
\end{align*}

Let us set $C_{5}:=C_{1}C_{3}C_{4}^{q_{2}^{-1}}L$. By combining the
estimate above with equation~(\ref{eq:NecessarySimplifiedConjugateExponentLemmaApplication}),
we conclude
\begin{equation}
\left\Vert \left(\frac{v_{j_{i}}}{w_{i}}\cdot\left|\det T_{i}\right|^{q_{2}^{-1}+p_{1}^{-1}-1}\left|\det S_{j_{i}}\right|^{1-p_{2}^{-1}-q_{2}^{-1}}\right)_{i\in I_{0}}\right\Vert _{\ell^{q_{2}\cdot\left(q_{1}/q_{2}\right)'}}\leq C_{5}\cdot\vertiii{\iota}.\label{eq:NecessaryConjugateExponentCoarseInFineAlmostDone}
\end{equation}

\medskip{}

With this preparation, we can deduce the actual claim of the theorem.
For technical reasons, however, we need to distinguish three cases:

\textbf{Case 1}: $p_{1}\in\left(0,1\right)$. In this case, we have
\[
\frac{1}{\SignedUpperExpo{p_{1}}}=\min\left\{ \frac{1}{p_{1}},1-\frac{1}{p_{1}}\right\} =1-\frac{1}{p_{1}}<0\quad\quad\text{and hence}\quad\quad s=\left(\frac{1}{q_{2}}-\smash{\frac{1}{\SignedUpperExpo{p_{1}}}}\right)_{+}=\frac{1}{q_{2}}+\frac{1}{p_{1}}-1\,.
\]
Thus,
\[
\left|\det T_{i}\right|^{s}\cdot\left|\det S_{j_{i}}\right|^{p_{1}^{-1}-p_{2}^{-1}-s}=\left|\det T_{i}\right|^{q_{2}^{-1}+p_{1}^{-1}-1}\cdot\left|\det S_{j_{i}}\right|^{1-p_{2}^{-1}-q_{2}^{-1}},
\]
so that the desired estimate~(\ref{eq:NecessarySimplifiedConjugateExponentCoarseInFine})
is a direct consequence of equation~(\ref{eq:NecessaryConjugateExponentCoarseInFineAlmostDone}).

\medskip{}

\textbf{Case 2}: $p_{1}\in\left[1,2\right]$ and $q_{2}\leq p_{1}'$.
Here, we have 
\[
\frac{1}{\SignedUpperExpo{p_{1}}}=\min\left\{ \frac{1}{p_{1}},1-\frac{1}{p_{1}}\right\} =1-\frac{1}{p_{1}}=\frac{1}{p_{1}'}\quad\quad\text{and hence}\quad\quad\frac{1}{q_{2}}-\frac{1}{\SignedUpperExpo{p_{1}}}=\frac{1}{q_{2}}-\frac{1}{p_{1}'}\geq0\,,
\]
which yields
\[
s=\left(\frac{1}{q_{2}}-\smash{\frac{1}{\SignedUpperExpo{p_{1}}}}\right)_{+}=\frac{1}{q_{2}}-\frac{1}{p_{1}'}=\frac{1}{q_{2}}+\frac{1}{p_{1}}-1.
\]
Thus, we see exactly as in the previous case that the desired estimate
is a direct consequence of equation~(\ref{eq:NecessaryConjugateExponentCoarseInFineAlmostDone}).

\medskip{}

\textbf{Case 3}: We have $p_{1}\in\left[1,2\right]$ and $q_{2}>p_{1}'$
or we have $p_{1}\in\left[2,\infty\right]$. In this case, we will
\emph{not} use equation~(\ref{eq:NecessaryConjugateExponentCoarseInFineAlmostDone}).
Instead, we will invoke Theorem~\ref{thm:BurnerNecessaryConditionCoarseInFine}.
As a preparation, set $t:=q_{2}\cdot\left(p_{1}/q_{2}\right)'$ and
note that the exponent $s=\left(\frac{1}{q_{2}}-\smash{\frac{1}{\SignedUpperExpo{p_{1}}}}\right)_{\!+}$
satisfies
\begin{equation}
s=\frac{1}{q_{2}\cdot\left(p_{1}/q_{2}\right)'}=\frac{1}{t}.\label{eq:NecessaryConditionForModerateCoarseInFineSpecialCaseExponentForm}
\end{equation}
To see this, we first note that equation~(\ref{eq:InverseOfSpecialExponent})
shows $\frac{1}{t}=\left(\frac{1}{q_{2}}-\frac{1}{p_{1}}\right)_{+}$.
To complete the proof of equation~(\ref{eq:NecessaryConditionForModerateCoarseInFineSpecialCaseExponentForm}),
we distinguish two sub-cases:

\begin{casenv}
\item We have $p_{1}\in\left[1,2\right]$ and $q_{2}>p_{1}'$. In this case,
note $\SignedUpperExpo{p_{1}}=\UpperExpo{p_{1}}=\max\left\{ p_{1},p_{1}'\right\} =p_{1}'$
and hence $\frac{1}{q_{2}}-\frac{1}{\SignedUpperExpo{p_{1}}}=\frac{1}{q_{2}}-\frac{1}{p_{1}'}<0$,
since $q_{2}>p_{1}'$. This implies $s=0$. But since $p_{1}\in\left[1,2\right]$,
we also have $p_{1}\leq p_{1}'$ and thus $\frac{1}{q_{2}}-\frac{1}{p_{1}}\leq\frac{1}{q_{2}}-\frac{1}{p_{1}'}<0$,
so that $s=0=\left(\frac{1}{q_{2}}-\frac{1}{p_{1}}\right)_{+}=\frac{1}{t}$,
as desired.
\item We have $p_{1}\in\left[2,\infty\right]$. In this case, note $\SignedUpperExpo{p_{1}}=\UpperExpo{p_{1}}=\max\left\{ p_{1},p_{1}'\right\} =p_{1}$,
which entails $s=\left(\frac{1}{q_{2}}-\smash{\frac{1}{\SignedUpperExpo{p_{1}}}}\right)_{+}=\left(\frac{1}{q_{2}}-\frac{1}{p_{1}}\right)_{+}=\frac{1}{t}$,
as claimed.
\end{casenv}
Recall from Lemma~\ref{lem:NecessaryConjugateCoarseInFineWithoutModerateness}
that $k=k\left(\smash{\CalP_{J_{0}}},\CalQ\right)$. Now, to show
that Theorem~\ref{thm:BurnerNecessaryConditionCoarseInFine} is applicable,
note for $j\in J_{0}$ that we have $\emptyset\neq P_{j}\subset Q_{i_{j}}^{k\ast}\subset\CalO$
for some $i_{j}\in I$. In particular, $P_{j}\cap Q_{i}\neq\emptyset$
for some $i\in I$, which implies $i\in I_{0}$ and (by Lemma~\ref{lem:SubordinatenessImpliesWeakSubordination})
that $P_{j}\subset Q_{i}^{\left(2k+2\right)\ast}\subset\overline{Q_{i}^{\left(2k+3\right)\ast}}$.
Thus, using the notation of Lemma~\ref{lem:NecessaryConjugateCoarseInFineWithoutModerateness},
we have 
\[
K=\bigcup_{i\in I_{0}}\overline{Q_{i}^{\left(2k+3\right)\ast}}\supset\bigcup_{j\in J_{0}}P_{j}.
\]

In view of the assumed boundedness of $\iota$ (see Lemma~\ref{lem:NecessaryConjugateCoarseInFineWithoutModerateness}),
we thus see that the prerequisites of Theorem~\ref{thm:BurnerNecessaryConditionCoarseInFine}
are satisfied (with $\vertiii{\iota_{0}}\leq\vertiii{\iota}$, where
$\iota_{0}$ denotes the embedding from Theorem~\ref{thm:BurnerNecessaryConditionCoarseInFine}).
Since $Z=\ell_{v}^{q_{2}}\left(J\right)$ satisfies the Fatou property,
since the triangle inequality of $Z$ only depends on $q_{2}$ and
since Lemma~\ref{lem:ModeratelyWeightedSpacesAreRegular} shows $\vertiii{\Gamma_{\CalP}}_{Z\to Z}\leq C_{v,\CalP}\cdot N_{\CalP}^{\max\left\{ 1,q_{2}^{-1}\right\} }$,
we conclude that there is a constant
\[
C_{6}=C_{6}\left(\dimension,k\left(\smash{\CalP_{J_{0}}},\CalQ\right),p_{1},p_{2},q_{2},\CalQ,\CalP,\varepsilon_{\CalP},C_{v,\CalP},C_{\CalQ,\Phi,p_{1}}\right)
\]
satisfying $\vertiii{\eta}\leq C_{6}\cdot\vertiii{\iota}$, for
\[
\eta:\ell_{w}^{q_{1}}\left(\left[\ell^{p_{1}}\left(J_{i}\cap J_{0}\right)\right]_{i\in I}\right)\hookrightarrow\ell_{v_{j}\cdot\left|\det S_{j}\right|^{p_{1}^{-1}-p_{2}^{-1}}}^{q_{2}}\left(J_{0}\right).
\]

Now, an application of Corollary~\ref{cor:EmbeddingCoarseIntoFineSimplification}
(with $r=p_{1}$, $u^{\left(1\right)}\equiv1$ and $u^{\left(2\right)}\equiv1$
and finally with $\bigl(v_{j}\cdot\left|\det S_{j}\right|^{p_{1}^{-1}-p_{2}^{-1}}\bigr)_{j\in J}$
instead of $v$) yields a constant $C_{7}=C_{7}\left(p_{1},q_{1},q_{2},\CalQ,C_{w,\CalQ},k\left(\smash{\CalP_{J_{0}}},\CalQ\right)\right)$
such that
\[
\left\Vert \left(w_{i}^{-1}\cdot\left\Vert \left(v_{j}\cdot\left|\det S_{j}\right|^{p_{1}^{-1}-p_{2}^{-1}}\right)_{j\in J_{0}\cap J_{i}}\right\Vert _{\ell^{q_{2}\cdot\left(p_{1}/q_{2}\right)'}}\right)_{i\in I}\right\Vert _{\ell^{q_{2}\cdot\left(q_{1}/q_{2}\right)'}}\leq C_{7}\cdot\vertiii{\eta}\leq C_{6}C_{7}\cdot\vertiii{\iota}.
\]
Now, the same arguments as in the previous cases—using the relative
$\CalQ$-moderateness of $\CalP_{J_{0}}$ and of $v|_{J_{0}}$ (see
in particular equation~(\ref{eq:NecessaryConjugateExponentCoarseInFineIntersectionCount}))—show
\begin{align*}
v_{j_{i}}\cdot\left|\det S_{j_{i}}\right|^{p_{1}^{-1}-p_{2}^{-1}-s}\cdot\left|\det T_{i}\right|^{s} & \leq C_{4}^{s}\cdot v_{j_{i}}\cdot\left|\det S_{j_{i}}\right|^{p_{1}^{-1}-p_{2}^{-1}}\cdot\left|J_{0}\cap J_{i}\right|^{s}\\
\left({\scriptstyle \text{since }s=1/t}\right) & =C_{4}^{1/t}\cdot v_{j_{i}}\cdot\left|\det S_{j_{i}}\right|^{p_{1}^{-1}-p_{2}^{-1}}\cdot\left|J_{0}\cap J_{i}\right|^{1/t}\\
 & \leq C_{4}^{1/t}C_{8}\cdot\left\Vert \left(v_{j}\cdot\left|\det S_{j}\right|^{p_{1}^{-1}-p_{2}^{-1}}\right)_{j\in J_{0}\cap J_{i}}\right\Vert _{\ell^{t}}\\
 & =C_{4}^{1/t}C_{8}\cdot\left\Vert \left(v_{j}\cdot\left|\det S_{j}\right|^{p_{1}^{-1}-p_{2}^{-1}}\right)_{j\in J_{0}\cap J_{i}}\right\Vert _{\ell^{q_{2}\cdot\left(p_{1}/q_{2}\right)'}}
\end{align*}
for some constant $C_{8}=C_{8}\left(p_{1},p_{2},C_{v|_{J_{0}},\CalP,\CalQ},C_{{\rm mod}}\left(\smash{\CalP_{J_{0}}},\CalQ\right)\right)=C_{8}\left(p_{1},p_{2},C_{2},L\right)$
and all $i\in I_{0}$.

All in all, we conclude as desired that
\begin{align*}
 & \left\Vert \left(\frac{v_{j_{i}}}{w_{i}}\cdot\left|\det S_{j_{i}}\right|^{p_{1}^{-1}-p_{2}^{-1}-s}\cdot\left|\det T_{i}\right|^{s}\right)_{i\in I_{0}}\right\Vert _{\ell^{q_{2}\cdot\left(q_{1}/q_{2}\right)'}}\\
\leq & C_{4}^{1/t}C_{8}\cdot\left\Vert \left(w_{i}^{-1}\cdot\left\Vert \left(v_{j}\cdot\left|\det S_{j}\right|^{p_{1}^{-1}-p_{2}^{-1}}\right)_{j\in J_{0}\cap J_{i}}\right\Vert _{\ell^{q_{2}\cdot\left(p_{1}/q_{2}\right)'}}\right)_{i\in I_{0}}\right\Vert _{\ell^{q_{2}\cdot\left(q_{1}/q_{2}\right)'}}\\
\leq & C_{4}^{1/t}C_{8}\cdot\left\Vert \left(w_{i}^{-1}\cdot\left\Vert \left(v_{j}\cdot\left|\det S_{j}\right|^{p_{1}^{-1}-p_{2}^{-1}}\right)_{j\in J_{0}\cap J_{i}}\right\Vert _{\ell^{q_{2}\cdot\left(p_{1}/q_{2}\right)'}}\right)_{i\in I}\right\Vert _{\ell^{q_{2}\cdot\left(q_{1}/q_{2}\right)'}}\\
\leq & C_{4}^{1/t}C_{6}C_{7}C_{8}\cdot\vertiii{\iota}\,\,.\qedhere
\end{align*}
\end{proof}
With similar techniques, we will now prove the analogous result for
the ``reverse'' case in which (a subfamily of) $\CalQ$ is almost
subordinate to and relatively moderate with respect to $\CalP$. But
as above, we first establish an auxiliary result for which we do \emph{not}
need to assume that $\CalQ$ is relatively moderate with respect to
$\CalP$. In contrast to Lemma~\ref{lem:NecessaryConjugateCoarseInFineWithoutModerateness},
however, we assume the ``global'' components $Y,Z$ to be weighted
$\ell^{q}$ spaces, since the proof will use (a slight variant of)
the ``reproducing'' property $\ell_{w}^{q_{1}}\left(I\right)=\ell^{q_{1}}\bigl(\left[\ell_{w}^{q_{1}}\left(I^{\left(k\right)}\right)\right]_{k\in K}\bigr)$,
which holds for arbitrary partitions $I=\biguplus_{k\in K}I^{\left(k\right)}$.
\begin{lem}
\label{lem:NecessaryConditionConjugateFineInCoarseWithoutModerateness}Let
$\emptyset\neq\CalO,\CalO'\subset\R^{\dimension}$ be open, let $\CalQ=\left(Q_{i}\right)_{i\in I}=\left(T_{i}Q_{i}'+b_{i}\right)_{i\in I}$
be a semi-structured $L^{p_{1}}$-decomposition covering of $\CalO$
and let $\CalP=\left(P_{j}\right)_{j\in J}=\left(\smash{S_{j}P_{j}'+c_{j}}\right)_{j\in J}$
be a semi-structured $L^{p_{2}}$-decomposition covering of $\CalO'$,
for certain $p_{1},p_{2}\in\left(0,\infty\right]$ . Let $w=\left(w_{i}\right)_{i\in I}$
and $v=\left(v_{j}\right)_{j\in J}$ be $\CalQ$-moderate and $\CalP$-moderate,
respectively. Finally, let $q_{1},q_{2}\in\left(0,\infty\right]$.

Choose an arbitrary subset $J_{0}\subset J$ such that there is\footnote{The easiest case in which this assumption is fulfilled is if $\CalP$
is tight with $P_{j}\subset\CalO$ for all $j\in J_{0}$. In this
case, one can simply take $\varepsilon=\varepsilon_{\CalP}$.} some $\varepsilon>0$ and for every $j\in J_{0}$ some $\xi_{j}\in\R^{\dimension}$
with $B_{\varepsilon}\left(\xi_{j}\right)\subset P_{j}'$ and $S_{j}\left[B_{\varepsilon}\left(\xi_{j}\right)\right]+c_{j}\subset\CalO$.

Define 
\[
I_{0}:=\left\{ i\in I\with J_{i}\cap J_{0}\neq\emptyset\right\} 
\]
and assume that $\CalQ_{I_{0}}:=\left(Q_{i}\right)_{i\in I_{0}}$
is almost subordinate to $\CalP$ and that\footnote{Note that each $j\in J_{0}$ satisfies $\emptyset\neq S_{j}\left[B_{\varepsilon}\left(\xi_{j}\right)\right]+c_{j}\subset P_{j}\cap\CalO$
by assumption. Hence, there is some $i\in I$ with $Q_{i}\cap P_{j}\neq\emptyset$,
which yields $i\in I_{0}\cap I_{j}\neq\emptyset$. In particular,
the quantity in equation~(\ref{eq:NecessaryConditionConjugateFineInCoarseSpecialAssumption})
is always positive. Finally, the equality in (\ref{eq:NecessaryConditionConjugateFineInCoarseSpecialAssumption})
is justified by equation~(\ref{eq:NecessaryConditionModerateFineInCoarseIntersectionSetContainedInI0})
below, since this equation yields $I_{j}\subset I_{0}$ and hence
$I_{0}\cap I_{j}=I_{j}$ for all $j\in J_{0}$.}
\begin{equation}
u_{j}:=\left\Vert \left(w_{i}\cdot\left|\det T_{i}\right|^{1-p_{1}^{-1}}\right)_{i\in I_{0}\cap I_{j}}\right\Vert _{\ell^{q_{1}}}=\left\Vert \left(w_{i}\cdot\left|\det T_{i}\right|^{1-p_{1}^{-1}}\right)_{i\in I_{j}}\right\Vert _{\ell^{q_{1}}}<\infty\qquad\forall\,j\in J_{0}\,.\label{eq:NecessaryConditionConjugateFineInCoarseSpecialAssumption}
\end{equation}

Finally, set
\[
K:=\bigcup_{i\in I_{0}}Q_{i}\subset\CalO\cap\CalO'
\]
and—with $\CalD_{K}:=\left\{ f\in\TestFunctionSpace{\R^{\dimension}}\with\supp f\subset K\right\} $—assume
that the identity map
\[
\iota:\left(\CalD_{K},\left\Vert \mybullet\right\Vert _{\FourierDecompSp{\CalQ}{p_{1}}{\ell_{w}^{q_{1}}}}\right)\hookrightarrow\FourierDecompSp{\CalP}{p_{2}}{\ell_{v}^{q_{2}}},f\mapsto f
\]
is bounded\footnote{Here, $\CalD_{K}\subset\TestFunctionSpace{\CalO\cap\CalO'}$. Note
that Lemma~\ref{lem:PartitionCoveringNecessary} shows that $\left(\varphi_{i}\right)_{i\in I}$
is a locally finite partition of unity on $\CalO$. Hence, for each
$g\in\CalD_{K}$, we can have $\varphi_{i}g\not\equiv0$ only for
finitely many $i\in I$. Since $\ell_{w}^{q_{1}}\left(I\right)$ contains
all finitely supported sequences, we easily see $\CalD_{K}\subset\FourierDecompSp{\CalQ}{p_{1}}{\ell_{w}^{q_{1}}}$.
An analogous argument shows $\CalD_{K}\subset\FourierDecompSp{\CalP}{p_{2}}{\ell_{v}^{q_{2}}}$,
so that $\iota$ is always well-defined, but not necessarily bounded.}. Then we have
\[
\left\Vert \left(v_{j}\cdot\left|\det S_{j}\right|^{1-p_{2}^{-1}}\big/\,u_{j}\right)_{j\in J_{0}}\right\Vert _{\ell^{q_{2}\cdot\left(q_{1}/q_{2}\right)'}}\leq C\cdot\vertiii{\iota}
\]
for some constant
\[
C=C\left(\dimension,p_{1},p_{2},q_{2},\varepsilon,k\left(\smash{\CalQ_{I_{0}}},\CalP\right),\CalQ,\CalP,C_{v,\CalP},C_{\CalQ,\Phi,p_{1}},C_{\CalP,\Psi,p_{2}}\right).
\]
Here, the $L^{p_{1}}$/$L^{p_{2}}$ BAPUs $\Phi=\left(\varphi_{i}\right)_{i\in I}$
and $\Psi=\left(\psi_{j}\right)_{j\in J}$ are those which are used
to calculate $\vertiii{\iota}$.
\end{lem}

\begin{proof}
As usual, we begin by setting up several auxiliary objects and by
deriving some of their properties. For brevity, set $k:=k\left(\smash{\CalQ_{I_{0}}},\CalP\right)$.

Using the inclusion $S_{j}\left[B_{\varepsilon}\left(\xi_{j}\right)\right]+c_{j}\subset P_{j}\cap\CalO$,
one can use exactly the same construction as in the proof of Lemma~\ref{lem:NormOfClusteredBAPUAndTestFunctionBuildingBlocks}
(take a nontrivial $\gamma\in\TestFunctionSpace{B_{1}\left(0\right)}$
and define $\gamma_{j}:=L_{c_{j}}\left(L_{\xi_{j}}\left[\gamma\left(\varepsilon^{-1}\mybullet\right)\right]\circ S_{j}^{-1}\right)$
for $j\in J_{0}$) to obtain a family $\left(\gamma_{j}\right)_{j\in J_{0}}$
satisfying $\gamma_{j}\in\TestFunctionSpace{P_{j}\cap\CalO}\subset\TestFunctionSpace{\CalO\cap\CalO'}$
and
\begin{equation}
\left\Vert \Fourier^{-1}\gamma_{j}\right\Vert _{L^{p}}=C_{1}^{\left(p\right)}\cdot\left|\det S_{j}\right|^{1-\frac{1}{p}}\qquad\forall\,j\in J_{0}\text{ and }p\in\left(0,\infty\right]\,,\label{eq:NecessaryConditionModerateFineInCoarseGammaNorm}
\end{equation}
with a constant  $C_{1}^{\left(p\right)}=C_{1}^{\left(p\right)}\left(\dimension,\varepsilon\right)>0$.

\medskip{}

We start with a few technical observations which will be useful later
on. We first note
\begin{equation}
K':=\bigcup_{j\in J_{0}}\left(P_{j}\cap\CalO\right)\subset K.\label{eq:NecessaryConditionModerateFineInCoarseTestSetInclusion}
\end{equation}
To see this, let $\xi\in K'$ be arbitrary. Hence, $\xi\in P_{j}\cap\CalO$
for some $j\in J_{0}$ and thus $\xi\in Q_{i}$ for some $i\in I$,
since $\CalQ$ covers $\CalO$. This entails $\xi\in Q_{i}\cap P_{j}\neq\emptyset$
and hence $j\in J_{i}\cap J_{0}\neq\emptyset$, which implies $i\in I_{0}$
and thus finally $\xi\in Q_{i}\subset\bigcup_{\ell\in I_{0}}Q_{\ell}=K$.

Furthermore, we have 
\begin{equation}
I_{j}\subset I_{0}\text{ for all }j\in J_{0},\label{eq:NecessaryConditionModerateFineInCoarseIntersectionSetContainedInI0}
\end{equation}
because for arbitrary $i\in I_{j}$ we have $j\in J_{i}\cap J_{0}\neq\emptyset$
and hence $i\in I_{0}$.

Next, we note 
\begin{equation}
\left|J_{i}\right|\leq N_{\CalP}^{k+1}\text{ for all }i\in I_{0}.\label{eq:NecessaryConditionModerateCoveringFineInCoarseJClusterBounded}
\end{equation}
To see that this is true, note that $j\in J_{i}$ entails $\emptyset\neq Q_{i}\cap P_{j}\subset P_{j_{i}}^{k\ast}\cap P_{j}$
for some $j_{i}\in J$ and thus $j\in j_{i}^{\left(k+1\right)\ast}$,
which means $J_{i}\subset j_{i}^{\left(k+1\right)\ast}$. Using Lemma~\ref{lem:SemiStructuredClusterInvariant},
this shows that equation~(\ref{eq:NecessaryConditionModerateCoveringFineInCoarseJClusterBounded})
is true.

Finally, set $r_{0}:=N_{\CalP}^{2\cdot1+1}=N_{\CalP}^{3}$, so that
Lemma~\ref{lem:DisjointizationPrinciple} ensures existence of a
partition $J=\biguplus_{r=1}^{r_{0}}J^{\left(r\right)}$ with $P_{j}^{\ast}\cap P_{m}^{\ast}=\emptyset$
for all $j,m\in J^{\left(r\right)}$ with $j\neq m$ and arbitrary
$r\in\underline{r_{0}}$. 

\medskip{}

Now we properly start the proof: Let $d=\left(d_{j}\right)_{j\in J_{0}}\in\ell_{0}\left(J_{0}\right)$
be arbitrary. Fix $r\in\underline{r_{0}}$ and define
\[
g^{\left(r\right)}:=\sum_{j\in J_{0}\cap J^{\left(r\right)}}d_{j}\cdot\gamma_{j}.
\]
Equation~(\ref{eq:NecessaryConditionModerateFineInCoarseTestSetInclusion})
and compactness of $\supp\gamma_{j}\subset P_{j}\cap\CalO\subset K'$
show that $\supp g^{\left(r\right)}\subset K\subset\CalO\cap\CalO'$
is compact, since the sum defining $g^{\left(r\right)}$ is finite.
This easily yields $g^{\left(r\right)}\in\CalD_{K}\subset\FourierDecompSp{\CalQ}{p_{1}}{\ell_{w}^{q_{1}}}$,
as well as $g^{\left(r\right)}\in\CalD_{K}\subset\FourierDecompSp{\CalP}{p_{2}}{\ell_{v}^{q_{2}}}$.
As usual, we now derive an upper bound on $\left\Vert g^{\left(r\right)}\right\Vert _{\FourierDecompSp{\CalQ}{p_{1}}{\ell_{w}^{q_{1}}}}$
and a lower bound on $\left\Vert g^{\left(r\right)}\right\Vert _{\FourierDecompSp{\CalP}{p_{2}}{\ell_{v}^{q_{2}}}}$.

\medskip{}

Let us begin with the lower bound: For $j\in J_{0}\cap J^{\left(r\right)}$
and $m\in J^{\left(r\right)}$ with $\psi_{m}^{\ast}\cdot\gamma_{j}\not\equiv0$,
we have 
\[
\emptyset\neq P_{m}^{\ast}\cap\supp\gamma_{j}\subset P_{m}^{\ast}\cap P_{j}\subset P_{m}^{\ast}\cap P_{j}^{\ast}
\]
and thus $m=j$ by choice of $J^{\left(r\right)}$. But Lemma~\ref{lem:PartitionCoveringNecessary}
shows $\psi_{m}^{\ast}\equiv1$ on $P_{m}=P_{j}\supset\supp\gamma_{j}$
and thus $\psi_{m}^{\ast}\cdot\gamma_{j}=\gamma_{j}=\gamma_{m}$.
In summary, we conclude
\begin{equation}
\psi_{m}^{\ast}\cdot g^{\left(r\right)}=d_{m}\cdot\gamma_{m}\qquad\forall\,m\in J_{0}\cap J^{\left(r\right)}.\label{eq:NecessaryConditionModerateFineInCoarseClusteredBAPUSelectsGammaM}
\end{equation}

Next, note that Remark~\ref{rem:ClusteredBAPUYIeldsBoundedControlSystem}
shows that $\Lambda:=\left(\psi_{m}^{\ast}\right)_{m\in J}$ is an
$L^{p_{2}}$-bounded (control) system for $\CalP$ with $C_{\CalP,\Lambda,p_{2}}\leq C_{2}=C_{2}\left(\dimension,p_{2},\CalP,C_{\CalP,\Psi,p_{2}}\right)$
and $\ell_{\Lambda,\CalP}=1$. In combination with Theorem~\ref{thm:BoundedControlSystemEquivalentQuasiNorm}
and with the estimate $\vertiii{\Gamma_{\CalP}}_{\ell_{v}^{q_{2}}\to\ell_{v}^{q_{2}}}\leq C_{v,\CalP}\cdot N_{\CalP}^{\max\left\{ 1,q_{2}^{-1}\right\} }$
from Lemma~\ref{lem:ModeratelyWeightedSpacesAreRegular}, we get
\begin{align}
C_{2}C_{3}\cdot\left\Vert \smash{g^{\left(r\right)}}\right\Vert _{\BAPUFourierDecompSp{\CalP}{p_{2}}{\ell_{v}^{q_{2}}}{\Psi}} & \geq\left\Vert \left(\left\Vert \Fourier^{-1}\left(\psi_{m}^{\ast}\cdot\smash{g^{\left(r\right)}}\right)\right\Vert _{L^{p_{2}}}\right)_{m\in J}\right\Vert _{\ell_{v}^{q_{2}}}\nonumber \\
 & \geq\left\Vert \left(\left\Vert \Fourier^{-1}\left(\psi_{m}^{\ast}\cdot\smash{g^{\left(r\right)}}\right)\right\Vert _{L^{p_{2}}}\right)_{m\in J_{0}\cap J^{\left(r\right)}}\right\Vert _{\ell_{v}^{q_{2}}}\nonumber \\
\left({\scriptstyle \text{eq. }\eqref{eq:NecessaryConditionModerateFineInCoarseClusteredBAPUSelectsGammaM}}\right) & =\left\Vert \left(\left\Vert \Fourier^{-1}\left(d_{m}\cdot\gamma_{m}\right)\right\Vert _{L^{p_{2}}}\right)_{m\in J_{0}\cap J^{\left(r\right)}}\right\Vert _{\ell_{v}^{q_{2}}}\nonumber \\
\left({\scriptstyle \text{eq. }\eqref{eq:NecessaryConditionModerateFineInCoarseGammaNorm}}\right) & =C_{1}^{\left(p_{2}\right)}\cdot\left\Vert \left(\left|\det S_{m}\right|^{1-p_{2}^{-1}}\cdot d_{m}\right)_{m\in J_{0}\cap J^{\left(r\right)}}\right\Vert _{\ell_{v}^{q_{2}}}\label{eq:NecessaryModerateFineInCoarseTestFunctionEstimateCoarseCovering}
\end{align}
for some constant $C_{3}=C_{3}\left(\dimension,q_{2},p_{2},\CalP,C_{v,\CalP}\right)>0$.

\medskip{}

Next, we want to establish an upper bound for $\left\Vert \smash{g^{\left(r\right)}}\right\Vert _{\BAPUFourierDecompSp{\CalQ}{p_{1}}{\ell_{w}^{q_{1}}}{\Phi}}$.
To this end, consider any $i\in I$ with
\[
0\not\equiv\varphi_{i}\cdot g^{\left(r\right)}=\sum_{j\in J_{0}\cap J^{\left(r\right)}}\left[d_{j}\cdot\varphi_{i}\cdot\gamma_{j}\right].
\]
Then, there is some $j\in J_{0}\cap J^{\left(r\right)}$ satisfying
$0\not\equiv\varphi_{i}\cdot\gamma_{j}$. But for \emph{any(!)\@}
such $j\in J_{0}\cap J^{\left(r\right)}$, we get $Q_{i}\cap P_{j}\neq\emptyset$,
which leads to $i\in I_{j}\subset I_{0}$ (see equation~(\ref{eq:NecessaryConditionModerateFineInCoarseIntersectionSetContainedInI0})),
or equivalently $i\in I_{0}$ and $j\in J_{i}$. Hence, $i\in I_{0}$,
and
\begin{equation}
\varphi_{i}\cdot g^{\left(r\right)}=\sum_{j\in J_{0}\cap J^{\left(r\right)}\cap J_{i}}\left[d_{j}\cdot\varphi_{i}\cdot\gamma_{j}\right]\,.\label{eq:NecessaryModerateFineInCoarseFineLocalization}
\end{equation}

Now, we distinguish the cases $p_{1}\in\left[1,\infty\right]$ and
$p_{1}\in\left(0,1\right)$. For $p_{1}\in\left[1,\infty\right]$,
we use equation~(\ref{eq:NecessaryModerateFineInCoarseFineLocalization}),
Young's inequality ($L^{1}\ast L^{p_{1}}\hookrightarrow L^{p_{1}}$)
and the triangle inequality for $\left\Vert \mybullet\right\Vert _{L^{p_{1}}}$,
to derive
\begin{align*}
w_{i}\cdot\left\Vert \Fourier^{-1}\left(\varphi_{i}\cdot\smash{g^{\left(r\right)}}\right)\right\Vert _{L^{p_{1}}} & =w_{i}\cdot\left\Vert \,\smash{\sum_{j\in J_{0}\cap J_{i}\cap J^{\left(r\right)}}}\vphantom{\sum}\Fourier^{-1}\left(d_{j}\cdot\varphi_{i}\cdot\gamma_{j}\right)\,\right\Vert _{L^{p_{1}}}\vphantom{\sum_{j\in J_{0}\cap J_{i}\cap J^{\left(r\right)}}}\\
 & \leq w_{i}\cdot\vphantom{\sum_{j\in J_{0}\cap J_{i}\cap J^{\left(r\right)}}^{T}}\sum_{j\in J_{0}\cap J_{i}\cap J^{\left(r\right)}}\left(\vphantom{I^{\left(r\right)}}\left|d_{j}\right|\cdot\left\Vert \Fourier^{-1}\gamma_{j}\right\Vert _{L^{1}}\cdot\left\Vert \Fourier^{-1}\varphi_{i}\right\Vert _{L^{p_{1}}}\right)\\
\left({\scriptstyle \text{equations }\eqref{eq:NecessaryConditionModerateFineInCoarseGammaNorm}\text{ and }\eqref{eq:UpperBAPUNormEstimate}}\right) & \leq C_{1}^{\left(1\right)}C_{4}\cdot\left|\det T_{i}\right|^{1-p_{1}^{-1}}\cdot w_{i}\cdot\vphantom{\sum_{j\in J_{0}\cap J_{i}\cap J^{\left(r\right)}}^{T}}\sum_{j\in J_{0}\cap J_{i}\cap J^{\left(r\right)}}\left|d_{j}\right|,
\end{align*}
where the left-hand side of the estimate vanishes for $i\in I\setminus I_{0}$,
as the discussion above showed. In the previous estimate, the last
step used that Lemma~\ref{lem:NormOfClusteredBAPUAndTestFunctionBuildingBlocks}
yields 
\begin{equation}
\left\Vert \Fourier^{-1}\varphi_{i}\right\Vert _{L^{p_{1}}}\leq C_{4}\cdot\left|\det T_{i}\right|^{1-p_{1}^{-1}}\qquad\forall\,i\in I\,,\label{eq:NecessaryModerateFineInCoarseFineBAPUNormEstimate}
\end{equation}
with a constant $C_{4}=C_{4}\left(\dimension,p_{1},\CalQ,C_{\CalQ,\Phi,p_{1}}\right)$.

In case of $p_{1}\in\left(0,1\right)$, we note that $\left\Vert \mybullet\right\Vert _{L^{p_{1}}}$
is a quasi-norm with triangle constant only depending on $p_{1}$.
Together with the uniform bound $\left|J_{0}\cap J_{i}\cap J^{\left(r\right)}\right|\leq\left|J_{i}\right|\leq N_{\CalP}^{k+1}$
from equation~(\ref{eq:NecessaryConditionModerateCoveringFineInCoarseJClusterBounded})—and
recalling the identity~(\ref{eq:NecessaryModerateFineInCoarseFineLocalization})—this
yields a constant $C_{5}=C_{5}\left(p_{1},k,N_{\CalP}\right)$ satisfying
\begin{align*}
w_{i}\cdot\left\Vert \Fourier^{-1}\left(\varphi_{i}\cdot\smash{g^{\left(r\right)}}\right)\right\Vert _{L^{p_{1}}} & =w_{i}\cdot\vphantom{\sum_{j\in J_{0}\cap J_{i}\cap J^{\left(r\right)}}}\left\Vert \,\smash{\sum_{j\in J_{0}\cap J_{i}\cap J^{\left(r\right)}}}\vphantom{\sum}\Fourier^{-1}\left(d_{j}\cdot\varphi_{i}\cdot\gamma_{j}\right)\,\right\Vert _{L^{p_{1}}}\\
 & \leq C_{5}\cdot w_{i}\cdot\vphantom{\sum_{j\in J_{0}\cap J_{i}\cap J^{\left(r\right)}}^{T}}\sum_{j\in J_{0}\cap J_{i}\cap J^{\left(r\right)}}\left|d_{j}\right|\cdot\left\Vert \Fourier^{-1}\left(\varphi_{i}\cdot\gamma_{j}\right)\right\Vert _{L^{p_{1}}}.
\end{align*}
Now, observe for $i\in I_{0}$ and $j\in J_{0}\cap J_{i}\cap J^{\left(r\right)}$
that we have $\supp\varphi_{i}\subset\overline{Q_{i}}\subset\overline{P_{j}^{\left(2k+2\right)\ast}}$
as well as $\supp\gamma_{j}\subset P_{j}\subset\overline{P_{j}^{\left(2k+2\right)\ast}}$.
Indeed, $Q_{i}\subset P_{j}^{\left(2k+2\right)\ast}$ is a consequence
of $Q_{i}\cap P_{j}\neq\emptyset$ (since $j\in J_{i}$) and of Lemma~\ref{lem:SubordinatenessImpliesWeakSubordination},
together with the fact that $\CalQ_{I_{0}}$ is almost subordinate
to $\CalP$ with $k=k\left(\smash{\CalQ_{I_{0}}},\CalP\right)$. Thus,
Corollary~\ref{cor:QuasiBanachConvolutionSemiStructured} yields
a constant $C_{6}=C_{6}\left(\dimension,p_{1},k,\CalP\right)>0$ with
\begin{align*}
\left\Vert \Fourier^{-1}\left(\varphi_{i}\cdot\gamma_{j}\right)\right\Vert _{L^{p_{1}}} & \leq C_{6}\cdot\left|\det S_{j}\right|^{p_{1}^{-1}-1}\cdot\left\Vert \Fourier^{-1}\varphi_{i}\right\Vert _{L^{p_{1}}}\cdot\left\Vert \Fourier^{-1}\gamma_{j}\right\Vert _{L^{p_{1}}}\\
\left({\scriptstyle \text{eq. }\eqref{eq:NecessaryConditionModerateFineInCoarseGammaNorm}}\right) & =C_{1}^{\left(p_{1}\right)}C_{6}\cdot\left\Vert \Fourier^{-1}\varphi_{i}\right\Vert _{L^{p_{1}}}\\
\left({\scriptstyle \text{def. of an }L^{p_{1}}\text{-BAPU}}\right) & \leq C_{1}^{\left(p_{1}\right)}C_{6}C_{\CalQ,\Phi,p_{1}}\cdot\left|\det T_{i}\right|^{1-p_{1}^{-1}}.
\end{align*}

All in all, if we set $C_{7}:=C_{1}^{\left(p_{1}\right)}C_{5}C_{6}C_{\CalQ,\Phi,p_{1}}$,
we arrive at
\begin{equation}
w_{i}\cdot\left\Vert \Fourier^{-1}\left(\varphi_{i}\cdot\smash{g^{\left(r\right)}}\right)\right\Vert _{L^{p_{1}}}\leq C_{7}\cdot\left|\det T_{i}\right|^{1-p_{1}^{-1}}\cdot w_{i}\cdot\sum_{j\in J_{0}\cap J_{i}\cap J^{\left(r\right)}}\left|d_{j}\right|,\label{eq:NecessaryModerateFineInCoarseFineDecompNormEstimatePrep}
\end{equation}
where the left-hand side vanishes for $i\notin I_{0}$. Together with
the case $p_{1}\in\left[1,\infty\right]$ considered above, we conclude
that this estimate holds for all $p_{1}\in\left(0,\infty\right]$,
with $C_{7}:=C_{1}^{\left(1\right)}C_{4}$ in case of $p_{1}\in\left[1,\infty\right]$.

Now, we can estimate $\left\Vert g^{\left(r\right)}\right\Vert _{\BAPUFourierDecompSp{\CalQ}{p_{1}}{\ell_{w}^{q_{1}}}{\Phi}}$.
We first assume $q_{1}<\infty$ and make use of equation~(\ref{eq:NecessaryConditionModerateCoveringFineInCoarseJClusterBounded}),
which yields $\left|J_{0}\cap J_{i}\cap J^{\left(r\right)}\right|\leq N_{\CalP}^{k+1}$,
to derive
\[
\vphantom{\sum_{j\in J_{0}\cap J_{i}\cap J^{\left(r\right)}}}\left(\,\smash{\sum_{j\in J_{0}\cap J_{i}\cap J^{\left(r\right)}}}\vphantom{\sum}\theta_{j}\,\right)^{q_{1}}\leq\left(\left|J_{0}\cap J_{i}\cap\smash{J^{\left(r\right)}}\right|\cdot\max_{j\in J_{0}\cap J_{i}\cap J^{\left(r\right)}}\theta_{j}\right)^{q_{1}}\leq\left(N_{\CalP}^{k+1}\right)^{q_{1}}\cdot\sum_{j\in J_{0}\cap J_{i}\cap J^{\left(r\right)}}\theta_{j}^{q_{1}},
\]
for arbitrary non-negative sequences $\left(\theta_{j}\right)_{j}$.
We then use equation~(\ref{eq:NecessaryModerateFineInCoarseFineDecompNormEstimatePrep})
and set $C_{8}:=C_{7}N_{\CalP}^{k+1}$ to deduce
\begin{align}
\left\Vert \smash{g^{\left(r\right)}}\right\Vert _{\BAPUFourierDecompSp{\CalQ}{p_{1}}{\ell_{w}^{q_{1}}}{\Phi}} & =\left\Vert \left(\left\Vert \Fourier^{-1}\left(\varphi_{i}\cdot\smash{g^{\left(r\right)}}\right)\right\Vert _{L^{p_{1}}}\right)_{i\in I}\right\Vert _{\ell_{w}^{q_{1}}}\nonumber \\
\left({\scriptstyle \text{l.h.s. of }\eqref{eq:NecessaryModerateFineInCoarseFineDecompNormEstimatePrep}\text{ vanishes for }i\notin I_{0}}\right) & \leq C_{8}\cdot\left(\sum_{i\in I_{0}}\left[\,\left(\left|\det T_{i}\right|^{1-p_{1}^{-1}}\cdot w_{i}\right)^{q_{1}}\cdot\smash{\sum_{j\in J_{0}\cap J_{i}\cap J^{\left(r\right)}}}\vphantom{\sum}\left|d_{j}\right|^{q_{1}}\,\right]\right)^{1/q_{1}}\vphantom{\sum_{j\in J_{0}\cap J_{i}\cap J^{\left(r\right)}}}\nonumber \\
 & =C_{8}\cdot\left(\sum_{j\in J_{0}\cap J^{\left(r\right)}}\left[\,\left|d_{j}\right|^{q_{1}}\cdot\smash{\sum_{\substack{i\in I_{0}\\
\text{with }j\in J_{i}
}
}}\vphantom{\sum_{i}}\left(\left|\det T_{i}\right|^{1-p_{1}^{-1}}\cdot w_{i}\right)^{q_{1}}\,\right]\right)^{1/q_{1}}\vphantom{\sum_{\substack{i\in I_{0}\\
\text{with }j\in J_{i}
}
}}\nonumber \\
\left({\scriptstyle j\in J_{i}\,\Longleftrightarrow\,i\in I_{j}}\right) & \leq C_{8}\cdot\!\left\Vert \left(d_{j}\cdot\left\Vert \left(\left|\det T_{i}\right|^{1-p_{1}^{-1}}\cdot w_{i}\right)_{i\in I_{0}\cap I_{j}}\right\Vert _{\ell^{q_{1}}}\right)_{\!j\in J_{0}}\right\Vert _{\ell^{q_{1}}}\!.\label{eq:NecessaryModerateFineInCoarseAlmostDoneQ1Finite}
\end{align}
Now, let us consider the case $q_{1}=\infty$. Here, we have
\begin{align}
\left\Vert \smash{g^{\left(r\right)}}\right\Vert _{\BAPUFourierDecompSp{\CalQ}{p_{1}}{\ell_{w}^{q_{1}}}{\Phi}} & =\sup_{i\in I}\left[w_{i}\cdot\left\Vert \Fourier^{-1}\left(\varphi_{i}\cdot\smash{g^{\left(r\right)}}\right)\right\Vert _{L^{p_{1}}}\right]\nonumber \\
\left({\scriptstyle \text{l.h.s. of }\eqref{eq:NecessaryModerateFineInCoarseFineDecompNormEstimatePrep}\text{ vanishes for }i\notin I_{0}}\right) & \leq C_{7}\cdot\sup_{i\in I_{0}}\left[\,\left|\det T_{i}\right|^{1-p_{1}^{-1}}\cdot w_{i}\cdot\smash{\sum_{j\in J_{0}\cap J_{i}\cap J^{\left(r\right)}}}\vphantom{\sum}\left|d_{j}\right|\,\right]\vphantom{\sum_{j\in J_{0}\cap J_{i}\cap J^{\left(r\right)}}}\nonumber \\
\left({\scriptstyle \text{since }\left|J_{0}\cap J_{i}\cap J^{\left(r\right)}\right|\leq N_{\CalP}^{k+1}}\right) & \leq C_{7}N_{\CalP}^{k+1}\cdot\sup_{i\in I_{0}}\left[\,\left|\det T_{i}\right|^{1-p_{1}^{-1}}\cdot w_{i}\cdot\smash{\sup_{j\in J_{0}\cap J_{i}\cap J^{\left(r\right)}}}\vphantom{\sum}\left|d_{j}\right|\,\right]\vphantom{\sum_{i}^{T}}\nonumber \\
 & =C_{7}N_{\CalP}^{k+1}\cdot\sup_{j\in J_{0}\cap J^{\left(r\right)}}\left[\,\left|d_{j}\right|\cdot\smash{\sup_{\substack{i\in I_{0}\\
\text{with }j\in J_{i}
}
}}\vphantom{\sum_{i}}\left|\det T_{i}\right|^{1-p_{1}^{-1}}\cdot w_{i}\,\right]\vphantom{\sup_{\substack{i\in I_{0}\\
\text{with }j\in J_{i}
}
}}\nonumber \\
\left({\scriptstyle j\in J_{i}\,\Longleftrightarrow\,i\in I_{j}\text{ and }q_{1}=\infty}\right) & \leq C_{8}\cdot\left\Vert \left(d_{j}\cdot\left\Vert \left(\left|\det T_{i}\right|^{1-p_{1}^{-1}}\cdot w_{i}\right)_{i\in I_{0}\cap I_{j}}\right\Vert _{\ell^{q_{1}}}\right)_{j\in J_{0}}\right\Vert _{\ell^{q_{1}}}.\label{eq:NecessaryModerateFineInCoarseAlmostDoneQ1Infinite}
\end{align}
Observe that the right-hand sides of equations (\ref{eq:NecessaryModerateFineInCoarseAlmostDoneQ1Finite})
and (\ref{eq:NecessaryModerateFineInCoarseAlmostDoneQ1Infinite})
are finite, since $\left(d_{j}\right)_{j\in J_{0}}\in\ell_{0}\left(J_{0}\right)$
and because of equation~(\ref{eq:NecessaryConditionConjugateFineInCoarseSpecialAssumption}).

All in all, using the weight $u_{j}$ defined in equation~(\ref{eq:NecessaryConditionConjugateFineInCoarseSpecialAssumption}),
we see that equations (\ref{eq:NecessaryModerateFineInCoarseAlmostDoneQ1Finite})
and (\ref{eq:NecessaryModerateFineInCoarseAlmostDoneQ1Infinite})
show
\begin{equation}
\left\Vert g^{\left(r\right)}\right\Vert _{\BAPUFourierDecompSp{\CalQ}{p_{1}}{\ell_{w}^{q_{1}}}{\Phi}}\leq C_{8}\cdot\left\Vert \left(u_{j}\cdot d_{j}\right)_{j\in J_{0}}\right\Vert _{\ell^{q_{1}}}\label{eq:NecessaryModerateFineInCoarseAlmostDone}
\end{equation}
for arbitrary $q_{1}\in\left(0,\infty\right]$.

\medskip{}

Now, we finish the proof: We use estimates (\ref{eq:NecessaryModerateFineInCoarseTestFunctionEstimateCoarseCovering})
and (\ref{eq:NecessaryModerateFineInCoarseAlmostDone}), and the boundedness
of $\iota$, to deduce
\begin{align}
 & \phantom{\leq\:}\left\Vert \left(d_{m}\cdot\left|\det S_{m}\right|^{1-p_{2}^{-1}}\right)_{m\in J_{0}\cap J^{\left(r\right)}}\right\Vert _{\ell_{v}^{q_{2}}}\nonumber \\
\left({\scriptstyle \text{eq. }\eqref{eq:NecessaryModerateFineInCoarseTestFunctionEstimateCoarseCovering}}\right) & \leq\frac{C_{2}C_{3}}{C_{1}^{\left(p_{2}\right)}}\cdot\left\Vert \smash{g^{\left(r\right)}}\right\Vert _{\BAPUFourierDecompSp{\CalP}{p_{2}}{\ell_{v}^{q_{2}}}{\Psi}}\nonumber \\
 & \leq\frac{C_{2}C_{3}}{C_{1}^{\left(p_{2}\right)}}\cdot\vertiii{\iota}\cdot\left\Vert \smash{g^{\left(r\right)}}\right\Vert _{\BAPUFourierDecompSp{\CalQ}{p_{1}}{\ell_{w}^{q_{1}}}{\Phi}}\nonumber \\
\left({\scriptstyle \text{eq. }\eqref{eq:NecessaryModerateFineInCoarseAlmostDone}}\right) & \leq\frac{C_{2}C_{3}C_{8}}{C_{1}^{\left(p_{2}\right)}}\cdot\vertiii{\iota}\cdot\left\Vert \left(u_{j}\cdot d_{j}\right)_{j\in J_{0}}\right\Vert _{\ell^{q_{1}}}\qquad\forall\,d=\left(d_{j}\right)_{j\in J_{0}}\in\ell_{0}\left(J_{0}\right)\,.\label{eq:NecessaryConditionConjugateFineInCoarseBeforePartitionSummation}
\end{align}

Since $\ell_{v}^{q_{2}}\left(J_{0}\right)$ is a quasi-normed space
with triangle constant only depending on $q_{2}$, and because of
$J_{0}=\biguplus_{r=1}^{r_{0}}\left(J^{\left(r\right)}\cap J_{0}\right)$,
there is a constant $C_{9}=C_{9}\left(q_{2},r_{0}\right)=C_{9}\left(q_{2},\CalP\right)>0$
with
\[
\left\Vert \left(e_{m}\right)_{m\in J_{0}}\right\Vert _{\ell_{v}^{q_{2}}}\leq C_{9}\cdot\sum_{r=1}^{r_{0}}\left\Vert \left(e_{m}\right)_{m\in J^{\left(r\right)}\cap J_{0}}\right\Vert _{\ell_{v}^{q_{2}}}\text{ for arbitrary sequences }\left(e_{m}\right)_{m\in J_{0}}\,.
\]
Combining this with estimate~(\ref{eq:NecessaryConditionConjugateFineInCoarseBeforePartitionSummation})
implies
\[
\left\Vert \left(d_{j}\cdot\left|\det S_{j}\right|^{1-p_{2}^{-1}}\right)_{j\in J_{0}}\right\Vert _{\ell_{v}^{q_{2}}}\leq\frac{C_{2}C_{3}C_{8}C_{9}r_{0}}{C_{1}^{\left(p_{2}\right)}}\cdot\vertiii{\iota}\cdot\left\Vert \left(u_{j}\cdot d_{j}\right)_{j\in J_{0}}\right\Vert _{\ell^{q_{1}}}\qquad\forall\,\left(d_{j}\right)_{j\in J_{0}}\in\ell_{0}\left(J_{0}\right)\,,
\]
which simply means that the embedding 
\[
\eta:\ell_{0}\left(J_{0}\right)\cap\ell_{u}^{q_{1}}\left(J_{0}\right)\hookrightarrow\ell_{\mu}^{q_{2}}\left(J_{0}\right)\qquad\text{with}\qquad\mu_{j}:=v_{j}\cdot\left|\det S_{j}\right|^{1-p_{2}^{-1}}\text{ for }j\in J
\]
is bounded, with $\vertiii{\eta}\leq C_{10}\cdot\vertiii{\iota}$
for $C_{10}:=\frac{C_{2}C_{3}C_{8}C_{9}r_{0}}{C_{1}^{\left(p_{2}\right)}}$.
But since $\ell_{\mu}^{q_{2}}\left(J_{0}\right)$ satisfies the Fatou
property, Lemmas \ref{lem:EmbeddingBetweenWeightedSequenceSpaces}
and \ref{lem:FinitelySupportedSequencesSufficeUnderFatouProperty}
show that this implies
\[
\left\Vert \left(v_{j}\cdot\left|\det S_{j}\right|^{1-p_{2}^{-1}}\big/\,u_{j}\right)_{j\in J_{0}}\right\Vert _{\ell^{q_{2}\cdot\left(q_{1}/q_{2}\right)'}}=\left\Vert \left(\mu_{j}/u_{j}\right)_{j\in J_{0}}\right\Vert _{\ell^{q_{2}\cdot\left(q_{1}/q_{2}\right)'}}=\vertiii{\eta}\leq C_{10}\cdot\vertiii{\iota}<\infty.
\]
To see that $J_{0}\subset J$ is countable (as required by Lemma~\ref{lem:FinitelySupportedSequencesSufficeUnderFatouProperty}),
one can argue as at the end of the proof of Theorem~\ref{thm:BurnerNecessaryConditionFineInCoarse}.
\end{proof}
As before, we now specialize the above lemma to the case where $\CalQ$
(or more precisely $\CalQ_{I_{0}}$) is relatively $\CalP$-moderate.
We remark that the following theorem is a slightly generalized version
of \cite[Theorem 5.3.14]{VoigtlaenderPhDThesis} from my PhD thesis.
\begin{thm}
\label{thm:NecessaryConditionForModerateCoveringFineInCoarse}Let
$\emptyset\neq\CalO,\CalO'\subset\R^{\dimension}$ be open, let $\CalQ=\left(Q_{i}\right)_{i\in I}=\left(T_{i}Q_{i}'+b_{i}\right)_{i\in I}$
be a \emph{tight} semi-structured $L^{p_{1}}$-decomposition covering
of $\CalO$ and let $\CalP=\left(P_{j}\right)_{j\in J}=\left(\smash{S_{j}P_{j}'+c_{j}}\right)_{j\in J}$
be a \emph{tight} semi-structured $L^{p_{2}}$-decomposition covering
of $\CalO'$, for certain $p_{1},p_{2}\in\left(0,\infty\right]$.
Finally, let $q_{1},q_{2}\in\left(0,\infty\right]$ and let $w=\left(w_{i}\right)_{i\in I}$
and $v=\left(v_{j}\right)_{j\in J}$ be $\CalQ$-moderate and $\CalP$-moderate,
respectively. 

Choose an arbitrary subset $J_{0}\subset J$, set 
\[
I_{0}:=\left\{ i\in I\with J_{i}\cap J_{0}\neq\emptyset\right\} 
\]
and assume that the following hold:

\begin{enumerate}
\item \label{enu:NecessaryConditionModerateCoveringFineInCoarseAlmostTightness}There
is $\varepsilon>0$ and for each $j\in J_{0}$ some $\xi_{j}\in\R^{\dimension}$
with $B_{\varepsilon}\left(\xi_{j}\right)\subset P_{j}'$ and $S_{j}\left[B_{\varepsilon}\left(\xi_{j}\right)\right]+c_{j}\subset\CalO$.
\item $\CalQ_{I_{0}}:=\left(Q_{i}\right)_{i\in I_{0}}$ is almost subordinate
to $\CalP$.
\item \label{enu:NecessaryConditionModerateCoveringFineInCoarseModerateness}$\CalQ_{I_{0}}$
is relatively $\CalP$-moderate.
\item \label{enu:NecessaryConditionModerateCoveringFineInCoarseWeightModerate}The
weight $w|_{I_{0}}$ is relatively $\CalP$-moderate.
\end{enumerate}
Furthermore, set
\[
K:=\bigcup_{i\in I_{0}}Q_{i}\subset\CalO\cap\CalO'
\]
and—with $\CalD_{K}:=\left\{ f\in\TestFunctionSpace{\R^{\dimension}}\with\supp f\subset K\right\} $—assume
that the identity map
\[
\iota:\left(\CalD_{K},\left\Vert \mybullet\right\Vert _{\FourierDecompSp{\CalQ}{p_{1}}{\ell_{w}^{q_{1}}}}\right)\rightarrow\FourierDecompSp{\CalP}{p_{2}}{\ell_{v}^{q_{2}}},f\mapsto f
\]
is bounded.

For each $j\in J_{0}$, let $i_{j}\in I_{j}$ be arbitrary\footnote{We have $S_{j}\left[B_{\varepsilon}\left(\xi_{j}\right)\right]+c_{j}\subset\CalO\cap P_{j}$,
so that there is some $i\in I$ with $Q_{i}\cap P_{j}\neq\emptyset$.
In particular, $I_{j}\neq\emptyset$, so that $i_{j}$ can be chosen
as desired.}. Then, we have
\begin{equation}
\left\Vert \!\left(\!\frac{v_{j}}{w_{i_{j}}}\cdot\left|\det\smash{T_{i_{j}}}\right|^{\frac{1}{p_{1}}-\left(\smash{\frac{1}{\LowerExpo{p_{2}}}}-\frac{1}{q_{1}}\right)_{+}-\frac{1}{p_{2}}}\cdot\left|\det S_{j}\right|^{\left(\smash{\frac{1}{\LowerExpo{p_{2}}}}-\frac{1}{q_{1}}\right)_{+}}\right)_{\!\!\!j\in J_{0}}\right\Vert _{\ell^{q_{2}\cdot\left(q_{1}/q_{2}\right)'}}\leq C\cdot\vertiii{\iota}\label{eq:NecessaryModerateFineInCoarseTargetEstimate}
\end{equation}
for some constant $C>0$ depending only on
\[
\dimension,p_{1},p_{2},q_{1},q_{2},k\left(\smash{\CalQ_{I_{0}}},\CalP\right),C_{{\rm mod}}\left(\smash{\CalQ_{I_{0}}},\CalP\right),\CalQ,\CalP,\varepsilon,\varepsilon_{\CalQ},\varepsilon_{\CalP},C_{w,\CalQ},C_{v,\CalP},C_{w|_{I_{0}},\CalQ,\CalP},C_{\CalQ,\Phi,p_{1}},C_{\CalP,\Psi,p_{2}}.
\]
Here, the $L^{p_{1}}$/$L^{p_{2}}$-BAPUs $\Phi=\left(\varphi_{i}\right)_{i\in I}$
and $\Psi=\left(\psi_{j}\right)_{j\in J}$ are those which are used
to calculate $\vertiii{\iota}$.
\end{thm}

\begin{proof}
For brevity, let $k:=k\left(\smash{\CalQ_{I_{0}}},\CalP\right)$.
Furthermore, set $L:=C_{w|_{I_{0}},\CalQ,\CalP}$, so that we have
$w_{i}\leq L\cdot w_{\ell}$ for all $i,\ell\in I_{0}\cap I_{j}$
and arbitrary $j\in J$. Moreover, set 
\[
P_{j}^{\left(0\right)}:=S_{j}\left[B_{\varepsilon}\left(\xi_{j}\right)\right]+c_{j}\subset P_{j}\cap\CalO\qquad\forall\,j\in J_{0},
\]
with $\varepsilon,\xi_{j}$ as in part~(\ref{enu:NecessaryConditionModerateCoveringFineInCoarseAlmostTightness})
of the prerequisites of the current theorem.

As in the proof of Lemma~\ref{lem:NecessaryConditionConjugateFineInCoarseWithoutModerateness}
(see equation~(\ref{eq:NecessaryConditionModerateFineInCoarseIntersectionSetContainedInI0})),
we note 
\[
\emptyset\neq I_{j}\subset I_{0}\qquad\forall\,j\in J_{0}.
\]
Indeed, we have $\emptyset\neq P_{j}^{\left(0\right)}\subset P_{j}\cap\CalO$
and thus $Q_{i}\cap P_{j}\neq\emptyset$ for some $i\in I$, i.e.\@
$I_{j}\neq\emptyset$. Furthermore, every $i\in I_{j}$ satisfies
$j\in J_{i}\cap J_{0}\neq\emptyset$ and thus $i\in I_{0}$.

\medskip{}

Our first goal is to invoke Lemma~\ref{lem:IntersectionCountForModerateCoverings}
to get $\left|I_{j}\right|=\left|I_{0}\cap I_{j}\right|\asymp\left|\det S_{j}\right|/\left|\det T_{i_{j}}\right|$
for all $j\in J_{0}$. Thus, we first need to verify the assumptions
of Lemma~\ref{lem:IntersectionCountForModerateCoverings}. Note that
$\CalQ_{I_{0}}$ is relatively $\CalP$-moderate and thus in particular
relatively $\CalP_{J_{0}}$-moderate. Hence, it remains to verify
part~(\ref{enu:IntersectionCountCoarseCoveringCoveredByFine}) of
the assumptions of Lemma~\ref{lem:IntersectionCountForModerateCoverings},
i.e., we need to show
\begin{equation}
\lambda\left(P_{j}\right)\leq C_{0}\cdot\lambda\left(\,\smash{\bigcup_{i\in I_{0}\cap I_{j}}}\vphantom{\bigcup}Q_{i}^{r\ast}\,\right)\vphantom{\bigcup_{i\in I_{0}\cap I_{j}}}\label{eq:IntersectionCountSpecialPrerequisite}
\end{equation}
for all $j\in J_{0}$ and suitable $r\in\N_{0}$ and $C_{0}>0$. Note
$\lambda\left(\smash{P_{j}^{\left(0\right)}}\right)=C_{1}\cdot\left|\det S_{j}\right|$
for all $j\in J_{0}$ and a suitable constant $C_{1}=C_{1}\left(\dimension,\varepsilon\right)$.
But Corollary~\ref{cor:SemiStructuredDifferenceSetsMeasureEstimate}
shows $\lambda\left(P_{j}\right)\leq C_{2}\cdot\left|\det S_{j}\right|=\frac{C_{2}}{C_{1}}\cdot\lambda\left(\smash{P_{j}^{\left(0\right)}}\right)$
for a constant $C_{2}=C_{2}\left(\dimension,\CalP\right)$. Furthermore,
we have $P_{j}^{\left(0\right)}\subset P_{j}\cap\CalO$ and hence
\[
P_{j}^{\left(0\right)}\subset\bigcup_{\substack{i\in I\\
\text{with }Q_{i}\cap P_{j}^{\left(0\right)}\neq\emptyset
}
}\!\!\!Q_{i}\subset\bigcup_{i\in I_{j}}Q_{i}\quad\overset{I_{j}\subset I_{0}}{=}\quad\bigcup_{i\in I_{0}\cap I_{j}}Q_{i}^{0\ast}\qquad\forall\,j\in J_{0}\,.
\]
Thus equation~(\ref{eq:IntersectionCountSpecialPrerequisite}) is
fulfilled for $r:=0$ and $C_{0}:=C_{2}/C_{1}$.

All in all, we have shown that Lemma~\ref{lem:IntersectionCountForModerateCoverings}
is applicable, so that there are constants 
\[
C_{3}=C_{3}\left(\dimension,k\left(\smash{\CalQ_{I_{0}}},\CalP\right),C_{{\rm mod}}\left(\smash{\CalQ_{I_{0}}},\CalP\right),\CalQ,\varepsilon_{\CalQ},\CalP\right)\quad\text{and}\quad C_{4}=C_{4}\left(\dimension,C_{{\rm mod}}\left(\smash{\CalQ_{I_{0}}},\CalP\right),\varepsilon,\CalQ,\CalP,\varepsilon_{\CalP}\right)
\]
which satisfy
\begin{equation}
\left|I_{0}\cap I_{j}\right|\leq C_{3}\cdot\left|\det S_{j}\right|\cdot\left|\det\smash{T_{i_{j}}}\right|^{-1}\qquad\forall\,j\in J_{0},\label{eq:NecessaryModerateCoveringFineInCoarseIntersectionEstimateAbove}
\end{equation}
and
\begin{equation}
\left|\det\smash{T_{i_{j}}}\right|^{-1}\cdot\left|\det S_{j}\right|\leq C_{4}\cdot\left|I_{0}\cap I_{j}\right|\qquad\forall\,j\in J_{0}\,.\label{eq:NecessaryModerateCoveringFIneInCoarseIntersectionSetEstimate}
\end{equation}
Here, we also used that $I_{0}\cap I_{j}=I_{j}$ is nonempty for all
$j\in J_{0}$, as seen above.

\medskip{}

In particular, equation~(\ref{eq:NecessaryModerateCoveringFineInCoarseIntersectionEstimateAbove})
shows that $I_{0}\cap I_{j}$ is finite for every $j\in J_{0}$, so
that equation~(\ref{eq:NecessaryConditionConjugateFineInCoarseSpecialAssumption})
from the prerequisites of Lemma~\ref{lem:NecessaryConditionConjugateFineInCoarseWithoutModerateness}
is satisfied. In fact, using the notation from that lemma, we get
\begin{align*}
u_{j}=\left\Vert \left(w_{i}\cdot\left|\det T_{i}\right|^{1-p_{1}^{-1}}\right)_{i\in I_{0}\cap I_{j}}\right\Vert _{\ell^{q_{1}}} & \leq LC_{5}\cdot w_{i_{j}}\cdot\left|\det T_{i_{j}}\right|^{1-p_{1}^{-1}}\cdot\left|I_{0}\cap I_{j}\right|^{q_{1}^{-1}}\\
 & \leq LC_{3}^{q_{1}^{-1}}\!C_{5}\cdot w_{i_{j}}\!\cdot\left|\det S_{j}\right|^{q_{1}^{-1}}\!\cdot\left|\det T_{i_{j}}\right|^{1-p_{1}^{-1}-q_{1}^{-1}}\!<\!\infty\quad\forall\,j\in J_{0}\,,
\end{align*}
for some constant $C_{5}=C_{5}\left(p_{1},C_{{\rm mod}}\left(\smash{\CalQ_{I_{0}}},\CalP\right)\right)$.
Note that this also holds in case of $1-p_{1}^{-1}<0$, since we have
$\left|\det T_{i}\right|\asymp\left|\det T_{i_{j}}\right|$ for all
$i\in I_{0}\cap I_{j}$. For brevity, set $C_{6}:=LC_{3}^{1/q_{1}}C_{5}$.

It is now easy to see that all assumptions of Lemma~\ref{lem:NecessaryConditionConjugateFineInCoarseWithoutModerateness}
are satisfied, so that we get a constant $C_{7}>0$, which only depends
on quantities mentioned in the statement of the theorem, and which
satisfies
\begin{align}
\left\Vert \left(\frac{v_{j}}{w_{i_{j}}}\!\cdot\!\left|\det S_{j}\right|^{1-q_{1}^{-1}-p_{2}^{-1}}\!\cdot\!\left|\det T_{i_{j}}\right|^{q_{1}^{-1}+p_{1}^{-1}-1}\right)_{\!j\in J_{0}}\right\Vert _{\ell^{q}} & \leq C_{6}\cdot\left\Vert \left(\left|\det S_{j}\right|^{1-p_{2}^{-1}}\cdot\frac{v_{j}}{u_{j}}\right)_{\!j\in J_{0}}\right\Vert _{\ell^{q}}\nonumber \\
 & \leq C_{6}C_{7}\cdot\vertiii{\iota},\label{eq:NecessaryModerateFineInCoarseDone}
\end{align}
where we defined $q:=q_{2}\cdot\left(q_{1}/q_{2}\right)'$ for brevity.

To see that this indeed implies the claim, we distinguish two cases:

\textbf{Case 1}: We have $p_{2}\in\left[2,\infty\right]$ and $q_{1}\geq p_{2}'$.
In this case, we get $\LowerExpo{p_{2}}=\min\left\{ p_{2},p_{2}'\right\} =p_{2}'$
and $\frac{1}{\LowerExpo{p_{2}}}-\frac{1}{q_{1}}=\frac{1}{p_{2}'}-\frac{1}{q_{1}}\geq0$.
Thus, $\left(\smash{\frac{1}{\LowerExpo{p_{2}}}}-\frac{1}{q_{1}}\right)_{+}=\frac{1}{p_{2}'}-\frac{1}{q_{1}}=1-\frac{1}{p_{2}}-\frac{1}{q_{1}}$
and hence
\begin{align*}
\left|\det\smash{T_{i_{j}}}\right|^{\frac{1}{p_{1}}-\left(\smash{\frac{1}{\LowerExpo{p_{2}}}}-\frac{1}{q_{1}}\right)_{+}-\frac{1}{p_{2}}}\cdot\left|\det S_{j}\right|^{\left(\smash{\frac{1}{\LowerExpo{p_{2}}}}-\frac{1}{q_{1}}\right)_{+}} & =\left|\det S_{j}\right|^{1-\frac{1}{p_{2}}-\frac{1}{q_{1}}}\cdot\left|\det\smash{T_{i_{j}}}\right|^{\frac{1}{p_{1}}-\left(1-\frac{1}{p_{2}}-\frac{1}{q_{1}}\right)-\frac{1}{p_{2}}}\\
 & =\left|\det S_{j}\right|^{1-\frac{1}{p_{2}}-\frac{1}{q_{1}}}\cdot\left|\det\smash{T_{i_{j}}}\right|^{\frac{1}{p_{1}}+\frac{1}{q_{1}}-1}.
\end{align*}
Plugging this into the target inequality (\ref{eq:NecessaryModerateFineInCoarseTargetEstimate}),
we see that equation~(\ref{eq:NecessaryModerateFineInCoarseDone})
implies the claim.

\medskip{}

\textbf{Case 2}: We have $p_{2}\in\left(0,2\right]$ or $q_{1}<p_{2}'$.
For brevity, define $s:=\left(\smash{\frac{1}{\LowerExpo{p_{2}}}}-\frac{1}{q_{1}}\right)_{+}$
and $t:=p_{2}\cdot\left(q_{1}/p_{2}\right)'$. Recall from equation~(\ref{eq:InverseOfSpecialExponent})
that
\[
\frac{1}{t}=\left(\frac{1}{p_{2}}-\frac{1}{q_{1}}\right)_{+}.
\]
Our first goal is to show $s=1/t$. In case of $p_{2}\in\left(0,2\right]$,
this is clear, since we have $\LowerExpo{p_{2}}=p_{2}$. Thus, let
us assume $p_{2}\in\left(2,\infty\right]$. In view of our case distinction,
this entails $q_{1}<p_{2}'=\LowerExpo{p_{2}}$ and hence $s=0$. But
since $p_{2}\in\left(2,\infty\right]$, we also have $p_{2}'<p_{2}$
and hence $q_{1}<p_{2}$, which entails $s=0=1/t$, as desired.

Now, it is easy to see that the assumptions of the current theorem
include those of Theorem~\ref{thm:BurnerNecessaryConditionFineInCoarse}.
Hence (using the estimate for $\vertiii{\Gamma_{\CalQ}}_{\ell_{w}^{q_{1}}\to\ell_{w}^{q_{1}}}$
and $\vertiii{\Gamma_{\CalP}}_{\ell_{v}^{q_{2}}\to\ell_{v}^{q_{2}}}$
from Lemma~\ref{lem:ModeratelyWeightedSpacesAreRegular}), we see
that there is a constant 
\[
C_{8}=C_{8}\left(\dimension,p_{1},p_{2},q_{1},q_{2},k\left(\smash{\CalQ_{I_{0}}},\CalP\right),\CalQ,\varepsilon_{\CalQ},\CalP,C_{w,\CalQ},C_{v,\CalP},C_{\CalQ,\Phi,p_{1}},C_{\CalP,\Psi,p_{2}}\right)
\]
which satisfies $\vertiii{\eta}\leq C_{8}\cdot\vertiii{\iota}$, with
\[
\eta:\ell_{\left(w_{i}\cdot\left|\det T_{i}\right|^{p_{2}^{-1}-p_{1}^{-1}}\right)_{i}}^{q_{1}}\left(I_{0}\right)\hookrightarrow\ell_{v}^{q_{2}}\left(\left[\ell^{p_{2}}\left(I_{0}\cap I_{j}\right)\right]_{j\in J}\right).
\]
In view of Corollary~\ref{cor:EmbeddingFineInCoarseSimplification}
(with $u\equiv1$ and $r=p_{2}$, as well as $J_{0}=J$), there is
hence a suitable constant $C_{9}=C_{9}\left(q_{1},q_{2},p_{2},k\left(\smash{\CalQ_{I_{0}}},\CalP\right),\CalP,C_{v,\CalP}\right)$
satisfying
\[
\left\Vert \left(v_{j}\cdot\left\Vert \left(\left|\det T_{i}\right|^{p_{1}^{-1}-p_{2}^{-1}}\big/w_{i}\right)_{i\in I_{0}\cap I_{j}}\right\Vert _{\ell^{p_{2}\cdot\left(q_{1}/p_{2}\right)'}}\right)_{j\in J}\right\Vert _{\ell^{q_{2}\cdot\left(q_{1}/q_{2}\right)'}}\leq C_{9}\cdot\vertiii{\eta}\leq C_{8}C_{9}\cdot\vertiii{\iota}.
\]
Hence,
\begin{align*}
 & \left\Vert \left(\frac{v_{j}}{w_{i_{j}}}\cdot\left|\det T_{i_{j}}\right|^{p_{1}^{-1}-p_{2}^{-1}-s}\cdot\left|\det S_{j}\right|^{s}\right)_{j\in J_{0}}\right\Vert _{\ell^{q_{2}\cdot\left(q_{1}/q_{2}\right)'}}\\
\left({\scriptstyle \text{since }s=1/t}\right) & =\left\Vert \left(\frac{v_{j}}{w_{i_{j}}}\cdot\left|\det T_{i_{j}}\right|^{p_{1}^{-1}-p_{2}^{-1}}\cdot\left|\det S_{j}\right|^{t^{-1}}\big/\left|\det T_{i_{j}}\right|^{t^{-1}}\right)_{j\in J_{0}}\right\Vert _{\ell^{q_{2}\cdot\left(q_{1}/q_{2}\right)'}}\\
\left({\scriptstyle \text{eq. }\eqref{eq:NecessaryModerateCoveringFIneInCoarseIntersectionSetEstimate}}\right) & \leq C_{4}^{1/t}\cdot\left\Vert \left(\frac{v_{j}}{w_{i_{j}}}\cdot\left|\det T_{i_{j}}\right|^{p_{1}^{-1}-p_{2}^{-1}}\cdot\left\Vert \left(1\right)_{i\in I_{0}\cap I_{j}}\right\Vert _{\ell^{t}}\right)_{j\in J_{0}}\right\Vert _{\ell^{q_{2}\cdot\left(q_{1}/q_{2}\right)'}}\\
 & \leq LC_{4}^{1/t}C_{10}\cdot\left\Vert \left(v_{j}\cdot\left\Vert \left(\left|\det T_{i}\right|^{p_{1}^{-1}-p_{2}^{-1}}\big/w_{i}\right)_{i\in I_{0}\cap I_{j}}\right\Vert _{\ell^{p_{2}\cdot\left(q_{1}/p_{2}\right)'}}\right)_{j\in J}\right\Vert _{\ell^{q_{2}\cdot\left(q_{1}/q_{2}\right)'}}\\
 & \leq LC_{4}^{1/t}C_{8}C_{9}C_{10}\cdot\vertiii{\iota}
\end{align*}
for some constant $C_{10}=C_{10}\left(C_{{\rm mod}}\left(\smash{\CalQ_{I_{0}}},\CalP\right),p_{1},p_{2}\right)$.
This is precisely the desired estimate.
\end{proof}

\section{An overview of the derived embedding results}

\label{sec:SummaryOfEmbeddingResults}In this section, we summarize
those embedding results from the present paper which are most convenient
to use. We will always be interested in sufficient or necessary conditions
for an embedding of the form
\[
\FourierDecompSp{\CalQ}{p_{1}}{\ell_{w}^{q_{1}}}\hookrightarrow\FourierDecompSp{\CalP}{p_{2}}{\ell_{v}^{q_{2}}},
\]
under varying assumptions on the coverings $\CalQ,\CalP$ of the open
sets $\CalO,\CalO'\subset\R^{\dimension}$.

In the first subsection, we state three theorems which provide sufficient
as well as necessary criteria for the existence of the above embedding:
Theorem~\ref{thm:SummaryCoarseIntoFine} is applicable when $\CalP$
is almost subordinate to $\CalQ$, while Theorem~\ref{thm:SummaryFineIntoCoarse}
applies in the ``converse'' case where $\CalQ$ is almost subordinate
to $\CalP$. Finally, Theorem~\ref{thm:SummaryMixedSubordinateness}
applies if we can write $\CalO\cap\CalO'=A\cup B$ where $\CalQ$
is almost subordinate to $\CalP$ ``near $A$'', while $\CalP$
is almost subordinate to $\CalQ$ ``near $B$''.

In the second subsection, we accompany these results with a ``user's
guide'' for their application. This user's guide describes a sequence
of steps that can be followed to determine the existence of an embedding
between two decomposition spaces.

While we strove for maximal generality in the preceding sections,
our aim in this section is ease of applicability, even if this makes
our results slightly less general. Thus, those readers who are interested
in embeddings between decomposition spaces which are highly ``incompatible''—for
instance if there is no subordinateness or no relative moderateness
or we have neither $\CalO\subset\CalO'$, nor $\CalO'\subset\CalO$—or
for which the ``global components'' are not just weighted $\ell^{q}$-spaces,
are encouraged to browse the results from the preceding sections and
not only this one. But for the most common cases, the results in this
section are easily applicable and yield optimal results.

\subsection{A collection of readily applicable embedding results}

To simplify notation, we will employ the following general assumptions
in this section:
\begin{assumption}
\label{assu:GeneralSummaryAssumptions}Let $\emptyset\neq\CalO,\CalO'\subset\R^{\dimension}$
be open and let $\CalQ=\left(T_{i}Q_{i}'+b_{i}\right)_{i\in I}$ and
$\CalP=\left(S_{j}P_{j}'+c_{j}\right)_{j\in J}$ be two open, tight,
semi-structured coverings of $\CalO$ and $\CalO'$, respectively.
Let $w=\left(w_{i}\right)_{i\in I}$ and $v=\left(v_{j}\right)_{j\in J}$
be $\CalQ$-moderate and $\CalP$-moderate, respectively and let $p_{1},p_{2},q_{1},q_{2}\in\left(0,\infty\right]$.

Furthermore, assume that $\Phi=\left(\varphi_{i}\right)_{i\in I}$
and $\Psi=\left(\psi_{j}\right)_{j\in J}$ are $L^{p}$-BAPUs for
$\CalQ$ or $\CalP$, respectively, simultaneously for all $p\in\left(0,\infty\right]$.

Finally, assume that the (quasi)-norms $\left\Vert \mybullet\right\Vert _{\FourierDecompSp{\CalQ}{p_{1}}{\ell_{w}^{q_{1}}}}$
and $\left\Vert \mybullet\right\Vert _{\FourierDecompSp{\CalP}{p_{2}}{\ell_{v}^{q_{2}}}}$
are calculated using the BAPUs $\Phi$ and $\Psi$, respectively.
\end{assumption}

\begin{rem*}
Note that the preceding assumptions—excluding the moderateness of
$w,v$—are always fulfilled if $\CalQ$ and $\CalP$ are \emph{almost
structured} coverings (see Definition~\ref{defn:DifferentTypesOfCoverings}
and the ensuing remark), for a proper choice of $\Phi,\Psi$ (see
Theorem~\ref{thm:AlmostStructuredAdmissibleAdmitsBAPU}).
\end{rem*}
We begin with the case in which $\CalP$ is almost subordinate to
$\CalQ$. We remark that the following result is a generalized version
of \cite[Theorem 5.4.1]{VoigtlaenderPhDThesis} from my PhD thesis.
\begin{thm}
\label{thm:SummaryCoarseIntoFine}In addition to our standing assumptions,
assume that $\CalP$ is almost subordinate to $\CalQ$. Note that
this entails $\CalO'\subset\CalO$.

For $i\in I$, let
\[
J_{i}:=\left\{ j\in J\with P_{j}\cap Q_{i}\neq\emptyset\right\} 
\]
and for $r\in\left(0,\infty\right]$ and $\ell\in\left\{ 1,2\right\} $,
define
\[
K_{r,\ell}:=\left\Vert \left(w_{i}^{-1}\cdot\left\Vert \left(v_{j}/\smash{u_{i,j}^{\left(\ell\right)}}\right)_{j\in J_{i}}\right\Vert _{\ell^{q_{2}\cdot\left(r/q_{2}\right)'}}\right)_{i\in I}\right\Vert _{\ell^{q_{2}\cdot\left(q_{1}/q_{2}\right)'}}\in\left[0,\infty\right],
\]
with
\[
u_{i,j}^{\left(1\right)}:=\left|\det S_{j}\right|^{p_{2}^{-1}-p_{1}^{-1}}\quad\!\text{and}\quad\!u_{i,j}^{\left(2\right)}:=\begin{cases}
\left|\det S_{j}\right|^{p_{2}^{-1}-1}\cdot\left|\det T_{i}\right|^{1-p_{1}^{-1}}, & \text{if }p_{1}<1,\\
\vphantom{\rule{0.1cm}{0.55cm}}\left|\det S_{j}\right|^{p_{2}^{-1}-p_{1}^{-1}}, & \text{if }p_{1}\geq1
\end{cases}\:\:\quad\:\:\text{for }i\in I\text{ and }j\in J.
\]

\noindent Then, the following hold:

\begin{enumerate}[leftmargin=0.85cm]
\item \label{enu:SummaryCoarseInFineSufficientNonModerate}If we have $p_{1}\leq p_{2}$
and if $K_{\UpperExpo{p_{1}},2}<\infty$, then
\[
\iota:\FourierDecompSp{\CalQ}{p_{1}}{\ell_{w}^{q_{1}}}\to\FourierDecompSp{\CalP}{p_{2}}{\ell_{v}^{q_{2}}},f\mapsto f|_{\TestFunctionSpace{\CalO'}}
\]
is well-defined and bounded, with $\vertiii{\iota}\leq C_{1}\cdot K_{\UpperExpo{p_{1}},2}$
for some constant
\[
C_{1}=C_{1}\left(\dimension,p_{1},p_{2},q_{1},q_{2},\CalQ,\CalP,k\left(\CalP,\CalQ\right),C_{\CalQ,\Phi,p_{1}},C_{\CalP,\Psi,p_{1}},C_{w,\CalQ},C_{v,\CalP}\right).
\]
\item \label{enu:SummaryCoarseInFineNecessaryNonModerate}Conversely, if
the identity map
\[
\theta:\left(\TestFunctionSpace{\CalO'},\left\Vert \mybullet\right\Vert _{\FourierDecompSp{\CalQ}{p_{1}}{\ell_{w}^{q_{1}}}}\right)\to\FourierDecompSp{\CalP}{p_{2}}{\ell_{v}^{q_{2}}},f\mapsto f
\]
is bounded, then we have $p_{1}\leq p_{2}$ and $K_{p_{1},1}\leq C_{2}\cdot\vertiii{\theta}<\infty$
for some constant
\[
C_{2}=C_{2}\left(\dimension,p_{1},p_{2},q_{1},q_{2},\CalQ,\CalP,\varepsilon_{\CalP},k\left(\CalP,\CalQ\right),C_{\CalQ,\Phi,p_{1}},C_{w,\CalQ},C_{v,\CalP}\right).
\]
\item \label{enu:SummaryCoarseInFineNecessaryKhinchin}Under the assumptions
of the preceding point, if $p_{1}=p_{2}$, then $K_{2,1}\leq C_{3}\cdot\vertiii{\theta}<\infty$
for some constant
\[
C_{3}=C_{3}\left(\dimension,p_{1},q_{1},q_{2},\CalQ,\CalP,k\left(\CalP,\CalQ\right),C_{\CalQ,\Phi,p_{1}},C_{\CalP,\Psi,p_{2}},C_{w,\CalQ},C_{v,\CalP}\right).
\]
\item \label{enu:SummaryCoarseInFineModerate}Finally, if we have $\CalO=\CalO'$
and if $\CalP$ and $v$ are relatively $\CalQ$-moderate, then we
have the following equivalence (with $\iota,\theta$ as above): Setting
$s:=\left(\frac{1}{q_{2}}-\smash{\frac{1}{\SignedUpperExpo{p_{1}}}}\right)_{+}$
and choosing for each $i\in I$ some $j_{i}\in J$ with $Q_{i}\cap P_{j_{i}}\neq\emptyset$,
we have:
\begin{align*}
 & \iota\text{ well-defined and bounded}\\
\Longleftrightarrow\: & \theta\text{ well-defined and bounded}\\
\Longleftrightarrow\: & p_{1}\leq p_{2}\text{ and }K:=\left\Vert \left(\frac{v_{j_{i}}}{w_{i}}\cdot\left|\det T_{i}\right|^{s}\cdot\left|\det\smash{S_{j_{i}}}\right|^{p_{1}^{-1}-p_{2}^{-1}-s}\right)_{i\in I}\right\Vert _{\ell^{q_{2}\cdot\left(q_{1}/q_{2}\right)'}}<\infty.
\end{align*}
Furthermore, given $p_{1}\leq p_{2}$, we have
\[
\vertiii{\iota}\asymp\vertiii{\theta}\asymp K,
\]
where the implied constants only depend on
\[
\dimension,p_{1},p_{2},q_{1},q_{2},\CalQ,\varepsilon_{\CalQ},\CalP,\varepsilon_{\CalP},k\left(\CalP,\CalQ\right),C_{{\rm mod}}\left(\CalP,\CalQ\right),C_{w,\CalQ},C_{v,\CalP},C_{v,\CalP,\CalQ},C_{\CalQ,\Phi,p_{1}},C_{\CalP,\Psi,p_{1}}.\qedhere
\]
\end{enumerate}
\end{thm}

\begin{proof}
In the whole proof, let us set $Y:=\ell_{w}^{q_{1}}\left(I\right)$
and $Z:=\ell_{v}^{q_{2}}\left(J\right)$, as well as $J_{0}:=J$.
We now prove each statement individually.

Ad (1): Remark~\ref{rem:SufficientCoarseIntoFineSimplification}
shows that the embedding 
\[
\eta:Y\Bigl(\Bigl[\ell_{u^{\left(2\right)}}^{\UpperExpo{p_{1}}}\left(J_{0}\cap J_{i}\right)\Bigr]_{i\in I}\Bigr)\hookrightarrow Z|_{J_{0}}
\]
satisfies $\vertiii{\eta}\asymp K_{\UpperExpo{p_{1}},2}$, where the
implied constant only depends on $p_{1},q_{1},q_{2},C_{w,\CalQ},\CalQ,k\left(\CalP,\CalQ\right)$.

Hence, the boundedness of $\iota$, with 
\[
\vertiii{\iota}\leq C\cdot\vertiii{\eta}\leq C'\cdot K_{\UpperExpo{p_{1}},2}<\infty
\]
for suitable constants $C,C'$ as in the statement of the theorem,
follows directly from Corollary~\ref{cor:EmbeddingCoarseIntoFine}.
Here, we also used Lemma~\ref{lem:ModeratelyWeightedSpacesAreRegular}
to estimate $\vertiii{\Gamma_{\CalQ}}_{\ell_{w}^{q_{1}}\to\ell_{w}^{q_{1}}}$
and $\vertiii{\Gamma_{\CalP}}_{\ell_{v}^{q_{2}}\to\ell_{v}^{q_{2}}}$
by quantities involving only $q_{1},q_{2},\CalQ,\CalP,C_{w,\CalQ},C_{v,\CalP}$.

\medskip{}

Ad (2): First, note that $\TestFunctionSpace{\CalO'}\subset\TestFunctionSpace{\CalO}$
and hence $\TestFunctionSpace{\CalO}\subset\FourierDecompSp{\CalQ}{p_{1}}{\ell_{w}^{q_{1}}}\cap\FourierDecompSp{\CalP}{p_{2}}{\ell_{v}^{q_{2}}}$,
so that $\theta$ is always well-defined, but not necessarily bounded.
Of course, in the present case, $\theta$ is bounded by assumption.

Now, it is easy to see that Lemma~\ref{lem:SimpleNecessaryCondition}
is applicable with $K:=\CalO'\subset\CalO\cap\CalO'$ . Picking any
$j\in J$ and any $\xi\in P_{j}\subset\CalO'\subset\CalO$, there
is some $i\in I$ with $\xi\in Q_{i}$. Hence, $\xi\in K^{\circ}\cap Q_{i}^{\circ}\cap P_{j}^{\circ}$,
since by our standing assumptions, $\CalQ$ and $\CalP$ are open
coverings. Finally, $\delta_{i}\in Y=\ell_{w}^{q_{1}}\left(I\right)$
for arbitrary $i\in I$, so that Lemma~\ref{lem:SimpleNecessaryCondition}
yields $p_{1}\leq p_{2}$, as claimed.

Next, note that our assumptions (in particular, boundedness of $\theta$)
imply that Theorem~\ref{thm:BurnerNecessaryConditionCoarseInFine}
(with $Y,Z,J_{0}$ as above) is applicable. Since $Z=\ell_{v}^{q_{2}}\left(J\right)$
satisfies the Fatou property, this yields boundedness of
\[
\eta:\:\ell_{w}^{q_{1}}\left(\left[\ell^{p_{1}}\left(J_{i}\right)\right]_{i\in I}\right)=Y\left(\left[\ell^{p_{1}}\left(J_{i}\cap J_{0}\right)\right]_{i\in I}\right)\hookrightarrow Z_{\left|\det S_{j}\right|^{p_{1}^{-1}-p_{2}^{-1}}}=\ell_{\bigl[v_{j}\cdot\left|\det S_{j}\right|^{p_{1}^{-1}-p_{2}^{-1}}\bigr]_{j}}^{q_{2}}\left(J\right),
\]
with $\vertiii{\eta}\leq C\cdot\vertiii{\theta}$ for some constant
\[
C=C\left(\dimension,p_{1},p_{2},q_{2},k\left(\CalP,\CalQ\right),\CalQ,\CalP,\varepsilon_{\CalP},C_{\CalQ,\Phi,p_{1}},C_{v,\CalP}\right).
\]
Here, we used that the triangle constant $C_{Z}$ for $Z=\ell_{v}^{q_{2}}\left(J\right)$
only depends on $q_{2}$ and that Lemma~\ref{lem:ModeratelyWeightedSpacesAreRegular}
allows to estimate $\vertiii{\Gamma_{\CalP}}_{Z\to Z}$ in terms of
$q_{2},\CalP,C_{v,\CalP}$.

But in view of Corollary~\ref{cor:EmbeddingCoarseIntoFineSimplification}
(with $r=p_{1}$, $u^{\left(1\right)}\equiv1$ and $u^{\left(2\right)}\equiv1$,
as well as $\bigl(v_{j}\cdot\left|\det S_{j}\right|^{p_{1}^{-1}-p_{2}^{-1}}\bigr)_{j\in J}$
instead of $v$), this entails
\[
\left\Vert \left(w_{i}^{-1}\cdot\left\Vert \left(v_{j}\cdot\left|\det S_{j}\right|^{p_{1}^{-1}-p_{2}^{-1}}\right)_{j\in J_{i}}\right\Vert _{\ell^{q_{2}\cdot\left(p_{1}/q_{2}\right)'}}\right)_{i\in I}\right\Vert _{\ell^{q_{2}\cdot\left(q_{1}/q_{2}\right)'}}\leq C'\cdot\vertiii{\eta}\leq CC'\cdot\vertiii{\theta}
\]
for some constant $C'=C'\left(p_{1},q_{1},q_{2},C_{w,\CalQ},N_{\CalQ},k\left(\CalP,\CalQ\right)\right)$.
But the left-hand side of this estimate is simply $K_{p_{1},1}$.

\medskip{}

Ad (3): In this case, Theorem~\ref{thm:KhinchinNecessaryCoarseInFine}
(and the ensuing remark) show that the identity map
\[
\eta:\ell_{w}^{q_{1}}\Bigl(\left[\ell^{2}\left(J_{i}\right)\right]_{i\in I}\Bigr)=Y\Bigl(\left[\ell^{2}\left(J_{0}\cap J_{i}\right)\right]_{i\in I}\Bigr)\hookrightarrow Z|_{J_{0}}=\ell_{v}^{q_{2}}\left(J\right)
\]
is well-defined and bounded, with $\vertiii{\eta}\leq C\cdot\vertiii{\theta}$
for some constant
\[
C=C\left(\dimension,p_{1},q_{1},q_{2},\CalQ,\CalP,C_{\CalQ,\Phi,p_{1}},C_{\CalP,\Psi,p_{2}},C_{w,\CalQ},C_{v,\CalP}\right).
\]
Here, we again used that the triangle constant $C_{Z}$ of $Z=\ell_{v}^{q_{2}}\left(J\right)$
only depends on $q_{2}$ and that Lemma~\ref{lem:ModeratelyWeightedSpacesAreRegular}
allows to estimate $\vertiii{\Gamma_{\CalQ}}_{Y\to Y}$ and $\vertiii{\Gamma_{\CalP}}_{Z\to Z}$
in terms of $q_{1},q_{2},\CalQ,\CalP,C_{w,\CalQ}$ and $C_{v,\CalP}$.

The remainder of the proof proceeds precisely as in the previous case
(with $r=2$ instead of $r=p_{1}$), noting that we have $\bigl(v_{j}\cdot\left|\det S_{j}\right|^{p_{1}^{-1}-p_{2}^{-1}}\bigr)_{j\in J}=v$,
since $p_{1}=p_{2}$.

\medskip{}

Ad (4): Let us first prove the estimate $K\lesssim\vertiii{\theta}$.
To do this, it suffices to verify that the assumptions of Theorem~\ref{thm:NecessaryConditionForModerateCoveringCoarseInFine}
are met; these consist of the following:

\begin{enumerate}
\item $\CalQ$ is tight, and $\CalP_{J_{0}}=\CalP$ is relatively $\CalQ$-moderate.
This holds by assumption.
\item The weight $v=v|_{J_{0}}$ is relatively $\CalQ$-moderate. Again,
this holds by assumptions.
\item There is $r\in\N_{0}$ and some $C_{0}>0$ satisfying
\[
\lambda\left(Q_{i}\right)\leq C_{0}\cdot\lambda\left(\,\smash{\bigcup_{j\in J_{0}\cap J_{i}}}\,\vphantom{\bigcup}P_{j}^{r\ast}\,\right)\vphantom{\bigcup_{j\in J_{0}\cap J_{i}}}\qquad\forall\,i\in I_{0}\text{ for a certain set }I_{0}\subset I\,.
\]
But since we are assuming $\CalO=\CalO'$, we have $Q_{i}\subset\CalO=\CalO'=\bigcup_{j\in J}P_{j}$.
Since we also have $J_{0}=J$, we see that the preceding estimate
is satisfied for $C_{0}=1$ and $r=0$.
\item The assumptions of Lemma~\ref{lem:NecessaryConjugateCoarseInFineWithoutModerateness}
are satisfied. This is indeed the case, as we will verify now.
\end{enumerate}
First, $\CalP=\CalP_{J_{0}}$ is almost subordinate to $\CalQ$, as
required in Lemma~\ref{lem:NecessaryConjugateCoarseInFineWithoutModerateness}.
Finally, we want to verify that $\theta$ satisfies the properties
of the map $\iota$ from that lemma. Hence, we have to verify
\[
\left\langle \theta f,\,\varphi\right\rangle _{\CalD'}=\left\langle f,\,\varphi\right\rangle _{\CalD'}\qquad\forall\,\varphi\in\TestFunctionSpace{\CalO\cap\CalO'}\text{ and }f\in\CalD_{K}^{\CalQ,p_{1},\ell_{w}^{q_{1}}\left(I\right)},
\]
for a certain set $K\subset\CalO$, which is defined in the statement
of Lemma~\ref{lem:NecessaryConjugateCoarseInFineWithoutModerateness},
but whose precise definition is unimportant for us. Indeed, we have
$K\subset\CalO=\CalO'$ and hence $\theta f=f$ for all $f\in\CalD_{K}^{\CalQ,p_{1},\ell_{w}^{q_{1}}\left(I\right)}\subset\TestFunctionSpace{\CalO'}$,
which easily implies the desired identity from above.

All in all, we see that Lemma~\ref{lem:NecessaryConjugateCoarseInFineWithoutModerateness},
and thus also Theorem~\ref{thm:NecessaryConditionForModerateCoveringCoarseInFine}
are applicable, so that we get $K\lesssim\vertiii{\theta}$, with
an implied constant as in the statement of the current theorem. In
view of the second part of the current theorem, it is clear that the
boundedness of $\theta$ also entails $p_{1}\leq p_{2}$.

Next, if $\iota$ is bounded, it is clear that $\iota|_{\TestFunctionSpace{\CalO'}}=\theta$,
so that we get $K\lesssim\vertiii{\,\iota|_{\TestFunctionSpace{\CalO'}}\,}\leq\vertiii{\iota}$
and $p_{1}\leq p_{2}$.

Hence, it remains to show that if $p_{1}\leq p_{2}$ and if $K$ is
finite, then $\iota$ is well-defined and bounded with $\vertiii{\iota}\lesssim K$.
In view of the first part of the theorem, it is thus sufficient to
verify $K_{\UpperExpo{p_{1}},2}\lesssim K$. But this is a direct
consequence of Remark~\ref{rem:SufficientCoarseIntoFineSimplification};
see in particular equation~(\ref{eq:SufficientCoarseIntoFineSimplificationRelativelyModerate}).
\end{proof}
One easy, but nevertheless frequently useful special case is when
$\CalQ$ and $\CalP$ \emph{coincide}.
\begin{cor}
\label{cor:SummarySameCovering}In addition to our standing assumptions,
assume that $\CalQ=\CalP$. Then, the identity map
\[
\iota:\FourierDecompSp{\CalQ}{p_{1}}{\ell_{w}^{q_{1}}}\to\FourierDecompSp{\CalQ}{p_{2}}{\ell_{v}^{q_{2}}},f\mapsto f
\]
is well-defined and bounded if and only if we have $p_{1}\leq p_{2}$
and if
\[
K:=\left\Vert \left(\left|\det T_{i}\right|^{p_{1}^{-1}-p_{2}^{-1}}\cdot\frac{v_{i}}{w_{i}}\right)_{i\in I}\right\Vert _{\ell^{q_{2}\cdot\left(q_{1}/q_{2}\right)'}}<\infty.
\]

More precisely, in case of $p_{1}\leq p_{2}$, we have $\vertiii{\iota}\asymp K$,
where the implied constants only depend on 
\[
\dimension,p_{1},p_{2},q_{1},q_{2},\CalQ,\varepsilon_{\CalQ},C_{w,\CalQ},C_{v,\CalQ},C_{\CalQ,\Phi,p_{1}}.\qedhere
\]
\end{cor}

\begin{proof}
We apply part~(\ref{enu:SummaryCoarseInFineModerate}) of Theorem~\ref{thm:SummaryCoarseIntoFine}
(with $\CalQ=\CalP$). First of all, $\CalP(=\CalQ)$ is almost subordinate
to $\CalQ$ with $k\left(\CalP,\CalQ\right)=0$, since  $P_{j}=Q_{j}\subset Q_{j}^{0\ast}$
for all $j\in J=I$. Next, note that $\CalP(=\CalQ)$ and $v$ are
clearly relatively $\CalQ$-moderate, with $C_{{\rm mod}}\left(\CalP,\CalQ\right)=C\left(\dimension,C_{\CalQ}\right)$
and $C_{v,\CalP,\CalQ}=C\left(C_{v,\CalQ}\right)$; see equation~(\ref{eq:DeterminantIsModerate}).
Furthermore, $\CalQ$ and $\CalP(=\CalQ)$ both cover the same set
$\CalO$, so that part~(\ref{enu:SummaryCoarseInFineModerate}) of
Theorem~\ref{thm:SummaryCoarseIntoFine} is indeed applicable.

Thus, we see that $\iota$ is well-defined and bounded if and only
if we have $p_{1}\leq p_{2}$ and if
\[
\widetilde{K}:=\left\Vert \left(\frac{v_{j_{i}}}{w_{i}}\cdot\left|\det T_{i}\right|^{s}\cdot\left|\det\smash{S_{j_{i}}}\right|^{p_{1}^{-1}-p_{2}^{-1}-s}\right)_{i\in I}\right\Vert _{\ell^{q_{2}\cdot\left(q_{1}/q_{2}\right)'}}
\]
is finite, where $s=\left(\frac{1}{q_{2}}-\smash{\frac{1}{\SignedUpperExpo{p_{1}}}}\right)_{+}$
and where for each $i\in I$, some $j_{i}\in J_{i}$ can be arbitrarily
selected. Because of $\CalQ=\CalP$, we can simply choose $j_{i}:=i$,
so that we get
\[
\widetilde{K}=\left\Vert \left(\left|\det T_{i}\right|^{p_{1}^{-1}-p_{2}^{-1}}\cdot\frac{v_{i}}{w_{i}}\right)_{i\in I}\right\Vert _{\ell^{q_{2}\cdot\left(q_{1}/q_{2}\right)'}}.
\]
Now, all claims follow from part~(\ref{enu:SummaryCoarseInFineModerate})
of Theorem~\ref{thm:SummaryCoarseIntoFine}.
\end{proof}
As our next result, compared to Theorem~\ref{thm:SummaryCoarseIntoFine},
we consider the ``reverse'' case, in which $\CalQ$ is almost subordinate
to $\CalP$. We remark that the following theorem is an improved version
of \cite[Theorem 5.4.4]{VoigtlaenderPhDThesis} from my PhD thesis.
\begin{thm}
\label{thm:SummaryFineIntoCoarse}In addition to our standing assumptions,
assume that $\CalQ$ is almost subordinate to $\CalP$. Note that
this entails $\CalO\subset\CalO'$.

For $j\in J$, let
\[
I_{j}:=\left\{ i\in I\with Q_{i}\cap P_{j}\neq\emptyset\right\} 
\]
and for $r\in\left(0,\infty\right]$, define
\begin{equation}
K_{r}:=\left\Vert \left(v_{j}\cdot\left\Vert \left(\left|\det T_{i}\right|^{p_{1}^{-1}-p_{2}^{-1}}/\,w_{i}\right)_{i\in I_{j}}\right\Vert _{\ell^{r\cdot\left(q_{1}/r\right)'}}\right)_{j\in J}\right\Vert _{\ell^{q_{2}\cdot\left(q_{1}/q_{2}\right)'}}\in\left[0,\infty\right].\label{eq:SummaryFineIntoCoarseCriticalQuantitiy}
\end{equation}
Then, the following hold:

\begin{enumerate}[leftmargin=0.7cm]
\item \label{enu:SummaryFineInCoarseSufficientNonModerate}If we have $p_{1}\leq p_{2}$
and if $K_{\LowerExpo{p_{2}}}<\infty$, then
\[
\iota:\FourierDecompSp{\CalQ}{p_{1}}{\ell_{w}^{q_{1}}}\to\FourierDecompSp{\CalP}{p_{2}}{\ell_{v}^{q_{2}}},f\mapsto\sum_{i\in I}\varphi_{i}f
\]
is well-defined and bounded, with $\vertiii{\iota}\leq C_{1}\cdot K_{\LowerExpo{p_{2}}}$
for some constant
\[
C_{1}=C_{1}\left(\dimension,p_{1},p_{2},q_{1},q_{2},\CalQ,\CalP,k\left(\CalQ,\CalP\right),C_{\CalQ,\Phi,p_{1}},C_{\CalP,\Psi,p_{2}},C_{v,\CalP}\right).
\]
Furthermore, the following hold:

\begin{enumerate}
\item \label{enu:SummaryFineIntoCoarseEmbeddingYieldsExtension}For each
$f\in\FourierDecompSp{\CalQ}{p_{1}}{\ell_{w}^{q_{1}}}\subset\DistributionSpace{\CalO}$,
the distribution $\iota f\in\DistributionSpace{\CalO'}$ is an extension
of $f$ onto $\TestFunctionSpace{\CalO'}$, i.e., $\left\langle \iota f,\,g\right\rangle _{\CalD'}=\left\langle f,\,g\right\rangle _{\CalD'}$
for all $g\in\TestFunctionSpace{\CalO}\subset\TestFunctionSpace{\CalO'}$.\vspace{0.1cm}
\item \label{enu:SummaryFineIntoCoarseEmbeddingConsistency}For each $f\in\TestFunctionSpace{\CalO}$,
we have $\iota f=f$.
\end{enumerate}
\item \label{enu:SummaryFineInCoarseNecessaryNonModerate}Conversely, if
the identity map
\[
\theta:\left(\TestFunctionSpace{\CalO},\left\Vert \mybullet\right\Vert _{\FourierDecompSp{\CalQ}{p_{1}}{\ell_{w}^{q_{1}}}}\right)\to\FourierDecompSp{\CalP}{p_{2}}{\ell_{v}^{q_{2}}},f\mapsto f
\]
is bounded, then we have $p_{1}\leq p_{2}$ and 
\[
K_{s}\leq C_{2}\cdot\vertiii{\theta}<\infty\qquad\text{ with }\qquad s:=\begin{cases}
p_{2}, & \text{if }p_{2}<\infty,\\
1=\LowerExpo{p_{2}}, & \text{if }p_{2}=\infty,
\end{cases}
\]
for some constant
\[
C_{2}=C_{2}\left(\dimension,p_{1},p_{2},q_{1},q_{2},\CalQ,\varepsilon_{\CalQ},\CalP,k\left(\CalQ,\CalP\right),C_{\CalQ,\Phi,p_{1}},C_{\CalP,\Psi,p_{2}},C_{w,\CalQ},C_{v,\CalP}\right).
\]
\item \label{enu:SummaryFineInCoarseNecessaryKhinchin}Under the assumptions
of the preceding point, if $p_{1}=p_{2}=:p$ and
\[
s:=\begin{cases}
1 & \text{if }p=\infty,\\
\min\left\{ 2,p\right\} , & \text{if }p<\infty,
\end{cases}
\]
then $K_{s}\leq C_{3}\cdot\vertiii{\theta}<\infty$ for some constant
\[
C_{3}=C_{3}\left(\dimension,p_{1},q_{1},q_{2},\CalQ,\CalP,k\left(\CalQ,\CalP\right),C_{\CalQ,\Phi,p_{1}},C_{\CalP,\Psi,p_{2}},C_{w,\CalQ},C_{v,\CalP}\right).
\]
\item \label{enu:SummaryFineInCoarseModerate}Finally, if we have $\CalO=\CalO'$
and if $\CalQ$ and $w$ are relatively $\CalP$-moderate, then we
have the following equivalence (with $\iota,\theta$ as above): Choosing
for each $j\in J$ some $i_{j}\in I$ with $P_{j}\cap Q_{i_{j}}\neq\emptyset$,
we have:
\begin{align*}
 & \iota\text{ well-defined and bounded}\\
\qquad\qquad\Longleftrightarrow\: & \theta\text{ well-defined and bounded}\\
\Longleftrightarrow\: & p_{1}\leq p_{2}\text{ and}\\
 & K:=\left\Vert \!\left(\!\frac{v_{j}}{w_{i_{j}}}\cdot\left|\det\smash{T_{i_{j}}}\right|^{p_{1}^{-1}-p_{2}^{-1}-s}\cdot\left|\det S_{j}\right|^{s}\right)_{\!\!\!j\in J}\right\Vert _{\ell^{q_{2}\cdot\left(q_{1}/q_{2}\right)'}}\!\!\!<\infty,\text{ where }s:=\left(\smash{\frac{1}{\LowerExpo{p_{2}}}}-\frac{1}{q_{1}}\right)_{+}.
\end{align*}
Furthermore, given $p_{1}\leq p_{2}$, we have
\[
\vertiii{\iota}\asymp\vertiii{\theta}\asymp K,
\]
where the implied constants only depend on
\[
\qquad\dimension,p_{1},p_{2},q_{1},q_{2},\CalQ,\CalP,\varepsilon_{\CalQ},\varepsilon_{\CalP},k\left(\CalQ,\CalP\right),C_{{\rm mod}}\left(\CalQ,\CalP\right),C_{w,\CalQ,\CalP},C_{w,\CalQ},C_{v,\CalP},C_{\CalQ,\Phi,p_{1}},C_{\CalP,\Psi,p_{2}}.\qedhere
\]
\end{enumerate}
\end{thm}

\begin{proof}
For the whole proof, set $Y:=\ell_{w}^{q_{1}}\left(I\right)$, $Z:=\ell_{v}^{q_{2}}\left(J\right)$
and $I_{0}:=I$. We now prove each statement individually.

Ad (1): Remark~\ref{rem:SufficientFineInCoarseSimplification} shows
that the embedding $\eta$ from Corollary~\ref{cor:EmbeddingFineIntoCoarse}
(with $I_{0}=I$) satisfies
\[
\vertiii{\eta}\asymp\left\Vert \left(v_{j}\cdot\left\Vert \left(\left|\det T_{i}\right|^{p_{1}^{-1}-p_{2}^{-1}}/\,w_{i}\right)_{i\in I_{j}}\right\Vert _{\ell^{\LowerExpo{p_{2}}\cdot\left(q_{1}/\LowerExpo{p_{2}}\right)'}}\right)_{j\in J}\right\Vert _{\ell^{q_{2}\cdot\left(q_{1}/q_{2}\right)'}}=K_{\LowerExpo{p_{2}}},
\]
where the implied constant only depends on $q_{1},q_{2},p_{2},N_{\CalP},k\left(\CalQ,\CalP\right),C_{v,\CalP}$.

Since we also assume $p_{1}\leq p_{2}$, Corollary~\ref{cor:EmbeddingFineIntoCoarse}
shows that the stated map $\iota$ is well-defined and bounded and
satisfies property (\ref{enu:SummaryFineIntoCoarseEmbeddingYieldsExtension}),
since $I_{0}=I$. Furthermore, since $\Phi=\left(\varphi_{i}\right)_{i\in I}$
is a locally finite partition of unity on $\CalO$, we have for $g\in\TestFunctionSpace{\CalO}$
that $g=\sum_{i\in I}\varphi_{i}g=\iota g$, where only finitely many
terms do not vanish.

Finally, Corollary~\ref{cor:EmbeddingFineIntoCoarse} even yields
$\vertiii{\iota}\lesssim\vertiii{\eta}$, where the implied constant
only depends on quantities mentioned in the current theorem—at least,
once we use Lemma~\ref{lem:ModeratelyWeightedSpacesAreRegular} to
estimate $\vertiii{\Gamma_{\CalP}}_{Z\to Z}$ in terms of $q_{2},\CalP,C_{v,\CalP}$.

\medskip{}

Ad (2): The present assumptions show that Lemma~\ref{lem:SimpleNecessaryCondition}
(with $K=\CalO\subset\CalO\cap\CalO'$) is applicable. Now, for $i\in I$
and arbitrary $\xi\in Q_{i}\subset\CalO\subset\CalO'$, we get $\xi\in Q_{i}\cap P_{j}\neq\emptyset$
for some $j\in J$. But since $\CalQ,\CalP$ are open coverings, we
get $K^{\circ}\cap Q_{i}^{\circ}\cap P_{j}^{\circ}\neq\emptyset$,
so that Lemma~\ref{lem:SimpleNecessaryCondition} yields $p_{1}\leq p_{2}$.

Next, note that $Z=\ell_{v}^{q_{2}}\left(J\right)$ satisfies the
Fatou property and that our assumptions imply that Theorem~\ref{thm:BurnerNecessaryConditionFineInCoarse}
is applicable, so that the embedding
\[
\eta:\ell_{\bigl(\left|\det T_{i}\right|^{p_{2}^{-1}-p_{1}^{-1}}\cdot w_{i}\bigr)_{i}}^{q_{1}}\left(I\right)=\left(Y|_{I_{0}}\right)_{\left|\det T_{i}\right|^{p_{2}^{-1}-p_{1}^{-1}}}\hookrightarrow Z\bigl(\left[\ell^{s}\left(I_{0}\cap I_{j}\right)\right]_{j\in J}\bigr)=\ell_{v}^{q_{2}}\bigl(\left[\ell^{s}\left(I_{j}\right)\right]_{j\in J}\bigr)
\]
is well-defined and bounded, with $\vertiii{\eta}\leq C\cdot\vertiii{\theta}$
for some constant
\[
C=C\left(\dimension,p_{1},p_{2},q_{1},q_{2},\CalQ,\varepsilon_{\CalQ},\CalP,k\left(\CalQ,\CalP\right),C_{\CalQ,\Phi,p_{1}},C_{\CalP,\Psi,p_{2}},C_{w,\CalQ},C_{v,\CalP}\right).
\]
Here, we again used Lemma~\ref{lem:ModeratelyWeightedSpacesAreRegular}
to estimate $\vertiii{\Gamma_{\CalQ}}_{Y\to Y}$ and $\vertiii{\Gamma_{\CalP}}_{Z\to Z}$
in terms of $q_{1},q_{2},\CalQ,\CalP,C_{w,\CalQ}$ and $C_{v,\CalP}$,
and we also used that the triangle constant $C_{Z}$ of $Z=\ell_{v}^{q_{2}}\left(J\right)$
only depends on $q_{2}$.

But in view of Corollary~\ref{cor:EmbeddingFineInCoarseSimplification}
(with $J_{0}=J$, $u\equiv1$, $r=s$ and $\bigl(\left|\det T_{i}\right|^{p_{2}^{-1}-p_{1}^{-1}}\cdot w_{i}\bigr)_{i\in I}$
instead of $w$), the boundedness of $\eta$ entails
\[
\left\Vert \left(v_{j}\cdot\left\Vert \left(\left|\det T_{i}\right|^{p_{1}^{-1}-p_{2}^{-1}}/\,w_{i}\right)_{i\in I_{j}}\right\Vert _{\ell^{s\cdot\left(q_{1}/s\right)'}}\right)_{j\in J}\right\Vert _{\ell^{q_{2}\cdot\left(q_{1}/q_{2}\right)'}}\leq C'\cdot\vertiii{\eta}\leq C'C\cdot\vertiii{\theta}
\]
for some constant $C'=C'\left(q_{1},q_{2},p_{2},\CalP,k\left(\CalQ,\CalP\right),C_{v,\CalP}\right)$.
The proof is completed by noting that the left-hand side of this estimate
is simply $K_{s}$.

\medskip{}

Ad (3): In this case, Theorem~\ref{thm:KhinchinNecessaryFineInCoarse}
shows that
\[
\eta:Y|_{I_{0}}\hookrightarrow Z\bigl(\left[\ell^{s}\left(I_{0}\cap I_{j}\right)\right]_{j\in J}\bigr)
\]
is well-defined and bounded, with $\vertiii{\eta}\leq C\cdot\vertiii{\theta}$
for some constant
\[
C=C\left(\dimension,p_{1},q_{1},q_{2},\CalQ,\CalP,k\left(\CalQ,\CalP\right),C_{\CalQ,\Phi,p_{1}},C_{\CalP,\Psi,p_{2}},C_{w,\CalQ},C_{v,\CalP}\right).
\]
Here, we again used that the triangle constant $C_{Z}$ of $Z=\ell_{v}^{q_{2}}\left(J\right)$
only depends on $q_{2}$ and that Lemma~\ref{lem:ModeratelyWeightedSpacesAreRegular}
allows to estimate $\vertiii{\Gamma_{\CalQ}}_{Y\to Y}$ and $\vertiii{\Gamma_{\CalP}}_{Z\to Z}$
in terms of $q_{1},q_{2},\CalQ,\CalP,C_{w,\CalQ}$ and $C_{v,\CalP}$.

The rest of the proof is as in the previous case, noting that we have
$\bigl(\left|\det T_{i}\right|^{p_{2}^{-1}-p_{1}^{-1}}\cdot w_{i}\bigr)_{i\in I}=w$,
since $p_{1}=p_{2}$.

\medskip{}

Ad (4): Under the present assumptions, let us first assume that $\theta$
is well-defined and bounded. As seen above, this yields $p_{1}\leq p_{2}$,
so that it suffices to show $K\lesssim\vertiii{\theta}$, with an
implied constant as stated. To see this, we want to apply Theorem~\ref{thm:NecessaryConditionForModerateCoveringFineInCoarse},
with $J_{0}:=J$. In the notation of that theorem, this implies
\begin{align*}
I_{0}=\left\{ i\in I\with J_{i}\cap J_{0}\neq\emptyset\right\}  & =\left\{ i\in I\with J_{i}\neq\emptyset\right\} \\
 & =\left\{ i\in I\with\exists\,j\in J:\,Q_{i}\cap P_{j}\neq\emptyset\right\} \\
\left({\scriptstyle \text{since }\CalP\text{ covers }\CalO'}\right) & =\left\{ i\in I\with Q_{i}\cap\CalO'\neq\emptyset\right\} =I,
\end{align*}
where the last step used that we have $\CalO'=\CalO$ and that each
$Q_{i}$ is nonempty. Now, let us verify assumptions (\ref{enu:NecessaryConditionModerateCoveringFineInCoarseAlmostTightness})–(\ref{enu:NecessaryConditionModerateCoveringFineInCoarseWeightModerate})
of Theorem~\ref{thm:NecessaryConditionForModerateCoveringFineInCoarse}:

\begin{enumerate}
\item Since $\CalP=\left(S_{j}P_{j}'+c_{j}\right)_{j\in J}$ is a tight
semi-structured covering, there is for $\varepsilon:=\varepsilon_{\CalP}$
for each $j\in J=J_{0}$ some $\xi_{j}=c_{j}\in\R^{\dimension}$ with
$B_{\varepsilon}\left(\xi_{j}\right)\subset P_{j}'$. Particularly,
$S_{j}\left[B_{\varepsilon}\left(\xi_{j}\right)\right]+c_{j}\subset P_{j}\subset\CalO'=\CalO$,
as desired.
\item By assumption, $\CalQ_{I_{0}}=\CalQ$ is almost subordinate to $\CalP$.
\item By assumption, $\CalQ_{I_{0}}=\CalQ$ is relatively $\CalP$-moderate.
\item By assumption, the weight $w=w|_{I_{0}}$ is relatively $\CalP$-moderate.
\end{enumerate}
Finally, in view of $\bigcup_{i\in I_{0}}Q_{i}=\CalO$ we see that
the map $\iota$ from Theorem~\ref{thm:NecessaryConditionForModerateCoveringFineInCoarse}
coincides with $\theta$ from the current theorem. In particular,
$\iota$ (as in Theorem~\ref{thm:NecessaryConditionForModerateCoveringFineInCoarse})
is bounded, so that Theorem~\ref{thm:NecessaryConditionForModerateCoveringFineInCoarse}
shows $K\leq C\vertiii{\theta}<\infty$ for some constant
\[
C=C\left(\dimension,p_{1},p_{2},q_{1},q_{2},\CalQ,\CalP,\varepsilon_{\CalQ},\varepsilon_{\CalP},k\left(\CalQ,\CalP\right),C_{{\rm mod}}\left(\CalQ,\CalP\right),C_{w,\CalQ},C_{v,\CalP},C_{w,\CalQ,\CalP},C_{\CalQ,\Phi,p_{1}},C_{\CalP,\Psi,p_{2}}\right),
\]
as desired.

Now, by the properties of $\iota$ shown in the first part of the
present theorem, we see that $\theta$ is a restriction of $\iota$.
In particular, if $\iota$ is bounded, then so is $\theta$, so that
we get $p_{1}\leq p_{2}$ and $K\lesssim\vertiii{\theta}\leq\vertiii{\iota}$.

Finally, we need to show that $p_{1}\leq p_{2}$ and $K<\infty$ together
imply that $\iota$ is bounded. In view of the first part of the current
theorem, it suffices to show $K_{\LowerExpo{p_{2}}}\lesssim K$. But
this is an immediate consequence of the second part of Corollary~\ref{cor:EmbeddingFineInCoarseSimplification}
(with $J_{0}=J$, $u_{i}=\left|\det T_{i}\right|^{p_{1}^{-1}-p_{2}^{-1}}$,
$r=\LowerExpo{p_{2}}$, and $I_{0}=I$), once we verify its assumptions.
But the only assumptions which are not obviously implied by the present
ones are the relative moderateness of $u=\bigl(\left|\det T_{i}\right|^{p_{1}^{-1}-p_{2}^{-1}}\bigr)_{i\in I}$—which
holds since $\CalQ$ is relatively $\CalP$-moderate—and the existence
of $s\in\N_{0}$ and $C_{0}>0$ satisfying
\[
\lambda\left(P_{j}\right)\leq C_{0}\cdot\lambda\left(\,\smash{\bigcup_{i\in I_{0}\cap I_{j}}}\,\vphantom{\bigcup}Q_{i}^{s\ast}\,\right)\vphantom{\bigcup_{i\in I_{0}\cap I_{j}}}\qquad\forall\,j\in J\,.
\]
But since we have $P_{j}\subset\CalO'=\CalO$ and since $\CalQ$ covers
$\CalO$, we have $P_{j}\subset\bigcup_{i\in I_{j}}Q_{i}$. Since
we have $I_{0}=I$, we can thus take $C_{0}=1$ and $s=0$.

Finally, the implicit constant in the estimate $K_{\LowerExpo{p_{2}}}\lesssim K$
provided by Corollary~\ref{cor:EmbeddingFineInCoarseSimplification}
only depends on
\[
\dimension,q_{1},q_{2},p_{1},p_{2},\CalQ,\CalP,\varepsilon_{\CalQ},\varepsilon_{\CalP},k\left(\CalQ,\CalP\right),C_{{\rm mod}}\left(\CalQ,\CalP\right),C_{w,\CalQ,\CalP},C_{v,\CalP}.\qedhere
\]
\end{proof}
Our final theorem in this summary is concerned with the ``intermediate''
case in which some part of $\CalQ$ is almost subordinate to $\CalP$,
while some part of $\CalP$ is almost subordinate to $\CalQ$.

Strictly speaking, the two preceding theorems are special cases of
this theorem. We nevertheless stated them separately, since their
application is \emph{much} more convenient than invoking the following
theorem, which is an improved version of \cite[Theorem 5.4.5]{VoigtlaenderPhDThesis}
from my PhD thesis.
\begin{thm}
\label{thm:SummaryMixedSubordinateness}In addition to our standing
assumptions, suppose that we have
\[
\emptyset\neq\CalO\cap\CalO'=A\cup B
\]
for certain sets $A,B\subset\R^{\dimension}$, and such that the following
additional properties are satisfied:

\begin{enumerate}
\item With $I_{A}:=\left\{ i\in I\with Q_{i}\cap A\neq\emptyset\right\} $,
the family $\CalQ_{I_{A}}:=\left(Q_{i}\right)_{i\in I_{A}}$ is almost
subordinate to $\CalP$.
\item With $J_{B}:=\left\{ j\in J\with P_{j}\cap B\neq\emptyset\right\} $,
the family $\CalP_{J_{B}}:=\left(P_{j}\right)_{j\in J_{B}}$ is almost
subordinate to $\CalQ$.
\end{enumerate}
For $r\in\left(0,\infty\right]$ and $\ell\in\left\{ 1,2\right\} $,
define
\begin{align*}
K_{r}^{\left(1\right)} & :=\left\Vert \left(v_{j}\cdot\left\Vert \left(\left|\det T_{i}\right|^{p_{1}^{-1}-p_{2}^{-1}}/\,w_{i}\right)_{i\in I_{j}\cap I_{A}}\right\Vert _{\ell^{r\cdot\left(q_{1}/r\right)'}}\right)_{j\in J}\right\Vert _{\ell^{q_{2}\cdot\left(q_{1}/q_{2}\right)'}}\in\left[0,\infty\right]\\
K_{r}^{\left(2,\ell\right)} & :=\left\Vert \left(w_{i}^{-1}\cdot\left\Vert \left(v_{j}\,/\,\smash{u_{i,j}^{\left(\ell\right)}}\,\right)_{j\in J_{i}\cap J_{B}}\right\Vert _{\ell^{q_{2}\cdot\left(r/q_{2}\right)'}}\right)_{i\in I}\right\Vert _{\ell^{q_{2}\cdot\left(q_{1}/q_{2}\right)'}}\in\left[0,\infty\right],
\end{align*}
with
\[
u_{i,j}^{\left(1\right)}:=\left|\det S_{j}\right|^{p_{2}^{-1}-p_{1}^{-1}}\qquad\text{ and }\qquad u_{i,j}^{\left(2\right)}:=\begin{cases}
\left|\det S_{j}\right|^{p_{2}^{-1}-1}\cdot\left|\det T_{i}\right|^{1-p_{1}^{-1}}, & \text{if }p_{1}<1,\\
\vphantom{\rule{0.1cm}{0.55cm}}\left|\det S_{j}\right|^{p_{2}^{-1}-p_{1}^{-1}}, & \text{if }p_{1}\geq1.
\end{cases}
\]
Then the following hold:

\begin{enumerate}[leftmargin=0.65cm]
\item \label{enu:SummaryMixedSubordinatenessSufficientNonModerate}If $p_{1}\leq p_{2}$
and if $K_{\LowerExpo{p_{2}}}^{\left(1\right)}<\infty$ and $K_{\UpperExpo{p_{1}}}^{\left(2,2\right)}<\infty$,
then there is a bounded linear map
\[
\iota:\FourierDecompSp{\CalQ}{p_{1}}{\ell_{w}^{q_{1}}}\to\FourierDecompSp{\CalP}{p_{2}}{\ell_{v}^{q_{2}}}
\]
with $\vertiii{\iota}\lesssim K_{\LowerExpo{p_{2}}}^{\left(1\right)}+K_{\UpperExpo{p_{1}}}^{\left(2,2\right)}$,
where the implied constant only depends on
\[
\dimension,p_{1},p_{2},q_{1},q_{2},\CalQ,\CalP,k\left(\smash{\CalQ_{I_{A}}},\CalP\right),k\left(\smash{\CalP_{J_{B}}},\CalQ\right),C_{w,\CalQ},C_{v,\CalP},C_{\CalQ,\Phi,p_{1}},C_{\CalP,\Psi,p_{1}}\:,
\]
and with the following additional properties:

\begin{enumerate}[leftmargin=0.6cm]
\item For $f\in\FourierDecompSp{\CalQ}{p_{1}}{\ell_{w}^{q_{1}}}\subset\DistributionSpace{\CalO}$,
we have
\[
\left\langle f,\,g\right\rangle _{\CalD'}=\left\langle \iota f,\,g\right\rangle _{\CalD'}\qquad\forall\,g\in\TestFunctionSpace{\CalO\cap\CalO'}.
\]
In particular, if $\CalO=\CalO'$, then $\iota f=f$ for all $f\in\FourierDecompSp{\CalQ}{p_{1}}{\ell_{w}^{q_{1}}}\subset\DistributionSpace{\CalO}$.
\item If $f\in\FourierDecompSp{\CalQ}{p_{1}}{\ell_{w}^{q_{1}}}$ is given
by (integration against) a measurable function $f:\R^{\dimension}\to\Compl$
with $f\in L_{{\rm loc}}^{1}\left(\CalO\cup\CalO'\right)$ and with
$f=0$ almost everywhere on $\CalO'\setminus\CalO$, then $\iota f=f$
as elements of $\DistributionSpace{\CalO'}$.
\item In particular, $\iota f=f$ for all $f\in\TestFunctionSpace{\CalO\cap\CalO'}$.
\item For $f\in\FourierDecompSp{\CalQ}{p_{1}}{\ell_{w}^{q_{1}}}$, we have
$\supp\iota f\subset\CalO'\cap\overline{\supp f}$, where the closure
$\overline{\supp f}$ is taken in $\R^{\dimension}$.
\end{enumerate}
\item \label{enu:SummaryMixedSubordinatenessNecessaryNonModerate}Conversely,
if the identity map
\[
\theta:\left(\TestFunctionSpace{\CalO\cap\CalO'},\left\Vert \mybullet\right\Vert _{\FourierDecompSp{\CalQ}{p_{1}}{\ell_{w}^{q_{1}}}}\right)\to\FourierDecompSp{\CalP}{p_{2}}{\ell_{v}^{q_{2}}},f\mapsto f
\]
is bounded, then the following hold:\vspace{0.1cm}

\begin{enumerate}
\item $p_{1}\leq p_{2}$,
\item $K_{s}^{\left(1\right)}\lesssim\vertiii{\theta}$ with $s=\begin{cases}
1, & \text{if }p_{2}=\infty,\\
p_{2}, & \text{otherwise},
\end{cases}$
\item $K_{p_{1}}^{\left(2,1\right)}\lesssim\vertiii{\theta}$.\vspace{0.15cm}
\end{enumerate}
Here, all implied constants only depend on
\[
\qquad\qquad\dimension,p_{1},p_{2},q_{1},q_{2},\CalQ,\CalP,\varepsilon_{\CalQ},\varepsilon_{\CalP},k\left(\smash{\CalQ_{I_{A}}},\CalP\right),k\left(\smash{\CalP_{J_{B}}},\CalQ\right),C_{w,\CalQ},C_{v,\CalP},C_{\CalQ,\Phi,p_{1}},C_{\CalP,\Psi,p_{2}}.
\]

\item \label{enu:SummaryMixedSubordinatenessNecessaryKhinchin}Under the
assumptions of the previous case, if additionally $p_{1}=p_{2}$,
we also have $K_{2}^{\left(1\right)}\lesssim\vertiii{\theta}$ and
$K_{2}^{\left(2,1\right)}\lesssim\vertiii{\theta}$, where the implied
constants only depend on 
\[
\dimension,p_{1},q_{1},q_{2},\CalQ,\CalP,k\left(\smash{\CalQ_{I_{A}}},\CalP\right),k\left(\smash{\CalP_{J_{B}}},\CalQ\right),C_{w,\CalQ},C_{v,\CalP},C_{\CalQ,\Phi,p_{1}},C_{\CalP,\Psi,p_{2}}.
\]
\item \label{enu:SummaryMixedSubordinatenessModerate}Finally, if

\begin{enumerate}
\item $\CalO=\CalO'$, and
\item $\CalQ_{I_{A}}$ and $w|_{I_{A}}$ are relatively $\CalP$-moderate,
and
\item $\CalP_{J_{B}}$ and $v|_{J_{B}}$ are relatively $\CalQ$-moderate,
\end{enumerate}
then the following equivalence holds:
\[
\qquad\theta\text{ bounded}\quad\Longleftrightarrow\quad\iota\text{ bounded}\quad\Longleftrightarrow\quad\left(p_{1}\leq p_{2}\quad\text{and}\quad K^{\left(1\right)}<\infty\quad\text{and}\quad K^{\left(2\right)}<\infty\right)\,,
\]
where for each $j\in J\setminus J_{B}$, some $i_{j}\in I_{j}$ and
for each $i\in I^{\left(B\right)}$, some $j_{i}\in J_{i}\cap J_{B}$
is selected, to define
\begin{align*}
\quad\quad\quad K^{\left(1\right)} & :=\left\Vert \!\left(\!\frac{v_{j}}{w_{i_{j}}}\cdot\left|\det\smash{T_{i_{j}}}\right|^{p_{1}^{-1}-p_{2}^{-1}-s_{1}}\cdot\left|\det S_{j}\right|^{s_{1}}\right)_{\!\!\!j\in J\setminus J_{B}}\right\Vert _{\ell^{q_{2}\cdot\left(q_{1}/q_{2}\right)'}}\text{ with }s_{1}:=\left(\smash{\frac{1}{\LowerExpo{p_{2}}}}-\frac{1}{q_{1}}\right)_{+},\\
\quad\quad\quad K^{\left(2\right)} & :=\left\Vert \left(\frac{v_{j_{i}}}{w_{i}}\cdot\left|\det T_{i}\right|^{s_{2}}\cdot\left|\det S_{j_{i}}\right|^{p_{1}^{-1}-p_{2}^{-1}-s_{2}}\right)_{i\in I^{\left(B\right)}}\right\Vert _{\ell^{q_{2}\cdot\left(q_{1}/q_{2}\right)'}}\text{ with }s_{2}:=\left(\frac{1}{q_{2}}-\smash{\frac{1}{\SignedUpperExpo{p_{1}}}}\right)_{+},
\end{align*}
where $I^{\left(B\right)}:=\left\{ i\in I\with J_{B}\cap J_{i}\neq\emptyset\right\} $.

Precisely, given $p_{1}\leq p_{2}$, we have $\vertiii{\theta}\asymp\vertiii{\iota}\asymp K^{\left(1\right)}+K^{\left(2\right)}$,
where the implied constant only depends on $d,p_{1},p_{2},q_{1},q_{2}$
and on
\begin{align*}
\CalQ,\CalP,\varepsilon_{\CalQ},\varepsilon_{\CalP},k\left(\smash{\CalQ_{I_{A}}},\CalP\right),k\left(\smash{\CalP_{J_{B}}},\CalQ\right),C_{\CalQ,\Phi,p_{1}},C_{\CalP,\Psi,p_{1}},\\
C_{{\rm mod}}\left(\smash{\CalQ_{I_{A}}},\CalP\right),C_{{\rm mod}}\left(\smash{\CalP_{J_{B}}},\CalQ\right),C_{w,\CalQ},C_{v,\CalP},C_{w|_{I_{A}},\CalQ,\CalP},C_{v|_{J_{B}},\CalP,\CalQ}. & \qedhere
\end{align*}

\end{enumerate}
\end{thm}

\begin{proof}
For the whole proof, let $Y:=\ell_{w}^{q_{1}}\left(I\right)$ and
$Z:=\ell_{v}^{q_{2}}\left(J\right)$. We will prove each statement
individually.

\medskip{}

Ad (1): This is an immediate consequence of Corollary~\ref{cor:MixedSubordinateness},
see in particular the last statement of that corollary \emph{and the
ensuing remark}. Also note that the triangle constant of $Z=\ell_{v}^{q_{2}}\left(J\right)$
only depends on $q_{2}$ and that Lemma~\ref{lem:ModeratelyWeightedSpacesAreRegular}
can be used to estimate $\vertiii{\Gamma_{\CalQ}}_{Y\to Y}$ and $\vertiii{\Gamma_{\CalP}}_{Z\to Z}$
only in terms of $q_{1},q_{2},\CalQ,\CalP,C_{w,\CalQ}$ and $C_{v,\CalP}$.

\medskip{}

Ad (2): First, note that we can apply Lemma~\ref{lem:SimpleNecessaryCondition}
with $K:=\CalO\cap\CalO'$. Thus, choose some $\xi\in\CalO\cap\CalO'$,
which is possible since $\CalO\cap\CalO'\neq\emptyset$ by assumption.
Since $\CalQ$ and $\CalP$ are open covers of $\CalO$ and $\CalO'$,
respectively, there are $i\in I$ and $j\in J$ with $\xi\in Q_{i}\cap P_{j}\cap K=Q_{i}^{\circ}\cap P_{j}^{\circ}\cap K^{\circ}$.
Furthermore, $\delta_{i}\in Y=\ell_{w}^{q_{1}}\left(I\right)$, so
that Lemma~\ref{lem:SimpleNecessaryCondition} yields $p_{1}\leq p_{2}$,
by boundedness of $\theta$.

Next, we apply Theorem~\ref{thm:BurnerNecessaryConditionCoarseInFine},
with $J_{0}:=J_{B}$. Note that $\CalP_{J_{0}}$ is almost subordinate
to $\CalQ$ by assumption. In particular, $K:=\bigcup_{j\in J_{0}}P_{j}$
satisfies $K\subset\CalO\cap\CalO'$, so that the map $\iota:=\theta|_{\CalD_{K}}$,
with $\CalD_{K}=\left\{ f\in\TestFunctionSpace{\R^{\dimension}}\with\supp f\subset K\right\} \subset\TestFunctionSpace{\CalO\cap\CalO'}\subset\FourierDecompSp{\CalQ}{p_{1}}{\ell_{w}^{q_{1}}}$,
i.e.\@ with $\CalD_{K}=\CalD_{K}^{\CalQ,p_{1},\ell_{w}^{q_{1}}\left(I\right)}$,
satisfies the assumptions of Theorem~\ref{thm:BurnerNecessaryConditionCoarseInFine}.
Since $Z=\ell_{v}^{q_{2}}\left(J\right)$ satisfies the Fatou property,
this implies that the embedding
\[
\eta:\ell_{w}^{q_{1}}\left(\left[\ell^{p_{1}}\left(J_{i}\cap J_{B}\right)\right]_{i\in I}\right)\hookrightarrow\left(\ell_{v}^{q_{2}}\left(J_{B}\right)\right)_{\left|\det S_{j}\right|^{p_{1}^{-1}-p_{2}^{-1}}}=\ell_{\bigl[\left|\det S_{j}\right|^{p_{1}^{-1}-p_{2}^{-1}}\cdot v_{j}\bigr]_{j}}^{q_{2}}\left(J_{B}\right)
\]
is bounded with $\vertiii{\eta}\lesssim\vertiii{\theta|_{\CalD_{K}}}\leq\vertiii{\theta}$,
where the implied constant only depends on
\[
\dimension,p_{1},p_{2},q_{2},k\left(\smash{\CalP_{J_{B}}},\CalQ\right),\CalQ,\CalP,\varepsilon_{\CalP},C_{\CalQ,\Phi,p_{1}},C_{v,\CalP}\:,
\]
since the triangle constant for $Z=\ell_{v}^{q_{2}}\left(J\right)$
only depends on $q_{2}$ and since Lemma~\ref{lem:ModeratelyWeightedSpacesAreRegular}
allows us to estimate $\vertiii{\Gamma_{\CalP}}_{Z\to Z}$ only in
terms of $q_{2},\CalP,C_{v,\CalP}$.

Finally, an application of Corollary~\ref{cor:EmbeddingCoarseIntoFineSimplification}
(with $r=p_{1}$, $u\equiv1$, $u^{\left(1\right)}\equiv1$ and $u^{\left(2\right)}\equiv1$)
yields $\vertiii{\eta}\asymp K_{p_{1}}^{\left(2,1\right)}$, where
the implied constant only depends on $p_{1},q_{1},q_{2},\CalQ,k\left(\smash{\CalP_{J_{B}}},\CalQ\right),C_{w,\CalQ}$.

To prove the statement involving $K_{s}^{\left(1\right)}$, we invoke
Theorem~\ref{thm:BurnerNecessaryConditionFineInCoarse}, with $I_{0}=I_{A}$.
To see that it is applicable, note that we have $K:=\bigcup_{i\in I_{0}}Q_{i}\subset\CalO\cap\CalO'$,
since $\CalQ_{I_{A}}$ is almost subordinate to $\CalP$. As above,
this implies that the map $\iota:=\theta|_{\CalD_{K}}$ satisfies
the assumptions of Theorem~\ref{thm:BurnerNecessaryConditionFineInCoarse}.
Since $Z=\ell_{v}^{q_{2}}\left(J\right)$ satisfies the Fatou property,
this yields boundedness of
\[
\eta:\ell_{\bigl(\left|\det T_{i}\right|^{p_{2}^{-1}-p_{1}^{-1}}\cdot w_{i}\bigr)_{i}}^{q_{1}}\left(I_{A}\right)=\left(Y|_{I_{0}}\right)_{\left|\det T_{i}\right|^{p_{2}^{-1}-p_{1}^{-1}}}\hookrightarrow Z\bigl(\left[\ell^{s}\left(I_{0}\cap I_{j}\right)\right]_{j\in J}\bigr)=\ell_{v}^{q_{2}}\bigl(\left[\ell^{s}\left(I_{A}\cap I_{j}\right)\right]_{j\in J}\bigr).
\]
Precisely, we get $\vertiii{\eta}\lesssim\vertiii{\theta|_{\CalD_{K}}}\leq\vertiii{\theta}$,
where the implied constant only depends on
\[
\dimension,p_{1},p_{2},q_{1},q_{2},\CalQ,\CalP,\varepsilon_{\CalQ},k\left(\smash{\CalQ_{I_{A}}},\CalP\right),C_{w,\CalQ},C_{v,\CalP},C_{\CalQ,\Phi,p_{1}},C_{\CalP,\Psi,p_{2}},
\]
by the usual arguments involving Lemma~\ref{lem:ModeratelyWeightedSpacesAreRegular}.
But by Corollary~\ref{cor:EmbeddingFineInCoarseSimplification} (with
$J_{0}=J$, $r=s$, $u\equiv1$ and $\bigl(\left|\det T_{i}\right|^{p_{2}^{-1}-p_{1}^{-1}}\cdot w_{i}\bigr)_{i}$
instead of $w$), we have $K_{s}^{\left(1\right)}\lesssim\vertiii{\eta}$,
where the implied constant only depends on $p_{2},q_{1},q_{2},\CalP,k\left(\smash{\CalQ_{I_{A}}},\CalP\right),C_{v,\CalP}$.

\medskip{}

Ad (3): The proof is similar to the previous case, but using Theorems
\ref{thm:KhinchinNecessaryCoarseInFine} and \ref{thm:KhinchinNecessaryFineInCoarse}
(and the ensuing remarks) instead of Theorems \ref{thm:BurnerNecessaryConditionCoarseInFine}
and \ref{thm:BurnerNecessaryConditionFineInCoarse}, respectively.

\medskip{}

Ad (4): Let us first assume that $\theta$ is bounded. As seen above,
this yields $p_{1}\leq p_{2}$. We want to show that also $K^{\left(1\right)}\lesssim\vertiii{\theta}$
and $K^{\left(2\right)}\lesssim\vertiii{\theta}$.

For the first part, we invoke Theorem~\ref{thm:NecessaryConditionForModerateCoveringFineInCoarse},
with $J_{0}:=J\setminus J_{B}$. Let us verify the prerequisites of
that theorem. First of all, we need to verify that $\CalQ_{I_{0}}$
is almost subordinate to and relatively moderate with respect to $\CalP$,
where
\[
I_{0}=\left\{ i\in I\with J_{i}\cap J_{0}\neq\emptyset\right\} .
\]
But for $i\in I_{0}$, there is some $j\in J_{0}\cap J_{i}$, which
yields existence of some $\xi\in Q_{i}\cap P_{j}\subset\CalO\cap\CalO'=A\cup B$.
In case of $\xi\in B$, we would have $j\in J_{B}$, in contradiction
to $j\in J_{0}$. Hence, $\xi\in A$, so that we get $i\in I_{A}$,
since $\xi\in A\cap Q_{i}$. Thus, we have shown $I_{0}\subset I_{A}$.

Since $\CalQ_{I_{A}}$ is almost subordinate to and relatively moderate
with respect to $\CalP$, so is $\CalQ_{I_{0}}$, with $k\left(\smash{\CalQ_{I_{0}}},\CalP\right)\leq k\left(\smash{\CalQ_{I_{A}}},\CalP\right)$
and $C_{{\rm mod}}\left(\smash{\CalQ_{I_{0}}},\CalP\right)\leq C_{{\rm mod}}\left(\smash{\CalQ_{I_{A}}},\CalP\right)$.
Furthermore, since $w|_{I_{A}}$ is relatively $\CalP$-moderate,
so is $w|_{I_{0}}$, with $C_{w|_{I_{0}},\CalQ,\CalP}\leq C_{w|_{I_{A}},\CalQ,\CalP}$.
Finally, since $\CalP$ is tight and since we have $\CalO=\CalO'$,
assumption (\ref{enu:NecessaryConditionModerateCoveringFineInCoarseAlmostTightness})
of Theorem~\ref{thm:NecessaryConditionForModerateCoveringFineInCoarse}
is also satisfied (with $\varepsilon=\varepsilon_{\CalP}$). All in
all, we have verified assumptions (\ref{enu:NecessaryConditionModerateCoveringFineInCoarseAlmostTightness})–(\ref{enu:NecessaryConditionModerateCoveringFineInCoarseWeightModerate})
of that theorem. Finally, note that we have $C_{\CalP,\Psi,p_{2}}\lesssim C_{\CalP,\Psi,p_{1}}$,
where the implied constant only depends on $\dimension,p_{1},\CalP$;
see Corollary~\ref{cor:LpBAPUsAreAlsoLqBAPUsForLargerq}.

Thus, it remains to verify boundedness of the map $\iota$ as defined
in Theorem~\ref{thm:NecessaryConditionForModerateCoveringFineInCoarse}.
But since the set $K$ defined in Theorem~\ref{thm:NecessaryConditionForModerateCoveringFineInCoarse}
satisfies $K\subset\CalO\cap\CalO'$, we see that this map $\iota$
satisfies $\iota=\theta|_{\CalD_{K}}$ and hence $\vertiii{\iota}\leq\vertiii{\theta}<\infty$.
Hence, Theorem~\ref{thm:NecessaryConditionForModerateCoveringFineInCoarse}
yields a constant $C>0$, depending only on
\[
\dimension,p_{1},p_{2},q_{1},q_{2},k\left(\smash{\CalQ_{I_{A}}},\CalP\right),C_{{\rm mod}}\left(\smash{\CalQ_{I_{A}}},\CalP\right),\CalQ,\CalP,\varepsilon_{\CalQ},\varepsilon_{\CalP},C_{w,\CalQ},C_{v,\CalP},C_{w|_{I_{A}},\CalQ,\CalP},C_{\CalQ,\Phi,p_{1}},C_{\CalP,\Psi,p_{1}},
\]
which satisfies $K^{\left(1\right)}\leq C\cdot\vertiii{\iota}\leq C\cdot\vertiii{\theta}<\infty$.

To show $K^{\left(2\right)}\lesssim\vertiii{\theta}$, we invoke Theorem~\ref{thm:NecessaryConditionForModerateCoveringCoarseInFine},
with $J_{0}:=J_{B}$. Let us first verify assumptions (\ref{enu:NecessaryConditionModerateCoarseInFineCoveringModerate})–(\ref{enu:NecessaryConditionModerateCoarseInFineReverseSubordinateness})
of that theorem. The first two of these are direct consequences of
our assumptions and our choice of $J_{0}=J_{B}$. For the last one,
we need to find $C_{0}>0$ and $r\in\N_{0}$ satisfying
\begin{equation}
\lambda\left(Q_{i}\right)\leq C_{0}\cdot\lambda\left(\,\smash{\bigcup_{j\in J_{0}\cap J_{i}}}\,\vphantom{\bigcup}P_{j}^{r\ast}\,\right)\vphantom{\bigcup_{j\in J_{0}\cap J_{i}}}\qquad\forall\,i\in I_{0}:=\left\{ i\in I\with J_{0}\cap J_{i}\neq\emptyset\right\} =I^{\left(B\right)}\,.\label{eq:SummaryMixedSubordinatenessSpecialMeasureEstimate}
\end{equation}
To this end, set $k:=k\left(\smash{\CalQ_{I_{A}}},\CalP\right)$ and
let $i\in I^{\left(B\right)}$ be arbitrary. We distinguish two cases:

\begin{casenv}
\item We have $i\in I_{A}$. Because of $i\in I_{0}=I^{\left(B\right)}$,
there is some $\ell\in J_{0}\cap J_{i}$, i.e.\@ with $Q_{i}\cap P_{\ell}\neq\emptyset$.
Since $i\in I_{A}$ and because $\CalQ_{I_{A}}$ is almost subordinate
to $\CalP$, Lemma~\ref{lem:SubordinatenessImpliesWeakSubordination}
yields the inclusion $Q_{i}\subset P_{\ell}^{\left(2k+2\right)\ast}\subset\bigcup_{j\in J_{0}\cap J_{i}}P_{j}^{\left(2k+2\right)\ast}$.
\item We have $i\notin I_{A}$, i.e.\@ $Q_{i}\cap A=\emptyset$. But we
have $Q_{i}\subset\CalO=\CalO\cap\CalO'=A\cup B$, so that we get
$Q_{i}\subset B$. Now, let $\xi\in Q_{i}\subset B\subset\CalO\cap\CalO'=\CalO'$
be arbitrary. Then there is some $\ell\in J$ with $\xi\in P_{\ell}$.
Since $\xi\in B$, we get $\xi\in P_{\ell}\cap B\neq\emptyset$ and
hence $\ell\in J_{B}=J_{0}$, i.e.\@ $\ell\in J_{0}\cap J_{i}$,
which implies $\xi\in P_{\ell}\subset\bigcup_{j\in J_{0}\cap J_{i}}P_{j}\subset\bigcup_{j\in J_{0}\cap J_{i}}P_{j}^{\left(2k+2\right)\ast}$.
Since $\xi\in Q_{i}$ was arbitrary, we get $Q_{i}\subset\bigcup_{j\in J_{0}\cap J_{i}}P_{j}^{\left(2k+2\right)\ast}$.
\end{casenv}
Together, the two cases easily show that we can choose $r=2k+2$ and
$C_{0}=1$.

We have thus verified all prerequisites formulated in Theorem~\ref{thm:NecessaryConditionForModerateCoveringCoarseInFine}
itself; but note that the prerequisites of the theorem also include
those of Lemma~\ref{lem:NecessaryConjugateCoarseInFineWithoutModerateness},
which we verify now. First of all, $\CalP_{J_{0}}=\CalP_{J_{B}}$
is almost subordinate to $\CalQ$, so that it remains to establish
the existence of a bounded linear map 
\[
\iota:\left(\CalD_{K}^{\CalQ,p_{1},\ell_{w}^{q_{1}}\left(I\right)},\left\Vert \mybullet\right\Vert _{\FourierDecompSp{\CalQ}{p_{1}}{\ell_{w}^{q_{1}}}}\right)\to\FourierDecompSp{\CalP}{p_{2}}{\ell_{v}^{q_{2}}}
\]
with $\left\langle \iota f,\,\varphi\right\rangle _{\CalD'}=\left\langle f,\,\varphi\right\rangle _{\CalD'}$
for all $\varphi\in\TestFunctionSpace{\CalO\cap\CalO'}$ and all $f\in\CalD_{K}^{\CalQ,p_{1},\ell_{w}^{q_{1}}\left(I\right)}$,
where 
\[
K=\bigcup_{i\in I_{0}}\overline{Q_{i}^{\left(2k+3\right)\ast}}\subset\CalO=\CalO\cap\CalO'.
\]
But given our present assumptions, it is easy to see that $\iota=\theta|_{\CalD_{K}^{\CalQ,p_{1},\ell_{w}^{q_{1}}\left(I\right)}}$
satisfies these properties.

All in all, we can thus apply Theorem~\ref{thm:NecessaryConditionForModerateCoveringCoarseInFine},
so that we get a constant $C>0$, depending only on
\[
\dimension,p_{1},p_{2},q_{1},q_{2},\CalQ,\varepsilon_{\CalQ},\CalP,\varepsilon_{\CalP},k\left(\smash{\CalQ_{I_{A}}},\CalP\right),k\left(\smash{\CalP_{J_{B}}},\CalQ\right),C_{{\rm mod}}\left(\smash{\CalP_{J_{B}}},\CalQ\right),C_{w,\CalQ},C_{v,\CalP},C_{v|_{J_{B}},\CalP,\CalQ},C_{\CalQ,\Phi,p_{1}},
\]
which satisfies $K^{\left(2\right)}\leq C\cdot\vertiii{\iota}\leq C\cdot\vertiii{\theta}<\infty$.

\medskip{}

Next, if $\iota$ is bounded, the properties of $\iota$ from part~(\ref{enu:SummaryMixedSubordinatenessSufficientNonModerate})
(together with $\CalO=\CalO'$) imply $\iota f=f$ for all $f\in\FourierDecompSp{\CalQ}{p_{1}}{\ell_{w}^{q_{1}}}$,
so that the map $\theta$ from part~(\ref{enu:SummaryMixedSubordinatenessNecessaryNonModerate})
satisfies $\theta=\iota|_{\TestFunctionSpace{\CalO}}$. In view of
what we just showed, we thus get $p_{1}\leq p_{2}$, $K^{\left(1\right)}\lesssim\vertiii{\theta}\leq\vertiii{\iota}$
and $K^{\left(2\right)}\lesssim\vertiii{\theta}\leq\vertiii{\iota}$.

\medskip{}

It thus remains to show that $\iota$ is bounded if $p_{1}\leq p_{2}$,
$K^{\left(1\right)}<\infty$ and $K^{\left(2\right)}<\infty$. To
this end, we will show that the assumptions of Corollary~\ref{cor:MixedSubordinateness}
are satisfied. All of these assumptions are trivially included in
those of the present theorem, with the exception of finiteness of
the right-hand sides of equations (\ref{eq:MixedSubordinatenessSimplificationBeta1Norm})
and (\ref{eq:MixedSubordinatenessSimplificationBeta2Norm}), which
we will verify now.

First, note that the second part of Corollary~\ref{cor:EmbeddingFineInCoarseSimplification}
(with $J_{0}:=J\setminus J_{B}$, with $I_{0}:=L:=\bigcup_{j\in J\setminus J_{B}}I_{j}$
and with $r=\LowerExpo{p_{2}}$ and $u_{i}:=\left|\det T_{i}\right|^{p_{1}^{-1}-p_{2}^{-1}}$)
yields
\begin{align}
\text{r.h.s. of eq. }\eqref{eq:MixedSubordinatenessSimplificationBeta1Norm} & =\left\Vert \left(v_{j}\cdot\left\Vert \left(\left|\det T_{i}\right|^{p_{1}^{-1}-p_{2}^{-1}}/w_{i}\right)_{i\in I_{j}}\right\Vert _{\ell^{\LowerExpo{p_{2}}\cdot\left(q_{1}/\LowerExpo{p_{2}}\right)'}}\right)_{j\in J\setminus J_{B}}\right\Vert _{\ell^{q_{2}\cdot\left(q_{1}/q_{2}\right)'}}\nonumber \\
\left({\scriptstyle \text{since }I_{j}\subset L=I_{0}\text{ for }j\in J_{0}=J\setminus J_{B}}\right) & =\left\Vert \left(v_{j}\cdot\left\Vert \left(u_{i}/w_{i}\right)_{i\in I_{0}\cap I_{j}}\right\Vert _{\ell^{\LowerExpo{p_{2}}\cdot\left(q_{1}/\LowerExpo{p_{2}}\right)'}}\right)_{j\in J_{0}}\right\Vert _{\ell^{q_{2}\cdot\left(q_{1}/q_{2}\right)'}}\nonumber \\
 & \asymp\left\Vert \left(v_{j}\cdot\frac{u_{i_{j}}}{w_{i_{j}}}\cdot\left[\left|\det S_{j}\right|\bigg/\left|\det T_{i_{j}}\right|\right]^{\left(\smash{\frac{1}{\LowerExpo{p_{2}}}}-\frac{1}{q_{1}}\right)_{+}}\right)_{j\in J^{\left(0\right)}}\right\Vert _{\ell^{q_{2}\cdot\left(q_{1}/q_{2}\right)'}}\nonumber \\
 & \overset{\left(\ast\right)}{=}K^{\left(1\right)}<\infty,\label{eq:SummaryMixedSubordinatenessModerateSufficiency}
\end{align}
with $J^{\left(0\right)}:=\left\{ j\in J_{0}\with I_{0}\cap I_{j}\neq\emptyset\right\} $
and with an implied constant only depending on 
\[
\dimension,p_{1},p_{2},q_{1},\CalQ,\CalP,\varepsilon_{\CalQ},\varepsilon_{\CalP},k\left(\smash{\CalQ_{I_{A}}},\CalP\right),C_{{\rm mod}}\left(\smash{\CalQ_{I_{A}}},\CalP\right),C_{w|_{I_{A}},\CalQ,\CalP}.
\]
The prerequisites for the application of Corollary~\ref{cor:EmbeddingFineInCoarseSimplification}
from above are not hard to verify: As seen in Corollary~\ref{cor:MixedSubordinateness},
we have $I_{0}=L\subset I_{A}$, so that $\CalQ_{I_{0}}$ is almost
subordinate to and relatively moderate with respect to $\CalP$, since
$\CalQ_{I_{A}}$ is. We also get $C_{{\rm mod}}\left(\smash{\CalQ_{I_{0}}},\CalP\right)\leq C_{{\rm mod}}\left(\smash{\CalQ_{I_{A}}},\CalP\right)$
and $k\left(\smash{\CalQ_{I_{0}}},\CalP\right)\leq k\left(\smash{\CalQ_{I_{A}}},\CalP\right)$.
Hence, the only nontrivial point is to establish existence of $C_{0}>0$
and $s\in\N_{0}$ satisfying
\[
\lambda\left(P_{j}\right)\leq C_{0}\cdot\lambda\left(\,\smash{\bigcup_{i\in I_{0}\cap I_{j}}}\,\vphantom{\bigcup}Q_{i}^{s\ast}\,\right)\vphantom{\bigcup_{i\in I_{0}\cap I_{j}}}\qquad\forall\,j\in J_{0}\,.
\]
But since we have $I_{0}\cap I_{j}=L\cap I_{j}=I_{j}$ for all $j\in J_{0}=J\setminus J_{B}$,
and because of $P_{j}\subset\CalO'=\CalO=\bigcup_{i\in I}Q_{i}$,
we easily see $P_{j}\subset\bigcup_{i\in I_{j}}Q_{i}=\bigcup_{i\in I_{0}\cap I_{j}}Q_{i}^{0\ast}$
for all $j\in J_{0}$, so that we can choose $C_{0}=1$ and $s=0$.

Finally, to justify the step marked with $\left(\ast\right)$ in estimate~(\ref{eq:SummaryMixedSubordinatenessModerateSufficiency}),
note that we have $u_{i_{j}}=\left|\det\smash{T_{i_{j}}}\right|^{p_{1}^{-1}-p_{2}^{-1}}$
and that we have $I_{0}\cap I_{j}=I_{j}\neq\emptyset$ for all $j\in J_{0}$,
since $\emptyset\neq P_{j}\subset\bigcup_{i\in I_{0}\cap I_{j}}Q_{i}$,
as we just saw. Hence, $J^{\left(0\right)}=J_{0}=J\setminus J_{B}$.

We have thus established finiteness of the right-hand side of equation~(\ref{eq:MixedSubordinatenessSimplificationBeta1Norm}).

\medskip{}

To establish finiteness of the right-hand side of equation~(\ref{eq:MixedSubordinatenessSimplificationBeta2Norm}),
we apply the second part of Corollary~\ref{cor:EmbeddingCoarseIntoFineSimplification},
with $r=\UpperExpo{p_{1}}$, $J_{0}=J_{B}$ and
\[
u_{i}^{\left(1\right)}:=\begin{cases}
\left|\det T_{i}\right|^{1-p_{1}^{-1}}, & \text{if }p_{1}<1,\\
1, & \text{if }p_{1}\geq1,
\end{cases}\qquad\text{ as well as }\qquad u_{j}^{\left(2\right)}:=\begin{cases}
\left|\det S_{j}\right|^{p_{2}^{-1}-1}, & \text{if }p_{1}<1,\\
\left|\det S_{j}\right|^{p_{2}^{-1}-p_{1}^{-1}}, & \text{if }p_{1}\geq1.
\end{cases}
\]
All of the prerequisites of Corollary~\ref{cor:EmbeddingCoarseIntoFineSimplification}
are easily seen to be fulfilled, with the possible exception of
\[
\lambda\left(Q_{i}\right)\leq C_{0}\cdot\lambda\left(\,\smash{\bigcup_{j\in J_{0}\cap J_{i}}}\,\vphantom{\bigcup}P_{j}^{s\ast}\,\right)\vphantom{\bigcup_{j\in J_{0}\cap J_{i}}}\qquad\forall\,i\in I^{\left(0\right)}=\left\{ i\in I\with J_{0}\cap J_{i}\neq\emptyset\right\} =I^{\left(B\right)}.
\]
But following equation~(\ref{eq:SummaryMixedSubordinatenessSpecialMeasureEstimate}),
we verified exactly this estimate, for $C_{0}=1$ and $s=2k+2$, for
$k=k\left(\smash{\CalQ_{I_{A}}},\CalP\right)$. All in all, Corollary~\ref{cor:EmbeddingCoarseIntoFineSimplification}
thus yields
\begin{align*}
\text{r.h.s. of eq. }\eqref{eq:MixedSubordinatenessSimplificationBeta2Norm} & =\left\Vert \left(w_{i}^{-1}\cdot\left\Vert \left(v_{j}/u_{i,j}\right)_{j\in J_{0}\cap J_{i}}\right\Vert _{\ell^{q_{2}\cdot\left(\UpperExpo{p_{1}}/q_{2}\right)'}}\right)_{i\in I}\right\Vert _{\ell^{q_{2}\cdot\left(q_{1}/q_{2}\right)'}}\\
 & \asymp\left\Vert \left(w_{i}^{-1}\cdot\frac{v_{j_{i}}}{u_{i,j_{i}}}\cdot\left[\left|\det T_{i}\right|\bigg/\left|\det S_{j_{i}}\right|\right]^{\left(\frac{1}{q_{2}}-\smash{\frac{1}{\UpperExpo{p_{1}}}}\right)_{+}}\right)_{i\in I^{\left(B\right)}}\right\Vert _{\ell^{q_{2}\cdot\left(q_{1}/q_{2}\right)'}}\\
\left({\scriptstyle \text{see Remark }\ref{rem:SufficientCoarseIntoFineSimplification}}\right) & =K^{\left(2\right)}<\infty,
\end{align*}
where the implied constant only depends on
\[
\dimension,p_{1},p_{2},q_{2},\CalQ,\CalP,\varepsilon_{\CalQ},\varepsilon_{\CalP},k\left(\smash{\CalQ_{I_{A}}},\CalP\right),k\left(\smash{\CalP_{J_{B}}},\CalQ\right),C_{{\rm mod}}\left(\smash{\CalP_{J_{B}}},\CalQ\right),C_{v|_{J_{B}},\CalP,\CalQ}.
\]

\medskip{}

Finally, we can invoke Corollary~\ref{cor:MixedSubordinateness}
to conclude that $\iota$ is bounded, with
\[
\vertiii{\iota}\leq C\cdot\left[\left(\text{r.h.s. of eq. }\eqref{eq:MixedSubordinatenessSimplificationBeta1Norm}\right)+\left(\text{r.h.s. of eq. }\eqref{eq:MixedSubordinatenessSimplificationBeta2Norm}\right)\right],
\]
for some constant
\[
C=C\left(\dimension,p_{1},p_{2},q_{1},q_{2},\CalQ,\CalP,k\left(\smash{\CalP_{J_{B}}},\CalQ\right),k\left(\smash{\CalQ_{I_{A}}},\CalP\right),C_{\CalQ,\Phi,p_{1}},C_{\CalP,\Psi,p_{1}},C_{w,\CalQ},C_{v,\CalP}\right),
\]
thereby completing the proof.
\end{proof}

\subsection{Embeddings between decomposition spaces: A user's guide}

We close this section with a user's guide which describes a convenient
workflow for deciding the existence of an embedding between decomposition
spaces. Concrete applications of (essentially) this workflow are given
in Section~\ref{sec:Applications}. The typical approach is as follows:
\begin{enumerate}[leftmargin=0.7cm]
\item Determine the coverings $\CalQ,\CalP$ (and the remaining parameters
$p_{1},p_{2},q_{1},q_{2},w,v$) which define the decomposition spaces
for which an embedding $\FourierDecompSp{\CalQ}{p_{1}}{\ell_{w}^{q_{1}}}\hookrightarrow\FourierDecompSp{\CalP}{p_{2}}{\ell_{v}^{q_{2}}}$
is desired.
\item If $p_{1}>p_{2}$, the embedding does not exist. Thus, assume in the
following that $p_{1}\leq p_{2}$.
\item Determine whether $\CalP$ is \textbf{almost subordinate} to $\CalQ$
or vice versa (or neither). Probably the most convenient way to do
this is as follows:

\begin{enumerate}[leftmargin=0.55cm]
\item Determine upper and lower bounds (in terms of set-inclusion) for
the \textbf{intersection sets}
\[
\qquad\qquad I_{j}:=\left\{ i\in I\with Q_{i}\cap P_{j}\neq\emptyset\right\} \qquad\text{ and }\qquad J_{i}:=\left\{ j\in J\with P_{j}\cap Q_{i}\neq\emptyset\right\} .
\]
This amounts to finding conditions on $i,j$ which ensure/prevent
$Q_{i}\cap P_{j}\neq\emptyset$. This step is the (only) one in which
the \emph{geometry} of the two coverings plays a significant role.
\item If we have $\sup_{j\in J}\left|I_{j}\right|<\infty$, then $\CalP$
is \textbf{weakly subordinate} to $\CalQ$. If $\CalP$ consists of
path-connected sets and if $\CalQ$ is an open covering, then this
implies that $\CalP$ is almost subordinate to $\CalQ$, see Lemma~\ref{lem:WeaksubordinatenessImpliesSubordinatenessIfConnected}.

The same holds if $\CalQ$ and $\CalP$ are exchanged everywhere and
if the condition $\sup_{j\in J}\left|I_{j}\right|<\infty$ is replaced
by $\sup_{i\in I}\left|J_{i}\right|<\infty$.
\end{enumerate}
\item A quick check: If $p_{1}=p_{2}$, try to prove/disprove that $v_{j}\lesssim w_{i}$
if $Q_{i}\cap P_{j}\neq\emptyset$. If this fails, the embedding does
\emph{not} exist.
\item Determine whether $\CalP$ and $v$ are relatively moderate with respect
to $\CalQ$ (or vice versa).

This amounts to showing that any two sets $P_{j},P_{\ell}$ with $j,\ell\in J_{i}$
have essentially the same measure, $\lambda\left(P_{j}\right)\asymp\lambda\left(P_{\ell}\right)$
and that $v_{j}\asymp v_{\ell}$ for $j,\ell\in J_{i}$. Of course,
this depends on having good upper bounds for the sets $\left(J_{i}\right)_{i\in I}$.
Usually, this step will yield sequences $\gamma_{i}^{\left(1\right)},\gamma_{i}^{\left(2\right)}$
satisfying
\begin{equation}
\lambda\left(P_{j}\right)\asymp\gamma_{i}^{\left(1\right)}\qquad\text{ and }\qquad v_{j}\asymp\gamma_{i}^{\left(2\right)}\qquad\text{ for }j\in J_{i}\text{ and }i\in I.\label{eq:IntroductionRelativeModeratenessProxyWeight}
\end{equation}

\item Choose the right criterion to use:

\begin{enumerate}[leftmargin=0.55cm]
\item If $\CalQ$ is almost subordinate to $\CalP$: Consider Theorem~\ref{thm:SummaryFineIntoCoarse}.
Specifically:

\begin{enumerate}[leftmargin=0.6cm]
\item If $\CalQ$ and $w$ are relatively $\CalP$-moderate: Use part~(\ref{enu:SummaryFineInCoarseModerate})
of Theorem~\ref{thm:SummaryFineIntoCoarse}. This requires only to
check whether a \emph{single norm} is finite, which can usually be
easily done (using the equivalents of the weights $\gamma^{\left(1\right)}$
and $\gamma^{\left(2\right)}$ from equation~(\ref{eq:IntroductionRelativeModeratenessProxyWeight})).
Furthermore, this yields a \emph{complete characterization} of the
existence of the embedding, i.e.\@ a \emph{yes/no answer}.
\item If $\CalQ$ or $w$ are \emph{not} relatively $\CalP$-moderate: Consider
the expression
\[
\qquad\qquad\qquad K\left(r\right):=\left\Vert \left(v_{j}\cdot\left\Vert \left(\left|\det T_{i}\right|^{p_{1}^{-1}-p_{2}^{-1}}/\,w_{i}\right)_{i\in I_{j}}\right\Vert _{\ell^{r\cdot\left(q_{1}/r\right)'}}\right)_{j\in J}\right\Vert _{\ell^{q_{2}\cdot\left(q_{1}/q_{2}\right)'}}\in\left[0,\infty\right]
\]
from Theorem~\ref{thm:SummaryFineIntoCoarse} (equation~(\ref{eq:SummaryFineIntoCoarseCriticalQuantitiy}))
for general $r\in\left(0,\infty\right]$. More precisely, determine
for which values of $r$ the expression $K\left(r\right)$ is finite.
For this, it could be helpful to note $\frac{1}{r\cdot\left(q_{1}/r\right)'}=\left(r^{-1}-q_{1}^{-1}\right)_{+}$.
In particular, $r\mapsto K\left(r\right)$ is non-increasing.\vspace{0.1cm}\\
The following hold:\vspace{0.1cm}
\begin{enumerate}
\item If $K\left(\smash{\LowerExpo{p_{2}}}\right)<\infty$, the embedding
exists.\vspace{0.1cm}
\item If $K\left(p_{2}\right)=\infty$ or if $p_{2}=\infty$ and $K\left(\smash{\LowerExpo{p_{2}}}\right)=\infty$,
the embedding does \emph{not} exist.\vspace{0.1cm}
\item If $p_{1}=p_{2}$ and $K\left(2\right)=\infty$, the embedding does
\emph{not} exist.\vspace{0.1cm}
\item If none of the cases above applies, the (convenient) criteria given
in this paper are inconclusive.
\end{enumerate}
\end{enumerate}
\item If $\CalP$ is almost subordinate to $\CalQ$: Consider Theorem~\ref{thm:SummaryCoarseIntoFine}.
Specifically:

\begin{enumerate}[leftmargin=0.6cm]
\item If $\CalP$ and $v$ are relatively $\CalQ$-moderate: Use part~(\ref{enu:SummaryCoarseInFineModerate})
of Theorem~\ref{thm:SummaryCoarseIntoFine}. This requires only to
check whether a \emph{single norm} is finite, which can usually be
easily done (using the weights $\gamma^{\left(1\right)}$ and $\gamma^{\left(2\right)}$
from equation~(\ref{eq:IntroductionRelativeModeratenessProxyWeight})).
Furthermore, this yields a \emph{complete characterization} of the
existence of the embedding, i.e.\@ a \emph{yes/no answer}.
\item If $\CalP$ or $v$ are not relatively $\CalQ$-moderate: Consider
the expression
\[
\qquad\qquad\qquad K\left(r,\alpha,\beta\right):=\left\Vert \left(w_{i}^{-1}\cdot\left|\det T_{i}\right|^{\alpha}\cdot\left\Vert \left(\left|\det S_{j}\right|^{\beta}\cdot v_{j}\right)_{j\in J_{i}}\right\Vert _{\ell^{q_{2}\cdot\left(r/q_{2}\right)'}}\right)_{i\in I}\right\Vert _{\ell^{q_{2}\cdot\left(q_{1}/q_{2}\right)'}}\in\left[0,\infty\right]
\]
for general $r\in\left(0,\infty\right]$ and $\alpha,\beta\in\R$.
Then,

\begin{enumerate}
\item If $p_{1}\geq1$ and if $K\left(\UpperExpo{p_{1}},0,p_{1}^{-1}-p_{2}^{-1}\right)<\infty$,
the embedding exists.
\item If $p_{1}<1$ and if $K\left(\UpperExpo{p_{1}},p_{1}^{-1}-1,1-p_{2}^{-1}\right)<\infty$,
the embedding exists.
\item If $K\left(p_{1},0,p_{1}^{-1}-p_{2}^{-1}\right)=\infty$, the embedding
does \emph{not} exist.
\item If $p_{1}=p_{2}$ and $K\left(2,0,0\right)=\infty$, the embedding
does \emph{not} exist.
\item If none of the cases above applies, the (convenient) criteria given
in this paper are inconclusive.
\end{enumerate}
\end{enumerate}
\item If neither $\CalQ$ is almost subordinate to $\CalP$, nor vice versa:

\begin{enumerate}[leftmargin=0.6cm]
\item Try to find sets $A,B$ satisfying $A\cup B=\CalO\cap\CalO'$ and
such that $\CalQ$ is ``smaller than'' $\CalP$ ``near'' $A$
and vice versa ``near'' $B$. Precisely, such that if
\[
\qquad\qquad I_{A}:=\left\{ i\in I\with Q_{i}\cap A\neq\emptyset\right\} \qquad\text{ and }\qquad J_{B}:=\left\{ j\in J\with P_{j}\cap B\neq\emptyset\right\} ,
\]
then $\CalQ_{I_{A}}:=\left(Q_{i}\right)_{i\in I_{A}}$ is almost subordinate
to $\CalP$ and $\CalP_{J_{B}}:=\left(P_{j}\right)_{j\in J_{B}}$
is almost subordinate to $\CalQ$.

If such sets can be found, use Theorem~\ref{thm:SummaryMixedSubordinateness}.
A more detailed explanation of this is outside of the scope of this
user's guide. See the mentioned theorem for more details.

If no sets $A,B$ as above can be found, see the next two points for
further options.
\item If a sufficient criterion is desired:

\begin{enumerate}
\item Try to find a covering $\CalR$ such that $\CalQ$ and $\CalP$ are
both almost subordinate to $\CalR$ or such that $\CalR$ is almost
subordinate to both $\CalQ$ and $\CalP$. Then, use the preceding
criteria to establish a chain of embeddings
\[
\qquad\qquad\qquad\FourierDecompSp{\CalQ}{p_{1}}{\ell_{w}^{q_{1}}}\hookrightarrow\FourierDecompSp{\CalR}{p_{3}}{\ell_{u}^{q_{3}}}\hookrightarrow\FourierDecompSp{\CalP}{p_{2}}{\ell_{v}^{q_{2}}}
\]
for suitably selected parameters $p_{3},q_{3}$ and a suitable weight
$u$.
\item Try to apply the (very general) Theorem~\ref{thm:NoSubordinatenessWithConsiderationOfOverlaps}.
This theorem does not require subordinateness of any kind, but its
prerequisites are quite involved and fall outside the scope of this
user's guide.
\end{enumerate}
\item If a necessary criterion is desired: In case of $p_{1}=p_{2}$, Theorem~\ref{thm:KhinchinNecessaryCoarseInFine}
might be helpful. In case of $1\leq p_{1}\leq p_{2}$, Theorem~\ref{thm:BurnerNecessaryConditionCoarseInFine}
and Remark~\ref{rem:BurnerNecessaryConditionCoarseInFine} yield
readily applicable necessary criteria. Finally, recall that the necessary
conditions from Section~\ref{subsec:ImprovedNecessaryConditions}
only need a \emph{subfamily} of $\CalQ$ to be almost subordinate
to $\CalP$, or vice versa. Thus, even though $\CalQ$ might not be
almost subordinate to $\CalP$, one can try to find a ``large''
(infinite) subset $I_{0}\subset I$ such that $\CalQ_{I_{0}}=\left(Q_{i}\right)_{i\in I_{0}}$
is almost subordinate to $\CalP$ (or a ``large'' subset $J_{0}\subset J$
such that $\CalP_{J_{0}}$ is almost subordinate to $\CalQ$). If
such a subset can be found, the necessary criteria from Section~\ref{subsec:ImprovedNecessaryConditions}
can be readily applied.

If all these criteria fail, then there is no readily applicable criterion
in this paper which applies in this generality. However, reading the
statements and proofs of our various necessary conditions might yield
inspiration for finding a ``custom-tailored'' necessary criterion.
The interested reader might want to study the proof of equation~(\ref{eq:InhomIntoHomBesovNecessarySmallFrequencies2})
in the proof of Theorem~\ref{thm:InhomogeneousIntoHomogeneousBesov}
as an example of such a criterion.
\end{enumerate}
\end{enumerate}
\end{enumerate}

\section{Decomposition spaces as spaces of tempered distributions}

\label{sec:DecompositionSpacesAsSpacesOfTemperedDistributions}In
this section we study embeddings of decomposition spaces into the
space of tempered distributions. This question is nontrivial, since
the (Fourier-side) decomposition spaces in this paper are defined
as subspaces of $\DistributionSpace{\CalO}$ instead of $\Schwartz'\left(\R^{\dimension}\right)$.

As noted in Section~\ref{sec:DecompositionSpaces}—in particular
in Remark~\ref{rem:TemperedDistributionsAsReservoirIncomplete} and
Example~\ref{exa:BorupNielsenDecompositionSpaceIncomplete}—this
has two main reasons:
\begin{enumerate}
\item We want to allow the case $\CalO\subsetneq\R^{\dimension}$. In this
case, one would have to factor out a certain subspace of $\Schwartz'\left(\R^{\dimension}\right)$
to obtain a positive definite quasi-norm $\left\Vert \mybullet\right\Vert _{\FourierDecompSp{\CalQ}pY}$.
\item Even for $\CalO=\R^{\dimension}$, using $\Schwartz'\left(\R^{\dimension}\right)$
as the reservoir can lead to incomplete spaces $\CalD_{\Schwartz'}\left(\CalQ,L^{p},Y\right)$,
even in case of $Y=\ell_{u}^{1}$ with a $\CalQ$-moderate weight
$u$.
\end{enumerate}
The second reason can be seen as just a technical nuisance in most
cases. That is why we now provide a readily verifiable sufficient
criterion which ensures that each element $f\in\FourierDecompSp{\CalQ}pY\subset\DistributionSpace{\CalO}$
of a (Fourier-side) decomposition space has an extension to a tempered
distribution $f_{\Schwartz'}\in\Schwartz'\left(\R^{\dimension}\right)$.
In case of $\CalO=\R^{\dimension}$, this also implies that every
element $f\in\DecompSp{\CalQ}pY\subset\SpaceReservoir{\CalO}$ of
the \emph{space-side} decomposition space admits an extension $\Fourier^{-1}\left(\Fourier f\right)_{\Schwartz'}\in\Schwartz'\left(\R^{\dimension}\right)$,
since $\Schwartz'\left(\R^{\dimension}\right)$ is invariant under
the (inverse) Fourier transform and because of $\Fourier f\in\FourierDecompSp{\CalQ}pY\hookrightarrow\Schwartz'\left(\R^{\dimension}\right)$,
see Remark~\ref{rem:TemperedDistributionsAsReservoirIncomplete}.
This implies that the two spaces $\DecompSp{\CalQ}pY$ and $\CalD_{\Schwartz'}\left(\CalQ,L^{p},Y\right)$
(as defined in Remark~\ref{rem:TemperedDistributionsAsReservoirIncomplete})
coincide up to obvious identifications.

Our criterion for ensuring $\FourierDecompSp{\CalQ}pY\hookrightarrow\Schwartz'\left(\R^{\dimension}\right)$
can be most conveniently formulated using the notions of so-called
\textbf{regular partitions of unity} and \textbf{regular coverings}.
These terms were first used informally by Borup and Nielsen\cite{BorupNielsenDecomposition}
and formally introduced in \cite[Definition 2.4]{DecompositionIntoSobolev}.
\begin{defn}
\label{def:RegularCoveringRegularPartitionOfUnity}Let $\CalQ=\left(T_{i}Q_{i}'+b_{i}\right)_{i\in I}$
be a semi-structured covering of an open set $\emptyset\neq\CalO\subset\R^{\dimension}$
and let $\Phi=\left(\varphi_{i}\right)_{i\in I}$ be a smooth partition
of unity subordinate to $\CalO$. For $i\in I$, define the \textbf{normalized
version} of $\varphi_{i}$ as
\[
\varphi_{i}^{\#}:\R^{\dimension}\to\Compl,\xi\mapsto\varphi_{i}\left(T_{i}\xi+b_{i}\right).
\]
We say that $\Phi$ is a \textbf{regular partition of unity} subordinate
to $\CalQ$ if $\varphi_{i}\in\TestFunctionSpace{\CalO}$ with $\varphi_{i}\equiv0$
on $\CalO\setminus Q_{i}$ for all $i\in I$ and if additionally the
following expression (then a constant) is finite for all $\alpha\in\N_{0}^{\dimension}$:
\[
C_{\Phi,\alpha}:=\sup_{i\in I}\left\Vert \partial^{\alpha}\smash{\varphi_{i}^{\#}}\right\Vert _{\sup}\:.
\]

$\CalQ$ is called a \textbf{regular covering} of $\CalO$ if there
is a regular partition of unity $\Phi$ subordinate to $\CalQ$.
\end{defn}

Although the notion of a regular covering seems somewhat restrictive,
every ``reasonable'' covering turns out to be regular. For readers
who want to recall the notion of an \textbf{almost structured covering},
we refer to Definition~\ref{defn:DifferentTypesOfCoverings} from
above.
\begin{thm}
\label{thm:AlmostStructuredCoveringsAreRegular}(cf.\@ \cite[Corollary 2.7 and Theorem 2.8]{DecompositionIntoSobolev};
see \cite[Proposition 1]{BorupNielsenDecomposition} for a similar
statement)

Every almost structured covering is regular. In particular, every
structured covering is regular.

Finally, every regular partition of unity subordinate to $\CalQ$
is an $L^{p}$-BAPU for $\CalQ$, for all $p\in\left(0,\infty\right]$.
\end{thm}

Using the notion of regular partitions of unity, we can now state
our criterion for embeddings into the space of tempered distributions.
One can probably tweak some of the parameters to obtain similar, but
different sufficient criteria; but for our purposes, the present version
suffices.
\begin{thm}
\label{thm:DecompositionSpacesAsTemperedDistributions}Assume that
$\CalQ=\left(Q_{i}\right)_{i\in I}=\left(T_{i}Q_{i}'+b_{i}\right)_{i\in I}$
is a regular covering of the open set $\emptyset\neq\CalO\subset\R^{\dimension}$
and let $\Phi=\left(\varphi_{i}\right)_{i\in I}$ be a regular partition
of unity for $\CalQ$. Let $Y\subset\Compl^{I}$ be $\CalQ$-regular
and let $p\in\left(0,\infty\right]$.

For each $N\in\N_{0}$, define $w^{\left(N\right)}=\left(\smash{w_{i}^{\left(N\right)}}\right)_{i\in I}$
by 
\[
w_{i}^{\left(N\right)}:=\left|\det T_{i}\right|^{1/p}\cdot\max\left\{ 1,\left\Vert \smash{T_{i}^{-1}}\right\Vert ^{\dimension+1}\right\} \cdot\left[\inf_{\xi\in Q_{i}^{\ast}}\left(1+\left|\xi\right|\right)\right]^{-N}.
\]

If there is some subset $I_{0}\subset I$ and some $N\in\N_{0}$ for
which we have 
\begin{equation}
w^{\left(N\right)}\cdot\Indicator_{I_{0}}\cdot\theta\in\ell^{1}\left(I\right)\qquad\forall\,\theta=\left(\theta_{i}\right)_{i\in I}\in Y,\label{eq:EmbeddingIntoTemperedDistributionsWeightAssumption}
\end{equation}
then
\[
\Phi:\FourierDecompSp{\CalQ}pY\rightarrow\Schwartz'\left(\R^{\dimension}\right),f\mapsto\left[\,g\mapsto\smash{\sum_{i\in I_{0}}}\,\vphantom{\sum}\left\langle \varphi_{i}f,\,g\right\rangle _{\Schwartz'}\,\right]\vphantom{\sum_{i\in I_{0}}}
\]
is well-defined, and continuous with respect to the weak-$\ast$-topology
on $\Schwartz'\left(\R^{\dimension}\right)$. The series defining
$\left\langle \Phi f,\,g\right\rangle _{\Schwartz'}$ converges absolutely
for every $g\in\Schwartz\left(\R^{\dimension}\right)$.

Finally, for $I_{0}=I$, the tempered distribution $\Phi f\in\Schwartz'\left(\R^{\dimension}\right)$
is an extension of $f:\TestFunctionSpace{\CalO}\rightarrow\Compl$
to $\Schwartz\left(\R^{\dimension}\right)$ for any $f\in\FourierDecompSp{\CalQ}pY\subset\DistributionSpace{\CalO}$,
i.e., $\left\langle \Phi f,\,g\right\rangle _{\Schwartz'}=\left\langle f,\,g\right\rangle _{\CalD'}$
for all $g\in\TestFunctionSpace{\CalO}$.
\end{thm}

\begin{rem*}
In case of $Y=\ell_{u}^{q}\left(I\right)$, it is not hard to see
that the assumption~(\ref{eq:EmbeddingIntoTemperedDistributionsWeightAssumption})
regarding the weight $w^{\left(N\right)}$ is equivalent to requiring
$w^{\left(N\right)}|_{I_{0}}\in\ell_{1/u}^{q'}\left(I_{0}\right)$.
\end{rem*}
\begin{proof}
We start the proof with some preliminary observations, estimates and
definitions.

By assumption, the linear map
\[
\Lambda:Y\to\ell^{1}\left(I\right),\theta\mapsto w^{\left(N\right)}\cdot\Indicator_{I_{0}}\cdot\theta
\]
is well-defined. Furthermore, since convergence in $Y$ and in $\ell^{1}\left(I\right)$
yields convergence in each coordinate (see Lemma~\ref{lem:SolidSequenceSpaceEmbedsIntoWeightedLInfty}),
it is easy to see that $\Lambda$ has a closed graph. Finally, as
discussed in Section~\ref{subsec:Notation}, the closed graph theorem
is valid for the quasi-Banach space $Y$. Therefore, $\Lambda$ is
continuous and hence bounded, i.e., there is a constant $C_{1}>0$
with
\begin{equation}
\left\Vert \left(\smash{w_{i}^{\left(N\right)}}\cdot\theta_{i}\right)_{i\in I_{0}}\right\Vert _{\ell^{1}\left(I_{0}\right)}=\left\Vert \,w^{\left(N\right)}\cdot\Indicator_{I_{0}}\cdot\theta\,\right\Vert _{\ell^{1}\left(I\right)}\leq C_{1}\cdot\left\Vert \theta\right\Vert _{Y}\qquad\forall\,\theta\in Y\,.\label{eq:DecompositionSpaceAsTemperedDistributionsQuantitativeSequenceEstimate}
\end{equation}

For brevity, let us set
\[
S_{i}:\R^{\dimension}\to\R^{\dimension},\xi\mapsto T_{i}\xi+b_{i}\qquad\text{ and }\qquad c_{i}:=\inf_{\xi\in Q_{i}^{\ast}}\left(1+\left|\xi\right|\right)\quad\text{ for }i\in I\,.
\]
Furthermore, let
\[
\left\Vert g\right\Vert _{k,N}:=\max_{\substack{\alpha\in\N_{0}^{\dimension}\\
\left|\alpha\right|\leq N
}
}\:\sup_{x\in\R^{\dimension}}\left[\left(1+\left|x\right|\right)^{k}\cdot\left|\left(\partial^{\alpha}g\right)\left(x\right)\right|\right]\quad\text{for }g\in\Schwartz\left(\R^{\dimension}\right)\text{ and }k,N\in\N_{0}\,.
\]

Next, Lemma~\ref{lem:LocalEmbeddingInHigherLpSpaces} furnishes a
constant $C_{2}=C_{2}\left(\dimension,p,\CalQ\right)>0$ such that
\begin{equation}
\quad\left\Vert \Fourier^{-1}\left(\varphi_{i}\,f\right)\right\Vert _{L^{\infty}}\leq C_{2}\cdot\left|\det T_{i}\right|^{1/p}\cdot\left\Vert \Fourier^{-1}\left(\varphi_{i}\,f\right)\right\Vert _{L^{p}}\quad\forall\,i\in I\text{ and }f\in\DistributionSpace{\CalO}\supset\FourierDecompSp{\CalQ}pY\,.\label{eq:DecompositionSpacesAsTemperedDistributionsLocalEmbeddingInLInfty}
\end{equation}
Furthermore, for $\ell\in\underline{\dimension}$, $k\in\N_{0}$,
and $j\in I$, the chain rule\footnote{See \cite[Lemma 2.6]{DecompositionIntoSobolev} for a more detailed
explanation of how the chain rule is applied here.} implies (with $\varphi_{j}^{\#}:=\varphi_{j}\circ S_{j}$ as in Definition~\ref{def:RegularCoveringRegularPartitionOfUnity})
that
\begin{align*}
\left|\left(\partial_{\ell}^{k}\varphi_{j}\right)\left(\xi\right)\right|=\left|\left(\partial_{\ell}^{k}\left(\varphi_{j}^{\#}\circ\smash{S_{j}^{-1}}\right)\right)\left(\xi\right)\right| & =\left|\left(\partial_{\ell}^{k}\left[\eta\mapsto\varphi_{j}^{\#}\left(T_{j}^{-1}\left(\eta-b_{j}\right)\right)\right]\right)\left(\xi\right)\right|\\
 & =\left|\left(\partial_{\ell}^{k}\left[\varphi_{j}^{\#}\circ T_{j}^{-1}\right]\right)\left(\xi-b_{j}\right)\right|\\
 & \leq C_{3}^{\left(k\right)}\cdot\left\Vert \smash{T_{j}^{-1}}\right\Vert ^{k}\cdot\max_{\left|\alpha\right|\leq k}\left\Vert \partial^{\alpha}\varphi_{j}^{\#}\right\Vert _{\sup}\\
 & \leq C_{3}^{\left(k\right)}\cdot\left\Vert \smash{T_{j}^{-1}}\right\Vert ^{k}\cdot\max_{\left|\alpha\right|\leq k}C_{\Phi,\alpha}\leq C_{4}^{\left(k\right)}\cdot\left\Vert \smash{T_{j}^{-1}}\right\Vert ^{k}
\end{align*}
for suitable constants $C_{3}^{\left(k\right)},C_{4}^{\left(k\right)}$
depending only on $\dimension\in\N$, on $k\in\N_{0}$ and on the
constants $\left(C_{\Phi,\beta}\right)_{\beta\in\N_{0}^{\dimension}}$
from the definition of a regular partition of unity (see Definition~\ref{def:RegularCoveringRegularPartitionOfUnity}).

If we additionally assume $j\in i^{\ast}$ for some $i\in I$, we
get
\[
\left\Vert \smash{T_{j}^{-1}}\right\Vert =\left\Vert \smash{T_{j}^{-1}}T_{i}\smash{T_{i}^{-1}}\right\Vert \leq\left\Vert \smash{T_{j}^{-1}}T_{i}\right\Vert \cdot\left\Vert \smash{T_{i}^{-1}}\right\Vert \leq C_{\CalQ}\cdot\left\Vert \smash{T_{i}^{-1}}\right\Vert ,
\]
which leads to
\[
\left|\left(\smash{\partial_{\ell}^{k}}\varphi_{i}^{\ast}\right)\left(\xi\right)\right|\leq\sum_{j\in i^{\ast}}\left|\left(\smash{\partial_{\ell}^{k}}\varphi_{j}\right)\left(\xi\right)\right|\leq C_{4}^{\left(k\right)}\cdot\sum_{j\in i^{\ast}}\left\Vert \smash{T_{j}^{-1}}\right\Vert ^{k}\leq C_{4}^{\left(k\right)}N_{\CalQ}C_{\CalQ}^{k}\cdot\left\Vert \smash{T_{i}^{-1}}\right\Vert ^{k}=:C_{5}^{\left(k\right)}\cdot\left\Vert \smash{T_{i}^{-1}}\right\Vert ^{k}.
\]
Because of $\varphi_{i}^{\ast}\equiv0$ on $\R^{\dimension}\setminus Q_{i}^{\ast}$,
we arrive at
\begin{equation}
\left|\left(\smash{\partial_{\ell}^{k}}\varphi_{i}^{\ast}\right)\left(\xi\right)\right|\leq C_{5}^{\left(k\right)}\cdot\left\Vert \smash{T_{i}^{-1}}\right\Vert ^{k}\cdot\Indicator_{\overline{Q_{i}^{\ast}}}\left(\xi\right)\;\text{ for all }\xi\in\R^{\dimension},\:i\in I,\text{ and }\ell\in\underline{\dimension},\,k\in\N_{0}\,.\label{eq:DecompositionSpacesAsTemperedDistributionsBAPUDerivativeControl}
\end{equation}

\medskip{}

Now, let $g\in\Schwartz\left(\R^{\dimension}\right)$ and $k\in\N_{0}$
be arbitrary. Using Leibniz's rule we derive
\[
\left|\left(\partial_{\ell}^{k}\left[\varphi_{i}^{\ast}\cdot g\right]\right)\left(\xi\right)\right|\leq\sum_{j=0}^{k}\left[\binom{k}{j}\cdot\left|\left(\smash{\partial_{\ell}^{k-j}}\varphi_{i}^{\ast}\right)\left(\xi\right)\right|\cdot\left|\left(\smash{\partial_{\ell}^{j}}g\right)\left(\xi\right)\right|\right].
\]
But using $c_{i}=\inf_{\xi\in Q_{i}^{\ast}}\left(1+\left|\xi\right|\right)=\min_{\xi\in\overline{Q_{i}^{\ast}}}\left(1+\left|\xi\right|\right)$
and estimate~(\ref{eq:DecompositionSpacesAsTemperedDistributionsBAPUDerivativeControl}),
we see 
\begin{align*}
\left(1+\left|\xi\right|\right)^{-N}\cdot\left|\left(\smash{\partial_{\ell}^{k-j}}\varphi_{i}^{\ast}\right)\left(\xi\right)\right| & \leq C_{5}^{\left(k-j\right)}\cdot\left\Vert \smash{T_{i}^{-1}}\right\Vert ^{k-j}\cdot\Indicator_{\overline{Q_{i}^{\ast}}}\left(\xi\right)\cdot\left(1+\left|\xi\right|\right)^{-N}\\
 & \leq C_{5}^{\left(k-j\right)}\cdot\left\Vert \smash{T_{i}^{-1}}\right\Vert ^{k-j}\cdot\left[\inf_{y\in Q_{i}^{\ast}}\left(1+\left|y\right|\right)\right]^{-N}\\
\left({\scriptstyle \text{distinguish }\left\Vert T_{i}^{-1}\right\Vert \leq1\text{ and }\left\Vert T_{i}^{-1}\right\Vert \geq1}\right) & \leq C_{5}^{\left(k-j\right)}\cdot c_{i}^{-N}\cdot\max\left\{ 1,\left\Vert \smash{T_{i}^{-1}}\right\Vert ^{k}\right\} \qquad\forall\,\xi\in\R^{\dimension}\,.
\end{align*}
Furthermore,
\[
\left(1+\left|\xi\right|\right)^{N+\dimension+1}\cdot\left|\left(\smash{\partial_{\ell}^{j}}g\right)\left(\xi\right)\right|\leq\left\Vert g\right\Vert _{N+\dimension+1,k}\qquad\forall\,\xi\in\R^{\dimension}\text{ and }j\in\left\{ 0,\dots,k\right\} \,.
\]

By setting $C_{6}^{\left(k\right)}:=\sum_{j=0}^{k}\binom{k}{j}C_{5}^{\left(k-j\right)}$,
a combination of the three estimates above shows
\[
\left|\left(\partial_{\ell}^{k}\left[\varphi_{i}^{\ast}\cdot g\right]\right)\left(\xi\right)\right|\leq C_{6}^{\left(k\right)}\left\Vert g\right\Vert _{N+\dimension+1,k}\cdot c_{i}^{-N}\cdot\max\left\{ 1,\left\Vert \smash{T_{i}^{-1}}\right\Vert ^{k}\right\} \cdot\left(1+\left|\xi\right|\right)^{-\dimension-1}
\]
for arbitrary $\xi\in\R^{\dimension}$, $\ell\in\underline{d}$ and
$k\in\N_{0}$. For $C_{7}^{\left(k\right)}:=C_{6}^{\left(k\right)}\cdot\left\Vert \left(1+\left|\mybullet\right|\right)^{-\dimension-1}\right\Vert _{L^{1}}$,
this implies
\begin{equation}
\left\Vert \partial_{\ell}^{k}\left(\varphi_{i}^{\ast}\cdot g\right)\right\Vert _{L^{1}}\leq C_{7}^{\left(k\right)}\cdot\left\Vert g\right\Vert _{N+\dimension+1,k}\cdot c_{i}^{-N}\cdot\max\left\{ 1,\left\Vert \smash{T_{i}^{-1}}\right\Vert ^{k}\right\} \qquad\forall\,i\in I,\,\ell\in\underline{\dimension}\text{ and }k\in\N_{0}.\label{eq:DecompositionSpaceAsTemperedDistributionsMainInequality}
\end{equation}

\medskip{}

Next, it is not hard to see that we have 
\[
\left(1+\left|\xi\right|\right)^{\dimension+1}\leq M_{1}\cdot\left(\,1+\smash{\sum_{\ell=1}^{\dimension}}\,\vphantom{\sum}\left|\xi_{\ell}\right|^{\dimension+1}\,\right)\vphantom{\sum_{\ell=1}^{\dimension}}\qquad\forall\,\xi\in\R^{\dimension}
\]
for a suitable constant $M_{1}>0$ which only depends on $\dimension\in\N$.
Employing the well-known identity (see for instance \cite[Theorem 8.22(e)]{FollandRA})
\[
\xi^{\alpha}\cdot\widehat{g}\left(\xi\right)=\left(2\pi i\right)^{-\left|\alpha\right|}\cdot\widehat{\partial^{\alpha}g}\left(\xi\right)\qquad\forall\,\xi\in\R^{\dimension}\text{ and }g\in\Schwartz\left(\R^{\dimension}\right)\,,
\]
and setting $M_{2}:=M_{1}\cdot\left\Vert \smash{\left(1+\left|\mybullet\right|\right)^{-\dimension-1}}\vphantom{T^{i}}\right\Vert _{L^{1}}$,
we thus derive
\begin{align}
\left\Vert \widehat{\varphi_{i}^{\ast}\cdot g}\right\Vert _{L^{1}} & =\left\Vert \left(1+\left|\mybullet\right|\right)^{-\dimension-1}\cdot\left(1+\left|\mybullet\right|\right)^{\dimension+1}\cdot\widehat{\varphi_{i}^{\ast}\cdot g}\right\Vert _{L^{1}}\nonumber \\
 & \leq M_{1}\cdot\left\Vert \left(1+\left|\mybullet\right|\right)^{-\dimension-1}\right\Vert _{L^{1}}\cdot\left\Vert \xi\mapsto\left(\,1+\smash{\sum_{\ell=1}^{\dimension}}\,\vphantom{\sum}\left|\xi_{\ell}\right|^{\dimension+1}\,\right)\cdot\widehat{\varphi_{i}^{\ast}\cdot g}\left(\xi\right)\right\Vert _{\sup}\vphantom{\sum_{\ell=1}^{\dimension}}\nonumber \\
 & \leq M_{2}\cdot\left[\left\Vert \widehat{\varphi_{i}^{\ast}\cdot g}\right\Vert _{\sup}+\sum_{\ell=1}^{\dimension}\left\Vert \xi\mapsto\xi_{\ell}^{\dimension+1}\cdot\widehat{\varphi_{i}^{\ast}\cdot g}\left(\xi\right)\right\Vert _{\sup}\right]\nonumber \\
 & =M_{2}\cdot\left[\left\Vert \widehat{\varphi_{i}^{\ast}\cdot g}\right\Vert _{\sup}+\sum_{\ell=1}^{\dimension}\left|\left(2\pi i\right)^{-\left(\dimension+1\right)}\right|\cdot\left\Vert \Fourier\left[\partial_{\ell}^{\dimension+1}\left(\varphi_{i}^{\ast}\cdot g\right)\right]\right\Vert _{\sup}\right]\nonumber \\
 & \leq M_{2}\cdot\left[\left\Vert \varphi_{i}^{\ast}\cdot g\right\Vert _{L^{1}}+\sum_{\ell=1}^{\dimension}\left\Vert \partial_{\ell}^{\dimension+1}\left(\varphi_{i}^{\ast}\cdot g\right)\right\Vert _{L^{1}}\right]\label{eq:EstimateFourierL1NormByL1NormOfDerivatives}\\
 & \leq\left(\dimension+1\right)M_{2}\cdot\max\left\{ C_{7}^{\left(0\right)},C_{7}^{\left(\dimension+1\right)}\right\} \cdot\left\Vert g\right\Vert _{N+\dimension+1,d+1}\cdot c_{i}^{-N}\cdot\max\left\{ 1,\left\Vert \smash{T_{i}^{-1}}\right\Vert ^{\dimension+1}\right\} ,\nonumber 
\end{align}
where the last step used estimate~(\ref{eq:DecompositionSpaceAsTemperedDistributionsMainInequality})
and that the quantity $\left\Vert g\right\Vert _{N+\dimension+1,k}\cdot\max\left\{ 1,\smash{\left\Vert T_{i}^{-1}\right\Vert ^{k}}\right\} $
is nondecreasing as a function of $k$ (because both factors are).
The step before used the boundedness of $\Fourier:L^{1}\left(\R^{\dimension}\right)\to C_{0}\left(\R^{\dimension}\right)$.

Setting $M_{3}:=\left(\dimension+1\right)M_{2}\cdot\max\left\{ \smash{C_{7}^{\left(0\right)}},\smash{C_{7}^{\left(\dimension+1\right)}}\right\} $
and $v_{i}^{\left(N\right)}:=c_{i}^{-N}\cdot\max\left\{ 1,\smash{\left\Vert T_{i}^{-1}\right\Vert ^{\dimension+1}}\right\} $
for each $i\in I$, we have thus shown
\begin{equation}
\left\Vert \widehat{\varphi_{i}^{\ast}\cdot g}\right\Vert _{L^{1}}\leq M_{3}\cdot\left\Vert g\right\Vert _{N+\dimension+1,\dimension+1}\cdot v_{i}^{\left(N\right)}\qquad\forall\,i\in I\text{ and }g\in\Schwartz\left(\R^{\dimension}\right).\label{eq:DecompositionSpaceAsTemperedDistributionMainSchwartzEstimate}
\end{equation}

\medskip{}

Next, recall from Lemma~\ref{lem:PartitionCoveringNecessary} that
$\varphi_{i}^{\ast}\equiv1$ on $Q_{i}$. Consequently, $\varphi_{i}^{\ast}\,\varphi_{i}=\varphi_{i}$.
Furthermore, note for $f\in\FourierDecompSp{\CalQ}pY\subset\DistributionSpace{\CalO}$
that $\varphi_{i}\,f$ is a distribution on $\R^{\dimension}$ with
compact support (since $\varphi_{i}\in\TestFunctionSpace{\CalO}$)
and hence $\varphi_{i}\,f\in\Schwartz'\left(\R^{\dimension}\right)$.
All in all, we thus arrive at 
\begin{align*}
\sum_{i\in I_{0}}\left|\left\langle \varphi_{i}\,f,\,g\right\rangle _{\Schwartz'}\right|=\sum_{i\in I_{0}}\left|\left\langle \varphi_{i}\,f,\,\varphi_{i}^{\ast}\,g\right\rangle _{\Schwartz'}\right| & =\sum_{i\in I_{0}}\left|\left\langle \Fourier^{-1}\left(\varphi_{i}\,f\right),\,\widehat{\varphi_{i}^{\ast}\,g}\right\rangle _{\Schwartz'}\right|\\
 & \leq\sum_{i\in I_{0}}\left[\left\Vert \Fourier^{-1}\left(\varphi_{i}\,f\right)\right\Vert _{L^{\infty}}\cdot\left\Vert \widehat{\varphi_{i}^{\ast}\,g}\right\Vert _{L^{1}}\right]\\
\left({\scriptstyle \text{eqs. }\eqref{eq:DecompositionSpacesAsTemperedDistributionsLocalEmbeddingInLInfty}\text{ and }\eqref{eq:DecompositionSpaceAsTemperedDistributionMainSchwartzEstimate}}\right) & \leq C_{2}M_{3}\cdot\left\Vert g\right\Vert _{N+\dimension+1,\dimension+1}\cdot\sum_{i\in I_{0}}\left[\left\Vert \Fourier^{-1}\left(\varphi_{i}\,f\right)\right\Vert _{L^{p}}\cdot\left|\det T_{i}\right|^{\frac{1}{p}}\cdot v_{i}^{\left(N\right)}\right]\\
 & =C_{2}M_{3}\cdot\left\Vert g\right\Vert _{N+\dimension+1,\dimension+1}\cdot\left\Vert \left(\left\Vert \Fourier^{-1}\left(\varphi_{i}\,f\right)\right\Vert _{L^{p}}\cdot w_{i}^{\left(N\right)}\right)_{i\in I_{0}}\right\Vert _{\ell^{1}}\\
\left({\scriptstyle \text{eq. }\eqref{eq:DecompositionSpaceAsTemperedDistributionsQuantitativeSequenceEstimate}}\right) & \leq C_{1}C_{2}M_{3}\cdot\left\Vert g\right\Vert _{N+\dimension+1,\dimension+1}\cdot\left\Vert \left(\left\Vert \Fourier^{-1}\left(\varphi_{i}\,f\right)\right\Vert _{L^{p}}\right)_{i\in I}\right\Vert _{Y}\\
 & =C_{1}C_{2}M_{3}\cdot\left\Vert g\right\Vert _{N+\dimension+1,\dimension+1}\cdot\left\Vert f\right\Vert _{\BAPUFourierDecompSp{\CalQ}pY{\Phi}}<\infty.
\end{align*}
This estimate proves that $\Phi f\in\Schwartz'\left(\R^{\dimension}\right)$
is well-defined with absolute convergence of the series defining $\left\langle \Phi f,\,g\right\rangle _{\Schwartz'}$
for every $g\in\Schwartz\left(\R^{\dimension}\right)$.

Additionally, we see $\left\langle \Phi f,\,g\right\rangle _{\Schwartz'}\rightarrow0$
for $f\rightarrow0$ in $\FourierDecompSp{\CalQ}pY$, for arbitrary
$g\in\Schwartz\left(\R^{\dimension}\right)$, so that $\Phi$ is continuous,
as claimed.

\medskip{}

Finally, Lemma~\ref{lem:PartitionCoveringNecessary} shows that $\left(\varphi_{i}\right)_{i\in I}$
is a locally finite partition of unity on $\CalO$, so that we get
$g=\sum_{i\in I}\varphi_{i}g$ for every $g\in\TestFunctionSpace{\CalO}$,
where only finitely many terms of the sum do not vanish. In case of
$I_{0}=I$, this implies
\[
\left\langle \Phi f,\,g\right\rangle _{\Schwartz'}=\sum_{i\in I}\left\langle \varphi_{i}\,f,\,g\right\rangle _{\Schwartz'}=\sum_{i\in I}\left\langle f,\,\varphi_{i}\,g\right\rangle _{\CalD'}=\vphantom{\sum_{i\in I}}\left\langle \,f,\,\smash{\sum_{i\in I}}\,\vphantom{\sum}\varphi_{i}\,g\,\right\rangle _{\CalD'}=\left\langle f,\,g\right\rangle _{\CalD'}
\]
for all $g\in\TestFunctionSpace{\CalO}$, so that $\Phi f$ is indeed
an extension of $f$.
\end{proof}

\section{Applications}

\label{sec:Applications}In this section, we demonstrate the power
and ease-of-use of our results by considering two classes of classical
spaces, for which our general results extend the state of the art. 

First, we \emph{completely characterize} the existence of embeddings
between $\alpha$-modulation spaces (formally defined below) with
general parameters. Formally, we characterize the existence of the
embedding $\AlphaModSpace{p_{1}}{q_{1}}{s_{1}}{\alpha_{1}}\left(\R^{\dimension}\right)\hookrightarrow\AlphaModSpace{p_{2}}{q_{2}}{s_{2}}{\alpha_{2}}\left(\R^{\dimension}\right)$
in terms of the ``smoothness parameters'' $s_{1},s_{2}\in\R$, the
integrability exponents $p_{1},p_{2},q_{1},q_{2}\in\left(0,\infty\right]$
and in terms of $\alpha_{1},\alpha_{2}\in\left[0,1\right]$.

In this generality, the derived embedding result greatly improves
the state of the art: Previously, the most general characterization
of embeddings between $\alpha$-modulation spaces was obtained by
Han and Wang in \cite{HanWangAlphaModulationEmbeddings}; but they
only allow $\left(p_{1},q_{1}\right)\neq\left(p_{2},q_{2}\right)$
in case of $\alpha_{1}=\alpha_{2}$, whereas they assume $\left(p_{1},q_{1}\right)=\left(p_{2},q_{2}\right)$
in case of $\alpha_{1}\neq\alpha_{2}$. Our results need no such restriction.

Incidentally, only a few days after the first version of the present
paper appeared on the arXiv, the preprint version of the paper \cite{GuoAlphaModulationEmbeddingCharacterization}
of Guo, Fan, Wu, and Zhao was submitted to the arXiv. In that paper,
the authors independently obtain \emph{the same} fully general characterization
of the existence of embeddings between different $\alpha$-modulation
spaces. The proof techniques used in \cite{GuoAlphaModulationEmbeddingCharacterization}
and in the present paper are both based in parts on the arguments
of Han and Wang in \cite{HanWangAlphaModulationEmbeddings}. Nevertheless,
the two proofs differ quite a bit: While we developed a general framework
for deciding the existence of embeddings between decomposition spaces,
which we then simply need to apply to the special case of $\alpha$-modulation
spaces, the authors of \cite{GuoAlphaModulationEmbeddingCharacterization}
work only in the special setting of $\alpha$-modulation spaces.

Precisely, they first show for $0\leq\alpha_{1}\leq\alpha_{2}\leq1$
that $\AlphaModSpace pqs{\alpha_{1}}\left(\R^{\dimension}\right)=\CalD_{{\rm general}}\left(\CalQ^{\left(\alpha_{2}\right)},\AlphaModSpace pq0{\alpha_{1}},\ell_{w^{\left(\alpha_{2},s\right)}}^{q}\right)$;
see \cite[Proposition 3.1]{GuoAlphaModulationEmbeddingCharacterization}.
Here, the \emph{generalized decomposition space} on the right-hand
side is to be understood as in the paper \cite{DecompositionSpaces1}
by Feichtinger and Gröbner—it is defined just as the decomposition
spaces considered in the present paper, the only difference being
that when computing the (quasi)-norm, the space $L^{p}$ is replaced
by the $\alpha_{1}$-modulation space $\AlphaModSpace pq0{\alpha_{1}}$.
The above characterization of $\AlphaModSpace pqs{\alpha_{1}}\left(\R^{\dimension}\right)$
shows that the $\alpha_{1}$-modulation space $\AlphaModSpace pqs{\alpha_{1}}\left(\R^{\dimension}\right)=\DecompSp{\CalQ^{\left(\alpha_{1}\right)}}p{\ell_{w^{\left(\alpha_{1},s\right)}}^{q}}$—which
is defined using the ``$\alpha_{1}$-modulation covering'' $\CalQ^{\left(\alpha_{1}\right)}$—can
be interpreted as a (generalized) decomposition space that uses the
very different covering $\CalQ^{\left(\alpha_{2}\right)}$; but this
requires changing the local component $L^{p}$ to the $\alpha_{1}$-modulation
space $\AlphaModSpace pq0{\alpha_{1}}\left(\R^{\dimension}\right)$.
We remark that the space $\ell_{w^{\left(\alpha,s\right)}}^{q}\left(\Z^{\dimension}\right)$
is the usual ``global component'' used to define the $\alpha$-modulation
spaces $\AlphaModSpace pqs{\alpha}\left(\R^{\dimension}\right)=\DecompSp{\CalQ^{\left(\alpha\right)}}p{\ell_{w^{\left(s,\alpha\right)}}^{q}}$;
see Definition~\ref{def:AlphaModulationFourierSpace} and Corollary~\ref{cor:FourierAlphaModulationAsTemperedDistributions}.

Guo et al.\@ then use this characterization of the $\alpha_{1}$-modulation
space $\AlphaModSpace pqs{\alpha_{1}}\left(\R^{\dimension}\right)$
to give an alternative description of the space $M_{\Fourier}\left(\AlphaModSpace{p_{1}}{q_{1}}{s_{1}}{\alpha_{1}},\AlphaModSpace{p_{2}}{q_{2}}{s_{2}}{\alpha_{2}}\right)$
of all Fourier multipliers that act boundedly as a map $\AlphaModSpace{p_{1}}{q_{1}}{s_{1}}{\alpha_{1}}\left(\R^{\dimension}\right)\to\AlphaModSpace{p_{2}}{q_{2}}{s_{2}}{\alpha_{2}}\left(\R^{\dimension}\right)$.
Precisely, \cite[Theorem 3.2]{GuoAlphaModulationEmbeddingCharacterization}
shows for $\alpha:=\max\left\{ \alpha_{1},\alpha_{2}\right\} $ that
\[
M_{\Fourier}\left(\AlphaModSpace{p_{1}}{q_{1}}{s_{1}}{\alpha_{1}},\AlphaModSpace{p_{2}}{q_{2}}{s_{2}}{\alpha_{2}}\right)=\CalD_{{\rm general}}\left(\CalQ^{\left(\alpha\right)},M_{\Fourier}\left(\AlphaModSpace{p_{1}}{q_{1}}0{\alpha_{1}},\AlphaModSpace{p_{2}}{q_{2}}0{\alpha_{2}}\right),M_{p}\left(\ell_{w^{\left(s_{1},\alpha\right)}}^{q_{1}},\ell_{w^{\left(s_{2},\alpha\right)}}^{q_{2}}\right)\right)\,.\tag{\ensuremath{\ast}}
\]
Here, the right-hand side is again a generalized decomposition space
as introduced in \cite{DecompositionSpaces1}; furthermore, $M_{p}\left(\ell_{w^{\left(s_{1},\alpha\right)}}^{q_{1}},\ell_{w^{\left(s_{2},\alpha\right)}}^{q_{2}}\right)$
denotes the space of all sequences $\left(m_{k}\right)_{k\in\Z^{\dimension}}$
for which the associated \emph{pointwise multiplication operator}
$\left(c_{k}\right)_{k\in\Z^{\dimension}}\mapsto\left(m_{k}\cdot c_{k}\right)_{k\in\Z^{\dimension}}$
defines a bounded operator from $\ell_{w^{\left(s_{1},\alpha\right)}}^{q}\left(\Z^{\dimension}\right)$
into $\ell_{w^{\left(s_{2},\alpha\right)}}^{q}\left(\Z^{\dimension}\right)$.

Finally, we have $\AlphaModSpace{p_{1}}{q_{1}}{s_{1}}{\alpha_{1}}\left(\R^{\dimension}\right)\hookrightarrow\AlphaModSpace{p_{2}}{q_{2}}{s_{2}}{\alpha_{2}}\left(\R^{\dimension}\right)$
if and only if $\left(\xi\mapsto1\right)\in M_{\Fourier}\left(\AlphaModSpace{p_{1}}{q_{1}}{s_{1}}{\alpha_{1}},\AlphaModSpace{p_{2}}{q_{2}}{s_{2}}{\alpha_{2}}\right)$,
which can be decided using the characterization $\left(\ast\right)$.
Note, however, that this is still highly nontrivial; see \cite[Section 4]{GuoAlphaModulationEmbeddingCharacterization}.
The arguments in that section are based on those of Han and Wang in
\cite{HanWangAlphaModulationEmbeddings}, and are essentially simplified
versions of those used in Sections~\ref{sec:SufficientConditions}
and \ref{sec:NecessaryConditions} of the present paper. Of course,
the arguments in the present paper are much more involved, since they
have to apply to general coverings, and not only to the relatively
simple case of the coverings $\CalQ^{\left(\alpha\right)}$ used to
define the $\alpha$-modulation spaces.

Overall, the approach in \cite{GuoAlphaModulationEmbeddingCharacterization}
is more elementary and thus also more accessible than the present
paper. But this is only true if one considers the entirety of the
present paper; if instead we consider Theorems~\ref{thm:SummaryCoarseIntoFine}
and \ref{thm:SummaryFineIntoCoarse} as given, then the present treatment
is much shorter. Furthermore, our approach clearly shows which properties
of the $\alpha$-modulation spaces determine the existence of the
embedding: As we will see, the crucial point is that the covering
$\CalQ^{\left(\alpha\right)}$ used to define the $\alpha$-modulation
spaces is almost subordinate to, and relatively moderate with respect
to the covering $\CalQ^{\left(\beta\right)}$ if $\alpha\leq\beta$.

Finally, we emphasize that—once the ``embedding framework'' is given—most
of the work is needed to obtain a proper understanding of the relation
between the coverings $\CalQ^{\left(\alpha\right)}$ and $\CalQ^{\left(\beta\right)}$
which are used to define the $\alpha$-modulation spaces. Once this
understanding is obtained, the actual embedding statements can be
derived with ease. In many other treatments, the properties of the
coverings $\CalQ^{\left(\alpha\right)}$ which we formally derive
here are essentially taken for granted, and only the existence of
the embedding is considered with care.

\medskip{}

As our second example, we consider the question of embeddings between
homogeneous and inhomogeneous Besov spaces. The relation between these
spaces is particularly interesting, since they live on slightly different
frequency domains: The covering $\CalB$ used to define the inhomogeneous
Besov spaces covers $\CalO=\R^{\dimension}$, while the covering $\dot{\CalB}$
for the homogeneous spaces only covers $\CalO'=\R^{\dimension}\setminus\left\{ 0\right\} $.
In particular, $\dot{\CalB}$ is not relatively moderate with respect
to $\CalB$ (see Lemma~\ref{lem:DifferentOrbitsPreventModerateness}).
Therefore, this example shows that the criteria developed in Sections
\ref{subsec:ElementaryNecessaryConditions}, \ref{subsec:ImprovedNecessaryConditions},
and \ref{subsec:NecessaryForP1EqualP2} are important, since there
are relevant cases in which the complete characterizations from Section~\ref{subsec:RelativelyModerateCase}
do not apply.

\medskip{}

We finally remark that a huge number of other examples—including ($\alpha$)-shearlet
smoothness spaces and a large class of wavelet-type coorbit spaces—can
be handled using the framework developed in this paper. But since
the present paper already has a considerable length, we postpone these
applications to later contributions. The interested reader can already
find some of these applications in my PhD thesis \cite{VoigtlaenderPhDThesis},
and in the recent preprint \cite{AlphaShearletSparsity}. The examples
considered in the present section are essentially taken from my thesis
\cite{VoigtlaenderPhDThesis}.

\subsection{Embeddings of \texorpdfstring{$\alpha$-modulation}{α-modulation}
spaces}

\label{subsec:AlphaModulationEmbedding}We begin our study of $\alpha$-modulation
spaces by describing the construction of the associated covering $\CalQ^{\left(\alpha\right)}$,
as given by Borup and Nielsen in \cite{BorupNielsenAlphaModulationSpaces}.
\begin{thm}
\label{thm:AlphaCoveringExistence}(cf. \cite[Theorem 2.6]{BorupNielsenAlphaModulationSpaces})
Let $\dimension\in\N$ and $0\leq\alpha<1$ be arbitrary and set $\alpha_{0}:=\frac{\alpha}{1-\alpha}$.
Then there is a constant $r_{0}=r_{0}\left(\dimension,\alpha\right)>0$
such that for every $r>r_{0}$, the family
\[
\CalQ^{\left(\alpha\right)}:=\CalQ_{r}^{\left(\alpha\right)}:=\left(Q_{r,k}^{\left(\alpha\right)}\right)_{k\in\Z^{\dimension}\setminus\left\{ 0\right\} }:=\left(\vphantom{Q_{r,k}^{\left(\alpha\right)}}B_{r\left|k\right|^{\alpha_{0}}}\left(\left|k\right|^{\alpha_{0}}k\right)\right)_{k\in\Z^{\dimension}\setminus\left\{ 0\right\} }
\]
is an admissible covering of $\R^{\dimension}$, called the \textbf{$\alpha$-covering
of $\R^{\dimension}$}.
\end{thm}

Our plan is to examine the (relative) geometry of the coverings $\CalQ^{\left(\alpha\right)}$
by showing that
\begin{enumerate}
\item the covering $\CalQ^{\left(\alpha\right)}$ is a \emph{structured}
admissible covering of $\R^{\dimension}$,
\item the weight $\Z^{\dimension}\setminus\left\{ 0\right\} \rightarrow\left(0,\infty\right),k\mapsto\left\langle k\right\rangle :=\sqrt{1+\left|k\right|^{2}}$
is moderate with respect to $\CalQ^{\left(\alpha\right)}$,
\item $\CalQ^{\left(\alpha\right)}$ is almost subordinate to and relatively
moderate with respect to $\CalQ^{\left(\beta\right)}$ for $\alpha\leq\beta$.
\end{enumerate}
We begin with establishing the second point.
\begin{lem}
\label{lem:AlphaModulationCoveringNormEstimate}Let $\dimension\in\N$,
$0\leq\alpha<1$ and $r>0$ so that $\CalQ^{\left(\alpha\right)}=\CalQ_{r}^{\left(\alpha\right)}$
is an admissible covering of $\R^{\dimension}$. For $\xi\in\R^{\dimension}$,
set $\left\langle \xi\right\rangle :=\left(1+\smash{\left|\xi\right|^{2}}\right)^{1/2}\vphantom{\left(\left|\xi\right|^{2}\right)^{1/2}}$.
We then have
\[
\left\langle \xi\right\rangle \asymp\left\langle k\right\rangle ^{\frac{1}{1-\alpha}}\qquad\text{ for all }\qquad k\in\Z^{\dimension}\setminus\left\{ 0\right\} \quad\text{and}\quad\xi\in Q_{r,k}^{\left(\alpha\right)},
\]
where the implied constants only depend on $r,\alpha$.

In particular, for any $\gamma\in\R$, the weight $w^{\left(\gamma\right)}:=\left(\left\langle k\right\rangle ^{\gamma}\right)_{k\in\Z^{\dimension}\setminus\left\{ 0\right\} }$
is moderate with respect to $\CalQ_{r}^{\left(\alpha\right)}$.
\end{lem}

\begin{proof}
Define $\alpha_{0}:=\frac{\alpha}{1-\alpha}\in\left[0,\infty\right)$.
For $k\in\Z^{\dimension}\setminus\left\{ 0\right\} $ and arbitrary
$\xi\in Q_{r,k}^{\left(\alpha\right)}=B_{r\left|k\right|^{\alpha_{0}}}\left(\left|k\right|^{\alpha_{0}}k\right)$,
we have
\begin{align*}
\left|\xi\right|\leq\bigl|\xi-\left|k\right|^{\alpha_{0}}k\bigr|+\bigl|\left|k\right|^{\alpha_{0}}k\bigr| & <r\left|k\right|^{\alpha_{0}}+\left|k\right|^{\alpha_{0}}\left|k\right|\\
 & \leq\left\langle k\right\rangle ^{\alpha_{0}}\cdot\left[r+\left\langle k\right\rangle \right]\\
\left({\scriptstyle \text{since }\left\langle k\right\rangle \geq1}\right) & \leq\left(1+r\right)\cdot\left\langle k\right\rangle ^{\alpha_{0}}\cdot\left\langle k\right\rangle =\left(1+r\right)\cdot\left\langle k\right\rangle ^{\frac{1}{1-\alpha}}.
\end{align*}
Because of $\frac{1}{1-\alpha}>0$, we have $\left\langle k\right\rangle ^{\frac{1}{1-\alpha}}\geq1$.
In combination, this yields ``$\lesssim$''; indeed,
\[
\left\langle \xi\right\rangle \leq1+\left|\xi\right|\leq\left\langle k\right\rangle ^{\frac{1}{1-\alpha}}+\left|\xi\right|\leq\left(2+r\right)\cdot\left\langle k\right\rangle ^{\frac{1}{1-\alpha}}.
\]
For the proof of ``$\gtrsim$'', we distinguish two cases:

\begin{casenv}
\item We have $\left|k\right|\geq2r$. This implies $\left|k\right|-r\geq\left|k\right|-\frac{\left|k\right|}{2}=\frac{\left|k\right|}{2}$
and hence
\begin{align*}
\qquad\qquad\left|\xi\right|=\left|\vphantom{\sum}\left|k\right|^{\alpha_{0}}k-\left(\left|k\right|^{\alpha_{0}}k-\xi\right)\right| & \geq\bigl|\left|k\right|^{\alpha_{0}}k\bigr|-\bigl|\left|k\right|^{\alpha_{0}}k-\xi\bigr|\\
 & >\left|k\right|^{\alpha_{0}}\left|k\right|-r\left|k\right|^{\alpha_{0}}\\
 & =\left|k\right|^{\alpha_{0}}\cdot\left(\left|k\right|-r\right)\geq\frac{1}{2}\cdot\left|k\right|^{\alpha_{0}}\left|k\right|=\frac{1}{2}\cdot\left|k\right|^{\frac{1}{1-\alpha}}.
\end{align*}
Because of $k\in\Z^{\dimension}\setminus\left\{ 0\right\} $, we have
$\left|k\right|\geq1$. But this implies $\left\langle k\right\rangle \leq1+\left|k\right|\leq2\left|k\right|$;
therefore,
\[
\left\langle k\right\rangle ^{\frac{1}{1-\alpha}}\leq2^{\frac{1}{1-\alpha}}\cdot\left|k\right|^{\frac{1}{1-\alpha}}\leq2^{1+\frac{1}{1-\alpha}}\cdot\left|\xi\right|\leq2^{\frac{2-\alpha}{1-\alpha}}\cdot\left\langle \xi\right\rangle .
\]
\item We have $\left|k\right|\leq2r$. In this case, we use the positivity
of $\frac{1}{1-\alpha}$ together with $\left\langle \xi\right\rangle \geq1$
and $\left\langle k\right\rangle \leq1+\left|k\right|\leq1+2r$ to
derive $\left\langle k\right\rangle ^{\frac{1}{1-\alpha}}\leq\left(1+2r\right)^{\frac{1}{1-\alpha}}\leq\left(1+2r\right)^{\frac{1}{1-\alpha}}\cdot\left\langle \xi\right\rangle $.\vspace{0.2cm}
\end{casenv}
Overall, this proves $\left\langle k\right\rangle ^{\frac{1}{1-\alpha}}\leq C_{0}\cdot\left\langle \xi\right\rangle $
with $C_{0}:=\max\left\{ \left(1+2r\right)^{\frac{1}{1-\alpha}},2^{\frac{2-\alpha}{1-\alpha}}\right\} $.
We have thus established the first part of the lemma.

To prove the $\CalQ_{r}^{\left(\alpha\right)}$-moderateness of $w^{\left(\gamma\right)}$,
note that $\xi\in Q_{r,k}^{\left(\alpha\right)}\cap Q_{r,\ell}^{\left(\alpha\right)}$
implies $\left\langle k\right\rangle \asymp\left\langle \xi\right\rangle ^{1-\alpha}\asymp\left\langle \ell\right\rangle $
and hence $1\lesssim\left\langle k\right\rangle /\left\langle \ell\right\rangle \lesssim1$,
which yields $1\lesssim\left(\left\langle k\right\rangle /\left\langle \ell\right\rangle \right)^{\gamma}\lesssim1$,
and thus finally $w_{k}^{\left(\gamma\right)}=\left\langle k\right\rangle ^{\gamma}\asymp\left\langle \ell\right\rangle ^{\gamma}=w_{\ell}^{\left(\gamma\right)}$,
where the implied constants only depend on $r,\alpha,\gamma$.
\end{proof}
Using the moderateness of the weight $w^{\left(\gamma\right)}$, we
can now show that $\CalQ_{r}^{\left(\alpha\right)}$ is a \emph{structured}
admissible covering of $\R^{\dimension}$.
\begin{lem}
\label{lem:AlphaModulationCoveringIsStructured}Let $\dimension\in\N$,
$\alpha\in\left[0,1\right)$ and set $\alpha_{0}:=\frac{\alpha}{1-\alpha}$.
Let $r_{0}=r_{0}\left(\dimension,\alpha\right)$ as in Theorem~\ref{thm:AlphaCoveringExistence}.
Then
\[
\CalQ_{r}^{\left(\alpha\right)}=\left(Q_{r,k}^{\left(\alpha\right)}\right)_{k\in\Z^{\dimension}\setminus\left\{ 0\right\} }=\left(T_{k}^{\left(\alpha\right)}Q^{\left(r\right)}+b_{k}\right)_{k\in\Z^{\dimension}\setminus\left\{ 0\right\} }
\]
is a structured admissible covering of $\R^{\dimension}$ for each
$r>r_{0}$, with
\[
T_{k}^{\left(\alpha\right)}:=\left|k\right|^{\alpha_{0}}\cdot{\rm id}\quad\text{and}\quad b_{k}:=\left|k\right|^{\alpha_{0}}k\quad\text{for}\quad k\in\Z^{\dimension},\quad\text{and with}\quad Q^{\left(r\right)}:=B_{r}\left(0\right)\,.\qedhere
\]
\end{lem}

\begin{proof}
Let $s:=\frac{r_{0}+r}{2}$, so that we have $r_{0}<s<r$. Define
$P:=Q^{\left(s\right)}:=B_{s}\left(0\right)$. This implies $\overline{P}=\overline{Q^{\left(s\right)}}\subset Q^{\left(r\right)}=:Q$,
where $P,Q$ are both open and bounded. Using Theorem~\ref{thm:AlphaCoveringExistence},
we see that $\CalQ_{s}^{\left(\alpha\right)}=\left(\smash{T_{k}^{\left(\alpha\right)}P+b_{k}}\right)_{k\in\Z^{\dimension}\setminus\left\{ 0\right\} }$
and $\CalQ_{r}^{\left(\alpha\right)}=\left(\smash{T_{k}^{\left(\alpha\right)}Q+b_{k}}\right)_{k\in\Z^{\dimension}\setminus\left\{ 0\right\} }$
are both admissible coverings of $\R^{\dimension}$.

It remains to show that $C_{\CalQ_{r}^{\left(\alpha\right)}}$ is
finite. To see this, let $k,\ell\in\Z^{\dimension}\setminus\left\{ 0\right\} $
with $\emptyset\neq Q_{r,k}^{\left(\alpha\right)}\cap Q_{r,\ell}^{\left(\alpha\right)}$.
Using $\left\langle k\right\rangle \leq1+\left|k\right|\leq2\left|k\right|$,
as well as $\left|\ell\right|\leq\left\langle \ell\right\rangle $
and $\alpha_{0}\geq0$, we see
\[
\left\Vert \bigl(\smash{T_{k}^{\left(\alpha\right)}}\bigr)^{-1}T_{\ell}^{\left(\alpha\right)}\right\Vert =\left(\frac{\left|\ell\right|}{\left|k\right|}\right)^{\alpha_{0}}\leq\left(\frac{\left\langle \ell\right\rangle }{\left|k\right|}\right)^{\alpha_{0}}\leq\left(2\frac{\left\langle \ell\right\rangle }{\left\langle k\right\rangle }\right)^{\alpha_{0}}=2^{\alpha_{0}}\cdot\frac{w_{\ell}^{\left(\smash{\alpha_{0}}\right)}}{w_{k}^{\left(\smash{\alpha_{0}}\right)}}\leq2^{\alpha_{0}}\cdot C_{w^{\left(\smash{\alpha_{0}}\right)},\CalQ_{r}^{\left(\alpha\right)}}.
\]
But $C_{w^{\left(\smash{\alpha_{0}}\right)},\CalQ_{r}^{\left(\alpha\right)}}$
is finite as a consequence of Lemma~\ref{lem:AlphaModulationCoveringNormEstimate}.
This completes the proof.
\end{proof}
Lemma~\ref{lem:AlphaModulationCoveringIsStructured} (together with
Theorem~\ref{thm:AlmostStructuredAdmissibleAdmitsBAPU}) implies
in particular that the $\alpha$-covering $\CalQ_{r}^{\left(\alpha\right)}$
admits a family $\left(\varphi_{k}\right)_{k\in\Z^{\dimension}\setminus\left\{ 0\right\} }$
which is an $L^{p}$-BAPU for $\CalQ_{r}^{\left(\alpha\right)}$ for
every $p\in\left(0,\infty\right]$. Furthermore, the weight $w^{\left(\gamma\right)}$
is $\CalQ_{r}^{\left(\alpha\right)}$-moderate by Lemma~\ref{lem:AlphaModulationCoveringNormEstimate}.
Hence, the associated decomposition spaces are well-defined.
\begin{defn}
\label{def:AlphaModulationFourierSpace}Let $\dimension\in\N$ and
$0\leq\alpha<1$ be arbitrary. Choose $r_{0}=r_{0}\left(\dimension,\alpha\right)$
as in Theorem~\ref{thm:AlphaCoveringExistence}. For each $r>r_{0}$,
and $p,q\in\left(0,\infty\right]$, as well as $\gamma\in\R$, we
define the \textbf{Fourier-side $\alpha$-modulation space} with integrability
exponents $p,q$ and weight-exponent $\gamma$ as
\[
\AlphaModSpace pq{\Fourier,\gamma}{\alpha}\left(\R^{\dimension}\right):=\FourierDecompSp{\CalQ_{r}^{\left(\alpha\right)}}p{\ell_{w^{\left(\gamma/\left(1-\alpha\right)\right)}}^{q}},
\]
where the weight $w^{\left(\gamma/\left(1-\alpha\right)\right)}=\left(\smash{\left\langle k\right\rangle ^{\gamma/\left(1-\alpha\right)}}\right)_{k\in\Z^{\dimension}\setminus\left\{ 0\right\} }$
is defined as in Lemma~\ref{lem:AlphaModulationCoveringNormEstimate}.
\end{defn}

\begin{rem*}
(1) We will show below (see Corollary~\ref{cor:FourierAlphaModulationAsTemperedDistributions})
that the Fourier transform restricts to an (isometric) isomorphism
\[
\Fourier:\AlphaModSpace pq{\gamma}{\alpha}\left(\smash{\R^{\dimension}}\right)\subset\Schwartz'\left(\R^{\dimension}\right)\rightarrow\AlphaModSpace pq{\Fourier,\gamma}{\alpha}\left(\R^{\dimension}\right)\subset\DistributionSpace{\R^{\dimension}},f\mapsto\widehat{f}|_{\TestFunctionSpace{\R^{\dimension}}}\,,
\]
where the (space-side) $\alpha$-modulation space on the left-hand
side is defined as in \cite[Definition 2.4]{BorupNielsenAlphaModulationSpaces}.
This justifies the name ``Fourier-side $\alpha$-modulation space''.

\smallskip{}

(2) The notation for the $\alpha$-modulation spaces is not completely
consistent in the literature. For example, the space that we denote
by $\AlphaModSpace pq{\gamma}{\alpha}\left(\R^{\dimension}\right)$
is denoted by $M_{p,q}^{\gamma,\alpha}\left(\R^{\dimension}\right)$
in \cite{BorupNielsenAlphaModulationSpaces}. The notation adopted
here is motivated by the usual notation $L^{p}\left(\R^{\dimension}\right)$
and $L^{p,q}\left(\R^{\dimension}\right)$ for the Lebesgue spaces
and Lorentz spaces: the integrability exponents are at the top.

\smallskip{}

(3) We observe that the parameter $r$ is suppressed on the left-hand
side of the definition above; but since any two coverings $\CalQ_{r}^{\left(\alpha\right)},\CalQ_{s}^{\left(\alpha\right)}$
(with $r,s>r_{0}\left(\dimension,\alpha\right)$) use the \emph{same}
parametrization, it follows that every $L^{p}$-BAPU $\Phi$ for $\CalQ_{r}^{\left(\alpha\right)}$
is also an $L^{p}$-BAPU for $\CalQ_{s}^{\left(\alpha\right)}$—at
least for $s\geq r$, which we can always assume by symmetry. With
this choice of the BAPUs for both spaces, we get
\[
\left\Vert \mybullet\right\Vert _{\BAPUFourierDecompSp{\smash{\CalQ_{r}^{\left(\alpha\right)}}}p{\smash{\ell_{w^{\left(\gamma/\left(1-\alpha\right)\right)}}^{q}}}{\Phi}}=\left\Vert \mybullet\right\Vert _{\BAPUFourierDecompSp{\smash{\CalQ_{s}^{\left(\alpha\right)}}}p{\smash{\ell_{w^{\left(\gamma/\left(1-\alpha\right)\right)}}^{q}}}{\Phi}},
\]
which implies that both spaces coincide. Since different choices of
the BAPU yield the same space with equivalent quasi-norms (see Corollary~\ref{cor:DecompositionSpaceWellDefined}),
we see that the right-hand side of the above definition is independent
of the choice of $r>r_{0}\left(\dimension,\alpha\right)$, with equivalent
quasi-norms for different choices.
\end{rem*}
As the next step in our program, we will estimate the cardinalities
of the ``intersection sets'' $J_{i}$, where $\CalQ=\left(Q_{i}\right)_{i\in I}$
and $\CalP=\left(P_{j}\right)_{j\in J}$ are both certain $\alpha$-coverings,
for different values of $\alpha$. As a byproduct, we will then see
(using Lemma~\ref{lem:WeaksubordinatenessImpliesSubordinatenessIfConnected})
that $\CalQ_{r}^{\left(\alpha\right)}$ is almost subordinate to $\CalQ_{s}^{\left(\beta\right)}$
for $\alpha\le\beta$.
\begin{lem}
\label{lem:AlphaModulationIntersectionCountSubordinateness}Let $\dimension\in\N$
and $\alpha,\beta\in\left[0,1\right)$ with $\alpha\leq\beta$. Let
$r_{1}=r_{0}\left(\dimension,\alpha\right)$ and $r_{2}=r_{0}\left(\dimension,\beta\right)$
with $r_{0}$ as in Theorem~\ref{thm:AlphaCoveringExistence}, and
let $r>r_{1}$ and $s>r_{2}$.

There is a constant $C=C\left(\dimension,\alpha,\beta,r,s\right)>0$
with
\begin{equation}
1\leq\left|\left\{ \ell\in\Z^{\dimension}\setminus\left\{ 0\right\} \with Q_{s,\ell}^{\left(\beta\right)}\cap Q_{r,k}^{\left(\alpha\right)}\neq\emptyset\right\} \right|\leq C\qquad\forall\,k\in\Z^{\dimension}\setminus\left\{ 0\right\} \,.\label{eq:AlphaModulationIntersectionCountSubordinateness}
\end{equation}
In particular, $\CalQ_{r}^{\left(\alpha\right)}$ is almost subordinate
to $\CalQ_{s}^{\left(\beta\right)}$.
\end{lem}

\begin{proof}
To be consistent with our usual notation, set $\CalQ=\left(Q_{i}\right)_{i\in I}:=\CalQ_{r}^{\left(\alpha\right)}$
and $\CalP=\left(P_{j}\right)_{j\in J}:=\CalQ_{s}^{\left(\beta\right)}$.
Then, the intersection set $J_{k}$ is given by
\[
J_{k}=\left\{ \ell\in\Z^{\dimension}\setminus\left\{ 0\right\} \with Q_{s,\ell}^{\left(\beta\right)}\cap Q_{r,k}^{\left(\alpha\right)}\neq\emptyset\right\} \quad\text{for any }k\in\Z^{\dimension}\setminus\left\{ 0\right\} \,.
\]
Once we have established estimate~(\ref{eq:AlphaModulationIntersectionCountSubordinateness}),
we have shown that $\CalQ$ is \emph{weakly} subordinate to $\CalP$,
see Definition~\ref{def:RelativeIndexClustersSubordinateCoveringsModerateCoverings}.
But since all sets in $\CalQ$ and $\CalP$ are open balls, which
are open and path-connected, Corollary~\ref{cor:WeakSubordinationImpliesSubordinationIfConnected}
easily implies that $\CalQ=\CalQ_{r}^{\left(\alpha\right)}$ is almost
subordinate to $\CalP=\CalQ_{s}^{\left(\beta\right)}$, since $\CalQ,\CalP$
are both coverings of all of $\R^{\dimension}$.

It thus remains to establish equation~(\ref{eq:AlphaModulationIntersectionCountSubordinateness}).
To this end, set $\alpha_{0}:=\frac{\alpha}{1-\alpha}$ and $\beta_{0}:=\frac{\beta}{1-\beta}$.
For $m\in\Z^{\dimension}\setminus\left\{ 0\right\} $ and $\gamma\in\left[0,1\right)$,
let us denote the \textbf{$\gamma$-normalized version} of $m$ by
$m^{\left(\gamma\right)}:=\left|m\right|^{\gamma_{0}}m$ with $\gamma_{0}:=\frac{\gamma}{1-\gamma}$.
By definition, we have $Q_{t,m}^{\left(\gamma\right)}=B_{t\left|m\right|^{\gamma_{0}}}\left(\smash{m^{\left(\gamma\right)}}\right)$
for all $t>0$ and $m\in\Z^{\dimension}\setminus\left\{ 0\right\} $.

Lemma~\ref{lem:AlphaModulationCoveringNormEstimate} implies
\begin{equation}
\begin{split}\left\langle k\right\rangle ^{\frac{1}{1-\alpha}}\asymp\left\langle \xi\right\rangle  & \qquad\forall\,k\in\Z^{\dimension}\setminus\left\{ 0\right\} \text{ and }\xi\in Q_{r,k}^{\left(\alpha\right)},\\
\left\langle \ell\right\rangle ^{\frac{1}{1-\beta}}\asymp\left\langle \xi\right\rangle  & \qquad\forall\,\ell\in\Z^{\dimension}\setminus\left\{ 0\right\} \text{ and }\xi\in Q_{s,\ell}^{\left(\beta\right)},
\end{split}
\label{eq:AlphaModulationIntersectionCountNormEstimate}
\end{equation}
where the implied constants only depend on $\dimension,\alpha,\beta,r,s$.

The lower estimate in equation~(\ref{eq:AlphaModulationIntersectionCountSubordinateness})
is trivial, since $\emptyset\neq Q_{r,k}^{\left(\alpha\right)}\subset\R^{\dimension}=\bigcup_{\ell\in\Z^{\dimension}\setminus\left\{ 0\right\} }Q_{s,\ell}^{\left(\beta\right)}$.

To prove the upper estimate, fix $k\in\Z^{\dimension}\setminus\left\{ 0\right\} $,
some $\ell_{0}\in J_{k}$ and some $\xi_{0}\in Q_{s,\ell_{0}}^{\left(\beta\right)}\cap Q_{r,k}^{\left(\alpha\right)}$.
Using equation~(\ref{eq:AlphaModulationIntersectionCountNormEstimate}),
this yields
\[
\left|k\right|^{\alpha_{0}}\leq\left(\left\langle k\right\rangle ^{\frac{1}{1-\alpha}}\right)^{\alpha}\lesssim\left\langle \xi_{0}\right\rangle ^{\alpha}\lesssim\left(\left\langle \ell_{0}\right\rangle ^{\frac{1}{1-\beta}}\right)^{\alpha}\,\overset{\left(\dagger\right)}{\leq}\,\left\langle \ell_{0}\right\rangle ^{\frac{\beta}{1-\beta}}=\left\langle \ell_{0}\right\rangle ^{\beta_{0}}
\]
where the step marked with $\left(\dagger\right)$ is justified by
$\left\langle \ell_{0}\right\rangle \geq1$ and $\alpha\leq\beta$.
As usual, the implied constants only depend on $\dimension,\alpha,\beta,r,s$.

Note that for each $\ell\in J_{k}$, there is some $\xi\in Q_{r,k}^{\left(\alpha\right)}\cap Q_{s,\ell}^{\left(\beta\right)}$.
Hence, $\xi,\xi_{0}\in Q_{r,k}^{\left(\alpha\right)}$, which yields
$\left\langle \xi\right\rangle \asymp\left\langle k\right\rangle ^{1/\left(1-\alpha\right)}\asymp\left\langle \xi_{0}\right\rangle $
by equation~(\ref{eq:AlphaModulationIntersectionCountNormEstimate}).
Now, another application of equation~(\ref{eq:AlphaModulationIntersectionCountNormEstimate}),
together with $\left\langle \ell\right\rangle \leq1+\left|\ell\right|\leq2\left|\ell\right|$,
leads to
\begin{equation}
\left\langle \ell_{0}\right\rangle ^{\beta_{0}}=\left\langle \ell_{0}\right\rangle ^{\frac{\beta}{1-\beta}}\lesssim\left\langle \xi_{0}\right\rangle ^{\beta}\lesssim\left\langle \xi\right\rangle ^{\beta}\lesssim\left(\left\langle \ell\right\rangle ^{\frac{1}{1-\beta}}\right)^{\beta}\lesssim\left|\ell\right|^{\beta_{0}}\label{eq:AlphaModulationNormEstimateApplication}
\end{equation}
and
\[
\left|\ell\right|^{\beta_{0}}\leq\left(\left\langle \ell\right\rangle ^{\frac{1}{1-\beta}}\right)^{\beta}\lesssim\left\langle \xi\right\rangle ^{\beta}\lesssim\left\langle \xi_{0}\right\rangle ^{\beta}\lesssim\left(\left\langle \ell_{0}\right\rangle ^{\frac{1}{1-\beta}}\right)^{\beta}=\left\langle \ell_{0}\right\rangle ^{\beta_{0}}.
\]

For arbitrary $\eta\in Q_{s,\ell}^{\left(\beta\right)}$, a combination
of the estimates from the previous two paragraphs shows
\begin{align*}
\left|\eta-\left|\ell_{0}\right|^{\beta_{0}}\ell_{0}\right| & \leq\left|\eta-\ell^{\left(\beta\right)}\right|+\left|\ell^{\left(\beta\right)}-\xi\right|+\left|\xi-k^{\left(\alpha\right)}\right|+\left|k^{\left(\alpha\right)}-\xi_{0}\right|+\left|\xi_{0}-\ell_{0}^{\left(\beta\right)}\right|\\
 & <s\cdot\left|\ell\right|^{\beta_{0}}+s\cdot\left|\ell\right|^{\beta_{0}}+r\cdot\left|k\right|^{\alpha_{0}}+r\cdot\left|k\right|^{\alpha_{0}}+s\cdot\left|\ell_{0}\right|^{\beta_{0}}\\
 & \lesssim\left(s+r\right)\cdot\left\langle \ell_{0}\right\rangle ^{\beta_{0}}\,,
\end{align*}
and thus $Q_{s,\ell}^{\left(\beta\right)}\subset B_{C_{1}\left\langle \ell_{0}\right\rangle ^{\beta_{0}}}\left(\smash{\left|\ell_{0}\right|^{\beta_{0}}\ell_{0}}\right)$
for all $\ell\in J_{k}$ and some constant $C_{1}=C_{1}\left(\dimension,\alpha,\beta,r,s\right)$.

But the admissibility of $\CalQ_{s}^{\left(\beta\right)}$ yields
a constant $C_{2}=C_{2}\left(\dimension,\beta,s\right)>0$ with $\sum_{\ell\in\Z^{\dimension}\setminus\left\{ 0\right\} }\Indicator_{Q_{s,\ell}^{\left(\beta\right)}}\leq C_{2}$.
In fact, we can take $C_{2}=N_{\CalQ_{s}^{\left(\beta\right)}}$;
see the remark after Definition~\ref{defn:AdmissibleCoveringModerateWeight}.
In combination with the above inclusion, we conclude
\begin{equation}
\sum_{\ell\in J_{k}}\Indicator_{Q_{s,\ell}^{\left(\beta\right)}}\leq C_{2}\cdot\Indicator_{B_{C_{1}\left\langle \ell_{0}\right\rangle ^{\beta_{0}}}\left(\smash{\left|\ell_{0}\right|^{\beta_{0}}\ell_{0}}\right)}.\label{eq:AlphaCoveringSubordinatePointwiseEstimate}
\end{equation}
Furthermore, equation~(\ref{eq:AlphaModulationNormEstimateApplication})
yields 
\[
\lambda\left(\,\smash{Q_{s,\ell}^{\left(\beta\right)}}\,\right)=\lambda\left(B_{1}\left(0\right)\right)\cdot s^{\dimension}\cdot\left|\ell\right|^{\dimension\beta_{0}}\gtrsim s^{\dimension}\cdot\left\langle \ell_{0}\right\rangle ^{\dimension\beta_{0}}\qquad\forall\,\ell\in J_{k}\,.
\]
Thus, integration of estimate~(\ref{eq:AlphaCoveringSubordinatePointwiseEstimate})
leads to
\begin{align*}
s^{\dimension}\cdot\left|J_{k}\right|\cdot\left\langle \ell_{0}\right\rangle ^{\dimension\beta_{0}}\lesssim\sum_{\ell\in J_{k}}\lambda\left(\,\smash{Q_{s,\ell}^{\left(\beta\right)}}\,\right) & =\int_{\R^{\dimension}}\sum_{\ell\in J_{k}}\Indicator_{Q_{s,\ell}^{\left(\beta\right)}}\left(\xi\right)\,\d\xi\\
\left({\scriptstyle \text{eq. }\eqref{eq:AlphaCoveringSubordinatePointwiseEstimate}}\right) & \leq C_{2}\cdot\lambda\left(B_{C_{1}\left\langle \ell_{0}\right\rangle ^{\beta_{0}}}\left(\smash{\left|\ell_{0}\right|^{\beta_{0}}\ell_{0}}\right)\right)\lesssim C_{2}C_{1}^{\dimension}\cdot\left\langle \ell_{0}\right\rangle ^{\dimension\beta_{0}}\,.
\end{align*}
This establishes the upper estimate in equation~(\ref{eq:AlphaModulationIntersectionCountSubordinateness})
and thereby completes the proof.
\end{proof}
As we just saw, the $\alpha$-covering $\CalQ_{r}^{\left(\alpha\right)}$
is almost subordinate to the $\beta$-covering $\CalQ_{s}^{\left(\beta\right)}$,
for $0\leq\alpha\leq\beta<1$. Even more is true:
\begin{lem}
\label{lem:AlphaModulationSubordinatenessModerateness}Let $\dimension\in\N$
and $\alpha,\beta\in\left[0,1\right)$ with $\alpha\leq\beta$ be
arbitrary. Choose $r_{1}=r_{0}\left(\dimension,\alpha\right)$ and
$r_{2}=r_{0}\left(\dimension,\beta\right)$ with $r_{0}$ as in Theorem~\ref{thm:AlphaCoveringExistence},
and let $r>r_{1}$, as well as $s>r_{2}$.

Then $\CalQ_{r}^{\left(\alpha\right)}$ is almost subordinate to $\CalQ_{s}^{\left(\beta\right)}$
and relatively $\CalQ_{s}^{\left(\beta\right)}$-moderate. More precisely,
\begin{equation}
C^{-1}\cdot\left\langle \ell\right\rangle ^{\dimension\cdot\frac{\alpha}{1-\beta}}\leq\left|\det\left(\smash{T_{k}^{\left(\alpha\right)}}\right)\right|\leq C\cdot\left\langle \ell\right\rangle ^{\dimension\cdot\frac{\alpha}{1-\beta}}\qquad\forall\,k,\ell\in\Z^{\dimension}\setminus\left\{ 0\right\} \text{ with }Q_{r,k}^{\left(\alpha\right)}\cap Q_{s,\ell}^{\left(\beta\right)}\neq\emptyset\,,\label{eq:AlphaModulationRelativeModerateness}
\end{equation}
for a suitable constant $C=C\left(\dimension,\alpha,\beta,r,s\right)\geq1$.

Finally, let $\gamma\in\R$ be arbitrary, and let $w^{\left(\gamma\right)}$
as defined in Lemma~\ref{lem:AlphaModulationCoveringNormEstimate}.
Then, there is a constant $L=L\left(\dimension,\alpha,\beta,\gamma,r,s\right)\geq1$
with
\begin{equation}
L^{-1}\cdot w_{\ell}^{\left(\gamma\cdot\frac{1-\alpha}{1-\beta}\right)}\leq w_{k}^{\left(\gamma\right)}\leq L\cdot w_{\ell}^{\left(\gamma\cdot\frac{1-\alpha}{1-\beta}\right)}\qquad\forall\,k,\ell\in\Z^{\dimension}\text{ with }Q_{r,k}^{\left(\alpha\right)}\cap Q_{s,\ell}^{\left(\beta\right)}\neq\emptyset\,.\qedhere\label{eq:AlphaModulationSpacesWeightEstimate}
\end{equation}
\end{lem}

\begin{proof}
Lemma~\ref{lem:AlphaModulationIntersectionCountSubordinateness}
shows that $\CalQ_{r}^{\left(\alpha\right)}$ is almost subordinate
to $\CalQ_{s}^{\left(\beta\right)}$, so that we only have to establish
estimates (\ref{eq:AlphaModulationRelativeModerateness}) and (\ref{eq:AlphaModulationSpacesWeightEstimate}).

To this end, let $k,\ell\in\Z^{\dimension}\setminus\left\{ 0\right\} $
with $\emptyset\neq Q_{r,k}^{\left(\alpha\right)}\cap Q_{s,\ell}^{\left(\beta\right)}\ni\xi$.
Lemma~\ref{lem:AlphaModulationCoveringNormEstimate} shows $\left\langle k\right\rangle ^{\frac{1}{1-\alpha}}\asymp\left\langle \xi\right\rangle \asymp\left\langle \ell\right\rangle ^{\frac{1}{1-\beta}}$,
where the implied constants only depend on $\dimension,\alpha,\beta,r,s$.
If implied constants also depend on $\gamma$, this is indicated by
$\asymp_{\gamma}$ and $\lesssim_{\gamma}$. By what we just saw,
$\left\langle k\right\rangle ^{\frac{1}{1-\alpha}}/\left\langle \ell\right\rangle ^{\frac{1}{1-\beta}}\asymp1$,
and hence
\[
1\lesssim_{\gamma}\left(\left\langle k\right\rangle ^{\frac{1}{1-\alpha}}\,\big/\,\left\langle \ell\right\rangle ^{\frac{1}{1-\beta}}\right)^{\left(1-\alpha\right)\gamma}\lesssim_{\gamma}1,
\]
which means $w_{k}^{\left(\gamma\right)}=\left\langle k\right\rangle ^{\gamma}\asymp_{\gamma}\left\langle \ell\right\rangle ^{\gamma\cdot\frac{1-\alpha}{1-\beta}}=w_{\ell}^{\left(\gamma\cdot\frac{1-\alpha}{1-\beta}\right)}$.
Hence, equation~(\ref{eq:AlphaModulationSpacesWeightEstimate}) is
established.

Finally, the considerations from above (with $\gamma=\dimension\cdot\alpha_{0}$)
also show that $\CalQ_{r}^{\left(\alpha\right)}=\left(\smash{T_{k}^{\left(\alpha\right)}Q^{\left(r\right)}+b_{k}}\right)_{k\in\Z^{\dimension}\setminus\left\{ 0\right\} }$
is relatively moderate with respect to $\CalQ_{s}^{\left(\beta\right)}$,
because of
\[
\left|\det\left(\smash{T_{k}^{\left(\alpha\right)}}\right)\right|=\left|\det\left(\left|k\right|^{\alpha_{0}}{\rm id}\right)\right|=\left|k\right|^{\dimension\cdot\alpha_{0}}\asymp\left\langle k\right\rangle ^{\dimension\cdot\alpha_{0}}=w_{k}^{\left(\dimension\cdot\alpha_{0}\right)}\asymp\left\langle \ell\right\rangle ^{\dimension\cdot\alpha_{0}\cdot\frac{1-\alpha}{1-\beta}}=\left\langle \ell\right\rangle ^{\dimension\cdot\frac{\alpha}{1-\beta}}
\]
for all $k\in\Z^{\dimension}\setminus\left\{ 0\right\} $ with $Q_{r,k}^{\left(\alpha\right)}\cap Q_{s,\ell}^{\left(\beta\right)}\neq\emptyset$.
This also establishes equation~(\ref{eq:AlphaModulationRelativeModerateness}).
\end{proof}
After the above analysis of the geometric relationship between the
different $\alpha$-coverings, it is now essentially straightforward
to derive the announced characterization of the existence of embeddings
between $\alpha$-modulation spaces with different values of $\alpha$.

The ``limit case'' $\alpha=1$—which corresponds to (inhomogeneous)
Besov spaces—will be treated below; see Theorem~\ref{thm:AlphaModulationBesovEmbedding}.
Finally, we remark that the present theorem is essentially identical
to \cite[Theorem 6.1.7]{VoigtlaenderPhDThesis} from my PhD thesis.
\begin{thm}
\label{thm:AlphaModulationEmbeddings}Let $\dimension\in\N$ and $\alpha,\beta\in\left[0,1\right)$
with $\alpha\leq\beta$. Finally, let $p_{1},p_{2},q_{1},q_{2}\in\left(0,\infty\right]$,
as well as $\gamma_{1},\gamma_{2}\in\R$. Define
\begin{align*}
\gamma^{\left(0\right)} & :=\alpha\left(\frac{1}{p_{2}}-\frac{1}{p_{1}}\right)+\left(\alpha-\beta\right)\left(\smash{\frac{1}{\LowerExpo{p_{2}}}}-\frac{1}{q_{1}}\right)_{+},\\
\gamma^{\left(1\right)} & :=\alpha\left(\frac{1}{p_{2}}-\frac{1}{p_{1}}\right)+\left(\alpha-\beta\right)\left(\frac{1}{q_{2}}-\smash{\frac{1}{\SignedUpperExpo{p_{1}}}}\right)_{\!+}.
\end{align*}

We have $\AlphaModSpace{p_{1}}{q_{1}}{\Fourier,\gamma_{1}}{\alpha}\left(\R^{\dimension}\right)\subset\AlphaModSpace{p_{2}}{q_{2}}{\Fourier,\gamma_{2}}{\beta}\left(\R^{\dimension}\right)$
if and only if\footnote{This equivalence uses that both (Fourier-side) $\alpha$-modulation
spaces embed into the Hausdorff space $\DistributionSpace{\R^{\dimension}}$
(see Theorem~\ref{thm:DecompositionSpaceComplete}), and that the
closed graph theorem is also valid for quasi-Banach spaces, as seen
in Section~\ref{subsec:Notation}.} $\AlphaModSpace{p_{1}}{q_{1}}{\Fourier,\gamma_{1}}{\alpha}\left(\R^{\dimension}\right)\hookrightarrow\AlphaModSpace{p_{2}}{q_{2}}{\Fourier,\gamma_{2}}{\beta}\left(\R^{\dimension}\right)$,
if and only if we have
\[
p_{1}\leq p_{2}\quad\text{and}\quad\begin{cases}
\gamma_{2}\leq\gamma_{1}+\dimension\cdot\gamma^{\left(0\right)}, & \text{if }q_{1}\leq q_{2},\\
\gamma_{2}<\gamma_{1}+\dimension\cdot\left(\gamma^{\left(0\right)}+\left(1-\beta\right)\left(q_{1}^{-1}-q_{2}^{-1}\right)\right), & \text{if }q_{1}>q_{2}.
\end{cases}
\]

Conversely, we have $\AlphaModSpace{p_{1}}{q_{1}}{\Fourier,\gamma_{1}}{\beta}\left(\R^{\dimension}\right)\subset\AlphaModSpace{p_{2}}{q_{2}}{\Fourier,\gamma_{2}}{\alpha}\left(\R^{\dimension}\right)$
if and only if $\AlphaModSpace{p_{1}}{q_{1}}{\Fourier,\gamma_{1}}{\beta}\left(\R^{\dimension}\right)\hookrightarrow\AlphaModSpace{p_{2}}{q_{2}}{\Fourier,\gamma_{2}}{\alpha}\left(\R^{\dimension}\right)$,
if and only if we have
\[
p_{1}\leq p_{2}\quad\text{and}\quad\begin{cases}
\gamma_{2}\leq\gamma_{1}+\dimension\cdot\gamma^{\left(1\right)}, & \text{if }q_{1}\leq q_{2},\\
\gamma_{2}<\gamma_{1}+\dimension\cdot\left(\gamma^{\left(1\right)}+\left(1-\beta\right)\left(q_{1}^{-1}-q_{2}^{-1}\right)\right), & \text{if }q_{1}>q_{2}.
\end{cases}\qedhere
\]
\end{thm}

\begin{rem}
\label{rem:AlphaModulationWellDefined}In case of $q_{1}=q_{2}=q$
and $p_{1}=p_{2}=p$, the first condition reduces to
\[
\gamma_{2}+\dimension\left(\beta-\alpha\right)\left(\smash{\frac{1}{\LowerExpo p}}-\frac{1}{q}\right)_{+}\leq\gamma_{1},
\]
which is identical to the condition given by Han and Wang in \cite[Theorem 4.1]{HanWangAlphaModulationEmbeddings}.
Analogously, the second condition reduces to
\[
\gamma_{2}+\dimension\left(\beta-\alpha\right)\left(\frac{1}{q}-\smash{\frac{1}{\SignedUpperExpo p}}\right)_{\!+}\leq\gamma_{1},
\]
which also coincides with the condition stated in \cite[Theorem 4.1]{HanWangAlphaModulationEmbeddings}.
The verification that the stated conditions are equivalent to the
ones given by Han and Wang is straightforward, but slightly tedious,
since Han and Wang use a different notation.

Finally, as already mentioned at the beginning of this section, we
remark that the result obtained here is identical to the one in \cite[Theorem 1.2]{GuoAlphaModulationEmbeddingCharacterization},
which was obtained independently and the preprint version of which
appeared on the arXiv only a few days after the first version of the
present paper appeared there. The result appeared even earlier, however,
in my PhD thesis \cite[Theorem 6.1.7]{VoigtlaenderPhDThesis}.
\end{rem}

\begin{proof}
Let $r_{1}=r_{0}\left(\dimension,\alpha\right)$ and $r_{2}=r_{0}\left(\dimension,\beta\right)$,
with $r_{0}$ as in Theorem~\ref{thm:AlphaCoveringExistence}; choose
$r>r_{1}$, as well as $s>r_{2}$ and recall the definition of the
spaces $\AlphaModSpace pq{\Fourier,\gamma}{\alpha}\left(\R^{\dimension}\right)=\FourierDecompSp{\smash{\CalQ_{r}^{\left(\alpha\right)}}}p{\smash{\ell_{w^{\left(\gamma/\left(1-\alpha\right)\right)}}^{q}}}\vphantom{\CalQ_{r}^{\left(\alpha\right)}}$
from Definition~\ref{def:AlphaModulationFourierSpace}.

We first analyze for which parameters the inclusion $\AlphaModSpace{p_{1}}{q_{1}}{\Fourier,\gamma_{1}}{\alpha}\left(\R^{\dimension}\right)\subset\AlphaModSpace{p_{2}}{q_{2}}{\Fourier,\gamma_{2}}{\beta}\left(\R^{\dimension}\right)$
holds. By the closed graph theorem, this holds iff the identity map
$\iota:\AlphaModSpace{p_{1}}{q_{1}}{\Fourier,\gamma_{1}}{\alpha}\left(\R^{\dimension}\right)\to\AlphaModSpace{p_{2}}{q_{2}}{\Fourier,\gamma_{2}}{\beta}\left(\R^{\dimension}\right),f\mapsto f$
is well-defined and bounded. To characterize when this is the case,
we define $\alpha_{0}:=\frac{\alpha}{1-\alpha}$ and $\beta_{0}:=\frac{\beta}{1-\beta}$,
and set
\[
\CalQ=\left(Q_{i}\right)_{i\in I}=\left(T_{i}Q_{i}'+b_{i}\right)_{i\in I}:=\CalQ_{r}^{\left(\alpha\right)}=\left(Q_{r,k}^{\left(\alpha\right)}\right)_{k\in\Z^{\dimension}\setminus\left\{ 0\right\} }=\left(\left(\left|k\right|^{\alpha_{0}}\cdot\identity\right)Q^{\left(r\right)}+\left|k\right|^{\alpha_{0}}k\right)_{k\in\Z^{\dimension}\setminus\left\{ 0\right\} }\,,
\]
and likewise
\[
\CalP=\left(P_{j}\right)_{j\in J}=\left(S_{j}P_{j}'+c_{j}\right)_{j\in J}:=\CalQ_{s}^{\left(\beta\right)}=\left(Q_{s,k}^{\left(\beta\right)}\right)_{k\in\Z^{\dimension}\setminus\left\{ 0\right\} }=\left(\left(\left|k\right|^{\beta_{0}}\cdot\identity\right)Q^{\left(s\right)}+\left|k\right|^{\beta_{0}}k\right)_{k\in\Z^{\dimension}\setminus\left\{ 0\right\} }\,;
\]
see Lemma~\ref{lem:AlphaModulationCoveringIsStructured} for the
notations $Q^{\left(r\right)}$ and $Q^{\left(s\right)}$. Furthermore,
with $w^{\left(\gamma\right)}$ as in Lemma~\ref{lem:AlphaModulationCoveringNormEstimate},
we set $w:=w^{\left(\gamma_{1}/\left(1-\alpha\right)\right)}$ and
$v:=w^{\left(\gamma_{2}/\left(1-\beta\right)\right)}$. These choices
ensure that $\AlphaModSpace{p_{1}}{q_{1}}{\Fourier,\gamma_{1}}{\alpha}\left(\R^{\dimension}\right)=\FourierDecompSp{\CalQ}{p_{1}}{\ell_{w}^{q_{1}}}$
and $\AlphaModSpace{p_{2}}{q_{2}}{\Fourier,\gamma_{2}}{\beta}\left(\R^{\dimension}\right)=\FourierDecompSp{\CalP}{p_{2}}{\ell_{v}^{q_{2}}}$.
Lemma~\ref{lem:AlphaModulationCoveringIsStructured} shows that $\CalQ,\CalP$
are structured admissible coverings of $\R^{\dimension}$, which thus
satisfy the standing assumptions from Section~\ref{sec:SummaryOfEmbeddingResults}
(i.e., Assumption~\ref{assu:GeneralSummaryAssumptions}). Furthermore,
$w$ is $\CalQ$-moderate and $v$ is $\CalP$-moderate, thanks to
Lemma~\ref{lem:AlphaModulationCoveringNormEstimate}, so that Assumption~\ref{assu:GeneralSummaryAssumptions}
is completely satisfied.

Finally, since we have $\alpha\leq\beta$, Lemma~\ref{lem:AlphaModulationSubordinatenessModerateness}
shows that $\CalQ$ is almost subordinate to $\CalP$ and that $\CalQ$
and $w$ are relatively $\CalP$-moderate. Thus, all assumptions of
part~(\ref{enu:SummaryFineInCoarseModerate}) of Theorem~\ref{thm:SummaryFineIntoCoarse}
are satisfied. Note that the embedding $\iota$ from Theorem~\ref{thm:SummaryFineIntoCoarse}
satisfies $\iota f=f$ for all $f\in\FourierDecompSp{\CalQ}{p_{1}}{\ell_{w}^{q_{1}}}$,
since $\CalQ,\CalP$ both cover the same set $\CalO=\R^{\dimension}=\CalO'$.
Thus, all in all, we conclude that $\iota$ is well-defined and bounded
if and only if we have $p_{1}\leq p_{2}$ and if
\[
K:=\left\Vert \!\left(\!\frac{v_{j}}{w_{i_{j}}}\cdot\left|\det\smash{T_{i_{j}}}\right|^{p_{1}^{-1}-p_{2}^{-1}-s}\cdot\left|\det S_{j}\right|^{s}\right)_{\!\!\!j\in J}\right\Vert _{\ell^{q_{2}\cdot\left(q_{1}/q_{2}\right)'}}
\]
is finite, where $s=\left(\smash{\frac{1}{\LowerExpo{p_{2}}}}-\frac{1}{q_{1}}\right)_{+}$
and where for each $j\in J=\Z^{\dimension}\setminus\left\{ 0\right\} $,
some $i_{j}\in I_{j}$, i.e.\@ with $Q_{r,i_{j}}^{\left(\alpha\right)}\cap Q_{s,j}^{\left(\beta\right)}\neq\emptyset$
is selected. According to Lemma~\ref{lem:AlphaModulationSubordinatenessModerateness},
we have
\[
\left|\det\smash{T_{i_{j}}}\right|\asymp\left\langle j\right\rangle ^{\dimension\frac{\alpha}{1-\beta}}\qquad\text{ and }\qquad w_{i_{j}}=w_{i_{j}}^{\left(\gamma_{1}/\left(1-\alpha\right)\right)}\asymp w_{j}^{\left(\frac{\gamma_{1}}{1-\alpha}\cdot\frac{1-\alpha}{1-\beta}\right)}=w_{j}^{\left(\gamma_{1}/\left(1-\beta\right)\right)}=\left\langle j\right\rangle ^{\frac{\gamma_{1}}{1-\beta}}.
\]
Since we also have $\left|\det S_{j}\right|=\left|j\right|^{\dimension\frac{\beta}{1-\beta}}$
and $v_{j}=w_{j}^{\left(\gamma_{2}/\left(1-\beta\right)\right)}=\left\langle j\right\rangle ^{\frac{\gamma_{2}}{1-\beta}}$,
we get
\begin{align*}
K & \asymp\left\Vert \left(\left\langle j\right\rangle ^{\frac{\gamma_{2}-\gamma_{1}}{1-\beta}}\cdot\left\langle j\right\rangle ^{\dimension\frac{\alpha}{1-\beta}\left(p_{1}^{-1}-p_{2}^{-1}-s\right)}\cdot\left|j\right|^{\dimension s\frac{\beta}{1-\beta}}\right)_{j\in\Z^{\dimension}\setminus\left\{ 0\right\} }\right\Vert _{\ell^{q_{2}\cdot\left(q_{1}/q_{2}\right)'}}\\
 & \asymp\left\Vert \left(\left\langle j\right\rangle ^{\frac{1}{1-\beta}\left[\left(\gamma_{2}-\gamma_{1}\right)+\dimension\alpha\left(p_{1}^{-1}-p_{2}^{-1}-s\right)+\dimension\beta s\right]}\right)_{j\in\Z^{\dimension}\setminus\left\{ 0\right\} }\right\Vert _{\ell^{q_{2}\cdot\left(q_{1}/q_{2}\right)'}}.
\end{align*}

Recall from equation~(\ref{eq:SpecialExponentFiniteness}) that $q_{2}\cdot\left(q_{1}/q_{2}\right)'$
is finite if and only if $q_{2}<q_{1}$. Hence, there are two cases:

\textbf{Case 1}: $q_{1}\leq q_{2}$. Here, we have $q_{2}\cdot\left(q_{1}/q_{2}\right)'=\infty$,
so that we get the following equivalence:
\begin{align*}
K<\infty & \Longleftrightarrow\frac{1}{1-\beta}\left[\left(\gamma_{2}-\gamma_{1}\right)+\dimension\alpha\left(p_{1}^{-1}-p_{2}^{-1}-s\right)+\dimension\beta s\right]\leq0\\
 & \Longleftrightarrow\gamma_{2}+\dimension\left[s\left(\beta-\alpha\right)+\alpha\left(p_{1}^{-1}-p_{2}^{-1}\right)\right]\leq\gamma_{1}\\
 & \Longleftrightarrow\gamma_{2}\leq\gamma_{1}+\dimension\cdot\gamma^{\left(0\right)}.
\end{align*}
This completes the proof for $q_{1}\leq q_{2}$.

\textbf{Case 2}: $q_{1}>q_{2}$. Here, we have $t:=q_{2}\cdot\left(q_{1}/q_{2}\right)'<\infty$.
But for arbitrary $\varrho\in\R$ and $t\in\left(0,\infty\right)$,
it is an elementary fact that
\[
\left\Vert \left(\left\langle j\right\rangle ^{\varrho}\right)_{j\in\Z^{\dimension}\setminus\left\{ 0\right\} }\right\Vert _{\ell^{t}}<\infty\qquad\Longleftrightarrow\qquad t\cdot\varrho<-\dimension\qquad\Longleftrightarrow\qquad\varrho<-\frac{\dimension}{t}.
\]
In the present setting, equation~(\ref{eq:InverseOfSpecialExponent})
yields $\frac{1}{t}=\left(\frac{1}{q_{2}}-\frac{1}{q_{1}}\right)_{+}=\frac{1}{q_{2}}-\frac{1}{q_{1}}$.
Furthermore,
\begin{align*}
\varrho & =\frac{1}{1-\beta}\left[\left(\gamma_{2}-\gamma_{1}\right)+\dimension\alpha\left(p_{1}^{-1}-p_{2}^{-1}-s\right)+\dimension\beta s\right]\\
 & =\frac{1}{1-\beta}\left[\left(\gamma_{2}-\gamma_{1}\right)+\dimension\left(s\left(\beta-\alpha\right)+\alpha\left(p_{1}^{-1}-p_{2}^{-1}\right)\right)\right],
\end{align*}
so that we get all in all that $K<\infty$ is equivalent to
\begin{align*}
 & \frac{1}{1-\beta}\left[\left(\gamma_{2}-\gamma_{1}\right)+\dimension\left(s\left(\beta-\alpha\right)+\alpha\left(p_{1}^{-1}-p_{2}^{-1}\right)\right)\right]\overset{!}{<}-\dimension\left(q_{2}^{-1}-q_{1}^{-1}\right)\\
\Longleftrightarrow\quad & \left(\gamma_{2}-\gamma_{1}\right)+\dimension\left(s\left(\beta-\alpha\right)+\alpha\left(p_{1}^{-1}-p_{2}^{-1}\right)\right)\overset{!}{<}\dimension\left(1-\beta\right)\left(q_{1}^{-1}-q_{2}^{-1}\right)\\
\Longleftrightarrow\quad & \gamma_{2}\overset{!}{<}\gamma_{1}+\dimension\left(s\left(\alpha-\beta\right)+\alpha\left(p_{2}^{-1}-p_{1}^{-1}\right)\right)+\dimension\left(1-\beta\right)\left(q_{1}^{-1}-q_{2}^{-1}\right)\\
\Longleftrightarrow\quad & \gamma_{2}\overset{!}{<}\gamma_{1}+\dimension\left(\gamma^{\left(0\right)}+\left(1-\beta\right)\left(q_{1}^{-1}-q_{2}^{-1}\right)\right).
\end{align*}
This completes the proof for $q_{1}>q_{2}$.

\medskip{}

Now, we analyze for which parameters the inclusion $\AlphaModSpace{p_{1}}{q_{1}}{\Fourier,\gamma_{1}}{\beta}\left(\R^{\dimension}\right)\subset\AlphaModSpace{p_{2}}{q_{2}}{\Fourier,\gamma_{2}}{\alpha}\left(\R^{\dimension}\right)$
holds. By the closed graph theorem, this holds iff the identity map
$\theta:\AlphaModSpace{p_{1}}{q_{1}}{\Fourier,\gamma_{1}}{\beta}\left(\R^{\dimension}\right)\to\AlphaModSpace{p_{2}}{q_{2}}{\Fourier,\gamma_{2}}{\alpha}\left(\R^{\dimension}\right),f\mapsto f$
is well-defined and bounded. To characterize when this is the case,
we define $\alpha_{0}:=\frac{\alpha}{1-\alpha}$ and $\beta_{0}:=\frac{\beta}{1-\beta}$,
and set
\[
\CalQ=\left(Q_{i}\right)_{i\in I}=\left(T_{i}Q_{i}'+b_{i}\right)_{i\in I}:=\CalQ_{s}^{\left(\beta\right)}=\left(Q_{s,k}^{\left(\beta\right)}\right)_{k\in\Z^{\dimension}\setminus\left\{ 0\right\} }=\left(\left(\left|k\right|^{\beta_{0}}\cdot\identity\right)Q^{\left(s\right)}+\left|k\right|^{\beta_{0}}k\right)_{k\in\Z^{\dimension}\setminus\left\{ 0\right\} }\,,
\]
and likewise
\[
\CalP=\left(P_{j}\right)_{j\in J}=\left(S_{j}P_{j}'+c_{j}\right)_{j\in J}:=\CalQ_{r}^{\left(\alpha\right)}=\left(Q_{r,k}^{\left(\alpha\right)}\right)_{k\in\Z^{\dimension}\setminus\left\{ 0\right\} }=\left(\left(\left|k\right|^{\alpha_{0}}\cdot\identity\right)Q^{\left(r\right)}+\left|k\right|^{\alpha_{0}}k\right)_{k\in\Z^{\dimension}\setminus\left\{ 0\right\} }\,,
\]
where again $Q^{\left(r\right)}$ and $Q^{\left(s\right)}$ are as
in Lemma~\ref{lem:AlphaModulationCoveringIsStructured}. Further,
we set $w:=w^{\left(\gamma_{1}/\left(1-\beta\right)\right)}$ and
$v:=w^{\left(\gamma_{2}/\left(1-\alpha\right)\right)}$, so that we
have $\AlphaModSpace{p_{1}}{q_{1}}{\Fourier,\gamma_{1}}{\beta}\left(\R^{\dimension}\right)=\FourierDecompSp{\CalQ}{p_{1}}{\ell_{w}^{q_{1}}}$
and $\AlphaModSpace{p_{2}}{q_{2}}{\Fourier,\gamma_{2}}{\alpha}\left(\R^{\dimension}\right)=\FourierDecompSp{\CalP}{p_{2}}{\ell_{v}^{q_{2}}}$.

Exactly as above, we see that Assumption~\ref{assu:GeneralSummaryAssumptions}
is fulfilled. Furthermore, since $\alpha\leq\beta$, Lemma~\ref{lem:AlphaModulationSubordinatenessModerateness}
shows that $\CalP=\CalQ_{r}^{\left(\alpha\right)}$ is almost subordinate
to $\CalQ=\CalQ_{s}^{\left(\beta\right)}$. The same lemma also shows
that $\CalP$ and $v$ are relatively $\CalQ$-moderate, so that all
assumptions of part~(\ref{enu:SummaryCoarseInFineModerate}) of Theorem~\ref{thm:SummaryCoarseIntoFine}
are satisfied.

As above, we thus see that $\theta$ is well-defined and bounded if
and only if we have $p_{1}\leq p_{2}$ and if
\[
K:=\left\Vert \left(\frac{v_{j_{i}}}{w_{i}}\cdot\left|\det T_{i}\right|^{s}\cdot\left|\det\smash{S_{j_{i}}}\right|^{p_{1}^{-1}-p_{2}^{-1}-s}\right)_{i\in I}\right\Vert _{\ell^{q_{2}\cdot\left(q_{1}/q_{2}\right)'}}<\infty
\]
holds, where $s=\left(\frac{1}{q_{2}}-\smash{\frac{1}{\SignedUpperExpo{p_{1}}}}\right)_{+}$
and where for each index $i\in I=\Z^{\dimension}\setminus\left\{ 0\right\} $,
some $j_{i}\in J_{i}$, i.e.\@ with $Q_{r,j_{i}}^{\left(\alpha\right)}\cap Q_{s,i}^{\left(\beta\right)}\neq\emptyset$
is selected. But in view of Lemma~\ref{lem:AlphaModulationSubordinatenessModerateness},
we have
\[
\left|\det\smash{S_{j_{i}}}\right|\asymp\left\langle i\right\rangle ^{\dimension\frac{\alpha}{1-\beta}}\qquad\text{ and }\qquad v_{j_{i}}=w_{j_{i}}^{\left(\gamma_{2}/\left(1-\alpha\right)\right)}\asymp w_{i}^{\left(\frac{\gamma_{2}}{1-\alpha}\cdot\frac{1-\alpha}{1-\beta}\right)}=w_{i}^{\left(\gamma_{2}/\left(1-\beta\right)\right)}=\left\langle i\right\rangle ^{\frac{\gamma_{2}}{1-\beta}},
\]
so that we get
\begin{align*}
K & \asymp\left\Vert \left(\left\langle i\right\rangle ^{\frac{\gamma_{2}-\gamma_{1}}{1-\beta}}\cdot\left\langle i\right\rangle ^{\dimension\frac{\alpha}{1-\beta}\left(p_{1}^{-1}-p_{2}^{-1}-s\right)}\cdot\left|i\right|^{\dimension s\frac{\beta}{1-\beta}}\right)_{i\in\Z^{\dimension}\setminus\left\{ 0\right\} }\right\Vert _{\ell^{q_{2}\cdot\left(q_{1}/q_{2}\right)'}}\\
 & \asymp\left\Vert \left(\left\langle i\right\rangle ^{\frac{1}{1-\beta}\left[\left(\gamma_{2}-\gamma_{1}\right)+\dimension\left(s\left(\beta-\alpha\right)+\alpha\left(p_{1}^{-1}-p_{2}^{-1}\right)\right)\right]}\right)_{i\in\Z^{\dimension}\setminus\left\{ 0\right\} }\right\Vert _{\ell^{q_{2}\cdot\left(q_{1}/q_{2}\right)'}}.
\end{align*}
The remainder of the proof is entirely analogous to the considerations
for the other inclusion, and is hence omitted.
\end{proof}
In summary, we have shown that the embedding results for $\alpha$-modulation
spaces obtained in \cite[Theorem 4.1]{HanWangAlphaModulationEmbeddings}
can be derived with ease using our more general approach. We remark
that our approach can also handle the case $\left(p_{1},q_{1}\right)\neq\left(p_{2},q_{2}\right)$
which is not covered by \cite[Theorem 4.1]{HanWangAlphaModulationEmbeddings}
(but which was independently considered in \cite{GuoAlphaModulationEmbeddingCharacterization}).
Furthermore, the most tedious parts of our derivation—namely the proof
that $\CalQ_{r}^{\left(\alpha\right)}$ is almost subordinate to $\CalQ_{s}^{\left(\beta\right)}$
for $\alpha\leq\beta$ and the proof of the $\CalQ_{s}^{\left(\beta\right)}$-relative
moderateness of $\CalQ_{r}^{\left(\alpha\right)}$ and $w^{\left(\gamma\right)}$—are
used in \cite{HanWangAlphaModulationEmbeddings,GuoAlphaModulationEmbeddingCharacterization}
without much justification; see for example \cite[equations (4.4) and (4.6)]{HanWangAlphaModulationEmbeddings}
and \cite[equations (4.12) and (4.29)]{GuoAlphaModulationEmbeddingCharacterization}.

\medskip{}

In the next subsection, we will begin our study of Besov spaces. Precisely,
we will extend Theorem~\ref{thm:AlphaModulationEmbeddings} to include
the case $\beta=1$, which shows that the results for embeddings between
$\alpha$-modulation spaces and inhomogeneous Besov spaces obtained
in \cite[Theorem 4.2]{HanWangAlphaModulationEmbeddings} are also
special cases of our approach; in fact our approach extends these
results to the case $\left(p_{1},q_{1}\right)\neq\left(p_{2},q_{2}\right)$.

\subsection{Embeddings between \texorpdfstring{$\alpha$-modulation}{α-modulation}
spaces and Besov spaces}

\label{subsec:AlphaModulationBesov}It is common to consider the (inhomogeneous)
Besov spaces $\BesovInhom pq{\gamma}\left(\R^{\dimension}\right)$
as the ``limit case'' of $\alpha$-modulation spaces for $\alpha\uparrow1$,
i.e.\@ to set $\AlphaModSpace pq{\gamma}1\left(\R^{\dimension}\right):=B_{\gamma}^{p,q}\left(\R^{\dimension}\right)$.
In the present subsection, we want to extend Theorem~\ref{thm:AlphaModulationEmbeddings}
to include the case $\beta=1$.

As a further result, we will show in the present subsection that the
Fourier transform restricts to an isomorphism of quasi-Banach spaces
\[
\Fourier:\AlphaModSpace pq{\gamma}{\alpha}\left(\R^{\dimension}\right)\to\FourierDecompSp{\CalQ_{r}^{\left(\alpha\right)}}p{\ell_{w^{\left(\gamma/\left(1-\alpha\right)\right)}}^{q}}=\AlphaModSpace pq{\Fourier,\gamma}{\alpha}\left(\R^{\dimension}\right),
\]
thereby justifying the name ``Fourier-side $\alpha$-modulation spaces''
for the spaces $\AlphaModSpace pq{\Fourier,\gamma}{\alpha}\left(\R^{\dimension}\right)$.
As a consequence, the characterization of the embeddings for the Fourier-side
$\alpha$-modulation spaces from Theorem~\ref{thm:AlphaModulationEmbeddings}
readily extends to the ``usual'' (space-side) $\alpha$-modulation
spaces.

\smallskip{}

But let us begin by defining the inhomogeneous Besov covering:
\begin{defn}
\label{def:InhomogeneousBesovCovering}For $\dimension\in\N$, define
\[
B_{0}:=B_{4}\left(0\right)\subset\R^{\dimension}\qquad\text{and}\qquad B_{n}:=B_{2^{n+2}}\left(0\right)\setminus\overline{B_{2^{n-2}}}\left(0\right)\subset\R^{\dimension}\qquad\text{for }n\in\N\,.
\]
The family $\CalB:=\left(B_{n}\right)_{n\in\N_{0}}$ is called the
\textbf{inhomogeneous Besov covering} of $\R^{\dimension}$.
\end{defn}

The following lemma (partly) justifies our nomenclature, by showing
that $\CalB$ is indeed an (almost structured) covering of $\R^{\dimension}$.
The full justification of our nomenclature will be given in Lemma~\ref{lem:FourierInhomogeneousBesovAsTemperedDistributions}.
\begin{lem}
\label{lem:InhomogeneousBesovCoveringAlmostStructured}The inhomogeneous
Besov covering $\CalB=\left(T_{n}B_{n}'+b_{n}\right)_{n\in\N_{0}}$
is an almost structured covering of $\R^{\dimension}$, with parametrization
given by 
\[
T_{n}:=2^{n}\cdot\identity,\qquad b_{n}:=0\qquad\text{ and }\qquad B_{n}':=\begin{cases}
B_{4}\left(0\right)\setminus\overline{B_{1/4}}\left(0\right), & \text{if }n\in\N,\\
B_{4}\left(0\right), & \text{if }n=0.
\end{cases}
\]

Furthermore, for any $\gamma\in\R$, the weight $v^{\left(\gamma\right)}:=\left(2^{\gamma n}\right)_{n\in\N_{0}}$
is $\CalB$-moderate.
\end{lem}

\begin{proof}
The given $T_{n},B_{n}'$ and $b_{n}$ indeed satisfy $B_{n}=T_{n}B_{n}'+b_{n}$
for all $n\in\N_{0}$. Furthermore, if we set
\[
P_{n}:=\begin{cases}
B_{2}\left(0\right)\setminus\overline{B_{1/2}}\left(0\right), & \text{if }n\in\N,\\
B_{2}\left(0\right), & \text{if }n=0,
\end{cases}
\]
then all $P_{n}$ and all $B_{n}'$ are open and bounded, the sets
$\left\{ P_{n}\with n\in\N_{0}\right\} $ and $\left\{ B_{n}'\with n\in\N_{0}\right\} $
are finite, and we have $\overline{P_{n}}\subset B_{n}'$ for all
$n\in\N_{0}$.

Further, $\R^{\dimension}=\bigcup_{n\in\N_{0}}\left(T_{n}P_{n}+b_{n}\right)$.
Indeed, for $\xi\in\R^{\dimension}$ with $\left|\xi\right|\geq1$,
choose $n\in\N_{0}$ maximal with $2^{n}\leq\left|\xi\right|$. This
implies $2^{n}\leq\left|\xi\right|<2^{n+1}$. If $n=0$, this yields
$\left|\xi\right|<2$ and hence $\xi\in B_{2}\left(0\right)=T_{0}P_{0}+b_{0}$.
Otherwise, for $n\in\N$, we have $2^{n-1}<2^{n}\leq\left|\xi\right|<2^{n+1}$,
i.e.\@ $\xi\in2^{n}\cdot\left(B_{2}\left(0\right)\setminus\overline{B_{1/2}}\left(0\right)\right)=T_{n}P_{n}+b_{n}$,
as desired. Finally, in case of $\left|\xi\right|<1$, we clearly
have $\xi\in P_{0}=T_{0}P_{0}+b_{0}$.

Thus, recalling the definition of an almost structured covering (Definition~\ref{defn:DifferentTypesOfCoverings}),
it remains to show that $\CalB$ is admissible and that $C_{\CalB}$
is finite. To this end, fix $n\in\N$ and assume that there is some
$m\in\N$ with $\emptyset\neq B_{n}\cap B_{m}$. If we choose $\xi\in B_{n}\cap B_{m}$,
we thus get
\[
2^{m-2}<\left|\xi\right|<2^{n+2}\qquad\text{ and }\qquad2^{n-2}<\left|\xi\right|<2^{m+2}.
\]
In particular, $m-2<n+2$ and $n-2<m+2$, i.e.\@ $n-4<m<n+4$. Since
all these quantities are integers, we get $n-3\leq m\leq n+3$. All
in all, we derive
\begin{equation}
n^{\ast}\subset\left\{ 0\right\} \cup\left\{ n-3,\dots,n+3\right\} \qquad\forall\,n\in\N\,,\label{eq:InhomogeneousBesovClusterEstimate}
\end{equation}
and hence $\left|n^{\ast}\right|\leq8$. Finally, for $m\in\N$ with
$\xi\in B_{0}\cap B_{m}\neq\emptyset$, we get $2^{m-2}<\left|\xi\right|<4=2^{2}$
and hence $m-2<2$, i.e.\@ $m\leq3$. We have thus shown $0^{\ast}\subset\left\{ 0,1,2,3\right\} $
and hence $\left|0^{\ast}\right|\leq4$. All in all, we get $N_{\CalB}=\sup_{n\in\N_{0}}\left|n^{\ast}\right|\leq8<\infty$,
so that $\CalB$ is an admissible covering of $\R^{\dimension}$.

It remains to show that 
\[
C_{\CalB}:=\sup_{n\in\N_{0}}\sup_{m\in n^{\ast}}\left\Vert T_{n}^{-1}T_{m}\right\Vert =\sup_{n\in\N_{0}}\sup_{m\in n^{\ast}}2^{m-n}
\]
is finite. This is indeed the case: For $n\in\N$, we saw in equation~(\ref{eq:InhomogeneousBesovClusterEstimate})
that $n^{\ast}\subset\left\{ 0,\dots,n+3\right\} $, so that every
$m\in n^{\ast}$ satisfies $2^{m-n}\leq2^{3}=8$. Finally, for $n=0$
and $m\in n^{\ast}=0^{\ast}\subset\left\{ 0,1,2,3\right\} $, we also
have $2^{m-n}\leq2^{3}=8$, so that $C_{\CalB}\leq8<\infty$. All
in all, we have thus shown that $\CalB$ is an almost structured covering
of $\R^{\dimension}$.

\medskip{}

As an almost structured covering, $\CalB$ is in particular a semi-structured
covering, so that equation~(\ref{eq:DeterminantIsModerate}) shows
that the weight $\left(\left|\det T_{n}\right|\right)_{n\in\N_{0}}=\left(2^{n\dimension}\right)_{n\in\N_{0}}$
is $\CalB$-moderate. But if $w=\left(w_{n}\right)_{n\in\N_{0}}$
is $\CalB$-moderate, it is not hard to see that the same holds for
$\left(w_{n}^{\varrho}\right)_{n\in\N_{0}}$, for arbitrary $\varrho\in\R$.
If we apply this with $\varrho=\frac{\gamma}{\dimension}$, we see
that $v^{\left(\gamma\right)}$ is $\CalB$-moderate, as claimed.
\end{proof}
Recall from Theorem~\ref{thm:AlmostStructuredAdmissibleAdmitsBAPU}
that every almost structured covering admits a subordinate partition
of unity which is an $L^{p}$-BAPU, simultaneously for all $p\in\left(0,\infty\right]$.
Thus, the decomposition spaces $\FourierDecompSp{\CalB}p{\ell_{v^{\left(\gamma\right)}}^{q}}$
are well-defined.
\begin{defn}
\label{def:InhomogeneousBesovFourier}For $\dimension\in\N$, $p,q\in\left(0,\infty\right]$
and $\gamma\in\R$, we define the \textbf{Fourier-side inhomogeneous
Besov space} with integrability exponents $p,q$ and smoothness parameter
$\gamma$ as
\[
\BesovInhom pq{\Fourier,\gamma}\left(\R^{\dimension}\right):=\AlphaModSpace pq{\Fourier,\gamma}1\left(\R^{\dimension}\right):=\FourierDecompSp{\CalB}p{\ell_{v^{\left(\gamma\right)}}^{q}},
\]
where the weight $v^{\left(\gamma\right)}=\left(2^{\gamma n}\right)_{n\in\N_{0}}$
is as in Lemma~\ref{lem:InhomogeneousBesovCoveringAlmostStructured}.
\end{defn}

\begin{rem*}
As we will see in Lemma~\ref{lem:FourierInhomogeneousBesovAsTemperedDistributions},
the Fourier transform yields an isomorphism 
\[
\Fourier:\BesovInhom pq{\gamma}\left(\R^{\dimension}\right)\subset\Schwartz'\left(\R^{\dimension}\right)\to\BesovInhom pq{\Fourier,\gamma}\left(\R^{\dimension}\right)\subset\DistributionSpace{\R^{\dimension}},f\mapsto\widehat{f}|_{\TestFunctionSpace{\R^{\dimension}}}
\]
between the usual (inhomogeneous) Besov space $\BesovInhom pq{\gamma}\left(\R^{\dimension}\right)$
and $\BesovInhom pq{\Fourier,\gamma}\left(\R^{\dimension}\right)$.
This will justify the name ``Fourier-side inhomogeneous Besov space.''
\end{rem*}
As in the previous subsection, the main step for establishing embedding
results is to obtain a suitable subordinateness and moderateness statement.
\begin{lem}
\label{lem:AlphaModulationBesovSubordinatenessModerateness}Let $\dimension\in\N$
and $\alpha\in\left[0,1\right)$ be arbitrary and choose $r>r_{0}\left(\dimension,\alpha\right)$,
with $r_{0}\left(\dimension,\alpha\right)$ as in Theorem~\ref{thm:AlphaCoveringExistence}.
Then, the covering $\CalQ_{r}^{\left(\alpha\right)}$ is almost subordinate
to $\CalB$.

Furthermore, we have
\begin{equation}
\left|k\right|\asymp\left\langle k\right\rangle \asymp2^{n\left(1-\alpha\right)}\qquad\forall\,k\in\Z^{\dimension}\setminus\left\{ 0\right\} \text{ and }n\in\N_{0}\text{ with }Q_{r,k}^{\left(\alpha\right)}\cap B_{n}\neq\emptyset,\label{eq:AlphaModulationBesovModerateness}
\end{equation}
where the implied constants only depend on $\dimension,r,\alpha$.

Thus, $\CalQ_{r}^{\left(\alpha\right)}$ and $w^{\left(\gamma\right)}$
are relatively $\CalB$-moderate, with $w^{\left(\gamma\right)}=\left(\left\langle k\right\rangle ^{\gamma}\right)_{k\in\Z^{\dimension}\setminus\left\{ 0\right\} }$
as in Lemma~\ref{lem:AlphaModulationCoveringNormEstimate}.
\end{lem}

\begin{proof}
Let us first show that $\CalQ_{r}^{\left(\alpha\right)}=\left(\smash{Q_{r,k}^{\left(\alpha\right)}}\right)_{k\in\Z^{\dimension}\setminus\left\{ 0\right\} }$
is almost subordinate to $\CalB=\left(B_{n}\right)_{n\in\N_{0}}$.
Observe that each of the sets $Q_{r,k}^{\left(\alpha\right)}$ is
convex and hence path-connected and that $\CalQ_{r}^{\left(\alpha\right)}$
and $\CalB$ are both admissible coverings of all of $\R^{\dimension}$.
Moreover, each of the sets $B_{n}$ is open, so that Corollary~\ref{cor:WeakSubordinationImpliesSubordinationIfConnected}
implies that it suffices to show that the cardinality of the sets
\[
J_{k}:=\left\{ n\in\N_{0}\with B_{n}\cap Q_{r,k}^{\left(\alpha\right)}\neq\emptyset\right\} 
\]
is uniformly bounded with respect to $k\in\Z^{\dimension}\setminus\left\{ 0\right\} $.

To see this, let $k\in\Z^{\dimension}\setminus\left\{ 0\right\} $
be arbitrary and choose $n\in J_{k}$. Hence, there is some $\xi\in B_{n}\cap Q_{r,k}^{\left(\alpha\right)}$.
Lemma~\ref{lem:AlphaModulationCoveringNormEstimate} implies
\[
\left\langle \xi\right\rangle \asymp\left\langle k\right\rangle ^{\frac{1}{1-\alpha}}
\]
where the implied constants only depend on $r,\alpha$.

Let us assume $n\geq2$ for the moment. This yields $\left|\xi\right|\geq2^{n-2}\geq1$
because of $\xi\in B_{n}$. Hence, $\left|\xi\right|\leq\left\langle \xi\right\rangle \leq1+\left|\xi\right|\leq2\left|\xi\right|$.
In combination with $2^{n-2}<\left|\xi\right|<2^{n+2}$, i.e. $\left\langle \xi\right\rangle \asymp\left|\xi\right|\asymp2^{n}$,
we arrive at
\begin{equation}
2^{n}\asymp\left\langle \xi\right\rangle \asymp\left\langle k\right\rangle ^{\frac{1}{1-\alpha}}\quad\text{ for }k\in\Z^{\dimension}\setminus\left\{ 0\right\} \text{ and }n\in J_{k}\cap\N_{\geq2}.\label{eq:AlphaModulationBesovModeratenessAwayFromOrigin}
\end{equation}
This yields 
\[
2^{\frac{\log_{2}\left\langle k\right\rangle }{1-\alpha}+\log_{2}C_{0}}=C_{0}\cdot\left\langle k\right\rangle ^{\frac{1}{1-\alpha}}\leq2^{n}\leq C_{1}\cdot\left\langle k\right\rangle ^{\frac{1}{1-\alpha}}=2^{\frac{\log_{2}\left\langle k\right\rangle }{1-\alpha}+\log_{2}C_{1}}
\]
for suitable constants $C_{0}\in\left(0,1\right)$ and $C_{1}\geq1$,
which only depend on $r,\alpha$. We conclude
\[
\log_{2}C_{0}\leq n-\frac{\log_{2}\left\langle k\right\rangle }{1-\alpha}\leq\log_{2}C_{1}\,,
\]
and hence $n\in\N_{0}\cap\overline{B_{R}}\left(\frac{\log_{2}\left\langle k\right\rangle }{1-\alpha}\right)$
for $R:=\max\left\{ -\log_{2}C_{0},\log_{2}C_{1}\right\} $.

By dropping the assumption $n\geq2$, we finally arrive at
\[
J_{k}\subset\left\{ 0,1\right\} \cup\left[\N_{0}\cap\overline{B_{R}}\left(\frac{\log_{2}\left\langle k\right\rangle }{1-\alpha}\right)\right],
\]
where the latter set has cardinality at most $2+2\left\lceil R\right\rceil +1=2\left\lceil R\right\rceil +3=:N$.
As observed above, this implies that $\CalQ_{r}^{\left(\alpha\right)}$
is almost subordinate to $\mathcal{B}$. More precisely, Lemma~\ref{lem:WeaksubordinatenessImpliesSubordinatenessIfConnected}
shows that $Q_{r,k}^{\left(\alpha\right)}\subset B_{n}^{N\ast}$ holds
for all $k\in\Z^{\dimension}\setminus\left\{ 0\right\} $ and all
$n\in J_{k}$.

\medskip{}

Now, we establish equation~(\ref{eq:AlphaModulationBesovModerateness}).
For $n\in\N_{\geq2}$, and $k\in\Z^{\dimension}\setminus\left\{ 0\right\} $
with $Q_{r,k}^{\left(\alpha\right)}\cap B_{n}\neq\emptyset$, we get
$\left\langle k\right\rangle \asymp2^{n\left(1-\alpha\right)}$ from
equation~(\ref{eq:AlphaModulationBesovModeratenessAwayFromOrigin}).
Furthermore, since $\left|k\right|\geq1$, we also have $\left|k\right|\leq\left\langle k\right\rangle \leq1+\left|k\right|\leq2\left|k\right|$
and hence $\left|k\right|\asymp\left\langle k\right\rangle \asymp2^{n\left(1-\alpha\right)}$.
This establishes equation~(\ref{eq:AlphaModulationBesovModerateness})
for $n\geq2$. In case of $n\leq1$, we note that Lemma~\ref{lem:AlphaModulationCoveringNormEstimate}
yields for arbitrary $\xi\in Q_{r,k}^{\left(\alpha\right)}\cap B_{n}$
that 
\[
1\leq\left|k\right|\leq\left\langle k\right\rangle \asymp\left\langle \xi\right\rangle ^{1-\alpha}\leq\left(1+\left|\xi\right|\right)^{1-\alpha}<\left(1+2^{n+2}\right)^{1-\alpha}\leq9^{1-\alpha}\lesssim1,
\]
as well as $1\leq2^{n\left(1-\alpha\right)}\leq2^{1-\alpha}\lesssim1$,
where all implied constants only depend on $r$ and on $\alpha\in\left[0,1\right)$.
Thus, equation~(\ref{eq:AlphaModulationBesovModerateness}) holds
in all cases.
\end{proof}
Now, it is again easy to establish sharp embeddings between the Fourier-side
$\alpha$-modulation spaces and the Fourier-side inhomogeneous Besov
spaces. We remark that the following theorem is essentially identical
to \cite[Theorem 6.2.8]{VoigtlaenderPhDThesis} from my PhD thesis.
\begin{thm}
\label{thm:AlphaModulationBesovEmbedding}Let $\dimension\in\N$,
$\alpha\in\left[0,1\right)$ and $p_{1},p_{2},q_{1},q_{2}\in\left(0,\infty\right]$,
as well as $\gamma_{1},\gamma_{2}\in\R$.

We have $\AlphaModSpace{p_{1}}{q_{1}}{\Fourier,\gamma_{1}}{\alpha}\left(\R^{\dimension}\right)\subset\AlphaModSpace{p_{2}}{q_{2}}{\Fourier,\gamma_{2}}1\left(\R^{\dimension}\right)$
if and only if\footnote{This equivalence uses that both spaces continuously embed into the
Hausdorff space $\DistributionSpace{\R^{\dimension}}$ (see Theorem~\ref{thm:DecompositionSpaceComplete}),
and that the closed graph theorem also applies to quasi-Banach spaces,
as seen in Section~\ref{subsec:Notation}.} $\AlphaModSpace{p_{1}}{q_{1}}{\Fourier,\gamma_{1}}{\alpha}\left(\R^{\dimension}\right)\hookrightarrow\AlphaModSpace{p_{2}}{q_{2}}{\Fourier,\gamma_{2}}1\left(\R^{\dimension}\right)$,
if and only if we have
\[
p_{1}\leq p_{2}\qquad\text{and}\qquad\begin{cases}
\gamma_{2}\leq\gamma_{1}+\alpha\dimension\left(\frac{1}{p_{2}}-\frac{1}{p_{1}}\right)+\dimension\left(\alpha-1\right)\left(\smash{\frac{1}{\LowerExpo{p_{2}}}}-\frac{1}{q_{1}}\right)_{+}, & \text{if }q_{1}\leq q_{2},\vspace{0.15cm}\\
\gamma_{2}<\gamma_{1}+\alpha\dimension\left(\frac{1}{p_{2}}-\frac{1}{p_{1}}\right)+\dimension\left(\alpha-1\right)\left(\smash{\frac{1}{\LowerExpo{p_{2}}}}-\frac{1}{q_{1}}\right)_{+}, & \text{if }q_{1}>q_{2}.
\end{cases}
\]

Conversely, we have $\AlphaModSpace{p_{1}}{q_{1}}{\Fourier,\gamma_{1}}1\left(\R^{\dimension}\right)\subset\AlphaModSpace{p_{2}}{q_{2}}{\Fourier,\gamma_{2}}{\alpha}\left(\R^{\dimension}\right)$
if and only if $\AlphaModSpace{p_{1}}{q_{1}}{\Fourier,\gamma_{1}}1\left(\R^{\dimension}\right)\hookrightarrow\AlphaModSpace{p_{2}}{q_{2}}{\Fourier,\gamma_{2}}{\alpha}\left(\R^{\dimension}\right)$,
if and only if we have 
\[
p_{1}\leq p_{2}\qquad\text{and}\qquad\begin{cases}
\gamma_{2}<\gamma_{1}+\alpha\dimension\left(\frac{1}{p_{2}}-\frac{1}{p_{1}}\right)+\dimension\left(\alpha-1\right)\left(\frac{1}{q_{2}}-\smash{\frac{1}{\SignedUpperExpo{p_{1}}}}\right)_{\!+}, & \text{if }q_{1}>q_{2},\vspace{0.15cm}\\
\gamma_{2}\leq\gamma_{1}+\alpha\dimension\left(\frac{1}{p_{2}}-\frac{1}{p_{1}}\right)+\dimension\left(\alpha-1\right)\left(\frac{1}{q_{2}}-\smash{\frac{1}{\SignedUpperExpo{p_{1}}}}\right)_{\!+}, & \text{if }q_{1}\leq q_{2}.
\end{cases}\qedhere
\]
\end{thm}

\begin{proof}
Let $r_{0}=r_{0}\left(\dimension,\alpha\right)$ as in Theorem~\ref{thm:AlphaCoveringExistence},
choose $r>r_{0}$ and recall the definition of the spaces $\AlphaModSpace pq{\Fourier,\gamma}{\alpha}\left(\R^{\dimension}\right)=\FourierDecompSp{\CalQ_{r}^{\left(\alpha\right)}}p{\ell_{w^{\left(\gamma/\left(1-\alpha\right)\right)}}^{q}}$
from Definition~\ref{def:AlphaModulationFourierSpace}.

We first analyze for which parameters the inclusion $\AlphaModSpace{p_{1}}{q_{1}}{\Fourier,\gamma_{1}}{\alpha}\left(\R^{\dimension}\right)\subset\AlphaModSpace{p_{2}}{q_{2}}{\Fourier,\gamma_{2}}1\left(\R^{\dimension}\right)$
holds. By the closed graph theorem, this holds if and only if the
identity map $\iota:\AlphaModSpace{p_{1}}{q_{1}}{\Fourier,\gamma_{1}}{\alpha}\left(\R^{\dimension}\right)\to\AlphaModSpace{p_{2}}{q_{2}}{\Fourier,\gamma_{2}}1\left(\R^{\dimension}\right),f\mapsto f$
is well-defined and bounded. To characterize when this is the case,
we set $\alpha_{0}:=\frac{\alpha}{1-\alpha}$, and define
\[
\CalQ=\left(Q_{i}\right)_{i\in I}=\left(T_{i}Q_{i}'+b_{i}\right)_{i\in I}:=\CalQ_{r}^{\left(\alpha\right)}=\left(Q_{r,k}^{\left(\alpha\right)}\right)_{k\in\Z^{\dimension}\setminus\left\{ 0\right\} }=\left(\left(\left|k\right|^{\alpha_{0}}\cdot\identity\right)Q^{\left(r\right)}+\left|k\right|^{\alpha_{0}}k\right)_{k\in\Z^{\dimension}\setminus\left\{ 0\right\} }\,,
\]
and
\[
\CalP=\left(P_{j}\right)_{j\in J}=\left(S_{j}P_{j}'+c_{j}\right)_{j\in J}:=\CalB=\left(B_{n}\right)_{n\in\N_{0}}=\bigl(\left(2^{n}\cdot\identity\right)B_{n}'\bigr)_{n\in\N_{0}}\,,
\]
where the sets $Q^{\left(r\right)}$ and $B_{n}'$ are defined as
in Lemmas~\ref{lem:AlphaModulationCoveringIsStructured} and \ref{lem:InhomogeneousBesovCoveringAlmostStructured},
respectively. Furthermore, we choose $w:=w^{\left(\gamma_{1}/\left(1-\alpha\right)\right)}$
and $v:=v^{\left(\gamma_{2}\right)}$, where $w^{\left(\gamma\right)}$
and $v^{\left(\gamma\right)}$ are defined as in Lemmas~\ref{lem:AlphaModulationCoveringNormEstimate}
and \ref{lem:InhomogeneousBesovCoveringAlmostStructured}, respectively.
With these choices, we have $\AlphaModSpace{p_{1}}{q_{1}}{\Fourier,\gamma_{1}}{\alpha}\left(\R^{\dimension}\right)=\FourierDecompSp{\CalQ}{p_{1}}{\ell_{w}^{q_{1}}}$
and $\AlphaModSpace{p_{2}}{q_{2}}{\Fourier,\gamma_{2}}1\left(\R^{\dimension}\right)=\FourierDecompSp{\CalP}{p_{2}}{\ell_{v}^{q_{2}}}$.

Lemmas~\ref{lem:AlphaModulationCoveringIsStructured}, \ref{lem:AlphaModulationCoveringNormEstimate}
and \ref{lem:InhomogeneousBesovCoveringAlmostStructured} imply that
both $\CalQ$ and $\CalP$ are almost structured coverings and that
$w$ and $v$ are $\CalQ$-moderate and $\CalP$-moderate, respectively.
Hence, the standing assumptions from Section~\ref{sec:SummaryOfEmbeddingResults}
(i.e., Assumption~\ref{assu:GeneralSummaryAssumptions}) are satisfied.

Finally, $\CalQ,\CalP$ are both coverings of the same set $\CalO=\CalO'=\R^{\dimension}$,
and Lemma~\ref{lem:AlphaModulationBesovSubordinatenessModerateness}
shows that $\CalQ$ is almost subordinate to $\CalP$ and that $\CalQ$
and $w$ are relatively $\CalP$-moderate. Thus, all assumptions of
part~(\ref{enu:SummaryFineInCoarseModerate}) of Theorem~\ref{thm:SummaryFineIntoCoarse}
are satisfied. Note that the embedding $\iota$ from Theorem~\ref{thm:SummaryFineIntoCoarse}
satisfies $\iota f=f$ for all $f\in\FourierDecompSp{\CalQ}{p_{1}}{\ell_{w}^{q_{1}}}$,
since both $\CalQ$ and $\CalP$ cover the same set $\CalO=\R^{\dimension}=\CalO'$.
Thus, the map $\iota$ from above coincides with $\iota$ from Theorem~\ref{thm:SummaryFineIntoCoarse}.
All in all, we conclude that $\iota$ is well-defined and bounded
if and only if we have
\[
p_{1}\leq p_{2}\qquad\text{and}\qquad K:=\left\Vert \!\left(\!\frac{v_{j}}{w_{i_{j}}}\cdot\left|\det\smash{T_{i_{j}}}\right|^{p_{1}^{-1}-p_{2}^{-1}-s}\cdot\left|\det S_{j}\right|^{s}\right)_{\!\!\!j\in J}\right\Vert _{\ell^{q_{2}\cdot\left(q_{1}/q_{2}\right)'}}<\infty,
\]
where $s:=\left(\smash{\frac{1}{\LowerExpo{p_{2}}}}-\frac{1}{q_{1}}\right)_{+}$
and where for each $j\in J=\N_{0}$, some $i_{j}\in I_{j}\subset\Z^{\dimension}\setminus\left\{ 0\right\} $,
i.e.\@ with $\emptyset\neq B_{j}\cap Q_{r,i_{j}}^{\left(\alpha\right)}$
is selected. But according to Lemma~\ref{lem:AlphaModulationBesovSubordinatenessModerateness},
we have
\[
w_{i_{j}}=w_{i_{j}}^{\left(\gamma_{1}/\left(1-\alpha\right)\right)}=\left\langle i_{j}\right\rangle ^{\frac{\gamma_{1}}{1-\alpha}}\asymp2^{j\left(1-\alpha\right)\frac{\gamma_{1}}{1-\alpha}}=2^{j\gamma_{1}}\,,
\]
and
\[
\left|\det\smash{T_{i_{j}}}\right|=\left|i_{j}\right|^{\dimension\alpha_{0}}\asymp\left(2^{j\left(1-\alpha\right)}\right)^{\dimension\frac{\alpha}{1-\alpha}}=2^{j\dimension\alpha},\quad\text{since}\quad\alpha_{0}=\frac{\alpha}{1-\alpha}\,.
\]

All in all, we get
\begin{align}
K & \asymp\left\Vert \left(2^{j\left(\gamma_{2}-\gamma_{1}\right)}\cdot2^{j\dimension\alpha\left(p_{1}^{-1}-p_{2}^{-1}-s\right)}\cdot2^{j\dimension s}\right)_{j\in\N_{0}}\right\Vert _{\ell^{q_{2}\cdot\left(q_{1}/q_{2}\right)'}}\nonumber \\
 & \asymp\left\Vert \left(2^{j\left(\gamma_{2}-\gamma_{1}+\dimension\left[s\left(1-\alpha\right)+\alpha\left(p_{1}^{-1}-p_{2}^{-1}\right)\right]\right)}\right)_{j\in\N_{0}}\right\Vert _{\ell^{q_{2}\cdot\left(q_{1}/q_{2}\right)'}}.\label{eq:BesovAlphaModulationEmbeddingSequenceNormAsymptotic}
\end{align}
Recall from equation~(\ref{eq:SpecialExponentFiniteness}) that $q_{2}\cdot\left(q_{1}/q_{2}\right)'$
is finite if and only if $q_{2}<q_{1}$. Hence, due to the exponential
nature of the sequence in equation~(\ref{eq:BesovAlphaModulationEmbeddingSequenceNormAsymptotic}),
we get
\begin{align*}
K<\infty & \Longleftrightarrow\begin{cases}
\gamma_{2}-\gamma_{1}+\dimension\left[s\left(1-\alpha\right)+\alpha\left(p_{1}^{-1}-p_{2}^{-1}\right)\right]\leq0, & \text{if }q_{1}\leq q_{2},\vspace{0.15cm}\\
\gamma_{2}-\gamma_{1}+\dimension\left[s\left(1-\alpha\right)+\alpha\left(p_{1}^{-1}-p_{2}^{-1}\right)\right]<0, & \text{if }q_{1}>q_{2}
\end{cases}\\
 & \Longleftrightarrow\begin{cases}
\gamma_{2}\leq\gamma_{1}+\alpha\dimension\left(p_{2}^{-1}-p_{1}^{-1}\right)+\dimension\left(\alpha-1\right)\left(\smash{\frac{1}{\LowerExpo{p_{2}}}}-\frac{1}{q_{1}}\right)_{+}, & \text{if }q_{1}\leq q_{2},\vspace{0.2cm}\\
\gamma_{2}<\gamma_{1}+\alpha\dimension\left(p_{2}^{-1}-p_{1}^{-1}\right)+\dimension\left(\alpha-1\right)\left(\smash{\frac{1}{\LowerExpo{p_{2}}}}-\frac{1}{q_{1}}\right)_{+}, & \text{if }q_{1}>q_{2}.
\end{cases}
\end{align*}
This completes the characterization of the inclusion $\AlphaModSpace{p_{1}}{q_{1}}{\Fourier,\gamma_{1}}{\alpha}\left(\R^{\dimension}\right)\subset\AlphaModSpace{p_{2}}{q_{2}}{\Fourier,\gamma_{2}}1\left(\R^{\dimension}\right)$.

\medskip{}

To characterize the cases in which the inclusion $\AlphaModSpace{p_{1}}{q_{1}}{\Fourier,\gamma_{1}}1\left(\R^{\dimension}\right)\subset\AlphaModSpace{p_{2}}{q_{2}}{\Fourier,\gamma_{2}}{\alpha}\left(\R^{\dimension}\right)$
holds, we note as above that this holds if and only if the identity
map $\theta:\AlphaModSpace{p_{1}}{q_{1}}{\Fourier,\gamma_{1}}1\left(\R^{\dimension}\right)\to\AlphaModSpace{p_{2}}{q_{2}}{\Fourier,\gamma_{2}}{\alpha}\left(\R^{\dimension}\right),f\mapsto f$
is well-defined and bounded. To characterize when this is the case,
we set
\[
\CalQ=\left(Q_{i}\right)_{i\in I}=\left(T_{i}Q_{i}'+b_{i}\right)_{i\in I}:=\CalB=\left(B_{n}\right)_{n\in\N_{0}}=\bigl(\left(2^{n}\cdot\identity\right)B_{n}'\bigr)_{n\in\N_{0}}\,,
\]
and
\[
\CalP=\left(P_{j}\right)_{j\in J}=\left(S_{j}P_{j}'+c_{j}\right)_{j\in J}:=\CalQ_{r}^{\left(\alpha\right)}=\left(Q_{r,k}^{\left(\alpha\right)}\right)_{k\in\Z^{\dimension}\setminus\left\{ 0\right\} }=\left(\left(\left|k\right|^{\alpha_{0}}\cdot\identity\right)Q^{\left(r\right)}+\left|k\right|^{\alpha_{0}}k\right)_{k\in\Z^{\dimension}\setminus\left\{ 0\right\} }\,,
\]
with $B_{n}'$, $Q^{\left(r\right)}$, and $\alpha_{0}$ as above.
Further, we set  $w:=v^{\left(\gamma_{1}\right)}$ and $v:=w^{\left(\gamma_{2}/\left(1-\alpha\right)\right)}$,
where $w^{\left(\gamma\right)}$ and $v^{\left(\gamma\right)}$ are
defined as in Lemmas \ref{lem:AlphaModulationCoveringNormEstimate}
and \ref{lem:InhomogeneousBesovCoveringAlmostStructured}, respectively.
These choices ensure $\AlphaModSpace{p_{1}}{q_{1}}{\Fourier,\gamma_{1}}1\left(\R^{\dimension}\right)=\FourierDecompSp{\CalQ}{p_{1}}{\ell_{w}^{q_{1}}}$
and $\AlphaModSpace{p_{2}}{q_{2}}{\Fourier,\gamma_{2}}{\alpha}\left(\R^{\dimension}\right)=\FourierDecompSp{\CalP}{p_{2}}{\ell_{v}^{q_{2}}}$.
Precisely as above, we see that Assumption~\ref{assu:GeneralSummaryAssumptions}
is fulfilled for these choices.

Finally, $\CalQ$ and $\CalP$ both cover the same set $\CalO=\R^{\dimension}=\CalO'$,
and Lemma~\ref{lem:AlphaModulationBesovSubordinatenessModerateness}
shows that $\CalP$ is almost subordinate to $\CalQ$ and that $\CalP$
and $v$ are relatively $\CalQ$-moderate. Thus, all assumptions of
part~(\ref{enu:SummaryCoarseInFineModerate}) of Theorem~\ref{thm:SummaryCoarseIntoFine}
are satisfied. As above, we see that the map $\iota$ from that theorem
satisfies $\iota f=f=\theta f$ for all $f\in\FourierDecompSp{\CalQ}{p_{1}}{\ell_{w}^{q_{1}}}$
with $\theta$ as above. Thus, all in all, we conclude that $\theta$
is well-defined and bounded if and only if we have
\[
p_{1}\leq p_{2}\qquad\text{and}\qquad K:=\left\Vert \left(\frac{v_{j_{i}}}{w_{i}}\cdot\left|\det T_{i}\right|^{s}\cdot\left|\det\smash{S_{j_{i}}}\right|^{p_{1}^{-1}-p_{2}^{-1}-s}\right)_{i\in I}\right\Vert _{\ell^{q_{2}\cdot\left(q_{1}/q_{2}\right)'}}<\infty,
\]
where $s:=\left(\frac{1}{q_{2}}-\smash{\frac{1}{\SignedUpperExpo{p_{1}}}}\right)_{+}$
and where for each $i\in I=\N_{0}$, some $j_{i}\in J_{i}\subset\Z^{\dimension}\setminus\left\{ 0\right\} $,
i.e.\@ with $\emptyset\neq B_{i}\cap Q_{r,j_{i}}^{\left(\alpha\right)}$
is selected.

But since $\alpha_{0}=\frac{\alpha}{1-\alpha}$ and in view of Lemma~\ref{lem:AlphaModulationBesovSubordinatenessModerateness},
we see—exactly as above—that
\[
v_{j_{i}}=\left\langle j_{i}\right\rangle ^{\frac{\gamma_{2}}{1-\alpha}}\asymp2^{i\left(1-\alpha\right)\frac{\gamma_{2}}{1-\alpha}}=2^{i\gamma_{2}}\qquad\text{ and }\qquad\left|\det\smash{S_{j_{i}}}\right|=\left|j_{i}\right|^{\dimension\alpha_{0}}\asymp2^{i\left(1-\alpha\right)\dimension\alpha_{0}}=2^{i\dimension\alpha}\,.
\]
Thus, setting $t:=q_{2}\cdot\left(q_{1}/q_{2}\right)'$, we get
\[
K\asymp\left\Vert \left(2^{i\left(\gamma_{2}-\gamma_{1}\right)}\cdot2^{i\dimension s}\cdot2^{i\dimension\alpha\left(p_{1}^{-1}-p_{2}^{-1}-s\right)}\right)_{i\in\N_{0}}\right\Vert _{\ell^{t}}\asymp\left\Vert \left(2^{i\left(\gamma_{2}-\gamma_{1}+\dimension s\left(1-\alpha\right)+\dimension\alpha\left(p_{1}^{-1}-p_{2}^{-1}\right)\right)}\right)_{i\in\N_{0}}\right\Vert _{\ell^{t}}\,.
\]
The remainder of the proof is exactly as above and hence omitted.
\end{proof}
We note that our embedding results for $\alpha$-modulation spaces
from Theorem~\ref{thm:AlphaModulationEmbeddings} apply also to the
case $\alpha=\beta$. In contrast, the preceding theorem requires
$\alpha\in\left[0,1\right)$, so that embeddings between two (different)
inhomogeneous Besov spaces are strictly speaking not covered by that
criterion. This motivates the following theorem—which is well-known
in the theory of Besov spaces, and only presented here to show that
it is a straightforward consequence of the developed criteria.
\begin{thm}
\label{thm:InhomBesovEmbeddings}Let $\dimension\in\N$, $p_{1},p_{2},q_{1},q_{2}\in\left(0,\infty\right]$
and $\gamma_{1},\gamma_{2}\in\R$ be arbitrary.

We have $\AlphaModSpace{p_{1}}{q_{1}}{\Fourier,\gamma_{1}}1\left(\R^{\dimension}\right)\subset\AlphaModSpace{p_{2}}{q_{2}}{\Fourier,\gamma_{2}}1\left(\R^{\dimension}\right)$
if and only if $\AlphaModSpace{p_{1}}{q_{1}}{\Fourier,\gamma_{1}}1\left(\R^{\dimension}\right)\hookrightarrow\AlphaModSpace{p_{2}}{q_{2}}{\Fourier,\gamma_{2}}1\left(\R^{\dimension}\right)$,
if and only if
\[
p_{1}\leq p_{2}\qquad\text{and}\qquad\begin{cases}
\gamma_{2}\leq\gamma_{1}+\dimension\left(p_{2}^{-1}-p_{1}^{-1}\right), & \text{if }q_{1}\leq q_{2},\vspace{0.1cm}\\
\gamma_{2}<\gamma_{1}+\dimension\left(p_{2}^{-1}-p_{1}^{-1}\right), & \text{if }q_{1}>q_{2}.
\end{cases}\qedhere
\]
\end{thm}

\begin{rem*}
In contrast to Theorems \ref{thm:AlphaModulationEmbeddings} and \ref{thm:AlphaModulationBesovEmbedding},
the preceding theorem is \emph{not} a new result. For example, \cite[Proposition 2.4]{HanWangAlphaModulationEmbeddings}
already shows that the stated conditions are sufficient for the existence
of the embedding. We state the theorem here mainly for completeness,
and to give an example application of Corollary~\ref{cor:SummarySameCovering}.
\end{rem*}
\begin{proof}
Let $\CalQ=\left(Q_{i}\right)_{i\in I}=\left(T_{i}Q_{i}'+b_{i}\right)_{i\in I}:=\CalB=\left(B_{n}\right)_{n\in\N_{0}}=\bigl(\left(2^{n}\cdot\identity\right)B_{n}'\bigr)_{n\in\N_{0}}$,
and $w:=v^{\left(\gamma_{1}\right)}$, as well as $v:=v^{\left(\gamma_{2}\right)}$.
Here, the sets $B_{n}'$ and the weight $v^{\left(\gamma\right)}$
are chosen as in Lemma~\ref{lem:InhomogeneousBesovCoveringAlmostStructured}.
We want to apply Corollary~\ref{cor:SummarySameCovering} with these
choices. Since $\CalB$ is an almost structured covering by Lemma~\ref{lem:InhomogeneousBesovCoveringAlmostStructured}
and since the same lemma shows that $v,w$ are $\CalB$-moderate,
the general assumptions of Section~\ref{sec:SummaryOfEmbeddingResults}
(i.e., Assumption~\ref{assu:GeneralSummaryAssumptions}) are satisfied.

Thus, Corollary~\ref{cor:SummarySameCovering} shows that the embedding
$\AlphaModSpace{p_{1}}{q_{1}}{\Fourier,\gamma_{1}}1\left(\R^{\dimension}\right)\hookrightarrow\AlphaModSpace{p_{2}}{q_{2}}{\Fourier,\gamma_{2}}1\left(\R^{\dimension}\right)$
holds if and only if we have $p_{1}\leq p_{2}$ and if
\[
K:=\left\Vert \left(\left|\det T_{i}\right|^{p_{1}^{-1}-p_{2}^{-1}}\cdot\frac{v_{i}}{w_{i}}\right)_{i\in I}\right\Vert _{\ell^{q_{2}\cdot\left(q_{1}/q_{2}\right)'}}=\left\Vert \left(2^{n\left[\gamma_{2}-\gamma_{1}+\dimension\left(p_{1}^{-1}-p_{2}^{-1}\right)\right]}\right)_{n\in\N_{0}}\right\Vert _{\ell^{q_{2}\cdot\left(q_{1}/q_{2}\right)'}}
\]
is finite. But equation~(\ref{eq:SpecialExponentFiniteness}) shows
that $q_{2}\cdot\left(q_{1}/q_{2}\right)'$ is finite if and only
if we have $q_{2}<q_{1}$. Due to the exponential nature of the weight
in the preceding equation, we thus see that $K$ is finite if and
only if we have
\[
\begin{cases}
\gamma_{2}-\gamma_{1}+\dimension\left(p_{1}^{-1}-p_{2}^{-1}\right)\leq0, & \text{if }q_{1}\leq q_{2},\vspace{0.1cm}\\
\gamma_{2}-\gamma_{1}+\dimension\left(p_{1}^{-1}-p_{2}^{-1}\right)<0, & \text{if }q_{1}>q_{2}.
\end{cases}
\]
Finally, the equivalence between the embedding $\AlphaModSpace{p_{1}}{q_{1}}{\Fourier,\gamma_{1}}1\left(\R^{\dimension}\right)\hookrightarrow\AlphaModSpace{p_{2}}{q_{2}}{\Fourier,\gamma_{2}}1\left(\R^{\dimension}\right)$
and the inclusion $\AlphaModSpace{p_{1}}{q_{1}}{\Fourier,\gamma_{1}}1\left(\R^{\dimension}\right)\subset\AlphaModSpace{p_{2}}{q_{2}}{\Fourier,\gamma_{2}}1\left(\R^{\dimension}\right)$
follows from the closed graph theorem (which applies also to quasi-Banach
spaces, see Section~\ref{subsec:Notation}), since both involved
spaces embed continuously into the Hausdorff space $\DistributionSpace{\R^{\dimension}}$,
as we saw in Theorem~\ref{thm:DecompositionSpaceComplete}.
\end{proof}
Now, we show that our ``Fourier-side'' inhomogeneous Besov space
$\BesovInhom pq{\Fourier,\gamma}\left(\R^{\dimension}\right)$ indeed
coincides with the Fourier transform of the usual inhomogeneous Besov
space $\BesovInhom pq{\gamma}\left(\R^{\dimension}\right)$:
\begin{lem}
\label{lem:FourierInhomogeneousBesovAsTemperedDistributions}Let $\dimension\in\N$,
$p,q\in\left(0,\infty\right]$ and $\gamma\in\R$ be arbitrary. Let
the inhomogeneous Besov space $\BesovInhom pq{\gamma}\left(\R^{\dimension}\right)\subset\Schwartz'\left(\R^{\dimension}\right)$
be as defined in\footnote{Note that the exact notation in \cite[Definition 2.2.1]{GrafakosModern}
is actually slightly different. What we call $\BesovInhom pq{\gamma}\left(\R^{\dimension}\right)$
here is written as $B_{p}^{\gamma,q}$ in \cite[Definition 2.2.1]{GrafakosModern}.
Just as for the $\alpha$-modulation spaces, the notation is not completely
consistent in the literature. Our choice of notation is motivated
by the Lebesgue spaces $L^{p}$ and the Lorentz spaces $L^{p,q}$,
where the integrability exponents are always at the top.} \cite[Definition 2.2.1]{GrafakosModern}. 

Then, the Fourier transform restricts to an isomorphism of (quasi)-Banach
spaces
\begin{equation}
\Fourier:\BesovInhom pq{\gamma}\left(\R^{\dimension}\right)\to\AlphaModSpace pq{\Fourier,\gamma}1\left(\R^{\dimension}\right),f\mapsto\widehat{f}|_{\TestFunctionSpace{\R^{\dimension}}}\:.\label{eq:BesovFourierIsomorphism}
\end{equation}
In particular, every $f\in\AlphaModSpace pq{\Fourier,\gamma}1\left(\R^{\dimension}\right)\subset\DistributionSpace{\R^{\dimension}}$
extends to a \emph{tempered} distribution.
\end{lem}

\begin{proof}
As in \cite[Sections 1.3.2 and 2.2.1]{GrafakosModern}, we fix a
radial Schwartz function $\Psi\in\Schwartz\left(\R^{\dimension}\right)$
such that $\widehat{\Psi}$ is nonnegative with $\supp\widehat{\Psi}\subset\overline{B_{2}}\left(0\right)\setminus B_{6/7}\left(0\right)$,
with $\widehat{\Psi}\equiv1$ on $\overline{B_{12/7}}\left(0\right)\setminus B_{1}\left(0\right)$,
and such that
\begin{equation}
\sum_{j=-\infty}^{\infty}\widehat{\Psi}\left(\smash{2^{-j}\xi}\right)=1\qquad\forall\,\xi\in\R^{\dimension}\setminus\left\{ 0\right\} .\label{eq:BesovPartitionOfUnityEquation}
\end{equation}
Let  $\Psi_{j}\left(x\right):=2^{\dimension j}\cdot\Psi\left(2^{j}x\right)$
for $j\in\Z$ and $x\in\R^{\dimension}$. We then have $\widehat{\Psi_{j}}\left(\xi\right)=\widehat{\Psi}\left(2^{-j}\xi\right)=:\psi_{j}\left(\xi\right)$
for all $\xi\in\R^{\dimension}$. From this, it follows easily that
\[
\left\Vert \Fourier^{-1}\psi_{j}\right\Vert _{L^{p}}=\left\Vert \Psi_{j}\right\Vert _{L^{p}}=2^{\dimension j}\cdot\left\Vert D_{2^{j}{\rm id}}\Psi\right\Vert _{L^{p}}=2^{\dimension j\left(1-\frac{1}{p}\right)}\cdot\left\Vert \Psi\right\Vert _{L^{p}}=\left|\det\left(\smash{2^{j}}{\rm id}\right)\right|^{1-\frac{1}{p}}\cdot\left\Vert \Psi\right\Vert _{L^{p}}
\]
for all $j\in\Z$ and $p\in\left(0,\infty\right]$. Furthermore, using
$\supp\widehat{\Psi}\subset\overline{B_{2}}\left(0\right)\setminus B_{6/7}\left(0\right)\subset B_{4}\left(0\right)\setminus\overline{B_{1/4}}\left(0\right)$,
we see $\supp\psi_{j}=2^{j}\cdot\supp\widehat{\Psi}\subset B_{j}$
for all $j\in\N$, where $\CalB=\left(B_{n}\right)_{n\in\N_{0}}$
denotes the inhomogeneous Besov covering from Definition~\ref{def:InhomogeneousBesovCovering}.

Grafakos (cf. \cite[equation (2.2.2)]{GrafakosModern}) states that
there is a Schwartz function $\Phi\in\Schwartz\left(\R^{\dimension}\right)$
satisfying
\begin{equation}
\widehat{\Phi}\left(\xi\right)=\begin{cases}
\sum_{j\leq0}\widehat{\Psi}\left(2^{-j}\xi\right)=\sum_{j\leq0}\psi_{j}\left(\xi\right), & \text{if }\xi\neq0,\\
1, & \text{if }\xi=0
\end{cases}\label{eq:GrafakosSpecialLowPassFunction}
\end{equation}
and $\widehat{\Phi}\left(\xi\right)=1$ for $\left|\xi\right|\leq1$,
as well as $\widehat{\Phi}\left(\xi\right)=0$ for $\left|\xi\right|\geq2$.

Set $\varphi_{i}:=\psi_{i}$ for $i\in\N$, and let $\varphi_{0}:=\widehat{\Phi}$.
Using equations (\ref{eq:BesovPartitionOfUnityEquation}) and (\ref{eq:GrafakosSpecialLowPassFunction}),
it is not hard to see $\sum_{i\in\N_{0}}\varphi_{i}\left(\xi\right)=1$
for all $\xi\in\R^{\dimension}$. Furthermore, we have $\varphi_{i}\in\TestFunctionSpace{\R^{\dimension}}$
with $\supp\varphi_{i}\subset B_{i}$ for all $i\in\N$ and with $\supp\varphi_{0}\subset\overline{B_{2}}\left(0\right)\subset B_{4}\left(0\right)=B_{0}$.
Finally, we have
\[
\left|\det\left(\smash{2^{i}}{\rm id}\right)\right|^{\frac{1}{p}-1}\cdot\left\Vert \Fourier^{-1}\varphi_{i}\right\Vert _{L^{p}}\leq\max\left\{ \left\Vert \Phi\right\Vert _{L^{p}},\,\sup_{j\in\Z}\left|\det\left(\smash{2^{j}}{\rm id}\right)\right|^{\frac{1}{p}-1}\cdot\left\Vert \Fourier^{-1}\psi_{j}\right\Vert _{L^{p}}\right\} =\max\left\{ \left\Vert \Phi\right\Vert _{L^{p}},\,\left\Vert \Psi\right\Vert _{L^{p}}\right\} 
\]
for all $i\in\N_{0}$, and thus $C_{\CalB,\Phi,p}<\infty$ for $\Phi:=\left(\varphi_{i}\right)_{i\in\N_{0}}$
and arbitrary $p\in\left(0,\infty\right]$. All in all, we see that
$\Phi$ is an $L^{p}$-BAPU for $\CalB$ for all $p\in\left(0,\infty\right]$.

Now, still following Grafakos (see \cite[Section 2.2.1]{GrafakosModern}),
we define
\[
\Delta_{j}^{\Psi}:\Schwartz'\left(\smash{\R^{\dimension}}\right)\rightarrow\Schwartz'\left(\smash{\R^{\dimension}}\right),f\mapsto\Fourier^{-1}\left(\,\smash{\widehat{\Psi_{j}}}\cdot\smash{\widehat{f}}\,\,\right)=\Fourier^{-1}\left(\,\psi_{j}\cdot\smash{\widehat{f}}\,\,\right)
\]
for $j\in\Z$, as well as 
\[
S_{0}:\Schwartz'\left(\smash{\R^{\dimension}}\right)\rightarrow\Schwartz'\left(\smash{\R^{\dimension}}\right),f\mapsto\Phi\ast f=\Fourier^{-1}\left(\,\smash{\widehat{\Phi}}\cdot\smash{\widehat{f}}\,\,\right)=\Fourier^{-1}\left(\,\varphi_{0}\cdot\smash{\widehat{f}}\,\,\right).
\]
Using these notations, Grafakos defines the (quasi)-norm on $\BesovInhom pq{\gamma}\left(\R^{\dimension}\right)$
as
\begin{align*}
\left\Vert f\right\Vert _{\BesovInhom pq{\gamma}} & :=\left\Vert S_{0}f\right\Vert _{L^{p}}+\left\Vert \left(2^{j\gamma}\cdot\left\Vert \Delta_{j}^{\Psi}f\right\Vert _{L^{p}}\right)_{j\in\N}\right\Vert _{\ell^{q}}\\
 & =\left\Vert \Fourier^{-1}\left(\varphi_{0}\cdot\smash{\widehat{f}}\,\,\right)\right\Vert _{L^{p}}+\left\Vert \left(\left\Vert \Fourier^{-1}\left(\varphi_{j}\cdot\smash{\widehat{f}}\,\,\right)\right\Vert _{L^{p}}\right)_{j\in\N}\right\Vert _{\ell_{v^{\left(\gamma\right)}}^{q}}\in\left[0,\infty\right]\quad\text{for }f\in\Schwartz'\left(\R^{\dimension}\right)\,.
\end{align*}
Note that we have $\widehat{f}\in\Schwartz'\left(\R^{\dimension}\right)\subset\DistributionSpace{\R^{\dimension}}$
for $f\in\Schwartz'\left(\R^{\dimension}\right)$. Now, it is not
hard to see
\begin{equation}
\left\Vert \widehat{f}|_{\TestFunctionSpace{\R^{\dimension}}}\right\Vert _{\BAPUFourierDecompSp{\CalB}p{\ell_{v^{\left(\gamma\right)}}^{q}}{\Phi}}=\left\Vert \left(\left\Vert \Fourier^{-1}\left(\varphi_{i}\,\smash{\widehat{f}}\,\,\right)\right\Vert _{L^{p}}\right)_{i\in\N_{0}}\right\Vert _{\ell_{v^{\left(\gamma\right)}}^{q}}\asymp\left\Vert f\right\Vert _{\BesovInhom pq{\gamma}}\qquad\forall\,f\in\Schwartz'\left(\R^{\dimension}\right),\label{eq:BesovFourierNormEquivalence}
\end{equation}
so that the map $\Fourier$ defined in equation~(\ref{eq:BesovFourierIsomorphism})
is a well-defined isomorphism of quasi-normed vector spaces, when
considered as a map onto a certain subspace of $\AlphaModSpace pq{\Fourier,\gamma}1\left(\R^{\dimension}\right)$
(namely onto its range).

\medskip{}

It remains to show that the map $\Fourier$ from equation~(\ref{eq:BesovFourierIsomorphism})
is surjective. To this end, it suffices to show that every $f\in\FourierDecompSp{\CalB}p{\ell_{v^{\left(\gamma\right)}}^{q}}$
extends to a tempered distribution $\widetilde{f}\in\Schwartz'\left(\R^{\dimension}\right)$.
Indeed, if this is the case, equation~(\ref{eq:BesovFourierNormEquivalence})
implies $\left\Vert \Fourier^{-1}\smash{\widetilde{f}}\,\right\Vert _{\BesovInhom pq{\gamma}}\asymp\left\Vert f\right\Vert _{\FourierDecompSp{\CalB}p{\smash{\ell_{v^{\left(\gamma\right)}}^{q}}}}<\infty$,
and hence ${g:=\Fourier^{-1}\widetilde{f}\in\BesovInhom pq{\gamma}\left(\R^{\dimension}\right)}$
with $\Fourier g=\widetilde{f}|_{\TestFunctionSpace{\R^{\dimension}}}=f$,
so that $\Fourier$ is surjective. By completeness of $\FourierDecompSp{\CalB}p{\ell_{v^{\left(\gamma\right)}}^{q}}$
(see Theorem~\ref{thm:DecompositionSpaceComplete}), this incidentally
proves the completeness of $\BesovInhom pq{\gamma}\left(\R^{\dimension}\right)$,
so that the map $\Fourier$ as defined in equation~(\ref{eq:BesovFourierIsomorphism})
is indeed an isomorphism of quasi-Banach spaces.

But Theorem~\ref{thm:DecompositionSpacesAsTemperedDistributions}
(with $\CalQ=\left(Q_{i}\right)_{i\in I}=\left(T_{i}Q_{i}'+b_{i}\right)_{i\in I}:=\CalB=\left(B_{n}\right)_{n\in\N_{0}}=\bigl(\left(2^{n}\cdot\identity\right)B_{n}'\bigr)_{n\in\N_{0}}$
and $I_{0}=I=\N_{0}$, with $B_{n}'$ as in Lemma~\ref{lem:InhomogeneousBesovCoveringAlmostStructured})
is well-suited for showing that every $f\in\FourierDecompSp{\CalB}p{\ell_{v^{\left(\gamma\right)}}^{q}}$
extends to a tempered distribution $\widetilde{f}\in\Schwartz'\left(\R^{\dimension}\right)$.

It is not hard to verify that the partition of unity $\Phi$ defined
above is regular (see Definition~\ref{def:RegularCoveringRegularPartitionOfUnity}).
Thus, we only need to prove existence of $N\in\N_{0}$ with $w^{\left(N\right)}\cdot c\in\ell^{1}\left(\N_{0}\right)$
for all ${c=\left(c_{n}\right)_{n\in\N_{0}}\in\ell_{v^{\left(\gamma\right)}}^{q}\left(\N_{0}\right)}$,
where
\[
w_{n}^{\left(N\right)}=\left|\det T_{n}\right|^{1/p}\cdot\max\left\{ 1,\left\Vert T_{n}^{-1}\right\Vert ^{\dimension+1}\right\} \cdot\left[\inf_{\xi\in B_{n}^{\ast}}\left(1+\left|\xi\right|\right)\right]^{-N}.
\]
Because of $\ell_{v^{\left(\gamma\right)}}^{q}\left(\N_{0}\right)\hookrightarrow\ell_{v^{\left(\gamma\right)}}^{\infty}\left(\N_{0}\right)$,
it suffices to show $w^{\left(N\right)}\cdot\left(v^{\left(\gamma\right)}\right)^{-1}\in\ell^{1}\left(\N_{0}\right)$
for some $N\in\N_{0}$.

Now, for $n\in\N$, we saw in equation~(\ref{eq:InhomogeneousBesovClusterEstimate})
that $n^{\ast}\subset\left\{ 0\right\} \cup\left\{ n-3,\dots,n+3\right\} $.
Furthermore, directly after that equation, we showed $0^{\ast}\subset\left\{ 0,1,2,3\right\} $,
which implies $0\notin n^{\ast}$ for $n\in\N_{>3}$. Thus, $n^{\ast}\subset\left\{ n-3,\dots,n+3\right\} $
for $n\in\N_{>3}$, which easily implies $\left|\xi\right|\geq2^{\left(n-3\right)-2}=2^{n-5}$
for all $\xi\in B_{n}^{\ast}$. Furthermore, we have $\left\Vert T_{n}^{-1}\right\Vert =\left\Vert 2^{-n}\cdot\identity\right\Vert =2^{-n}\leq1$
for all $n\in\N_{0}$, so that we get
\[
w_{n}^{\left(N\right)}=2^{\dimension n/p}\cdot\left[\inf_{\xi\in B_{n}^{\ast}}\left(1+\left|\xi\right|\right)\right]^{-N}\leq2^{\dimension n/p}\cdot\left(2^{n-5}\right)^{-N}=2^{5N}\cdot2^{n\left(\frac{\dimension}{p}-N\right)}\qquad\forall\,n\in\N_{>3}\,.
\]
Thus,
\[
0\leq w_{n}^{\left(N\right)}\,\big/\,v_{n}^{\left(\gamma\right)}\leq2^{-\gamma n}\cdot2^{5N}\cdot2^{n\left(\frac{\dimension}{p}-N\right)}=2^{5N}\cdot2^{n\left(\frac{\dimension}{p}-\gamma-N\right)}\qquad\forall\,n\in\N_{>3},
\]
which implies $\left(v^{\left(\gamma\right)}\right)^{-1}\cdot w^{\left(N\right)}\in\ell^{1}\left(\N_{0}\right)$,
as soon as $\frac{\dimension}{p}-\gamma-N<0$, i.e.\@ as soon as
$N>\frac{\dimension}{p}-\gamma$.
\end{proof}
As a corollary, we can now show that the Fourier transform also yields
an isomorphism between the classical ``space-side'' $\alpha$-modulation
spaces and our ``Fourier-side'' variants of these spaces.
\begin{cor}
\label{cor:FourierAlphaModulationAsTemperedDistributions}Let $\dimension\in\N$,
$p,q\in\left(0,\infty\right]$, $\alpha\in\left[0,1\right]$ and $\gamma\in\R$.
Then, the map
\[
\Fourier:\AlphaModSpace pq{\gamma}{\alpha}\left(\R^{\dimension}\right)\to\FourierDecompSp{\smash{\CalQ^{\left(\alpha\right)}}}p{\smash{\ell_{w^{\left(\gamma/\left(1-\alpha\right)\right)}}^{q}}}=\AlphaModSpace pq{\Fourier,\gamma}{\alpha}\left(\R^{\dimension}\right),f\mapsto\widehat{f}|_{\TestFunctionSpace{\R^{\dimension}}}
\]
yields an isomorphism of (quasi)-Banach spaces, where the $\alpha$-modulation
space $\AlphaModSpace pq{\gamma}{\alpha}\left(\R^{\dimension}\right)$
is as defined in \cite[Definition 2.4]{BorupNielsenAlphaModulationSpaces}.

In particular, Theorems \ref{thm:AlphaModulationEmbeddings}, \ref{thm:AlphaModulationBesovEmbedding},
and \ref{thm:InhomBesovEmbeddings} concerning the embeddings between
the ``Fourier-side'' $\alpha$-modulation spaces apply just as well
to the classical ``space-side'' $\alpha$-modulation spaces.
\end{cor}

\begin{proof}
For $\alpha=1$, the claim was just shown in Lemma~\ref{lem:FourierInhomogeneousBesovAsTemperedDistributions}.
Thus, we can assume $\alpha\in\left[0,1\right)$.

Note that we have (by definition) that $\AlphaModSpace pq{\Fourier,\gamma}{\alpha}\left(\R^{\dimension}\right)=\FourierDecompSp{\smash{\CalQ_{r}^{\left(\alpha\right)}}}p{\smash{\ell_{w^{\left(\gamma/\left(1-\alpha\right)\right)}}^{q}}}\vphantom{\CalQ_{r}^{\left(\alpha\right)}}$,
for $r>r_{0}\left(\dimension,\alpha\right)$, with $r_{0}\left(\dimension,\alpha\right)$
as in Theorem~\ref{thm:AlphaCoveringExistence}. But the covering
$\CalQ_{r}^{\left(\alpha\right)}$ is exactly the same covering which
is used in \cite[Definition 2.4]{BorupNielsenAlphaModulationSpaces}
to define the (space-side) $\alpha$-modulation spaces $\AlphaModSpace pq{\gamma}{\alpha}\left(\R^{\dimension}\right)$.
Thus, let $\Phi=\left(\varphi_{k}\right)_{k\in\Z^{\dimension}\setminus\left\{ 0\right\} }$
be an $L^{p}$-BAPU for $\CalQ_{r}^{\left(\alpha\right)}$, jointly
for all $p\in\left(0,\infty\right]$. Existence of $\Phi$ is ensured
by Theorem~\ref{thm:AlmostStructuredAdmissibleAdmitsBAPU}, since
$\CalQ_{r}^{\left(\alpha\right)}$ is an (almost) structured covering
of $\R^{\dimension}$, by Lemma~\ref{lem:AlphaModulationCoveringIsStructured}.
Using such a system $\Phi$, the (quasi)-norm on $\AlphaModSpace pq{\gamma}{\alpha}\left(\R^{\dimension}\right)$,
as defined in \cite[Definition 2.4]{BorupNielsenAlphaModulationSpaces},
is given by
\[
\left\Vert f\right\Vert _{\AlphaModSpace pq{\gamma}{\alpha}}=\left\Vert \left(\left\langle \xi_{k}\right\rangle ^{\gamma}\cdot\left\Vert \Fourier^{-1}\left(\varphi_{k}\cdot\smash{\widehat{f}}\,\,\right)\right\Vert _{L^{p}}\right)_{k\in\Z^{\dimension}\setminus\left\{ 0\right\} }\right\Vert _{\ell^{q}},
\]
where $\xi_{k}\in Q_{r,k}^{\left(\alpha\right)}$ can be selected
arbitrarily. Note that we have $\xi_{k}:=\left|k\right|^{\alpha_{0}}k\in Q_{r,k}^{\left(\alpha\right)}$
(with $\alpha_{0}=\frac{\alpha}{1-\alpha}$), and that $\left|\xi_{k}\right|=\left|k\right|^{\alpha_{0}+1}=\left|k\right|^{\frac{1}{1-\alpha}}\asymp\left\langle k\right\rangle ^{\frac{1}{1-\alpha}}$
(since $k\in\Z^{\dimension}\setminus\left\{ 0\right\} $). In particular,
$\left|\xi_{k}\right|\gtrsim1$, and hence $\left\langle \xi_{k}\right\rangle \asymp\left|\xi_{k}\right|\asymp\left\langle k\right\rangle ^{1/\left(1-\alpha\right)}$,
so that we get $\left\langle \xi_{k}\right\rangle ^{\gamma}\asymp\left\langle k\right\rangle ^{\gamma/\left(1-\alpha\right)}\asymp w_{k}^{\left(\gamma/\left(1-\alpha\right)\right)}$.
All in all, it is now not hard to see that we have
\begin{equation}
\left\Vert f\right\Vert _{\AlphaModSpace pq{\gamma}{\alpha}}\asymp\left\Vert \left(\left\Vert \Fourier^{-1}\left(\varphi_{k}\cdot\smash{\widehat{f}}\,\,\right)\right\Vert _{L^{p}}\right)_{k\in\Z^{\dimension}\setminus\left\{ 0\right\} }\right\Vert _{\ell_{w^{\left(\gamma\right)}}^{q}}\!\!\!=\left\Vert \widehat{f}|_{\TestFunctionSpace{\R^{\dimension}}}\right\Vert _{\FourierDecompSp{\CalQ_{r}^{\left(\alpha\right)}}p{\ell_{w^{\left(\gamma/\left(1-\alpha\right)\right)}}^{q}}}\label{eq:FourierAlphaModulationNormEquivalence}
\end{equation}
for all $f\in\Schwartz'\left(\R^{\dimension}\right)$, which shows
that the map $\Fourier$ defined in the statement of the corollary
is a well-defined isomorphism of quasi-normed vector spaces onto a
certain subspace of $\FourierDecompSp{\smash{\CalQ_{r}^{\left(\alpha\right)}}}p{\smash{\ell_{w^{\left(\gamma/\left(1-\alpha\right)\right)}}^{q}}}\vphantom{\CalQ_{r}^{\left(\alpha\right)}}$
(namely, onto the range of $\Fourier$).

Thus, to complete the proof, we only have to show that $\Fourier$
is surjective. To this end, it suffices to show that each distribution
$f\in\FourierDecompSp{\smash{\CalQ_{r}^{\left(\alpha\right)}}}p{\smash{\ell_{w^{\left(\gamma/\left(1-\alpha\right)\right)}}^{q}}}\vphantom{\CalQ_{r}^{\left(\alpha\right)}}=\AlphaModSpace pq{\Fourier,\gamma}{\alpha}\left(\R^{\dimension}\right)$
can be extended to a \emph{tempered} distribution $\widetilde{f}\in\Schwartz'\left(\R^{\dimension}\right)$.
Indeed, if this is the case, equation~(\ref{eq:FourierAlphaModulationNormEquivalence})
shows 
\[
\left\Vert \Fourier^{-1}\widetilde{f}\,\right\Vert _{\AlphaModSpace pq{\gamma}{\alpha}}\asymp\left\Vert \widetilde{f}|_{\TestFunctionSpace{\R^{\dimension}}}\right\Vert _{\AlphaModSpace pq{\Fourier,\gamma}{\alpha}}=\left\Vert f\right\Vert _{\AlphaModSpace pq{\Fourier,\gamma}{\alpha}}<\infty\,,
\]
and hence $g:=\Fourier^{-1}\widetilde{f}\in\AlphaModSpace pq{\gamma}{\alpha}\left(\R^{\dimension}\right)$
with $\Fourier g=\widetilde{f}|_{\TestFunctionSpace{\R^{\dimension}}}=f$,
so that the map $\Fourier$ (as defined in the present corollary)
is surjective.

But Theorem~\ref{thm:AlphaModulationBesovEmbedding} shows that we
have $\AlphaModSpace pq{\Fourier,\gamma}{\alpha}\left(\R^{\dimension}\right)\hookrightarrow\AlphaModSpace pq{\Fourier,\sigma}1\left(\R^{\dimension}\right)$
for a suitable $\sigma\in\R$, so that every $f\in\AlphaModSpace pq{\Fourier,\gamma}{\alpha}\left(\R^{\dimension}\right)$
also satisfies $f\in\AlphaModSpace pq{\Fourier,\sigma}1\left(\R^{\dimension}\right)$.
But in view of Lemma~\ref{lem:FourierInhomogeneousBesovAsTemperedDistributions},
this implies that $f$ extends to a tempered distribution, as desired.
\end{proof}

\subsection{Embeddings between homogeneous and inhomogeneous Besov spaces}

\label{subsec:HomogeneousInhomogeneousBesovEmbeddings}In the previous
subsections, the coverings which we considered were very \emph{compatible}:
They all covered the same set (namely, $\R^{\dimension}$) and $\CalQ$
was almost subordinate to $\CalP$ and relatively $\CalP$-moderate,
or vice versa. Thus, we could completely characterize the existence
of the associated embeddings, but we could not indicate in how far
our theory also applies in cases where the two coverings are more
incompatible to each other.

Therefore, in this subsection, we discuss one such example of more
incompatible coverings: We consider embeddings between inhomogeneous
and homogeneous (Fourier-side) Besov spaces. The inhomogeneous Besov
covering was introduced in Definition~\ref{def:InhomogeneousBesovCovering};
it consists of the dyadic annuli $B_{n}=B_{2^{n+2}}\left(0\right)\setminus\overline{B_{2^{n-2}}}\left(0\right)$
for $n\in\N$, but the low-frequency part is covered by one large
ball $B_{0}=B_{4}\left(0\right)$. In contrast, the homogeneous Besov
covering only consists of dyadic annuli:
\begin{defn}
\label{def:HomogeneousBesovCovering}For $\dimension\in\N$ and $n\in\Z$,
the \textbf{$n$-th dyadic ring} is $\dot{B}_{n}:=B_{2^{n+2}}\left(0\right)\setminus\overline{B_{2^{n-2}}}\left(0\right)\subset\R^{\dimension}$.
We define the ($\dimension$-dimensional) \textbf{homogeneous Besov
covering} as $\dot{\CalB}:=\left(\smash{\dot{B}_{n}}\right)_{n\in\Z}$.
\end{defn}

Our first aim is to show that $\dot{\CalB}$ is an almost structured
covering of $\R^{\dimension}\setminus\left\{ 0\right\} $. In fact,
$\dot{\CalB}$ is a structured covering:
\begin{lem}
\label{lem:HomogeneousBesovCoveringStructured}Let $\dimension\in\N$.
Then $\dot{\CalB}=\left(T_{n}Q+b_{n}\right)_{n\in\Z}$ is a structured
admissible covering of $\R^{\dimension}\setminus\left\{ 0\right\} $,
with parametrization given by
\[
T_{n}:=2^{n}\cdot\identity\quad\text{and}\quad b_{n}:=0\quad\text{for }n\in\Z,\qquad\text{and}\qquad Q:=B_{4}\left(0\right)\setminus\overline{B_{1/4}}\left(0\right)\,.
\]
Furthermore, for any $\gamma\in\R$, the weight $u^{\left(\gamma\right)}:=\left(2^{n\gamma}\right)_{n\in\Z}$
is $\dot{\CalB}$-moderate.
\end{lem}

\begin{proof}
With $T_{n},b_{n},Q$ as in the statement of the lemma, we clearly
have $\dot{B}_{n}=T_{n}Q+b_{n}$ for all $n\in\Z$.

Now, let $P:=B_{2}\left(0\right)\setminus\overline{B_{1/2}}\left(0\right)$.
It is clear that $P,Q\subset\R^{\dimension}$ are open and bounded,
with $\overline{P}\subset Q$. Furthermore, 
\[
\CalP:=\left(P_{n}\right)_{n\in\Z}:=\left(T_{n}P+b_{n}\right)_{n\in\Z}=\left(B_{2^{n+1}}\left(0\right)\setminus\overline{B_{2^{n-1}}}\left(0\right)\right)_{n\in\Z}
\]
is a covering of $\R^{\dimension}\setminus\left\{ 0\right\} $. Indeed,
clearly $P_{n}\subset\R^{\dimension}\setminus\left\{ 0\right\} $
for all $n\in\Z$. Conversely, for arbitrary $\xi\in\R^{\dimension}\setminus\left\{ 0\right\} $,
we can choose $n\in\Z$ maximal with $2^{n}\leq\left|\xi\right|$,
since the set of those $n\in\Z$ is nonempty (this uses $\left|\xi\right|>0$
and $2^{n}\xrightarrow[n\to-\infty]{}0$), closed and bounded from
above. Then $2^{n-1}<2^{n}\leq\left|\xi\right|<2^{n+1}$, which means
$\xi\in P_{n}$.

Since we have $P_{n}\subset\dot{B}_{n}\subset\R^{\dimension}\setminus\left\{ 0\right\} $
for all $n\in\Z$, we see that $\dot{\CalB}$ is also a covering of
$\R^{\dimension}\setminus\left\{ 0\right\} $. To show that $\dot{\CalB}$
is a structured admissible covering of $\R^{\dimension}\setminus\left\{ 0\right\} $,
it remains to show that $\dot{\CalB}$ is admissible and that 
\begin{equation}
C_{\dot{\CalB}}:=\sup_{n\in\Z}\,\,\sup_{k\in n^{\ast}}\left\Vert T_{n}^{-1}T_{k}\right\Vert =\sup_{n\in\Z}\,\,\sup_{k\in n^{\ast}}2^{k-n}\label{eq:BesovHomogeneousStructuredEstimate}
\end{equation}
is finite. But for $n\in\Z$ and $k\in n^{\ast}$, there is some $\xi\in\dot{B}_{k}\cap\dot{B}_{n}$,
which implies $2^{n-2}<\left|\xi\right|<2^{k+2}$ and $2^{k-2}<\left|\xi\right|<2^{n+2}$.
Hence, $n-2<k+2$ and $k-2<n+2$, which implies $n-4<k<n+4$. Since
all of these quantities are integers, we conclude $k\in\left\{ n-3,\dots,n+3\right\} $
and thus
\begin{equation}
n^{\ast}\subset\left\{ n-3,\dots,n+3\right\} ,\label{eq:BesovHomogeneousClusterSet}
\end{equation}
which yields $\left|n^{\ast}\right|\leq7$. Since this holds for all
$n\in\Z$, $\dot{\CalB}$ is admissible.

But we also get $2^{k-n}\leq2^{\left(n+3\right)-n}=2^{3}=8$. Thanks
to equation~(\ref{eq:BesovHomogeneousStructuredEstimate}), this
yields $C_{\dot{\CalB}}\leq8<\infty$, as desired.

\medskip{}

Finally, we also get the $\dot{\CalB}$-moderateness of $u^{\left(\gamma\right)}$:
For $n\in\Z$ and $k\in n^{\ast}$, we saw above that $\left|k-n\right|\leq3$.
Therefore, $\bigl|u_{k}^{\left(\gamma\right)}\,\big/\,u_{n}^{\left(\gamma\right)}\bigr|=2^{\gamma\left(k-n\right)}\leq2^{\left|\gamma\left(k-n\right)\right|}\leq2^{3\left|\gamma\right|}$,
so that we get $C_{u^{\left(\gamma\right)},\dot{\CalB}}\leq2^{3\left|\gamma\right|}<\infty$.
\end{proof}
In view of the preceding lemma, the following decomposition spaces
are well-defined:
\begin{defn}
\label{def:HomogeneousFouriersideBesovSpaces}For $\dimension\in\N$,
$p,q\in\left(0,\infty\right]$ and $\gamma\in\R$, we define the \textbf{homogeneous
(Fourier-side) Besov space} with exponents $p,q$ and smoothness parameter
$\gamma$ as
\[
\BesovHom pq{\Fourier,\gamma}\left(\R^{\dimension}\right):=\FourierDecompSp{\smash{\dot{\CalB}}}p{\ell_{u^{\left(\gamma\right)}}^{q}}\subset\DistributionSpace{\R^{\dimension}\setminus\left\{ 0\right\} },\qquad\text{with}\qquad u^{\left(\gamma\right)}=\left(2^{k\gamma}\right)_{k\in\Z}\,.\qedhere
\]
\end{defn}

\begin{rem*}
In addition to the homogeneous Fourier-side Besov spaces that we defined
above, there are other function spaces that one can define using the
homogeneous Besov covering $\dot{\CalB}$. All of these fit into the
(slightly generalized) framework of \emph{general decomposition} spaces
as introduced in \cite{DecompositionSpaces1}. We mention here the
three most important examples of such spaces: If $\Phi=\left(\varphi_{k}\right)_{k\in\Z}$
is a suitable partition of unity subordinate to $\dot{\CalB}$, if
$u^{\left(\gamma\right)}=\left(2^{k\gamma}\right)_{k\in\Z}$ is as
above, and if $p,q\in\left(0,\infty\right]$, then we define
\[
\begin{alignedat}{2}\left\Vert f\right\Vert _{\BesovHom pq{\gamma}} & :=\!\!\left\Vert \left(\left\Vert \Fourier^{-1}\left(\varphi_{k}\cdot\smash{\widehat{f}}\:\right)\right\Vert _{L^{p}}\right)_{k\in\Z}\right\Vert _{\ell_{u^{\left(\gamma\right)}}^{q}}\in\left[0,\infty\right] & \qquad & \text{ for }f\in\Schwartz'\left(\R^{\dimension}\right)\,\big/\,\CalP,\\
\left\Vert f\right\Vert _{\BesovHom pq{\Fourier,\gamma}} & :=\!\!\left\Vert \left(\left\Vert \Fourier^{-1}\left(\varphi_{k}\cdot f\:\right)\right\Vert _{L^{p}}\right)_{k\in\Z}\right\Vert _{\ell_{u^{\left(\gamma\right)}}^{q}}\in\left[0,\infty\right] & \qquad & \text{ for }f\in\DistributionSpace{\R^{\dimension}\setminus\left\{ 0\right\} },\\
\left\Vert f\right\Vert _{\dot{K}_{\gamma}^{p,q}} & :=\!\!\left\Vert \left(\left\Vert \varphi_{k}\cdot f\,\right\Vert _{L^{p}}\right)_{k\in\Z}\right\Vert _{\ell_{u^{\left(\gamma\right)}}^{q}}\!\!\in\left[0,\infty\right] & \qquad & \text{ for }f:\R^{\dimension}\to\Compl\text{ measurable}.
\end{alignedat}
\]
Here, $\CalP$ denotes the space of polynomials. The corresponding
function spaces $\BesovHom pq{\gamma}\left(\R^{\dimension}\right)$,
$\BesovHom pq{\Fourier,\gamma}\left(\R^{\dimension}\right)$, and
$\dot{K}_{\gamma}^{p,q}\left(\R^{\dimension}\right)$ then simply
consist of those functions/distributions for which the respective
quasi-norm is finite.

The spaces $\BesovHom pq{\gamma}\left(\R^{\dimension}\right)$ are
the \textbf{homogeneous Besov spaces}; they are studied in greater
detail for example in \cite{PeetreNewThoughtsOnBesovSpaces}, \cite[Chapter 5]{TriebelTheoryOfFunctionSpaces},
and \cite{HerzLipschitzSpaces}. The last mentioned paper \cite{HerzLipschitzSpaces}
by Herz also introduces the spaces $\dot{K}_{\gamma}^{p,q}\left(\R^{\dimension}\right)$,
which are nowadays called \textbf{(homogeneous) Beurling-Herz spaces};
see also \cite[Chapter 3, Example 7]{PeetreNewThoughtsOnBesovSpaces}.
We mention that with the definition given in \cite{HerzLipschitzSpaces},
it is not quite trivial to see that the spaces $K_{p,q}^{\gamma}$
defined there really coincide with the spaces $\dot{K}_{\gamma}^{p,q}\left(\R^{\dimension}\right)$
as defined above. One can either trust Peetre (see \cite[Chapter 3, Example 7]{PeetreNewThoughtsOnBesovSpaces})
on this, or one can use \cite[Proposition 2.4]{HerzLipschitzSpaces}
with the kernel $\widehat{\kappa}\left(\xi;h\right):=\Indicator_{\left(2^{-1},2\right)}\left(\left|\xi\right|\cdot\left|h\right|\right)$
to see that this is indeed the case.

Of course, all of these spaces are related: As shown in \cite[Lemma 6.2.2]{VoigtlaenderPhDThesis},
the Fourier transform induces an isomorphism of (quasi)-Banach spaces
\[
\begin{split}\Fourier:\BesovHom pq{\gamma}\left(\R^{\dimension}\right)\subset\Schwartz'\left(\R^{\dimension}\right)/\mathcal{P} & \to\FourierDecompSp{\smash{\dot{\CalB}}}p{\ell_{v^{\left(\gamma\right)}}^{q}}=\BesovHom pq{\Fourier,\gamma}\left(\R^{\dimension}\right)\subset\DistributionSpace{\R^{\dimension}\setminus\left\{ 0\right\} },\\
f+\mathcal{P} & \mapsto\widehat{f}|_{\TestFunctionSpace{\R^{\dimension}\setminus\left\{ 0\right\} }}
\end{split}
\tag{\ensuremath{\ast}}
\]
between the homogeneous Besov space $\BesovHom pq{\gamma}\left(\R^{\dimension}\right)$
and the homogeneous \emph{Fourier-side} Besov space $\BesovHom pq{\Fourier,\gamma}\left(\R^{\dimension}\right)$.
This justifies our name for these spaces. Furthermore, one of the
main results of the paper \cite{HerzLipschitzSpaces} by Herz is \cite[Proposition 3.1]{HerzLipschitzSpaces},
where it is shown for $q\in\left[1,\infty\right]$, $p\in\left[1,2\right]$
and $\gamma\in\R$ that the Fourier transform maps $\BesovHom pq{\gamma}\left(\R^{\dimension}\right)$
boundedly into $\dot{K}_{\gamma}^{p',q}\left(\R^{\dimension}\right)$.
In view of $\left(\ast\right)$, this is equivalent to the embedding
$\BesovHom pq{\Fourier,\gamma}\left(\R^{\dimension}\right)\hookrightarrow\dot{K}_{\gamma}^{p',q}\left(\R^{\dimension}\right)$.
From a modern perspective, this statement is very easy to obtain:
The Hausdorff-Young inequality shows 
\[
\left\Vert \varphi_{k}\cdot\Fourier f\right\Vert _{L^{p'}}=\left\Vert \Fourier\Fourier^{-1}\left[\varphi_{k}\cdot\Fourier f\right]\right\Vert _{L^{p'}}\leq\left\Vert \Fourier^{-1}\left[\varphi_{k}\cdot\smash{\widehat{f}}\,\,\right]\right\Vert _{L^{p}}\:.
\]
With this observation, the boundedness of $\Fourier:\BesovHom pq{\gamma}\left(\R^{\dimension}\right)\to\dot{K}_{\gamma}^{p',q}\left(\R^{\dimension}\right)$
is a direct consequence of the definitions of the involved spaces.
We emphasize, however, that the proof given in \cite{HerzLipschitzSpaces}
is very different and much more involved. The reason for this is that
Herz uses a definition of the homogeneous Besov spaces $\BesovHom pq{\gamma}\left(\R^{\dimension}\right)$
in terms of difference quotients; the resulting spaces are called
(homogeneous) \textbf{Lipschitz spaces} and denoted by $\dot{\Lambda}_{\gamma}^{p,q}\left(\R^{\dimension}\right)$.
The equality $\BesovHom pq{\gamma}\left(\R^{\dimension}\right)=\dot{\Lambda}_{\gamma}^{p,q}\left(\R^{\dimension}\right)$,
i.e., the equivalence of the definition using difference quotients
with the Fourier-analytic definition of $\BesovHom pq{\gamma}\left(\R^{\dimension}\right)$
from above—while well-known nowadays—is not trivial. It can be found
for example in \cite[Theorem 2 in Section 5.2.3]{TriebelTheoryOfFunctionSpaces}.
We also mention the paper \cite{JohnsonTemperatures} as an early
source in which embeddings between homogeneous Besov spaces (again
under the name Lipschitz spaces) have been studied; see \cite[Theorem 8]{JohnsonTemperatures}.

Finally, we emphasize that we are mainly interested in the \emph{space-side
homogeneous} Besov spaces $\BesovHom pq{\gamma}\left(\R^{\dimension}\right)$,
and their inclusion relations with the \emph{space-side inhomogeneous}
Besov spaces $\BesovInhom pq{\gamma}\left(\R^{\dimension}\right)$.
But in view of equation~$\left(\ast\right)$, this is equivalent
to studying embeddings between the \emph{Fourier-side} counterparts
of these spaces, which are somewhat more convenient to handle. Thus,
we are not really interested in the Fourier-side spaces themselves,
but only use them as a convenient means of studying their space-side
counterparts.
\end{rem*}
Our next step is to verify that the homogeneous Besov covering $\dot{\CalB}$
is almost subordinate to the inhomogeneous Besov covering $\CalB$.
\begin{lem}
\label{lem:HomogeneousInhomogeneousBesovSubordinateness}Let $\dimension\in\N$
be arbitrary. The homogeneous Besov covering $\dot{\CalB}=\left(\smash{\dot{B}_{k}}\right)_{k\in\Z}$
is almost subordinate to the inhomogeneous Besov covering $\CalB=\left(B_{n}\right)_{n\in\N_{0}}$,
and we have
\begin{equation}
\left\{ n\right\} \subset\left\{ k\in\Z\with\smash{\dot{B}_{k}}\cap B_{n}\neq\emptyset\right\} \subset\left\{ n-3,\dots,n+3\right\} \qquad\forall\,n\in\N,\label{eq:BesovHomogeneousInhomogeneousIntersectionHighFrequency}
\end{equation}
and
\begin{equation}
\Z_{\leq1}\subset\left\{ k\in\Z\with\smash{\dot{B}_{k}}\cap B_{0}\neq\emptyset\right\} \subset\Z_{\leq3}.\qedhere\label{eq:BesovHomogeneousInhomogeneousIntersectionLowFrequency}
\end{equation}
\end{lem}

\begin{proof}
For $k\in\N$, we simply have $\dot{B}_{k}=B_{k}=B_{k}^{0\ast}$.
But for $k\in\Z_{\leq0}$, we have 
\[
\dot{B}_{k}=B_{2^{k+2}}\left(0\right)\setminus\overline{B_{2^{k-2}}}\left(0\right)\subset B_{2^{2}}\left(0\right)=B_{4}\left(0\right)=B_{0}=B_{0}^{0\ast}.
\]
All in all, we have shown that $\dot{\CalB}$ is almost subordinate
to $\CalB$ with $k\left(\smash{\dot{\CalB},\CalB}\right)=0$.

Now, for $n\in\N$ and $k\in\Z$ with $\emptyset\neq\dot{B}_{k}\cap B_{n}=\dot{B}_{k}\cap\dot{B}_{n}$,
equation~(\ref{eq:BesovHomogeneousClusterSet}) shows that we have
$k\in\left\{ n-3,\dots,n+3\right\} $, which establishes equation~(\ref{eq:BesovHomogeneousInhomogeneousIntersectionHighFrequency}).

Finally, for $k\in\Z$ with $\dot{B}_{k}\cap B_{0}\neq\emptyset$,
there is some $\xi\in B_{0}\cap\dot{B}_{k}$, which implies $2^{k-2}<\left|\xi\right|<4=2^{2}$
and hence $k<4$, i.e.\@ $k\leq3$, which establishes the right inclusion
in equation~(\ref{eq:BesovHomogeneousInhomogeneousIntersectionLowFrequency}).

Conversely, for $k\in\Z_{\leq1}$, we have $\xi:=2^{k}\cdot\left(1,0,\dots,0\right)\in\dot{B}_{k}$
and $\left|\xi\right|\leq2^{1}<4$, which implies $\xi\in B_{0}\cap\dot{B}_{k}\neq\emptyset$.
This establishes the left inclusion in equation~(\ref{eq:BesovHomogeneousInhomogeneousIntersectionLowFrequency}).
\end{proof}
Now, we are in a position to consider embeddings between the homogeneous
and the inhomogeneous (Fourier-side) Besov spaces. We remark that
the following theorem already appeared in my PhD thesis, \cite[Theorem 6.2.6]{VoigtlaenderPhDThesis}.
\begin{thm}
\label{thm:InhomogeneousIntoHomogeneousBesov}Let $\dimension\in\N$,
$\gamma_{1},\gamma_{2}\in\R$ and $p_{1},p_{2},q_{1},q_{2}\in\left(0,\infty\right]$.
Then, the map
\[
\iota:\BesovInhom{p_{1}}{q_{1}}{\Fourier,\gamma_{1}}\left(\R^{\dimension}\right)\to\BesovHom{p_{2}}{q_{2}}{\Fourier,\gamma_{2}}\left(\R^{\dimension}\right),f\mapsto f|_{\TestFunctionSpace{\R^{\dimension}\setminus\left\{ 0\right\} }}
\]
is well-defined and bounded if we have $p_{1}\leq p_{2}$ and
\begin{equation}
\begin{cases}
\gamma_{2}\leq\gamma_{1}+\dimension\left(p_{2}^{-1}-p_{1}^{-1}\right), & \text{if }q_{1}\leq q_{2},\vspace{0.1cm}\\
\gamma_{2}<\gamma_{1}+\dimension\left(p_{1}^{-1}-p_{2}^{-1}\right), & \text{if }q_{1}>q_{2},
\end{cases}\label{eq:InhomIntoHomBesovSufficientLargeFrequencies}
\end{equation}
as well as
\begin{equation}
\begin{cases}
\gamma_{2}\geq\dimension\left(p_{2}^{-1}-\min\left\{ 1,p_{1}^{-1}\right\} \right), & \text{if }\UpperExpo{p_{1}}\leq q_{2},\vspace{0.1cm}\\
\gamma_{2}>\dimension\left(p_{2}^{-1}-\min\left\{ 1,p_{1}^{-1}\right\} \right), & \text{if }\UpperExpo{p_{1}}>q_{2}.
\end{cases}\label{eq:InhomIntoHomBesovSufficientSmallFrequencies}
\end{equation}

Conversely, if $\iota$ is bounded, then we have $p_{1}\leq p_{2}$,
condition~(\ref{eq:InhomIntoHomBesovSufficientLargeFrequencies})
is fulfilled and we have
\begin{equation}
\begin{cases}
\gamma_{2}\geq\dimension\left(p_{2}^{-1}-p_{1}^{-1}\right), & \text{if }p_{1}\leq q_{2},\vspace{0.1cm}\\
\gamma_{2}>\dimension\left(p_{2}^{-1}-p_{1}^{-1}\right), & \text{if }p_{1}>q_{2}.
\end{cases}\label{eq:InhomIntoHomBesovNecessarySmallFrequencies1}
\end{equation}
Furthermore, we also have
\begin{equation}
\begin{cases}
\gamma_{2}\geq\dimension\left(p_{2}^{-1}-1\right), & \text{if }q_{2}=\infty,\vspace{0.1cm}\\
\gamma_{2}>\dimension\left(p_{2}^{-1}-1\right), & \text{if }q_{2}<\infty.
\end{cases}\label{eq:InhomIntoHomBesovNecessarySmallFrequencies2}
\end{equation}

Finally, in case of $p_{1}=p_{2}$, we also have
\begin{equation}
\begin{cases}
\gamma_{2}\geq0, & \text{if }2\leq q_{2},\\
\gamma_{2}>0, & \text{if }2>q_{2}.
\end{cases}\qedhere\label{eq:InhomIntoHomBesovNecessarySmallFrequenciesKhinchin}
\end{equation}
\end{thm}

\begin{rem*}
Note that the conditions $p_{1}\leq p_{2}$ and (\ref{eq:InhomIntoHomBesovSufficientLargeFrequencies})
belong to the sufficient conditions as well as to the necessary conditions
derived above.

Thus, the only difference between the sufficient and the necessary
conditions is that condition~(\ref{eq:InhomIntoHomBesovSufficientSmallFrequencies})
is not fully equivalent to the conjunction of conditions~(\ref{eq:InhomIntoHomBesovNecessarySmallFrequencies1})
and (\ref{eq:InhomIntoHomBesovNecessarySmallFrequencies2}). There
are, however, many cases in which the equivalence \emph{does} hold:
For $p_{1}\in\left[2,\infty\right]$, we have $\UpperExpo{p_{1}}=p_{1}$
and $\min\left\{ 1,p_{1}^{-1}\right\} =p_{1}^{-1}$, so that conditions~(\ref{eq:InhomIntoHomBesovSufficientSmallFrequencies})
and (\ref{eq:InhomIntoHomBesovNecessarySmallFrequencies1}) are equivalent.
Furthermore, for $p_{1}\in\left(0,1\right]$, we have $\UpperExpo{p_{1}}=\infty$
and $\min\left\{ 1,p_{1}^{-1}\right\} =1$, so that conditions~(\ref{eq:InhomIntoHomBesovSufficientSmallFrequencies})
and (\ref{eq:InhomIntoHomBesovNecessarySmallFrequencies2}) are equivalent.

Thus, the only case in which we do \emph{not} get a complete characterization
is if we have
\begin{equation}
p_{1}\leq p_{2},\qquad p_{1}\in\left(1,2\right),\qquad\gamma_{2}=\dimension\left(p_{2}^{-1}-p_{1}^{-1}\right)\qquad\text{ and finally }\qquad p_{1}\leq q_{2}<p_{1}'.\label{eq:InhomIntoHomBesovEquivalenceGap}
\end{equation}
In this case, the necessary conditions are fulfilled, but the sufficient
conditions are not, so that our criteria cannot decide whether $\iota$
is or is not bounded.

In case of $p_{1}=p_{2}=:p$, this ``gap'' between necessary and
sufficient criteria can be narrowed further. In this case, condition~(\ref{eq:InhomIntoHomBesovNecessarySmallFrequenciesKhinchin})
shows that the only case in which we do \emph{not} get a complete
characterization is if we have
\begin{equation}
p\in\left(1,2\right),\qquad\gamma_{2}=0\qquad\text{ and finally }\qquad2\leq q_{2}<p'.\label{eq:InhomIntoHomBesovEquivalenceGapKhinchin}
\end{equation}

If conditions~(\ref{eq:InhomIntoHomBesovEquivalenceGap}) or (\ref{eq:InhomIntoHomBesovEquivalenceGapKhinchin})
are fulfilled, our criteria are inconclusive. But I do not know of
any more comprehensive results in this direction which could clarify
this ambiguity. In fact, the only statement about embeddings of inhomogeneous
into homogeneous Besov spaces of which I know is \cite[Theorem in §2.3.3]{TriebelTheoryOfFunctionSpaces2},
which shows that
\[
f\mapsto\left\Vert f\right\Vert _{L^{p}}+\left\Vert f\right\Vert _{\BesovHom pq{\gamma}}
\]
is an equivalent norm on the space $\BesovInhom pq{\gamma}\left(\R^{\dimension}\right)$,
\emph{provided that} $\gamma>\sigma_{p}:=\dimension\left(p^{-1}-1\right)_{+}$.
In particular, this yields $\BesovInhom pq{\gamma}\left(\R^{\dimension}\right)\hookrightarrow\BesovHom pq{\gamma}\left(\R^{\dimension}\right)$
for $\gamma>\dimension\left(p^{-1}-1\right)_{+}$. Note, however,
that our results from above show that this embedding is even true
for $\gamma\geq\dimension\left(p^{-1}-1\right)_{+}$, as long as $q\geq\UpperExpo p$.
In this sense, our results improve that given by Triebel in \cite{TriebelTheoryOfFunctionSpaces2}.

I would be thrilled about any further progress in closing the gap
described in equation~(\ref{eq:InhomIntoHomBesovEquivalenceGap}).
\end{rem*}
\begin{proof}[Proof of Theorem~\ref{thm:InhomogeneousIntoHomogeneousBesov}]
Define
\[
\CalQ=\left(Q_{i}\right)_{i\in I}=\left(T_{i}Q_{i}'+b_{i}\right)_{i\in I}:=\CalB=\left(B_{n}\right)_{n\in\N_{0}}=\bigl(\left(2^{n}\cdot\identity\right)B_{n}'\bigr)_{n\in\N_{0}}\,,
\]
and
\[
\CalP=\left(P_{j}\right)_{j\in J}=\left(S_{j}P_{j}'+c_{j}\right)_{j\in J}:=\dot{\CalB}=\left(\smash{\dot{B}_{k}}\right)_{k\in\Z}=\bigl(\left(\smash{2^{k}}\cdot\identity\right)Q\bigr)_{k\in\Z}\,,
\]
with $B_{n}'$ and $Q$ as in Lemmas \ref{lem:InhomogeneousBesovCoveringAlmostStructured}
and \ref{lem:HomogeneousBesovCoveringStructured}, respectively.
Furthermore, set $w:=v^{\left(\gamma_{1}\right)}=\left(2^{n\gamma_{1}}\right)_{n\in\N_{0}}$
and $v:=u^{\left(\gamma_{2}\right)}=\left(2^{k\gamma_{2}}\right)_{k\in\Z}$.
In view of Lemmas \ref{lem:InhomogeneousBesovCoveringAlmostStructured}
and \ref{lem:HomogeneousBesovCoveringStructured}, we see that $\CalQ$
and $\CalP$ are both almost structured coverings of $\CalO:=\R^{\dimension}$
and $\CalO':=\R^{\dimension}\setminus\left\{ 0\right\} $, respectively.
Furthermore, the same lemmas also show that $w$ and $v$ are $\CalQ$-moderate
and $\CalP$-moderate, respectively, so that the standing assumptions
of Section~\ref{sec:SummaryOfEmbeddingResults} (i.e., Assumption~\ref{assu:GeneralSummaryAssumptions})
are satisfied.

Finally, thanks to Lemma~\ref{lem:HomogeneousInhomogeneousBesovSubordinateness},
we know that $\CalP$ is almost subordinate to $\CalQ$ (but \emph{not}
relatively $\CalQ$-moderate). For the sufficient condition, we can
thus apply part~(\ref{enu:SummaryCoarseInFineSufficientNonModerate})
of Theorem~\ref{thm:SummaryCoarseIntoFine} to conclude that $\iota$
is well-defined and bounded if we have $p_{1}\leq p_{2}$ and if $K_{\UpperExpo{p_{1}},2}<\infty$
(in the notation of Theorem~\ref{thm:SummaryCoarseIntoFine}).

To verify this, we more generally investigate finiteness of
\begin{align*}
K\left(r,\alpha,\beta\right) & :=\left\Vert \left(w_{i}^{-1}\cdot\left|\det T_{i}\right|^{\alpha}\cdot\left\Vert \left(\smash{\left|\det S_{j}\right|^{\beta}}\cdot v_{j}\right)_{j\in J_{i}}\right\Vert _{\ell^{q_{2}\cdot\left(r/q_{2}\right)'}}\right)_{i\in I}\right\Vert _{\ell^{q_{2}\cdot\left(q_{1}/q_{2}\right)'}}\\
 & =\left\Vert \left(2^{n\dimension\alpha}\cdot2^{-n\gamma_{1}}\cdot\left\Vert \left(2^{k\dimension\beta}\cdot2^{k\gamma_{2}}\right)_{k\in\Z^{\left(n\right)}}\right\Vert _{\ell^{q_{2}\cdot\left(r/q_{2}\right)'}}\right)_{n\in\N_{0}}\right\Vert _{\ell^{q_{2}\cdot\left(q_{1}/q_{2}\right)'}}
\end{align*}
for arbitrary $r\in\left(0,\infty\right]$ and $\alpha,\beta\in\R$
and with
\begin{equation}
\Z^{\left(n\right)}:=\left\{ k\in\Z\with\dot{B_{k}}\cap B_{n}\neq\emptyset\right\} .\label{eq:InhomIntoHomBesovIntersectionSet}
\end{equation}

We first observe that $K\left(r,\alpha,\beta\right)$ is finite if
and only if we have 
\[
K_{0}\left(r,\alpha,\beta\right):=\left\Vert \left(2^{k\dimension\beta}\cdot2^{k\gamma_{2}}\right)_{k\in\Z^{\left(0\right)}}\right\Vert _{\ell^{q_{2}\cdot\left(r/q_{2}\right)'}}<\infty
\]
and if
\begin{align*}
\widetilde{K}\left(r,\alpha,\beta\right) & :=\left\Vert \left(2^{n\left(\dimension\alpha-\gamma_{1}\right)}\cdot\left\Vert \left(2^{k\left(\gamma_{2}+\dimension\beta\right)}\right)_{k\in\Z^{\left(n\right)}}\right\Vert _{\ell^{q_{2}\cdot\left(r/q_{2}\right)'}}\right)_{n\in\N}\right\Vert _{\ell^{q_{2}\cdot\left(q_{1}/q_{2}\right)'}}\\
\left({\scriptstyle \text{equation }\eqref{eq:BesovHomogeneousInhomogeneousIntersectionHighFrequency}}\right) & \asymp\left\Vert \left(2^{n\left(\dimension\alpha-\gamma_{1}\right)}\cdot2^{n\left(\gamma_{2}+\dimension\beta\right)}\right)_{n\in\N}\right\Vert _{\ell^{q_{2}\cdot\left(q_{1}/q_{2}\right)'}}=\left\Vert \left(2^{n\left(\gamma_{2}-\gamma_{1}+\dimension\left(\alpha+\beta\right)\right)}\right)_{n\in\N}\right\Vert _{\ell^{q_{2}\cdot\left(q_{1}/q_{2}\right)'}}
\end{align*}
is finite. Since equation~(\ref{eq:SpecialExponentFiniteness}) shows
that $q_{2}\cdot\left(q_{1}/q_{2}\right)'$ is finite if and only
if $q_{1}>q_{2}$, we see that $\widetilde{K}\left(r,\alpha,\beta\right)$
is finite if and only if we have
\begin{equation}
\begin{cases}
\gamma_{2}-\gamma_{1}+\dimension\left(\alpha+\beta\right)\leq0, & \text{if }q_{1}\leq q_{2},\\
\gamma_{2}-\gamma_{1}+\dimension\left(\alpha+\beta\right)<0, & \text{if }q_{1}>q_{2}.
\end{cases}\label{eq:BesovInhomogeneousIntoHomogeneousLargeFrequencyConditionGeneral}
\end{equation}

Furthermore, equation~(\ref{eq:BesovHomogeneousInhomogeneousIntersectionLowFrequency})
from Lemma~\ref{lem:HomogeneousInhomogeneousBesovSubordinateness}
yields $\Z_{\leq1}\subset\Z^{\left(0\right)}\subset\Z_{\leq3}$. Additionally,
by another application of equation~(\ref{eq:SpecialExponentFiniteness}),
we see that $q_{2}\cdot\left(r/q_{2}\right)'$ is finite if and only
if $r>q_{2}$, so that $K_{0}\left(r,\alpha,\beta\right)$ is finite
if and only if we have $\left\Vert \left(2^{k\left(\gamma_{2}+\dimension\beta\right)}\right)_{k\in\Z_{\leq1}}\right\Vert _{\ell^{q_{2}\cdot\left(r/q_{2}\right)'}}<\infty$
and thus if and only if
\begin{equation}
\begin{cases}
\gamma_{2}+\dimension\beta\geq0, & \text{if }r\leq q_{2},\\
\gamma_{2}+\dimension\beta>0, & \text{if }r>q_{2}.
\end{cases}\label{eq:BesovInhomogeneousIntoHomogeneousSmallFrequencyConditionGeneral}
\end{equation}
All in all, we see that $K\left(r,\alpha,\beta\right)$ is finite
if and only if conditions~(\ref{eq:BesovInhomogeneousIntoHomogeneousLargeFrequencyConditionGeneral})
and (\ref{eq:BesovInhomogeneousIntoHomogeneousSmallFrequencyConditionGeneral})
are satisfied.

Now, there are two cases. In case of $p_{1}\geq1$, the quantity $K_{\UpperExpo{p_{1}},2}$
from Theorem~\ref{thm:SummaryCoarseIntoFine} satisfies $K_{\UpperExpo{p_{1}},2}=K_{\UpperExpo{p_{1}},1}=K\left(\smash{\UpperExpo{p_{1}}},0,p_{1}^{-1}-p_{2}^{-1}\right)$.
Plugging these values into conditions~(\ref{eq:BesovInhomogeneousIntoHomogeneousLargeFrequencyConditionGeneral})
and (\ref{eq:BesovInhomogeneousIntoHomogeneousSmallFrequencyConditionGeneral}),
we see that these conditions are satisfied, thanks to our assumptions
(\ref{eq:InhomIntoHomBesovSufficientLargeFrequencies}) and (\ref{eq:InhomIntoHomBesovSufficientSmallFrequencies}).
To see this, note that we have $\min\left\{ 1,p_{1}^{-1}\right\} =p_{1}^{-1}$,
since $p_{1}\geq1$.

In case of $p_{1}<1$, the quantity $K_{\UpperExpo{p_{1}},2}$ from
Theorem~\ref{thm:SummaryCoarseIntoFine} satisfies $K_{\UpperExpo{p_{1}},2}=K\left(\smash{\UpperExpo{p_{1}}},p_{1}^{-1}-1,1-p_{2}^{-1}\right)$.
Plugging these values into conditions~(\ref{eq:BesovInhomogeneousIntoHomogeneousLargeFrequencyConditionGeneral})
and (\ref{eq:BesovInhomogeneousIntoHomogeneousSmallFrequencyConditionGeneral}),
we again see that these conditions are satisfied, thanks to our assumptions
(\ref{eq:InhomIntoHomBesovSufficientLargeFrequencies}) and (\ref{eq:InhomIntoHomBesovSufficientSmallFrequencies}).
To see this, note in this case that $\min\left\{ 1,p_{1}^{-1}\right\} =1$,
since $p_{1}\in\left(0,1\right)$.

\medskip{}

Now, let us derive the stated necessary conditions. To this end, assume
that $\iota$ is bounded. This easily implies boundedness of the map
$\theta$ from part~(\ref{enu:SummaryCoarseInFineNecessaryNonModerate})
of Theorem~\ref{thm:SummaryCoarseIntoFine}, since $\iota$ maps
each function $f\in\TestFunctionSpace{\CalO'}=\TestFunctionSpace{\R^{\dimension}\setminus\left\{ 0\right\} }$
to itself (interpreted as a distribution on $\R^{\dimension}\setminus\left\{ 0\right\} $).
Hence, part~(\ref{enu:SummaryCoarseInFineNecessaryNonModerate})
of Theorem~\ref{thm:SummaryCoarseIntoFine} yields $p_{1}\leq p_{2}$
and shows that we have $K\left(p_{1},0,p_{1}^{-1}-p_{2}^{-1}\right)=K_{p_{1},1}<\infty$.
In view of conditions (\ref{eq:BesovInhomogeneousIntoHomogeneousLargeFrequencyConditionGeneral})
and (\ref{eq:BesovInhomogeneousIntoHomogeneousSmallFrequencyConditionGeneral}),
this shows that conditions (\ref{eq:InhomIntoHomBesovSufficientLargeFrequencies})
and (\ref{eq:InhomIntoHomBesovNecessarySmallFrequencies1}) from the
statement of the theorem are indeed fulfilled.

\medskip{}

Next, if we have $p_{1}=p_{2}$, we can apply part~(\ref{enu:SummaryCoarseInFineNecessaryKhinchin})
of Theorem~\ref{thm:SummaryCoarseIntoFine}, which shows that $K_{2,1}=K\left(2,0,0\right)$
is finite. In view of condition~(\ref{eq:BesovInhomogeneousIntoHomogeneousSmallFrequencyConditionGeneral}),
we thus see that condition~(\ref{eq:InhomIntoHomBesovNecessarySmallFrequenciesKhinchin})
from the statement of the theorem is indeed fulfilled.

\medskip{}

It remains to show that condition~(\ref{eq:InhomIntoHomBesovNecessarySmallFrequencies2})
is fulfilled, regardless of whether $p_{1}=p_{2}$ or $p_{1}\neq p_{2}$.
This appears to be difficult, since we have exhausted the necessary
criteria provided by Theorem~\ref{thm:SummaryCoarseIntoFine}. But
in this case, it helps to invoke Lemma~\ref{lem:NecessaryConjugateCoarseInFineWithoutModerateness},
with our choices $\CalQ=\CalB$ and $\CalP=\dot{\CalB}$ from above
and with $J_{0}=J=\Z$. Indeed, as required in Lemma~\ref{lem:NecessaryConjugateCoarseInFineWithoutModerateness},
$\CalP_{J_{0}}$ is almost subordinate to $\CalQ$.

Furthermore, Lemma~\ref{lem:NecessaryConjugateCoarseInFineWithoutModerateness}
requires existence of a bounded map $\iota$ (which we will write
as $\iota_{0}$ from now on to avoid confusion with the map $\iota$
from the current theorem) which satisfies $\left\langle \iota_{0}f,\,\varphi\right\rangle _{\CalD'}=\left\langle f,\,\varphi\right\rangle _{\CalD'}$
for all $\varphi\in\TestFunctionSpace{\CalO\cap\CalO'}=\TestFunctionSpace{\R^{\dimension}\setminus\left\{ 0\right\} }$
and all $f\in\CalD_{K}^{\CalQ,p_{1},\ell_{w}^{q_{1}}}$, where $K$
is some subset of $\CalO=\R^{\dimension}$ (defined in Lemma~\ref{lem:NecessaryConjugateCoarseInFineWithoutModerateness}),
whose precise definition is immaterial for us.

In our present case, we can simply select $\iota_{0}:=\iota|_{\CalD_{K}^{\CalQ,p_{1},\ell_{w}^{q_{1}}}}$.
Indeed, for $f\in\CalD_{K}^{\CalQ,p_{1},\ell_{w}^{q_{1}}}\subset\TestFunctionSpace{\R^{\dimension}}$
(interpreted as a distribution on $\R^{\dimension}$) and $\varphi\in\TestFunctionSpace{\R^{\dimension}\setminus\left\{ 0\right\} }$,
we have
\[
\left\langle \iota f,\,\varphi\right\rangle _{\CalD'}=\left\langle f|_{\TestFunctionSpace{\R^{\dimension}\setminus\left\{ 0\right\} }},\,\varphi\right\rangle _{\CalD'}=\left\langle f,\,\varphi\right\rangle _{\CalD'},
\]
as required. Now, note that $B_{0}=B_{4}\left(0\right)$ intersects
$\R^{\dimension}\setminus\left\{ 0\right\} $, so that $B_{0}\cap\dot{B}_{k}\neq\emptyset$
for some $k\in\Z=J_{0}$; i.e., $0\in I_{0}$ with $I_{0}$ as in
Lemma~\ref{lem:NecessaryConjugateCoarseInFineWithoutModerateness}.
Since we also have $\delta_{0}\in\ell_{w}^{q_{1}}\left(\N_{0}\right)$,
Lemma~\ref{lem:NecessaryConjugateCoarseInFineWithoutModerateness}
shows (for $i=0$) that we have
\begin{align*}
\infty>u_{i} & =\left|\det T_{i}\right|^{p_{1}^{-1}-1}\cdot\left\Vert \left(v_{j}\cdot\left|\det S_{j}\right|^{1-p_{2}^{-1}}\right)_{j\in J_{0}\cap J_{i}}\right\Vert _{\ell^{q_{2}}}\\
\left({\scriptstyle \text{since }i=0,\text{ since }J_{0}=\Z\text{ and since }J_{i}=\Z^{\left(i\right)}}\right) & =\left\Vert \left(2^{j\gamma_{2}}\cdot2^{j\dimension\left(1-p_{2}^{-1}\right)}\right)_{j\in\Z^{\left(0\right)}}\right\Vert _{\ell^{q_{2}}},
\end{align*}
with $\Z^{\left(0\right)}$ as in equation~(\ref{eq:InhomIntoHomBesovIntersectionSet}).
But because of $\Z_{\leq1}\subset\Z^{\left(0\right)}\subset\Z_{\leq3}$,
we see that $u_{0}$ is finite if and only if condition~(\ref{eq:InhomIntoHomBesovNecessarySmallFrequencies2})
is satisfied.
\end{proof}
As our final result, we study the ``reverse'' direction of the preceding
theorem, i.e.\@ the existence of embeddings of the homogeneous into
the inhomogeneous Besov spaces. We remark that the following theorem
already appeared in a slightly different form in my PhD thesis, as
\cite[Theorem 6.2.3]{VoigtlaenderPhDThesis}.
\begin{thm}
\label{thm:HomogeneousIntoInhomogeneousBesov}Let $\dimension\in\N$,
$p_{1},p_{2},q_{1},q_{2}\in\left(0,\infty\right]$ and $\gamma_{1},\gamma_{2}\in\R$.

If we have $p_{1}\leq p_{2}$,
\begin{equation}
\begin{cases}
\gamma_{2}\leq\gamma_{1}+\dimension\left(p_{2}^{-1}-p_{1}^{-1}\right), & \text{if }q_{1}\leq q_{2},\vspace{0.1cm}\\
\gamma_{2}<\gamma_{1}+\dimension\left(p_{2}^{-1}-p_{1}^{-1}\right), & \text{if }q_{1}>q_{2}.
\end{cases}\label{eq:BesovHomIntoInhomHighFrequencyConditionSufficient}
\end{equation}
and
\begin{equation}
\begin{cases}
\gamma_{1}\leq\dimension\left(p_{1}^{-1}-p_{2}^{-1}\right), & \text{if }q_{1}\leq\LowerExpo{p_{2}},\vspace{0.1cm}\\
\gamma_{1}<\dimension\left(p_{1}^{-1}-p_{2}^{-1}\right), & \text{if }q_{1}>\LowerExpo{p_{2}},
\end{cases}\label{eq:BesovHomIntoInhomLowFrequencyConditionSufficient}
\end{equation}
then there is a bounded linear map
\[
\iota:\BesovHom{p_{1}}{q_{1}}{\Fourier,\gamma_{1}}\left(\R^{\dimension}\right)\hookrightarrow\BesovInhom{p_{2}}{q_{2}}{\Fourier,\gamma_{2}}\left(\R^{\dimension}\right)
\]
satisfying
\[
\left\langle \iota f,\,\varphi\right\rangle _{\CalD'}=\left\langle f,\,\varphi\right\rangle _{\CalD'}\qquad\forall\,\varphi\in\TestFunctionSpace{\R^{\dimension}\setminus\left\{ 0\right\} }\text{ and }f\in\BesovHom{p_{1}}{q_{1}}{\Fourier,\gamma_{1}}\left(\R^{\dimension}\right)
\]
and $\iota f=f$ for all $f\in\TestFunctionSpace{\R^{\dimension}\setminus\left\{ 0\right\} }$.

Conversely, if a map with these properties exist, then we necessarily
have $p_{1}\leq p_{2}$, condition~(\ref{eq:BesovHomIntoInhomHighFrequencyConditionSufficient})
is satisfied, and we have
\begin{equation}
\begin{cases}
\gamma_{1}\leq\dimension\left(p_{1}^{-1}-p_{2}^{-1}\right), & \text{if }q_{1}\leq p_{2}<\infty,\vspace{0.1cm}\\
\gamma_{1}<\dimension\left(p_{1}^{-1}-p_{2}^{-1}\right), & \text{if }q_{1}>p_{2},\vspace{0.1cm}\\
\gamma_{1}<\dimension\left(p_{1}^{-1}-p_{2}^{-1}\right), & \text{if }p_{2}=\infty\text{ and }q_{1}>1=\LowerExpo{p_{2}},\vspace{0.1cm}\\
\gamma_{1}\leq\dimension\left(p_{1}^{-1}-p_{2}^{-1}\right), & \text{if }p_{2}=\infty\text{ and }q_{1}\leq1=\LowerExpo{p_{2}}.
\end{cases}\label{eq:BesovHomIntoInhomLowFrequencyConditionNecessary}
\end{equation}
Finally, in case of $p_{1}=p_{2}$, we also have
\begin{equation}
\begin{cases}
\gamma_{1}\leq0, & \text{if }q_{1}\leq2,\\
\gamma_{1}<0, & \text{if }q_{1}>2.
\end{cases}\qedhere\label{eq:BesovHomIntoInhomLowFrequencyKhinchinConditionNecessary}
\end{equation}
\end{thm}

\begin{rem*}
As for the previous theorem, we precisely determine the case in which
our criteria are inconclusive. The condition $p_{1}\leq p_{2}$ and
condition~(\ref{eq:BesovHomIntoInhomHighFrequencyConditionSufficient})
appear in both the necessary and the sufficient conditions. Hence,
the only difference is between conditions (\ref{eq:BesovHomIntoInhomLowFrequencyConditionSufficient})
and (\ref{eq:BesovHomIntoInhomLowFrequencyConditionNecessary}). In
case of $p_{2}\in\left(0,2\right]$, we have $\LowerExpo{p_{2}}=p_{2}$,
so that both conditions coincide. They also coincide for $p_{2}=\infty$.
Thus, the \emph{only} case in which our criteria are inconclusive
is if
\begin{equation}
p_{1}\leq p_{2},\qquad p_{2}\in\left(2,\infty\right),\qquad\gamma_{1}=\dimension\left(p_{1}^{-1}-p_{2}^{-1}\right),\qquad\text{ and finally }\qquad p_{2}'<q_{1}\leq p_{2}.\label{eq:BesovHomIntoInhomInconclusiveRegion1}
\end{equation}
In this case, the necessary conditions are fulfilled, but the sufficient
conditions are not.

In case of $p_{1}=p_{2}=:p$, we can narrow this region down even
further, using condition~(\ref{eq:BesovHomIntoInhomLowFrequencyKhinchinConditionNecessary}).
In this case, the only case in which our criteria are inconclusive
is if
\begin{equation}
p\in\left(2,\infty\right),\qquad\gamma=0\qquad\text{ and finally }\qquad p'<q_{1}\leq2.\label{eq:BesovHomIntoInhomInconclusiveRegion2}
\end{equation}

Finally, note that our necessary conditions \emph{always} entail $\gamma_{1}\leq\dimension\left(p_{1}^{-1}-p_{2}^{-1}\right)$
and 
\[
\gamma_{2}\leq\gamma_{1}+\dimension\left(p_{2}^{-1}-p_{1}^{-1}\right)\leq\dimension\left(p_{1}^{-1}-p_{2}^{-1}\right)+\dimension\left(p_{2}^{-1}-p_{1}^{-1}\right)=0,
\]
so that the homogeneous Besov spaces can only ever embed into inhomogeneous
Besov spaces of \emph{non-positive} smoothness.

\medskip{}

We also remark that the preceding theorem is connected to the paper
\cite{BourdaudRealizationsOfHomogenousBesovSpaces}. In that paper
Bourdaud studies the existence of \textbf{translation commuting realizations}
of the homogeneous Besov spaces $\BesovHom pq{\gamma}\left(\R^{\dimension}\right)$.
More precisely, Bourdaud characterizes the conditions under which
there is a continuous linear map $\eta:\BesovHom pq{\gamma}\left(\R^{\dimension}\right)\to\Schwartz'\left(\R^{\dimension}\right)/\CalP_{\nu}$,
where $\CalP_{\nu}$ denotes the set of polynomials of degree \emph{strictly
less} than $\nu$ with $\nu\in\N_{0}$ fixed, such that $\eta\left(f+\CalP\right)+\CalP=f+\CalP$
and $\eta\left(L_{x}f+\CalP\right)=\eta\left(f+\CalP\right)$ for
all $f+\CalP\in\BesovHom pq{\gamma}\left(\R^{\dimension}\right)\subset\Schwartz'\left(\R^{\dimension}\right)/\CalP$
and $x\in\R^{\dimension}$. One can indeed show that the map $\iota$
constructed above is translation commuting. Thus, for the case $\nu=0$,
our result gives another proof of the fact that there exists a translation
commuting realization $\eta:\BesovHom pq{\gamma}\left(\R^{\dimension}\right)\to\Schwartz'\left(\R^{\dimension}\right)$
if $\gamma\leq d/p$ (for $q\leq1$) or if $\gamma<d/p$ (for $q>1$).
For more details, we refer to \cite[Remark 6.2.5]{VoigtlaenderPhDThesis}.
\end{rem*}
\begin{proof}[Proof of Theorem~\ref{thm:HomogeneousIntoInhomogeneousBesov}]
Define 
\[
\CalQ=\left(Q_{i}\right)_{i\in I}=\left(T_{i}Q_{i}'+b_{i}\right)_{i\in I}:=\dot{\CalB}=\left(\smash{\dot{B}_{k}}\right)_{k\in\Z}=\bigl(\left(\smash{2^{k}}\cdot\identity\right)Q\bigr)_{k\in\Z}
\]
and
\[
\CalP=\left(P_{j}\right)_{j\in J}=\left(S_{j}P_{j}'+c_{j}\right)_{j\in J}:=\CalB=\left(B_{n}\right)_{n\in\N_{0}}=\bigl(\left(2^{n}\cdot\identity\right)B_{n}'\bigr)_{n\in\N_{0}}\,,
\]
with $Q$ and $B_{n}'$ as in Lemmas~\ref{lem:HomogeneousBesovCoveringStructured}
and \ref{lem:InhomogeneousBesovCoveringAlmostStructured}, respectively.Furthermore,
set $w:=u^{\left(\gamma_{1}\right)}=\left(2^{k\gamma_{1}}\right)_{k\in\Z}$
and $v:=v^{\left(\gamma_{2}\right)}=\left(v^{n\gamma_{2}}\right)_{n\in\N_{0}}$.
In view of Lemmas~\ref{lem:InhomogeneousBesovCoveringAlmostStructured}
and \ref{lem:HomogeneousBesovCoveringStructured}, we see that $\CalQ$
and $\CalP$ are both almost structured coverings of $\CalO:=\R^{\dimension}\setminus\left\{ 0\right\} $
and of $\CalO':=\R^{\dimension}$, respectively. Furthermore, the
same lemmas also show that $w$ and $v$ are $\CalQ$-moderate and
$\CalP$-moderate, respectively, so that the standing assumptions
of Section~\ref{sec:SummaryOfEmbeddingResults} (i.e., Assumption~\ref{assu:GeneralSummaryAssumptions})
are satisfied.

Finally, thanks to Lemma~\ref{lem:HomogeneousInhomogeneousBesovSubordinateness},
we know that $\CalQ$ is almost subordinate to $\CalP$ (but \emph{not}
relatively $\CalP$-moderate). For the sufficient condition, we can
thus apply part~(\ref{enu:SummaryFineInCoarseSufficientNonModerate})
of Theorem~\ref{thm:SummaryFineIntoCoarse} to conclude that a map
$\iota$ with properties as in the current theorem exists if we have
$p_{1}\leq p_{2}$ and if $K_{\LowerExpo{p_{2}}}<\infty$ (in the
notation of Theorem~\ref{thm:SummaryFineIntoCoarse}). To see this,
note that the map $\iota$ given by Theorem~\ref{thm:SummaryFineIntoCoarse}
satisfies $\left\langle \iota f,\,\varphi\right\rangle _{\CalD'}=\left\langle f,\,\varphi\right\rangle _{\CalD'}$
for all $\varphi\in\TestFunctionSpace{\R^{\dimension}\setminus\left\{ 0\right\} }$
and $f\in\FourierDecompSp{\CalQ}{p_{1}}{\ell_{w}^{q_{1}}}=\BesovHom{p_{1}}{q_{1}}{\Fourier,\gamma_{1}}\left(\R^{\dimension}\right)$.
Furthermore, the map $\iota$ from Theorem~\ref{thm:SummaryFineIntoCoarse}
also fulfills $\iota f=f$ for all $f\in\TestFunctionSpace{\CalO}=\TestFunctionSpace{\R^{\dimension}\setminus\left\{ 0\right\} }$,
as desired.

To verify $K_{\LowerExpo{p_{2}}}<\infty$, we more generally investigate
the finiteness of
\begin{align*}
K\left(r\right):=K_{r} & =\left\Vert \left(v_{j}\cdot\left\Vert \left(\left|\det T_{i}\right|^{p_{1}^{-1}-p_{2}^{-1}}/w_{i}\right)_{i\in I_{j}}\right\Vert _{\ell^{r\cdot\left(q_{1}/r\right)'}}\right)_{j\in J}\right\Vert _{\ell^{q_{2}\cdot\left(q_{1}/q_{2}\right)'}}\\
 & =\left\Vert \left(2^{n\gamma_{2}}\cdot\left\Vert \left(2^{k\dimension\left(p_{1}^{-1}-p_{2}^{-1}\right)}/2^{k\gamma_{1}}\right)_{k\in\Z^{\left(n\right)}}\right\Vert _{\ell^{r\cdot\left(q_{1}/r\right)'}}\right)_{n\in\N_{0}}\right\Vert _{\ell^{q_{2}\cdot\left(q_{1}/q_{2}\right)'}}
\end{align*}
for general $r\in\left(0,\infty\right]$, where
\[
\Z^{\left(n\right)}:=\left\{ k\in\Z\with\smash{\dot{B}_{k}}\cap B_{n}\neq\emptyset\right\} 
\]
is defined as in equation~(\ref{eq:InhomIntoHomBesovIntersectionSet}).
Recall from Lemma~\ref{lem:HomogeneousInhomogeneousBesovSubordinateness}
that we have $\left\{ n\right\} \subset\Z^{\left(n\right)}\subset\left\{ n-3,\dots,n+3\right\} $
for $n\in\N$ and $\Z_{\leq1}\subset\Z^{\left(0\right)}\subset\Z_{\leq3}$.

Now, we observe that $K\left(r\right)$ is finite if and only if
\[
K_{0}\left(r\right):=\left\Vert \left(2^{k\dimension\left(p_{1}^{-1}-p_{2}^{-1}\right)}/2^{k\gamma_{1}}\right)_{k\in\Z^{\left(0\right)}}\right\Vert _{\ell^{r\cdot\left(q_{1}/r\right)'}}=\left\Vert \left(2^{k\left[-\gamma_{1}+\dimension\left(p_{1}^{-1}-p_{2}^{-1}\right)\right]}\right)_{k\in\Z^{\left(0\right)}}\right\Vert _{\ell^{r\cdot\left(q_{1}/r\right)'}}
\]
is finite and if furthermore
\begin{align*}
\widetilde{K}\left(r\right) & :=\left\Vert \left(2^{n\gamma_{2}}\cdot\left\Vert \left(2^{k\dimension\left(p_{1}^{-1}-p_{2}^{-1}\right)}/2^{k\gamma_{1}}\right)_{k\in\Z^{\left(n\right)}}\right\Vert _{\ell^{r\cdot\left(q_{1}/r\right)'}}\right)_{n\in\N}\right\Vert _{\ell^{q_{2}\cdot\left(q_{1}/q_{2}\right)'}}\\
\left({\scriptstyle \text{equation }\eqref{eq:BesovHomogeneousInhomogeneousIntersectionHighFrequency}}\right) & \asymp\left\Vert \left(2^{n\gamma_{2}}\cdot2^{n\dimension\left(p_{1}^{-1}-p_{2}^{-1}\right)}/2^{n\gamma_{1}}\right)_{n\in\N}\right\Vert _{\ell^{q_{2}\cdot\left(q_{1}/q_{2}\right)'}}=\left\Vert \left(2^{n\left[\gamma_{2}-\gamma_{1}+\dimension\left(p_{1}^{-1}-p_{2}^{-1}\right)\right]}\right)_{n\in\N}\right\Vert _{\ell^{q_{2}\cdot\left(q_{1}/q_{2}\right)'}}
\end{align*}
is finite. But equation~(\ref{eq:SpecialExponentFiniteness}) shows
that the exponent $q_{2}\cdot\left(q_{1}/q_{2}\right)'$ is finite
if and only if $q_{2}<q_{1}$. Hence, due to the exponential nature
of the sequence above, we conclude
\begin{equation}
\widetilde{K}\left(r\right)<\infty\qquad\Longleftrightarrow\qquad\begin{cases}
\gamma_{2}-\gamma_{1}+\dimension\left(p_{1}^{-1}-p_{2}^{-1}\right)\leq0, & \text{if }q_{1}\leq q_{2},\vspace{0.1cm}\\
\gamma_{2}-\gamma_{1}+\dimension\left(p_{1}^{-1}-p_{2}^{-1}\right)<0, & \text{if }q_{1}>q_{2}.
\end{cases}\label{eq:BesovHomogeneousIntoInhomogeneousHighFrequencyCondition}
\end{equation}
Note that the right-hand side is \emph{independent} of $r$.

Finally, recalling $\Z_{\leq1}\subset\Z^{\left(0\right)}\subset\Z_{\leq3}$
and noting that $r\cdot\left(q_{1}/r\right)'$ is finite if and only
if $r<q_{1}$, we also conclude
\begin{equation}
K_{0}\left(r\right)<\infty\qquad\Longleftrightarrow\qquad\begin{cases}
-\gamma_{1}+\dimension\left(p_{1}^{-1}-p_{2}^{-1}\right)\geq0, & \text{if }q_{1}\leq r,\vspace{0.1cm}\\
-\gamma_{1}+\dimension\left(p_{1}^{-1}-p_{2}^{-1}\right)>0, & \text{if }q_{1}>r.
\end{cases}\label{eq:BesovHomogeneousIntoInhomogeneousLowFrequencyCondition}
\end{equation}
If we now recall from above that a map $\iota$ with the properties
claimed in the current theorem exists if we have $p_{1}\leq p_{2}$
and if $K_{\LowerExpo{p_{2}}}<\infty$, we see that the sufficient
conditions stated in the current theorem are indeed sufficient for
the existence of $\iota$.

\medskip{}

Now, we consider the necessary conditions. To this end, assume that
a map $\iota$ as in the statement of the current theorem exists.
Since $\iota$ satisfies $\iota f=f$ for all $f\in\TestFunctionSpace{\R^{\dimension}\setminus\left\{ 0\right\} }=\TestFunctionSpace{\CalO}$,
we see that the map $\theta$ from part~(\ref{enu:SummaryFineInCoarseNecessaryNonModerate})
of Theorem~\ref{thm:SummaryFineIntoCoarse} is well-defined and bounded.
Hence, that part of Theorem~\ref{thm:SummaryFineIntoCoarse} shows
that we have $p_{1}\leq p_{2}$ and $K_{s}<\infty$, with $s=p_{2}$
for $p_{2}\in\left(0,\infty\right)$ and with $s=1=\LowerExpo{p_{2}}$
for $p_{2}=\infty$. In view of conditions~(\ref{eq:BesovHomogeneousIntoInhomogeneousHighFrequencyCondition})
and (\ref{eq:BesovHomogeneousIntoInhomogeneousLowFrequencyCondition}),
we thus see that conditions~(\ref{eq:BesovHomIntoInhomHighFrequencyConditionSufficient})
and (\ref{eq:BesovHomIntoInhomLowFrequencyConditionNecessary}) from
the statement of the current theorem are indeed fulfilled.

Finally, if we have $p_{1}=p_{2}$, we can apply part~(\ref{enu:SummaryFineInCoarseNecessaryKhinchin})
of Theorem~\ref{thm:SummaryFineIntoCoarse}, which yields $K_{1}<\infty$
in case of $p=\infty$, and $K_{\min\left\{ 2,p\right\} }<\infty$
if $p<\infty$. Just as above, this easily implies that condition~(\ref{eq:BesovHomIntoInhomLowFrequencyKhinchinConditionNecessary})
is indeed fulfilled.
\end{proof}
\begin{rem}
One can show—similar to Lemma~\ref{lem:FourierInhomogeneousBesovAsTemperedDistributions}—that
the map
\begin{equation}
\BesovHom pq{\gamma}\left(\R^{\dimension}\right)\subset\Schwartz'\left(\R^{\dimension}\right)/\mathcal{P}\to\FourierDecompSp{\smash{\dot{\CalB}}}p{\ell_{u^{\left(\gamma\right)}}^{q}}=\BesovHom pq{\Fourier,\gamma}\left(\R^{\dimension}\right),f+\mathcal{P}\mapsto\widehat{f}|_{\TestFunctionSpace{\R^{\dimension}\setminus\left\{ 0\right\} }}\label{eq:HomogeneousBesovFourierSpaceIsomorphism}
\end{equation}
is an isomorphism of (quasi)-Banach spaces; see \cite[Lemma 6.2.2]{VoigtlaenderPhDThesis}.
This justifies the name ``homogeneous (Fourier-side) Besov space''
of the space $\BesovHom pq{\Fourier,\gamma}\left(\R^{\dimension}\right)$.
Furthermore, using this isomorphism, it is easy to see that the embedding
results from Theorems \ref{thm:InhomogeneousIntoHomogeneousBesov}
and \ref{thm:HomogeneousIntoInhomogeneousBesov} also yield corresponding
embedding results for the ``usual, space-side'' homogeneous and
inhomogeneous Besov spaces.

The main step in showing that (\ref{eq:HomogeneousBesovFourierSpaceIsomorphism})
indeed yields an isomorphism is to show that each distribution $f\in\FourierDecompSp{\smash{\dot{\CalB}}}p{\ell_{u^{\left(\gamma\right)}}^{q}}\subset\DistributionSpace{\R^{\dimension}\setminus\left\{ 0\right\} }$
can be extended to a \emph{tempered} distribution $\tilde{f}\in\Schwartz'\left(\R^{\dimension}\right)$.
To prove this, however, we cannot use Theorem~\ref{thm:DecompositionSpacesAsTemperedDistributions},
since this theorem relies on the quantity 
\[
c_{i}:=\inf_{\xi\in Q_{i}^{\ast}}\left(1+\left|\xi\right|\right)
\]
being large for most $i\in I$. But for the homogeneous Besov covering,
we have $c_{i}\asymp1$ for all $i\in\Z_{\leq0}$. The actual proof
given in \cite[Lemma 6.2.2]{VoigtlaenderPhDThesis} is quite technical
and does not yield a further illustration of our embedding results.
Hence, we omit the proof.
\end{rem}

\let\section\origsection\addtocontents{toc}{\SkipTocEntry}

\section*{Concluding remarks}

In this paper, we developed a comprehensive and easy to use framework
for deciding the existence of embeddings between decomposition spaces.
This framework essentially allows to reduce the study of embeddings
between the Fourier analytic decomposition spaces to \emph{purely
combinatorial conditions}. To check whether these conditions are satisfied,
\emph{no knowledge of Fourier analysis is required} anymore. Instead,
all one needs to understand is the relation between the two coverings
under consideration and the associated discrete sequence spaces.

The applications given in Section~\ref{sec:Applications} underline
this claim: Once we understood the relation between the different
``$\alpha$-modulation coverings'', or between the homogeneous and
the inhomogeneous Besov coverings, the application of the criteria
produced by our embedding framework was straightforward.

Note, however, that we were unable to \emph{completely} characterize
the existence of embeddings between homogeneous and inhomogeneous
Besov spaces. This shows that for highly ``incompatible'' coverings,
the developed theory is still incomplete—although it is more comprehensive
than any previously known results. The desire to close this gap provides
a fascinating topic for further research.\addtocontents{toc}{\SkipTocEntry}

\section*{Acknowledgments}

I warmly thank Hartmut Führ, Hans Feichtinger, and David Rottensteiner
for fruitful discussions and inspiring comments regarding the topics
in this paper. Jordy van Velthoven heavily contributed to this paper
by pushing me to finally finish revising it.

I especially thank the anonymous referee for several valuable suggestions
that led to improvements in the presentation of the paper, and also
for suggesting important references that were missing in the earlier
version of the paper.

This research was funded by the Excellence Initiative of the German
federal and state governments, and by the German Research Foundation
(DFG), under the contract FU 402/5-1. While this research was mainly
carried out at RWTH Aachen University, I am currently employed at
the KU Eichstätt–Ingolstadt. Finally, this manuscript was written
to a large extent at TU Berlin. I therefore gratefully acknowledge
support from the European Commission through DEDALE (contract no.~665044)
within the H2020 Framework Program.\vspace{-0.1cm}

\bibliographystyle{plain}
\bibliography{felixbib}

\end{document}